\documentclass[12pt,leqno]{article}
\usepackage{fancyhdr}
\usepackage{amsmath}
\usepackage{amssymb}

\newtheorem{definition}{Definition}[section] \rm
\newtheorem{lemma}{Lemma}[section]
\newtheorem{theorem}{Theorem}[section]
\newtheorem{proposition}{Proposition}[section]
\newtheorem{corollary}{Corollary}[section]
\newenvironment{proof}{{\noindent\it Proof.}}\rm

\begin{document}




\textwidth 16.59cm
\textheight 21.94cm

\renewcommand\theequation{\thesection.\arabic{equation}}


\def\oo{\omega}
\def\C{\mathbb{C}}
\def\R{\mathbb{R}}
\def\Z{\mathbb{Z}}
\def\bfe{{\mathbf{E}}}
\def\[{\left[}
\def\]{\right]}
\def\({\left(}
\def\){\right)}
\def\wt{\widetilde}
\def\a{\alpha}
\def\b{\beta}
\def\de{\delta}
\def\lbb{\lambda}
\def\bfe{\mathbf{E}}
\def\bfa{\mathbf{A}}
\def\bfb{\mathbf{B}}
\def\bff{\mathbf{F}}
\def\bfr{\mathbf{R}}
\def\bfm{\mathbf{M}}
\def\bft{\mathbf{T}}
\def\bfc{\mathbf{C}}
\def\bfj{\mathbf{J}}
\def\bfi{\mathbf{I}}
\def\dd{\displaystyle}
\def\g{\gamma}
\def\G{\Gamma}
\def\9{\infty}
\def\pp{\partial}
\def\vare{\varepsilon}
\def\ms{\medskip}
\def\F{\mathbf{F}}
\def\M{\mathbf{M}}
\def\H{\mathbf{H}}
\def\J{\mathbf{J}}

\newif\iftitlepage \titlepagetrue
\newtoks\titlepagehead \titlepagehead={\hfil}
\newtoks\titlepagefoot \titlepagefoot={\hfil}
\newtoks\runningauthor \runningauthor={\hfil}
\newtoks\runningtitle \runningtitle={\hfil}
\newtoks\evenpagehead \newtoks\oddpagehead
\evenpagehead={\hfil\the\runningauthor\hfil}
\oddpagehead={\hfil\the\runningtitle\hfil}
\newtoks\evenpagefoot \evenpagefoot={\hfil}
\newtoks\oddpagefoot \oddpagefoot={\hfil}
\thispagestyle{empty}

\centerline{\footnotesize{``Alexandru Ioan Cuza'' University of Ia\c{s}i}}

\vspace*{1cm}

\centerline{\bf Mihai ANASTASIEI}

\vspace*{6cm}

\centerline{\Large{\bf SELECTED PAPERS}}

\vspace*{8cm}

\centerline{IA\c{S}I, 2012}

\newpage

\centerline{\bf \Large PREFACE}

\medskip

Professor Doctor Mihai Anastasiei is a distinguished member of the
Faculty of Mathematics  from the University Alexandru Ioan Cuza of
Ia\c si. Born in 1946, in Boto\c sani county, he passed the entrance
admission for the Faculty of Mathematics and Mechanics, University Alexandru Ioan Cuza
 of
Ia\c si, in 1964. He graduated in 1969, by receiving a Diploma of Merit in
Mathematics with the specialization in Geometry. After completing
the post-graduate studies, between 1972 and 1976, he received, in 1977, the Ph. D
title with the Thesis {\it Remarkable geometric structures on infinite
dimensional manifolds}, under the supervision of Professor Doctor Radu
Miron.

After graduating the University, in 1969, he was appointed as
an assistant and some  years latter was promoted as a lecturer in the
Department of Geometry, Faculty of Mathematics. In 1990 he became an
associate professor and in 1995 was promoted to a full professor.

Professor Mihai Anastasiei has devoted his lifelong academic
career to fulfill his teaching, research and administrative duties as well. As
a teacher, he guided, with love for mathematics and especially
geometry, many generations of students. The courses he taught, including analytic
geometry, differential geometry of  curves and surfaces, differential
geometry of manifolds, Riemannian geometry, theory of $G$-structures,
Finsler and Lagrange geometry, where always rigorous, well prepared,
and characterized by a personal style.

During his career, Professor Mihai Anastasiei has served the academic
community as a Head of the Department of Geometry, Chancellor and vice-Dean  of the
Faculty of Mathematics, and in the last four years as the Head of an interdisciplinary
(Biology, Chemistry, Geography, Geology, Informatics, Mathematics, Physics), Department of
Science at the University Alexandru Ioan Cuza of Ia\c si.

The research activity of Professor Mihai Anastasiei is very
substantial and had a good impact in the scientific community. He
coauthored three monographs, edited nine volumes, and
published more than ninety scientific papers in prestigious
journals. For his research activity he was awarded by the Romanian
Academy, in 1987, the Prize ``Gheorghe \c Ti\c teica''. He
participated to more than twenty international conferences
(ICM 1990, Kyoto; SUA, Canada, Greece, Italy, Hungary). In his
research activity he had fruitful collaborations with his teacher, Prof. Dr. Radu Miron,
now a member of the Romanian Academy, and with other Romanian professors as well as  with researchers from Japan, Canada,
Hungary, the US. He had the status of visiting professor at prestigious
universities around the world: Italy, Bari, 1986; Japan, Tsukuba,
1990; Canada, Edmonton, 1993; Germany, Munich, 1998; Japan, Sapporo,
1998; Belgium, Leuven, 2001; Poland, Wroclaw, 2003.

The present volume contains a selection of thirty two of the articles published by
Professor Mihai Anastasiei, in prestigious journals or Proceedings that made
a certain impact in the scientific world. These articles cover the
following research topics: geometry of conformal and spin structures
on Hilbert manifolds, geometry of tangent bundle of a Finsler or a
Lagrange manifold, applications of techniques from Finsler geometry to
Mechanics and Physics, geometry of total space of a vector bundle, applications of Lie algebroids to Mechanics.

Ia\c si, May 15, 2012

\hfill Prof. Dr. Radu Miron, Member of the Romanian Academy.

\
\hfill Asso. Prof. Dr. Ioan Bucataru

\newpage

\centerline {\bf \Large AUTHOR'S NOTES}

\medskip

I began my undergraduate studies in 1964, a time when the political atmosphere in Romania had begun to emerge from a strong freezing caused by the events that followed the Second World War.

The trend continued till about 1972. In the meantime I graduated and I was appointed an Assistant at the Alexandru Ioan Cuza University. In 1972 I became a Ph D  student under the supervision of Prof. Dr. Radu Miron, later a member of the Romanian Academy. Then I began to systematically frequent the library of the Mathematical Seminar Alexandru Myller, a special library that contains mainly mathematical books. Its premises also offered me the possibility of meeting all the other mathematicians in Ia\c si.

In spite of the iron curtain, our professors succeeded in (re)establishing contacts with Western mathematicians and some of them were able to travel to Europe to attend various professional meetings. Consequently, the library of the Mathematical Seminar Alexandu Myller begun receiving as a gift or in exchange with our mathematical journal ``Analele \c Stiintifice ale Universit\u a\c tii Alexandru Ioan Cuza, Ia\c si, s I a Matematica'', important books published in Europe, the US, and Japan. Our Geometry Chair benefited of the book on foundations of differential geometry authored by S. Kobayashi and K. Nomizu. This was the main reference in my doctoral studies. Around 1960, some important studies signed by J. Eells, R.S. Palais and others on infinite dimensional manifolds with applications to the space functions had appeared. They were mainly concerned with the topology of such manifolds but the differential geometry of Hilbert or Banach manifolds looked also interesting. A group of mathematicians from Ia\c si, namely: M. Craioveanu, L. Maxim-Raileanu, A. Bejancu, N. Papaghiuc, Al. Neagu, I. Burdujan started studying the topology and the differential geometry of Hilbert and Banach manifolds. I joined this group and I benefited a lot from the references they accessed, from the discussions and from reports we systematically presented in regular meetings at the Geometry Chair. My general task was to study various geometrical structures on Banach or Hilbert manifolds. Finally I concentrated on two such structures: conformal structures and spin structures. Concerning these structures I elaborated the papers [1-6, {\it Papers}]\footnote[1]{The numbers in the brackets refer to {\it List of Publications}} based on which I wrote my Ph D Thesis and delivered it in June 1977. I learned about spin structures in the finite dimensional setting from Iulian Popovici, a senior researcher at the Mathematical Institute of the Romanian Academy in Bucharest. He was a very kind person, well knowledgeable in differential geometry and in connection with several Western mathematicians. Later on we published three joint papers [8, 9 and 12, {\it Papers}]. Unfortunately, he died suddenly in 1981.

In the autumn of 1978 I started learning Finsler geometry as it was presented at that time by the French and the Japanese schools. Thus I joined a group led by Prof. Dr. Radu Miron who organized the First National Seminar on Finsler Geometry at the University of Bra\c sov in February 1980. There I met some Hungarian geometers (L. Tam\'{a}ssy and his group: Z. Szab\'{o}, J. Szilasi, L. Kozma, S. Bacso, Tran Quoc Binh), and I remained in contact with them to this day. This Seminar has been held every 2 years to this day, with some changes in its topics and the location it was held in.

Studying various papers on Finsler spaces I was intrigued by the objects that were depending on points on base manifold as well as on the directions in those points but were behaving as objects on the base manifold. I call them Finsler geometric objects and I included them in a general theory of geometrical objects by constructing some fiber bundles whose sections were Finsler geometric objects. Their existence was therefore a consequence of some general results about the existence of sections in fiber bundles [11, 13-16, {\it Papers}]. In particular, I was interested in the existence of Finsler metrics. Here I compared the notion of Finsler structure on Banach manifolds introduced by R.S. Palais and that used in the finite dimensional setting. I came to the conclusion that the so-called positive definiteness condition was not adequate to Banach setting. It would imply the reduction of the Finsler structure to a Riemannian one.

In our Geometry Chair seminars I listened and then carefully studied the geometry of the total space of a vector bundle that Prof. Dr. Radu Miron developed using some techniques from Finsler geometry. His connection with Prof. Dr. Makoto Matsumoto, established in 1976 proved very useful. Two members of the Japanese school: M. Hashiguchi and Y. Ichijio were our guests in Ia\c si and attended our seminars. At these meetings a Romanian –Japanese Colloquium was planned and was held in 1984 at the University Alexandru Ioan Cuza of Ia\c si and at the University of Bra\c sov. Sixteen Japanese (two accompanying persons) were our guests for ten days (one day field trip to Suceava County). Back then, I had the great opportunity of meeting some very nice persons. I remained in touch with all of them for many years. I met some of them again and again in Romania, Japan, the US, Canada, and Hungary. In the following years, some 10 members of the Japanese school completed doctoral studies in Ia\c si under the guidance of Prof. Dr. Radu Miron. One member of this school (M. Kitayama) studied under my supervision and was granted Doctor's degree in 2001. My Japanese colleagues grew older and older. Some of them passed away. I keep good memories of them all.

Being so far much involved in the study of Finsler geometry I understood there were three essential modern ways to present it: using vector bundles, using principal bundles or using the total space of the tangent bundle. All these include the old–fashioned Finsler geometry as P.Finsler, Y. Thomas, \'{E}. Cartan developed and for which these appear as a model or a paradigm, [22, {\it Papers}]. I reported this at the University of Bari (1986) where I was invited by Biaggio Casciaro for five weeks. It was my first visit to a Western country (three days  in Venezia and in Rome). Later on, I had attended various meetings or fulfilled stages of visiting researcher in Belgium, Bulgaria, Canada, Germany, Greece, Hungary, Japan,  Poland, Yugoslavia, the US.

My contributions to Finsler geometry were at that time (1986) numerous and valuable. Consequently I was deemed an appropriate candidate as to co-author when Prof. Dr. Radu Miron decided to invite someone to co-write a book (in Romanian) on Finsler geometry. In the meantime he also suggested some applications in the Theory of Relativity, and as aforementioned, the main construction from Finsler geometry was extended to the total space of a vector bundle. This fact justified the book's title, and the chapter by S. Ikeda in which the arguments from Physics were provided. The book [1, {\it Monograph}] appeared in 1987 at the Publishing House of the Romanian Academy. It influenced very much the young generation and many young (at that time) Romanian geometers embarked in the study of Finsler geometry and its generalizations. Some obtained very nice results. To me it brought much intellectual satisfaction since later on based on it I received the ``Gh. Titeica'' prize awarded by the Romanian Academy.

In the next few years I was focused on applications in Mechanics and Physics and on techniques from the Finsler geometry and from the geometry of the total space of a vector bundle. Thus I wrote a kind of a general Einstein equation for vector bundles in [21, {\it Papers}]. I showed that the case of one dimensional fibers provided a type of Kaluza –Klein theory. Then I studied the geometry of Lagrangians that are explicitly time dependent. These are living on the manifold  $\mathbf{R}\times TM $regarded  as the total space of a convenient vector bundle. I worked on this subject together with H. Kawaguchi and I presented the main results at the International Congress of Mathematician (ICM) that was held on August 1990 in Kyoto. An overview of all these results can be found in [36, {\it Papers}]. Based on these results a Ph D Thesis was delivered in Ia\c si and generalizations of them were recently presented by V. Balan and M. Neagu in a book published by the John Wiley \& Sons (2011).

Various applications in Mechanics and Physics developed by Prof. Dr. Radu Miron and myself were included in our joint book [2, {\it Monographs}] published in 1994. In the meantime (1990-1992) S.S. Chern, D. Bao and Z. Shen brought a fresh air in Finsler geometry by approaching some old global problems using the so-called Chern connection. They contacted us and we exchanged some ideas so that in 1995 when they obtained funds for a a Joint Summer Research Conference on Finsler Geometry that held in Seattle in July 16-20,1995, Prof. Dr. Radu Miron and I were invited as speakers. The Proceedings were published the following year in Contemporary Mathematics. Upon my return I found out that the Chern connection coincides in fact with a well- known connection bearing the name of H. Rund. I wrote the details and the paper was published in the same volume of Proceedings see [40-42, {\it Papers}]. I continued to stay in contact with D. Bao and especially with Z. Shen, so I met again the latter in Edmonton and I invited him to Romania to attend a Seminar on Finsler Geometry that was held in Bac\u au. They gave me as a gift their famous textbook on Finsler geometry published by Springer in 2000.

At the Conference in Seattle I stressed upon the importance of the structures defined by a generalized Finsler or Lagrange metric. Back home, I continued to study the generalized Lagrange spaces and I published the papers [44-46, 49, 50, {\it Papers}]. In meantime, R.G. Beil, an American Physicist who was my guest in 1997, used in a kind of unified theory of fields a generalized Lagrange metric that appeared as a deformation of a Finsler metric or a Riemannian metric. I studied such deformations in the most general case in co-operation with H. Shimada [50, 57, {\it Papers}] and by myself  some particular cases [52, 54, 55, {\it Papers}].

Learning some symplectic geometry from I. Vaisman I connected it with the geometry of Lagrangians ([59-61, {\it Papers}].

At the time when Prof. Dr. L. Tam\'{a}ssy turned 80 I was invited to dedicate him an article. Thus I approached the question of metrizability he studied in different settings, in my case the metrizability of a linear connection in a vector bundle. This problem suggested me new ones and so I published [63, 66-71, {\it Papers}].

Beside my purely theoretical studies, I always kept an eye open to applications in Mechanics, Physics and Biology.  Thus I published in a team or by myself the papers [62, 65, 72, 74, 77]. The most extensive and interesting applications of Finsler and Lagrange geometry as well as of techniques suggested by these geometries were worked out by  Dr. Sergiu  Vacaru. He has done a lot for promoting the Finsler and Lagrange geometries in the world of Physicists. We published two joint papers [83, 86, {\it Papers}].

My studies from geometry of vector bundles led me to Lie algebroids. I studied these structures keeping in mind possible applications to Mechanics and Systems Theory. The results are in the papers [73, 74, 76-78, 84, {\it Papers}]. Then I reverted to the setting that I had worked in my youth, namely I studied the Lie algebroids in the category of Banach vector bundles in the paper [87, {\it Papers}]. So a circle was closed.

I wrote in cooperation or by myself shorter or larger overviews upon the topics that I studied, [40, 53, 56, 64, {\it Papers}]. I was an Editor or co-Editor for nine volumes, [1-9, {\it Edited Volumes}]. I gave lectures on many subjects from Geometry and I published five textbook, [1-5, {\it Textbooks}].

\medskip

\noindent {\bf Acknowledgements.} The preparation of this volume was partially supported by a grant of the Romanian National Authority for Scientific Research, CNCS-UEFISCDI project number PN -II- ID-PCE-2011-3-0256. The author is grateful to Professors Radu Miron and Ioan Bucataru for their Preface and to Professor  Marian Munteanu who kindly went throw the text and pointed out some misprints.

{\footnotesize{\it May, 2012}}

\begin{flushright}
Prof. Dr. Mihai Anastasiei
\end{flushright}

\newpage

\centerline{\bf\Large LIST OF PUBLICATIONS}

\medskip

\centerline{\bf PROF. DR. Mihai ANASTASIEI}
\vskip .25cm

\begin{itemize}
\item[{\bf I.}] {\bf Monographs}
\begin{enumerate}
\item Miron, Radu; Anastasiei, Mihai, {\it Fibrate
      vectoriale. Spa\c tii Lagrange. Aplica\c tii în teoria relativit\u a\c
      tii}. (Romanian) [Vector bundles. Lagrange spaces. Applications to the
      theory of relativity] With a chapter in English by Satoshi Ikeda. Editura
      Academiei Republicii Socialiste România, Bucharest, 1987. 216 p.
 \item Miron, Radu; Anastasiei, Mihai, {\it The geometry of
      Lagrange spaces: theory and applications}. Fundamental Theories of Physics,
      59. Kluwer Academic Publishers Group, Dordrecht, 1994. xiv+285 pp.
\item Miron, Radu; Anastasiei, Mihai, {\it Vector bundles
      and Lagrange spaces with applications to relativity}. With a chapter by
      Satoshi Ikeda. Translated from the 1987 Romanian original. Balkan Society
      of Geometers, Monographs and Textbooks, no. 1. Geometry Balkan Press,
      Bucharest, 1997. iv + 219 pp.
      \end{enumerate}

\item[{\bf II.}] {\bf Textbooks}

\begin{enumerate}
\item Anastasiei, Mihai, {\it Metodica predarii matematicii}.(Romanian). \\ Methods for Teaching Mathematics. Univ.
``Al. I. Cuza'', Ia\c si, 1985, 138 p.
\item Anastasiei, Mihai et al., {\it Metodica Predarii Matematicii}
    (Romanian). Methods for teaching Mathematics, Vol. III.
    Editura ``Lumina'', Chi\c{s}in\u au, 1996, 300 p.
\item Anastasiei, M., {\it Geometrie: Curbe \c si suprafete} In
      Romanian. [Geometry: Curves and Surfaces] Editura CERMI,
      Iasi, Romania, 2003
\item Anastasiei M., Crasmareanu M., {\it Lec\c tii de
    geometrie (Curbe \c si suprafe\c te)}. In Romanian [Lectures
    on geometry (Curves and surfaces)] Editura Tehnopress, Ia\c si,Romania,
    2005.
\item Anastasiei, M., {\it Geometrie superioar\u a \c si aplica\c tii.} (pentru sec\c tia ID)(in Romanian). Editura Universit\u a\c tii ``Alexandru Ioan Cuza'' Ia\c si, Ia\c si, 2005.
\end{enumerate}
\item[{\bf III.}]{\bf Edited Volumes}
\begin{enumerate}

\item {\it Proceedings of the 25th National Conference
      on Geometry and Topology,} September 18-23, 1995, Iasi,
      Romania. Edited by V. Cruceanu and M. Anastasiei. An.
      Stint.Univ. Al.I.Cuza Iasi, 42, Suppl. 264 p. 1996.
 \item {\it Lagrange and Finsler geometry. Applications to
      physics and biology}. Edited by P. L. Antonelli and R. Miron in cooperation
      with M. Anastasiei and Gh. Zet. Fundamental Theories of Physics, no.76.
      Kluwer Academic Publishers Group, Dordrecht, 1996. x + 279 pp.
\item Anastasiei, M., {\it The Mathematician RADU MIRON. His Work
      and Life}. (Editor). Geometry Balkan Press, Bucharest, 1998, 170 p.
 \item {\it Proceedings of the 10th National Conference
      on Finsler, Lagrange and Hamilton Geometry}. Selected papers from the
      conference held at the University of Craiova, Craiova, February 1998.
      Edited by R. Miron and M. Anástásiei. Algebras Groups Geom. 16 (1999), no.
      1. Hadronic Press, Inc., Palm Harbor, FL, 1999. pp. i--vi and 1--126.
\item {\it Proceedings of the 11th National Conference
      on Finsler, Lagrange and Hamilton Geometry}. Held in Craiova, February
      2000. Edited by Radu Miron, Mihai Anastasiei and Victor Bl\u anu\c t\u a.
      Algebras, Groups and Geometries. 17 (2000), no. 3. Hadronic Press, Inc., Palm Harbor,
      FL, 2000. pp. i--iv and 253--382.
\item {\it Finsler and Lagrange Geometries}.Proceedings of
   a Conference held on August 26-31, 2001, in Ia\c si, Romania.
   Edited by M. Anastasiei and P.L. Antonelli. Kluwer Academic
   Publishers, Dordrecht, 2003.
\item {\it The Mathematician Radu Miron. His work and life at
    75 anniversary.}Edited by M. Anastasiei. ``Al. I. Cuza''
    University Press, Ia\c si, Romania, 2003.

\item Cruceanu Vasile, {\it Selected papers}. Edited by M. Anastasiei. Editura PIM, Ia\c si, 2006.

\item {\it Proceedings of the Conference on Differential Geometry: Lagrange and hamilton spaces, september 3-8, 2007, Ia\c si, Romania. dedicated to Academician Radu Miron at his 80th anniversary}. Edited by M. Anastasiei as a supplement of T. LIII of An. \c Stiin\c t. Univ. ``Al. I. Cuza'' Ia\c si,s I a, Matematic\u a, 2007, 355p.

\end{enumerate}

\item[{\bf IV.}]{\bf Papers}

\begin{enumerate}
\item Anastasiei, Mihai, {\it Conformal structures on Banach vector bundles}. An. \c Sti. Univ. ``Al. I. Cuza'' Ia\c si
      Sec\c t. I a Mat. (N.S.) 20 (1974), no. 2, 351--358.
\item Anastasiei, Mihai, {\it Affine trensformations on Banach manifolds}. Lincei - Rend. Sc. mat. e nat. Vol. LX, marzo 1976, 4p.
\item Anastasiei, Mihai, {\it Structures spinorielles sur les
      vari\'{e}t\'{e}s hilbertiennes}. C. R. Acad. Sci. Paris S\'{e}r. A-B 284 (1977), no.
      16, A943--A946.
\item Anastasiei, Mihai, {\it Riemannian $P$-structures on
      vector bundle}. Papers from the National Colloquium on Geometry and
      Topology (Univ. Timi\c soara, Timi\c soara, 1977), pp. 55--62, Univ.
      Timi\c soara, Timi\c soara, 1977.
\item Anastasiei, Mihai, {\it Spin structures on Hilbert manifolds}. An. \c Stiin\c t. Univ. ``Al. I. Cuza'' Ia\c si Sec\c t. I a Mat.
      (N.S.) 24 (1978), no. 2, 367--373.
\item Anastasiei, Mihai, {\it Connexions sur les fibr\'{e}s
      spinoriels}. (French) An. \c Stiin\c t. Univ. ``Al. I. Cuza'' Ia\c si Sec\c
      t. I a Mat. (N.S.) 25 (1979), no. 2, 337--344.
\item Anastasiei, Mihai; Nicolescu, Liviu, {\it Directions
      caractéristiques associées à une paire de connexions lin\'{e}aires sur une
      variété de Banach}. (French) An. \c Stiin\c t. Univ. "Al. I. Cuza" Ia\c si
      Sec\c t. I a Mat. (N.S.) 25 (1979), no. 2, 349--354.
\item Anastasiei, Mihai; Popovici, Iulian, {\it Constant
      linear connections on Banach manifolds}. Rev. Roumaine Math. Pures Appl. 25
      (1980), no. 1, 3--11.
\item Popovici, Iulian; Anastasiei, Mihai, {\it Sur les
      bases de la géométrie finslerienne}. (French) C. R. Acad. Sci. Paris Sér.
      A-B 290 (1980), no. 17, A807--A810.
item Anastasiei, Mihai, {\it Generalized affine
      connections on Banach manifolds}. An. \c Stiin\c t. Univ. "Al. I. Cuza"
      Ia\c si Sec\c t. I a Mat. (N.S.) 26 (1980), no. 2, 381--388.
\item Anastasiei, Mihai, {\it Finsler geometric objects and
      their Lie derivative}. The Proceedings of the National Seminar on Finsler
      Spaces (Bra\c sov, 1980), pp. 11--25, Univ. Timi\c soara, Timi\c soara,
      1981.
\item Anastasiei, Mihai; Popovici, Iulian, {\it An
      intrinsic characterization of Finsler connections}. The Proceedings of the
      National Seminar on Finsler Spaces (Bra\c sov, 1980), pp. 27--39, Univ.
      Timi\c soara, Timi\c soara, 1981.
\item Anastasiei, Mihai, {\it Some tensorial Finsler
      structures on the tangent bundle}. An. \c Stiin\c t. Univ. ``Al. I. Cuza''
      Ia\c si Sec\c t. I a Mat. (N.S.) 27 (1981),suppl., 9--16.
\item Miron, Radu; Anastasiei, Mihai, {\it On the notion of
      Finsler geometric object}. Mem. Sec\c t. \c Stiin\c t. Acad. Repub. Soc.
      România Ser. IV 4 (1981), no. 1, 25--31 .
\item Anastasiei, Mihai, {\it Some existence theorems in
      Finsler geometry}. The Proceedings of the National Seminar on Finsler
      Spaces (Bra\c sov, 1982), 25--33, Soc. \c Stiin\c te Mat. R.S. România,
      Bucharest, 1982.
\item Anastasiei, Mihai, {\it Some existence theorems in
      Finsler geometry}. An. \c Stiin\c t. Univ. Al. I. Cuza Ia\c si Sec\c t. I a
      Mat. 29 (1983), no. 1, 79--84.
\item Anastasiei, Mihai, {\it On Miron's theory of
      subspaces in Finsler spaces}. An. \c Stiin\c t. Univ. Al. I. Cuza Ia\c si
      Sec\c t. I a Mat. 30 (1984), no. 4, 11--14.
\item Anastasiei, Mihai, {\it Finsler geometry on principal
      bundles}. Proceedings of the national conference on geometry and topology\\
      (Piatra Neam\c t, 1983), 70--77, Univ. Al. I. Cuza, Ia\c si, 1984.
\item Anastasiei, Mihai, {\it On metrical Finsler
      connections}. Proceedings of the third national seminar on Finsler spaces
      (Bra\c sov, 1984), 19--28, Soc. \c Stiin\c te Mat. R.S. România,
      Bucharest, 1984.
\item Anastasiei, Mihai, {\it On subspaces in Finsler
      spaces}. Proceedings of the Romanian-Japanese colloquium on Finsler
      geometry (Ia\c si -Bra\c sov, 1984), 17--30, Univ. Bra\c sov, Bra\c sov,
      1984.
\item Anastasiei, Mihai, {\it Vector bundles. Einstein
      equations}. An. \c Stiin\c t. Univ. Al. I. Cuza Ia\c si Sec\c t. I a Mat.
      32 (1986), no. 3, 17--24.
\item Anastasiei, Mihai, {\it Models of Finsler and
      Lagrange geometry}. The Proceedings of the Fourth National Seminar on
      Finsler and Lagrange Spaces (Bra\c sov, 1986), 43--56, Soc. \c Stiin\c te
      Mat. R.S. România, Bucharest, 1986.
\item Miron, Radu; Anastasiei, Mihai, {\it Invariant theory of
      subspaces in generalized Lagrange spaces}. Progr. Math. (Varanasi) 21
      (1987), no. 2, 59--74.
\item Anastasiei, Mihai, {\it Some models in geometry of
      Hamilton spaces}. The XVIIIth National Conference on Geometry and Topology
      (Oradea, 1987), 11--14, Preprint, 88-2, Univ. ``Babe\c s-Bolyai'',
      Cluj-Napoca, 1988.
\item Anastasiei, Mihai, {\it Nonlinear connections in
      Hamilton spaces}. The Proceedings of the Fifth National Seminar of Finsler
      and Lagrange Spaces (Bra\c sov, 1988), 47--57, Soc. \c Stiin\c te Mat.
      R.S. România, Bucharest, 1989.
\item Miron, Radu; Kirkovits, Magdalen Sz.;
      Anastasiei, Mihai, {\it A geometrical model for variational problems of multiple
      integrals}. Differential geometry and its applications (Dubrovnik, 1988),
      209--216, Univ. Novi Sad, Novi Sad, 1989.
\item Miron, Radu; Ianu\c s, Stere; Anastasiei, Mihai,
      {\it The geometry of the dual of a vector bundle}. Publ. Inst. Math. (Beograd)
      (N.S.) 46(60) (1989), 145--162.
\item Anastasiei, Mihai; Kawaguchi, Hiroaki, {\it A
      geometrical theory of time dependent Lagrangians. I. Nonlinear
      connections}. Tensor (N.S.) 48 (1989), no. 3, 273--282.
\item Anastasiei, Mihai; Kawaguchi, Hiroaki, {\it A
      geometrical theory of time dependent Lagrangians. II. $M$-connections}.
      Tensor (N.S.) 48 (1989), no. 3, 283--293.
\item Anastasiei, Mihai; Kawaguchi, Hiroaki, {\it A
      geometrical theory of time dependent Lagrangians. III. Applications}.
      Tensor (N.S.) 49 (1990), no. 3, 296--304.
\item Anastasiei, Mihai, {\it Cross-section submanifolds of
      cotangent bundle over an Hamilton space}. Rend. Sem. Fac. Sci. Univ.
      Cagliari 60 (1990), no. 1, 13--21.
\item Miron, Radu; Ro\c sca, Radu; Anastasiei, Mihai;
      Buchner, Klaus, {\it New aspects of Lagrangian relativity}. Found. Phys. Lett. 5
      (1992), no. 2, 141--171.
\item Anastasiei, M., {\it Cross-section submanifolds in
      vector bundles}. Differential geometry and its applications (Eger, 1989),
      23--31, Colloq. Math. Soc. János Bolyai, 56, North-Holland, Amsterdam,
      1992.
\item Anastasiei, Mihai; Kawaguchi, Hiroaki, {\it Geometry
      of multiparame\-trized Lagrangians}. Publ. Math. Debrecen 42 (1993), no. 1-2,
      29--37.
\item Anastasiei, Mihai; Kawaguchi, Hiroaki, {\it Absolute
      energy of a Finsler space}. International Conference on Differential
      Geometry and its Applications (Bucharest, 1992). Tensor (N.S.) 53 (1993),
      Commemoration Volume I, 108--113.
\item Anastasiei, Mihai, {\it The geometry of time-dependent
      Lagrangians}. Lagrange geometry, Finsler spaces and noise applied in
      Biology and Physics. Math. Comput. Modelling 20 (1994), no. 4-5, 67--81.
\item Anastasiei, Mihai, {\it On divergence in Finsler spaces},Proceedings of the Workshop on Global Analysis,
      Diff. Geometry and Lie Algebras, 1994, 1-8.
\item Anastasiei, Mihai, {\it On deflection tensor field in
      Lagrange geometries}. Lagrange and Finsler geometry, 1--14, Fund. Theories
      Phys., 76, Kluwer Acad. Publ., Dordrecht, 1996.
\item Anastasiei, Mihai; Antonelli, Peter Louis, {\it The
      differential geometry of Lagrangians which generate sprays}. Lagrange and
      Finsler geometry, 15--34, Fund. Theories Phys., 76, Kluwer Acad. Publ.,
      Dordrecht, 1996.
\item Anastasiei, Mihai; Miron, Radu, {\it Preface for
      ``Generalized Finsler metrics}.\\ Finsler geometry (Seattle, WA, 1995),
      157--159, Contemp. Math., 196, Amer. Math. Soc., Providence, RI, 1996.
\item Anastasiei, Mihai, {\it Certain generalizations of
      Finsler metrics}. Finsler geometry (Seattle, WA, 1995), 161--169, Contemp.
      Math., 196, Amer. Math. Soc., Providence, RI, 1996.
\item Anastasiei, Mihai, {\it A historical remark on the
      connections of Chern and Rund}. Finsler geometry (Seattle, WA, 1995),
      171--176, Contemp. Math., 196, Amer. Math. Soc., Providence, RI, 1996.
\item Anastasiei, Mihai, {\it Geometry of higher order
      sprays}. New frontiers in algebras, groups and geometries (Monteroduni,
      1995), 169--178, Ser. New Front. Adv. Math. Ist. Ric. Base, Hadronic
      Press, Palm Harbor, FL, 1996.
\item Anastasiei, Mihai, {\it A class of generalized
      Lagrange spaces}. Proceedings of the 25th National Conference on Geometry
      and Topology (Ia\c si, 1995). An. \c Stiin\c t. Univ. Al. I. Cuza Ia\c si.
      Mat. (N.S.) 42 (1996), suppl., 259--264.
\item Anastasiei, Mihai, {\it Finsler connections in
      generalized Lagrange spaces}. Balkan J. Geom. Appl. 1 (1996), no. 1, 1--9.
\item Anastasiei, Mihai, {\it Gradient, divergence and
      Laplacian in generalized Lagrange spaces}. Mem. Sec\c t. \c Stiin\c t.
      Acad. Rom\^{a}n\u a Ser. IV 19 (1996), 115--120 (1998).
\item Anastasiei, Mihai; Bucataru, Ioan, {\it A notable
      submersion in the higher order geometry}. Proceedings of the Workshop on
      Global Analysis, Differential Geometry and Lie Algebras (Thessaloniki,
      1995), 1--10, BSG Proc., 1, Geom. Balkan Press, Bucharest, 1997.
\item Anastasiei, Mihai; \v Comi\'c, Irena, {\it Geometry
      of $k$-Lagrange spaces of second order}. 11th Yugoslav Geometrical Seminar
      (Div\v cibare, 1996). Mat. Vesnik 49 (1997), no. 1, 15--22.
\item Anastasiei, Mihai; Buc\u{a}taru, Ioan, {\it Jacobi
      fields in generalized Lagrange spaces}. Collection of papers in honor of
      Academician Radu Miron on his 70th birthday. Rev. Roumaine Math. Pures
      Appl. 42 (1997), no. 9-10, 689--695.
\item Anastasiei, Mihai; Shimada, Hideo, {\it Beil metrics
      associated to a Finsler space}. Balkan J. Geom. Appl. 3 (1998), no. 2,
      1--16.
\item Anastasiei, Mihai; Shimada, Hideo; Sab\u{a}u,
      Sorin, {\it On the nonlinear connection of a second order Finsler space}.
      Proceedings of the 10th National Conference on Finsler, Lagrange and
      Hamilton Geometry (Craiova, 1998). Algebras Groups Geom. 16 (1999), no. 1,
      1--10.
\item Anastasiei, Mihai, {\it Locally conformal Kaehler structures on tangent manifold of a space form}. Libertas Math. 19
      (1999), 71--76.
\item Miron Radu, Anastasiei Mihai, Bucataru Ioan, {\it On the geometry of
    higher order Lagrange spaces}. APH N.S. Heavy Ion Physics
    11)(2000), 1-10.
\item Anastasiei, Mihai, {\it Some Riemannian almost product structures on tangent manifold}. Proceedings of the 11th
    National Conference on Finsler, Lagrange and Hamilton Geometry (Craiova,
    2000). Algebras Groups Geom. 17 (2000), no. 3, 253--262.
\item Anastasiei, Mihai, {\it A framed $f$-structure on tangent manifold of a Finsler space}. Proceedings of the
      Centennial "G. Vr\u anceanu" and the Annual Meeting of the Faculty of Mathematics(Bucharest, 2000). An. Univ. Bucure\c sti Mat.
      Inform. 49 (2000), no. 2, 3--9.
\item Anastasiei, Mihai; Hrimiuc, Drago\c s, {\it Generalizations of Finsler geometry}. Finslerian geometries (Edmonton, AB, 1998), 3--7,
    Fund. Theories Phys., 109, Kluwer Acad. Publ., Dordrecht, 2000.
\item Anastasiei, Mihai; Shimada, Hideo, {\it Deformations of Finsler metrics}. Finslerian geometries (Edmonton, AB,
   1998), 53--65, Fund. Theories Phys., 109, Kluwer Acad. Publ., Dordrecht, 2000.
\item Anastasiei, Mihai,  {\it Distributions on spray spaces}. Balkan J. Geom. Appl. 6 (2001), no. 1, 1--6.
\item Anastasiei, Mihai, {\it Symplectic connections in Lagrange
      geometry}. First Conference of the Mathematical Society of the Republic of
      Moldova (Chi\c sin\u au, 2001). Bul. Acad. \c Stiin\c te Repub. Mold. Mat.
      2001, no. 3, 57--66.
\item Anastasiei, Mihai, {\it Symplectic structures in fibre bundles}. Stud. Cercet. \c Stiin\c t. Ser. Mat.
    Univ. Bac\u au no. 11 (2001), 7--12 (2002).
\item Anastasiei, Mihai Symplectic structures and
Lagrange geometry. Finsler and Lagrange geometries (Ia\c si,
2001), 9--16, Kluwer Acad. Publ., Dordrecht, 2003.
\item Anastasiei, Mihai;Ciobanu, Gabriela; Gottlieb, Ioan, Contraforms on pseudo-Riemannian
manifolds. Finsler and Lagrange geometries (Ia\c si, 2001),
249--258, Kluwer Acad. Publ., Dordrecht, 2003.
\item Anastasiei, Mihai Metrizable linear
connections in vector bundles. Dedicated to Professor Lajos
Tamássy on the occasion of his 80th birthday. Publ. Math. Debrecen
62 (2003), no. 3-4, 277--287.
\item Miron, Radu; Anastasiei, Mihai;
Bucataru, Ioan The geometry of Lagrange spaces. Handbook of
Finsler geometry. Vol. 1, 2, 969--1122, Kluwer Acad. Publ.,
Dordrecht, 2003.
\item Anastasiei, Mihai On P. L. Antonelli
works in mathematical biology. Sci. Ann. Univ. Agric. Sci. Vet.
Med. 46 (2003), no. 2, 3--8.
\item Anastasiei Mihai,{\it Minkowskian G-structures in
vector bundles}, R. Miron (ed.) Lagrange and Hamilton Geometries and
Applications (in memory of Gr. Tsagas). Ed. Fair Partners,
Bucharest, 2004, 1-10.
\item Anastasiei, Mihai, Finsler vector bundles--metrizable
connections. Period. Math. Hungar. 48 (2004), no. 1-2, 83--91.
\item Anastasiei, Mihai, Geometry of Berwald vector bundles. Algebras
Groups Geom. 21 (2004), no. 3, 251--262.
\item Anastasiei, Mihai Geometry of Berwald-Cartan spaces. Conference ``Applied
Differential Geometry: General Relativity'' - Workshop ``Global
Analysis, Differential Geometry, Lie Algebras'', 1--9, BSG Proc.,
11, Geom. Balkan Press, Bucharest, 2004.
\item Anastasiei, Mihai , Metrizable nonlinear
connections. Sci. Ann. Univ. Agric. Sci. Vet. Med. 47 (2004), no.
2, 71--78.
\item Anastasiei, Mihai, New properties of
Berwald-Cartan spaces. Libertas Math. 24 (2004), 3--10.
\item Ikeda, Satoshi; Anastasiei, Mihai, Some modified
connection structures associated with the Finslerian gravitational
field. Balkan J. Geom. Appl. 11 (2006), no. 1, 81--86.
\item  Anastasiei, Mihai Geometry of Lagrangians and semisprays on Lie
algebroids. Proceedings of the 5th Conference of Balkan Society of
Geometers, 10--17, BSG Proc., 13, Geom. Balkan Press, Bucharest,
2006.
\item Anastasiei, Mihai, Mechanical systems on Lie algebroids.
Algebras Groups Geom. 23 (2006), no. 3, 235--245.
\item \v{C}omi\'c, Irena; Anastasiei,Mihai, Curvature theory of generalized connection in $J\sp 2\sb
kM$. Novi Sad J. Math. 36 (2006), no. 1, 141--151.
\item Anastasiei, Mihai; Nimine\c t, Valer Prolongations of Lie
algebroids. Stud. Cercet. \c Stiin\c t. Ser. Mat. Univ. Bac\u au
no. 16, suppl. (2006), 29--37.
\item Anastasiei, Mihai Dynamics on Lie algebroids. Sci. Ann. Univ. Agric. Sci. Vet. Med. 49 (2006), no. 2, 137--143.
\item  Anastasiei, Mihai Metrical aspects in theory of Lie algebroids. Bull. Transilv. Univ. Brasov Ser. B (N.S.) 14(49) (2007), suppl., 11--21.
\item Anastasiei, Mihai Lagrange geometry on fibre bundles and related Finsler geometry. I: Lagrangians on bundles. (English)
Bull. Soc. Sci. Lett. Lódz, Sér. Rech. Déform. 57 (54), 33-41 (2007).
\item Anastasiei, Mihai Lagrange geometry on fibre bundles and related Finsler geometry. II: Semisprays, nonlinear connections and metrical connections in Lagrange spaces. (English) Bull. Soc. Sci. Lett. Lódz, Sér. Rech. Déform. 57 (54), 43-61 (2007).
\item Anastasiei, Mihai, A generalization of Myers theorem. An. \c Siin\c{}. Univ. Al. I. Cuza Iasi, sI a, Matematica, T. LIII,2007, Supliment  p. 33-40.
\item Anastasiei, Mihai Semisprays on Lie algebroids. Applications. Tensor N.S., vol. 69(2008),190-198
\item  Anastasiei, Mihai; Gheorghe, Marinela, Some examples of Randers spaces. Acta Math. Acad. Paedagog. Nyházi. (N.S.) 24 (2008), no. 1, 15--23.
\item Anastasiei, Mihai; Vacaru, Sergiu I. Fedosov quantization of Lagrange - Finsler and Hamilton-Cartan spaces and Einstein gravity lifts on (co) tangent bundles. J. Math. Phys. 50 (2009), no. 1, 013510, 23 pp.
\item Anastasiei, Mihai Metrizable linear connections in a Lie algebroid. J. Adv. Math. Stud. 3 (2010), no. 1, 9--18.
\item Anastasiei, Mihai; G\^{\i}r\c{t}u, Manuela,  Indicatrix of a Finsler vector bundle. Sci. Stud. Res. Ser. Math. Inform. 20 (2010), no. 2, 21-–28.
\item Anastasiei Mihai ; Vacaru Sergiu, Nonholonomic Black Ring and Solitonic Solutions in Finsler and Extra Dimension Gravity Theories.Internat.  J. Theoret. Phys 49(2010), no. 8, 1788-1804
\item Anastasiei, M., Banach Lie algebroids. An.{\c S}tiin{\c t}.Univ. ``Al.~I.~Cuza'' Ia{\c s}i (S.N.).Matematica, Tomul LVII, 2011,f.2, 409--416

\end{enumerate}

\item[{\bf V.}]{\bf Others}
\begin{enumerate}
\item Miron, Radu; Anastasiei, Mihai; Yawata, Makoto; Kawaguchi, Tomoaki, {\it On the state of Mathematics Education in Romania}.
      Report of Chiba Inst. of Tech. nr.41, 1994, 179-182.
\item Anastasiei, M.; Atanasiu, Gheorghe, {\it Professor Radu
      Miron, on his 60th birthday}. The Proceedings of the Fifth National Seminar
      of Finsler and Lagrange Spaces (Bra\c sov, 1988), 9--31, Soc. \c Stiin\c
      te Mat. R.S. România, Bucharest, 1989.
\item Anastasiei, Mihai, {\it Academician Radu Miron on his
      70th birthday}. Collection of papers in honor of Academician Radu Miron on
      his 70th birthday. Rev. Roumaine Math. Pures Appl. 42 (1997), no. 9-10,
      673--687.
\item Miron, Radu;Anastasiei, Mihai, {\it Professor Dr. Peter Louis
    Antonelli at sixty.} In the volume Finsler and Lagrange
    Geometries edited by M. Anastasiei and P.L.Antonelli, Kluwer
    Academic Publishers, Dordrecht, 2003.
\item Anastasiei, Mihai Academician Radu Miron at eighty. An. \c Stiin\c t. Univ. Al. I. Cuza Ia\c si. Mat. (N.S.) 53 (2007), suppl. 1, 3--10.
\item Anastasiei, Mihai; Cioban, Mitrofan; Soltan, Petru, Academician Radu Miron---eighty years of life and sixty years of efforts. Bul. Acad. \c Stiin\c te Repub. Mold. Mat. 2008, no. 2, 130--136.

\end{enumerate}

\end{itemize}

\newpage

\runningauthor={M. ANASTASIEI}
\runningtitle={CONFORMAL STRUCTURES ON BANACH VECTOR BUNDLES}
\noindent
\baselineskip 8pt
\noindent{\footnotesize{ANALELE {\c S}TIIN{\c T}IFICE ALE UNIVERSIT{\u A}{\c T}II ``AL.~I.~CUZA" IA{\c S}I}}
\hfill\break
{\footnotesize{Tomul XX, s.I a, Matematic\u a, 1974, f.2., p. 351--358}}
\vskip 2cm
\baselineskip 11.5pt plus .15pt
\centerline{\bf\Large CONFORMAL STRUCTURES}
\vskip .2cm
\centerline{\bf\Large ON BANACH VECTOR BUNDLES}
\vskip .5cm
\centerline{\bf {\footnotesize{BY}}}
\vskip .5cm
\centerline{\bf {\footnotesize{M. ANASTASIEI}}}
\vskip 1cm

In this paper, conformal structures, in particular Weyl structures on Banach vector bundles are defined. We prove that there is a one-to-one correspondence between the set of conformal structures on a Banach vector bundle and the set of reductions of its structural group to the conformal group.

The existence and the uniqueness of a connection without torsion, compatible with a Weyl structure on a Banach manifold are  proved.

\section{Linear conformal space. Conformal group} Let $\mathbf{E}$ be a real, infinite-dimensional linear space.

{\bf Definition 1.1.} A conformal structure on $\mathbf{E}$ is a set $C(\mathbf{E})$ of scalar products on $\mathbf{E}$, denoted by $(\ ,\ )_a$, $a\in{\cal I}$ which satisfy
\begin{equation}
(\ ,\ )_a=\lbb_{ab}(\ ,\ )_b\ \ a,b\in{\cal I}
\end{equation}
where $\lbb_{ab}$ is a positive real number, and ${\cal I}$ a set of indices.

\ms

{\bf Definition 1.2.} The linear space ${\mathbf{E}}$ with conformal structure $C(\bfe)$ is called a {\it conformal space}.

\ms

{\bf Remark 1.1.} In the conformal space $\bfe$ the angle between two vectors can be defined by
\begin{equation}
\cos(u,v)=\dfrac{(u,v)_a}{\|u\|_a\|v\|_a},\ \forall u,v\in\bfe,\ a\in{\cal I},
\end{equation}
where $\|u\|_a=\sqrt{(u,u)_a}$, the ratio of their lengths is well defined but their absolute lengths are not defined.

If $(\ ,\ )$ is a fixed element of $C(\bfe)$, (1.1) can be replaced by
\begin{equation}
(\ ,\ )_a=\lbb_a(\ , \ )\ \ \forall a\in{\cal I}.
\end{equation}
This implies that all norms $u\to\|u\|_a$ are equivalent to the fixed norm $u\to\|u\|=\sqrt{(u,u)}$. In the following, the space $\bfe$ with the norm $\|\cdot\|$ will be supposed to be complete.

\ms

{\bf Remark 1.2.} Let $R_i(\bfe)$ be the set of all scalar products on $\bfe$. We have
\begin{equation}
C(\bfe)\subseteq R_i(\bfe)\subseteq L^2_s(\bfe).
\end{equation}
(Here $L_s^2(\bfe)$ is the linear space of bilinear and symmetric maps $s:\bfe\times\bfe\to\mathbf{R}$.)

Let $L(\bfe)$ be the linear space of linear bounded operators on $\bfe$.

\ms

{\bf Definition 1.3.} We say that $A\in L(\bfe)$ preserves the conformal structure of $\bfe$ if there is a unique $a\in{\cal J}$ such that $$(Au,Av)=(u,v)_a,\ \ u,v\in\bfe.$$

If $A^*$ denotes the adjoint operator of $A$ with respect to the scalar product $(\ ,\ )$, then we have

\ms

{\bf Theorem 1.1.} {\it Let $\bfe$ be a conformed space and $A\in L(\bfe)$. The following conditions are equivalent

1) $A$ preserves the conformal structure of $\bfe$;

2) There is a unique real number $k_A >0$ such that
$$A^* A = k_A I\ \ (I\ {\rm{identity\ operator}});$$

3) $A$ preserves the angle between vectors of $\bfe$.}

\ms

\noindent{\it Proof.} Obvious.

\ms

{\bf Definition 1.4.} An operator $A\in L(\bfe)$ which satisfies one of the conditions of Theorem 1.1 will be called a {\it conformal operator}.

\ms

Let $CO(\bfe)$ be the set of invertible conformal operators and $O(\bfe)$ the set of invertible operators which satisfy $A^*A = I$. We immediately obtain the following

\ms

{\bf Theorem 1.2.} {\it The sets $CO(\bfe)$ and $O(\bfe)$ are subgroups of the group $GL(\bfe)$ of all invertible operators.}

\ms

We call $CO(\bfe)$ and $O(\bfe)$ the conformal group and the orthogonal group of $\bfe$, respectively.

\ms

{\bf Theorem 1.3.} {\it There is an isomorphism $$\a:CO(\bfe)\to O(\bfe)\times R^*_+,$$ where $R^*_+$ is the positive real multiplicative group.}

\ms

\noindent{\it Proof.} For $A\in CO(\bfe)$ and $A^*A = k_A I$, we put $\a(A)=\left(\dfrac{1}{\sqrt{k_A}}A,k_A\right)$ and for $(B,1)\in O(\bfe)\times R^*_+$, $\a^{-1}(B,1)=\sqrt{1}B$.

Let $GL(\bfe)$ be endowed with the topology induced by the norm topology of $L(\bfe)$. We identify the group $R^*_+$ with the group of homotheties of $\bfe$. The following result is obvious.

\ms

{\bf Theorem 1.4.} {\it The subgroups $CO(\bfe)$, $O(\bfe)$ and $R^*_+$ are closed, and the map $\a$ from Theorem 1.3 is a topological isomorphism.}

\section{Conformal structures on a Banach vector bundle}
\setcounter{equation}{0}

All vector bundles, manifolds and maps considered in the following sections will be assumed of class $C^\9$.

Let $E$ and $M$ be manifolds, modeled on Banach spaces, and suppose $M$ connected. Let $\pi : E \to M$ be a vector bundle with fibre the conformal space $\bfe$. We denote by $Ri(\pi)$ the set of Riemannian metrics on $\pi$, and we define the following equivalence relation: $$\Lambda:g\sim g'\Leftrightarrow g'=e^\lbb\cdot g, \ \ g,g'\in Ri(\pi),$$ where $\lbb$ is a smooth function on $M$.
(The use of exponential function is a handy way of ensuring positivity).\ms

{\bf Definition 2.1.} A conformal structure on $\pi$ is an equivalence class $\bfc$ with respect to $\Lambda$ of Riemannian metrics on $\pi$.

\ms

{\bf Remark 2.1.} If the equivalence class $\bfc$ contains only one element, we obtain a Riemannian structure on $\pi$.

\ms

The proofs of the two following theorems are standard. (See [3, Ch. 7] for the particular case of the Riemannian structure).

\ms

{\bf Theorem 2.1.} {\it Let $\pi : E\to M$ be a vector bundle with fibre $\bfe$ and suppose $M$ admits a partition of unity. Then the vector bundle $\pi$ admits a conformal structure.}

\ms

{\bf Theorem 2.2.} {\it Let $\pi : E\to M$ and $\pi': E'\to M'$ be vector bundles and let $f: E'\to E$ be a bundle morphism such that the map $f_{p'}: E'_{p'}\to E_{f(p')}$, where $E'_{p'}=\pi'^{-1}(p')$ and $E_{f(p')}=\pi^{-1}(f(p'))$, is injective and such that $f_{p'}(E'_{p'})$ has a complementary closed subspace in $E_{f(p')}$. Then a conformal structure on $\pi$ canonically induces a conformal structure on $\pi'$.}

\ms

{\bf Definition 2.2.} The vector bundle $\pi: E\to M$ with fibre $\bfe$ admits a reduction of its structural group to $CO(\bfe)$, if and only if there exists a bundle atlas $(U_i,\tau_i)_{i\in I}$, such that the maps	$(\tau_j\circ\tau^{-1}_i)_p:\bfe\to\bfe$ for each $p$ of $U_i\cap U_j$ belong to the group $CO(\bfe)$.

\ms

{\bf Theorem 2.3.} {\it Let $\pi: E \to M$ be a vector bundle with fibre $\bfe$. There exists a one-to-one correspondence between the set of reductions of the structural group to the conformal group and the set of conformal structures on $\pi$.}

\ms

\noindent{\it Proof.} Every reduction of $\pi$ to the conformal group $CO(\bfe)$ determines the conformal structure of $\pi$. Indeed, we define
$$g_{a,p}(v,w)=(\tau_{i,p},v,\tau_{i,p},w)_a,\  \forall v,w\in \bfe {\rm{\ and\ }} a\in{\cal J}.$$ The maps $g_a:p\to g_{a,p}$ define the sections of vector bundle $L^2_s(\pi)$ (see [3, Ch. 3] for the definition of this vector bundle) and the set $\{g_a\}$ is a conformal structure on $\pi$.

Conversely, let $\{g_j\}_{j\in \bfj}$ be a conformal structure on $\pi$ and let $(U_i,\tau_i)_{i\in \bfi}$ be a bundle atlas for $\pi$. We consider an arbitrary map $\varepsilon : {\mathbf I}\to{\mathbf J}$ and let $g^{\varepsilon(i)}_i$ be the induced metric by $g_j$ with $j=\vare(i)$, on $U_i\times\bfe$ by the isomorphism $\tau_i$. There exists a positive definite symmetric operator $A^{\vare(i)}_{i,p}$ such that $$g^{\vare(i)}_{i,p}(v,w)=(A^{\vare(i)}_{i,p}v,w),\ \forall p\in U_i,\ v,w\in\bfe.$$

We denote $B_{i,p} = \sqrt{A^{\vare(i)}_{i,p}}$ and we put $\sigma_i = B_{i,p}\circ\tau_{i,p}$. Then $(U_i,\sigma_i)$ is the bundle atlas we looked for. It is sufficient to prove that $B_i:U_i\times\bfe\to U_i\times\bfe$ which is defined on fibres by $B_{i,p}$ map $g^{\vare (i)}_i$ on the scalar product $(\ ,\ )$ of $\bfe$. But we have
$$(B_{ip}v, B_{i,p}w) = (A^{\vare(i)}_{i,p}v,w) = g^{\vare (i)}_i(v,w),$$ since $B_{i,p}$ is symmetric.

\ms

{\bf Definition 2.3.} Let $\pi:E\to M$ be a vector bundle with a conformal structure $\C$; $\C$ is called a Weyl structure if and only if there exists a map $W:\C\to C^\9(T^* M)$ which satisfies
$$W(e^\lbb\cdot g)= W(g)-d\lbb,$$ where $C^\9(T^*M)$ denotes the set of sections in cotangent bundle of $M$.

\ms

\noindent{\bf Remark 2.2.} A Riemannian metric $g$ and a $1-$form $\eta$ on $M$ determine a Weyl structure, namely $W: \C\to C^\9(T^* M)$, where $\C$ is the equivalence class of $g$ and $W(e^\lbb g)=\eta-d\lbb$.

\ms

{\bf Theorem 2.4.} {\it Let $\pi : E\to M$ and $\pi' : E'\to M'$ be vector bundles with conformal structures and $f: E'\to E$ a bundle morphism compatible with this conformal structure (in the sense of Theorem 2.2). Every Weyl structure on $\pi$ canonically induces a Weyl structure on $\pi'$.}

\ms

\noindent {\it Proof.} If $(g,\eta)$ defines Weyl structure on $\pi$, then $(f^*g,f^*\eta)$ defines a Weyl structure on $\pi'$.

\ms

{\bf Theorem 2.5.} {\it Let $\pi:E\to M$ be a vector bundle with a conformal structure, where $M$ admits a partition of unity. Then $\pi$ admits a Weyl structure.}

\ms

\noindent{\it Proof.} It is sufficient to prove that there is a global 1-form on $M$. But this follows from [1, Lemma 1.3].

\section{Connections compatibles with\\ conformal structures}
\setcounter{equation}{0}

In this section we give some results from the connections theory on Banach vector bundles which we shall use later.

\ms

{\bf Theorem 3.1.} {\it Let $\pi: E\to M$ be a vector bundle and $M$ admits a partition of unity; then

i) there exists a connection map $K$ for $\pi$,

ii) there exists a canonical bijective map from the set of connection maps on $\pi$ to the set of covariant derivatives on $\pi$ given by
\begin{equation}
K\circ T\xi = \nabla\xi,\ \ \forall\xi\in{\cal X}_E(M),
\end{equation}
where ${\cal X}_E(M)$ is the set of the sections in $\pi$.}

\ms

The proof is given in [1, Theorem 2.2].

\ms

{\bf Remark 3.1.} a) $\nabla\xi$ is considered as a section in $L(\tau,\pi):L(TM,E)\to M$ where $\tau : TM\to M$ is tangent bundle.

b) The implication $K\to\nabla$ is given by (3.1) without the hypothesis of existence of the partition of unity on $M$.

Let $c:[0,1]\to M$ be a piecewise differentiable curve on $M$. We denote by $P_c|_{[t,t_0]}$, where $t,t_0\in(0,1)$, the parallel displacement from $E_{c(t)}$ to $E_{c(t_0)}$ defined by the connection $\nabla,K$. The map $\wt{Q}_c : c^*E\to (0,1)\times E_{c(t_0)}$ defined by
\begin{equation}
\wt{Q}_c(t,v) = (t,P_c|_{[t,t_0]}v)\ \  (t,v)\in c^*E.
\end{equation}
is a vector bundles isomorphism. See [1, Theorem 3.5].

Let ${\cal X}_E(c)$ be the vector space of section in $\pi$ along the curve $c$ and let $C^\9((0,1),E_{c(t_0)})$ be the vector space of maps of class $C^\9$ from $(0,1)$, to $E_{c(t_0)}$. We consider the map $Q_c: {\cal X}_E(c)\to C^\9((0,1), E_{c(t_0)})$ defined by
\begin{equation}
Y\to Q_c Y=pr_2\circ\wt{Q}_c(t,Y(t))\ \ \forall Y\in{\cal X}_E(c),\ t\in(0,1).
\end{equation}

\ms

{\bf Theorem 3.2.} {\it a) $Q_c$ is a vector space isomorphism;

b) $\dfrac{d}{dt}(Q_c Y) = Q_c(\nabla_c,Y)$ where $\nabla_c Y$ is covariant differentiation of section $Y$ along curve $c$.}

\ms

For proof see [1, Theorem 3.6].

In [4], the holonomy group of connection $\nabla,K$ with reference point $p$, denoted by $\Phi(p)$, is defined. For each $p$ of $M$, the group $\Phi(p)$ can be realized as a subgroup of the structural group of $\pi$. Theorem 2.11 of [4] suggests the following

\ms

{\bf Definition 3.1.} Let $\pi: E\to M$ be a vector bundle with a conformal structure $\C$. The connection $\nabla, K$ is compatible with the conformal structure $\C$ if and only if
\begin{equation}
\Phi(p)\subseteq CO(\bfe),\ \  \forall p\in M,
\end{equation}
where $\bfe$ is the fibre of $\pi$.

In the interesting case when $\C$ is at the same time a Weyl structure, we will use the following

\ms

{\bf Definition 3.2.} (See [2]). Let $\pi : E\to M$ be a vector bundle with a Weyl structure $(g,\eta)$. The connection $\nabla, K$ is compatible with the Weyl structure $(g,\eta)$ if and only if along every curve $c:[0,1]\to M$ and for at least one $g$ from $\C$,
\begin{equation}
g_p(Q^t_C Y,Q^t_C Z)=\exp\left[\int^t_0 c^*\eta\right]g_{c(t)}(Y_t,Z_t),
\end{equation}
where $Q^t_c=P_C|_{[t,t_0]}, p=c(0)$ and $Y_t,Z_t\in E_{c(t)}$.

\ms

{\bf Remark 3.2.} If condition (3.5) is satisfied by one $g\in\C$, it will be satisfied by every $g'=e^\lbb\cdot g$ belonging to $\C$.

\ms

{\bf Theorem 3.3.} {\it Let $\pi:E\to M$ be a vector bundle with the Weyl structure $(g,\eta)$ and let $\nabla, K$ be a connection on $\pi$.

The following assertions are equivalent:

$1)$	The connection $\nabla, K$ is compatible with the Weyl structure $(g,\eta)$;\\

$2)$	$\dfrac{d}{dt}g_c(Y,Z)=g_c(\nabla_c Y,Z)+g_c(Y, \nabla_c Z)- c^*\eta\cdot g_c(Y,Z)$, $\forall c:[0,1]\to M$ and $Y,Z\in{\cal X}_E(M)$;\\

$3)$	$Xg(Y,Z)=g(\nabla_X Y,Z)+g(Y,\nabla_X Z)-\eta(X)g(Y,Z)$ $\forall X\in{\cal X}_{TM}(M)$ and $Y,Z\in{\cal X}_E(M)$}

\ms

\noindent{\it Proof.} 1) $\to$ 2). For each curve $c$ with $c(0) = p$, we shall denote $Y_{c(t)}=Y_t$, $Z_{c(t)}=Z_t$, $Y_p=Y$, $Z_p=Z$ for $Y,Z\in{\cal X}_E(M)$. It follows from (3.4) and b) of Theorem 3.2 that
$$\begin{array}{c}\dfrac{d}{dt}g_c(Y,Z)=\dd\lim_{t\to 0}\dfrac1t(g_{c(t)}(Y_t,Z_t)-g_p(Y,Z))=\\ \\
=\dd\lim_{t\to 0}\dfrac1t\left(\exp\[-\int^t_0c^*\eta\]g_p(Q^t_c Y_t,Q^t_c Z_t)-g_p(Y,Z)\right)=\\ \\
=\dd\lim_{t\to 0}\dfrac1t\exp\[-\int^t_0c^*\eta\](g_p(Q^t_c Y_t,Q^t_c Z_t)-g_p(Y,Z))+\\ \\
+g_{p}(Y,Z)\dd\lim_{t\to 0}\dfrac1t\left(\exp\[-\int^t_0c^*\eta\]_{-1}=g_p(\lim_{t\to 0}\dfrac{1}{t}(Q_c Y_t-Y),Z)+\right.\\ \\
+g_p(Y,\dd\lim_{t\to 0}\dfrac{1}{t}(Q^t_c Z_t-Z))+g_p(Y,Z)\dfrac{d}{dt}\exp\[-\int^t_0c^*\eta\]=\\ \\
\left.=g_p(\nabla_c Y,Z)+g_p(Y,\nabla_c Z)-c^*\eta\cdot g_p(Y,Z)\ {\rm{i.e.}} 2\right).\end{array}$$

2) $\to$ 1). Let $Y,Z$ be parallel sections in $\pi$ along $c$, i.e $Q^t_c Y_t = Y_p$ and $Q^t_c Z_t=Z_p$. Assertion 2) of Theorem becomes
\begin{equation}
\dfrac{d}{dt}g_{c(t)}(Y_t,Z_t)=-c^*\eta\cdot g_{c(t)}(Y_t,Z_t),
\end{equation}
and we get (3.5) by integration.

The proof of 2) $\leftrightarrow$ 3) can be obtained in the same way as in the Riemannian case. See [1, Theorem 3.8].

\ms

{\bf Definition 3.3.} We say that a manifold $M$ modeled by the conformal space $\bfm$ is endowed with a conformal structure if there exists a collection of charts $(U_i,\varphi_i)$, covering $M$ and satisfying
\begin{equation}
D(\varphi_j\circ\varphi^{-1}_i)_{\varphi_i(p)}\in CO(\bfm)\ {\rm{for\ all\ }}i,j\ {\rm{and\ }}p\in U_i\cap U_j,
\end{equation}
where $D$ denotes the differentiation operator.

\ms

{\bf Theorem 3.4.} {\it A manifold $M$ modeled by a conformal space $\bfm$ admits a conformal structure if and only if the tangent bundle $TM$ admits a conformal structure.}

\ms

\noindent{\it Proof.} Let $(U_i, \varphi_i)$ be the collection of charts which defines the conformal structure on $M$. The transition maps of $TM$ are $D(\varphi_j\circ \varphi^{-1})_{\varphi_i(p)}$ and belong to $CO(\bfm)$ i.e. $TM$ admits a reduction to conformal group, therefore a conformal structure by Theorem 2.3.

Conversely, a conformal structure on $TM$ induces a reduction of this vector bundle to the conformal group $CO(\bfm)$ i.e. the maps $D(\varphi_j\circ \varphi^{-1}_i)_{\varphi_i(p)}$ belong to $CO(\bfm)$.

\ms

{\bf Definition 3.4.} A conformal manifold $M$ is called a Weyl manifold if and only if the conformal structure of $TM$ is a Weyl structure.

\ms

{\bf Theorem 3.5.} {\it Let $M$ be a Weyl manifold, modeled by the conformal space $\bfm$. There exists a unique connection $\nabla, K$, such that

\noindent i) $Xg(Y,Z)=g(\nabla_X Y,Z)+g(Y,\nabla_X Z)-\eta(X)g(Y,Z)$ for $X,Y,Z\in{\cal X}_{TM}(M)$

\noindent ii) $T(X,Y)\overset{def}{=}\nabla_X Y-\nabla_Y X-[X,Y]=0$, $\forall X,Y\in{\cal X}_{TM}(M)$
where $(g,\eta)$ is the Weyl structure of $TM$.}

\ms

\noindent{\it Proof. Existence.} Let $(U,\varphi)$ be a chart for $M$ at $p$. We consider the following equation with Fr\'{e}chet derivatives
\begin{equation}
\begin{array}{c}
2g_\varphi(\Gamma_{\varphi(p)}((u,v),w))=Dg_\varphi|_{\varphi(p)}(u)(v,w)+\\ \\
+Dg_\varphi|_{\varphi(p)}(v)(u,w)-Dg_\varphi|_{\varphi(p)}(w)(u,v)+\eta_\varphi(u)g_\varphi(v,w)+\\ \\
+\eta_\varphi(v)g_\varphi(u,w)-\eta_\varphi(w)g_\varphi(u,v),\ \ \forall u,v,w\in \bfm,
\end{array}
\end{equation}
where $g_\varphi$ and $\eta_\varphi$ are local representatives of $g$ and $\eta$, respectively. This equation defines a map $\Gamma_{\varphi(p)}\in L^2_s(M,M)$, such that $\varphi(p)\to \Gamma_{\varphi(p)}$ is of class $C^\9$. As the $\Gamma_{\varphi(p)}$ satisfies the usual transformation formula of a local connector, under change of trivialization, it defines a connection on $M$. The connection obtained in such a way satisfies i) and ii) of Theorem. Indeed, ii) has the following local expression $$T(X,Y)_{\varphi(p)}=\Gamma_{\varphi(p)}(X_{\varphi(p)},Y_{\varphi(p)})-\Gamma_{\varphi(p)}(Y_{\varphi(p)},X_{\varphi(p)})=0.$$

This equality is satisfied because $\Gamma_{\varphi(p)}$ is a bilinear symmetric map. The local expression of condition i) is
$$\begin{array}{c}Dg_\varphi|_{\varphi(p)}(X_{\varphi(p)})(Y_{\varphi(p)},Z_{\varphi(p)})=g_{\varphi(p)}(\Gamma_{\varphi(p)}(X_{\varphi(p)}),Z_{\varphi(p)})+\\ \\
g_{\varphi(p)}(Y_{\varphi(p)},\Gamma_{\varphi(p)}(X_{\varphi(p)},Z_{\varphi(p)}))-\eta_{\varphi(p)}(X_{\varphi(p)})g_{\varphi(p)}(Y_{\varphi(p)},Z_{\varphi(p)}).\end{array}$$

This equality can be easily verified using (3.8).

\ms

{\it Uniqueness.} Let $\Gamma'_{\varphi(p)}$ be another local connector which fulfils i) and ii). It follows that $\Gamma'_{\varphi(p)}$ must satisfy equation (3.8) i.e, $\Gamma'_{\varphi(p)}=\Gamma_{\varphi(p)}$.

\newpage

{\bf BIBLIOGRAPHY}

\begin{enumerate}
  \item {{Flaschel, P., Klingenberg W.}}, {\it Riemannsche Hilbertmannigfaltigkei\-ten, Periodische Geod\"{a}tische}. Lecture Notes in Mathematics, Vol. 282, Springer - Verlag, 1972.
  \item {{Folland, G. B.}}, {\it Weyl manifolds}, J. Differential geometry 4(1970) p. 145--153.
  \item {{Lang, S.}}, {\it Introduction to differentiable manifolds}, Interscience, New York, 1962.
  \item {{Maxim, L.}}, {\it Connections compatibles with Fredholm structures on Banach manifolds}, An. St. Univ. Iasi, Matematica, T. XVIII, fasc. 2 (1972), p. 389--399.
  \item {{Miron, R.}}, {\it Espaces \`{a} structure conforme presque simplectique}, Colloquium Mathematicum, T. XXVI (1972), p. 207--215.
\end{enumerate}

\bigskip

{\footnotesize{Received 30.XI.1973}}

\newpage





\def\G{\Gamma}

\runningauthor={M. ANASTASIEI}
\runningtitle={AFFINE TRANSFORMATIONS ON BANACH MANIFOLDS }
\noindent
\baselineskip 8pt
\noindent{\footnotesize{Lincei - Rend. Sc. mat. e nat.\hfill\break
Vol. LX, marzo 1976, 4p.}}
\vskip 2cm
\baselineskip 11.5pt plus .15pt
\centerline{\bf\Large GEOMETRIA DIFFERENZIALE}
\vskip .25cm
\centerline{\bf\Large Affine transformations}
\vskip .2cm
\centerline{\bf\Large on Banach manifolds}
\vskip .5cm
\centerline{\bf \footnotesize{Nota di M. ANASTASIEI}}
\vskip .1cm
\centerline{\bf \footnotesize{Presentata dal Socio B. Segre}}
\vskip 1cm

\begin{abstract}
{\bf RIASSUNTO.} In questo lavoro si dimostra che, in determinate condizioni, il gruppo delle trasformazioni affini di una variet\`{a} riemanniana di dimensione infinita coincide col gruppo delle isometrie. Un risultato di questo tipo, nel caso della dimensione finita, \`{e} stato precedentemente ottenuto da S. Kabayashi [2].

\smallskip

In this paper we prove that, under certain conditions, the group of affine transformations of a Riemannian infinite-dimensional manifold $M$ is equal to the group of isometries of $M$. A result of the same type, in the finite-dimensional case, has been obtained by S. Kobayashi [2].
\end{abstract}

\setcounter{section}{0}
\section{Affine morphisms of Banach ma\-nifolds}
We work in the category of infinite-dimensional manifolds of class $C^\9$. Let $M$ be a Banach manifold. We suppose the existence of a connection map $K:T^2M\to TM$ and denote by $\nabla$ the covariant differentiation associated to it, see [1]. For $X,Y$ in $\chi(M)$, the $F(M)-$module of vector fields on $M$,we set $$\nabla_X Y=K\circ TY(X),\ \ \dfrac{Dc}{dt}=K\circ T\dot{c},\leqno(1.1)$$ where $TY$ is the tangent map of $Y:M\to TM$ and $c:[0,1]\to M$ is a curve on $M$. The holonomy groups, denoted by $\Phi(p)$, for $p$ in $M$, were introduced and studied in [4].

{\bf Definition 1.1.} A Banach manifold $M$, endowed with a connection map, is said to be irreducible if $\Phi(p)$ does not have any trivial invariant subspace. Otherwise, it is called reducible.

{\bf Definition 1.2.} Let $M$ and $M'$ be endowed with the connection maps $K$ and $K'$, respectively. A morphism $f:M\to M'$ is called affine if and only if, $$Tf\circ K=K' \circ T^2 f.\leqno(1.2)$$
If $M=M'$ and $f$ is a diffeomorphism, we say that $f$ is an affine transformation.

In the following theorem we collect some facts about affine morphisms, needed in the next section; for the proof see [5].

{\bf Theorem 1.1.} {\it Let $M$ and $M'$ be Banach manifolds with the connection maps $K$ and $K'$, respectively. Suppose $f:M\to M'$ is an affine diffeomorphism. Then:

$a)$ $Tf\circ \tau_c=\tau'_{f\circ c}\circ Tf$ for every curve $c$, where $\tau_c$ (resp. $\tau'_{f\circ c}$) denotes the parallel displacement along the curve $c$ (resp. $f\circ c$);

$b)$ $Tf(\nabla_X Y)=\nabla'_{TfX}Tf Y$, for all $X,Y$ in $\chi(M)$;

$c)$ $Tf\circ R(X,Y)Z=R'(TfX,TfX)TfZ$, for all $X,Y,Z$ in ${\chi}(M)$,

\noindent where $R$ (resp. $R'$) denotes the curvature tensor field associated with $K$ (resp. $K'$).}

Let $(M,g)$ be a Riemannian manifold. As in the finite dimensional case, the sectional curvature for a $2-$plane $\sigma=\{X,Y\}$ in $T_pM$ (the tangent space at $p$ in $M$) is defined by $$K_p(\sigma)=\dfrac{g(R(X,Y)Y,X)}{g(X,Y)g(Y,Y)-g^2(X,Y)}.\leqno(1.3)$$

{\bf Definition 1.3.} Let $(M,g)$ and $(M',g')$ be Riemannian manifolds. A morphism $f:M\to M$ is called a homothety if there exists $c\in\mathbf{R}$ such that $$g'(TfX,TfY)=c^2g(X,Y)\ \textrm{for any } X,Y\textrm{ in } \chi(M).\leqno(1.4)$$
If in (1.4) $c=1$, then $f$ is an isometry.

It is proved in [1, p.38] that every isometry is an affine morphism (with respect to the unique connections without torsion defined by $g$ and $g'$, respectively).

In particular, the group of isometries of $M$ is a subgroup of the group of affine transformations of $M$.

\section{The main results}
The purpose of this section is to prove Theorems 2.1 and 2.2.

{\bf Theorem 2.1.} {\it Let $(M,g)$ be an irreducible Riemannian manifold, with bounded and non-identically zero sectional curvature. Then, the group of affine transformations of $(M,g)$ is equal to the group of isometries of $(M,g)$.}

\ms

\noindent {\it Proof.} The proof will be given in three steps.

{\it Step 1.} Every homothety is an affine transformation. Using a homothety $f$, we define a new Riemannian metric on $M$ by $g'(X,Y)=g(TfX,TfY)=c^2g(X,Y)$. Obviously, $f:(M,g')\to (M,g)$ is an isometry, hence an affine transformation. But, from the definitions of the Riemannian connection [1, p. 36], it follows that the connection defined by $g'$ and $g$ coincide; therefore $f:(M,g)\to(M,g)$ is an affine transformation.

{\it Step 2.} If $(M,g)$ is irreducible, every affine transformation is a homothety. For this we need the following

{\bf Lemma.} {\it Let $H$ be a real Hilbert space, $O(H)$ the orthogonal group and $S$ a subgroup of $O(H)$ which acts irreducibly on $H$. If $g$ is a symmetric and bilinear form on $H$, invariant under the action of $S$, then there is a constant $c$ such that $g(u,v)=c(u,v)$ for all $u,v$ in $H$, $(\ ,\ )$ being the standard inner product of $H$.}

\ms

\noindent{\it Proof of Lemma.} There exists a symmetric operator $A$ such that $g(u,v)=(Au,v)$. Let $s$ be an element of $S$. From $g(su,sv)=g(u,v)$ (invariance of $g$) it follows $As=sA$ for all $s$ in $S$ and from Theorem 6, Appendix II of [3], it follows that there exists a constant $c$, such that $A=cI$ (where $I$ is the identity operator) and therefore $g(u,v)=c(u,v)$. We remark that, if $g$ is positive definite, the constant $c$ must be positive.

We give now the proof of {\it Step 2}.

For $p$ in $M$ there are two inner products $g_p$ and $g'_p$ on $T_pM$, where $g'_p(X,Y)=g(T_pfX,T_pfY)$. As $f$ is an affine transformation $g$ is invariant under the action of $\Phi(p)$ which is a subgroup of the orthogonal group $O(T_pM)$ (with respect to the inner product $g$). We are in position to apply the Lemma and we obtain $g'_p=c^2_pg_p$. But $g$ and $g'$ are the parallel tensor fields with respect to the Riemannian connection defined by $g$, therefore $c_p$ does not depend on $p$ i.e. $f$ is a homothety.

{\it Step 3.} In the hypothesis of Theorem 2.1, every affine transformation is an isometry. Let $f$ be an affine transformation of $M$. By Step 2, $f$ is a homothety. If $c=1$, the proof is complete. Suppose $c<1$, otherwise we may use $f^{-1}$ and denote by $K<+\9$ the bound of the sectional curvature. For every $p$ in $M$ and the $2-$plane $\sigma=\{X,Y\}$ in $T_pM$ we have $$|K(\sigma)|=c^{2m}|K_{f^m(p)}((T_pf)^m X,\ (T_pf)^m Y|\leq c^{2m}\cdot K,$$ and, for $m\to\9$, we obtain $K_p(X,Y)\equiv 0$ which is a contradiction.

In the case of $M$ irreducible and complete, the hypothesis ``bounded sectional curvature'' can be {\it weakened}. Firstly, we prove

{\bf Lemma 2.1.} {\it Let $(M,g)$ be a complete Riemannian manifold. Every strict homothety (i.e. with $c\neq 1$) of $M$, has a fixed point.}

\ms

\noindent{\it Proof.} $(M,g)$ is a complete metric space with respect to the metric $d(p,q)=\dd\inf_b\left\{\dd\int^1_0 g(b,b)^{\frac12}dt\right\}$ for all curves $b$ on $M$, with $b(0)=p$ and $b(1)=q$, see [5].

Let $f$ be a homothety with $c<1$, otherwise we may take $f^{-1}$. We have $$d(f(p),f(q))=\inf\left\{\dd\int^1_0g(\dot{\stackrel\frown{f\circ b}},\dot{\stackrel\frown{f\circ b}})^{\frac12}dt\right\}\leq c\inf\left\{\dd\int^1_0 g(b,b)^{\frac12}dt\right\}\leq cd(p,q),$$ therefore $f$ is a contraction map. It follows that $f$ has a fixed point.

Now we give the following

{\bf Definition 2.1.} The Riemannian manifold $(M,g)$ is said to be with locally bounded sectional curvature if any $p$ in $M$ admits a closed neighborhood on which the sectional curvature is bounded.

\ms

{\bf Theorem 2.2.} {\it Let $(M,g)$ be a complete and irreducible Riemannian manifold with locally bounded and non-identically zero sectional curvature. Then, the group of affine transformations of $M$ is equal to the group of isometries of $M$.}

\ms

\noindent{\it Proof.} By Step 2 of the proof of Theorem 2.1, every affine transformation $f$ is a homothety and therefore by Lemma 2.1, has a fixed point, denoted by $p_0$. Let $U$ be a closed neighborhood of $p_0$ on which the sectional curvature is bounded by $K<+\9$. Suppose $c<1$ and we have $$d(p_0,f^m(p))=d(f^m(p_0),f^m(p))\leq c^m d(p_0,p),$$ for all $p$ in $M$; hence there exists an $m_0>m$ such that for $m_0>m$, $f^m(p)$ belongs to $U$. From $$|K_p(X,Y)|=c^{2m}|K_{f^m(p)}((T_pf)^mX,(T_pf)^mY)|\leq c^{2m}\cdot K$$ it follows, when $m\to\9$, $K_p(X,Y)\equiv 0$ which is a contradiction.

{\bf Remark 2 .1.} The hypothesis of Theorem 2.1 are satisfied by a $\delta-$pinched Riemannian manifold (i.e. there exists a constant $0<\delta<1$ such that $\delta<K_p<1$, for every $p$ in $M$).

{\bf Remark 2.2.} When $M$ is a finite-dimensional Riemannian manifold, every $p$ in $M$ admits a neighborhood such that $\overline{U}$ is compact. As $K_p$ is a continuous function, it follows that it is bounded; therefore $M$ has locally bounded sectional curvature. In the case when $M$ is complete and irreducible, we obtain the theorem by S. Kobayashi ([2]).

\noindent{\footnotesize{\it Received December 16, 1975}}

\newpage




\renewcommand{\thefootnote}{\fnsymbol{footnote}}

\runningauthor{M. ANASTASIEI}
\runningtitle{STRUCTURES SPINORIELLES SUR LES VARI\'{E}T\'{E}S HILBERTIENNES}

\noindent
\baselineskip 8pt
\noindent{\footnotesize{C.R. Acad. Sc. Paris\hfill\break
t. 284 (25 avril 1977)\hfill\break S\'{e}rie A 943--946}}

\vskip 2cm
\baselineskip 11.5pt plus .15pt
\centerline{\bf\Large ANALYSE MATHEMATIQUE}
\vskip .3cm
\centerline{\bf\large Structures spinorielles sur les vari\'{e}t\'{e}s hilbertiennes}
\vskip .2cm
\centerline{\bf\footnotesize{Note (\footnote{S\'{e}ance du 7 f\'{e}vrier 1977}) de Mihai Anastasiei}}
\vskip .1cm
\centerline{\bf \footnotesize{pr\'{e}sent\'{e} par M. Andr\'{e} Lichnerowicz}}
\vskip 1cm

\begin{abstract}
On donne quelques d\'{e}finitions et propri\'{e}t\'{e}s des structures spinorielles sur les vari\'{e}t\'{e}s model\'{e}es par des espaces de Hilbert.

{\it Some definitions and properties of the spin structures on the manifolds modeled by Hilbert spaces are given.}
\end{abstract}

\setcounter{section}{0}
\section{Introduction} Soit $H$ un espace de Hilbert r\'{e}el, s\'{e}parable et de dimension infinie. Nous notons par $GL(H)$ le group g\'{e}n\'{e}ral lin\'{e}aire de $H$ et par $O(H)$ le group orthogonal de $H$. Soit $P(H)$ un {\it classe de perturbation} pour l'anneau $L(H)$ des op\'{e}rateurs lin\'{e}aires born\'{e}s sur $H$ et soit $GL_p(H)=\{X\in GL(H)$, $X$ congruent \`{a} $I$ modulo $P(H)\}$, ou $I$ d\'{e}signe l'op\'{e}rateur identit\'{e} sur $H$. Le sous-groupe $O(H)_p=O(H)\cap GL_p(H)$ a deux composantes connexes; soit $SO(H)_p$ sa composante connexe de l'identit\'{e}. Si $P(H)$ coincide avec l'id\'{e}al des op\'{e}rateurs nucl\'{e}aires [resp. de Hilbert-Schmidt], les groupes $O(H)_p$, $SO(H)_p$ seront not\'{e}s par $O(H)_1$, $SO(H)_1$ [resp. par $O(H)_2$, $SO(H)_2$].

Pierre de la Harpe a donn\'{e} dans [3] une construction explicite du rev\^{e}tement universel Spin$(H)_1$ de $SO(H)_1$. Ult\'{e}rieurement, R.J. Plymen et R.F. Streater ont donn\'{e} dans [5] la construction explicite du rev\^{e}tement universel Spin$(H)_2$ de $SO(H)_2$. Les groupes Spin$(H)_1$ at Spin$(H)_2$ s'appellent les groupes spinoriels. Posons Spin$(H)=$Spin$(H)_1$ ou bien Spin$(H)_2$ et $SO(H)=SO(H)_1$ ou bien $SO(H)_2$ et notons par $\rho:$Spin$(H)\to SO(H)$, l'homomorphisme correspondant de rev\^{e}tement.

Dans la suite, nous allons d\'{e}finir les structures spinorielles en utilisant les groupes Spin$(H)$ et nous allons donner quelques propri\'{e}t\'{e}s de ces structures. Nous allons consid\'{e}rer aussi fibr\'{e}s de Clifford.

\section{D\'{e}finitions des structures spinorielles} Nous supposerons toujours la diff\'{e}rentiabilit\'{e} de classe $C^\infty$. Soient $M$ une vari\'{e}t\'{e} diff\'{e}rentiable, conexe, model\'{e}e sur un espace de Banach et soit\\ $\pi:E\to M$ un fibr\'{e} vectoriel (en abr\'{e}g\'{e} f.v.) de fibre type $H$.

\ms

{\bf D\'{e}finition 2.1.} Nous appelons $P-$structure riemannienne sur le f.v. $\pi$, une r\'{e}duction du groupe structural de $\pi$ au groupe $O(H)_p$.

\ms

Remarquons que ces structures existent toujours, $GL(H)$ \'{e}tant contractible (le th\'{e}or\`{e}me de Kuiper). Nous les \'{e}tudierons dans un autre travail. Pour le groupe $SO(H)_1$ [resp. $SO(H)_2$] nous obtenons la structure riemannienne nucl\'{e}aire orient\'{e}e [resp. riemannienne de Hilbert--Schmidt orient\'{e}e].

En supposant que le f.v. $\pi$ a une r\'{e}duction au groupe $SO(H)$ soit\\ $P(M,$ $SO(H))$ son fibr\'{e} de rep\`{e}res [1], qui est un fibr\'{e} principal (en abr\'{e}g\'{e} f.p.) de base $M$ et de groupe structural $SO(H)$.

\ms

{\bf D\'{e}finition 2.2.} Une structure spinorielle sur le f.v. $\pi$ avec une r\'{e}duction au groupe $SO(H)$ (ou sur le f.p. $P(M,SO(H))$) est une extension du f.p. $P(M,SO(H))$ associ\'{e}e \`{a} l'homomorphisme de rev\^{e}tement $\rho:$Spin$(H)\to SO(H)$.

\ms

Nous notons par $\Sigma(M,{\rm{Spin}})$ une telle extension et par $$\wt{\rho}:\Sigma(M,{\rm{Spin}}(H))\to P(M,SO(H))$$ l'homomorphisme qui correspond \`{a} l'homomorphisme $\rho$.

\ms

{\bf Remarque 2.1.} La d\'{e}finition 2.2 est \'{e}quivalente \`{a} la d\'{e}finition donn\'{e}e par A. Lichnerowicz [6].

\ms

{\bf Remarque 2.2.} Comme le f. p. $P(M,SO(H))$ est d\'{e}termin\'{e}, \`{a} un isomorphisme pr\`{e}s, par un recouvrement ouvert $\{U_i\}$ et un cocycle $g_{ij}:U_i\cap U_j\to SO(H)$ (1), le f. p. $\sum(M, {\rm{Spin}}(H))$ (s'il existe) est d\'{e}termin\'{e} par un cocycle $\wt{g}_{ij}: U_i\cap U_j\to {\rm{Spin}}(H)$ tel que $\rho(\wt{g}_{ij})=g_{ij}$. Cette remarque, r\'{e}unie avec la possibilit\'{e} d'identifier la classe d'isomorphie du f. p. $P(M,SO (H))$ avec un \'{e}l\'{e}ment de l'ensemble de cohomologie $H^1(M,SO(H))$ est tr\`{e}s utile.

\ms

{\bf Th\'{e}or\`{e}me 2.1.} {\it Le f. p. $P(M,SO(H))$ admet une structure spinorielle si et seulement s'il existe un \'{e}lement non nul $\sigma\in H^1(P,Z_2)$ tel que $\sigma$ restreint a chaque fibre soit non trivial.}

\ms

\noindent{\it Esquisse de preuve.} Nous consid\'{e}rons\\ $\sigma$ comme un homomorphisme $\sigma:H_1(P)\to Z_2$ et nous d\'{e}finissons un homomorphisme $\sigma\circ \varphi_1:\pi_1(P)\to Z_2$ ou $\varphi_1$ est l'homomorphisme de Hurewicz. Donc, ker$(\sigma\circ\varphi_1)$ est un sous-groupe d'ordre deux dans $\pi_1(P)$. Comme $P$ est localement contractible, il existe un rev\^{e}tement d'ordre deux $\Sigma$ de $P$, qui est l'espace total d'une extension de $P(M,SO (H))$ associ\'{e}e \`{a} $\rho$. La necessit\'{e} est imm\'{e}diate.

La d\'{e}finition d'une structure spinorielle qui decoule du th\'{e}or\`{e}me 2.1 est tr\`{e}s utile pour les d\'{e}monstrations des th\'{e}or\`{e}mes suivantes [voir [7] pour la dimension finie ou pour le cas topologique].

\ms

{\bf Th\'{e}or\`{e}me 2.2.} {\it L'ensemble des structures spinorielle, s'il n'est pas vide, est en bijection modulo l'isomorphisme de fibr\'{e}s principaux avec $H^1(M,Z_2)$.}

\ms

{\bf Th\'{e}or\`{e}me 2.3.} {\it Soient les f.v. $\pi_1$, et $\pi_2$, avec l'espace de base $M$ et $\pi_1\oplus\pi_2$, leur somme de Whitney. Si deux de ces fibres ont des structures spinorielles, le troisi\`{e}me est muni aussi d'une structure spinorielle.}

\ms

La suite exacte $$1\to Z_2\to{\rm{Spin}}(H)\to SO(H)\to 1$$ induit, une suite exacte de groupes et d'ensemble de cohomologie [4]:
$$H^1(M,Z_2)\longrightarrow H^1(M,{\rm{Spin}}(H))\longrightarrow H^1(M,SO(H))\underset{v}{\longrightarrow}H^2(M,Z_2).$$

Il existe une structure spinorielle sur $P(M,SO(H))$ si et seutlement si $v(P)=0$. Dans le cas $SO(H)=SO(H)_1$, $v(P)=w_2(P)$, $[3]$, la deuxi\`{e}me classe de Stiefel-Whitney de $P(M,SO(H)_1)$.

Soit $M'$ une autre vari\'{e}t\'{e} et soit $f:M'\to M$ un morphisme. Nous notons par $Pf$ le f. p. sur $M'$ induit par $P(M,SO(H))$ et $f$.  Si le f. p. $P(M,SO (H))$ admet une structure spinorielle, alors $Pf$ admet une structure spinorielle. Une propri\'{e}t\'{e} r\'{e}ciproque est donn\'{e}e par le

\ms

{\bf Th\'{e}or\`{e}me 2.4.} {\it Soit $f:M'\to M$ un f. p. de groupe structural $G$. Alors, $G$ op\`{e}re naturellement sur $Pf$ et soit $u'\cdot G$ l'orbite de $u'\in Pf$. Si le f. p. $Pf$ admet une structure spinorielle $\Sigma(M',{\rm{Spin}}(H))$ et si le groupe $G$ op\`{e}re sur $\Sigma$ tel que la projection $\Sigma\to\Sigma/G$ est un f. p. de groupe structural $G$ et $\wt{\rho}(u\cdot G)=\wt{\rho}(u)\cdot G$, o\`{u} $u\cdot G$ est l'orbite de $u\in\Sigma$, alors $P(M,SO(H))$ admet une structure spinorielle.}

\ms

\noindent{\it D\'{e}monstration.} La vari\'{e}t\'{e} $\Sigma/G$ est muni d'une structure naturelle de f.p. de base $M$ et de groupe structural Spin$(H)$ par l'action $(u\cdot G)a=ua\cdot G$, pour $a\in$Spin$(H)$, et projection $u\cdot G\to f(\a(u))$ avec $\a$ la projection $\Sigma\to M'$. L'homomorphisme canonique $\Sigma/G\to P$ est de la forme $u\cdot G\to f(f^*(u))$, o\`{u} $f^*:Pf\to P$ est l'homomorphisme induit par $f$.

Dans le cas $G=Z_2$, nous obtenons un r\'{e}sultat qui a \'{e}t\'{e} d\'{e}montr\'{e} par I. Popovici [9] pour la dimension finie et le cas non orientable. Une structure spinorielle sur $M$ sera par d\'{e}finition une structure spinorielle sur le fibr\'{e} tangent $TM$ (en supposant que $M$ est mod\'{e}l\'{e}e sur l'espace de Hilbert $H$).

\ms

{\it Exemples.} (a) Le f. p. trivial $M\times SO(H)$ admet toujours une structure spinorielle, unique si $M$ est simplement connexe.

(b) Soit l'espace de Hilbert: $$\mathcal{1}_2=\{x=(x_1,x_2,...)\ |\ x_i\in\R,\dd\sum^{\9}_{i=1}x^2_i<\9\},$$ avec la base $c_1,c_2,...$ Le tore hilbertien $T=\mathcal{1}_2/\dd\sum^\9_{i=1}\Z c_1$ est un groupe de Lie-Hilbert. Chaque structure riemannienne nucl\'{e}aire sur $T$ est orientable car $w_1(T)=0$ et il existe une structure spinorielle [r\'{e}latif \`{a} Spin$(\mathcal{1}_2)_1$] sur $T$, car $w_2(T)=0$ [5].

(c) Des exemples plus sofistiqu\'{e}s d\'{e}coulent de [5].

\section{Fibr\'{e}s de Clifford}\setcounter{equation}{0} Maintenant, nous allons limiter nos consid\'{e}rations au groupe $SO(H)_1$. Soit $C(H)_1$ l'alg\`{e}bre de Clifford avec sa structure de $C^*-$alg\`{e}bre [3]. Chaque $*-$au\-to\-mor\-phisme $\psi$ de $C(H)_1$ pour lequel $\psi(H)=H$, s'appelle de Bogoliubov. Le groupe de ces automorphismes est not\'{e} par Bog$(C(H)_1)$. Il existe un isomorphisme $C:O(H)\to$ Bog$(C(H)_1)$ qui est aussi un isomorphisme de groupes de Lie-Banach.

Soit un f.p. $P(M,O(H))$. Son fibr\'{e} associ\'{e} avec la fibre $C(H)_1$ est un fibr\'{e} en alg\`{e}bres localement trivial, nomm\'{e} fibr\'{e} de Clifford. Nous notons par Bog$(C(H)_1)_i$ les automorphismes de Bogoliubov int\'{e}rieurs et par Bog$(C(H)_1)_{ip}$ les automorphismes de Bogoliubov int\'{e}rieures et paires.

L'isomorphisme $O(H)_1\simeq$ Bog$(C(H)_1)_i$ implique le:

\ms

{\bf Th\'{e}or\`{e}me 3.1.} {\it Il existe une structure riemannienne nucl\'{e}aire sur f.v. riemannienne si et seulement si le fibr\'{e} de Clifford admet une r\'{e}duction au groupe ${\rm{Bog}}(C(H)_1)_i$.}

\ms

L'isomorphisme $SO(H)_1\simeq$ Bog$(C(H)_1)_{ip}$ [3] implique le

\ms

{\bf Th\'{e}or\`{e}me 3.2.} {\it Une structure riemannienne nucl\'{e}aire sur un f.v. riemannienne est orientable si et seulement si le fibr\'{e} de Clifford admet une r\'{e}duction au groupe ${\rm{Bog}}(C(H)_1)_{ip}$.}

\ms

Dans le contexte des structures spinorielles, on peut \'{e}tudier les champs spinoriels et les connexions spinorielles en dimension infinie. Ces aspects seront abord\'{e}s dans un autre travail.

\ms

I. Pop et I. Popovici ont tr\`{e}s utilement discut\'{e} sur ce travail.

\ms
\noindent {\bf BIBLIOGRAPHIE}
\begin{enumerate}
\item{ N. Bourbaki}, {\it Vari\'{e}t\'{e} diff\'{e}rentielles et analytiques}, Fascicule des r\'{e}sultats, Hermann, Paris, 1967.
\item{ K.D. Elworthy et A.J. Tromba}, {\it Proc. Sym. in Pure Math.}, XV, A.M.S., 1970, p. 45--94.
\item{ P. de la Harpe}, {\it Composition Math.}, 25, 1972, p. 245--261.
\item{ F. Hirzebruch}, {\it Topological Methods in Algebric Geometry}, Springer, 1966, p. 33.
\item{ U. Koschorke}, {\it Proc. Sym. in Pure Math.}, XV, A.M.S., 1970, p. 95-133.
\item{ A. Lichnerowicz}, {\it Bull. Soc. Math. Fr.}, 92, 1964, p. 11--100.
\item{ J. Milnor}, {\it Differential and Combinatorial Topology} in {\it Honor of M. Morse}, Princeton Univ. Press, 1965, p. 55--62.
\item{ R.J. Plymen et R.S. Streater}, {\it Bull. London Math. Soc.}, 7, 1975, p. 283--288.
\item{ I. Popovici}, {\it Rend. Circ. Mat. Palermo}, 23, 1974, p. 113-134.

\end{enumerate}

\newpage




\def\bfe{{\mathbf{E}}}

\runningauthor={M. ANASTASIEI}
\runningtitle={RIEMANNIAN P-STRUCTURES ON VECTOR BUNDLE}
\noindent
\baselineskip 8pt
\noindent{\footnotesize{Proceedings of the National Colloquium on\hfill\break
Geometry and Topology, Timi\c{s}oara (Romania)\hfill\break April 26--28, 1977}}

\vskip 2cm
\baselineskip 11.5pt plus .15pt
\centerline{\bf\Large RIEMANNIAN P-STRUCTURES}
\vskip .2cm
\centerline{\bf\Large ON VECTOR BUNDLE}
\vskip .5cm
\centerline{\bf\footnotesize{BY}}
\vskip .5cm
\centerline{\bf \footnotesize{M. ANASTASIEI}}
\vskip 1cm

\setcounter{section}{0}
\section*{Introduction} Let $H$ be a separable, real Hilbert space. Denote by $L(H)$ the algebra of bounded linear operators on $H$ and by $\phi(H)$ the set of Fredholm operators on $H$. A two-sided ideal $P(H)$ of $L(H)$ is said to be a $\phi$-perturbation class if $\phi(H) + P(H) =\phi(H)$ and $F(H)\subset P(H)$, where $F(H)$ is the two-sided ideal of the finite rank operators on $H$. Let $GL_P(H)$ be the group of those invertible operators on $H$, which can be written as $I+X$ with $X$ in $P(H)$, where $I$ is the identity operator. Denote by $O(H)$ the group of orthogonal operators on $H$ and we put $O(H)_P = GL_P(H)\cap O(H)$. The group $O(H)_P$ has two connected components. Denote by $SO(H)_P$ the connected component of $I$.

Let $M$ be a connected and paracompact manifold, locally diffeomorphic to a Banach space, and let $\pi : E\to M$ be a vector bundle having $H$ as the type fibre. A $P-$structure on $\pi$ is a reduction of its structural group to $GL_P(H)$. A vector bundle with a $P-$structure is called a $P-$bundle (see [3]). A reduction of the structural group of $\pi$ to the group $O(H)_P$ will be called a Riemannian $P-$structure and a vector bundle with a Riemannian $P-$structure will be called a $PR-$bundle.

We are going to discuss the reduction of a $P-$bundle to a $PR-$bundle and to describe the morphisms of the $PR-$bundles.

\section{On the reduction of a $P$-bundle to a $PR$-bundle} Let $\pi:E\to M$ be a vector bundle, as above. One says that $\pi$ admits a reduction of its structural group to a subgroup $G$ of $GL(H)$ if there exists a maximal collection of trivializations $(U_j,\phi_j)_{j\in J}$ with $U_j$ open in $M$ and $$\phi_j:\pi^{-1}(U_j)\to U_j\times H$$ such that the maps $$\phi_k\circ\phi^{-1}_j:U_j\cap U_k\to GL(H)$$ take their values in $G$.

Let $g$ be a Riemannian metric on $\pi$. Using	$\phi_j$ we can transport the restriction of the Riemannian metric $g$ to $\pi^{-1}(U_j)$ on $U_j\times H$ and for a fixed point $p$ in $U_j$ we obtain a symmetric, bilinear and positive defined form on $H$ whose corresponding operator (symmetric and positive) will be denoted by $A_{jp}$. The map $U_j\to L(H)$ given by $p\to A_{jp}$ is a morphism.

\ms

{\bf Definition 1.1.} Let $\pi:E\to M$ be a vector bundle with a $P-$structure. A Riemannian metric $g$ on $\pi$ is said to be adapted to the $P-$structure of $\pi$ if for every trivialization $(U_j,\phi_j)$ and for every $p$ in $U_j$, there exists $X_{jp}$ in $P(H)$ so that $A_{jp}=I+X_{jp}$.

\ms

{\bf Theorem 1.1.} {\it Let $\pi:E\to M$ be a vector bundle with a $P-$structure. Then $\pi$ admits a Riemannian $P-$structure if there exists a Riemannian metric $g$ adapted to its $P-$structure.}

\ms

\noindent{\it Proof.} Let $(U_j,\phi_j)_{j\in J}$ be the maximal collection of trivialization of $\pi$, such that $\phi_k\circ\phi_j^{-1}$ are $GL(H)-$valued. Denote by $g_j$, the metric on $U_j\times H$ obtained from the restriction of $g$ to $\pi^{-1}(U_j)$ and for every $p$ in $U_j$ we put	 $g_{jp}(v,w)=(A_{jp}v,w)$, where $(\ ,\ )$ is the inner product on $H$. Let $A$ be in $L(H)$; we agree to note by $\sqrt{A}$ an operator $B=\lim_n B_n$ where $B_n$ is a sequence inductively defined by $$B_{n+1}=\dfrac12(B_n+B^{-1}_n A),\ \ B_1 = I.$$ We define new trivializations for $\pi$ by $\Psi_{jp}=B_{jp}\circ\phi_{jp}$, where $B_{jp}=\sqrt{A_{jp}}$ and $\phi_{jp}=\phi_j |_{\pi^{-1}(p)}$. Since for every $v$ and $w$ in $H$, we have $$(B_{jp}v,B_{jp}w)=(B^2_{jp}v,w)=(A_{jp}v,w)=g_{jp}(v,w),$$ $B_{jp}$ is an isometric map with respect to inner product on $H$ and $g_{jp}$, hence $(\psi_k\circ \psi^{-1}_j)$ $(p)\in O(H)$. Since $A_{jp}=1+X_{jp}$ with $X_{jp}$ in $P(H)$, it is easy to see that $\sqrt{A_{jp}}=I+Y_{jp}$ with $Y_{jp}$ in $P(H)$.

Using this expression of $\sqrt{A_{jp}}$, we obtain: $$\begin{array}{l}(\psi_k\circ \psi_j^{-1})(p)=B_{kp}\circ\phi_{kp}\circ\phi^{-1}_{jp}\circ B^{-1}_{jp}=B_{kp}(I+Z_p)B^{-1}_{jp}=\\ \\
= B_{kp}\circ B^{-1}_{jp}+B_{kp}\circ Z_p\circ B^{-1}_{jp}=(I+ Y_{kp})(I+V_{jp})+B_{kp}\circ Z_p\circ B^{-1}_{jp}=\\ \\
=I+X_p\end{array}$$ with $X_p$ in $P(H)$, since $P(H)$ is a two-sided ideal. Therefore, $(U_j,\psi_j)_{i\in J}$ is the expected collection of trivializations of $\pi$.

\ms

{\bf Remark 1.1.} By the Theorem 1.1, a $P-$structure on $\pi$ and a Riemannian metric adapted to it, determine a Riemannian $P-$structure. Conversely, a Riemannian $P-$structure determines a $P-$structure (itself viewed as $P-$structure) and a Riemannian metric by $g_p(\xi,\eta)=(\phi_{jp}\xi,\phi_{jp}\eta)$ for $p$ in $M$ (the definition is correct because $\phi_{kp}\circ\phi^{-1}_{jp}\in O(H)$) whose associated operators $A_{jp}$ are all equal to I.

\ms

{\bf Definition 1.2.} A $PR-$bundle is orientable if it admits a reduction to the group $SO(H)_p$.

\ms

We have proved a criterion for the orientability of a $PR-$bundle in [l].

Consider for $F(H)$ the following $q-$norm ($1\leq q\leq\9$): $\|X\|_q=$ (trace $(\sqrt{X^*X})^q)^{1/q}$	for $1\leq q<\9$ and the usual norm for $q=\9$. The closure of $F(H)$ in this $q-$norm will be denoted by $F_q(H)$. Each set $F_q(H)$ is a $\phi-$perturbation class and it corresponds to it a group denoted by $O(H)_q$. Some structures of great importance in the study of the spin structures on Hilbert manifolds are the reductions of the structural group of a vector bundle to $O(H)_1$, $O(H)_2$ respectively, named in [2], Riemannian nuclear structure, Riemannian Hilbert-Schmidt structure respectively. In the same Note there exists a condition for the reduction of a Riemannian vector bundle to a nuclear vector bundle; another criterion for the orientability of a Riemannian nuclear vector bundle is also given.

\section{Morphisms of $PR-$bundles}
\setcounter{equation}{0}

Let $\phi_0(H)$ be the set of Fredholm operators of index 0.

\ms

{\bf Lemma 2.1.} [4] {\it Every $T\in \phi_0(H)$ can be written as: $T=S+a$, where $S\in GL(H)$ and $a\in F(H)$.}

\ms

{\bf Definition 2.1.} An operator $T\in\phi_0(H)$ is said to be an $O\phi_0-$operator if $a^*S+S^*a+a^*a=0$, where $S$ and $a$ are as in Lemma 2.1.

Let $\pi'$ be another vector bundle over $M$ having the type fibre $H$. A morphism $f:\pi\to\pi'$ is called a $\phi_0-$morphism if it is a $\phi_0-$operator on each fibre.

\ms

{\bf Definition 2.2.} Let $\pi'$ be a Riemannian vector bundle. A morphism $f:\pi\to\pi'$ will be called an $O\phi_0-$morphism if for every trivialization $(U,\phi)$ and $(V,\psi)$ with $f(U)\subset V$ of $\pi$ and $\pi'$ respectively, we have $$(\psi\circ f\circ\phi^{-1})(x,v)=(f(x),f_1(x)v)$$ with $f_1(x)$ an $O\phi_0-$operator for every $x$ in $U$.

\ms

{\bf Theorem 2.1.} {\it Let $\pi$ be a vector bundle and let $\pi'$ be a $PR-$bundle. An $O\phi_0-$morphism $f:\pi\to\pi'$ induces a unique Riemannian $P-$structure on $\pi$, such that one has $$\psi\circ f\circ \phi^{-1}(x,v)=(f(x),v+a_0(x)v)$$ with $a_0(x)$ in $F(H)$ and $I+a_0(x)$ in $O(H)$, whenever $(U,\phi)$ and $(V,\psi)$ with $f(U)\subset V$ are the trivializations of these Riemannian $P-$structures.}

\ms

\noindent {\it Proof.} Let be $(U,\phi)$ a trivialization of $\pi$ and $(V_0,\psi_0)$ with $V_0\subset f(U)$ a trivialization of $\pi'$. Therefore we have $$(\psi\circ f\circ\phi^{-1})(x,v)=(f(x),f_1(x)v),$$ where $f_1(x)$ is an $O\phi_0-$operator for every $x$ in $U$. By Lemma 2.1, $f_1(x)=S(x)+a(x)$, where $S(x)\in GL(H)$ and $a(x)\in F(H)$. Since $GL(H)$ is an open subset of $L(H)$, there exists a neighborhood $U_0$ of $x$, such that $S(U_0)\subset GL(H)$.

Define a new trivialization of $\pi$, $\phi_0:\pi^{-1}(U_0)\to U_0\times H$ by $$\phi_0\circ\phi^{-1}(x,v)=(f(x),S(x)v).$$ We have, using the definition of the $O\phi_0-$operators, $$\begin{array}{l}\psi_0\circ f\circ\phi^{-1}_0(x,v)=\psi_0\circ f\circ\phi^{-1}\circ\phi\circ\phi^{-1}_0(x,v)=\\ \\
=\psi_0\circ f\circ\phi^{-1}(x,S^{-1}(x)v)=(f(x),f_1(x)S^{-1}(x)v)=\\ \\
= (f(x),(S(x)+a(x))S^{-1}(x)v)=(f(x),v+a(x)S^{-1}(x)v)=\\ \\ = (f(x),v+a_0(x)v),\end{array}$$ with $a_0(x)$ in $F(H)$ and $I+a_0(x)$ in $O(H)$.

Let	$(V_1,\psi_1)$ be another trivialization of $\pi'$ and $(U_1,\phi_1)$ the trivialization of $\pi$ associated to it by the above construction. Therefore we have $$\psi_1\circ f\circ\phi^{-1}_1(x,v)=(f(x),v+a_1(x)v)$$ with $a_1(x)$ in $F(H)$ and $I+a_1(x)$ in $O(H)$. We put $\psi_1\circ\psi^{-1}_0(x,v)=(f(x),B'(x)v)$ with $B'(x)$ in $O(H)_p$ and $\phi_0\circ\phi^{-1}_1(x,v)=(f(x),B(x)v)$. From $\psi_1\circ f\circ\phi^{-1}_1=\psi_1\circ\psi^{-1}_0\circ\psi_0\circ f\circ\phi^{-1}_0\circ\phi\circ\phi^{-1}_1$ it follows
$$B'(x)\circ B(x)v+B'(x)a_0(x)B(x)v=v+a_1(x)v$$ hence, if we omit $x$, $$B'B+B'a_0 B=I+a_1\leqno(*)$$ or equivalently $$B'(I+a_0)B=I+a_1.\leqno(**)$$ Therefore $B$ is an orthogonal operator.

If we put $B'=I+b'$, from $(**)$ it follows $$B=I+a_1-b'\ B-a_0B-b'a_0B=I+a$$ with $a$ in $F(H)$ since $F(H)$ is a two-sided ideal. Hence $(\phi_0\circ\phi_1^{-1})(x)\in O(H)_p$ and the proof is complete.

\ms

{\bf Corollary 2.1.} {\it Let $\pi:E\to M$ be a vector bundle with fibre $H$. An $O\phi_0-$morphism $f:E\to M\times H$ induces a unique Riemannian $P-$structure on $E$, so that for any trivialization $(U,\phi)$ of $E$, we have $$f\circ \phi^{-1}(x,v)=(f(x), v+a(x)v),$$ with $a(x)$ in $F(H)$, and $I+a(x)$ in $O(H)$.}

\ms

\noindent{\it Proof.} One considers the trivial Riemannian $P-$structure on $M\times H$ and one applies the Theorem 2.1.

\ms

{\bf Corollary 2.2.} {\it Let $f:\pi\to\pi'$ be an isomorphism of vector bundles. If $\pi'$ admits an (orientable) Riemannian $P-$structure, there exists a unique (orientable) Riemannian $P-$structure on $\pi$ such that for every trivialization $(U,\phi)$ and $(V,\psi)$ with $f(U)\subset V$ of these Riemannian $P-$structures, we have $$(\psi\circ f\circ \phi^{-1})(x,v)=(x,v).$$}

\ms

\noindent{\it Proof.} One repeats the proof of Theorem 2.1 with $a(x)=0$ because $f$ is an isomorphism. The relation $(*)$ becomes $B'B=I$, hence $B'\in O(H)_p$ (resp. $SO(H)_p$) implies $B\in O(H)_p$ (resp. $SO(H)_p$).

\ms

{\bf Corollary 2.3.} {\it Let $N$ and $N'$ be manifolds modeled by $H$. Suppose that $N'$ admits an (orientable) Riemannian $P-$structure. A diffeomorphism $h:N\to N'$ induces an (orientable) Riemannian $P-$structure on $N$.}

\ms

{\bf Remark 2.1.} Let $f:\pi\to\pi'$ be an $O\phi_0-$morphism, where $\pi'$ is a $PR-$bundle. Suppose that the Riemannian $P-$structure of $\pi'$ is obtained from a $P-$structure and a Riemannian metric $g$. The Riemannian $P-$structure induced by $f$ on $\pi$ is not obtained from the $P-$structure induced by $f$ and $f^*g$, since $f^*g$ is not a Riemannian metric. However this happens in the context of the Corollary 2.2.

\bigskip

{\bf REFERENCES}

\begin{enumerate}
  \item {{M. Anastasiei}}, {\it Spin structures on Hilbert manifolds} (to appear in the Proceedings of the Institute of Mathematics Iasi, 1977).
  \item {{M. Anastasiei}}, {\it Structures spinorielles sur les vari\'{e}t\'{e}s hilbertiennes} (to appear in C.R. Acad. Sci. Paris, S\'{e}r. A, 1977).
  \item {{K.D. Elworthy}} and {{A.J. Tromba}}, {\it Differential structures and Fredholm maps on Banach manifolds}, Proc. Sym. in Pure Math., XV, (1970), 45-94.
  \item {{N. Moulis}}, {\it Structures de Fredholm sur les vari\'{e}t\'{e}s hilbertiennes}, Lecture. Notes in Math., 259, Springer-Verlag, 1972, p.4.
\end{enumerate}

\bigskip

\begin{flushleft}
    \footnotesize{Mathematical Seminar ``A. Myller'' Iasi (Romania)}
\end{flushleft}


\newpage





\def\G{\Gamma}

\runningauthor={M. ANASTASIEI}
\runningtitle={SPIN STRUCTURES ON HILBERT MANIFOLDS}
\noindent
\baselineskip 8pt
\noindent{\footnotesize{Analele \c{S}tiin\c{t}ifice ale Universit\u{a}\c{t}ii ``AI. I. Cuza'' Ia\c{s}i\hfill\break
Tomul XXIV, f. 2, s. Ia, 1978, 367-373}}
\vskip 2cm
\baselineskip 11.5pt plus .15pt
\centerline{\bf\Large SPIN STRUCTURES ON}
\vskip .2cm
\centerline{\bf\Large HILBERT MANIFOLDS}
\vskip .5cm
\centerline{\bf \footnotesize{BY}}
\vskip .5cm
\centerline{\bf \footnotesize{M. ANASTASIEI}}
\vskip 1cm

\setcounter{section}{0}
\section{Introduction}

Let $H$ be a separable, real Hilbert space. We denote by $L(H)$ the algebra of bounded linear operators on $H$, by $O(H)$ the orthogonal operators on $H$ and by $I$ the identity operator on $H$. Let $GL_P(H)$ be the group of those invertible operators on $H$ which can be written as $I+A$, where $A$ belongs to a ``perturbation class'' $P(H)$ [2, p. 46] of $L(H)$. The group $O(H)_P=O(H)\cap GL_P(H)$ is doubly connected and we denote by $SO(H)_P$ the connected component of $I$.

Let $F(H)$ be the ideal of finite rank operators on $H$. We denote by $F(H)_p$ the closure of the ideal $F(H)$ in the $p-$norm defined by $\|X\|_p=$(trace$(X^*X)^{\frac{p}{2}})^{\frac{1}{p}}$ for $1\leq p<\9$ and by the usual norm, for $p=\9$. The ideals $F(H)_p$ are ``perturbation classes'' for $L(H)$. For $p=1$ (resp. $p=2$) we obtain the ideal of nuclear operators (resp. the ideal of Hilbert-Schmid operators) and in this case the groups $O(H)_P$, $SO(H)_P$ will be denoted by $O(H)_1$, $SO(H)_1$ (resp. $O(H)_2$, $SO(H)_2$). It follows, from general principles, that the universal covering of $SO(H)_P$ is a Banach-Lie group and that the covering map is 2-sheeted. An explicit construction of the universal covering group Spin$(H)_1$ of $SO(H)_1$ has been given by P. de la Harpe [3]. Later, R.J. Plymen and R.F. Streater [9] gave an explicit construction of the universal covering group Spin$(H)_2$ of $SO(H)_2$. Both groups Spin$(H)_1$ and Spin$(H)_2$ will be called spinor groups and will be denoted by Spin$(H)$.

We denote by $SO(H)$ both the groups $SO(H)_1$ and $SO(H)_2$, and by $\rho:$Spin$(H)\to SO(H)$ the corresponding covering maps. In the following, we define the spin structures using the groups Spin$(H)$ and we give some properties of these structures. Some results about the Riemannian $P-$structures are given, too.

\section{Definitions of the spin structures}

All bundles, manifolds and morphisms considered in the following will be assumed of class $C^\9$. Let $M$ be a connected and paracompact manifold, modeled on a Banach space and let $\xi:E\to M$ be a vector bundle over $M$ with fibre $H$.

\ms

{\bf Definition 2.1.} A Riemannian $P-$structure on the vector bundle $\xi$ is a reduction of its structural group to the group $O(H)_{P}$.

\ms

{\bf Remark 2.1.} The existence of the Riemannian $P-$structures is a direct consequence of the fact that $GL(H)$ is contractible (Kuiper's theorem).

\ms

{\bf Definition 2.2.} A Riemannian $P-$structure on the vector bundle $\xi$ is said to be orientable if $\xi$ admits a further reduction of its structural group to the group $SO(H)_{P}$.

\ms

{\bf Theorem 2.1.} {\it A Riemannian $P-$structure on the vector bundle $\xi$ is orientable if and only if the first Stiefel-Whitney class $w_1(\xi)$, vanishes.}

\ms

\noindent{\it Proof.} The proof of proposition 6.2 from [6] can be repeated using the homomorphism $O(H)_P\to O(H)_P/SO(H)_P$. We give an alternative proof. The exact sequence $$1\to SO(H)_P\to O(H)_P\overset{p}{\to} Z_2\to 1\leqno(2.1)$$ induces an exact sequence of the cohomology groups ([5]) $$0\to H^1(M,SO(H)_P)\to H^1(M,O(H)_P)\overset{p^*}{\to}H^1(M,Z_2).\leqno(2.2)$$

We denote by $L$ the principal bundle of linear frames of $\xi$ (for definition see Bourbaki [l]), thought as an element of $H^1(M,O(H)_P)$. From exactness of the sequence (2.2) it follows that the {\it Riemannian $P-$structure of $\xi$ is orientable iff} $p^*(L)=0$. Now we prove that $p^*(L)=w_1(\xi)$. By the naturally property of the characteristic classes, it is sufficient to do so when $M$ is the classifying space $BO$ of the group $O(H)_P$. But $H^1(BO,Z_2)=Z_2$  [6], hence the map $p^*$ is either  identically zero, or is the class $w_1$. The first alternative is not possible, because there exists at least a vector bundle with a non-orientable Riemannian $P-$structure (see Exemple 4 from [2]). For the theory of the characteristic classes considered here, see U. Koschorke [6].

\ms

{\bf Corollary 2.1.} {\it A connected and paracompact manifold $N$, modeled on $H$, endowed with a Riemannian $P-$structure is orientable with respect to this structure iff $w_1(N)=0$.}

\ms

{\bf Theorem 2.2.} {\it A connected and paracompact manifold $N$, modeled on $H$, is orientable with respect to all Riemannian $P-$structures which are compatible with its manifold structure if $H^1(N,Z_2)=0$.}

\ms

\noindent{\it Proof.} It follows from a result of U. Koschorke [6, Proposition 6.3].

A reduction of the structural group of the vector bundle $\xi$ to the group $O(H)_1$ (resp. $O(H)_2)$ will be called a {\it Riemannian nuclear structure} (resp. a {\it Riemannian Hilbert-Schmidt structure}).

Let $G$ be a Banach-Lie group and let $P(M,\pi,G)$ (where $\pi:P\to M$), be a principal bundle over $M$ with group $G$. Let $G'$ be another Banach-Lie group.

\ms

{\bf Definition 2.3.} A principal bundle $P'(M,\pi',G')$ where $\pi':P'\to M$ is said to be an extension of the principal bundle $P(M,\pi,G)$, associated to the homomorphism $\varphi:G'\to G$ if there exists a morphism $\wt{\varphi}:P'\to P$ such that $(\wt{\varphi},\varphi)$ is a  morphism of principal bundles.
We suppose that the vector bundle $\xi$ is endowed with a reduction of its structural group to $SO(H)$ and we denote by $P(M,\pi,SO(H))$ the principal bundle of linear frames of it.

\ms

{\bf Definition 2.4.} A spin structure on the vector bundle $\xi$, endowed with a reduction of its structural group to $SO(H)$, is an extension of principal bundle $P(M,\pi,SO(H))$, associated to the covering map $\rho:$Spin$(H)\to SO(H)$.

\ms

Such an extension will be denoted by $\Sigma(M,\pi$,Spin$(H))$ and will be called a {\it spin structure on} $P(M,\pi,SO(H))$ or a {\it spin structure on $M$ with respect to} $P(M,\pi,SO(H))$, too.

The morphism $\wt{\rho}:\Sigma\to P$ is a 2-sheeted covering map and its restriction to fibres are 2-sheeted covering maps. For every $a$ in Spin$(H)$ and $u$ in $\Sigma$, we have $\wt{\rho}(ua)=\wt{\rho}(u)\rho(a)$ and $\pi(\wt{\rho}(u))=\pi'(u)$. It follows that the Definition 2.4 is equivalent to the following definition, given by A. Lichnerowicz [7] in a different context.

\ms

{\bf Definition 2.5.} A spin structure on the vector bundle $\xi$, endowed with a reduction of its structural group to $SO(H)$, is a principal bundle $\Sigma(M,\pi,{\rm{Spin}}(H))$ such that $\Sigma$ is a 2-fold covering of $P$, the restriction of the covering map $\wt{\rho}:\Sigma\to P$ to fibres are 2-sheeted covering maps and $\wt{\rho}(ua)=\wt{\rho}(u)\rho(a)$, $\pi(\wt{\rho}(u))=\pi'(u)$ hold, for every $a\in$Spin$(H)$ and $u\in \Sigma$.

\ms

{\bf Remark 2.2.} By a general result (see Bourbaki [1]), the principal bundle $P(M,\pi,SO(H))$ is determined (up to an isomorphism) by an open covering $\{U_i\}$ of $M$ and a cocycle $g_{ij}:U_i\cap U_j\to SO(H)$. From Definition 2.4 it follows that, the principal bundle $\Sigma(M,\pi$,Spin$(H))$, when it exists, is determined by a cocycle $\wt{g}_{ij}:U_i\cap U_j\to$Spin$(H)$ such that $\rho(\wt{g}_{ij})=g_{ij}$.

The following theorem holds.

\ms

{\bf Theorem 2.3.} {\it Let $P(M,\pi,SO(H))$ be the principal bundle of linear frames of $\xi$. The vector bundle $\xi$ admits a spin structure if and only if there exists a cohomology class $\sigma\in H^1(P,Z_2)$ whose restriction to each fibre is non-zero.}

\noindent{\it Proof.} From the isomorphism $H^1(P,Z_2)\simeq$Hom$(H_1(P),Z_2)$ it follows that, there exists a not trivial homomorphism $\sigma:H_1(P)\to Z_2$. If $\varphi_1:\pi(P)\to H_1(P)$ denotes Hurewicz's homomorphism, $\sigma\circ\varphi_1:\pi_1(P)\to Z_2$ is an epimorphism, hence ker$(\sigma\circ\varphi_1)$ is a subgroup of index 2 in $\pi_1(P)$. Consequently, since the manifold $P$ is locally contractible, there exists a covering space $\Sigma$ of $P$ such that $\wt{\rho}_*(\pi_1(\Sigma))=$ker$(\sigma\circ\varphi_1)$, where $\wt{\rho}:\Sigma\to P$ is the covering map. The covering space $\Sigma$ can be taken as the total space of a principal bundle over $M$ with group Spin$(H)$ which is an extension of $P(M,\pi,SO(H))$, therefore a spin structure of $\xi$.

Conversely, if $P(M,\pi,SO(H))$ admits a spin structure, $\Sigma(M,\wt{\pi}$,Spin$(H))$, the total space $\Sigma$ is a two-fold covering of $P$. Let $s_0$ and $s_1$ be two points in $\wt{\rho}^{-1}(u)$, where $u$ is a fixed point in $P$. Denote by $c$ a loop about $u$ and by $\hat{c}$ its lift to $\Sigma$ with $\hat{c}(0)=s_0$. The endpoint $\hat{c}(1)$ depends on $[c]\in\pi_1(P,u)$ the homotopy class of $c$. Define the homomorphism $\tau:\pi_1(P,u)\to Z_2$ by $\tau([c])=0$ if $\hat{c}(1)=s_0$ and $\tau([c])=1$ if $\hat{c}(1)=s_1$. Since $Z_2$ is commutative, $\tau$ vanishes on the commutator subgroup $[\pi_1(P),\pi_1(P)]$ of $\pi_1(P,u)$, therefore $\tau$ induces a homomorphism $\sigma:\pi_1(P,u)]/[\pi_1(P,u),\pi_1(P,u)]\to Z_2$. We can identify $\sigma$ with an element of $H^1(P,Z_2)$ via the isomorphisms $$\pi_1(P,u)/[\pi_1(P),\pi_1(P)]\simeq H_1(P),$$ $$Hom,(H_1(P),Z_2)\simeq H^1(P,Z_2).$$

The exact sequence of groups $$1\to Z_2\to{\rm{Spin}}(H)\to SO(H)\to 1\leqno(2.3)$$ induces an exact sequence of the cohomology groups $$0\to H^1(M,Z_2)\to H^1(M,{\rm{Spin}}(H))\to H^1(M,SO(H))\overset{v}{\to} H^2(M,Z_2).\leqno(2.4)$$

Denote by $P$ the element of $H^1(M,SO(H))$ determined by $P(M,\pi,SO(H))$. Using the Remark 2.2. and the exactness of the
sequence (2.4), we obtain

{\bf Theorem 2.4.} {\it The vector bundle $\xi$ admits a spin structure if and only if $v(P)=0$.}

When $SO(H)=SO(H)_1$, P. de la Harpe [3] has proved that $v(P)=w_2(\xi)$, where $w_2(\xi)$ is the second Stiefel-Whitney class of $\xi$.

The following exact sequence $$0\to H^1(M,Z_2)\overset{i^*}{\to}H^1(SO(H),Z_2)\to H^2(M,Z_2),\leqno(2.5)$$
where $i$ is the natural inclusion of the fibre in the total space, can be obtained from spectral sequence associated to the principal bundle $P(M,\pi,SO(H))$ (see J.-P. Serre [11] p. 456).

If $\xi$ admits a spin structure $\sigma\in H^1(P,Z_2)$, then $\sigma+\pi^*(b)$ where $b\in H^1(M,Z_2)$ is {\it the most general spin structure of} $\xi$. It follows that, {there is a bijection between the set of isomorphism classes of spin structures of $\xi$ and} $H^1(M,Z_2)$. Consequently, a spin structure of $\xi$ is unique (up to an isomorphism) iff $H^1(M,Z_2)=0$.

{\bf Theorem 2.5.} {\it Let $\xi_1\oplus\xi_2$ be the Whitney sum of the vector bundles $\xi_1$ and $\xi_2$ over $M$. Given spin structures on two of the three vector bundles $\xi_1,\xi_2,\xi_1\oplus\xi_2$, there is a uniquely determined spin structure on the third.}

\noindent{\it Proof.} As in J. Milnor [8].

Let $M'$ be another manifold and let $f:M'\to M$ be a morphism of manifolds. Denote by $Pf(M',\pi',SO(H))$ the principal bundle induced from $P(M,\pi,SO(H))$ by the map $f$. This principal bundle is determined (up to an isomorphism) by the cocycle $g_{ij}\circ f$ associated to the open covering $\{f^{-1}(U_i)\}$ where $g_{ij}$ is the cocycle of $P(M,\pi,SO(H))$ associated to an open covering $\{U_i\}$. Using the Remark 2.2 we obtain

{\bf Theorem 2.6.} {\it Let $f:M'\to M$ be a morphism of manifolds. If $M$ admits a spin structure $\Sigma(M,\wt{\pi},{\rm{Spin}}(H))$ with respect to $P(M,\pi,SO(H))$, then $\Sigma f(M',\pi',{\rm{Spin}}(H))$ is a spin structure on $M'$ with respect to $Pf(M',\pi',SO(H))$.}

Suppose now that $f:M'\to M$ is a principal bundle with  group $G$. It follows that, there is an action of $G$ on $Pf$ defined by $(p',u)g=(p'g,u)$ for $(p',u)\in Pf$ and $g\in G$.

{\bf Theorem 2.7.} {\it Let $f:M'\to M$ be a principal bundle with group $G$. Then the following conditions are equivalent:

$(1)$ There exists a spin structure $\Sigma(M,\wt{\pi},{\rm{Spin}}(H))$ on $P(M,\pi,SO(H))$.

$(2)$ There exists a spin structure $\Sigma'(M',\wt{\pi}',{\rm{Spin}}(H))$ on $Pf(M',\pi',SO(H))$ with the following properties:

$a)$ $G$ acts on $\Sigma'$ such that $\Sigma'/G$ is a manifold and the projection $\Sigma'\to\Sigma'/G$ is a submersion,

$b)$ The action of $G$ on $\Sigma'$ commutes with the action of ${\rm{Spin}}(H)$ on $\Sigma'$.

$c)$ $\wt{\rho}'(wg)=\wt{\rho}'(w)g$ holds, for every $g\in G$ and $w\in \Sigma'$, where $\wt{\rho}':\Sigma'\to Pf$ is the covering map.}

\noindent{\it Proof.} $(1) \Rightarrow (2)$. By the Theorem 2.6, $\Sigma f(M',\pi',{\rm{Spin}}(H))$ is a spin structure on $Pf(M',\Pi',SO(H))$. Define an action of $G$ on $\Sigma'=\Sigma f$ by $(p',v)g=(p'g,v)$, where $(p',v)\in\Sigma'$ and $g\in G$. This action is proper and free. Moreover, the map $g\to(p'g,v)$ is an immersion of $G$ in $\Sigma'$, since $f:M'\to M$ is a principal bundle. It follows that $\Sigma'\to \Sigma'/G$ is  a principal bundle (see Bourbaki [1]), hence the property a) is verified. From $(p',v)b=(p',vb)$ for $b\in$Spin$(H)$, it follows $(p'g,vb)=(p',v)gb=(p',v)bg$ i.e. the property b).

For $w=(p',v)\in\Sigma$, we have $\wt{\rho}'(wg)=(p'g,\wt{\rho}(v))=\wt{\rho}'(w)g$, i.e. the property c).

$(2) \Rightarrow (1)$ Define an action of Spin$(H)$ on $\Sigma'/G$ by $wGb=wbG$, where $w\in\Sigma'$, $b\in$Spin$(H)$ and $wG$ is the orbit of $w$, and a surjection $h:\Sigma'/G\to M$ by $h(wG)=f(\pi'(w))$. The local isomorphism $t:f^{-1}(U)\times$Spin$(H)\to\Sigma'$, where $U$ is an open subset of $M$, defines a local isomorphism $s:U\times$Spin$(H)\to\Sigma'/G$ by $$s(p,b)=t(p'G,b)=t(p',b)G,\ \textrm{where\ } f(p')=p.\leqno(2.6)$$
The last equality from (2.6) is a consequence of $$\left\{\begin{array}{l}\pi'(\wt{\rho}'(wg))=\pi'(\wt{\rho}'(w))g\\ \\ f^*(\wt{\rho}'(wg))=f^*\wt{\rho}'(w),\ \ g\in G,\ w\in\Sigma',\end{array}\right.\leqno(2.7)$$ where $f^*:Pf\to P$ is the morphism induced by $f$. But (2.7) is equivalent to the property c). It is not difficult to see that $h:\Sigma'/G\to M$ is a principal bundle with group Spin$(H)$. The morphism $\wt{\rho}:\Sigma'/G\to P$ defined by ${\wt{\rho}}(wG)=f^*(\wt{\rho}'(w))$ satisfies $\pi\circ\wt{\rho}=h$ and $\wt{\rho}(wGb)=\wt{\rho}(wG)\rho(b)$, therefore $\Sigma'/G(M,h,$Spin$(H))$ is a spin structure on $P(M,\pi,SO(H))$.

Suppose that $f:M'\to M$ is a principal bundle with group $Z_2$. Let $\eta$ denotes the involution of $Pf$ defined by the action of $Z_2$ on it.

{\bf Corollary 2.7.} {\it Let $f:M'\to M$ be a principal bundle with group $Z_2$. The following conditions are equivalent:

$(1)$ There exists a spin structure $\Sigma(M,\pi,{\rm{Spin}}(H))$ on $P(M,\pi,SO(H))$,

$(2)$ There exists a spin structure $\Sigma'(M',\wt{\pi},{\rm{Spin}}(H))$ on $Pf(M',\pi',SO(H))$ endowed with an involution $\eta'$ which corresponds to identity on {\rm{Spin}}$(H)$ and which commutes with the involution $\eta$.}

\noindent{\it Proof.} In order to apply the Theorem 2.7 it is sufficient to remark that $\eta'$ defines an action of $Z_2$ on $\Sigma'$, which commutes with the action of Spin$(H)$, such that $\Sigma'\to \Sigma'/Z_2$ is a submersion and $\rho'(wZ_2)=\wt{\rho}'(w)Z_2$ holds, for every $w\in\Sigma'$.

{\bf Remark 2.3.} The above corollary has been obtained by I. Popovici [10] in non-orientable and finite dimensional case.

{\bf Definition 2.6.} Let $N$ be a manifold modeled on $H$, with an oriented Riemannian nuclear structure (resp. an oriented Riemannian Hilbert-Schmidt structure). A spin structure on $N$ is a spin structure on the tangent bundle $TN$.

\section{Examples} a) The trivial bundle $M\times SO(H)$ admits a spin structure, unique if $M$ is simply connected.

b) Let $\left\{x=(x_1,x_2,...) / x_i\in R,\ \dd\sum^\9_{i=1}x^2_i<\9\right\}$ be the Hilbert space $l_2$. The Hilbert torus $T=l_2/\dd\sum^\9_{i=1}Ze_i$, where $e_1,e_2,...$ is a base of $l_2$, is a Hilbert-Lie group modeled on $l_2$. It admits a canonical analytic atlas such that the derivatives of the coordinate changes is always the identity on $l_2$. It follows that the corresponding Riemannian nuclear structure is orientable and $T$ with this structure, admits a spin structure (with respect to Spin$(l_2)_1)$ since $w_2(T)=0$.

c) Let $M\subset H$, be a smoothly imbedded manifold with an oriented Riemannian nuclear structure. From Theorem 2.6, it follows that a spin structure on $M$ determines a spin structure on the normal bundle to $M$. The converse is also valid.

\section{Clifford bundles} In this section we limit our considerations to the group $SO(H)_1$. Let $Cl(H)$ be the Clifford algebra of $H$ viewed as a $C^*-$algebra [3]. A $^*-$automorphism $\psi$ of $Cl(H)$ which satisfies $\psi(H)=H$ is called a Bogoliubov automorphism. Denote by Bog$(Cl(H))$ the group of Bogoliubov automorphisms and by $C:O(H)\to$Bog$(Cl(H))$ the canonical isomporhism  described in [3]. A Bogoliubov automorphism $\psi$ is said to be {\it inner} if there exists $u\in Cl(H)$ such that $\psi(v)=uvu^{-1}$ for every $v\in Cl(H)$ and is said to be {\it inner and even} if $u$ is even. Let $S(M,O(H))$ be a principal bundle. Its associated fibre bundle with fibre $Cl(H)$ (the action of the group $O(H)$ on $Cl(H)$ is given by $C$) is an algebric  bundle, called the Clifford bundle. We can obtain the Clifford bundle in another way. For this, let $\xi:E\to M$ be a Riemannian vector bundle. The fibres of $\xi$ are Hilbert spaces. Let $E_p$ be the fibre of $\xi$ in $p\in M$ and let $Cl(E_p)$ be the Clifford algebra of $E_p$. The set $\dd\bigcup_{p\in M}Cl(E_p)$ and the projection $Cl(E_p)\to p$ can be taken as the total space and projection of the Clifford bundle. We remark that $\xi$ can be identified with a subbundle of the Clifford bundle. Let Bog$(Cl(H))_i$ be the group of inner Bogoliubov automorphisms. The isomorphism $O(H)_1\simeq$Bog$(Cl(H)_i$ (see [4]) implies the following

{\bf Theorem 4.1.} {\it There exists a Riemannian nuclear structure on the Riemannian vector bundle $\xi:E\to M$ iff the Clifford bundle admits a reduction to ${\rm{Bog}}(Cl(H))_i$.}

From the isomophism $SO(H)_1\simeq$Bog$(Cl(H))_{ie}$ (see [3]), where Bog$(Cl(H))_{ie}$ is the group of inner and even Bogoliubov automorphism, it follows

{\bf Theorem 4.2.} {\it A Riemannian nuclear structure on the Riemannian vector bundle $\xi:E\to M$ is orientable iff the Clifford bundle admits a reduction to ${\rm{Bog}}(Cl(H))_{ie}$.}

\noindent{\footnotesize{\it Received February 27, 1978}}
\newpage




\def\bfe{{\mathbf{E}}}

\runningauthor={M. ANASTASIEI}
\runningtitle={CONNEXIONS SUR LES FIBRES SPINORIELS}
\noindent
\baselineskip 8pt
\noindent{\footnotesize{An. \c{S}tiin\c{t}.Univ. ``Al. I. Cuza''
\hfill\break Ia\c{s}i, Vol. XXV, s. Ia, Mat.,\hfill\break
1979, f.2., 337--344}}

\vskip 2cm
\baselineskip 11.5pt plus .15pt
\centerline{\bf\Large CONNEXIONS SUR LES FIBRES SPINORIELS}
\vskip .5cm
\centerline{\bf \footnotesize{PAR}}
\vskip .5cm
\centerline{\bf \footnotesize{MIHAI ANASTASIEI}}
\vskip 1cm

Dans [1] on a d\'{e}fini et \'{e}tudi\'{e} les structures spinorielles sur les vari\'{e}t\'{e}s hilbertiennes. Dans le pr\'{e}sent article on consid\`{e}re les connexions adapt\'{e}es aux structures spinorielles (les connexions spinorielles). Apr\`{e}s quelques r\'{e}sultats relatifs aux connexions sur les fibr\'{e}s principaux banachique (\S1), on d\'{e}finit les connexions spinorielles (\S2). En utilisant la d\'{e}riv\'{e}e covariante des spineurs, on met en \'{e}vidence une famille d'op\'{e}rateurs diff\'{e}rentiels du premier ordre (au sens de [7] p. 91) sur les fibr\'{e}s spinoriels. Si la vari\'{e}t\'{e} est a dimension finie l'operat\'{e}ur de Dirac peut \^{e}tre d\'{e}duit de cette famille des op\'{e}rateurs.

\setcounter{section}{0}

\section{Connexions sur les fibr\'{e}s principaux banachiques} Soient $M$ une vari\'{e}t\'{e} diff\'{e}rentiable de classe $C^\9$ model\'{e}e sur un espace de Banach ${\mathbf{M}}$, $G$ un groupe de Lie-Banach r\'{e}el et $P(M,\pi_P,G)$ un fibr\'{e} principal (abr\'{e}g\'{e} f.p.) de base $M$, de groupe structural $G$ et de projection $\pi_P$. Nous notons par $R_g:u\to ug$, $u\in P$, $g\in G$ l'action de $G$ sur $P$ et par $\sigma_u:G\to P$, $u\in P$, l'application $g\to ug$, $g\in G$. $\sigma_u$ est un diffeomorphisme de classe $C^\9$ entre le groupe $G$ et la fibre au-dessus du point $p=\pi_P(u)$. Le foncteur tangent sera not\'{e} par $T$. Consid\'{e}rons la suite exacte de fibr\'{e}s vectoriels au-dessus de $P$:
\begin{equation}
0\to P\times {\mathbf{G}}\overset{I}\rightarrow TP\overset{T\pi_P!}\longrightarrow \pi^*_P TM\to 0
\end{equation}
o\`{u}:
- $\mathbf{G}$ est l'alg\`{e}bre de Lie-Banach de $G$,

- $\pi^*_P TM$ est le fibr\'{e} image r\'{e}ciproque de $TM$ par $\pi_P$,

- $T\pi_P!=(\tau_P,T\pi_P)$, $\tau_P:TP\to P$,

- $I(u,A)=T\tau_u(A)$, $u\in P$, $A\in \mathbf{G}$.

Un $G-$fibr\'{e} vectoriel est un fibr\'{e} vectoriel sur lequel $G$ op\`{e}re par des automorphismes de fibr\'{e} vectoriel. Si nous consid\'{e}rons les actions naturelles de $G$ sur $TP$ et $\pi^*_P TM$ ainsi que l'action de $G$ sur $P\times \mathbf{G}$, $(u,A)g=(ug, ad(g{-1})A)$, $u\in P$, $g\in G$, $A\in \mathbf{G}$, alors la suite (1.1) devient une suite exacte de $G-$fibr\'{e}s vectoriels.

\ms

{\bf D\'{e}finition 1.1.} [8] Une connexion (infinit\'{e}simale) sur le f.p. $P(M,\pi_P,G)$ est une scission de la suite exacte $(1.1)$ de $G-$fibr\'{e}s vectoriels.

\ms

Soit $\G:\pi^*_P TM\to TP$ une telle scission. Alors, il existe un morphisme unique $\oo:TP\to P\times \mathbf{G}$, tel que $\oo\circ I={\rm{id}}|_{P\times \mathbf{G}}$. Pour chaque point $u\in P$, il existe une d\'{e}composition unique en somme directe de $T_uP$
\begin{equation}
T_uP=I(P\times G)_u\oplus\G(\pi^*_P TM)_u.
\end{equation}
La d\'{e}composition (1.2) definit de mani\`{e}re evidente deux projecteurs sur $T_uP$, qui seront not\'{e}s par $h$ (de noyau $I(P\times \mathbf{G})_u)$ et $v$ (de noyau $\G(\pi^*_P TM)_u$). La compatibilit\'{e} de $\G$ avec les actions de $G$ sur $\pi^*_P TM$ et $TP$ implique
\begin{equation}
\G(ug,Z)=TR_g\G(u,Z),\ u\in P, g\in G, Z\in TM,
\end{equation}
\begin{equation}
\oo(TR_g X_u)=(u,\oo_u(X_u))g,\ u\in P,\ g\in G,\ X_u\in T_uP,
\end{equation}
o\`{u} nous avons not\'{e} par $\oo_u$ la restriction de $\oo$ \`{a} $T_uP$ et nous avons pos\'{e} $\oo(X_u)=(u,\oo_u (X_u))$. L'application lin\'{e}aire $\oo_u:T_uP\to \mathbf{G}$, $u\in P$, a les deux propri\'{e}t\'{e}s suivantes
\begin{equation}
\oo_u(\sigma_u(A))=A,\ A\in \mathbf{G},
\end{equation}
\begin{equation}
(R^*_g\oo)(X_u)\overset{def}{=}\oo_{ug}(TR_g X_u)=ad(g^{-1})\oo_u(X_u),\ u\in P,g\in G.
\end{equation}

Soit $F_G(P)$ l'ensemble des fonctions diff\'{e}rentiables de classe $C^\9$ d\'{e}finies sur $P$ \`{a} valeurs dans $\mathbf{G}$ et soit $\oo:{\cal X}(P)\to F_G(P)$ une 1-forme de classe $C^\9$ sur $P$ \`{a} valeurs dans $\mathbf{G}$, d\'{e}finie par:
\begin{equation}
\oo(X)_u=\oo_u(X_u),\ u\in P, X\in{\cal X}(P),
\end{equation}
o\`{u} ${\cal X}(P)$ est le module des champs de vecteurs sur $P$. Evidemment, (1.5) et (1.6) impliquent:
\begin{equation}
\oo(\sigma(A))=A,\ A\in \mathbf{G},
\end{equation}
\begin{equation}
R^*_g\oo=ad(g^{-1})\oo,\ g\in G,
\end{equation}
o\`{u} $\sigma(A)$ est le champ vectoriel $u\to \sigma(A)$.

R\'{e}ciproquement, une 1-forme sur $P$ \`{a} valeurs dans $\mathbf{G}$, d\'{e}finit une application inverse \`{a} gauche pour $I$, compatible avec les actions de $G$ sur $P\times \mathbf{G}$ et $TP$, gr\^{a}ce aux formules (1.8) et (1.9). Vu que la suite (1.1) est exacte, cette inverse d\'{e}finit une connexion sur $P(M,\pi_P,G)$. Donc, nous avons \'{e}tabli l'equivalence entre la d\'{e}finition 1.1 et la

\ms

{\bf D\'{e}finition 1.2.} Une connexion (infinit\'{e}simale) sur le f.p. $P(M,\pi_P,G)$ est une 1-forme sur $P$ \`{a} valeurs dans $G$ avec les propri\'{e}t\'{e}s (1.8) et (1.9).

\ms

Ainsi, nous avons recup\'{e}r\'{e}, pour nos buts, une d\'{e}finition bien connue en dimension finie, des connexions (infinit\'{e}simales) sur un f.p. Une autre caract\'{e}risation est donn\'{e}e par le

\ms

{\bf Theoreme 1.1.} {\it L'existence sur un f.p. $P(M,\pi_P,G)$ d'une connexion (infinit\'{e}simale) \'{e}quivaut \`{a} l'existence d'un projecteur $h:TP\to TP$ $(h\circ h=h)$ avec les propri\'{e}t\'{e}s:}
\begin{equation}
\ker\ h=I(P\times \mathbf{G}),
\end{equation}
\begin{equation}
TR_g\circ h=h\circ TR_g,\ g\in G.
\end{equation}

\ms

La preuve de ce th\'{e}or\`{e}me est imm\'{e}diate si nous remarquons que (1.3) \'{e}quivaut \`{a} (1.11).

Le morphisme $F=v-h$ de $TP$ d\'{e}finit une structure presque-produit sur $P$, associ\'{e}e d'une mani\`{e}re naturelle \`{a} la scission $\G$. L'espace des vecteurs propres correspondant \`{a} la valeur propre 1 de l'op\'{e}rateur $F_u:T_uP\to T_uP$ est $I(P\times \mathbf{G})_u$. Vu que le projecteur $v$ a la propri\'{e}t\'{e}
\begin{equation}
TR_g\circ v=v\circ TR_g,\ g\in G,
\end{equation}
il r\'{e}sulte que (1.11) \'{e}quivaut \`{a}
\begin{equation}
TR_g\circ F=F\circ TR_g,\ g\in G.
\end{equation}

En utilisant le th\'{e}or\`{e}me 1.1 on obtient le

\ms

{\bf Th\'{e}or\`{e}me 1.2.} {\it Il existe une connexion (infinit\'{e}simale) sur $P(M,\pi_P,G)$ si et seulement si, il existe une structure presque-produit $F$ sur $P$, avec les propri\'{e}t\'{e}s}

a) $F_u(X_u)=X_u\Leftrightarrow X_u\in I(P\times G)_u$, $X_u\in T_uP$,

b) $TR_g\circ F=F\circ TR_g$, $g\in G$.

\ms

{\it Remarque.} Les th\'{e}or\`{e}mes 1.1 et 1.2 ont \'{e}t\'{e} \'{e}tablis en dimension finie, par V. Cruceanu dans [2] et [3].

\ms

Soit $(f,\varphi_0,h_0):P(M,\pi_P,G)\to P'(M',\pi'_{P'},G')$ o\`{u} $f:P\to P'$, $\varphi_0:G\to G'$, $h_0:M\to M'$ un homomorphisme de f.p., c'est-\`{a}-dire:
\begin{equation}
\pi_{P'}\circ f=h_0\circ \pi_P,\ f(ug)=f(u)\varphi_0(g),\ u\in P, g\in G.
\end{equation}

\ms

{\bf Definition 1.3.} Soient les f.p. $P(M,\pi_P,G)$ et $P'(M',\pi'_{P'},G')$ munis des connexions (infinit\'{e}simales) $\G$ resp. $\G'$. Nous dirons que l'homomorphisme $(f, \varphi_0,h_0)$ est compatible avec les connexions $\G$ et $\G'$ si nous avons
\begin{equation}
Tf\circ \G=\G'\circ(f\times Th_0).
\end{equation}

\ms

{\it Remarque.} La relation (1.15) \'{e}quivaut \`{a}
\begin{equation}
\oo'\circ Tf=(f\times T\varphi_0)\circ\oo.
\end{equation}

\ms

La d\'{e}monstration du th\'{e}or\`{e}me suivant se conduit comme en dimension finie (voir [6] p. 79--82).

\ms

{\bf Th\'{e}or\`{e}me 1.3.} {\it Soit $(f,\varphi_0,h_0):P(M,\pi_P,G)\to P'(M',\pi_{P'},G')$ un homomorphisme de fibr\'{e}s principaux, avec $h_0$ un diff\'{e}omorphisme.}

{\it a)	Soit $\G$ une connexion (infinit\'{e}simale) sur $P(M,\pi'_P,G)$. Alors, il existe une connexion (infinit\'{e}simale) unique $\G'$ sur $P'(M',\pi'_{P},G')$ de mani\`{e}re que l'homomorphisme $(f,\varphi_0,h_0)$ soit compatible avec les connexions (infinit\'{e}simales) $\G$ et $\G'$.}

{\it b)	En supposant que $\varphi_0$ est un diff\'{e}omorphisme local, soit $\G'$ une connexion (infinit\'{e}simale) sur $P'(M',\pi_{P'},G')$. Alors, il existe une connexion (infinit\'{e}simale) unique $\G$ sur $P(M,\pi_P,G)$ telle que l'homomorphisme $(f,\varphi_0,h_0)$ soit compatible avec les connexions (infinit\'{e}simales) $\G$ et $\G'$.}

\ms

Soit $\F$ un espace de Banach. En supposant que $G$ op\`{e}re sur $\F$ par un homomorphisme $\psi:G\to GL(\F)$, soit $\pi:E\to M$ le fibr\'{e} vectoriel associ\'{e} \`{a} $P(M,\pi_P,G)$ de fibre type $\F$. Par une modification l\'{e}g\`{e}re d'une preuve de J-P. Penot (voir [8]) on peut montrer que toute connexion (infinit\'{e}simale) sur $P(M,\pi_P,G)$ induit une connexion vectorielle unique sur $\pi:E\to M$ c'est-\`{a}-dire il existe un morphisme de fibr\'{e}s vectoriels $K:TE\to E$, tel que pour chaque carte vectorielle $(U,\varphi,\Phi)$ de $\pi$, nous avons
\begin{equation}
(\Phi\circ K\circ T\Phi^{-1})=(x,\xi,y,\eta)=(x,\eta+\G_\varphi(x))(y,\xi),\ x,y\in \M
\end{equation}
$\xi,\eta\in \F$, o\`{u} $\G_\varphi(x)\in L^2(M,\F;\F)$ correspond aux symboles de Christoffel usuels. Notons par ${\cal X}_E(M)$ le module des sections sur $M$ dans le fibr\'{e} vectoriel $\pi:E\to M$ et possons ${\cal X}_{TM}(M)={\cal X}(M)$. Il existe (voir [4], p. 17) une d\'{e}rivation covariante unique $\nabla_X$, $X\in{\cal X}(M)$, associ\'{e}e naturellement \`{a} l'application $K$, donn\'{e}e dans une carte vectorielle quelconque $(U,\varphi,\Phi)$ par
\begin{equation}
\nabla_XS|_{\varphi(p)}=\pp S_{\varphi}|_{\varphi(p)}(X_\varphi)+\G^{(x)}_{\varphi}(X_\varphi,S_{\varphi}),\ p\in U, X\in{\cal X}(M),
\end{equation}
$$x=\varphi(p),\ S\in{\cal X}_E(M),$$
o\`{u} $X_\varphi=T_\varphi\circ X$, $S_\varphi=\Phi\circ S$ et $\pp$ est le symbole de diff\'{e}rentiation de Fr\'{e}chet. Supposons que $M$ admet une partition de l'unit\'{e}. D'apr\`{e}s le lemme 3.1 de [4] et la formule (1.18) nous pouvons d\'{e}finir l'application $\nabla:T_pM\times {\cal X}_E(U)\to {\cal X}_E(U)$, avec $U$ un ouvert de $M$ par $(X_p,S)\to\nabla_{X_p}S=\nabla_X S$, o\`{u} $X$ est un champ arbitraire de vecteurs qui co\"{\i}ncide au point $p$ avec $X_p$ et $S\in {\cal X}_E(U)$. Cette application est $R-$lin\'{e}aire relativement \`{a} $X_p$ et
\begin{equation}
\nabla_{X_p}(fS)=T_pf(X_p)S+f(p)\nabla_{X_p}S,\ \ S\in{\cal X}_E(M),
\end{equation}
o\`{u} $f$ est une fonction r\'{e}elle quelconque sur $M$. Nous consid\'{e}rons l'op\'{e}rateur de diff\'{e}rentiation covariante $\nabla:{\cal X}_E(M)\to{\cal X}_{L(TM,E)}(M)(S\to\nabla S)$ donne par
\begin{equation}
(\nabla S)(p)=(\nabla_{X_p}S)(p),\ p\in M, X_p\in T_pM,\ S\in{\cal X}_E(M).
\end{equation}
En utilisant (1.19) on obtient le (voir aussi [4], p.6)

\ms

{\bf Th\'{e}or\`{e}me 1.1.} {\it L'op\'{e}rateur de diff\'{e}rentiation covariante $\nabla$ est un op\'{e}rateur diff\'{e}rentiel du premier ordre.}

\section{Structures spinorielles et connexions spi\-no\-ri\-elles}
\setcounter{equation}{0}

Soit $\H$ un espace de Hilbert r\'{e}el, s\'{e}parable et de dimension infinie. Nous notons par $Cl(\H)_\9$ l'alg\`{e}bre de Clifford sur $\H$ r\'{e}lative \`{a} la forme quadratique $Q(x)=\|x\|^2=(x,x)$, $x\in \H$, structur\'{e} comme une $C^*-$algebre (voir [5]). Soit $\J$ une structure complexe sur $\H$, c'est-\`{a}-dire un op\'{e}rateur orthogonal sur $\H$ avec $\J^2=-$id. Si nous posons i $x=\J x$ et $\langle x,y\rangle=(x,y)+(\J x,y)$, $x,y\in \H$, l'espace $H$ devient un espace de Hilbert complexe, qui sera not\'{e} par $\H_c$. Soient $\Lambda^n \H_C$ la puissance ext\'{e}rieure des $n$ exemplaires de $\H_C$ et $\Lambda(\H_C)=\subset_{n\geq 0}\oplus \Lambda^n \H_C$ avec la structure naturelle de l'espace de Hilbert complexe.

L'alg\`{e}bre ext\'{e}rieure $\Lambda(\H_C)$ a une $Z_2-$graduation naturelle, $\Lambda(\H_C)=\Lambda^0\oplus\Lambda^1$ o\`{u} $\Lambda^0=\dd\oplus_{k\geq 0}\Lambda^{2k}\H_C$ et $\Lambda^1=\dd\oplus_{k\geq 0}\Lambda^{2k+1}\H_C$. Il existe une repr\'{e}sentation fid\`{e}le et irr\'{e}ductibile $F$ de $Cl(\H)_\9$ sur l'espace de Hilbert $\Lambda(\H_C)$ (voir par exemple [9]). Soient $Cl(\H)^*_\9$ le groupe multiplicatif des \'{e}l\'{e}ments inversibles de $Cl(\H)_\9$ et Spin$(\H)_\9=\{u\in Cl(\H)^*_\9\ |\ u \H u^{-1}=\H, u\b(u)=\b(u)u=1, \a(u)=u\}$ o\`{u} $\a$ est l'involution et $\b$ est l'antiinvolution principale de $Cl(\H)_\9$. P. de la Harpe a montr\'{e} dans [5] que le groupe Spin$(\H)_\9$ est groupe de Lie-Banach.

Soit $O(\H)_1$ le groupe des op\'{e}rateurs orthogonaux sur $\H$ de la forme $id+A$, o\`{u} $A$ est un op\'{e}rateur nucl\'{e}aire. Le groupe de Lie-Banach $O(\H)_1$ a deux composantes connexes. Le rev\^{e}tement universel de la composante connexe de l'identit\'{e} $SO(\H)_1$, est exactement Spin$(\H)_\9$ (voir [5]). Soit $\Delta=F\ |\ $ Spin$(\H)_\9$. Les espaces $\Lambda^0$ et $\Lambda^1$ sont invariants par $\Delta$ (voir [9]) et ils d\'{e}finissent deux sous-repr\'{e}sentations de $\Delta$ qui seront not\'{e}es par $\Delta^0$ et $\Delta^1$, respectivement. Ces deux repr\'{e}sentations $\Delta^0$ et $\Delta^1$ sont continues, injectives et irr\'{e}ductibles (voir [9]). L'espace $\Lambda(\H_C)$ (en abr\'{e}g\'{e} $\Lambda$) s'appelle l'espace de spineurs relatif \`{a} $SO(\H)_1$ et les espaces $\Delta^0$ et $\Delta^1$ s'appellent les espaces de semi-spineurs relatifs \`{a} $SO(\H)_1$.

Soit $P(M,SO(\H)_1)$ un f.p. de base $M$ et de groupe structural $SO(\H)_1$. Un f.p. $\Sigma(M,$Spin$(\H)_\9)$ qui est l'extension de $P(M,SO(\H)_1)$ associ\'{e}e \`{a} l'ho\-mo\-mor\-phisme de rev\^{e}tement $\rho:$Spin$(\H)_\9\to SO(\H)_1$, s'appelle structure spi\-no\-ri\-elle sur $P(M,SO(\H)_1)$ ou structure spinorielle sur $M$ relative \`{a} $P(M,SO(\H)_1)$ (voir [1]). Nous notons par $(\wt{\rho},\rho):\Sigma(M,$Spin$(\H)_\9)\to P(M,SO(\H)_1)$ l'ho\-mo\-mor\-phisme d'extension.

\ms

{\bf Definition 2.1.} On appele {\it connexion spinorielle} une connexion (infinit\'{e}simale) sur $\Sigma(M,$Spin$(\H)_\9)$.

\ms

Soient $C(P)$ et $C(\Sigma)$ les ensembles de connexions sur $P(M,SO(\H)_1)$ et $\Sigma(M$, Spin$(\H)_\9)$, respectivement. Comme $\rho$ est un diff\'{e}omorphisme local il r\'{e}sulte en vertu du th\'{e}or\`{e}me 1.3, qu'il existe une bijection $\wt{\rho}:C(\Sigma)\to C(P)$. Si nous notons par $\wt{\oo}$ et $\oo$ les 1-formes des deux connexions correspondantes par $\wt{\rho}$, nous avons:
\begin{equation}
\oo_{\wt{\rho}(\wt{u})}(T\wt{\rho}X_{\wt{u}})=T\rho\wt{\oo}_{\wt{u}}(X_{\wt{u}}),\ \wt{u}\in\Sigma,\ X_{\wt{u}}\in T_{\wt{u}}\Sigma.
\end{equation}
$\wt{\oo}_{\wt{u}}(X_{\wt{u}})$ et $\oo_{\wt{\rho}(\wt{u})}(T\wt{\rho}X_{\wt{u}})$ sont dans les alg\`{e}bres de Lie-Banach des groupes Spin$(H)_\9$ et $SO(H)_1$, respectivement.

En vertu de la proposition 12 [5] nous obtenons
\begin{equation}
\|\oo_{\wt{\rho}(\wt{u})}(T\wt{\rho}X_{\wt{u}})\|_1=4\|\wt{oo}_{\wt{u}}(X_{\wt{u}})\|_\9,\ \ \wt{u}\in\Sigma,X_{\wt{u}}\in T_{\wt{u}}\Sigma
\end{equation}
o\`{u} $\|\ \|_1$ est la norme nucl\'{e}aire et $\|\ \|_\9$ est la norme de $C^*-$alg\`{e}bre sur $Cl(\H)_\9$. En dimension finie la relation (2.2) coincide avec la relation (5.4) de [10].

Vu que Spin$(\H)_\9$ op\`{e}re sur $\Lambda$ nous pouvons consid\'{e}rer le fibr\'{e} associ\'{e} \`{a} $P(M, $Spin$(\H)_\9)$ avec la fibre type $\Lambda$ qui sera nomm\'{e} fibr\'{e} spinoriel. Les fibr\'{e}s associ\'{e}s \`{a} $P(M,$ Spin$(\H)_\9)$ avec les fibr\'{e}s type $\Lambda^0$ et $\Lambda^1$, respectivement, seront nomm\'{e}s fibr\'{e}s semi-spinoriels. Une section du fibr\'{e} spinoriel sera nomm\'{e}e champ spinoriel.

Une connexion spinorielle induit une d\'{e}rivation covariante dans le fibr\'{e} spinoriel qui sera nom\'{e}e d\'{e}rivation spinorielle. L'existence de la bijection $\wt{\rho}$ implique que toute connexion sur $P(M,SO(\H)_1)$ induit une d\'{e}rivation spinorielle.

Une r\'{e}duction de groupe structural de $TM$ ($M$ model\'{e}e sur $\H$) au groupe $SO(\H)_1$ s'appelle structure riemanniene nucl\'{e}aire orient\'{e}e sur $M$. Si $M$ est munie d'une structure riemannienne nucl\'{e}aire orient\'{e}e, le fibr\'{e} de rep\`{e}res de $TM$ est un f.p. de base $M$ et de groupe structural $SO(\H)_1$ qui sera not\'{e} par $R(M)$. Nous supposons que $M$ a une structure spinorielle, c'est-\`{a}-dire il existe un f.p. $\Sigma(M)$ de base $M$ et de groupe structural Spin$(\H)_\9$, l'extension de $R(M)$ par $\rho$.

Dans la suite nous allons mettre en \'{e}vidence une classe d'op\'{e}rateurs diff\'{e}rentiels du premier ordre sur une vari\'{e}t\'{e} munie d'une structure spinorielle. Soit $\Lambda(M)$ le fibr\'{e} spinoriel. Vu que $H\subset Cl(\H)_\9$, pour tout $x\in \H$, on a une application lin\'{e}aire $F(x):\Lambda\to\Lambda$. Pour $b\in$Spin$(\H)_\9$ et $s\in\Lambda$ nous avons
\begin{equation}
\begin{array}{c}
\Delta(b)(F(x)s)=F(b)(F(x)s)=F(bx)s=F(bxb^{-1}b)s=\\ \\
=F(bxb^{-1})F(b)s=F(\rho(b)x)F(b)s.
\end{array}
\end{equation}

Soit $(U_i,\varphi_i)$ un atlas de la vari\'{e}t\'{e} $M$ et $\{\tau_i:\tau_s^{-1}(U_i)\to U_i\times\Lambda\}$ les cartes du fibr\'{e} spinoriel $\tau_s:\Lambda(M)\to M$. Les applications $\tau_j\circ\tau^{-1}_i:U_i\cap U_j\to L(\Lambda)$ ont leur images dans $\Delta($Spin$(H)_\9)$. Vu que $\Delta$ est injective, $\tau_j\circ\tau^{-1}_i(p)$, $p\in U_i\cap U_j$, s'identifie \`{a} son image dans Spin$(H)_\9$. Soient $\{\Phi_i:\tau^{-1}(U_i)\to U_i\times H\}$ les cartes du fibr\'{e} tangent $\tau:TM\to M$. Les applications $\Phi_j\circ\Phi_i^{-1}:U_i\cap U_j\to L(H)$ ont leur images dans $SO(H)_1$; comme $\rho$ est surjectif il existe $b\in$Spin$(H)_\9$ tel que $\rho(b)=\Phi_j\circ\Phi_i^{-1}(p)=\pp(\varphi_j\circ\varphi_i^{-1})(p)$, $p\in U_i\cap U_j$. Nous d\'{e}finissons maintenant une application $\Psi:TM\times\Lambda(M)\to\Lambda(M)$ par
\begin{equation}
\Psi(X_p,s_p)(p)=\tau^{-1}_{i,p}(F(X_{\varphi_i}))(s_{\tau_i}),\  p\in M, X_p\in T_pM, s_i\in\Lambda_p M
\end{equation}
o\`{u} $X_{\varphi_i}=\Phi_{i,p}(X_p)$, $s_{\tau_i}=\tau_{i,p}(s_p)$.

D'apr\`{e}s la relation (2.1) il r\'{e}sulte que la d\'{e}finition de $\Psi$ ne d\'{e}pend pas des cartes locales choisies. Une connexion lin\'{e}aire sur $M$ induit une connexion sur $\Lambda(M)$ et donc une d\'{e}rivation covariante $\nabla$. Pour tout $X\in{\cal X}(M)$ nous d\'{e}finissons un op\'{e}rateur $D_X:{\cal X}_{\Lambda(M)}(M)\to{\cal X}_{\Lambda(M)}(M)$ par
\begin{equation}
(D_XS)(p)=\Psi(X_p,\nabla_{X_p}S),\ p\in M,\ S\in{\cal X}_{\Lambda(M)}(M)
\end{equation}

\ms

{\bf Th\'{e}or\`{e}m 2.1.} {\it a) L'op\'{e}rateur $D_X$ est pour tout $X\in{\cal X}(M)$ un op\'{e}rateur diff\'{e}rentiel du premier ordre.

b) Le symbole de l'op\'{e}rateur $D_X$ est donn\'{e} par
\begin{equation}
\sigma_1(D_X)(v_p)=v_p(X_p)\Psi(X_{p,\cdot})
\end{equation}
o\`{u} $v_p$ est une 1-forme non nulle sur $T_p M$ pour chaque $p\in M$.}

\ms

\noindent{\it D\'{e}monstration.} a) Soit $S\in{\cal X}_{\Lambda(M)}(M)$ tel que le jet d'ordre 1, $(j^1S)(p)=0$. Vu que $\nabla$ est un op\'{e}rateur diff\'{e}rentiel du premier ordre, il r\'{e}sulte que $\nabla_{X_p}S=0$, donc $(D_XS)(p)=0$.

b) Soient $f$ une fonction r\'{e}elle sur $M$ telle que $f(p)=0$ et $T_p f=v_p$. Soit $S\in{\cal X}_{\Lambda(M)}(M)$ tel que $S(p)=s$. Nous avons $$\sigma_1(D_X)(v_p)(s)=D_X(fS)(p)=\Psi(X_p,\nabla_{X_p}(fS))=\Psi(X_p,T_pf(X_p)S(p)+$$ $$+f(p)\nabla_{X_p}S)=\Psi(X_p,v_p(X_p)s) {\rm{\ donc\ }}\sigma_1(D_X)(v_p)=v_p(X_p)\Psi(X_{p,\cdot}).$$

Un calcul direct montre que $$\sigma_1(D_X)(v_p)\circ\sigma_1(D_X)(v_p)=\a[v_p(X_p)]^2{\rm{id}}$$ o\`{u} $\a$ est un nombre r\'{e}el non nul et $X_p$ est non nul. Il r\'{e}sulte que le symbole $\sigma_1(D_X)(v_p)$ est injectif seulement si $v_p$ est injective, donc l'op\'{e}rateur $D_X$ n'est pas elliptique au sens de [7].

Supposons que $M$ est de dimension finie \'{e}gale \`{a} $n$ et introduisons l'op\'{e}ra\-teur de Dirac $D$ \`{a} l'aide des op\'{e}rateurs $D_X$ (voir aussi [7]). Soit $U$ un voisinage ouvert de $p\in M$ muni d'un champ de rep\`{e}res orthonorm\'{e}s $\{X_1, X_2,...,X_n\}$. D\'{e}finissons d'abord un op\'{e}rateur diff\'{e}rentiel du premier ordre $D_U$ sur $\Lambda(M)$ par
\begin{equation}
(D_U S)(p)=\dd\sum^n_{k=1}(D_{X_k}S)(p),\ \ S\in{\cal X}_{\Lambda(M)}(M).
\end{equation}

La d\'{e}finition de $D_U$ ne d\'{e}pend pas du champ $\{X_1,...,X_n\}$. Si l'on consid\`{e}re une famille d'op\'{e}rateurs $D_{U_i}$ de la forme (2.5) avec $\{U_i)$ un recouvrement ouvert de $M$, deux op\'{e}rateurs arbitraires $D_{U_i}$ et $D_{U_j}$ co\"{\i}ncident sur $U_i\cap U_j$. Il r\'{e}sulte que la famille d'op\'{e}rateurs $D_{U_i}$ d\'{e}finit un op\'{e}rateur diff\'{e}rentiel du premier ordre unique $D$ sur $M$ tel que la restriction de $D$ \`{a} chaque $U_i$ co\"{\i}ncide avec $D_{U_i}$. Le symbole de l'op\'{e}rateur $D$ est donn\'{e} par
\begin{equation}
\sigma_1(D)(v_p)=\dd\sum^n_{k=1}v_p(X_{k,p})\Psi(X_{k,p,\cdot}).
\end{equation}

Parce que le champ $\{X_1,...,X_n\}$ est un champ de rep\`{e}res orthonorm\'{e}s nous obtenons
\begin{equation}
\sigma_1(D)(v_p)\circ\sigma_1(D)(v_p)=\dd\sum^n_{k=1}[v_p(X_{k,p})]^2{\rm{id}}.
\end{equation}

Par suite l'op\'{e}rateur $D$ est elliptique.

{\bf BIBLIOGRAPHIE}

\begin{enumerate}
  \item Anastasiei M., {\it Structures spinorielles sur les vari\'{e}t\'{e}s hilbertiennes}, C.R. Acad. Sci. Paris, T. 284, S. A. p. 943--946.
  \item Cruceanu V., {\it Sur la structure presque-produit associ\'{e}e \`{a} une connexion sur un espace fibr\'{e}}, An \c{S}t. Univ. Iasi, Matematica, T. XV, fasc. 1 (1969), p. 159--167.
  \item Cruceanu V., {\it Sur l'ensemble des connexions sur un espace fibr\'{e}}, Bull. Math. Soc. Sci. Math. de Roumanie, T. XII (1969), p. 27--34.
  \item Flaschel P., Klingenberg W., {\it Riemannsche Hilbertmannigfaltigkeiten. Periodische Geoda

  tische}, Lecture Notes in Mathematics, T. 282, Springer-Verlag, 1972.
  \item P. de la Harpe, {\it The Clifford algebra and the spinor group of a Hilbert space}, Compositio Math., 25 (1972), p. 245--261.
  \item Kobayashi S.,Nomizu K.,  {\it Foundations of Differential Geometry} I, Interscience Publisher, New York-London, 1963.
  \item Palais R. S., {\it Seminar on the Atiah-Singer index theorem}, Princeton Univ. Press 1965.
  \item Penot J-P., {\it Connexion lin\'{e}aire d\'{e}duite d'une famille de connexions lin\'{e}aires par un foncteur multilin\'{e}aire}, C. R. Acad. Sci. Paris, T. 208, S. A., p. 100--103.
  \item Plymen R. J., {\it Spinors in Hilbert Spaces}, Math. Proc. Cambridge Phil. Soc. 80 (1970), p. 337--347.
  \item Popovici I., {\it Repr\'{e}sentation irr\'{e}ductibles des fibr\'{e}s de Clifford}, Ann. Inst. Henri Poincar\'{e}, S. A. T. XXV, 1976, p. 35--58.
\end{enumerate}

\bigskip

\noindent{\footnotesize{\it Re\c{c}u le 8.II.1978}\hfill {\it Facult\'{e} de Math\'{e}matique}\\}
\vspace*{-.3cm}
\hfill{\footnotesize{\it Universit\'{e} de Jassy, R.S. Rom\^ania}}


\newpage





\def\G{\Gamma}

\runningauthor={M. ANASTASIEI}
\runningtitle={CONSTANT LINEAR CONNECTIONS ON BANACH MANIFOLDS}
\noindent
\baselineskip 8pt
\noindent{\footnotesize{REVUE ROUMAINE DE MATH\'{E}MATIQUES PURES ET APPLIQU\'{E}ES\hfill\break
Tome XXV, No. 1, p. 3-11, Bucharest, 1980}}
\vskip 2cm
\baselineskip 11.5pt plus .15pt
\centerline{\bf\Large CONSTANT LINEAR CONNECTIONS}
\vskip .2cm
\centerline{\bf\Large ON BANACH MANIFOLDS}
\vskip .5cm
\centerline{\bf \footnotesize{BY}}
\vskip .5cm
\centerline{\bf \footnotesize{M. ANASTASIEI and IULIAN POPOVICI}}
\vskip 1cm

The notion of constant linear connection was studied by G. Vranceanu and others from various points of view. An approach of the second author to this subject is used to define and study this notion in the category of analytic Banach manifolds.

\setcounter{section}{0}
\section*{Introduction}

An affine connection on an open subset $U$ of $\R^n$ is well-determined by $n^3$ real functions $\Gamma^k_{ij}$ $(i,j,k=1,...,n)$, defined on $U$. G. Vranceanu has considered the affine connection defined by $\Gamma^k_{ij}=$ constant on $U$ and called it constant affine connection. By a remark of G. Vranceanu and Gr.C. Moisil, in this case $\Gamma^k_{ij}$ define on an $n-$dimensional vector space a structure of an $n-$dimensional algebra $\bfa$ and conversely, the constants of structure of such an algebra $\bfa$ define a constant affine connection $\nabla$ on an open subset of $\bfr^n$ [7]. In this way there appears a correspondence $\nabla\to\bfa$ studied in detail by G. Vranceanu [7], [8] and others. We quote the following interesting result: the constant affine connection $\nabla$ is flat if and only if the algebra $\bfa$ is commutative and associative.

The second author of this paper succeeded to give a global form of this notion and to obtain the global forms of the old results and some new results [5], [6]. His approach to this subject can be used to extend the notion of constant connection to Banach manifolds. This is the purpose of the present paper.

Firstly, some facts about Banach manifolds and linear connections on such manifolds are given.

The notion of constant linear connection is defined for Banach manifolds of class $C^\9$. A theorem which shows that the natural place of this concept is the category of analytic Banach manifolds is proved.

Finally, a splitting of the Banach manifolds with constant linear connections is given and a generalization of the result quoted above is proved.

\section{PRELIMINARIES AND NOTATIONS}
\setcounter{equation}{0}

Let $M$ be a paracompact Banach manifold of class $C^\9$ modeled by the Banach space $\bfm$. Assume that the norm of $\bfm$ is of class $C^\9$ on $\bfm-\{0\}$. It follows that $\bfm$ admits a $C^\9$partition of unity. We remark that the assumption concerning the norm of $\bfm$ is fulfilled if it originates in an inner product on $\bfm$, therefore if $\bfm$ is a Hilbert space. Let us denote by ${\cal F}(M)$ the ring of real functions of class $C^\9$ on $M$ and by ${\cal X}(M)$ the ${\cal F}(M)-$module of sections of class $C^\9$ of the tangent bundle $TM$.

Let $X\in{\cal X}(M)$ be a vector field on $M$ and let $(U,\varphi)$ be a local chart around of $p\in M$. The local section $X|_U$ is well-defined by a $C^\9-$map $X_\varphi$: $\varphi(U)\to\bfm$, called the local representation of $X$. We put $X_{\varphi(p)}=X_\varphi(\varphi(p))$. For another local chart $(V,\psi)$ around of $p$, the local representation $X_\psi$ of $X$ is given by
\begin{equation}
X_{\psi(p)}=D_{\varphi(p)}(\psi\circ\varphi^{-1})(X_{\varphi(p)}),
\end{equation}
where $D_{\varphi(p)}(\psi\circ\varphi^{-1})$ is the Fr\'{e}chet derivative of $\psi\circ\varphi^{-1}$ in the point $\varphi(p)\in\varphi(U\cap V)$. Let $Y$ be another vector field on $M$. The bracket $[X,Y]$ is a vector field whose local representation is given by (see [4])
\begin{equation}
[X,Y]_{\varphi(p)}=D_{\varphi(p)}X_\varphi(Y_{\varphi(p)})-D_{\varphi(p)}Y_\varphi(X_{\varphi(p)}).
\end{equation}

Given a local chart $(U,\varphi)$, we denote by $K(U,\varphi)$ the set of those vector fields on $U$, whose local representations are constant. The set $K(U,\varphi)$ has the following two properties.

\noindent(1.3) \qquad The map $K(U,\varphi)\to T_pM$ given by $X\to X_p$ is an isomorphism of vector spaces for every $p\in U$.

\noindent(1.4) \qquad The bracket $[X,Y]=0$ for every $X,Y\in K(U,\varphi)$.

Every $X\in{\cal X}(M)$ generates a local 1-parameter group $\a_t$, of diffeomorphisms of $M$. A vector field $Y$ is said to be invariant by $X$ if $\a_{t,*}Y=Y$. The Lie derivative of $Y$ with respect to $X$ is given by $(L_XY)_p=\dd\lim_{t\to 0}(Y_p—(\a_{t,*}Y_p))/t$, therefore $Y$ is invariant by $X$ if and only if $L_XY=[X,Y]=0$. We can say that $K(U,\varphi)$ is a set of vector fields on $U$ which are invariant by each other.

Let $\{(U_i,\varphi_i)\}$ be the complete atlas of $M$. By a linear connection $\Gamma$ on $M$ we shall understand (see also [2]) a local connector on $M$, i.e. a collection of $C^\9-$maps $\Gamma_{\varphi_i}:\varphi_i(U_i)\to L^2(\bfm;\bfm)$ such that
\begin{equation}
\begin{array}{l}
\Gamma_{\varphi_j(p)}=D_{\varphi_i(p)}(\varphi_j\circ\varphi^{-1}_i)\circ[D^2_{\varphi_j(p)}(\varphi_i\circ\varphi^{-1}_j)+\\
\\ \hspace*{30mm} + \Gamma_{\varphi_i(p)}(D\varphi_{j(p)}(\varphi_\circ\varphi^{-1}_j),
D_{\varphi_j(p)}(\varphi_i\circ\varphi^{-1}_j)]
\end{array}
\end{equation}
holds for $p\in U_i\cap U_j\neq \emptyset$, where $\Gamma_{\varphi_i(p)}=\Gamma_{\varphi_i}(\varphi_i(p))$.

Given $X,T\in{\cal X}(M)$, condition (1.5) assures that
\begin{equation}
\nabla_X Y\overset{def}= D_{\varphi_i(p)}Y_{\varphi_i}(X_{\varphi_i(p)})+\Gamma_{\varphi_i(p)}(X_{\varphi_i(p)},Y_{\varphi_i(p)})
\end{equation}
defines a new vector field on $M$ which is denoted by $\nabla_XY$ and is called the covariant derivative of $Y$ in the direction of $X$. The map $\nabla:{\cal X}(M)\times{\cal X}(M)\to{\cal X}(M)$ given by $(X,Y)\to\nabla_X Y$ is linear in the first variable and satisfies
\begin{equation}
\nabla_X(Y+Z)=\nabla_X Y+\nabla_X Z,\ \  X,Y,Z\in{\cal X}(M),
\end{equation}
and
\begin{equation}
\nabla_X(fY)=X(f)Y+f\nabla_X Y,\ \ f\in{\cal F}(M),
\end{equation}
therefore it is a covariant differentiation on $M$.

The torsion and the curvature of $\G$ are given by
\begin{equation}
T(X,Y)=\nabla_X Y-\nabla_Y X-[X,Y],
\end{equation}
and
\begin{equation}
R(X,Y)Z=\nabla_X\nabla_Y Z-\nabla_Y\nabla_X Z-\nabla_{[X,Y]}Z,\ X,Y,Z\in {\cal X}(M)
\end{equation}
respectively.

S. Kobayashi and K. Nomizu have defined generalized affine connections as connections in principal fibre bundle of affine frames over $M$. They proved that there exists a one-to-one correspondence between the set of generalized affine connections and the set of pairs $(\G,K)$, where $\G$ is a linear connection and $K$ is a tensor field of type $(1,1)$ (see Ch. III, \S3 of [3]). The generalized affine connection which corresponds to $(\G,I)$, where $I$ is the tensor of Kronecker, was called the affine connection associated to $\G$. This logical distinction between a linear connection and an affine connection can also be made in our context (see [1]). Moreover, the theorem which says that an affine connection is flat (cf. Ch. II, \S9 of [3]) if and only if $R=0$ and $T=0$ is still true.

Let $C$ be a vector field on $M$. A linear connection $\G$ is said to be invariant by $C$ if
\begin{equation}
[C,\nabla_XY]=\nabla_X[C,Y]+\nabla_{[C,X]}Y,\ {\rm{holds\ for\ every\ }} X,Y\in{\cal X}(M).
\end{equation}
It is easy to verify that a linear connection $\G$ is invariant by $C$ if and only if $\nabla_XY$ is invariant by $C$ when $X$ and $Y$ are invariant by $C$.

Finally, we remark that, with minor changes, the results what follow are true without hypothesis of paracompactness of manifolds and even in the case ``no Hausdorf''.

\section{CONSTANT LINEAR CONNECTIONS}\setcounter{equation}{0}

{\bf Definition 2.1.} Let $M$ be a Banach manifold of class $C^\9$. A linear connection $\G$ on $M$ is said to be {\it constant with respect to the local chart} $(U,\varphi)$ of $M$, if $\G$ is invariant by every vector field from $K(U,\varphi)$ i.e.

\begin{equation}
\nabla_XY\in K(U,\varphi)\ {\rm{for\ all\ }}X,Y\in K(U,\varphi).
\end{equation}

Suppose that $\G$ is constant with respect to $(U,\varphi)$. Then a new operation can be defined on $K(U,\varphi)$ by $XY=\nabla_XY$, $X,Y\in K(U,\varphi)$ or
\begin{equation}
X_{\varphi(p)}Y_{\varphi(p)}=\G_{\varphi(p)}(X_{\varphi(p)},Y_{\varphi(p)}).
\end{equation}
It follows from continuity of the bilinear map $\G_{\varphi(p)}$ that $$\|X_{\varphi(p)}Y_{\varphi(p)}\|\leq |\G_{\varphi(p)}|\ \|X_{\varphi(p)}\|\ \|Y_{\varphi(p)}\|,$$ where $\|\cdot\|$ is the norm on $\bfm$ and the $|\cdot|$ is the norm on $L^2(\bfm;\bfm)$, therefore $K(U,\varphi)$ with the product defined by (2.2) is a Banach algebra isomorphic to $\bfm$ as normed linear spaces via the isomorphism $T_pM\cong\bfm$. We denote this Banach algebra by $A(\Gamma,K(U,\varphi))$.

\ms

{\bf Definition 2.2.} Let $\bfa$ be a Banach algebra. A linear connection $\Gamma$ is said to be {\it constant} on a manifold $M$ if there exists an atlas $\a=\{(U_i,\varphi_i)\}$ on $M$ such that $\Gamma$ is constant with respect to each $(U_i,\varphi_i)$ and $A(\Gamma,K(U_i,\varphi_i))$ is isomorphic to $\bfa$ as Banach algebras. The triplet $(M,\Gamma,\bfa)$ will be called a {\it Vranceanu's space}. The atlas $\a$ will be called an {\it atlas adapted} to $(M,\Gamma,\bfa)$ and the atlas $\a$ completed with all local charts with the properties required above will be denoted by $\a^*$ and will be called {\it complete atlas adapted} to $(M,\Gamma,\bfa)$.

\ms

{\bf Proposition 2.1.} {\it Let $(M,\Gamma,\bfa)$ be a Vranceanu's space. The linear connection $\Gamma$ is symmetric $(T=0)$ if and only if $\bfa$ is commutative.}

\ms

\noindent{\it Proof.} The local representation of the torsion $T$ in a local chart $(U,\varphi)$ is $T(X,Y)_{\varphi(p)}=\Gamma_{\varphi(p)}(X_{\varphi(p)},Y_{\varphi(p)})-\Gamma_{\varphi(p)}(Y_{\varphi(p)},X_{\varphi(p)})$, therefore $\Gamma$ is symmetric if and only if $\Gamma_{\varphi(p)}$ are symmetric maps. It follows that $\Gamma$ is symmetric if and only if $A(\Gamma,K(U,\varphi))$ is commutative, therefore if and only if $\bfa$ is commutative. Q.E.D.

In Definition 2.2 the manifold $M$ was assumed of class $C^\9$. The following theorem shows that a manifold of class $C^\9$ with a constant linear connection has a structure of analytic manifold.

\ms

{\bf Theorem 2.1.} {\it Let $M$ be a Banach manifold of class $C^\9$ and let $\Gamma$ be a linear connection on $M$ which is constant with respect to an atlas $\a=\{(U_i,\varphi_i)\}$ and with a Banach algebra $\bfa$. Then $M$ has a structure of analytic manifold $'M$ given by $\a$ and $\Gamma$ induces on $'M$ an analytic connection $'\Gamma$.}

\ms

\noindent{\it Proof.} If $\Gamma$ is constant on $M$ with respect to $\a$ and $\bfa$, the maps $\Gamma_{\varphi_i}$ are necessarily constant maps. Relation (1.5) can be written as follows
\begin{equation}
\begin{array}{c}
D^2_{\varphi_j(p)}(\varphi_i\circ\varphi^{-1}_{j})\hspace*{-0.5mm}=\hspace*{-0.5mm}D_{\varphi_j(p)}(\varphi_i\circ\varphi_j^{-1})\circ\Gamma_{\varphi_i(p)}\hspace*{-0.5mm}-
\hspace*{-0.5mm}\Gamma_{\varphi_i(p)}(D_{\varphi_j(p)}(\varphi_i\circ\varphi^{-1}_j),\\ \\
\hspace*{20mm} D_{\varphi_j(p)}(\varphi_i\circ\varphi^{-1}_j)).
\end{array}
\end{equation}
Let us note $u=D(\varphi_i\circ\varphi^{-1}_j):\varphi_j(U_i\cap U_j)\to L(\bfm;\bfm)$. Then (2.3) becomes
\begin{equation}
D_{\varphi_j(p)}u(h)=u_{\varphi_i(p)}(\G_{\varphi_j(p)}(h,h))-\G_{\varphi_i(p)}(u_{\varphi_j(p)}(h),u_{\varphi_j(p)}(h)),\ h\in\bfm
\end{equation}
where $u_{\varphi_j(p)}=u(\varphi_j(p))$. We put $x=\varphi_j(p)$. The map $u$ is of class $C^\9$ by hypothesis of the theorem. Let us consider the Taylor series of $u$ in a neighborhood of $x$ (see [4], Ch. I, \S4)
\begin{equation}
u(x)+D_x uh+D^2_x uh^2+...+D^n_x uh^n+...,
\end{equation}
where $h^k=(h,h,...,h)\in\bfm^k$.

We shall prove that series (2.5) converges in a small neighborhood of $x$. Firstly, on differentiating by $n-$times the function $u$ and using (2.4) we obtain
\begin{equation}
\begin{array}{l}
D^n_x u(h_1,...,h_n)=D^{n-1}_x u(h_1,...,h_{n-2},\G_{\varphi_j(p)}(h_{n-1},h_n))-\\ \\
-\dd\sum_{k=0}^{n-1}\dfrac{(n-1)!}{k!(n-k-1)!}\G_{\varphi_i(p)}(D^k_x u(h_1,...,h_k),D^{n-k-1}_x u(h_{k+1},...,h_n)),
\end{array}
\end{equation}
where $(h_1,...,h_n)\in \bfm^n$. In what follows all norms will be denoted by $|\cdot|$ and the index $x$ will be omitted since all derivatives are in the point $x$. Equation (2.4) leads to
\begin{equation}
|Duh|\leq |u|\ |\G_{\varphi_j(p)}|\ |h|^2+|\G_{\varphi_i(p)}|\ |u|^2 |h|^2,\ h\in\bfm.
\end{equation}

By a well-known definition $|Duh|=\dd\sup_h\{|Duh|,\ |h|\leq 1\}$. Using (2.7) we arrive at
\begin{equation}
|Du|\leq |u|\ |\G_{\varphi_j(p)}| + |\G_{\varphi_i(p)}|\ |u|^2,
\end{equation}
\begin{equation}
|Du|\leq \lbb|u|\ {\rm{where\ }} \lbb=|\G_{\varphi_j(p)}| + |u|\ |\G_{\varphi_i(p)}|.
\end{equation}
Now we prove by mathematical induction
\begin{equation}
\dfrac{1}{n!}|D^n u|\leq \lbb^n|u|.
\end{equation}
From (2.6) it follows
\begin{equation}
\begin{array}{l}
|D^nu(h_1,...,h_n)|\leq |D^{n-1}u|\ |\G_{\varphi_j(p)}|\ |h_1|....|h_n| +\\ \\
+\dd\sum^{n-1}_{k=0}\dfrac{(n-1)!}{k!(n-k-1)!}|\G_{\varphi_i(p)}|\ |D^k u|\ |D^{n-k-1}u|\ |h_1|...|h_n|.
\end{array}
\end{equation}
Using $|D^nu|=\dd\sup_{(h_1...h_n)}\{|D^nu(h_1,...h_n),\ |h_i|\leq 1,i=1,...,n\}$ and, the inductive hypothesis we obtain $$\begin{array}{c}|D^n u|\leq|D^{n-1}u|\ |\G_{\varphi_j(p)}|+|\G_{\varphi_i(p)}|\cdot\dd\sum_{k=0}^{n-1}\dfrac{(n-1)!}{k!(n-k-1)!}|D^ku|\ |D^{n-k-1}u|\leq\\ \\
\leq (n-1)!\lbb^{n-1}|u|\ |\G_{\varphi_j(p)}|+|\G_{\varphi_i(p)}|\cdot\dd\sum^{n-1}_{k=0}\dfrac{(n-1)!}{k!(n-k-1)!}k!(n-k-1)!\lbb^{n-1}|u|^2,
\end{array}$$ hence $$\dfrac{1}{n!}|D^n u|\leq|u|\lbb^{n-1}\(\dfrac{1}{n}|\G_{\varphi_j(p)}|+|\G_{\varphi_i(p)}|\ |u|\)\leq\lbb^n|u|.$$
The series (2.5) converges if and only if the series $\dd\sum_{n\geq 0}\dfrac1{n!}|D^n uh^n|$ converges. Using (2.10) we obtain
\begin{equation}
dfrac1{n!}|D^n uh^n|\leq\dfrac1{n!}|D^n u|\ |h|^n\leq|u|(\lbb|h|)^n.
\end{equation}
Consequently, by comparison test, series (2.5) converges for $|h|\leq 1/\lbb$. It follows that $u$ is analytic i.e. $\{(U_i,\varphi_i)\}$ defines on $M$ a structure of analytic manifold. Q.E.D.	

\section{SOME TYPES of VRANCEANU'S SPACES}\setcounter{equation}{0}
Let $\bfa$ be a Banach algebra. The law of product on $\bfa$ defines an element $B\in L^2(\bfa;\bfa)$ putting $xy=B(x,y)$, $x,y\in\bfa$. Conversely, every element of $L^2(\bfa;\bfa)$ defines a law of product on the Banach space $\bfa$ which changes it into a Banach algebra. Let $(M,\G,\bfa)$, where $M$ is an analytic manifold, be a Vranceanu's space and let $\{(U_i,\varphi_i)\}$ be an atlas adapted to it. The isomorphism of $A(\G,K(U_i,\varphi_i))$ to $\bfa$ leads via the isomorphisms $A(\G,K(U_i,\varphi_i))\cong T_pM$ and $T_pM\cong\bfm$, to an isomorphism of normed linear spaces $\theta_i:\bfm\to\bfa$ such that $B(\theta_i u,\theta_i v)=\theta_i\G_{\varphi_i}(u,v)$ for $u,v\in\bfm$. The isomorphism $\theta_i$ depends on $(U_i,\varphi_i)$ but it does not depend on the points of $U_i$.

Now, let $\G$ be a certain linear connection on $M$ and let $\{(V_j,\psi_j)\}$, be an analytic atlas of $M$. Suppose that for each chart $(V_j,\psi_j)$ there exists an isomorphism of normed linear spaces $\theta_j:\bfm\to\bfa$ ($\theta_j$ does not depend on points of $V_j$) such that $B(\theta_ju,\theta_jv)=\theta_j\G_{\psi_j(p)}(u,v)$, $u,v\in\bfm$. Then $\G_{\psi_j}$ is constant on $V_j$ and $A(\G,K(V_j,\psi_j))$ is isomorphic to $\bfa$, i.e. the triplet $(M,\G,\bfa)$ is a Vranceanu's space. Therefore we have proved

\ms

{\bf Proposition 3.1.} {\it A triplet $(M,\G,\bfa)$ is a Vranceanu's space if and only if there exists an atlas $\{(V_j,\psi_j)\}$ on $M$ such that for each $(V_j,\psi_j)$ there exists an isomorphism $\theta_j:\bfm\to\bfa$ satisfying
\begin{equation}
B(\theta_ju,\theta_jv)=\theta_j\G_{\psi_j(p)}(u,v),\ \ u,v\in\bfm.
\end{equation}}

\ms

{\it Remarks.} 1) The isomorphisms $\theta_j$ are determined up to an isomorphism of the Banach algebra $\bfa$, i.e. if $\theta_j$ satisfies (3.1), then $h\circ\theta_j$, where $h$ is an isomorphism of $\bfa$, satisfies (3.1), too.

2)	The atlas $\beta=\{(V_j,\psi_j)\}$ from the above proposition is an atlas adapted to $(M,\G,\bfa)$. It will be called $\theta-$atlas and completed with all charts which satisfy (3.1) will be denoted by $\beta^*$ and will be called the complete $\theta-$atlas of $(M,\G,\bfa)$.

3)	If $(V_j,\psi_j)$ is a chart from $\beta^*$, then $(V_j,g\circ\psi_j)$, where $g\in GL(\bfm)$, satisfies (3.1) because we can write $\theta_{gj}=\theta_{gj}\circ\theta^{-1}_j\circ\theta_j$, where $\theta_{gj}:\bfm\to\bfa$ corresponds to $(V_j,g\circ\psi_j)$ and $\theta_{gj}\circ\theta^{-1}_j$ is an isomorphism of $\bfa$. It follows that $\G_{g\circ\psi_j}$ is constant on $V_j$. This shows that if we add to $\beta^*$ all charts of the form $(V_j,g\circ\psi_j)$ with $g$ from $GL(\bfm)$ and $(V_j,\psi_j)$ from $\beta^*$, we obtain the complete atlas adapted to $(M,\G,\bfa)$ (denoted above by $\a^*$).

Using Proposition 3.1 we obtain the following corollary of Theorem 2.1.

\ms

{\bf Corollary 3.1.} {\it Let $M$ be a Banach manifold of class $C^\9$ equipped with a linear connection $\G$ and let $\bfa$ be a Banach algebra. If there exists an atlas $\{(V_j,\psi_j)\}$ with the property that for each chart $(V_j,\psi_j)$ there exists an isomorphism of normed linear spaces $\theta_j:\bfm\to\bfa$ such that (3.1) to be true, then $\{(V_j,\psi_j)\}$ gives to $M$ a structure of analytic manifold.}

\ms

{\bf Definition 3.1.} An atlas $\{(U_i,\psi_i)\}$ of the analytic manifold $M$ is said to be {\it affine} if $\varphi_j\circ\varphi_i^{-1}:\varphi_i(U_i\cap U_j)\to\varphi_j(U_i\cap U_j)$ has the form
\begin{equation}
(\varphi_j\circ\varphi^{-1}_i)(u)=S(u)+u_0,\ {\rm\ where\ } u_0,u\in\bfm\ {\rm{and\ }}S\in GL(\bfm)),
\end{equation}
for all pairs $(i,j)$ with $U_i\cap U_j\neq\emptyset$.

We remark that in this case $D(\varphi_j\circ\varphi_i^{-1})=S$ and $D^2(\varphi_j\circ\varphi_i^{-1})=0$. Conversely, an atlas $\{(U_i,\varphi_i)\}$ on $M$ which satisfies
\begin{equation}
D^2(\varphi_j\circ\varphi^{-1}_i)=0\ {\rm{on\ }}U_i\cap U_j\neq\emptyset\ {\rm{for\ all\ pairs\ }}(i,j),
\end{equation}
is affine because the general solution of equation (3.3) is (3.2).

Let $(M,\G,\bfa)$ be a Vranceanu's space and let $\{(V_i,\psi_i)\}$ be a $\theta-$atlas. Suppose that $\{(V_i,\psi_i)\}$ is affine, therefore $D(\psi_j\circ\psi^{-1}_i)=S_{ji}\in GL(\bfm)$. Then (3.2) becomes
\begin{equation}
D_{\psi_i(p)}(\psi_j\circ\psi^{-1}_i)(\G_{\psi_i}(u,v))=\G_{\psi_j}(D_{\psi_i(p)}(\psi_j\circ\psi_i^{-1})u,D_{\psi_i(p)}(\psi_j\circ\psi_i^{-1})v),
\end{equation}
or
\begin{equation}
S_{ji}(\Gamma_{\psi_i}(u,v))=\G_{\psi_j}(S_{ji}u,S_{ji}v),\ \ u,v\in\bfm.
\end{equation}
Using (3.1) we obtain
\begin{equation}
S_{ji}\theta^{-1}_i B(\theta_i u,\theta_i v)=\theta^{-1}_j B(\theta_j S_{ji}u,\theta_j S_{ji}v), \ \ u,v\in\bfm.
\end{equation}
If we put $u=\theta_i^{-1}u'$, $v=\theta_i^{-1}v'$, in (3.6) we arrive at
\begin{equation}
\theta_j S_{ji}\theta^{-1}_i B(u',v')=B(\theta_j S_{ji}\theta^{-1}_i u',\theta_j S_{ji}\theta_i^{-1}v')\ \ u',v'\in\bfa,
\end{equation}
therefore $\overline{S}_{ji}=\theta_j S_{ji}\theta^{-1}_i$ is an isomorphism of $\bfa$.

We denote by $G$ the subset of $GL(\bfm)$ whose elements are of the form
\begin{equation}
S_{ji}=\theta_j^{-1}\overline{S}_{ji}\theta_i,\ {\rm{where\ }}\overline{S}_{ji}\ {\rm{is\ an\ isomorphism\ of\ }}\bfa.
\end{equation}

\ms

{\bf Theorem 3.1.} {\it Let $M$ be an analytic Banach manifold endowed with an affine atlas $\beta=\{(U_i,\varphi_i)\}$ and let $\bfa$ be a Banach algebra. If assume that for each chart $(U_i,\varphi_i)$ there exists an isomorphism of normed spaces $\theta_i:M\to A$ which does not depend on points from $U_i$ then the following statements are equivalent:

a)	there exists a linear connection $\Gamma$ such that $(M,\G,\bfa)$ is a Vranceanu's space with $\beta$ as $\theta-$atlas.

b)	every change of charts from $\beta$ is the composition of a translation on $\bfm$ and an element of $G$.}

\ms

\noindent{\it Proof.} Assuming a), since $\beta$ is affine, for $(U_i,\varphi_i)$ and $(U_j,\varphi_j)$ with $U_i\cap U_j\neq\emptyset$ we have $(\varphi_j\circ\varphi^{-1}_i)(u)=S_{ji}(u)+u_0$, therefore $\varphi_j\circ\varphi_i^{-1}=T_{u_0}\circ S_{ji}$, where $T_{u_0}(u)=u+u_0$ is the translation by $u_0$. By the considerations made above, $S_{ji}\in G$, therefore b) follows.

Let us suppose b). For each $(U_i,\varphi_i)$, we define $\G_{\varphi_i}$ by
\begin{equation}
\G_{\varphi_i}(u,v)=\theta^{-1}_i B(\theta_i u,\theta_i v)\ \ u,v\in\bfm.
\end{equation}
Let us prove that $\{\G_{\varphi_i}\}$ is a local connector. Since $\beta$ is affine we must verify (3.5), where $S_{ji}\in G$. If we replace $\G_{\varphi_i}$ and $\G_{\varphi_j}$ given by (3.9) in (3.5) we obtain (3.6) which is equivalent to (3.7). But (3.7) is true because $S_{ji}\in G$. From (3.9) and Proposition 3.1 it follows that $(M,\G,\bfa)$ is a Vranceanu's space with $\b$ as $\theta-$atlas. Q.E.D.

\ms

{\it Remark.} The connection $\G$ from a) of Theorem 3.1 is unique by Proposition 3.1.

\ms

{\bf Definition 3.2.} A Vranceanu's space $(M,\G,\bfa)$ will be called of the {\it first kind}, {\it second kind} or {\it third kind} if the complete atlas $\a^*$ adapted to it, satisfies the following conditions, respectively:

1)	$\a^*$ is affine,

2)	$\a^*$ is not affine but contains an affine atlas of $M$,

3)	$a^*$ does not contain any affine atlas of $M$.

Using Proposition 3.1 and the remark which follows it we obtain

\ms

{\bf Proposition 3.2.} {\it Let $\b$ be a $\theta-$atlas of the Vranceanu's space $(M,\G,\bfa)$. Then $(M,\G,\bfa)$ is of the first kind, second kind or third kind if $\b$ satisfies 1), 2) or 3) from Definition 3.2, respectively.}

On $\bfa$ we can consider a new structure of Banach algebra given by
\begin{equation}
B^s(x,y)=\dfrac12(B(x,y)+B(y,x)),\ \ x,y\in\bfa.
\end{equation}
We denote this new Banach algebra by $^s\bfa$ and we remark that $^s\bfa$ is commutative.

Let $\G$ be a certain linear connection on $M$ and let $\{\G_{\varphi_i}\}$ be its local connector. For each $\varphi_i$ let $^s\G_{\varphi_i}$ be given by
\begin{equation}
^s\G_{\varphi_i}(u,v)=\dfrac12[\G_{\varphi_i}(u,v)+\G_{\varphi_i}(v,u)]\ \ u,v\in M.
\end{equation}
It is easy to check that $\{^s\G_{\varphi_i}\}$ is a local connector. We denote by $^s\G$ the linear connection given by $\{^s\G_{\varphi_i}\}$ and by $^s\nabla$ the covariant differentiation associated to $^s\G$. It follows easily
\begin{equation}
2^s\nabla_X Y=\nabla_X Y+\nabla_Y X-[X,Y]\ \ X,Y\in{\cal X}(M).
\end{equation}

\ms

{\bf Proposition 3.3.} {\it If the triplet $(M,\G,\bfa)$ is a Vranceanu's space then $(M,^s\G,^s\bfa)$ is also a Vranceanu's space. Moreover, we have

1)	If $(M,^s\G,^s\bfa)$ is of the first kind (third kind) then $(M,\G,\bfa)$ is also of the first kind (third kind).

2)	If $(M,\G,\bfa)$ is of the second kind, then $(M,^s\G,^s\bfa)$ is of the second kind.}

\ms

\noindent{\it Proof.} Let $\b=\{(U_i,\varphi_i)\}$ a $\theta-$atlas of $(M,\G,\bfa)$, therefore $\theta_i\G_{\varphi_i}(u,v)=B(\theta_i u,\theta_i v)$ for every $i$ and $u,v\in\bfm$. It follows easily
that $\theta_i\ ^s\G_{\varphi_i}(u,v)=B^s(\theta_i u,\theta_i v)$, therefore $(M,^s\G,^s\bfa)$ is a Vranceanu's space with $\b$ as $\theta-$atlas. Let $\b^*$ and $^s\b^*$ be the complete $\theta-$atlas of $(M,\G,\bfa)$ and $(M,^s\G,^s\bfa)$, respectively. From $\b^*\subset ^s\b^*$ and Proposition 3.1 follow easily 1) and 2). Q.E.D.

\ms

{\it Remark.} The inclusion $\b^*\subset ^s\b^*$ shows also that if $(M,\G,\bfa)$ is of the first kind, then $(M,^s\G,^s\bfa)$ is of the first kind or of the second kind. Also, if $(M,\G,\bfa)$ is of the third kind, then $(M,^s\G,^s\bfa)$ if of the second kind or of the third kind.

The local representation of the curvature tensor of a linear connection $\G$ on $M$ is given by
\begin{equation}
\begin{array}{c}
R_x(u,v)w=D_x\G_{\varphi_i}(u)(v,w)-D_x\G_{\varphi_i}(v)(u,w)+\\ \\
+\G_x(u,\G_x(v,w))-\G_x(v,\G_x(u,w))\\ \\
u,v,w\in\bfm,\ x=\varphi_i(p),\ p\in M.
\end{array}
\end{equation}

\ms

{\bf Proposition 3.4.} {\it Let $(M,\G,\bfa)$ be a Vranceanu's space. The affine connection $\G'$ associated to $\G$ is plate if and only if $\bfa$ is associative and commutative.}

\ms

\noindent{\it Proof.} Assume that $\G'$ is flat. This is equivalent to $T=0$ and $R=0$. From $T=0$ it follows $\G_x(u,v)=\G_x(v,u)$ for $u,v\in\bfm$. From $R=0$ it follows
\begin{equation}
\G_x(u,\G_x(v,w))=\G_x(v,\G_x(u,w))\ u,v,w\in\bfm.
\end{equation}
Using (3.1) we obtain $B(u', v')=B(v',u')$ and $B(u',B(v', w'))=B(v',B(u',w'))$, where $u'=\theta_i u$, $v'=\theta_i v$, $w'=\theta_i w$. Using the first, which says that $\bfa$ is commutative, in the second we obtain $B(u',B(w',v'))=B(B(u',w'),v')$ i.e. $\bfa$ is associative. Conversely, if $\bfa$ is commutative and associative, using (3.1) we obtain easily that $\G_{\varphi_i}$ are symmetric and (3.14) i.e. $\G'$ is flat. Q.E.D.

\ms

{\it Remark.} Proposition 3.4 is the generalization of a result due to G. Vranceanu [8].

\bigskip

{\bf REFERENCES}

\begin{enumerate}
  \item M. Anastasiei, {\it Generalized affine connections on Banach manifolds} (to appear).
  \item P. Flaschel and W. Klingenberg, {\it Riemannische Hilbert Mannigfaltigkeiten. Periodische Geod\"{a}tische}, Lecture Notes in Math., 282, Springer-Verlag, 1972.
  \item S. Kobayashi and K. Nomizu, {\it Foundations of differential geometry} I, II. John Wiley \& Sons, New York - London, 1963, 1969.
  \item S. Lang, {\it Introduction aux vari\'{e}t\'{e}s diff\'{e}rentiables}, Dunod, Paris, 1967.
  \item I. Popovici, {\it Contributions a l'\'{e}tude des espaces \`{a} connexion constante}, Rev. Roum. Math. Pures et Appl., 33, 1978, 1211-1225.
  \item I. Popovici, {\it Sur les espaces \`{a} connexion constante attach\'{e}s aux alg\`{e}bres de Jordan simples et centrales} I, An. st. Univ. ``Al. I. Cuza'' Iasi, Mat., 25 (1979) 367-382.
  \item G. Vranceanu, {\it Le\c{c}ons de g\'{e}om\'{e}trie diff\'{e}rentielle} II, Edition de l'Acad. de R. P. Roumanie, Bucarest, 1957.
  \item G. Vranceanu, {\it Propri\'{e}t\'{e}s globales des espaces \`{a} connexion affine}, Bull. Math. Soc. Sci. Math. Phys., R. P. Roumanie, 1958, 2, 475-478.
\end{enumerate}

\noindent{\footnotesize{\it Received December 16, 1977}}\hfill {\it\footnotesize{Faculty of Mathematics}}\\ \hspace*{8cm}{\footnotesize{\it ``Al. I. Cuza'' University of Iasi}}\\ \hspace*{10cm}{\footnotesize{\it and}} \\  \hspace*{8.8cm}{\footnotesize{\it Institute of Mathematics}}\\ \hspace*{9.7cm}{\footnotesize{\it Bucharest}}

\newpage





\runningauthor={M. ANASTASIEI}
\runningtitle={GENERALIZED AFFINE CONNECTIONS ON BANACH MANIFOLDS}
\noindent
\baselineskip 8pt
\noindent{\footnotesize{ANALELE {\c S}TIIN{\c T}IFICE ALE UNIVERSIT{\u A}{\c T}II
``AL.~I.~CUZA" IA{\c S}I\hfill\break Tomul XXVI, s.I a, Matematic\u a, 1980, f.2., 381--388}}
\vskip 2cm
\baselineskip 11.5pt plus .15pt
\centerline{\bf\Large GENERALIZED AFFINE CONNECTIONS}
\vskip .2cm
\centerline{\bf\Large ON BANACH MANIFOLDS\footnote[1]{Communicated at the National Symposium on Theory of Relativity, April 25-28, 1979, Ia\c{s}i, Rom\^{a}nia}}
\vskip .5cm
\centerline{\bf\footnotesize{BY}}
\vskip .5cm
\centerline{\bf M. ANASTASIEI}
\vskip 1cm

The theory of nonlinear connections in the category of Banach vector bundles has been developed by J. Vilms ([7], [8]). A class of nonlinear connections, called {\it homogeneous connections}, is of great importance in the theory of Finsler connections ([3], [5]).

The purpose of this paper is the study of another class of nonlinear connections, called {\it generalized affine connections} (g.a.c., for short). The term agrees with the one used in [4, p. 127]. In the first section some new results regarding the nonlinear connections are given. The second section contains the definition of g.a.c. and some of their properties (associated linear connections, geodesics and others). The flat g.a.c. are studied in the third section.

\setcounter{section}{0}
\section{Nonlinear connections}\setcounter{equation}{0} Let $M$ be a paracompact manifold of class $C^\9$ (smooth), modeled by the Banach space $\bfm$ and let $p:E\to M$ be a smooth vector bundle of fiber type a Banach space $\bfe$. Denote by $p^{-1}TM\to E$ the pull-back by $p$ of the tangent bundle $\pi:TM\to M$ and by $p!=(Tp,\tau)$, where $Tp$ is the tangent map to $p$ and $\tau:TE\to E$ is the tangent bundle to the manifold $E$. The map $Tp:TE\to TM$ gives to $TE$ a second (different) structure of vector bundle.

A smooth nonlinear connection is a smooth splitting of the following exact sequence
\begin{equation}
0\to VE\overset{i}\longrightarrow TE\overset{p!}\longrightarrow p^{-1} TM\to 0
\end{equation}
of vector bundles over $E$. Here $VE:=\ker(p!)=\ker(Tp)$ denotes the vertical subbundle of $TE$ and $i$ is the inclusion map.

The vertical subbundle $VE\to E$ is canonically isomorphic to $p^{-1}E\to E$ (the pull-back of $E$ by $p$). Hence, there exists a canonical morphism (over $P$) $r:VE\to E$ of vector bundles, isomorphic on the fibres.

A splitting of the exact sequence (1.1), i.e. a nonlinear connection is given by a smooth morphism $V:TE\to VE$, such that $V\circ i={\rm{id}}|VE$, or equivalently, a smooth morphism $W:p^{-1}TM\to TE$ such that $p!\circ W={\rm{id}}|p^{-1}TM$. Moreover, we have $i\circ V+W\circ p!={\rm{id}}|TE$. This implies $TE=VE+HE$, where $HE=\ker V={\rm{im}}W$. Obviously, $HE$ is isomorphic to $p^{-1}TM$ as vector bundles. The morphism (over $p$) $K:=r\circ V:TE\to E$ is called the connection map and $v=i\circ V$, $h=W\circ p!$ are called vertical and horizontal projections, respectively. The morphism $J=i\circ p!$ of $TE$ satisfies $J^2=0$ since $p!\circ i=0$, therefore $J$ defines an almost tangent structure on $E$. Obviously, $JV(E)=0$ and Im$J=VE$. The morphism $\g=2h—I$, where $I$ is the identity on $TE$, satisfies
\begin{equation}
J\circ\g=J,\ \g\circ J=-J.
\end{equation}

Conversely, a morphism $\g$ satisfying (1.2) determines a unique splitting of the exact sequence (1.1), i.e a nonlinear connection on $p:E\to M$. Indeed, let $W'$ be any right splitting map of the sequence (1.1) ($W'$ exists if $M$ admits smooth partitions of unity). We put $W-hW'$, where $2h=I+\g$. The morphism $W$ does not depend on $W'$ and it is easy to check, using (1.2), that $p!\circ W={\rm{id}}_{p^{-1}TM}$. Therefore, we have the following definition of the nonlinear connections, equivalently to that previously given.

{\bf Definition 1.1.} A nonlinear connection on $p:E\to M$ is a smooth morphism $\g$ of $TE$ (over id$|_E$) satisfying (1.2).

The Definition 1.1 generalizes a definition of nonlinear connections on finite dimensional manifolds given by J. Grifone [3]. As in finite dimensional case (see [3]), one can prove the following.

{\bf Theorem 1.1.} {\it A smooth morphism $\g$ of $TE$ is a nonlinear connection on $p:E\to M$ if and only if it defines an almost product structure on $E(\g\circ\g=I)$ such that for every $u\in E$, the eigenspace of $\g_u$ (the restriction of $\g$ to $p^{-1}u)$) which corresponds to the eigenvalue $—1$ be $V_uE$.}

\section{Generalized affine connections}
\setcounter{equation}{0}

Let $\bff$ be a Banach space. The map $\to :\bff\times\bff\to\bff$ given by $(u,v)\to\overrightarrow{uv}=v-u$ defines the so-called canonical affine structure on $\bff$. Every vector bundle can be considered as an affine bundle if one considers its fibers with the canonical affine structure.

Let $\bff'$ be another Banach space. A map $t:\bff\to\bff'$ is said to be affine if $t(u)=T(u)+t(0)$ for every $u\in\bff$, where $T:\bff\to\bff'$ is a linear map. If we regard $\bff$ and $\bff'$ as affine spaces, the map $t$ is affine if and only if it is an affine morphism.

Given two vector bundles $E\to M$ and $E'\to M'$, a map $h:E\to E'$ which preserves the fibers is said to be affine if it is smooth and its restrictions to fibers are affine. Of course, $h$ can be considered as a morphism in the category of affine bundles.

{\bf Definition 2.1.} A nonlinear connection on $p:E\to M$ will be called generalized affine connection (briefly g.a.c.) if its connection map, denoted above by $K$, is an affine map with respect to the structure of vector bundle of $TE$ given by $Tp:TE\to TM$.

An examination of the local situation will be suitable to lead us to the essential properties of g.a.c. Let $(U,\varphi)$ be a local chart on $M$. We identify $U$ with $\varphi(U)$ and, restricting $U$ if necessary, suppose that there exists a bundle chart $U\times\bfe\cong\bfe|_U$. Then the tangent map gives a local chart $U\times\bfe\times\bfm\times\bfe\cong TE|_U$ and the sequence (1.1) restricted to $U$ becomes
\begin{equation}
0\to U\times\bfe\times 0\times\bfe\overset{i}\longrightarrow U\times\bfe\times\bfm\times\bfe\overset{p!}\longrightarrow U\times\bfe\times\bfm\to 0,
\end{equation}
where $p!(x,a,\lbb,b)=(x,a,\lbb)$, $x\in U$, $\lbb\in\bfm$, $a,b\in\bfe$.

The map $T_p$ is locally given by $T_p(x,a,\lbb,b)=(x,\lbb)$. Therefore the fibers of bundle $T_p:TE\to TM$ are isomorphic to $x\times\bfe\times\lbb\times\bfe\cong\bfe^2$. J. Vilms has proved (see [7]) the following

{\bf Lemma.} {\it A morphism (over $p$) $K:TE\to E$ is the connection map of a nonlinear connection on $p:E\to M$, if and only if it is locally given by
\begin{equation}
K(x,a,\lbb,b)=(x,b+\oo(x,a)\lbb),\ x\in U,\ \lbb\in\bfm,\ a,b\in\bfe,
\end{equation}
where $\oo:U\times\bfe\to L(\bfm,\bfe)$ is smooth.}

Further we prove

{\bf Lemma 2.1.} {\it A morphism (over $p$) $K:TE\to E$ is the connection map of a g.a.c. if and only if it is locally given by
\begin{equation}
K(x,a,\lbb,b)=(x,b+\G(x)(a,\lbb)+A(x)\lbb),
\end{equation}
where $\G:U\to L^2(\bfe,\bfm;\bfe)$ and $A:U\to L(\bfm,\bfe)$ are smooth maps.}

\noindent{\it Proof.} Let $K$ be the connection map of a g.a.c. By Definition 2.1. the map $(x,a,\lbb,b)\to(x,b+\oo(x,a)\lbb)$ must be affine on $Tp-$fibers. Consequently, the map $(a,b)\to b+\oo(x,a)\lbb$ of $\bfe\times\lbb\times\bfe\to\bfe$ must be affine with respect to both variables. Being linear, hence affine with respect to $b$, it remains to be affine with respect to $a$. This happens if and only if there exists a smooth map $\wt{\oo}:U\to L(\bfe,L(\bfm,\bfe))$ such that $\oo(x,a)=\wt{\oo}(x)(a)+\oo(x,0)$. We put $A(x)=\oo(x,0)$. Since $L(\bfe, L(\bfm,\bfe))\cong L^2(\bfe,\bfm;\bfe)$, $\wt{\oo}$ determines a unique smooth map $\G:U\to L^2(\bfe,\bfm;\bfe)$ such that $\wt{\oo}(x)(a)\lbb=\G(x)(a,\lbb)$. Therefore, $\oo(x,a)\lbb=\G(x)(a,\lbb)+A(x)\lbb$ and (2.3) follows from (2.2).

The maps $\G$ and $A$ will be called local components of the g.a.c..

{\it Remarks.} If the connection map of a nonlinear connection is linear on $Tp-$fibers, the connection becomes a linear connection. A g.a.c. is linear if and only if $A$ vanishes on $U$.

When $\oo(x,a)$ is 1-homogeneous with respect to $a$, or equivalently, $K$ is 1-homogeneous on $Tp-$fibers, the nonlinear connection is called homogeneous connection. The class of homogeneous connections is used in the theory of Finsler connections (see [3], [5]). In the definition of a homogeneous connections one needs a greater generality, namely the smoothness of it must be assumed only on $E-0$, otherwise it becomes linear. Our considerations from the first section remain true in such a generality (with the appropriate modifications), but it is not necessary for the theory of g.a.c.

Let $(U,\varphi)$ and $(V,\psi)$ two local charts on $M$, such that $U\cap V\neq\emptyset$. We put $f=\psi\circ\varphi^{-1}$. If $\Phi:p^{-1}(U)\to U\times\bfe$ and $\Psi:p^{-1}(V)\to V\times\bfe$ are bundle local charts, we denote by $B:U\cap V\to L(\bfe,\bfe)$ the map $x\to B(x)=\Psi\circ\Phi^{-1}$. In this notations the change of bundle local charts on $E$ can be written as $(x,a)\to(f(x), B(x)a)$, $a\in E$ and the change of bundle local charts on $Tp:TE\to TM$ induced by it, is given by $(x,a,\lbb,b)\to (f(x)),B(x)a,\pp f(x),\pp B(x)(\lbb)a+B(x)b)$, $x,\lbb\in\bfe$, $a,b\in\bfe$, where $\pp$ means Fr\'{e}chet differentiation. Let us denote by $\overline{\G}$ and $\overline{A}$ the local components of g.a.c. with respect to the local chart $(V,\psi)$. Using (2.3) and the expressions of changes of bundle local charts given above, we find the following transformation rule for the local components of a g.a.c.
\begin{equation}
\begin{array}{c}
\overline{\G}(f(x))(B(x)a,\pp f(x)\lbb)+A(f(x))\pp f(x)=B(x)\G(x)(a,\lbb)+\\ \\
+B(x)A(x)\lbb-\pp B(x)(\lbb)a,\ \ x\in U\cap V,\ a,b\in\bfe,\ \lbb\in\bfm.
\end{array}
\end{equation}

For $a=0$, the relation (2.4) becomes
\begin{equation}
A(f(x))\pp f(x)\lbb=B(x)A(x)\lbb,\ \ x\in U\cap V,\ \ \lbb\in\bfm,
\end{equation}
which, used in (2.4), leads to
\begin{equation}
\overline{\G}(f(x))(B(x)a,\pp f(x)\lbb)=B(x)\G(x)(a,\lbb)-\pp B(x)(\lbb)a.
\end{equation}

The relation (2.5) shows that $A$ is the local part of a section, denoted also by $A$, of the vector bundle $L,(TM,E)\to M$ (of fiber $L(T_xM,E_x)$, $x\in M$). The relation (2.6) is just the transformation rule of the local connector of a linear connection on $E$. Therefore, $\G$ defines a linear connection on $E$, which will be denoted also by $\G$. Conversely, a section $A$ of the vector bundle $L(TM,E)\to M$ and a linear connection $\G$ on $E$ determine a unique g.a.c. via their local parts. So, we have proved

{\bf Theorem 2.1.} {\it Let $p:E\to M$ be a Banach vector bundle. There exists a one-to-one correspondence between the set of g.a.c. on $p:E\to M$ and the pairs $(\G, A)$, where $\G$ is a linear connection on $p:E\to M$ and $A$ is a section of $L(TM,E)\to M$.}

Let us denote by ${\cal X}_E(M)$ the set of smooth sections of $p:E\to M$ and let us put ${\cal X}(M)={\cal X}_{TM}(M)$. Now, let us regard $p:E\to M$ as an affine bundle. Its fiber in $x\in M$ will be denoted by $^aE_x$ and $x$, identified with zero of $E_x$ will be called the {\it contact point} of $M$ with $^aE_x$. A map $P:M\to E$ given by $x\to P_x\in ^aE_x$ is by definition of class $C^\9$, if the map $a:M\to E$ defined by $x\to a_x=\overrightarrow{x}P_x\in E_x$ is of class $C^\9$.

Such a map $P$ of class $C^\9$ will be called a point field. We denote by ${\cal P}_E(M)$ the set of point fields of class $C^\9$ and we put ${\cal P}(M)={\cal P}_{TM}(M)$.

By Theorem 2.1 a g.a.c. is well determined by the pair $(\G,A)$. But a linear connection, defines a covariant differentiation i.e. a map $\nabla:{\cal X}(M)\times{\cal X}_E(M)\to{\cal X}_E(M)$ with the following properties:
\begin{equation}
\nabla(X,\a a)=\a\nabla(X,a)+X(\a)\cdot a,
\end{equation}
\begin{equation}
\nabla(X,a+b)=\nabla(X,a)+\nabla(X,b),
\end{equation}
\begin{equation}
\nabla(\a X+\b Y,a)=\a\nabla(X,a)+\b\nabla(Y,a)\ \ X,Y\in{\cal X}(M),\ a,b\in{\cal X}_E(M).
\end{equation}
and $\a,\b\in{\cal F}(M)$ — the module of real functions defined on $M$.

Using $\nabla$ and $A$ we shall define an analogue of $V$ for a g.a.c., namely: $D:{\cal P}(M)\times{\cal P}_E(M)\to{\cal P}_E(M)$ given by
\begin{equation}
\overrightarrow{QD(P,Q)}=\nabla(X,a)+A(X),\ {\rm{where\ }}X=\overrightarrow{xP}\ {\rm{and\ }} a=\overrightarrow{xQ}.
\end{equation}

{\bf Theorem 2.2.} {\it The map $D$ associated to a g.a.c. as above, has the following properties:
\begin{equation}
D(P,\a Q+\b R)=\a D(P,Q)+\b D(P,R)+X(\a)P+X(\b)Q
\end{equation}
\begin{equation}
D(\a P+\b P',Q)=\a D(P,Q)+\b D(P',Q),
\end{equation}
\begin{equation}
D(x,Q)=Q,\ {\rm{for\ }}\a,\b\in{\cal F}(M),\ \a+\b=1,\ P\in{\cal P}(M),\ Q,R\in{\cal P}_E(M)
\end{equation}
and $X=\overrightarrow{xP}$ where $x$ is the contact point field.}

\noindent{\it Proof.} To prove (2.11) we remark that it is equivalent to
$$\overrightarrow{xD(P,\a Q+\b R)}=\a\overrightarrow{xD(P,Q)}+\b\overrightarrow{xD(P,R)}+X(\a)\overrightarrow{xP}+X(\b)\overrightarrow{xQ},\leqno(*)$$ or
$$\overrightarrow{(\a Q+\b R)D(P,\a Q+\b R)}=\a\overrightarrow{QD(P,Q)}+\b\overrightarrow{RD(P,R)}+X(\a)\overrightarrow{xP}+X(\b)\overrightarrow{xQ}.
$$
Now we can use (2.10) to obtain
$$\nabla(X,\a a+\b b)+A(X)=
\a\nabla(X,a)+\a A(X)+\b\nabla(X,b)+\b A(X)+X(\a)a+X(\b)b,\leqno(**)$$ where $a=\overrightarrow{xQ},\ b=\overrightarrow{xR}$.

Since $\alpha+\beta=1$, (**) is true by virtue of (2.7) and (2.8). The proof of (2.11) follows the pattern of the previous proof. The property (2.13) is equivalent to $\overrightarrow{xD(x,Q)}=\overrightarrow{xQ}$, or $\overrightarrow{QD(x,Q}=0$. Using again (2.10) we obtain $V(0,a)+A(0)=0$, which is obviously true. Conversely, given a map $D:{\cal P}(M)\times{\cal P}_E(M)\to {\cal P}_E(M)$ which satisfies (2.11)-(2.13) we can derive from it a covariant differentiation $\nabla$ and a section of $L(TM,E)\to M$, as follows:
\begin{equation}
\nabla(X,a)=\overrightarrow{QD(P,Q)}-\overrightarrow{xD(P,x)},\ A(X)=\overrightarrow{xD(P,x)},
\end{equation}
where $$X=\overrightarrow{xP},\ \ a=\overrightarrow{xQ}.$$

But the covariant differentiation $\nabla$ does not define $\G$, such that in our framework the map $D$ does not determine a g.a.c. This happens when the dimension of $M$, as well as of $E$ is finite (see [2]). It is also easy to prove, using (2.14), the following

{\bf Theorem 2.3.} {\it A generalized affine connections on $M$ is affine if and only if $D(P,x)=P$ for every $P\in{\cal P}(M)$.}

We obtain a g.a.c. on $M$ when $E=TM$. Every section of $L(TM,TM)\to M$ is a tensor field of type (1,1). Therefore, we have:

{\bf Corollary 2.1.} {\it Let $M$ be a Banach manifold. There exists a one-to-one correspondence between the set of g.a.c. on $M$ and the set of pairs consisting of a linear connection on $M$ and a tensor field of type (1,1) on $M$.}

The g.a.c. $^a\G$ which corresponds to $(\G, I)$, where $I$ is the Kronecker tensor field, will be called affine connection.

Let $\wt{\G}$ be a g.a.c. on the Banach manifold $M$ and let be $K:TTM\to TM$ its connection map. The map $D$ defined above, induces a covariant differentiation $\wt{\nabla}:{\cal X}(M)\times{\cal X}(M)\to{\cal X}(M)$, which can be expressed as $\wt{\nabla}(X,Y)=K\circ TY(X)$, where $TY$ is the tangent map to $Y:M\to TM$. Indeed, the local parts of $\wt{\nabla}$ and $D$ are given by the right part of the equality
\begin{equation}
\wt{\nabla}_X Y|_\varphi=\pp Y_\varphi(X_\varphi)+\G_\varphi(Y_\varphi,X_\varphi)+A(X_\varphi),
\end{equation}
where $X_\varphi,Y_\varphi$ are the local parts of $X$ and $Y$, respectively and $\G_\varphi$ and $A$ are the local components of g.a.c.

Let $c:(a,b)\subset R\to M$ be a smooth curve on $M$ and let $Tc:(a,b)\times R\to TM$ be its tangent map. We denote by $\dot{c}$ the vector field on $c(a,b)\subset M$ given by $c(t)\to\dot{c}(t)$, where $t\in(a,b)$ and $\dot{c}(t) = Tc(t,1)$. In a local chart $(U,\varphi)$ with $U\cap c(a,b)\neq\emptyset$, the $\dot{c}(t)$ is given by
\begin{equation}
\dot{c}(t)=(c(t),\ \pp c_\varphi(t)(1)),
\end{equation}
where $c_\varphi=\varphi\circ c$.

A curve $c$ will be called a geodesic of the g.a.c. $\wt{\G}$ if $\wt{\nabla}_{\dot{c}}\dot{c}=0$. The local component $c_\varphi$ of a geodesic of $\wt{\G}$ satisfies the following differential equation
\begin{equation}
\pp^2 c_\varphi(t)+\G_\varphi(\pp c_\varphi(t), \pp c_\varphi(t))+A(\pp c_\varphi(t))=0.
\end{equation}

From the theory of differential equations, it follows the local existence and the uniqueness of a geodesic with the initial conditions $c_\varphi(t_0)=c_0\in M$ and $\pp c_\varphi(t_0)(1)=u_0\in M$.

In the following, we shall prove that the well-known relationship between geodesics and sprays holds within the general context. A vector field $S$ on $TM$, smooth on $TM-0$, is said to be a spray on $M$ if $T\pi\circ S=$id$|_{TM}$. Let $C$ be the canonical vector field on $TM$ defined locally by $C(x,a)=(x,a,0,a)$, $x\in U$, $a\in\bfm$.

{\bf Lemma 2.2.} {\it A vector field $S$ on $TM$, smooth on $TM-0$, is a spray on $M$ if and only if $J\circ S=C$, where $J$ is the natural almost tangent structure on $TM$.}

\noindent{\it Proof.} A vector field $S$ on $TM$ can be written locally as follows $$S(x,a)=(x,a,S_1(x,a),\ S_\varphi(x,a)),\ x\in U, a\in \bfm.$$ The condition $T\pi\circ S=$id$|_{TM}$ implies $S_1(x,a)=a$, therefore $S(x,a)=$ $(x,a,a,$ $S_\varphi(x,a))$, where $S_\varphi$ is smooth on $U\subset\bfm-0$. It follows easily $J\circ S=C$, because $J(x,a,b,c)=(x,a,0,b)$, $x\in U$, $a,b,c\in M$.

Conversely, given $S$ as above, the condition $J\circ S=C$ implies $S_1(x,a)=a$, hence $T\pi\circ S=$id$|_{TM}$.

{\bf Lemma 2.3.} {\it Let $K$ be the connection map of a nonlinear connection on $M$. There exists a unique spray on $M$ such that $K\circ S=0$, called geodesic spray.}

\noindent{\it Proof.} Locally, every spray $S$ can be written as follows: $$S(x,a)=(x,a,a,S_\varphi(x,a)).$$ We obtain the geodesics spray if we take $S(x,a)=—\oo(x,a)$ $a,x\in U$, $a\in\bfm$.

The geodesics of a nonlinear connection are the solutions of the following differential equation.
\begin{equation}
\pp^2 c_\varphi(t)+\oo(c_\varphi(t),\pp c_\varphi(t))\pp c_\varphi(t)=0.
\end{equation}

{\bf Theorem 2.4.} {\it A curve $c:(a,b)\to M$ is a geodesic of a nonlinear connection $N$ if and only if there exists an integral curve $\wt{c}:(a,b)\to TM$ of the geodesic spray $S$ of $N$, such that $\pi\circ \wt{c}=c$.}

\noindent{\it Proof.} The curve $\wt{c}$ on $TM$ is an integral curve of $S$ if $\wt{\dot{c}}=S$, therefore in a local chart $(U,\varphi)$ we have $\pp c_\varphi(t)=S_\varphi(\pi\circ\wt{c}_\varphi(t),\wt{c}_\varphi(t))$. Differentiating $\pi\circ\wt{c}=c$, we obtain $\wt{c}_\varphi(t)=\pp c_\varphi(t)$, therefore $\pp^2 c_\varphi(t)=S_\varphi(c_\varphi(t), \pp c_\varphi(t))=-\oo(c_\varphi(t),\pp c_\varphi(t))\pp c_\varphi(t)$, or $\pp^2 c_\varphi(t)+\oo(c_\varphi(t),\pp c_\varphi(t))\pp c_\varphi(t)=0$, i.e. $c$ is a geodesic of $N$. Conversely, if $c$ is a geodesic of $N$, then the curve $\wt{c}$ on $TM$ given by $\wt{c}(t)=\dot{c}(t)$ is an integral curve of the geodesic spray of $N$ and $\pi\circ\wt{c}=c$.

{\it Remarks.} As a corollary of the Theorem 2.4 one obtains again the local existence and uniqueness of a geodesic with given initial conditions. The geodesic spray of a g.a.c. is locally given by $^a S_\varphi(x,a)=(x,a,a,-\G(x)(a,a)-A(x)a)$. Using (2.18) one obtains again the equation (2.17) for the geodesics of a g.a.c. Let us suppose that $M$ has finite dimension. Then a curve $c$ can be written as follows: $x^i=x^i(t)$, $t\in(a,b)$, $i=1,2,...,m=$dim$M$ and the equation (2.17) becomes
\begin{equation}
\dfrac{d^2 x^i}{dt^2}+\G^i_{jk}\dfrac{dx^j}{dt}\dfrac{dx^k}{dt}+A^i_j\dfrac{dx^j}{dt}=0,
\end{equation}
where $\G^i_{jk}$ are the Christoffel symbols of the linear connection associated to the g.a.c. and $A^i_j$ are the components of a tensor of type $(1.1)$ on $M$. The solutions of the equation (2.19), called holomorphically planar curves have been used to give some geometrical meanings in the geometry of complex manifolds [6].

When $M$ is the space-time manifold of the general theory of relativity, the solutions of (2.19) are the trajectories of a charged particle moving in an electromagnetic field [1].

\section{Flat generalized affine connections}
\setcounter{equation}{0}
Let $V:TE\to VE$ be the splitting map which defines a g.a.c. on $p:E\to M$. The map $V$ can be viewed as a 1-form $VE$ valued. On the other hand, the linear connection $\G$ defined by $\wt{\G}$, induces a linear connection $\G_v$ on $VE\to E$.

{\bf Definition 3.1.} The exterior differential $dV$ of the 1-form $VE-$valued $V$, accounted using the linear connection $\G_v$ on $VE\to E$ will be called the curvature form of $\wt{\G}$.

The local component of $V$, denoted also by $V:U\times\bfe\to L(\bfm,\bfe,\bfe)$ is $V(x,a)(\lbb,b)=b+\G(x)(a,\lbb)+A(x)\lbb$, $x\in U$, $\lbb\in M$, $a,b\in\bfe$. After a calculus rather long but not difficult, one obtains the following local expression for the curvature of form $\wt{\G}$:
\begin{equation}
\begin{array}{c}
dV(z,a)((\lbb,b),(\mu,c))=R(x)(\lbb,\mu)+\pp A(x)(\lbb,\mu)-\\ \\
-\pp A(x)(\mu,\lbb)+\G(x)(A(x)\mu,\lbb)-\G(x)(A(x)\lbb\mu),
\end{array}
\end{equation}
$x\in U$, $\lbb,\mu\in\bfm$, $a,b,c\in\bfe$, where $R(x)$ is the local component of the curvature tensor of $\G$.

From (3.1) it follows that $dV$ vanishes when it is applied to a vertical vector field ($\lbb=0$ or $\mu=0$), therefore $dV$ is an horizontal 2-form i.e.
\begin{equation}
dV(\bfa,\bfb)=dV(h\bfa, h\bfb)
\end{equation}
holds for every vector fields $\bfa,\bfb$ on $E$. Using $dV(\bfa,\bfb)=^V\nabla_\bfa V\bfb-^V\nabla_\bfb V\bfa-V[\bfa,\bfb]$, where $^V\nabla$ is the covariant differentiation associated to $\G_v$ and (3.2), one obtains.
\begin{equation}
dV(\bfa,\bfb)=—V[h\bfa,h\bfb]\ ({\rm{the\ structure\ equation\ of\ }}\overline{\G}).
\end{equation}

{\bf Definition 3.2.} A g.a.c. will be called flat if its horizontal distribution is involutive i.e. the bracket of two horizontal vector fields is again a horizontal vector field.

From the structure equation (3.3) it follows

{\bf Theorem 3.1.} {\it The g.a.c. $\wt{\G}$ is flat if and only if $dV=0$.}

Taking $dV=0$ and $a=0$ in (5.1) one obtains
\begin{equation}
\pp A(x)(\lbb,\mu)-\pp A(x)(\mu,\lbb)+\G(x)(A(x)\mu,\lbb)-\G(x)(A(x)\lbb,\mu)=0.
\end{equation}
Taking again $dV=0$ in (3.1) and using (3.4) one obtains
\begin{equation}
R(x)(\lbb,\mu)a=0.
\end{equation}
Conversely, if (3.4) and (3.5) hold, then $dV=0$, therefore we have proved the following

{\bf Theorem 3.2.} {\it The g.a.c. $\wt{\G}=(\G,A)$ is flat if and only if the curvature tensor $R$ of $\G$ vanishes identically and (3.4) holds.}

When $p=\pi:TM\to M$, the conditions (3.4) is equivalent to the vanishing of the following tensor of type $(1, 2)$ on $M$
\begin{equation}
^aT(X,Y)=\nabla_XA(Y)-\nabla_YA(X)-A[X,Y],
\end{equation}
where $X,Y$ are vector fields on $M$, which will be called the torsion tensor of g.a.c. $\wt{\G}$. The Theorem 3.2 has the following

{\bf Corollary 3.2.} {\it A g.a.c. $\wt{\G}=(\G,A)$ on the manifold $M$ is flat if and only if $R=0$ and $^aT=0$.}

When $A=I$ (the Kronecker tensor), $^aT$ becomes the well-known torsion tensor of an a affine connection, therefore we have

{\bf Corollary 3.3.} {\it An affine connection on a Banach manifold $M$ is flat if and only if its curvature tensor and torsion tensor are both identically zero.}

\bigskip

{\bf REFERENCES}

\begin{enumerate}
  \item Cattaneo Gasparini I, {\it Sulle connessioni infinitesimale nello spazio fibrato dei riferimenti affini di una $V_n$}, Rend. Mat. Roma, 17 (1958), p. 326-404.
  \item Cruceanu V., {\it Sur les connexions affines g\'{e}n\'{e}rales}, An. \c{S}t. Univ. Iasi, sect. I a Matematica, XVII (1971), p. 421 - 435.
  \item Grifone J., {\it Structure presque-tangente et connexions} I, II. Ann. Inst. Fourier, Grenoble, t. 22, 1 p. 291 - 334 and t. 22, 3 p. 291 - 338.
  \item Kobayashi S. and Nomizu K., {\it Foundations of differential geometry} I, John Wiley and Sons, New-York-London, 1963.
  \item Miron R., {\it On transformations group of Finsler connections}, (to appear).
  \item Tachibana S. and Ishihara S., {\it On infinitesimal holomorphically projective transformations in K\"{a}hlerian manifolds}, Tohoku Math. J. 12 (I960), p. 77-101.
  \item Vilms J., {\it Connections on tangent bundles}, J. Differential Geometry 1 (1967), p. 235-243.
  \item Vilms J., {\it Curvature of nonlinear connections}, Proceedings AMS, 19 (1968), p. 1125-1129.
\end{enumerate}

\bigskip

\noindent{\footnotesize{\it Received 19.X.1979}}\hfill {\it University of Iasi}\\ \hspace*{9cm}{\footnotesize{\it Faculty of Mathematics}}


\newpage




\def\bfe{{\mathbf{E}}}

\runningauthor={M. ANASTASIEI}
\runningtitle={SOME EXISTENCE THEOREMS IN FINSLER GEOMETRY}
\noindent
\baselineskip 8pt
\noindent{\footnotesize{ANALELE {\c S}TIIN{\c T}IFICE ALE UNIVERSIT{\u A}{\c T}II
``AL.~I.~CUZA" IA{\c S}I\hfill\break
Tomul XXIX, s.I a, Matematic\u a, 1983, 79--84}}

\vskip 2cm
\baselineskip 11.5pt plus .15pt
\centerline{\bf\Large SOME EXISTENCE THEOREMS}
\vskip .2cm
\centerline{\bf\Large IN FINSLER GEOMETRY}
\vskip .5cm
\centerline{\bf\footnotesize{BY}}
\vskip .5cm
\centerline{\bf M. ANASTASIEI}
\vskip 1cm

Let $M$ be a differentiable manifolds of class $C^{k}$ ($k\geq 3$) and let $p:TM\to M$ be the tangent bundle to it. A positive real valued function $L:TM\to R^+$, with properties:

1) $L$ is differentiable of class $C^{k-1}$ on $TM\setminus O$ and continuous on the image of the null-section of $p$,

2) Its local representation in every chart	$(p^{-1}(U_i),x^i,y^i)$ on $TM$ induced by the chart $(U,x^i)$ on $M$, denoted by $L(x^i,y^i)$ or for brevity by $L(x,y)$ is (1) $p-$homogeneous i.e. $L(x,sy)=sL(x,y)$ for every $s>0$,

3) The matrix $(g_{ij}(x,y))=\(\dfrac{\pp^2(L^2(x,y)/2)}{\pp y^i\pp y^j}\)$ is invertible and the qua\-dra\-tic form associated to it is positive definite,

\noindent is called {\it fundamental Finsler function} and the pair $(M,L)$ is called a Finsler space.

The matrix $(g_{ij}(x,y))$ changes, when the local chart changes, as the components of a tensor of type $(0,2)$ on $M$ but it depends on direction (given by $y^i$) so $g_{ij}(x,y)$ defines a Finsler tensor field of type $(0,2)$ (see [1] for a general definition of Finsler geometric objects). New fields of Finsler geometric objects (tensors, connections) can be derived from $L$. Obviously, their existence is assured by the existence of $L$.

As it was pointed out in [3, p. 81], the conditions imposed on $L$ are too restrictive. It was an idea of R. Miron to eliminate the function $L$ and to define a Finsler space as a pair $(M,g)$, where $g$ is a symmetric Finsler tensor field of type $(0,2)$, nondegenerate, positive definite or not. He called such a $g$ a metrical Finsler structure on $M$. The fields of Finsler geometric objects can be defined independent on $L$ or $g$. So the problem of their existence, in particular that of the existence of $L$ and $g$ is quite natural. The aim of this paper is to prove the global existence of Finsler tensor fields, of linear Finsler connections and of metrical Finsler structures in the hypothesis $M$ paracompact, modeled by a separable Hilbert space. In a second section we make more explicit a proof due to S. Kashiwabara [2] of the global existence of a fundamental function $L$, in the hypothesis $M$ finite dimensional and paracompact. Some comments regarding the existence of $L$ when $M$ is infinite dimensional are made.

\setcounter{section}{0}
\section{The global existence of Finsler geometric objects} \setcounter{equation}{0}

Let us state again the hypothesis on $M$ in this section: differentiable of class $C^k$ $(k\geq3)$, modeled by a separable Hilbert space $H$ and paracompact i.e it is separate Hausdorff and every open covering of it admits a locally finite refinement. We recall that a $C^k-$partition of unity on a manifold $X$ is an indexed family of $C^k$ real valued functions $\{f_j\}_{j\in J}$ on $X$ such that

1) $f_j\geq 0$,

2) $\{{\rm{supp}}(f_j)\}_{j\in J}$ is locally finite and

3) $\dd\sum_j f_j(x)=1$, $x\in X$.

\noindent The partition of unity $\{f_j\}_{j\in J}$ is said to be subordinate to the covering $\{U_i\}_{i\in I}$ of $X$ if $\{{\rm{supp}}(f_j)\}_{j\in J}$ refines $\{U_i\}_{i\in I}$.

As it was proved in [6, p.57-60], every separable Hilbert space admits $C^k-$partition of unity and a necessary and sufficient condition that a $C^k-$ma\-ni\-fold $X$ admit $C^k-$partition of unity is that $X$ be paracompact and that each of its tangent spaces admit $C^k-$partition of unity. Consequently, our manifold $M$ admits $C^k-$partition of unity.

Let $A=\{(U_i,\varphi_i)\}_{i\in I}$ be an atlas on $M$. Suppose that $A$ is maximal i.e. it contains all charts compatible with it. Then $\{U_i\}$ is a basis for the topology of $M$. Let $x$ be a point of $M$ and let $(U_i,\varphi_i)$, $(U_j,\varphi_j)$ be two local charts around $x$. The triads $(U_i,\varphi_i,u)$ and $(U_j,\varphi_j,v)$, where $u,v\in H$, are called equivalent if $D_{\varphi_i}(x)(\varphi_j\circ\varphi_i^{-1})(u)=v$, where $D$ means the Fr\'{e}chet differentiation. This is indeed an equivalence relation on the set of such a triad $s$ and the class of equivalence $[(U_i,\varphi_i,u)]$ is called vector tangent to $M$ in $x$. Thus every chart $(U_i,\varphi_i)$ defines a map $\theta_{i,x}:T_xM\to H$, $\theta_{i,x}([U_i,\varphi_i,u)]=u$. We set $TM=\dd\cup_{x\in M}T_xM$ and $p:TM\to M$ projects $T_xM$ on $x$. The topology and the differentiable structure of $TM$ are induced by those of $M$. As a basis for the topology of $TM$ is taken $(p^{-1}(U_i))_{i\in I}$ and $p$ becomes a continuous map. One verifies easily that $(p^{-1}(U_i),h_i)$, where $h_i:p^{-1}(U_i)\to H\times H$, $h_i(z)=(\varphi_i(p(z))$, $\theta_{i,p(z)}(z)$ is a $C^{k-1}-$atlas (not maximal) on $TM$. Usually the differentiable structure defined by this atlas is considered. The map $p$ becomes a $C^{k-1}-$submersion.

\ms

{\bf Theorem 1.1.} {\it The manifold $TM$ is paracompact.}

\ms

\noindent{\it Proof.} Let $z_1,z_2\in TM$ and let us denote $x_1=p(z_1)$, $z_2=p(z_2)$. There exist open sets $D_1$ and $D_2$ such that $x_1\in D_1$, $x\in D_2$ and $D_1\cap D_2=\emptyset$. Putting $D_1=\dd\cup_{j\in I_1}U_j$ and $D_2=\dd\cup_{j\in I_2}U_j$, it follows that there exist $j_1\in I_1$, $j_2\in I_2$ such that $x_1\in U_{j_1}$, $x_2\in U_{j_2}$ and $U_{j_1}\cap U_{j_2}=\emptyset$. Then $p^{-1}(U_{j_1})\cap p^{-1}(U_{j_2})=\emptyset$ and $z_1\in p^{-1}(U_{j_1})$, $z_2\in p^{-1}(U_{j_2})$, therefore $TM$ is separate Hausdorff.

Let $\{D_j\}_{j\in J}$ be an open covering of $TM$. We may write $D_j=\dd\cup_{i\in I}p^{-1}(U_i)$. The open covering $\{U_i\}_{i\in I}$ admits an open locally finite refinement $\{V_k\}_{k\in K}$ i.e. for every $k\in K$ there exists $i(k)\in I$ such that $V_k\subset U_{i(k)}$. It follows $p^{-1}(V_k)\subset p^{-1}(U_{i(k)})$ and obviously there exists $D_{j(k)}\supset p^{-1}(U_{i(k)})$, therefore $(p^{-1}(V_k))_{k\in K}$ is an open refinement of the covering $\{D_j\}$. We prove that it is locally finite. If $z\in TM$ and $x=p(z)$, there exists an open neighborhood $U$ of $x$ which intersects only a finite number of $V$'s say $V_1,...,V_n$. It follows by reductio ad absurdum that $p^{-1}(U)$ intersects only $p^{-1}(V_1),...,p^{-1}(V_n)$. The proof is complete.

\ms

{\bf Corollary 1.1.} {\it Let $M$ be a paracompact manifold modeled by a separable Hilbert space $H$. Then the manifold $TM$ admits $C^{k-1}-$partition of unity.}

\ms

\noindent{\it Proof.} The space $H\times H$ being the product of two separable Hilbert spaces is itself a separable Hilbert space. The manifold $TM$ being paracompact, the proof follows via the above-mentioned theorems.

\ms

{\bf Corollary 1.2.} {\it Let $M$ be a finite dimensional manifold of class $C^k$ paracompact. Then $TM$ is a paracompact manifold of class $C^{k-1}$ and it admits a $C^{k-1}-$partition of unity.}

\ms

\noindent{\it Proof.} Obvious.

\ms

A partition of unity for $TM$ can be obtained from a partition of unity for $M$ as follows:

\ms

{\bf Theorem 1.2.} {\it Let	$\{f_j\}_{j\in J}$ be a $C^k-$partition of unity on $M$ subordinated to the covering $\{U_i\}_{i\in I}$. Then $\{f^V_j=f_i\circ p\}_{jeJ}$ is a $C^{k-1}-$partition of unity on $TM$ which is subordinated to the covering $\{p^{-1}(U_i)\}_{i\in I}$.}

\ms

\noindent{\it Proof.} Obviously, $f^V_j\geq 0$ for every $j\in J$. Then carr$f^V_j=\{z\in TM| f_j(p(z))\neq\emptyset\}\subset p^{-1}({\rm{supp}}f_j)$, therefore supp$f^V_j\subset p^{-1}({\rm{supp}}f_j)\subset p^{-1}(U_{i(j)})$ because supp$f_j$ is closed. Since $\{{\rm{supp}}f\}_{j\in J}$ is locally finite, so is $\{{\rm{supp}}f^v_j\}_{j\in J}$. The equalities $\sum f^V_i(z)=\sum f_i(p(z))=1$ end the proof.

The fields of Finsler geometric objects can be obtained as cross-sections of a convenient fibre bundle over $TM$. We recall briefly the construction of that bundle $FO^k\to TM$ (see [5]). Let us denote by $L_k(H)$ the set of $k-$jets of source $0\in H$ of the local diffeomorphism of $H$ which preserves $0\in H$. The composition of $k-$jets gives a group structure on $L_k(H)$. The set $P^k(M)$ of all $k-$jets of source $0\in H$ of the local diffeomorphisms of $H$ to $M$ can be structured as a principal fibre bundle over $M$ with structural group $L_k(H)$. The pull-back by $p$ of this bundle will be denoted by $F^k(M)\to TM$ (this is the Finsler bundle of order $k$).

A pair $(F,m)$, where $F$ is a manifold and $m$ a differentiable action of $L_k(H)$ on $F$, is called a manifold of geometric objects. Usually $F$ is taken a linear space or an open subset of a linear space. The fibre bundle associated to $F^k(M)$ of type fiber $(F,m)$ is denoted by $FO^k\to TM$ and is called the bundle of Finsler geometric objects. Its cross-sections are called fields of Finsler geometric objects on $M$. If $F$ is a linear space and $m$ is a linear action on $F$, the cross-sections of the bundle of Finsler geometric objects are called fields of linear Finsler geometric objects.

Now to state a result proved in [6, p. 62], some definitions are necessary. A subset $C\subseteq E$, where $E$ is the total space of a $C^k-$bundle $q:E\to M$, is said to be convex if for each $x\in M$, $C_x=C\cap E_x$ is a nonvoid and convex set. (Here $E_x$ denotes the fiber in $x$.) One says $C$ admits local $C^k-$sections if given $c_0\in C$ with $q(c_0)=x_0$, there is an open neighborhood $U$ of $x_0$ and a $C^k-$section $s$ over $U$ with $s(x_0)=c_0$ and $s(U)\subseteq C$.

\ms

{\bf Theorem 1.3.} [6] {\it Let $q:E\to M$ be a $C^k-$bundle and let $C$ be a convex subset of $E$ which admits local $C^k-$sections. There exists a $C^k-$section $S$ of $q$ over $M$ such that $S(x)\in C$ for every $x\in M$.}

\ms

The bundle $FO^k\to TM$ associated to $P^k(M)\to TM$ with $F$ a linear space and $m$ a linear action admits local $C^k-$sections, because it is a locally trivial vector bundle (every bundle chart of it defines a local $C^k-$section). Using the Theorem 1.3 one obtains

\ms

{\bf Theorem 1.4.} {\it Let $M$ be a paracompact manifold modeled by a separable Hilbert space. There exist global fields of linear Finsler geometric objects on $M$ i.e. cross-sections over $TM$ of vector bundle $FO^k\to TM$.}

\ms

The Finsler tensor fields of type $(0,r)$ or $(1,r)$, $r\geq 1$ can be considered as linear Finsler fields of geometric objects by adjusting, in an obvious manner, the corresponding  definition from finite dimensional case (see [1]). So, we have

\ms

{\bf Corollary 1.5.} {\it Under the hypothesis of the above theorem, there exist globally on $M$, Finsler vector fields and Finsler tensor fields of type $(0,r)$ and $(1,r)$.}

\ms

{\bf Corollary 1.6.} {\it Let $M$ be a paracompact finite dimensional manifold. Then there exists global Finsler tensor fields of any type on $M$.}

\ms

Let us take $F=L_2(H,H)$ i.e. the linear space of linear maps $H\times H\to H$ and let the action $m$ be denoted by $m_c$ and defined by $m_c:L_2(H)\times L_2(H,H)\to L_2(H,H)$, $m_c(h,K)=AK(A^{-1},A^{-1})-A_1(A^{-1},A^{-1})$ if $h=(A,A_1)\in L_2(H)$. With this choice of $F$ and $m$, the cross-sections of $p_c:FO^2\to TM$ are called linear Finsler connections on $M$. The fiber $p_c^{-1}(z)$, $z\in TM$, is not a linear space although its elements can be added and multiplied by reals, because $m_c$ is not linear. However it is a convex set because if $a+b=1$, $a,b\in R$, $m_c(h,aK_1+bK_2)=am_c(h,K_1)+bm_c(h,K_2)$ holds good. The bundle $p_c$ is locally trivial hence it admits local $C^k-$sections. The Theorem 1.3 applies and leads to

\ms

{\bf Theorem 1.7.} {\it On every paracompact manifold modeled by a separable Hilbert space there exist global Finsler linear connections.}

\ms

{\bf Corollary 1.8.} {\it On every paracompact finite dimensional manifold there exist global Finsler linear connections.}

\ms

Let us now take $F=L^s_2(H,R)$, the linear space of real valued bilinear and symmetric maps on $H$ and $m:GL(H)\times L^s_2(H,R)\to L^s_2(H,R)$ given by $m(A,g)=g(A^{-1},A^{-1})$, where $A$ belongs to the general linear group $GL(H)$ of $H$. The cross-sections of the bundle $p_1:FO^1\to TM$ obtained with such a choice of $F$ and $m$ are symmetric Finsler tensor fields of type $(0,2)$. Every $g\in L^s_2(H,R)$ defines a linear operator $\wt{g}:H\to H^*$. If $g$ is invertible, $g$ is called nondegenerate. The set ${\cal R}(H,R)$ of all nondegenerate symmetric bilinear maps on $H$ is an open and convex subset of $L^s_2(H,R)$. Obviously, $m$ leaves invariant the subset ${\cal R}(H,R)$. By applying the general construction sketched above taking ${\cal R}(H,R)$ as $F$, one obtains a fibre bundle over $TM$ whose total space ${\cal R}FO^1$ is a subset of $FO^1$ which is clearly open and convex. Being open it admits local $C^k-$sections, therefore by applying the Theorem 1.3 one obtains

\ms

{\bf Theorem 1.9.} {\it On every paracompact manifold modeled by a separable Hilbert space, there exist global metrical Finsler structures i.e. cross-sections $g:TM\to FO^1$ such that $g(TM)\subset{\cal R}FO^1$.}

\ms

{\bf Corollary 1.10.} {\it On every paracompact finite dimensional manifold there exist global metrical Finsler structures.}

\ms

When $F$ is a linear space and $m$ a linear action, the bundle $FO^k\to TM$ can be identified (is isomorphic) to a vector bundle obtained from the vertical subbundle $V$ of $TTM\to TM$ by algebraic operations. So, the Finsler vector fields appear as sections of $V\to TM$, the Finsler tensor fields of type $(0,r)$ appear as sections of $L(V,...,V;R)\to TM$, the Finsler tensor fields of type $(1,r)$ are sections of $L(V,...,V;V)\to TM$ and so on.

A Finsler almost product structure on $M$ is a section $P:TM\to L(V,V)$ which satisfies $P(z)\circ P(z)=I$, where $I$ is the identity map on fiber $V_z$. The global existence of such a structure on $M$ is a consequence of the following fact: if $V'$ is a subbundle of $V$, there exists a subbundle $V''$ of $V$ such that $V'\oplus V''=V$. Its proof is standard, using a partition of unity on $TM$. So, if $s_z=s'_z+s''_z$ where $s'_z\in V'_z$ and $s''_z\in V''_z$ we may define $P(s'_z)=s'_z$, $P(s''_z)=—s''_z$ and it follows that $P\circ P=I$.

\section{The global existence of a fundamental Finsler function}\setcounter{equation}{0} Let us suppose that $M$ is a separate Hausdorff, finite dimensional manifold, satisfying the second axiom of countability. It follows that $M$ is paracompact, therefore it admits $C^k-$partition of unity. We shall prove the existence of a real valued function $L$ on $TM$ verifying the conditions to be a fundamental Finsler function. Firstly, we prove the following

\ms

{\bf Lemma 2.1.} {\it Let $n$ be the dimension of $M$. There exists a continuous function $f:R^n\to R^+$ which is

a) 1 $(p)-$homogeneous,

b) differentiable at least of class $C^3$ on the complement of the origin and the quadratic form $\dfrac{\pp^2(f^2(y)/2)}{\pp y^i\pp y^j}z^i z^j$ is positive definite for all values of $z^i\neq 0$, where $y=(y^1,...,y^n)$, $z=(z_1,...,z^n)$ belong to $R^n$.
Furthermore, $f$ is a norm on $R^n$.}

\ms

\noindent {\it Proof.} Let $h:R^n\to R$ be a continuous function $1(p)-$homogeneous, differentiable at least of class $C^3$ on the complement of the origin and $h(0) =0$. Such a function always exists. For instance we may take $h(y)=\(\dd\sum^{n}_{i=1}(y^i)^p\)^{1/p}$ with $p\geq 1$, $p\neq 2$. Let us put $f(y)=\(\dd\sum^n_{i=1}(y^i)+\vare h^2(y)\)^{1/2}$, where $\vare$ is a positive real number. Obviously, $f$ is $1(p)-$homogeneous and differentiable at least of class $C^3$ on the complement of the origin. The matrix $A=\(\dfrac{\pp^2(f^2/2)}{\pp y^i\pp y^j}\)$ is given as follows: $A=I+B$, where $I$ is the unity matrix and $B=\(\dfrac{\pp h}{\pp y^i}\cdot\dfrac{\pp h}{\pp y^j}+h(y)\dfrac{\pp^2h}{\pp y^i\pp y^j}\)$. Choosing $\vare<1/\|B\|$, where $\|B\|$ means a norm on the space of matrices, $A$ becomes an invertible matrix and furthermore, the quadratic form associated to $A$ becomes positive definite. It is obvious that $f(y)\geq 0$, the equality sign occurring only if $y=0$. The condition $f(sy)=|s|f(y)$ is clearly satisfied. A proof of the inequality $f(x+y)\leq f(x)+f(y)$, $x,y\in R^n$ can be performed using a method of H. Rund [7, p. 18--20].

\ms

{\bf Theorem 2.1.} {\it Let $M$ be a finite dimensional paracompact manifold. There exists a fundamental Finsler function on $TM$.}

\ms

\noindent{\it Proof.} Let $(a_i)_{i\in I}$ be a $C^k-$partition of unity $(k\geq 1)$ on $M$ subordinate to a covering $\{U_i\}_{i\in I}$ of $M$, where $U_i$, is the domain of a bundle chart $\varphi_i:p^{-1}$ $(U_i)\to U_i\times R^n$. Define $L_i:U_i\times R^n\to R$ by $L_i(p,u)=f(u)$, where $f$ is given by the Lemma 2.1. The function $L$ defined by $L(v_p)=\dd\sum_{i}a_i(p)L_i(\varphi_i(v_p))$, $v_p\in TM$, satisfies all requirements to be a fundamental Finsler function.

In his lectures given at Brandeis University in 1965, R. S. Palais had considered what he called Finsler structures on a Banach bundle, in particular on a Banach manifold. Following his definition, a Finsler structure on M is a function $L:TM\to R^+$ such that for every $p_0\in M$ there exists a bundle chart $\varphi:p^{-1}(U)\to U\times H$ such that $L\circ \varphi^{-1}$ verifies

1)	$(L\circ \varphi^{-1})(p_0):H\to R^+$ is an admissible norm on $H$,

2)	There exists a neighborhood $U_0\subset U$ of $p_0$ such that $(L\circ\varphi^{-1})(p)$ be an equivalent norm to $(L\circ \varphi^{-1})(p_0)$ for every $p\in U_0$.

\noindent This notion is less restrictive then the usual notion of Finsler structure in finite dimensions, where the second condition becomes trivial. The existence of a function $L$ satisfying 1) and 2) was proved by R. S. Palais by means of a partition of unity on $M$. The basic tool in his proof was the so-called flat Finsler structure given by the map $N:M\times H\to R^+$, $N(p,u)=\|u\|$, for every $p\in M$ and $u\in H$, where $\|\cdot\|$ is the norm induced by the inner product of $H$. Such a Finsler structure satisfies a third condition

3) $D^2_u(L\circ \varphi^{-1})^2/2$ is for every $u\in H$ an isomorphism of $H$ to $H^*$, but it is essential a Riemannian one. Here $D^2$ means the Fr\'{e}chet differentiation of the second order.

When the norm on $H$ is not Hilbertian i.e. $H$ is a separable Banach space which admits a partition of unity, the procedure of R.~S. Palais leads to a proper Finsler structure, but if the condition 3) is imposed, it becomes a Riemannian one. This happens since the condition 3) implies that the norm of H is equivalent to a norm induced by an inner product of $H$. In the following we give a proof of this assertion. Let $p_0$ be a fixed point of $M$ and let us put $g=(L\circ \varphi^{-1})^2/2=N^2/2=\|\cdot\|^2/2$ and $T=D^2_u g$. The map $T$ can be viewed as a continuous and symmetric bilinear form on $H$. Differentiating $g$, one obtains $D^2_u g(v,v)=(D_uN(v))^2+\|u\|^2D^2_u N(v,v)\geq 0$ for every $v\in H$ since $D^2_u\|\cdot\|(v,v)\geq 0$ (see [7]). Therefore $T$ is also a positive bilinear form on $H$, hence the inequality of Cauchy--Schwarz $||T(u,v)||\leq p(u)p(v)$, where $p(u)=T(u,u)^{1/2}$, holds. If $p(v)=0$ it follows that $T(u,v)=0$ for every $u\in H$ or $(T(v))(u)=0$ for every $u\in H$ and $v=0$ because $T$ is an isomorphism. So, $T$ is an inner product on $H$. The inequality	 $p(v)\leq\|T\|^{1/2}\|v\|$ is obvious. Let us take $S=\{v\in H|p(v)\leq 1\}$. Then the inequality $|T(u)(v)|\leq p(v)$ holds for every $u\in S$. It follows there exists $c>0$ such that $\|T(u)\|\leq c$ for every $u\in S$. We have $\|u\|=\|(T^{-1}\circ T)(u)\|\leq\|T^{-1}\|\ \|T(u)\|\leq\|T^{-1}\|c$ for every $u\in S$. So, $\|u\|\leq c\|T^{-1}\|p(u)$ for every $u\in H$. Therefore the norms $\|\cdot\|$ and $p$ are equivalent.

The above considerations show that in the framework of Banach manifolds the definition of Finsler structures given by R.~S. Palais is the most convenient.

\bigskip

{\bf REFERENCES}

\begin{enumerate}
  \item Anastasiei, M., {\it Finsler geometric objects and their Lie derivative}, The Proceedings of the National Seminar on Finsler Spaces (University of Brasov, february 1980), Timisoara, 1981, 11-25.
  \item Kashiwabara, S., {\it On Euclidean connections in a Finsler manifolds}, Tohoku Math. J., 10 (1958), 69-80.
  \item Matsumoto, M., {\it Foundations of Finsler geometry and special Finsler spaces}, to appear from VEB Deutsch. Verlag. Wiss. Berlin.
  \item Miron, R. and Hashiguchi, M., {\it Metrical Finsler connections}, Rep. Fac. Sci. Kagoshima Univ. (Math. Phys. Chem.) 12 (1979), 21-35.
  \item Miron, R. and Anastasiei, M., {\it On the notion of Finsler geometric objects} (to appear in Mem. Acad. R.S.R. 4 (1981)).
  \item Palais, R. S., {\it Lectures on the differential topology of infinite dimensional manifolds}, Mimeo Notes at Brandeis Univ. by S. Greenfield, 1964--1965.
  \item Sundaresan, K., {\it Smooth Banach spaces}, Math. Annalen, 173 (1967), 191--199.
  \item Rund, H., {\it The differential geometry of Finsler spaces}, Springer-Verlag, Berlin, 1959.
\end{enumerate}

\bigskip

\noindent{\footnotesize{\it Received 22.VII.1981}\hfill {\it Faculty of Mathematics}}\\ \hspace*{9cm}{\footnotesize{\it University of Iasi}}


\newpage



\def\ovcirc{\overset{\circ}}



\runningauthor={M. ANASTASIEI}
\runningtitle={VECTOR BUNDLES. EINSTEIN EQUATIONS}
\noindent
\baselineskip 8pt
\noindent{\footnotesize{ANALELE {\c S}TIIN{\c T}IFICE ALE UNIVERSIT{\u A}{\c T}II
``AL.~I.~CUZA" IA{\c S}I\hfill\break
Tomul XXXII, s.I a, Matematic\u a, 1986, f. 3, 17-24}}

\vskip 2cm
\baselineskip 11.5pt plus .15pt
\centerline{\bf\Large VECTOR BUNDLES. EINSTEIN EQUATIONS}
\vskip .5cm
\centerline{\bf\footnotesize{BY}}
\vskip .5cm
\centerline{\bf M. ANASTASIEI}
\vskip 1cm

In the last years a Finslerian theory of relativity was built from various standpoints [2], [4], [7]. Recently, R. Miron has completed a more generalized version of this theory, which was called a Lagrangian theory of relativity in [6]. Some physical aspects of this theory were considered by S. Ikeda in [3]. His considerations show that the geometry of the total space of a vector bundle is useful from the view point of Physics.

In this paper the Einstein equations and the conservation law on the total space of a vector bundle are written. If the vector bundle is just the tangent bundle to the base manifold we recover the Einstein equations established by R. Miron in [6] as well as a new kind of Einstein's equations whose physical meaning remains to be found. If the vector bundle has 1-dimensional fibres, we obtain a geometrical framework for an unitary projective theory.

The author is indebted to Prof. Radu Miron who suggested him the subject of this work.

\setcounter{section}{0}
\section{Vector bundles} Let $\xi=(E,p,M)$, $p:E\to M$, be a vector bundle of paracompact base $M$ and finite dimensional type fibre $\bff$. We set $n=$dim$M$ and $m=$dim$\bff$. Let us denote, by $(x^i,y^a)$ the local coordinates on $p^{-1}(U)\subset E$, where $U\subset M$. In what follows we use $i,j,k,h...=1,2,...,n$ and $a,b,c,...= 1,2,...,m$.

The law of transformation of the local coordinates is the following:
\begin{equation}
x^{i}=x^i(x^1,...,x^n),\ y^{a'}=S^{a'}_a(x^1,...,x^n)y^a.
\end{equation}
If the vector bundle is endowed with a nonlinear connection, then, for every $u\in E$, we have $T_uE=H_uE\oplus V_uE$, where $V_uE$ is the vertical part and $H_uE$ is the horizontal part. A basis of $T_uE$ adapted to this decomposition is $(\de_i,\pp_a)$, where $\de_i=\pp_i-N^a_i(x,y)\pp_a$. Here $(N^a_i(x,y))$ are the local coefficients of the nonlinear connection and $\pp_i$ and $\pp_a$ stand for $\pp/\pp x^i$ and $\pp/\pp y^a$, respectively. The basis dual to it is $(dx^i,\de y^a)$, where $\de y^a=dx^a+N^a_i dx^i$.

\ms

{\bf Definition 1.1.} A linear connection $D$ on the manifold $E$ is said to be a $d-$connection if it preserves by parallel displacement the horizontal distribution $u\to H_uE$ and vertical distribution $u\to V_uE$.

\ms

If we set:
\begin{equation}
\left\{
\begin{array}{l}
D_{\de_k}\de_j=F^i_{jk}(x,y)\de_i,\ D_{\de_k}\pp_b=L^a_{bk}(x,y)\pp_a,\\ \\
D_{\pp_a}\de_j=M^i_{ja}(x,y)\de_i,\ D_{\pp_c}\pp_b=C^a_{bc}(x,y)\pp_a,
\end{array}\right.
\end{equation}
then $F^i_{jk}(x,y)$ and $L^a_{bk}(x,y)$ change like the local coefficients of a connection on $M$, respectively on $\xi$, and $M^i_{ja}(x,y)$, $C^a_{bc}(x,y)$ are tensor fields on $E$. A $d-$connection is completely determined by $F\G=(F^i_{jk}, L^a_{bk},M^i_{ja}, C^a_{bc})$. (See also [5].)

There exist $d-$connections on $E$. For instance, if $F^i_{jk}(x)$ are the local coefficients of a linear connection on $M$ (there exists such a connection because $M$ is paracompact), then $(F^i_{jk}(x),\pp_b N^a_j,0,0)$ is a $d-$connection on $E$.

A pair of linear connections on $M$ and $\xi$ defines a $d-$connection on $E$. Indeed, if $L^a_{bk}(x)$ are the local coefficients of a linear connection on $\xi$, then $N^a_i(x,y)=L^a_{bk}y^b$ are the local coefficients of a nonlinear connection on $\xi$ and $(F^i_{jk}(x), L^a_{bk}(x),0,0)$ is a $d-$connection.

A $d-$connection $F\G$ is called a Berwald connection if $$L^a_{bk}=\pp_b N^a_k(x,y),\ M^i_{ja}(x,y)=0.$$

The $d-$connections showed above are Berwald connections. We shall denote by $_|$ and $|$ the $h-$ and $v-$covariant derivative, respectively, associated to the $d-$connection $D$.

The Ricci identities introduce five torsions:
\begin{equation}
\left\{
\begin{array}{l}
T^i_{jk}=F^i_{jk}-F^i_{jk},\ R^a_{jk}=\de_k N^a_j-\de_j N^a_h,\ \overset{1}{P^a_{jb}}=\pp_b N^a_j-L^a_{bj},\\ \\
\overset{2}{P^i_{jb}}=M^i_{jb},\ S^a_{bc}=C^a_{bc}-C^a_{cb},
\end{array}
\right.
\end{equation}
and six curvatures:
\begin{equation}
\left\{
\begin{array}{l}
R^i_{jkh}=\de_h F^i_{jk}+F^l_{jk}F^l_{ih}-k|h+M^i_{ja}R^a_{kh},\\ \\
R^a_{bkh}=\de_h L^a_{bk}+L^c_{bk}L^a_{ch}-k|h+C^a_{bc}R^c_{kh},\\ \\
\overset{\hspace*{-2mm}1}{P^a_{bkc}}=\pp_c L^a_{bk}-C^a_{bc|k}+C^a_{bd}P^s_{kc},\\ \\
\overset{\hspace*{-2mm}2}{P^i_{jkc}}=\pp_c F^i_{jk}-M^i_{jc|k}+M^i_{jb}\overset{\hspace*{-2mm}1}{P^b_{kc}},\\ \\
M^i_{jbc}=\pp_c M^i_{jb}+M^h_{jb}M^i_{hc}-b|c,\\ \\
S^a_{bcd}=\pp_d C^a_{bc}+C^e_{bc}C^a_{ed}-c|d,
\end{array}
\right.
\end{equation}
for a $d-$connection $F\G$. Here and in the following $-k|h$ means the substraction of the previous terms after having changed the indices one to another one.

\section{Metrical structures on $E$. Metrical $d-$con\-nec\-tions} \setcounter{equation}{0} A metrical structure on $E$ is a tensor field $G$ on $E$ of type $(0,2)$, symmetric and nondegenerate. If such a metrical structure $G$ is given, then there exists a canonical nonlinear connection on $\xi$ defined by the orthogonal distribution to the vertical distribution with respect to $G$. In what follows we shall refer only to this nonlinear connection. It is obvious that, with respect to the adapted frame to this nonlinear connection, $G$ can be written as follows:
\begin{equation}
G=g_{ij}(x,y)dx\otimes dx^j+h_{ab}(x,y)\de y^a\otimes \de y^b.
\end{equation}

\ms

{\bf Definition 2.1.} A $d-$connection on $E$ is said to be metrical if
\begin{equation}
g_{ij|k}=0,\ g_{ij}|_a=0,\ h_{ab|k}=0,\ h_{ab}|_c=0,
\end{equation}
hold.

There exist metrical $d-$connections. Indeed, if $F\overset{\circ}{\G}=(\overset{\circ}{F}^{i}_{jk},\overset{\circ}{L}^a_{bk},\overset{\circ}{M}^i_{ja},\overset{\circ}{C}^a_{bc})$ is any $d-$connection on $E$, then the $d-$connection whose local coefficients are given below is metrical:
\begin{equation}
\left\{
\begin{array}{l}
F^i_{jk}=\overset{\circ}{F}^i_{jk}+\dfrac12 g^{ih}g_{hk\ovcirc{|}j}\\ \\
L^a_{bj}=\overset{\circ}{L}^a_{bj}+\dfrac12 h^{ac}hg_{cb\ovcirc{|}a}\\ \\
M^i_{jb}=\overset{\circ}{M}^i_{jb}+\dfrac12 g^{ih}g_{hi\ovcirc{|}b}\\ \\
C^a_{bc}=\overset{\circ}{C}^a_{bc}+\dfrac12 h^{ad}g_{db\ovcirc{|}c},
\end{array}
\right.
\end{equation}
where $\ovcirc{_|}$ and $\ovcirc{|}$ denotes the $h-$ and $v-$covariant derivative, respectively, associated to $F\ovcirc{\G}$.

The formulas (2.3) can be thought of as a process of metrization of any $d-$connection. This process will be called Kawaguchi metrization.

We say that a $d-$connection is $h-v-$metrical with respect to $G$ given by (2.1), if $g_{ij|k}=0$ and $h_{ab}|_c=0$. We remark that there exist $h-v-$metrical connections which are not metrical. Indeed, it is easy to check that the following Berwald connection
\begin{equation}
\left\{
\begin{array}{l}
\wt{F}^i_{jk}=\dfrac12 g^{ih}(\de_k g_{hj}+\de_j g_{hk}-\de_h g_{jk})\\ \\
\wt{L}^a_{bj}=\pp_b N^a_{j}\\ \\
\wt{M}^i_{jb}=0\\ \\
\wt{C}^a_{bc}=\dfrac12 h^{ad}(\pp_b h_{dc}+\pp_c h_{db}-\pp_d h_{bc})
\end{array}
\right.
\end{equation}
is $h-v-$metrical but it is not metrical.

\ms

{\bf Theorem 2.1.} {\it If two skew-symmetrical tensor fields $T^i_{jk}$ and $S^a_{bc}$ are given, then there exists a unique Berwald connection which is $h-v-$metric and has $h(hh)-$ and $v(vv)-$torsions the tensor fields $T^i_{jk}$ and $S^a_{bc}$, respectively. Its local coefficients are as follows:
\begin{equation}
\left\{
\begin{array}{l}
\widehat{F}^i_{jk}=\wt{F}^i_{jk}+\dfrac12 g^{ih}(g_{hr}T^r_{jk}-g_{jr}T^r_{hk}+g_{kr}T^r_{jh})\\ \\
\widehat{L}^a_{bk}=\pp_b N^a_k\\ \\
\widehat{M}^i_{ja}=0\\ \\
\widehat{C}^a_{bc}=\wt{C}^a_{bc}+\dfrac12 h^{ad}(h_{de}S^e_{bc}-h_{be}S^e_{dc}+h_{ce}S^e_{bd}).
\end{array}
\right.
\end{equation}}

\ms

\noindent{\it Proof.} All Berwald connections have the form $(F^i_{jk}+\tau^i_{jk},\ \pp_b N^a_k,0,C^a_{bc}+\wt{\tau}^a_{bc})$, where $\tau^i_{jk}$ and $\wt{\tau}^a_{bc}$ are arbitrary tensor fields. Imposing that such a connection be $h-v-$metrical and its $h(hh)-$ and $v(vv)-$torsions to be exactly $T^i_{jk}$ and
$S^a_{bc}$, respectively, one obtains that $\tau^i_{jk}$ and $\wt{\tau}^a_{bc}$ are uniquely determined and they have the expressions from (2.5), q.e.d.

\ms

{\bf Theorem 2.2.} {\it There exists a unique metrical d-connection with $h(hh)-$ and $v(vv)-$torsions $T^i_{jk}$ and $S^a_{bc}$ prescribed, obtained by Kawaguchi metrization of a $h-v-$metrical Berwald connection. Its local coefficients are as follows:}
\begin{equation}
\left\{
\begin{array}{l}
F^i_{jk}=\wt{F}^i_{jk}\\ \\
L^a_{bj}=\pp_b N^a_{k}+\dfrac12 h^{ac}[\de_k g_{bc}-(\pp_b N^d_k)k_{dc}-(\pp_c N^d_k)h_{db}]\\ \\
M^i_{jb}=\dfrac12 g^{ik}\pp_b g_{jk}\\ \\
C^a_{bc}=\wt{C}^a_{bc}.
\end{array}
\right.
\end{equation}

\ms

\noindent{\it Proof.} By the Kawaguchi metrization of the unique Berwald $h-v-$metrical connection given by (2.5) one obtains (2.6), q.e.d.

\section{Einstein equations on $E$} \setcounter{equation}{0} Let $E$ be the total space of the vector bundle $(E,p,M)$. Suppose that $E$ is endowed with a metrical structure $G$ and denote by $D$ the metrical $d-$connection having $h(hh)-$and $v(vv)-$torsions prescribed, given locally by (2.6).

We associate to $D$ the following Einstein equation
\begin{equation}
{\rm{Ric}}(D)-(1/2){\bfr}G=\varkappa \bft,
\end{equation}
where Ric$(D)$ and $\bfr$ denote the Ricci tensor and the scalar curvature of $D$, respectively, $\varkappa$ is a constant and $\bft$ is a tensor field of type $(0,2)$ called the energy-momentum tensor.	

\ms

{\it Remark 3.1.} The tensor field from the left hand of the eq. (3.1) is neither symmetric nor free divergence since $D$ has torsion.

\ms

To express (3.1) by using the curvature of the $d-$connection $D$, let us put $X_a=\{\de_i, \pp_a\}$. Then we have:
\begin{equation}
D_{X_\g}X_\b=\G^\a_{\b\g}X_\a,\ \a,\b,\g,\de ...=1,2,...,n+m,
\end{equation}
\begin{equation}
\bft^\a_{\b\g}=\G^\a_{\b\g}-\G^\a_{\g\b}+W^\a_{\b\g}, {\rm{\ where\ }}[X_\a,X_\b]=W^\g_{\a\b}X_\g,
\end{equation}
\begin{equation}
\bfr^\a_{\b\g\de}=X_\de\G^\a_{\b\g}+\G^\varphi_{\b\g}\G^\a_{\varphi\de}-\g|\de+\G^\a_{\b\varphi}W^\varphi_{\g\de},
\end{equation}
\begin{equation}
{\rm{Ric}}(D)=\bfr_{\b\g}=\bfr^\alpha_{\b\g \a},
\end{equation}
\begin{equation}
\bfr=G^{\a\b}\bfr_{\a\b},
\end{equation}
and the eq. (3.1) becomes:
\begin{equation}
\bfr_{\a\b}-\dfrac12\bfr G_{\a\b}=\varkappa\bft_{\a\b}.
\end{equation}
It results that it is equivalent to the following equations:
\begin{equation}
\begin{array}{l}
\left\{\begin{array}{l}R_{ij}-\dfrac12(R+S)g_{ij}=\varkappa\bft_{ij},\ \overset{1}{P}_{ai}=\varkappa\bft_{ai},\ \overset{2}{P}_{ia}=-\varkappa\bft_{ia}\\ \\
S_{ab}-\dfrac12(R+S)h_{ab}=\varkappa\bft_{ab},\ {\rm{where:\ }}\overset{2}{P}_{ia}=\overset{2}{P}^k_{i\ ka},
\end{array}\right.\\
R_{ij}=R^h_{i\ jh},\ \overset{1}{P}_{ai}=\overset{1}{P}^b_{a\ ib},\ S_{ab}=S^c_{abc},\ R=g^{ij}R_{ij},\ S=h^{ab}S_{ab}.
\end{array}
\end{equation}

All tensor fields from (3.8) are distinguished tensor fields on $E$ i.e. in their laws of transformations to a change of local coordinates, $y^\a$ does not appear explicitly.

The conservation law $D_{X_\a}(\bfr^\a_\b-(1/2)\bfr\de^\a_\b)=0$, where $\bfr^\a_\b=G^{\a\g}\bfr_{\g\b}$ can be written as follows:
\begin{equation}
\left\{
\begin{array}{l}
\[R^i_j-\dfrac12(R+S)\de^i_j\]_{|i}+\overset{1}{P}^a_j|_a=0,\\ \\
\[S^a_b-\dfrac12(R+S)\de^a_b\]|_a-\overset{2}{P}^i_{b|i}=0,
\end{array}
\right.
\end{equation}
where $R^i_j=g^{ik}R_{kj}$, $$S^a_b=g^{ac}S_{cb},\overset{1}{P}^a_j=g^{ab}\overset{1}{P}_{bj},\ \overset{2}{P}^j_b=g^{ij}\overset{2}{P}_{jb}.$$

Generally, the eqs. (3.9) are not identically satisfied; so it appears that the energy-momentum tensor is not conservative.

\ms

{\bf Definition 3.1.} The eqs. $(3.8)$ will be called the Einstein equations on the total space $E$ of the vector bundle $\xi$.

\section{Some particular cases} \setcounter{equation}{0} a) Let us take $\xi=(TM,\tau,M)$, where $M$ is a generalized Lagrange space i.e. $M=(M^n,g_{ij}(x,y))$ (cf. R. Miron [6]). If some additional conditions on $g_{ij}(x,y)$ are fulfilled (see R. Miron [6]) then $g_{ij}(x,y)$ determines an unique nonlinear connection on $(TM,\tau,M)$. Let $(N^i_j(x,y))$ be its local coefficients, $i,j,k,...= 1,2,...,n$, and let $(\de_i,\pp_i)$ be the frame adapted to it. The following Riemannian metric on $TM$ appears as natural:
\begin{equation}
\begin{array}{l}
G=g_{ij}(x,y)dx^i\otimes dx^j+g_{ij}(x,y)\de y^i\otimes\de y^j,{\rm{\ where}}\\ \\
\de y^i=dy^i+N^i_j dx^j.
\end{array}
\end{equation}

As in the general case, a linear $d-$connection on $TM$ is completely determined by a set of functions on $TM$, let say $$F\G=(F^i_{jk},L^i_{jk},M^i_{jk},C^i_{jk}).$$

Let $J$ be the natural almost tangent structure on $TM$ i.e. $$J(\de_i)=\pp_i,\ \ J(\pp_i)=0.$$

\ms

{\bf Definition 4.1.} A linear $d-$connection $D$ on $TM$ is said to be normal if $DJ=0$.

\ms

It is easy to see that a normal linear $d-$connection is characterized by $L^i_{jk}=F^i_{jk}$ and $M^i_{jk}=C^i_{jk}$, so a normal linear $d-$connection is completely determined by $F\G=(F^i_{jk}, C^i_{jk})$, where $F^i_{jk}$ and $C^i_{jk}$ have the laws of transformation like a linear connection and a tensor on $M$, respectively, if the local coordinates are changed.

\ms

{\bf Theorem 4.1.} {\it Given two skew-symmetric tensor fields $T^i_{jk}$ and $S^i_{jk}$, there exists a unique metrical normal linear d-connection on $TM$ which has $T^i_{jk}$ and $S^i_{jk}$ as $h(hh)-$ and $v(vv)-$torsions, respectively. Its local coefficients are as follows:}
\begin{equation}
\left\{
\begin{array}{l}
F^i_{jk}=\dfrac12 g^{ih}(\de_j g_{hk}+\de_k g_{hj}-\de_h g_{jk}+g_{hr}T^r_{jk}-g_{jr}T^r_{hk}+g_{kr}T^r_{jh}),\\ \\
C^i_{jk}=\dfrac12 g^{ih}(\pp_j g_{hk}+\pp_k g_{hj}-\pp_h g_{jk}+g_{hr}S^r_{jk}-g_{jr}S^r_{hk}+g_{kr}S^r_{jh}).
\end{array}
\right.
\end{equation}

\ms

\noindent{\it Proof.} Taking any linear $d-$connection $F\G=(F^i_{jk}, C^i_{jk})$ and imposing the conditions $g_{ij|k}=0$, $g_{ij}|_k=0$, $F^i_{jk}-F^i_{kj}=T^i_{jk}$ and $C^i_{jk}-C^i_{hj}=S^i_{jk}$, one gets that $F^i_{jk}$ and $C^i_{jk}$ are uniquely determined as in (4.2), q.e.d.

Einstein equations associated to the metrical, normal, linear $d-$connection given by the Theorem 4.1 are just the Einstein equations obtained by R. Miron in [6].

\medskip

b) Preserving the hypothesis from a) we only change the metric $G$ as follows:
\begin{equation}
G=g_{ij}(x,y)dx^i\otimes dx^j-g_{ij}(x,y)\de y^i\otimes \de y^j.
\end{equation}
This $G$ is nondegenerate but it is nondefinite. However, the Theorem 4.1 is still true. The Einstein equations associated to the metrical, normal, linear $d-$connection stated by it, written with respect to the adapted frame $(\de_i,\pp_i)$ are as follows:
\begin{equation}
\left\{
\begin{array}{l}
R_{ij}-\dfrac12 (R-S)g_{ij}=\varkappa\bft_{ij},\ \ S_{ij}-\dfrac12(R-S)g_{ij}=\varkappa\bft_{(i)(j)}\\ \\
\overset{1}{P}_{ij}=\varkappa\bft_{(i)j},\ \ \overset{2}{P}_{ij}=-\varkappa\bft_{i(j)},
\end{array}
\right.
\end{equation}
where $R_{ij}=R^k_{ijk}$, $R=g^{ij}R_{ij}$.

$S_{ij}=S^k_{ijk}$, $S=g^{ij}S_{ij}$ and in the right hand appear the components of the energy-momentum tensor with respect to the adapted frame.

\ms

{\it Remark 4.1.} The eqs. (4.4) could be also interesting for physicists because $G$ is nondefinite. Its signature is always $(n,n)$.

\ms

{\it Remark 4.2.} Setting $P(\de_i)=-\de_i$, $P(\pp_i)=—\de_i$ one obtains an almost product structure on $TM$ which satisfies $G(PX,PY)=—G(X,Y)$ for any vector fields $X$ and $Y$ on $TM$ and $G$ given by (4.3). Therefore, $(TM,P,G)$ is an almost hyperbolic manifold.

\medskip

c) Now, let us take $\xi=(E,p,M)$ with dim$\bff=1$. If $(U,\varphi)$ is a local chart on $M$, let $(x^1,...,x^n,x^0)$ the local coordinates of a point $u\in p^{-1}(U)$. These coordinates change as follows (cf. (1.1)):
\begin{equation}
x^{i'}=x^{i'}(x^1,...,x^n),\ x^{'0}=f(x^1,...,x^n)\cdot x^0,
\end{equation}
where $f$ is a real function locally defined on $M$, $f\neq 0$.

The formulas (4.5) show that the manifold $E$ is the most general framework for a unitary projective theory (cf. [8], p. 233).

A nonlinear connection on $\xi$ will be defined by a set of functions $(N_i)$ on $E$ such that $\de_i=\pp_i-N_i\pp_0$ verify $\de_{i'}=(\pp'_i,x^i)\de_i$, where $\pp_0=\pp/\pp x^0$ and $\de_{i'}=\de/\de x^i$.

A linear $d-$connection will be completely determined by the following set of functions on $E$, $F\G=(F^i_{jk},L_k,M^i_j,C)$ where $L_k=L^0_{0k}$, $M^i_j=M^i_{j0}$, $C=C^0_{00}$.

Such a connection has four torsions:
\begin{equation}
T^i_{jk}=F^i_{jk}-F^i_{hj},\ R_{kh}=\de_k N_k-\de_k N_h,\ \overset{1}{P}_j=\pp_0 N_j-L_j,\ \overset{2}{P}^i_j=M^i_j
\end{equation}
and four curvatures:
\begin{equation}
\left\{
\begin{array}{l}
R^i_{jkh}=\de_h F^i_{jk}+F^l_{ik}F^i_{lh}-h|k+M^i_j R_{kh}\\ \\
\wt{R}_{kh}=\de_h L_k-\de_k L_h+CR_{kh}\\ \\
\widehat{P}_k=\pp_0 L_k-\pp_k C+\pp_0(CN_k)\\ \\
\overset{2}{P}^i_{jk}=\pp_0 F^i_{jk}-\pp_k M^i_j+\pp_0(N_k M^i_j)+M^i_h F^h_{jk}-M^h_j F^i_{hk}.
\end{array}
\right.
\end{equation}

The $h-$ and $v-$covariant derivatives are defined as in the general case.

Let be $G=g_{ij}(x^1,...,x^n,x^0)dx^i\otimes dx^i+g_{00}(x^1,...,x^n,x^0)(\de x^0)$ where $\de x^0=dx^0+N_i dx^i$, a Riemannian metric on $E$.

There exists a metrical $d-$connection with $T^i_{jk}$ prescribed. Its local coefficients are as follows:
\begin{equation}
\left\{
\begin{array}{l}
F^{i}_{jk}=\dfrac12 g^{ih}(\de_j g_{hk}+\de_k g_{hj}-\de_k g_{hj}+g_{hr}T^r_{jk}-g_{jr}T^r_{hk}+g_{kr}T^r_{jh})\\ \\
L_k=\dfrac12g^{-1}_{00}\de_k g_{00},\ \ M^i_j=\dfrac12 g^{ik}\pp_0 g_{kj},\ \ C=\dfrac12 g^{-1}_{00}\pp_0 g_{00}.
\end{array}
\right.
\end{equation}

Einstein's equations associated to the metrical $d-$connection given by (4.8), written with respect to the adapted frame, are as follows:
\begin{equation}
\left\{
\begin{array}{l}
R_{ij}-\dfrac12 Rg_{ij}=\varkappa\bft_{ij}\\ \\
\overset{1}{P}_{i}=\varkappa\bft_{i0},\ \ \overset{2}{P}^k_{jk}=-\varkappa\bft_{0j},\ \ Rg_{00}=-2\bft_{00},
\end{array}
\right.
\end{equation}
where $R=g^{ij}R_{ij}$.

The conservation law looks as follows:
\begin{equation}
\left\{
\begin{array}{l}
R^i_{j|i}-\dfrac12 R_{|j}+\pp_0(g^{-1}_{00}\overset{1}{P}_j)-g^{-1}_{00}M^i_j\overset{1}{P}_i=0,\\ \\
\dfrac12\pp_0 R+g^{ij}\overset{2}{P}_{j|i}=0.
\end{array}
\right.
\end{equation}

\ms

{\it Remark 4.3.} If $(M,g_{ij}(x))$ is a Lorentz manifold and we set $G=g_{ij}(x^1,$ $x^2,$ $x^3,$ $x^4)dx^i\otimes dx^j+g_{00}(x^1,...,x^4)(\de x^0)^2$, $i,j=1,...,4$, then the first group from eqs. (4.9) are precisely the Einstein equations for $(M,g_{ij}(x))$. The following two groups can be thought of or interpreted as Maxwell equations. Therefore, a way for developing a unitary projective theory has appeared.

\bigskip

{\bf REFERENCES}

\begin{enumerate}
  \item Anastasiei, M., {\it Some tensorial Finsler structures on the tangent bundle}, An. \c{S}t. Univ. Iasi, s. I a, Supliment T. XXVII, (1981), 9-16.
  \item Ikeda, S., {\it On the Finslerian metrical structures of the gravitational field}, An. \c{S}t. Univ. Iasi, s. I a, T. XXX, (1984), 35-38.
  \item Ikeda, S., {\it Some physical aspects underlying the Lagrangian theory of relativity}, (to appear).
  \item Ishikawa, H., {\it Einstein equation in lifted Finsler space}, II Nuovo cimento, vol. 56 B(2), (1980), 252-262.
  \item Miron, R., {\it Vector bundles. Finsler geometry}, Proc. of the National Semin. on Finsler Spaces. Univ. Brasov, 1982, 147-188.
  \item Miron, R., {\it A Lagrangian theory of relativity}, Preprint. Univ. Timisoara (Romania), 1985, 53 p.
  \item Takano, Y., {\it On the theory of fields in Finsler spaces}, Proc. of International Symposium on Relativity and Unified Field Theory (1975-76). pp. 17-26. (S. N. Bose Inst. of Physical Scs., Calcutta).
  \item Tonnelat, M. A., {\it Les th\'{e}ories unitaires de l'\'{e}l\'{e}ctromagnetisme et de la gravitation}, Gauthier-Villars, Paris, 1965.
\end{enumerate}

\bigskip

\noindent{\footnotesize{\it Received 6.XII.1985}\hfill {\it Faculty of Mathematics}}\\ \hspace*{8cm}{\footnotesize{\it ``Al.~I.~Cuza'' University of Iasi}}


\newpage





\runningauthor={M. ANASTASIEI}
\runningtitle={MODELS OF FINSLER AND LAGRANGE GEOMETRY}
\noindent
\baselineskip 8pt
\noindent{\footnotesize{THE PROCEEDINGS OF THE FOURTH NATIONAL SEMINAR\hfill\break
ON FINSLER AND LAGRANGE SPACES\hfill\break Bra\c{s}ov - 1986, 43--56}}
\vskip 2cm
\baselineskip 11.5pt plus .15pt
\centerline{\bf\Large MODELS OF FINSLER AND}
\vskip .2cm
\centerline{\bf\Large LAGRANGE GEOMETRY}
\vskip .5cm
\centerline{\bf\footnotesize{BY}}
\vskip .5cm
\centerline{\bf\footnotesize{M. ANASTASIEI}}
\vskip .5cm

\setcounter{section}{0}
\section{Introduction} Although closely related to the Riemann geometry, the Finsler geometry had a more slow and sinuous development. Two reasons can he pointed out. Its foundation which is not so firm as of Riemann geometry (a prejudice!) and its too complicated character owing to a lot of differential invariants (a true!). The beginnings were stated by B. Riemann in 1854 (cf. M. Matsumoto [8]). Until 1960 almost all its results had a local character. Since 1960 up to now many efforts to modernize this geometry were made. The theory of connections in fibre bundles has been applied to this aim. In this period the studies in Finsler geometry have progressed very much mainly because three quite distinct models of this geometry were created. These models added to the model ``space of line elements'' introduced by E. Cartan, enriched considerably the area of researches in Finsler geometry. Our aim is to describe these models. We shall begin by giving a definition of Finsler, as well as of Lagrange geometry. Some historical facts which motivate these definitions are pointed out. The necessity and the usefulness of the models in studying Finsler and Lagrange geometry are explained. The model ``space of line elements'' will be only sketched since now it is of historical interest. The models which we call ``principal Finsler bundle'', ``vector Finsler bundle'' and ``almost Hermitian'' will be presented with some details.

It is not our purpose to establish accurately the history of appearance and development of these models. Our lecture is mainly an invitation for studying Finsler and Lagrange geometry by using one of these models. The author is indebted to Prof. Dr. Radu Miron for his helpful advices during the preparation of this lecture.

\section{A definition of Finsler and Lagrange geometry}
\setcounter{equation}{0}

We start with some historical facts (see M. Matsumoto [83]). In a famous lecture (1854), B. Riemann proposed the study of manifolds endowed with the so-called Riemannian metric $ds=\sqrt{g_{ij}(x)dx^i dx^j}$. Before arriving at this metric, he was concerned with the concept of generalized metric $ds=L(x^1,...,x^n,$ $dx^1,...,$ $dx^n)$, shortly $ds=L(x,dx)$, which gives the distance between two points $x$ and $x+dx$. B. Riemann himself gives a concrete example of generalized metric $L(x,y)=((y^1)^4+...+(y^n)^4)^{1/4}$ which satisfies the conditions imposed by him:

\noindent (L1) $\quad\quad$ $L(x,y)>0$ for any $y\neq 0$.

\noindent (L2) $\quad\quad$ $L(x,ay)=aL(x,y)$ for any $a>0$.

\noindent (L3) $\quad\quad$ $L(x,-y)=L(x,y)$.

Then, the notion of generalized metric space completely had been forgotten for almost 60 years. It was rediscovered in a geometrical treatment of the variational calculus about the beginning of this century. The dissertation of P. Finsler from 1918 is remarkable in this respect. He introduced the so-called fundamental tensor $g_{ij}(x,y)=(\pp^2 L^2/\pp y^i\pp y^j)/2$ and the $C-$tensor $C_{ijk}(x,y)=(\pp g_{ij}(x,y)/\pp y^k)$. The equations $C_{ijk}(x,y)=0$ characterize Riemannian metrics among Finslerian metrics. Finsler added a fourth condition on $L$:

\noindent (L4) $\quad\quad$ $g_{ij}(x,y)u^i u^j>0$ for any $u=(u^i)\neq 0$

\noindent This so-called positive definiteness condition is called the regularity condition in the calculus of variations.

A pair $(M,L)$ of an $n-$dimensional manifold $M$ and a general metric $L$ satisfying (L2) and

\noindent (L5) $\quad\quad$ $\det(g_{ij})\neq 0$,

\noindent is called a Finsler space. If $L$ satisfies only (L5), the pair $(M,L)$ is called a Lagrange space. This notion was recently introduced by J. Kern [6]. The other conditions (L1, 3, 4) are necessary in some theorems and in some geometrical theories or applications.

The variables $y^i$, $i=1,2,...n$ from $L(x,y)$ define, from the historical point of view, a direction in the point $(x)$. But these variables can also be thought as parameters and can be taken in a number different by $n$. This fact is a support for the study of the so-called Finsler geometry of vector bundles (R. Miron [10]).

The entities $g_{ij}(x,y)$ and $C_{ijk}(x,y)$ are not tensors in an usual sense because of their dependence upon $y$. However, if a change of local coordinates
\begin{equation}
x^{i'}=x^{i'}(x^1,...,x^n),
\end{equation}
is performed and if suppose that the law of transformation of $y$ is
\begin{equation}
y^{i'}=\dfrac{\pp x^{i'}}{\pp x^i}y^i,
\end{equation}
these entities have laws of transformation similar to that of a tensor of type $(0,2)$, respectively $(0,3)$, on manifold $M$. Such entities are called Finsler tensors. It is noteworthy that if $T(x,y)$ is a Finsler tensor then $\pp T(x,y)/\pp y^i$ is also a Finsler tensor.

The geodesics of a Finsler space $(M,L)$ are the extremals of the variation problem $\de\dd\int^b_a L(x(t);dx/dt)dt=0$. These are the solutions of the differential system
\begin{equation}
\dfrac{d^2 x^i}{dt^2}+\g^i_{jk}\(x,\dfrac{dx}{ds}\)\dfrac{dx^j}{ds}\dfrac{dx^k}{ds}=0,
\end{equation}
where $\g^i_{jk}$ are the Christoffel symbols constructed from $g_{ij}(x,y)$ with respect to $x^i$.

L. Berwald has put $G^i=\g^i_{jk} y^j y^k$ and has considered $G^i_j=\dfrac{\pp G^i}{\pp y^j}$ and $G^i_{jk}=\dfrac{\pp G^i_j}{\pp y^k}$. The laws of transformation of $G^i_j$ and $G^i_{jk}$ are as follows:
\begin{equation}
G^{i'}_{j'}(x',y')=\dfrac{\pp x^{i'}}{\pp x^i}\dfrac{\pp x^j}{\pp x^{j'}}G^i_j+\dfrac{\pp x^{i'}}{\pp x^i}\dfrac{\pp^2 x^i}{\pp x^{k'}\pp x^{j'}}y^{k'},
\end{equation}
\begin{equation}
G^{i'}_{j'k'}(x',y')=\dfrac{\pp x^{i'}}{\pp x^i}\dfrac{\pp x^j}{\pp x^{j'}}\dfrac{\pp x^k}{\pp x^{k'}}G^i_{jk}+\dfrac{\pp x^{i'}}{\pp x^i}\dfrac{\pp^2 x^i}{\pp x^{j'}\pp x^{k'}}.
\end{equation}
So $(G^i_{jk}(x,y))$ changes as a linear connection although it depends by $y$. The above consideration lead to the following

{\bf Definition 2.1.} A set of functions of $(x,y)$ whose law of transformation is like that of a geometric object on $M$ is called a Finsler geometric object.

Therefore, $g_{ij}(x,y)$, $C_{ijk}(x,y)$, $G_{jk}^i(x,y)$ are Finsler geometric objects while $G^i_j(x,y)$ are not so. It is now clear what means a field of Finsler geometric objects. These fields can be also defined as cross-sections in convenient fibre bundles (M. Anastasiei [2], R. Miron and M. Anastasiei [15]). But such an abstract definition has a little use without some interpretations of these fields. So it appears a necessity to construct the models in which a field of Finsler geometric objects to get a convenient interpretation.

A Finsler vector field is a set of functions $(X^i(x,y))$ with the following law of transformation:
\begin{equation}
X^{i'}(x',y')=\dfrac{\pp x^{i'}}{\pp x^i}X^{i}(x,y).
\end{equation}
It is easy to see that $\(\dfrac{\pp X^i}{\pp x^j}\)$ does not define a Finsler geometric object. So it is necessary to consider a derivation of $X^i$ which leads to a Finsler object. Such a (covariant) derivative has been defined by L. Berwald:
\begin{equation}
X^i_{;k}=\dfrac{\pp X^i}{\pp x^k}-G^i_k\dfrac{\pp X^i}{\pp y^j}+G^i_{jk}X^j
\end{equation}
and is called $h-$covariant derivative. If one puts	$\dfrac{\de}{\de x^k}=\dfrac{\pp}{\pp x^k}-G^j_k\dfrac{\pp}{\pp y^j}$, then $X^i_{;k}=\dfrac{\de X^i}{\de x^k}+G^i_{jk}X^j$, equality which reminds an usual formula for a covariant derivative. So it is natural to define: $$g_{ij;k}=\dfrac{\de g_{ij}}{\de x^k}-G^s_{ik}g_{sj}-G^s_{jk}g_{is}.$$ If one puts $X^i_{,j}=\dfrac{\pp X^i}{\pp y^j}$ one obtains a new covariant derivative called $v-$covariant derivative. The triad $(G^i_j,G^i_{jk},0)$ is called the {\it Berwald connection}. This connection is not metrical because $g_{ij;k}\neq 0$ and $g_{ij,k}\neq 0$, too. To look for a metrical connection one must modify the definition of $h-$ and $v-$covariant derivatives. One defines:
\begin{equation}
X^i_{|k}=\dfrac{\de X^i}{\de x^k}+F^i_{kj}X^j,\ X^i|_k=\dfrac{\pp X^i}{\pp y^k}+C^i_{kj}X^j,
\end{equation}
where $\dfrac{\de}{\de x^k}=\dfrac{\pp}{\pp x^k}-N^j_k\dfrac{\pp}{\pp y^j}$, $N^j_k=F^j_{ik}y^k$, $F^i_{jk}$ is a set of functions of $(x,y)$ which changes like a linear connection on $M$ and $C^i_{jk}$ is a Finsler tensor of type $(1,2)$ on $M$.

If the equalities $C^i_{jk}=C^i_{kj}$ and $F^i_{jk}=F^i_{kj}$ are assumed, then from the conditions $g_{ij|k}=0$ and $g_{ij}|_k=0$ by a ``Christoffel process'' one obtains
\begin{equation}
\left\{\begin{array}{l}N^i_j=G^i_j\\ \\ F^i_{jk}=\dfrac12 g^{is}\(\dfrac{\de g_{sk}}{\de x^j}+\dfrac{\de g_{sj}}{\de x^k}-\dfrac{\de g_{jk}}{\de x^s}\)\\ \\ C^i_{jk}=\dfrac12 g^{is}\(\dfrac{\pp g_{sk}}{\pp y^j}+\dfrac{\pp g_{sj}}{\pp y^k}-\dfrac{\pp g_{jk}}{\pp y^s}\)=\dfrac12 g^{is}\dfrac{\pp g_{sj}}{\pp y^k}.\end{array}\right.
\end{equation}

The triad $(G^i_j,F^i_{jk},C^i_{jk})$ whose elements are given by (2.9) is called the {\it Cartan connection} of the Finsler space $(M,L)$. More general, a triad $(N^i_j(x,y),$ $F^i_{jk}(x,y),$ $C^i_{jk}(x,y))$, where $N^i_j$ has a law of transformation similar to $G^i_j$,  $F^i_{jk}$ has a law of transformation similar to $G^i_{jk}$ and $C^i_{jk}$ is a Finsler tensor field, is called a {\it Finsler connection}. The definition of $v-$ and $h-$covariant derivative is similar to (2.8). The commutation formulae lead to the five torsion Finsler tensors:
\begin{equation}
\begin{array}{l}
T^i_{jk}=F^i_{jk}-F^i_{kj},\ R^i_{jk}=\dfrac{\de N^i_j}{\de x^k}-\dfrac{\de N^i_k}{\de x^j},\ P^i_{jk}=\dfrac{\pp N^i_j}{\pp y^k}-F^i_{kj},\\ \\
\ S^i_{jk}=C^i_{jk}-C^i_{kj} and C^i_{jk}.
\end{array}
\end{equation}
and three curvature Finsler tensors:
\begin{equation}
\begin{array}{l}
R^i_{h'jk}=\dfrac{\de F^i_{hj}}{\de x^k}-\dfrac{\de F^i_{hk}}{\de x^j}+ F^r_{hj}F^i_{rk}-F^r_{hk}F^i_{rj}+C^i_{hr}R^r_{jk},\\ \\
P^i_{hjk}=\dfrac{\pp F^i_{hj}}{\pp y^k}-C^i_{hk|j}+C^i_{hr}P^r_{jk},\\ \\
S^{\ i}_{h\ jk}=\dfrac{\pp C^i_{hj}}{\pp y^k}-\dfrac{\pp C^i_{hk}}{\pp y^j}+C^r_{hj}C^i_{rk}-C^r_{hk}C^i_{rj}.
\end{array}
\end{equation}
Then Bianchi identities can be established and the Finsler spaces with special properties can be studied. We conclude with the following definition of Finsler (Lagrange) geometry.

{\bf Definition 2.2.} We call Finsler (Lagrange) geometry the study of Finsler geometric objects on a manifold $M$ endowed with a general homogeneous (non-homogeneous) metric $L$.

The content of the classical Finsler geometry has mainly been obtained by using the methods discussed above. It is well represented by H. Rund's book [19].

\section{The model ``space of line elements''. Nonlinear connections}\setcounter{equation}{0}

E. Cartan has arrived at his connection by creating a model of Finsler geometry. As we have seen, the quantities appearing in Finsler geometry depend on $2n$ variables $x=(x^i)$ and $y=(y^i)$. E. Cartan calls the pair $(x,y)$ the supporting element of these quantities and considers the set $M'$ of all the supporting elements. Owing to the homogeneity, $y$ of a supporting element $(x,y)$ is an oriented direction in $x$ so $M'$ is a $(2n-1)-$dimensional manifold called the space of line elements. E. Cartan considers the Finsler geometry as the geometry of the manifold $M'$ and identifies a Finsler connection to an Euclidian connection on $M'$. By using four axioms he determines what is now called the Cartan connection i.e. a metrical Finsler connection completely determined by the fundamental function $L$ (see (2.9)). In this model the Finsler geometric objects are geometric objects defined on $M'$. Later, instead of $M'$ was considered the total space $TM$ of the tangent bundle over $M$ or $TM-0$ when the homogeneity is taken into account for.

The models of Finsler geometry created after 1960 have had mainly two purposes. The first one was to give a clear meaning to the notion of Finsler geometric object and the second one was to establish a global definition for connections in Finsler spaces and to re-examine E. Cartan'system of axioms. In all three models which we shall discuss in the following the notion of nonlinear connection appears, in two of them this notion being in a central place. The importance of the nonlinear connections was late recognized. Now there exist a lot of equivalent definitions of this notion. We shall give some here (cf. R. Miron, M. Anastasiei [43]).

1. Let $\tau:TM\to M$ be the tangent bundle to $M$, $\tau'$ its tangent map and $VTM=$ker$\tau'$. A nonlinear connection on $TM$ is a subbundle $HTM\subset TTM$ such that $TTM=HTM\oplus VTM$.

2. Let $\tau^{-1}(TM)\to TM$ be the pull-back of $TM$ by $\tau$. Let us denote by $\pi:TTM\to TM$ the projection map. The following sequence of vector bundles over $TM$ is exact: $$0\to VTM\overset{i}\longrightarrow TTM\overset{l=(\pi,\tau')}\longrightarrow\tau^{-1}(TM)\to 0.$$ A splitting $\G:\tau^{-1}(TM)\to TTM$ of this exact sequence is a
nonlinear connection on $TM$.

3. Let $J$ be the natural almost tangent structure on $TM$. A nonlinear connection on $TM$ is a tensor field $P$ of type $(1,1)$ on $TM$ such that $PJ=-J$ and $JP=J$ hold.

4. A nonlinear connection on $TM$ is an almost product structure $P$ on $TM$ which satisfies $P(X)=-X$ for any vertical vector field $X$.

5. A nonlinear connection is a set of functions $(N^i_j)$ on $TM$ which has the law of transformation like $G^i_j$ in (2.4).

A nonlinear connection always exists if the paracompactness of $M$ is assumed. If $(M,L)$ is a Finsler or a Lagrange space, then there exists a canonical nonlinear connection determined by $L$ only (cf. Section 2).

\section{The ``principal Finsler bundle'' model} \setcounter{equation}{0}

As it is well known, a linear connection on a manifold $M$ can be defined as a certain distribution on the linear frame bundle $L(M)$. A point of $L(M)$ is regarded as a pair $(x,z)$ of a point $x$ of $M$ and a frame $z$ in $x$. The Finsler geometric objects are special functions on $TM$ so they depend on $(x,y)$, a pair of a point $x$ of $M$ and a tangent vector $y$ at $x$. Therefore, the set of triads $(x,y,z)$ may be a good foundation for Finsler geometry. Such a set is obtained as follows. Let $\tau^{-1}(L(M))\to TM$ be the pull-back of $L(M)$ by $\tau$. This is a principal bundle over $TM$ with structural group $GL(n,R)$. It is called the Finsler bundle of $M$ and it will be denoted by $F(M)$. Its total space is $F=\{(y,z)\in TM\times L(M),\tau(y)=\pi(z)\}$. A right translation $\b_g$ of $F,g\in GL(n,R)$ is given as: $u=(y,z)\to ug=(y,zg)$. The bundle $F(M)$ was introduced by L. Auslander ([4]). It was also considered by H. Akbar-Zadeh [1]. $F(M)$ was called a Finsler bundle of $M$ by M. Matsumoto. He also used it systematically and efficiently as a model of Finsler geometry (see his monograph [7]). Let $(R^n)^r_s$ be the space of tensors of type $(r,s)$ over $R^n$. A Finsler tensor field of type $(r,s)$ is defined as a map $K:F\to (R^n)^r_s$ which satisfies a condition $K\circ \b_g=g^{-1}K$ for any $g\in GL(n,R)$. This definition is equivalent to a classical one (M. Matsumoto [7], p. 49). A Finsler connection $F$ on a manifold $M$ is a pair $(\G,N)$ of a connection $\G$ in $F(M)$ and a nonlinear connection $M$ on the tangent bundle $TM$. Such a definition is very general because $\G$ and $N$ are not yet related. A second definition of a Finsler connection ([7], p. 63) leads to a definition of $v-$ and $h-$covariant derivatives quite similar to that with respect to a linear connection. Torsions and curvatures are obtained by a suitable generalization of the structure equations of a linear connection to a Finsler connection ([7], p. 70-76).

All Bianchi identities are obtained from some general identities. If $M$ is endowed with a fundamental function $L$, then the nonlinear connection is completely determined by $L$. The following theorem of M. Matsumoto holds:

{\bf Theorem. } {\it The Cartan connection $C\G$ of a Finsler space $(M,L)$ is uniquely determined by the five axioms as follows:

(1)	$h-$metrical: $g_{ij\ |k}=0$

(2)	without $h-$torsion: $T=0$ $(F^i_{jk}=F^i_{kj})$

(3)	$v-$metrical: $g_{ij} |_{k}=0$

(4)	without $v-$torsion: $S^1=0$ $(C^i_{jk}=C^i_{kj})$

(5)	$D^i_{\ j}=N^i_j-F^i_{kj}y^k=0$ (the deflection tensor vanishes).}

The ``principal Finsler bundle'' model leads to a clear definition of the notion of Finsler geometric object, to a global definition of Finsler connections from which all classical Finsler connections are derived and to an interesting theory of transformations of Finsler spaces. Of course, there are some problems which can not be solved nor attacked by using this model. For instance the theory of Finsler spaces with constant curvature or the theory of subspaces in Finsler spaces.

\section{The ``vector Finsler bundle'' model}\setcounter{equation}{0}

In this model the base manifold is $TM$, too. The pull-back $\tau^{-1}(TM)\to TM$ of $TM$ by $\tau$ will be called the vector Finsler bundle of $M$. If $(x^i)$ is a coordinate system on $M$ and $(x^i,y^i)$ is the coordinate system on $TM$ induced by it, then $\(\dfrac{\pp}{\pp x^i},\dfrac{\pp}{\pp y^i}\)$ is a basis in $T_u TM$, $u\in TM$, and $\(\dfrac{\pp}{\pp y^i}\)$ is a basis in $\tau^{-1}(TM)$, the fiber of $\tau^{-1}(TM)$ over $u$. It follows easily that the vector Finsler bundle is isomorphic to the vertical subbundle. We denote by $v$ its inclusion map in $TTM$. A cross-section $\overline{X}$ of the vector Finsler bundle has the local form $\overline{X}=\overline{X}^i\dfrac{\pp}{\pp y^i}$. Since $\dfrac{\pp}{\pp y^i}=\dfrac{\pp x^{i'}}{\pp x^i}\dfrac{\pp}{\pp y^{i'}}$ it follows that the set of functions $(\overline{X}^i(x,y))$ defines a Finsler vector field. More general, the tensorial algebra on the vector Finsler bundle is a model for the algebra of Finsler tensor fields.

A Finsler connection is defined as a regular connection in the vector Finsler bundle. Let be $\nabla:{\cal X}(TM)\times S(\tau^{-1}(TM))\to S(\tau^{-1}(TM))$, $(X,Y)\to\nabla_X Y$ a linear connection in the vector Finsler bundle. Let be $C=y^i\dfrac{\pp}{\pp y^i}$ the canonical field (Liouville) on $TM$. A vector field $X$ on $TM$ is called horizontal if $\nabla_X C=0$. Let be $H_u$ the subspace of horizontal vectors and $V_u$ the subspace of vertical vectors. The connection $\nabla$ is called regular if $T_uTM=H_u\oplus V_u$ for any $u\in TM$. Such a decomposition of $T_uTM$ defines a splitting of the exact sequence from Section 3, hence a nonlinear connection on $TM$. If $\nabla$ is a regular connection, then $\tau'$ is an isomorphism on $H_u$. Let us denote by $h_u$ the map $(\tau'/H_u)^{-1}:T_{\tau(u)}M\to H_u$ and let us put $h_u\(\dfrac{\pp}{\pp x^i}\)=\dfrac{\de}{\de x^i}$. It results $\tau'_u\(\dfrac{\de}{\de x^i}-\dfrac{\pp}{\pp x^i}\)=0$ because of $\tau'_u\circ h_u=$identity and of $\tau'\(\dfrac{\pp}{\pp x^i}\)=\dfrac{\pp}{\pp x^i}$. Therefore, $\dfrac{\de}{\de x^i}-\dfrac{\pp}{\pp x^i}$ are vertical vector fields. We may write $\dfrac{\de}{\de x^i}=\dfrac{\pp}{\pp x^i}-N^j_i\dfrac{\pp}{\pp y^j}$ because $V_u$ is spanned by $\(\dfrac{\pp}{\pp y^j}\)$. A linear connection on the vector Finsler bundle is locally given as follows: $$\nabla_{\frac{\pp}{\pp x^k}}\dfrac{\pp}{\pp y^j}=\G^i_{jk}\dfrac{\pp}{\pp y^i},\ \nabla_{\frac{\pp}{\pp y^k}}\dfrac{\pp}{\pp y^j}=C^i_{jk}\dfrac{\pp}{\pp y^i}.$$ From the equality $\nabla_{\frac{\pp}{\pp y^i}}C=(\de^j_i+y^k C^j_{ki})$ it follows that $\nabla$ is regular if and only if the matrix $(\de^j_i+y^k C^j_{ki})$ is regular for any $y$. The condition $\nabla_{\frac{\de}{\de x^i}}C=0$ is equivalent to $N^k_{i}(\de^i_k+y^s C^j_{sk})=\G^j_{si}y^s$. It follows again that the regularity condition on $\nabla$ allows the determination of a nonlinear connection $(N^i_j)$. If we put $F^i_{jk}=\G^i_{jk}-N^p_k C^i_{jp}$ then it results that $(N^i_j,F^i_{jk},C^i_{jk})$ defines a Finsler connection in the classical sense. Therefore, any regular connection in the vector Finsler bundle is a model of a Finsler connection. The curvatures and torsions of such a Finsler connection are obtained as follows.

Let be $\wt{R}(X,Y)Z=\nabla_X\nabla_Y Z-\nabla_Y\nabla_X Z-\nabla_{[X,Y]}Z$ the curvature of $\nabla$. The following tensor fields: $R(\bar{X},\bar{Y})\bar{Z}=\wt{R}(h\bar{X},h\bar{Y})\bar{Z}$, $P(\bar{X},\bar{Y})\bar{Z}=\wt{R}(v\bar{X},h\bar{Y})\bar{Z}$, $S(\bar{X},\bar{Y})\bar{Z}=\wt{R}(v\bar{X},v\bar{Y})\bar{Z}$ are Finsler tensor fields and they are models for the curvatures of a Finsler connection. The tensor field $T(X,Y)=\nabla_X\ell Y-\nabla_Y\ell X-\ell[X,Y]$ , where $\ell$ is the horizontal lift, is called the torsion of $\nabla$. It results that $T(hX,hY)$, $T(hX,vY)$ and $T(vX,vY)$ are Finsler tensor fields, models for the three torsions of a Finsler connection. The others are $R^i_{jk}$ given by $\[\dfrac{\de}{\de x^j},\dfrac{\de}{\de x^k}\]=R^i_{kj}\dfrac{\pp}{\pp y^i}$ and  $C^i_{jk}$.

The model of a metrical Finsler connection is a regular connection $\nabla$ which verifies $\nabla g=0$. We remark that as a model of Finsler geometry can also serve the vertical subbundle which is isomorphic to the vectorFinsler bundle. This model appears in a paper by V. Oproiu [17]. The ``vectorFinsler bundle'' model was systematically used by B.T. Hassan [5]. A generalization of it was treated by D. Opris [16].

\section{The ``almost Hermitian'' model} \setcounter{equation}{0}

This model was pointed out by R. Miron and it was used by him for an interesting theory of Finslerian relativity (R. Miron [12]). The base manifold is $TM$ furnished with a nonlinear connection determined by the fundamental function or not. The tensorial Finsler fields have as models the elements of the tensorial algebra of the bundle $H\oplus V\to TM$. Let $\(\dfrac{\de}{\de x^i},\dfrac{\pp}{\pp y^i}\)$ be the local frame adapted to the decomposition $T_uTM=H_u\oplus V_u$. Putting $F\(\dfrac{\de}{\de x^i}\)=-\dfrac{\pp}{\pp y^i}$, $F\(\dfrac{\pp}{\pp y^i}\)=\dfrac{\de}{\de x^i}$, one obtains an almost complex structure $(F^2=-I)$ on $TM$. A model for a Finsler connection is a linear connection $D$ on $TM$ which satisfies the following two conditions:

(1)	$D$ preserves by parallelism the horizontal distribution $u\to H_u$ as well as the vertical distribution $u\to V_u$.

(2)	$DP=0$.

Indeed, the first condition leads to the following local form of $D$:	$$\begin{array}{l}D_{\frac{\de}{\de x^k}}\dfrac{\de}{\de x^j}=F^i_{jk}\dfrac{\de}{\de x^j},\ D_{\frac{\pp}{\pp y^k}}\dfrac{\de}{\de x^j}=\wt{C}^i_{jk}\dfrac{\de}{\de x^i},\\ \\ D_{\frac{\de}{\de x^k}}\dfrac{\pp}{\pp y^j}=\wt{F}^i_{jk}\dfrac{\pp}{\pp y^i},\ D_{\frac{\pp}{\pp y^k}}\dfrac{\pp}{\pp y^j}=C^i_{jk}\dfrac{\pp}{\pp y^i},\end{array}$$ and the second one gives $\wt{F}^i_{jk}=F^i_{jk}$, $\wt{C}^i_{jk}=C^i_{jk}$. So $D$, which satisfies (1) and (2), is a model for the Finsler connection $(N^i_j,F^i_{jk},C^i_{jk})$. As a model of the fundamental tensor $(g_{ij})$ is taken the tensor field $G=g_{ij}dx^i\otimes dx^j+g_{ij}\de y^i\otimes \de y^j$, where $\de y^i=dy^i+N^i_k dx^k$. This $G$ is called the $N-$lift of $(g_{ij})$. By a direct calculation it follows that $G(FX,FY)=G(X,X)$ for any vector fields $X,Y$ on $TM$. Therefore $(F,G)$ defines an almost Hermitian structure on $TM$. The triad $H^{2n}=(TM,G,F)$ is called the ``almost Hermitian'' model of a Finsler space $(M,L)$ or a Lagrange space $(M,{\cal L})$. The term is also justified by the following theorem:

{\it Let $G$ be a Riemannian metric on $TM$ of rank $n$ on the vertical distribution and let $N$ be the distribution supplementary and orthogonal to the vertical distribution. Let $F$ be the almost complex structure determined by $N$. If $(G,F)$ is an almost Hermitian structure then there exists an unique fundamental tensor $(g_{ij})$ whose $N-$lift is $G$.}

\noindent{\it Proof.} The distribution $N$ spanned by $\dfrac{\de}{\de x^i}$ is determined from the equations $G\(\dfrac{\de}{\de x^i},\dfrac{\pp}{\pp y^j}\)=0$. Locally, $G$ is as follows: $G=h_{ij}dx^i\otimes dx^j+g_{ij}\de y^i\otimes \de y^j$. From $G(FX,FY)=G(X,Y)$ it results $h_{ij}=g_{ij}$ is obvious that the $N-$lift of $g_{ij}$ is just $G$.

The following theorem holds (R. Miron [12]):

{\it The model $H^{2n}$ is Hermitian if and only if $$R^i_{jk}=\dfrac{\de N^i_k}{\de x^j}-\dfrac{\de N^i_j}{\de x^k}=0,\ t^i_{jk}=\dfrac{\pp N^i_k}{\pp y^j}-\dfrac{\pp N^i_j}{\pp y^k}=0.$$}

A metrical Finsler connection is a linear connection $D$ on $TM$ which satisfies and the third condition:

(3) $DG=0$.

There exists a unique metrical Finsler connection such that $T(hX,hY)=0$ and $T(vX,yY)=0$, where $T$ is its torsion.

This connection coincides to the Cartan connection when the nonlinear connection $N$ is determined by the fundamental function $L$. On $TM$ there exists also a symplectic structure defined by $\phi(X,Y)=G(X,FY)$. It is easy to see that a metrical Finsler connection is also a symplectic one $(D\phi=0)$. If $M$ is a Finsler space its model $H^{2n}$ is almost K\"{a}hler i.e. $d\phi=0$ (M.Matsumoto [9]). This result is also valid if $M$ is a Lagrange space (V. Oproiu [18]).

The ``almost Hermitian'' model suggests at least two generalizations studied until now. The first one is the considering of the linear connections $D$ on $TM$ which preserve the horizontal and vertical distributions but do not verify $DF=0$. Locally, such a connection has four distinct components. A Lagrangian theory of relativity by using such a connection was developed (R. Miron, S. Watanabe, S. Ikeda [13]). The second one is the considering of the geometry of the total space of a vector bundle (R. Miron [10]) or a principal bundle (M. Anastasiei [3]).

Quite recent it was observed the importance of the study of the pair $(M,g_{ij}(x,y))$, where $g_{ij}(x,y)$ is not provided by a fundamental function $L$ (R. Miron [10]). The models described above can also be used for studying such a spaces $(M,g_{ij}(x,y))$ called generalized Lagrange spaces (R. Miron [12]).

The ``almost Hermitian'' model is very complex so it allows to obtain more information about Finsler and Lagrange spaces.

The models which we just described rise a lot of new problems for Finsler and Lagrange geometry and allow the solving of the elder problems which could not be solved in a classical treatment. These models suggest also various generalizations and applications of Finsler and Lagrange geometry.

\medskip

{\bf REFERENCES}

\begin{enumerate}
  \item Akbar-Zadeh H., {\it Les espaces de Finsler et certain de leurs g\'{e}n\'{e}ralisations}, Ann. Sc. Ecole Norm. Sup. (3), 80, 1963, p. 1-79.
  \item Anastasiei M., {\it Finsler geometric objects and their Lie derivative}, Proc. Nation. Sem. on Finsler spaces, 1981, p. 11-25.
  \item Anastasiei M., {\it Finsler geometry on principal bundles}, Lucr. Conf. Nation. Geom. Topologie, Iasi, 1984, p.70-77.
  \item Auslander L., {\it On curvature in Finsler geometry}, Trans. Amer. Math. Soc. 79(1955), p. 378-388.
  \item Hassan B.T.M., {\it The theory of geodesics in Finsler spaces}. Thesis, Southampton, 1967, 108 p.
  \item Kern J., {\it Lagrange geometry}, Arch. Math. 25, 1974, p. 438-443.
  \item Matsumoto M., {\it Foundations of Finsler geometry and special Finsler spaces} (unpublished), Kyoto, 1982, 373 p.
  \item Matsumoto M., {\it A history of Finsler geometry}, Appendix to 7.
  \item Matsumoto M., {\it Connections, metrics and almost complex structures of the tangent bundles}, J. Math. Kyoto Univ., 5, 1966, p. 251-258.
  \item Miron R., {\it Vector bundles Finsler geometry}, Proc. Nation. Semin. on Finsler spaces, Brasov (Romania), 1982, p. 147-188.
  \item Miron R., {\it Metrical Finsler structures and metrical Finsler connections}, J. Math. Kyoto Univ., 23-2, 1983, p. 219-224.
  \item Miron R., {\it A Lagrangian theory of relativity}, Preprint nr. 84, Univ. Ti\-mi\-soa\-ra (Romania),1985, 53 p.
  \item Miron R., Watanabe S., Ikeda S., {\it A theory of relativity} (to appear).
  \item Miron R., Anastasiei M., {\it Fibrate vectoriale. Spa\c{t}ii Lagrange. Aplica\c{t}ii \^in teoria relativit\u{a}\c{t}ii}, to appear in Editura Acad. R.S. Romania, 300 p.
  \item Miron R., Anastasiei M., {\it On the notion of Finsler geometric object}, Mem. Sect. Stiint. Acad. R.S. Romania, s.IV, t. IV, 1, 1981, 25-31.
  \item Opris D., {\it Fibr\'{e}s vectoriels de Finsler et connexions associ\'{e}es}, Proc. Nation. Semin. on Finsler spaces, 1980, p. 184-193.
  \item Oproiu V., {\it Some properties of the tangent bundle related to the Finsler geometry}, Proc. Nation. Semin. on Finsler spaces. Brasov (R.S. Romania), 1980, p. 195-207.
  \item Oproiu V., {\it A Riemannian structure in Lagrange geometry}, to appear.
  \item Rund H., {\it The differential geometry of Finsler spaces}, Springer-Verlag, Berlin, 1959.
\end{enumerate}

\medskip

\hspace*{8.1cm}{\footnotesize{{\it Faculty of Mathematics}}}\\ \hspace*{8cm}{\footnotesize{\it ``Al.~I.~Cuza'' University of Iasi}}


\newpage





\runningauthor={M. ANASTASIEI}
\runningtitle={CONSERVATION LAWS IN THE $\{V,H\}-$BUNDLE MODEL OF RELATIVITY}
\noindent
\baselineskip 8pt
\noindent{\footnotesize{TENSOR, N.S.\hfill\break Vol. 46 (1987), 323--328}}
\vskip 2cm
\baselineskip 11.5pt plus .15pt
\centerline{\bf\Large CONSERVATION LAWS IN THE}
\vskip .2cm
\centerline{\bf\Large $\{V,H\}-$BUNDLE MODEL OF RELATIVITY\footnotetext[1]{Received January 27, 1987}}
\vskip .5cm
\centerline{\bf\footnotesize{Dedicated to the Memory of}}
\vskip .1cm
\centerline{\bf\footnotesize{Professor Dr. Akitsugu Kawaguchi, Founder of the Tensor Society}}
\vskip .5cm
\centerline{\bf by Mihai ANASTASIEI}
\vskip 1cm

\renewcommand{\thefootnote}{\arabic{footnote}}

\setcounter{section}{0}
\section{Introduction}

In a recent paper [1]\footnote[1]{Number in brackets refer to the references at the end of the paper} we have considered the Einstein equations on the total space of a vector bundle in order to obtain a Finslerian unitary projective theory as an extension of the Finslerian theory of relativity developed by R. Miron [5]. We have written the corresponding conservation laws which here, generally, are not identities, i. e., it is of interest to find vector bundles for which such a conservation laws are identically verified. In this paper we take, into considerations, the vector bundles whose type fibers are finite dimensional normed spaces $V$ and, moreover, admit the reductions to a subgroup $H$ of the group of the automorphisms of $V$ which preserve the norm, shortly, $\{V,H\}-$bundles. This class of vector bundles, which contains the tangent bundles to the $\{V,H\}-$manifolds of Y. Ichijyo [2], may be of own interest so we treat it with some details in \S3. In \S2 we give necessary preliminaries from the geometry of the total space of a vector bundle (cf. [7], [8]). The conservation laws on the $\{V,H\}-$bundles are discussed in \S4.

\section{Vector bundles. Metrical $d-$connections}

Let $M$ be an $n-$dimensional differentiable (of class $C^\9$) manifold and $\xi=(E,p,M)$, $p:E\to M$, a vector bundle whose type fiber is a vector space $V$ isomorphic to $R^m$. Let $\{(U_a,\wt{\phi}_a,R^n)\}$ be an atlas on $M$ and let $\{(U_a,\phi_a,V)\}$ be a bundle atlas of $\xi$ i.e. $\phi_a:p^{-l}(U_a)\to U_a\times V$ are bijective mappings such that $p\phi^{-1}_a(x,v)=x$ for $x\in M$ and $v\in V$, and the applications $g_{\b,x}(x)=\phi_{\b\a}\circ\phi^{-1}_{\a,x}:U_\a\cap U_\b\to V$ are differentiable. The manifold structure of $E$ is defined by the differentiable atlas $\{(p^{-1}(U_\a),h_\a)\}$, where $h_\a:p^{-l}(U_\a)\to R^n\times V$ is given as $h_\a(u)=(\wt{\phi}_\a(u),\phi_{\a,p(u)}(u))$. If we set $h_\a(u)=(x^h, y^a)$ and $h_\b(u)=(x^{k'},y^{a'})$, $k,k'=1,...,n$; $a,a'=1,...,m$, then $h_\b\circ h_\a^{-1}$ is as follows:
\begin{equation}
\left\{
\begin{array}{l}
x^{k'}=x^{k'}(x^1,...,x^n),\det\(\dfrac{\pp x^{k'}}{\pp x^k}\)\neq 0,\\ \\
y^{a'}=y^a\cdot S^{a'}_a(x^1,...,x^n),\ \|S^{a'}_a(x)\|\in GL(m,R).
\end{array}
\right.
\end{equation}
The transformation law (2.1) of the local coordinates on $E$ shows that $E$, for $m=1$ or $m=2$, is the most general framework for an unitary projective theory (cf. [10] p. 233).

Let $p^T:TE\to TM$ be the differential of $p$. Then ker$p^T=VE$ is a subbundle of $TE\to TM$ called the {\it vertical subbundle}.

{\bf Definition 2.1.} A nonlinear connection on $\xi$ is a subbundle $HE$ of $TE\to TM$ such that $TE=HE\oplus VE$.

The local fiber $H_uE$ of the vector bundle $HE\to TM$ is spanned by $(\de_k)$ given as follows:
\begin{equation}
\de_k=\pp_k-N_k^a(x,y)\pp_a.
\end{equation}
Here $\pp_k$ and $\pp_a$ stand for $\dfrac{\pp}{\pp x^k}$ and $\dfrac{\pp}{\pp y^a}$, respectively. The functions $N^a_k(x,y)$ are called the {\it local coefficients of the nonlinear connection}. This set of functions has appeared, for the particular case $E=TM$, early in the development of Finsler geometry, but the first who recognized its importance and treated it as defining a nonlinear connection was A. Kawaguchi ([3], [4]). The mapping $u\to H_uE$ $(u\to V_uE)$ is called the horizontal (vertical) distribution. The local frame $(\de_k,\pp_a)$ is adapted to the horizontal and the vertical distributions. Its dual is $(dx^k,\de y^a)$, where $\de y^a=dy^a+N^a_k dx^h$. The tensorial algebra spanned by $1,\de_k,\pp_a, d x^k,\de y^a$ is called the {\it algebra of d-tensor fields} ([6], [7]).

{\bf Definition 2.2.} A linear connection $D$ on $E$ is a $d-$connection if it preserves over parallel displacement of the horizontal and the vertical distributions.

Every $d-$connection $D$ on $E$ can be locally given as follows:
\begin{equation}
\left\{
\begin{array}{l}
D_{\de_k}\de_j=F^i_{jk}(x,y)\de_i,\ \ D_{\de k}\pp_b=L^a_{bk}(x,y)\pp_a,\\ \\
D_{\pp_a}\de_j=M^i_{ja}(x,y)\de_i,\ \ D_{\pp_c}\pp_b=C^a_{bc}(x,y)\pp_a.
\end{array}
\right.
\end{equation}
The coefficients $F^i_{jk}(x,y)$ and $L^a_{bk}(x,y)$ change under (2.1) like the coefficients of a linear connection on $M$ and on $\xi$, respectively, although they depend on $y^a$; $M^i_{ja}(x,y)$ and $C^a_{bc}(x,y)$ are defining tensor fields on $E$. Conversely, a set of coefficients
$L\G=(F^i_{jk}(x,y)$, $L^a_{bk}(x,y)$, $M^i_{ja}(x,y)$, $C^a_{bc}(x,y))$ which change under (2.1) as the above determines an unique $d-$connection on $E$. We shall denote by $_|$ and $|$ the $h-$ and
$v-$covariant derivative associated with the $d-$connection $D$ (see [6]). The commutation or Ricci formulae introduce for a $d-$connection five torsion $d-$tensor fields:
\begin{equation}
\left\{
\begin{array}{l}
T^i_{\ jk}=F^i_{jk}-F^i_{kj},\ \ R^a_{\ jk}=\de_k N^a_j-\de_j N^a_k,\\ \\
\overset{1}{P}^a_{\ jb}=\pp_b N^a_j-L^a_{bj},\ \ \overset{2}{P}^i_{\ jb}=M^i_{jb},\ S^a_{\ bc}=C^a_{bc}-C^a_{cb},
\end{array}
\right.
\end{equation}
and six curvature $d-$tensor fields:
\begin{equation}
\left\{
\begin{array}{l}
R^{\ i}_{j\ kh}=\de_h F^i_{jk}+F^l_{jk}F^i_{lh}-k/h+M^i_{ja}R^a_{\ kh},\\ \\
\wt{R}^{\ a}_{b\ kh}=\de_h L^a_{bk}+L^c_{bk}L^a_{ch}-k/h+C^a_{bc}R^c_{kh},\\ \\
\overset{1}{P}^{\ a}_{b\ kc}=\pp_c L^a_{bk}-C^a_{bc|k}+C^a_{bd}P^d_{\ kc},\\ \\
\overset{2}{P}^{\ i}_{j\ kc}=\pp_c F^i_{jk}-M^i_{jc|k}+M^i_{jb}P^d_{\ kc},\\ \\
M^i_{j\ bc}=\pp_c M^i_{jb}+M^h_{jb}M^i_{hc}-b/c,\\ \\
S^{\ a}_{b\ cd}=\pp_d C^a_{bc}+C^e_{bc}C^a_{ed}-c/d,
\end{array}
\right.
\end{equation}
where $k/h, b/c,c/d$ mean the interchange of indices in the foregoing terms.

A metrical structure on $E$ is a tensor field $G$ of type $(0,2)$ on $E$, symmetric and nondegenerate. It determines an unique nonlinear connection on $\xi$ by taking into consideration the distribution which is orthogonal to the vertical distribution, with respect to it. In the frames adapted to this nonlinear connection, $G$ can be expressed as follows:
\begin{equation}
G(x,y)=g_{ij}(x,y)dx^i dx^j+\wt{g}_{ab}(x,y)\de y^a\de y^b.
\end{equation}
A $d-$connection $D$ on $E$ is metrical with respect to $G$ if $DG=0$. It is easy to prove that a $d-$connection $D$ is metrical if and only if
\begin{equation}
g_{ij|k}=0,\ \ g_{ij}|_a=0,\ \ \wt{g}_{ab}|_c=0,\ \ \wt{g}_{ab|k}=0
\end{equation}
hold. There exists a metrical $d-$connection which has the torsion fields $T^i_{jk}$ and $S^a_{bc}$ prescribed and which is unique in a certain sense [1]. Another metrical $d-$connection will be constructed in \S4.

\section{$\{V,H\}-$bundles}\setcounter{equation}{0} Let $\xi$ be the vector bundle from \S2. Suppose that its fiber $V$ is endowed with a norm $\|\cdot\|:V\to R_+$, i.e. $V$ is a Minkowski space. If $v=v^a e_a$, where $(e_a)$ is a basis of $V$, we set $\|v\|=f(v^l,...,v^m)=f(v^a)$ and suppose that $f$ is differentiable at least of class $C^3$ for $v\neq 0$. The set $\{T| T\in GL(m,R),\ \|Tv\|=\|v\|, v\in V\}$ is a Lie group. Let $H$ be a subgroup of it.

{\bf Definition 3.1.} A vector bundle $\xi=(E,p,M)$ is said to be a $\{V,H\}-$bun\-dle if there exists a bundle atlas $\{(p^{-l}(U_\a),\phi_\a,V)\}$ such that the mappings $\psi_{\b,x}\circ\psi^{-1}_{\a,x}$ belongs to $H$ for every $x\in U_\a\cap U_\b\neq\psi$. We also say that $\xi$ admits an $H-$structure.

{\bf Proposition 3.1.} {\it If $\xi$ is a $\{V,H\}-$bundle, then its local fibers are Minkowski spaces isomorphic and isometric each to others.}

\noindent{\it Proof.} If $u\in E_x$, we set $\|u\|=f(\psi_{\a,x}(u))$ and obtain a norm on $E_x$ which does not depend on $\psi_{\a,x}$ because $\xi$ admits an $H-$structure. Namely $\psi_{\a,x}$ is also an isometry of $E_x$ and $V$ for every $x\in M$. Therefore the local fibers are isomorphic and isometric each to others.\hfill Q. E. D.

{\bf Examples} a) If $V$ is a Euclidian space then $O(m)$ leaves invariant its norm. Then $\xi$ is a $\{V,O(m)\}-$bun\-dle if and only if it is a Riemannian bundle.

b) If $\xi=(E=TM,\tau,M)$ and $M$ is modeled by $V$, then $\xi$ is a $\{V,H\}-$bun\-dle if and only if $M$ is a $\{V,H\}-$manifold in Ichijyo's sense [2].

If $\{(p^{-1}(U_\a),\phi_\a,V)\}$ is any bundle atlas on $\xi$, the cross-section $s_{\a,a}(x)=\phi^{-1}(x,e_a)$ define a frame in $E_x$ and the fiber coordinates $(y^a)$ are introduced by the equality $u_x=y^a s_{\a,a}(x)$. Setting $\sigma_{\a,a}(x)=\psi^{-1}_\a(x,e_a)$ we obtain a new frame in $E_x$, so that $u_x=u^a\sigma_{\a,a}(x)$. Taking $\sigma_{\a,a}(x)=\lbb^b_a(x)s_{\a,b}(x)$, it follows $y^a=\lbb^a_b(x)u^b$ or $u^a=\mu^a_b(x)y^b$, where $(\mu^a_b)$ is the inverse of the matrix $(\lbb^a_b)$. Now, $\|u_x\|=f(\psi_{\a,x}(u))=f(\psi_{a,x}(u^a\sigma_{\a,a}(x))=f(u^a)=f(\mu_b^a(x)y^b)$. Now we have a function $F:E\to R_+$, given locally by $F(x,y)=f(\mu_b^a(x)y^b)$, which is $(1)-$homogeneous and differentiable at least of class $C^3$ for $y\neq 0$. Moreover, as $F$ is provided by a norm, the matrix $(h_{ab}(x,y))=\(\dfrac12 \pp_a\pp_b(F^2)\)$ is nonsingular, and the quadratic form $h_{ab}\eta^a\eta^b$ is positive defined (see [9] p. 21 for a proof). We say that $F$ is a {\it fundamental Finsler function} on $E$. Therefore we have proved

{\bf Theorem 3.1.} {\it If a vector bundle $(E,p,M)$ admits an $H-$structure, then there exists on $E$ a fundamental Finsler function of the form $F(x,y)=f(\mu_b^a(x)y^b)$.}

{\bf Definition 3.2.} A linear connection $\nabla$ on a $\{V,H\}-$bundle is said to be an $H-${\it connection} if its parallel displacement preserves the Minkowski norms of fibers.

Let us set $\nabla_{\pp_k}s_a=\G^b_{ak}(x)s_b$. We have

{\bf Theorem 3.2.} {\it If $\nabla$ is an $H-$connection on the $\{V,H\}-$bundle $\xi$, then $\ovcirc{\de}_k F=0$, where $\ovcirc{\de}_k=\pp_k-\G^a_{bk}(x)y^b\pp_a$.}

\noindent{\it Proof.} Let $C=\{x(t),\ t\in [0,1]\}$ be a curve on $M$ and $S(x(t))=S^a(x(t))s_a(x(t))$ a cross-section of $\xi$ along $C$. It is parallel along $C$ with respect to $\nabla$ if and only if $\nabla_{\dot{x}(t)})S=0$, i.e. $\dfrac{dS^a}{dt}+\G^a_{bk}(x)S^b\dfrac{dx^k}{dt}=0$. If $\nabla$ is an $H-$connection, then $\dfrac{d\|S(t)\|}{dt}=0$. But $\|S(t)\|= F(x(t),S^a(x(t)))$ so we obtain $0=\pp_k F\dfrac{dx^k}{dt}+\pp_a F\dfrac{dS^a}{dt}=(\pp_kF-\G^a_{bk}(x)S^b)\dfrac{dx^k}{dt}=\ovcirc{\de}_k F\dfrac{dx^k}{dt}$. Since $C$ is arbitrary, it results $\ovcirc{\de}_kF=0$.\hfill Q. E. D.

\section{Einstein equations and the conservation laws}\setcounter{equation}{0} If $\G^a_{bk}(x)$ are the coefficients of a linear connection $\nabla$ on $\xi$, then $\ovcirc{N}^a_k(x,y)=\G^a_{bk}(x)y^b$ define a nonlinear connection on $\xi$. We consider a ``deformation'' of this nonlinear connection, i.e. $N^a_k(x,y)=\ovcirc{N}^a_k(x,y)+A^a_k(x)$, where $A^a_k(x)$ defines a $d-$tensor field on $E$, depending only on $x$. Using Theorem 3.2 we easily obtain

{\bf Proposition 4.1.} {\it If $\nabla$ given by $\G^a_{bk}(x)$ is an $H-$connection on the $\{V,H\}-$bundle $\xi$ then, $\de_kF=\pp_kF-N^a_k\pp_aF=0$ if and only if
\begin{equation}
A^a_k(x)\pp_a F=0,
\end{equation}
holds good.}

In what follows we assume that $\nabla$ is an $H-$connection on $\xi$ and that $A^a_k(x)$ satisfies (4.1). Now, let $g_{ij}(x)$ be a metric on $M$ and let $\G^i_{jk}(x)$ be the corresponding Christoffel symbols. The definition of $F$ shows that in the adapted frames to the $H-$structure on $\xi$, the functions $h_{ab}$ depend only on $y$. The following natural metrical structure on $E$ can be considered:
\begin{equation}
{\cal G}(x,y)=g_{ij}(x)dx^idx^j+h_{ab}(y)\de y^a\de y^b.
\end{equation}
This is an example of Riemann-Minkowski metric on $E$. A study of the Riemann-Minkowski metrics on $TM$ is given in [8]. Let $C^a_{bc}(y)$ be the Christoffel's symbols associated with $h_{ab}(y)$. Then it is clear that $E\G=(\G^i_{jk}$, $\G^a_{bk}(x)$, $0$, $C^a_{bc}(y))$ is a $d-$connection on $E$. Moreover, we have

{\bf Theorem 4.1.} {\it The $d-$connection $E\G$ is metrical with respect to the metric ${\cal G}$.}

\noindent{\it Proof.} The first three equalities from (2.7) hold by virtue of the definition of $E\G$. To prove the last one we remark that $\de_k F=0$ is the same with $F_{|k}=0$ and we note that $|k\circ\pp_a=\pp_a\circ |k$. So, $h_{ab|k}=\pp_a\pp_b(F_{|k})=0$.\hfill Q.E.D.

An easy computation shows that
\begin{equation}
R^a_{\ jk}=R^{\ a}_{b\ jk}(x)y^b+(\pp_k A^a_j+\G^a_{bk}A^b_j-j/k),
\end{equation}
where $R_{b\ jk}^{\ a}(x)$ is the curvature of $\nabla$. The others torsions of $E\G$ are vanishing. The $d-$tensor field $A^a_k(x)$ can be viewed as defining an 1-form $A$ on $M$, valued in $\xi$. If $\ovcirc{\nabla}$ denotes the Levi-Civita connection of $g_{ij}$ then the covariant differential of $A$ is the 2-form $$(\wt{\nabla}A)(X,Y)=\nabla_Y A(X)-A(\ovcirc{\nabla}_Y X),\ {\rm{for\ }} X, Y\in{\cal X}(M).$$

{\bf Theorem 4.2.} {\it The metrical $d-$connection $E\G$ coincide with the Levi-Civita connection of ${\cal G}$ if and only if a) $\nabla$ is flat, and b) $\wt{\nabla}A$ is symmetric.}

\noindent{\it Proof.} By annihilating $R^a_{\ jk}$ from (4.3) we obtain that $\nabla$ is flat, and $\wt{\nabla}A$ is symmetric. The converse is clear.\hfill Q. E. D.

Particularizing (2.5) one obtains

{\bf Proposition 4.2.} {\it The curvatures of $E\G$ are as follows: $R^{\ i}_{j\ kh}$ is the curvature of $\overset{\circ}{\nabla}$, $\wt{R}^{\ a}_{b\ kh}=R^{\ a}_{b\ kh}+C^a_{bc}R^c_{\ kh}$, $\overset{1}{P}^{\ a}_{b\ kc}=-C^a_{bc|k}=0$, $\overset{2}{P}^{\ i}_{j\ kc}=0$, $S^{\ a}_{b\ cd}$ has the general form.}

\noindent Note that $C^a_{bc|k}=0$ results from $C^a_{bc}=h^{ad}\pp_b\pp_c F^2/4$ and $F_{|k}=0$.

{\bf Corollary 4.1.} {\it The metrical $d-$connection $E\G$ has no curvature if a) $\ovcirc{\nabla}$ is flat, b) $\nabla$ is flat, c) $\wt{\nabla}A$ is symmetric, d) $S^{\ a}_{b\ cd}=0$.}

As to the metrical $d-$connection $E$ we are associated with the Einstein equation
\begin{equation}
{\rm{Ric}}(E\G)-{\cal RG}=\varkappa {\cal T},
\end{equation}
where Ric$(E\G)$ and ${\cal R}$ denote the Ricci tensor and the scalar curvature of $E\G$, respectively; $\varkappa$ is a constant; ${\cal T}$ is the energy momentum tensor. With respect to the adapted frame $(\de_k,\pp_a)$, equation (4.4) decomposes as follows (cf. [1]):
\begin{equation}
\left\{
\begin{array}{l}
R_{ij}-\dfrac12 (R+S)g_{ij}=\varkappa{\cal T}_{ij},\ 0={\cal T}_{ai},\ 0={\cal T}_{ia}\\ \\
S_{ab}-\dfrac12(S+R)h_{ab}=\varkappa{\cal T}_{ab},
\end{array}
\right.
\end{equation}
where $R_{ij}=R^{\ k}_{i\ jk}$, $S_{ab}=S^{\ c}_{a\ bc}$, $R=g^{ij}R_{ij}$, $S=h^{ab}S_{ab}$; in the right members appear the components of ${\cal T}$, two of them must be taken zero because the curvatures $\overset{1}P$ and $\overset{2}P$ of $E\G$ are vanishing.

Equation (4.5) will be called the {\it Einstein equations} on $E$. The conservation law is obtained by annihilating the divergence of the tensor which appear as the first member of (4.4), called the Einstein tensor. In the adapted frame one obtains as conservation laws:
\begin{equation}
\left\{
\begin{array}{l}
\(R^i_j-\dfrac12(R+S)\de^i_j\)_{|i}=0,\\ \\
(S^a_b-(R+S)\de^a_b)|_a=0.
\end{array}
\right.
\end{equation}
As is well known the divergence of the Einstein tensor associated with the Levi-Civita connection identically vanishes. By using Theorem 4.2 one obtains

{\bf Theorem 4.3.} {\it The conservation laws on $E$ with respect to $E\G$ are identities if: i) $\nabla$ is flat and ii) $\wt{\nabla}A$ is symmetric.}

The conditions ii) and (4.1) on $A$ can be easily satisfied taking, for instance, $A=0$. The condition  i) is a strong one because if $M$ is simply connected, then $E=M\times V$. Examining (4.6) we shall find algebraic conditions on $A$ under which the conservation laws are identities. The second equality (4.6) reduces to $\(S^a_b—\dfrac12 S\de^a_b\)|_a=0$ which is an identity by virtue of Bianchi identities. The first is reduced to $\(R^i_j—\dfrac R\de^i_j\)—S_{|i}=0$, and by virtue of the Bianchi identities it become an identity if and only if $S_{|i}=0$. Namely we have proved

{\bf Theorem 4.4.} {\it The conservation laws on $E$ with respect to $E\G$ are identities if
\begin{equation}
(\G^a_{bk}y^b+A^a_k)\pp_aS=0,
\end{equation}
holds good.}

The conditions (4.1) and (4.7) form together an algebraic system of $2n$ equations with $nm$ unknowns $A^a_k$. If $m=1$ the first $n$ equations give $A^1_k=0$, $k=1,...,n$ and the last $n$ equations are verified if $S$ is constant or if $\G^{\ 1}_{1k}=0$. For $m=2$ the determinant of the system is $-(\pp_1 F\pp_2 S-\pp_2 F\pp_1 S)^n$ which generally is different from zero. For $m>2$ some unknowns can be arbitrarily taken.

To conclude we must say that the total space $E$ of a $\{V,H\}-$bundle, whose base is a Lorentz manifold, $(M,g_{ij})$ can be endowed with a metrical $d-$connection with torsion for which the conservation laws are verified. We think this is a basic facts to an unitary theory of Finslerian type.

{\it Acknowledgements.} The author expresses many thanks to Professor Dr. Radu Miron for numerous discussions on the subject of this paper. He is also indebted to Professor Dr. Y. Ichijyo for some comments on \S3.

{\bf REFERENCES}

\begin{enumerate}
  \item M. Anastasiei, {\it Vector bundles. Einstein equations}, An. St. Univ. ``Al. I. Cuza'' Iasi, s. I-a Mat., 32(1986), 1-8.
  \item Y. Ichijyo, {\it Finsler manifolds modeled on a Minkowski space}, J. Math., Kyoto Univ., 16 (1976), 639-652.
  \item A. Kawaguchi, {\it On the theory of non-linear connections I. Introduction to the theory of general non-linear connection}, Tensor, N. S., 2 (1952), 123-142.
  \item A. Kawaguchi, {\it On the theory of non-linear connection II. Theory of Min\-kowski spaces and of nonlinear connections in a Finsler space}, Tensor, N. S., 6 (1956), 165-199.
  \item R. Miron, {\it A Lagrangian theory of relativity}, Preprint, Sem. Geom. Top., Timiqoara, R. S. Romania, (1985), 1-53.
  \item R. Miron, {\it Techniques of Finsler geometry in the theory of vector bundles}, Acta Sci. Math., 49 (1985), 119-129.
  \item R. Miron and M. Anastasiei, {\it Vector bundles. Lagrange geometry. Applications in Relativity} (in Romanian), Editura Academiei R. S. Romania, Bucuresti, 1987.
  \item R. Miron, S. Watanabe and S. Ikeda, {\it Some connections on tangent bundle and their applications to the general relativity}, Tensor, N. S., 46 (1987), 8-22.
  \item H. Rund, {\it The differential geometry of Finsler spaces}, Springer, 1957.
  \item M.-A. Tonnelat, {\it Les th\'{e}ories unitaries de l'\'{e}lectromagn\'{e}tisme et de la gravitation}, Gauthier-Villars, Paris, 1965.
\end{enumerate}

\bigskip

\hspace*{8.1cm}{\footnotesize{{\it Faculty of Mathematics}}}\\ \hspace*{8cm}{\footnotesize{\it ``Al.~I.~Cuza'' University of Iasi}}

\newpage





\def\G{\Gamma}

\runningauthor={M. ANASTASIEI}
\runningtitle={THE GEOMETRY OF TIME-DEPENDENT LAGRANGIANS}
\noindent
\baselineskip 8pt
\noindent{\footnotesize{Pergamon Mathl. Comput. Modelling\hfill\break
Vol. 20, No. 4/5, pp. 67-81, 1994}}
\vskip 2cm
\baselineskip 11.5pt plus .15pt
\centerline{\bf\Large THE GEOMETRY OF}
\vskip .2cm
\centerline{\bf\Large TIME-DEPENDENT LAGRANGIANS}
\vskip .5cm
\centerline{\bf \footnotesize{BY}}
\vskip .5cm
\centerline{\bf \footnotesize{M. ANASTASIEI}\footnote{I am very indebted to P. L. Antonelli for his interest in this subject. I would like to thank R. Miron for his comments and valuable suggestions during the preparation of this paper.}}
\vskip 1cm

\begin{abstract}
A generalization of Lagrange geometry appropriate for time-dependent Lagrangians arising in physics and biology, called rheonomic Lagrange geometry, is developed. Nonlinear and linear connections, their torsions, curvatures and deflections are explicitly given. Almost contact structures in rheonomic Lagrange spaces are characterized. Maxwell's equations, for a given Lagrangian, determined by the deflection tensors, are derived.
\end{abstract}

\setcounter{section}{0}
\section{Introduction}

Variational principles are basic for most mathematical models in mechanics, physics, ecology, physiology, and so on. These involve Lagrangians or Hamiltonians from which the Euler-Lagrange or Hamilton equations are derived, the theory being then centered on the later. From a geometrical point of view, the most general framework for such a theory is provided by differentiable (smooth) fibre bundles. It means that, for instance, a Lagrangian is a smooth real valued function on the total space $TM$ of the tangent bundle $(TM,\tau,M)$ over a smooth manifold $M$. For a geometrization of such Lagrangians, we refer to [1-5].

There exist certain mathematical models, as for instance those for the three-body problem [6, p. 206] and those concerning ecological systems due to Antonelli [7,8]in which an explicit dependence on time of the Lagrangian (Hamiltonian) is required. A time-dependent Lagrangian is smooth and real valued on $\R\times TM$, where $\R$ is the field of real numbers.

It is our aim to present a geometrization of time-dependent Lagrangians using as a pattern the geometry of Lagrange spaces developed by Miron [3-5,9]. The reader is invited to compare this geometrization to those of [10,11].

We begin with some facts (almost tangent structures, nonlinear connections) from the geometry of manifold $\R\times TM$ fibered over $\R\times M$ by $\pi(t,v)=(t,\tau(v))$, $t\in\R$, $v\in TM$. Then we associate with any nonlinear connection $N$ on $E=\R\times TM$ a semispray on $E$ whose integral curves coincide with the paths of $N$. Regular time-dependent Lagrangians are introduced in Section 3. It is shown that any such Lagrangian $L$ induces a canonical nonlinear connection $N_L$ on $E$. This nonlinear connection $N_L$ is derived from the Euler-Lagrange equations resulting from a variational problem involving $L$. Then a metrical almost contact structure on $E$ depending on $L$ only is exhibited. The geometry of $L$ is based on this structure and may be thought of as the counterpart of the almost K\"{a}hlerian model used in the geometry of time {\it independent} Lagrangians [3]. As a first step, the linear connections, which are compatible with $N_L$ as well as with the almost contact structure on $E$, called $N-$linear connections, are studied. The metrical $N-$linear connections are studied too. The existence of a canonical one is shown. Finally, some remarkable time-dependent Lagrangians are considered. Thus, we investigate homogeneous time-dependent Lagrangians.

Similarities with Finsler geometry are emphasized. A second class of time-dependent Lagrangians which we consider contains Lagrangians used in electrodynamics. It is shown that their geometry supports a theory of electromagnetism based on $N-$linear connections.

\subsection{On the Geometry of $\R\times TM$} Let $M$ be a smooth manifold of dimension $n$. It will be assumed Hausdorff connected, and paracompact. We assume $\R\times M$ is coordinated by $(t,x^i)\equiv (t,x)$. The indices $i,j,k,$..., will run over $1,2,..., n$, and the Einstein summation convention will be used. The coordinates in the fibres of the submersion $\pi:\R\times TM\to\R\times M$ are $(y^i)\equiv (y)$, introduced by $u_{(t,x)}=(t,v_x)=\(t,y^i\(\dfrac{\pp}{\pp x^i}\)_x\)\in\pi^{-1}(t,x)$, with $\(\dfrac{\pp}{\pp x^i}\)_x$ the natural basis in the tangent space $T_xM$, in $x\in M$. Thus, the manifold $\R\times TM$ is coordinated by $(t,x^i,y^i)\equiv(t,x,y)$ and $\pi$ takes the form $(t,x,y)\to (t,x)$. A change of local coordinates $(t,x,y)\to(\wt{t},\wt{x},\wt{y})$ on $E=\R\times TM$ has the following form
$$\wt{t}=t, \wt{x}=\wt{x}^i(x^1,...,x^n),\ \wt{y}^i=\dfrac{\pp\wt{x}^i}{\pp x^k}y^k,\leqno(1.1)$$ with rank $\(\dfrac{\pp\wt{x}^i}{\pp x^k}\)=n$.

The natural basis $\(\dfrac{\pp}{\pp t},\dfrac{\pp}{\pp x^i},\dfrac{\pp}{\pp y^i}\)$ transforms under (1.1) as follows: $$\begin{array}{l}\dfrac{\pp}{\pp t}=\dfrac{\pp}{\pp\wt{t}},\\ \\
\dfrac{\pp}{\pp x^j}=\dfrac{\pp\wt{x}^i}{\pp x^j}\dfrac{\pp}{\pp\wt{x}^i}+\dfrac{\pp \wt{y}^i}{\pp x^j}\dfrac{\pp}{\pp\wt{y}^i},\\ \\ \dfrac{\pp}{\pp y^j}=\dfrac{\pp\wt{y}^i}{\pp y^j}\dfrac{\pp}{\pp\wt{y}^i}.\end{array}\leqno(2.1)$$ The kernel of the Jacobian map $D\pi$ supplies a distribution $u\to V_u E$, $u\in E$, on $E$ which will be called the vertical distribution on $E$. A local basis of the vertical distribution is given by the local vector fields $\(\dfrac{\pp}{\pp y^i}\)$ denoted in what is to follow as $(\dot\pp_i)$. From (1.1) and (1.2), it follows that $C=y^i\dot\pp_i$ is a global vector field on $E$. This may be used in order to express the homogeneity with respect to $(y^i)$ of various geometrical objects on $E$. With the help of (1.2), one may check that setting $$J\(\dfrac{\pp}{\pp t}\)=0,\ \ J(\pp_i)=\dot{\pp}_i, J(\dot{\pp}_i)=0,\leqno(1.3)$$ where $\pp_i$ stands for $\dfrac{\pp}{\pp x^i}$ and requiring the linearity of $J$ one obtains a well-defined $(1,1)-$tensor field on $E$. Moreover, we have $J^2=0$, and the Nijenhuis tensor field $N_J(X,Y)=[JX,JY]+J^2[X,Y]-J[JX,Y]-J[X,JY]$, $X,Y\in \chi(E)$, the module of vector fields on $E$, identically vanishes. Thus $J$ defines an almost tangent structure on $E$. Sometimes it is convenient to put $t=x^0$ and to use the Greek indices $\a,\b,\g$,..., ranging over $0,1,2,...,n$.

A {\it nonlinear connection} on $E$ is a distribution, called horizontal, $u\to H_u E$, $u\in E$, which is supplementary to the vertical distribution on $E$. Such a distribution can be given by $(n+1)$ local vector fields, say $\delta_\a$. Choosing $\delta_\a$ such that they are projected by $D\pi$ to $\dfrac{\pp}{\pp x^\a}$ one gets $$\de_\a=\pp_\a-N^i_\a(t,x,y)\dot{\pp}_i,\leqno(1.4)$$ where $\pp_\a$ stands for $\dfrac{\pp}{\pp x^\a}$ and the minus sign is taken for convenience.

The invariance under (1.1) of the horizontal subspaces requires the condition $$\de_\a=\dfrac{\pp\wt{x}^\b}{\pp x^\a}\wt{\de}_\b.\leqno(1.5)$$ In turn, equation (1.5) implies the following law of transformation for the coefficients $N^i_\a$: $$\wt{N}^i_\a\dfrac{\pp\wt{x}^\a}{\pp x^\b}=\dfrac{\pp\wt{x}^i}{\pp x^k}N^k_\b-\dfrac{\pp^2\wt{x}^i}{\pp x^\b \pp x^k}y^k.\leqno(1.6)$$ If one rewrites (1.4) in the form $$\de_0=\dfrac{\pp}{\pp t}-N^i_0(t,x,y)\dot{\pp}_i,\ \de_i=\pp_i-N^k_i(t,x,y)\dot{\pp}_k,\leqno(1.7)$$ one may state the following theorem.

{\bf Theorem 1.1.} {\it To give a nonlinear connection on $E$ is equivalent to giving a set of functions $(N^k_0,N^k_i)$ defined in each coordinate chart on $E$, which transform under $(1.1)$ as follows: $$\wt{N}^k_0(\wt{t},\wt{x},\wt{y})=\dfrac{\pp\wt{x}^k}{\pp x^h}N^h_0(t,x,y),\leqno(1.8)$$ $$\wt{N}^k_i(\wt{t},\wt{x},\wt{y})\dfrac{\pp\wt{x}^i}{\pp x^j}=\dfrac{\pp\wt{x}^k}{\pp x^i}N^i_j(t,x,y)-\dfrac{\pp\wt{y}^y}{\pp x^j}.\leqno(1.9)$$}

\ms

\noindent{\it Proof.} If a nonlinear connection on $E$ is given by the local coefficients $(N^i_\a)$ satisfying (1.6), taking into account (1.5) and (1.7), it comes out that equation (1.6) is equivalent to (1.8) and (1.9). Conversely, a set of functions defined in each coordinate chart on $E$ verifying (1.6) on overlaps, provides, according to (1.4) and (1.5), a nonlinear connection on $E$.\hfill$\blacksquare$

The local coefficients $(N^i_0(t,x,y))$ transform under (1.1) like the components of a vector field on $M$, although they depend on $t,x$ and $y$. We shall say $(N^i_0)$ define a distinguished vector field on $E$, briefly a $d-$vector field. More generally, an $(r,s)-$tensor field on $E$ whose local components transform like those of an $(r,s)-$tensor field on $M$, ignoring their dependence on $t,x$, and $y$, will be called a $d-$tensor field of type $(r,s)$. A similar situation appears in [12], where a $d$-tensor field is called a Finsler tensor field, as well as in [7, p. 131], where a $d-$tensor field is called a Douglas tensor field.

The local coefficients $(N^k_i(t,x,y))$ transform under (1.1) like those of a nonlinear connection on $TM$ [3]. When these local coefficients do not depend on $t$, they really define a nonlinear connection on $TM$. Conversely, a nonlinear connection on $TM$ paired with a $d-$vector field on $E$ defines a nonlinear connection on $E$.

The decomposition $T_uE=H_uE\oplus V_uE$ gives rise to two projectors, an horizontal one denoted by $h$ and a vertical one denoted by $v$, as well as to an almost product structure $P=h-v$. All these depend smoothly on $u\in E$ and thus induce $(1,1)-$tensor fields on $E$, which will be denoted again by $h,v$, and $P$, respectively.

There exist many ways for introducing the curvature of a nonlinear connection. We choose the following formal one since it allows us to relate quickly the curvature to the integrability of the horizontal distribution. Namely, the {\it curvature} $\Omega$ of a nonlinear connection is defined as the Nijenhuis tensor field $N_h$ of the horizontal projector $h$, that is $\Omega=N_h$. In a coordinate chart $\Omega$, it is given as follows: $$\Omega(\de_\a,\de_\beta)=R^i_{\ \a\beta}\dot{\pp}_i,\ \Omega(\pp_\a,\dot{\pp}_i)=0,\ \Omega(\dot{\pp}_i,\dot{\pp}_j)=0,\leqno(1.10)$$ $$R^i_{\a\beta}=\de_\beta N^i_\a-\de_\a N^i_\beta=\pp_\beta N^i_\a-\pp_\a N^i_\beta+N^k_\a\dot{\pp}_k N^i_\beta-N^k_\beta\dot{\pp}_k N^i_\a.\leqno(1.11)$$ On the other hand, we have $$[\de_\a,\de_\beta]=R^i_{\a\beta}\dot{\pp}_i,\ [\de_\a,\dot{\pp}_i]=\dot{\pp}_i N^i_\a \dot{\pp}_j.\leqno(1.12)$$ Thus, the horizontal distribution on $E$ is integrable if and only if $\Omega=0$, or equivalently, $R^i_{\a\beta}=0$. We notice that $R^i_{jk}$ and $R^i_{ok}$ define $d-$tensor fields on $E$ of type $(1,2)$ and $(1,1)$, respectively.

Let $B^k_{j\a}=\dot{\pp}_j N^k_{\a}$. Differentiating with respect to $y^j$ both sides (1.8), one finds that $B^k_{\ j0}$ defines a $d-$tensor field of type $(1,1)$. Proceeding similarly with (1.9), one gets $$\wt{B}^k_{\ rs}\dfrac{\pp\wt{x}^r}{\pp x^j}\dfrac{\pp\wt{x}^s}{\pp x^i}=\dfrac{\pp\wt{x}^k}{\pp x^h}B^h_{\ ji}-\dfrac{\pp^2\wt{x}^k}{\pp x^i\pp x^j}.\leqno(1.13)$$ Thus, the functions $B^k_{\ ji}(t,x,y)$ transforms under (1.1) as the local coefficients of a classical linear connection, although they depend on $t,x,y$. It is said in [7, p. 131] that $B^h_{\ ji}$ define a Douglas connection. We have used the letter $B$ since in [12] these functions are related to the so-called Berwald connection.

A nonlinear connection $N(N^i_\a)$ is homogeneous (resp. linear) if the functions $N^i_0(t,x,y)$ and $N^i_k(t,x,y)$ are homogeneous of degree one (resp. linear) with respect to $(y^i)$. Of course, in order to speak about homogeneous connections we must delete from $E$ the points $(t,x,0)$ because any homogeneous real function of class $C^1$ at the origin becomes linear.

When a linear connection $(N^i_0,N^i_j)$ is given, the equalities $N^i_0(t,x,y)=K^i_j(t,x)y^j$, $N^i_j(t,x,y)=\Gamma^i_{jk}(t,x)y^k$ provide a pair $(K^i_j(t,x),\Gamma^i_{jk}(t,x))$ which may be thought of as a general affine connection on $\R\times M$ in the sense of [13].

\section{Semisprays and nonlinear connections}
A {\it time-dependent vector field} on $TM$ is a smooth map $X^0:\R\times TM\to T(TM)$, $(t,u)\to X^0(t,u)\in T_u(TM)$, $u\in TM$. It induces a vector field $X$ on $\R\times TM$ by setting $X(t,u)=(1,X^0(t,u))$, and we have also $X=\dfrac{\pp}{\pp t}+X^0$ [2].

A {\it time-dependent semispray} (second order differential equation on $M$) is a time-dependent vector field $S^0$ on $TM$ which satisfies $$D\tau\circ S^0(u)=u,\ \ \textrm{ for all }u\in TM.\leqno(2.1)$$ The vector field $S$ induced on $\R\times TM$ by a time-dependent semispray $S^0$, that is $$S=\dfrac{\pp}{\pp t}+S^0,\leqno(2.2)$$ will be called a {\it semispray}.

According to (2.1) and (2.2) a semispray in a coordinate chart on $E$ takes the form $$S=\dfrac{\pp}{\pp t}+y^i\pp_i+S^i(t,x,y)\dot{\pp}_i,\leqno(2.3)$$ with $(S^i)$ verifying on overlaps $$\wt{S}^i=\dfrac{\pp\wt{x}^i}{\pp x^k}S^k+\dfrac{\pp^2\wt{x}^i}{\pp x^j\pp x^k}y^j y^k.\leqno(2.4)$$ Conversely, a vector field $S$ which in each coordinate chart has the form (2.3) such that equation (2.4) is fulfilled, is a semispray.

A direct calculation gives the following result.

{\bf Proposition 2.1.} {\it A vector field $S$ on $E$ is a semispray if and only if $dt(S)=1$, $\psi^i(S)=0$, with $\psi^i=dx^i-y^i dt$.}

A relationship between semisprays and nonlinear connections is given by the following two theorems.

{\bf Theorem 2.1.} {\it Let $N$ be a nonlinear connection given by the local coefficients $(N^i_0(t,x,y),N^i_k(t,x,y))$. Then $S=\dfrac{\pp}{\pp t}+y^i\pp_i-(N^i_0+N^i_k y^k)\dot{\pp}_i$ is a semispray.}

{\bf Theorem 2.2.} {\it Let $S(S^i)$ be a semispray. Then $\(\dfrac{\pp S^i}{\pp t},-\dfrac{1}{2}\dfrac{\pp S^i}{\pp y^j}\)$ are local coefficients for a nonlinear connection on $E$.}

The proofs follow by showing that equations (1.8) and (1.9) imply (2.3), and conversely.

A time-dependent semispray is said to be a spray if it is invariant under (a gauge transformation, dilatation or contraction) {\it similarity} on $TM$ and a semispray will be called a spray if it is provided by a time dependent spray $S^0$. It is immediate that a semispray is a spray if and only if the functions $(S^i(t,x,y))$ are homogeneous of degree 2 in $(y^i)$. The later condition is clearly compatible with (2.4).

If $S(S^i)$ is a spray, then $(0,-\dfrac{1}{2}\dot{\pp}_j S^i)$ are the local coefficients of a homogeneous connection. Conversely, a homogeneous connection defines a spray $S^i=-N^i_k y^k$.

Let $c:\R\to M$ be a smooth curve on $M$ and $\dot{c}:\R\to TM$ its tangent vector field. Then $\sigma(t)=(t,\dot{c}(t))$ defines a smooth curve on $\R\times TM$. We say this curve is an integral curve of a semispray $S$ if $$\dot{\sigma}(t)=S(\sigma(t)),\ \ t\in\R.\leqno(2.5)$$ If we assume that $c(t)$ belongs to a coordinate chart for all $t\in\R$, and we take $x^i=x^i(t)$, $t\in\R$, as the equations of the curve $c$, then equation (2.5) is equivalent to $$\dfrac{d^2 x^i}{dt^2}=S^i(t,x,\dot{x}),\ \ \dot{x}=\dfrac{dx}{dt},\leqno(2.6)$$ because of $\dot{\sigma}(t)=\dfrac{\pp}{\pp t}+\dfrac{dx^i}{dt}\dot{\pp}_i+\dfrac{d^2 x^i}{dt^2}\dot{\pp}_i.$

A curve $c:t\to c(t)$, $t\in\R$, on $M$ is said to be a path for a nonlinear connection $N$ if the curve $\sigma:t\to(t,\dot{c}(t))$ is horizontal with respect to $N$, that is, its tangent vector field belongs to the horizontal distribution on $E$.

In a coordinate chart containing $\sigma(t)$, $t\in\R$, if $(\de_\a,\dot{\pp}_i)$ is the adapted basis introduced before and $(dx^\a,\de y^i)$, with $\de y^i=dy^i+N^i_\a dx^\a$ is its dual, it appears as obvious that $\sigma$ is an horizontal curve if and only if $\de y^i(\dot{\sigma})=0$ for every $i=1,..., n$. Writing down these equations, one obtains the following theorem.

{\bf Theorem 2.3.} {\it A curve $c:t\to x^i(t)$, $t\in\R$, on $M$ is a path for a nonlinear connection $(N^i_0,N^i_j)$ if and only if $$\dfrac{d^2 x^i}{dt}+N^i_j(t,x,\dot{x})\dfrac{dx^j}{dt}+N^i_0(t,x,\dot{x})=0,\ \ x^i=\dfrac{dx^i}{dt}.\leqno(2.7)$$}

Looking at the semispray associated with a nonlinear connection one immediately gets the next result.

{\bf Theorem 2.4.} {\it The paths of a nonlinear connection coincide with the integral curves of the semispray associated with it.}

We notice that the systems of differential equations (2.6) and (2.7) do not remain in the same form if an arbitrary change of parameter is performed. They keep their form only if one  sets $\tilde{t}=\pm t+a$, with $a\in\R$. Thus, $t$ plays the role of an affine parameter. We conclude that the solutions of these systems have to be considered together with the parameters in which they are given. In other words, the curve $c$ in the above has to be thought of as a parameterized curve.

\section{Time-dependent lagrangians}
Now, we shall point out that a regular time-dependent Lagrangian defines a nonlinear connection on $E=\R\times TM$ and, thus, a semispray on $E$.

A smooth function $L:\R\times TM\to\R$, $(t,v)\to L(t,v)$, is called a {\it time-dependent} Lagrangian on $M$. It is said $L$ is {\it regular} if the matrix with the entries $$g_{ij}(t,x,y)=\dfrac12\dfrac{\pp^2 L}{\pp y^i\pp y^j},\leqno(3.1)$$ is of rank $n$ on $E$.

The condition for $L$ regular does not depend on the coordinate chart involved.

{\bf Definition 3.1.} A pair $RL^n=(M,L(t,x,y))$ in which $L$ is a regular time-dependent Lagrangian such that the quadratic form with the coefficients $g_{ij}$ from (3.1) has constant signature, will be called a {\it rheonomic Lagrange space}.

Let $c:t\to c(t)$, $t\in\R$ be a parameterized curve on $M$ as before. If its image is in a coordinate chart, one may take $x^i=x^i(t)$, $t\in\R$ as its local representation, and then its tangent vector field $\dot{c}$ is locally represented as $(x^i(t),\dot{x}^i(t))$. When a regular time dependent Lagrangian $L$ on $M$ is given, one may define a functional $${\cal L}:c\to{\cal L}(c)=\dd\int^{t_1}_{t_0}L(t,x(t),\dot{x}(t))dt,$$ which suggests the following variational problem: find those curves, called extremals, which afford extremal values for ${\cal L}$. Looking for such an extremal in the space of all curves with fixed end points, one finds [1, p. 153; 8, p. 58], that it is among curves which are solutions of the Euler-Lagrange equations $$\dfrac{d}{dt}\(\dfrac{\pp L}{\pp \dot{x}^i}\)-\dfrac{\pp L}{\pp x^i}=0.\leqno(3.2)$$ Now, if one considers the curve $\wt{c}=(t,c(t))$, $t_0\leq t\leq t_1$, on ${\R}\times M$, it comes out [8, p. 58] that $\wt{c}$ is an extremal of the  functional $${\cal L}(\wt{c})=\dd\int^{t_1}_{t_0}L(t,x(t),\dot{x}(t))dt,$$ on the space of curves joining $(t_0,x_0)$ and $(t_1,x_1)$  if the Euler-Lagrange equations are satisfied along $\wt{c}$.

Expanding the derivative with respect to $t$, the equations (3.2) can be put in the form $$2g_{ij}\dfrac{d^2 x^j}{dt^2}+\(\dfrac{\pp^2 L}{\pp\dot{x}^i\pp x^j}\dot{x}^j-\dfrac{\pp L}{\pp x^i}\)+\dfrac{\pp^2 L}{\pp t\pp\dot{x}^i}=0.\leqno(3.3)$$ Using the inverse $(g^{ki})$ of the matrix $(g_{ij})$, one resolves (3.3) with respect to $\dfrac{d^2 x^k}{dt^2}$ as follows: $$\dfrac{d^2 x^k}{dt^2}+2G^k(t,x,\dot{x})+N^k_0(t,x,\dot{x})=0,\leqno(3.4)$$ in which the following notations were used $$N_0^k(t,x,y)=\dfrac12 g^{ki}\dfrac{\pp^2 L}{\pp t\pp y^i},\leqno(3.5)$$ $$G^k(t,x,y)=\dfrac14 g^{ki}\(\dfrac{\pp^2 L}{\pp y^i\pp x^j}y^j-\dfrac{\pp L}{\pp x^i}\),\ \ y^i=\dot{x}^i.\leqno(3.6)$$ Now, we state the following result:

{\bf Theorem 3.1.} {\it The functions $N_L=(N^k_0(t,x,y)$, $N^k_i(t,x,y))$, where $N^k_0$ is given by $(3.5)$ and $N^k_i$ by $$N^k_i(t,x,y)=\dfrac{\pp G^k(t,x,y)}{\pp y^i},\leqno(3.7)$$ are local coefficients of a nonlinear connection on $E$ completely determined by $L$.}

\ms

\noindent{\it Proof.} Under the coordinate transformation $(t,x,y)\to(\wt{t},\wt{x},\wt{y})$, given by (1.1), $\dfrac{\pp L}{\pp t}$ is invariant; $\dfrac{\pp L}{\pp y^i}$ transform like a covector on $M$, i.e., they define a $d-$covector field, and because $(g^{kj})$ transforms like the components of a $d-$tensor field of type $(2,0)$, it follows that $(N^k_0)$ from (3.5) are the components of a $d-$vector field on $E$. The partial derivatives of $L$ take, under (1.1), the following form:
$$\dfrac{\pp L}{\pp x^i}=\dfrac{\pp L}{\pp\wt{x}^k}\dfrac{\pp\wt{x}^k}{\pp x^i}+\dfrac{\pp L}{\pp\wt{y}^k}\dfrac{\pp^2\wt{x}^k}{\pp x^i\pp x^j}y^j,$$ $$\dfrac{\pp^2 L}{\pp y^i\pp x^k}=\dfrac{\pp^2 L}{\pp\wt{y}^j\pp\wt{y}^h}\dfrac{\pp\wt{x}^j}{\pp x^i}\dfrac{\pp\wt{x}^h}{\pp x^k}+\dfrac{\pp L}{\pp\wt{y}^j}\dfrac{\pp^2\wt{x}^j}{\pp x^i\pp x^k}+2\wt{g}_{jh}\dfrac{\pp\wt{x}^j}{\pp x^i}\dfrac{\pp^2\wt{x}^h}{\pp x^k\pp x^s}y^s.$$ Multiplying the second equality by $y^k$, we introduce the result in the form of (3.6): $$4g_{ij}G^j=\dfrac{\pp^2 L}{\pp y^i\pp x^k}y^k-\dfrac{\pp L}{\pp x^i},$$ which thus become $$4g_{ij}G^j=\(\dfrac{\pp^2 L}{\pp\wt{y}^j\pp\wt{x}^k}\wt{y}^k-\dfrac{\pp L}{\pp\wt{x}^j}\)\dfrac{\pp\wt{x}^j}{\pp x^i}+2g_{jk}\dfrac{\pp\wt{x}^j}{\pp x^i}\dfrac{\pp^2\wt{x}^k}{\pp x^r\pp x^s}y^r y^s,$$ since two terms cancel each other. Hence, we obtain the transformation law of $G^i$ as follows: $$\dfrac{\pp\wt{x}^k}{\pp x^i}G^i=\wt{G}^k+\dfrac12\dfrac{\pp^2\wt{x}^k}{\pp x^i\pp x^j}y^i y^j,$$ on account of rank$(g_{jj})=$rank$\(\dfrac{x^j}{x^i}\)=n$. Differentiating both sides of the last equality with respect to $y^j$, one gets $N^i_k$, from equation (3.7) has the transformation law (1.9), and the proof is complete.\hfill$\blacksquare$

As we have seen in Section 2, $N_L$ defines two semisprays on $E$ given by $$S^i_1=-\dfrac12 g^{ik}\dfrac{\pp^2 L}{\pp t\pp y^k}-\dfrac{\pp G^i}{\pp y^k}y^k\ \textrm{ and }S^i_0=-\dfrac{\pp G^i}{\pp y^k}y^k,$$ respectively. They coincide if, for instance, $\dfrac{\pp L}{\pp y^k}$ do not depend on $t$. Note that these semisprays are determined by $L$ only.

{\it Remark 3.1.} The nonlinear connection $N_L$ is without torsion. Let us consider an 1-form $\oo=\dfrac{\pp L}{\pp y^i}dx^i+\(L-\dfrac{\pp L}{\pp y^i}y^i\)dt$ on $E$ and let $$\theta=d\oo=\left[d\(\dfrac{\pp L}{\pp y^i}\)-\dfrac{\pp L}{\pp x^i}dt\right]\wedge(dx^i-y^i dt).$$ Thus, $\theta$ defines a contact structure on $E$. A vector field $X$ on $E$ is said to be characteristic for $\theta$ if the inner product of $\theta$ by $X$ vanishes, that is $X\cdot\theta=0$.

A curve on $E$ is said to be characteristic for $\theta$ if its tangent vector field is characteristic for $\theta$. We get the following theorem.

{\bf Theorem 3.2.} {\it If a curve $c:t\to c(t)$ on $M$ is an extremal for ${\cal L}$, then the curve $\sigma(t)=(t,\dot{c}(t))$ is a characteristic curve for $\theta$.}

\ms

\noindent{\it Proof.} As we have seen before, $\theta=\varphi_i\wedge\psi^i$, where $\varphi_i=d\(\dfrac{\pp L}{\pp y^i}\)—\dfrac{\pp L}{\pp x^i}dt$ and $\psi^i=dx^i-y^idt$.

Next, $(\dot{\sigma}(t)\cdot\theta)(Y)=\varphi_i(\dot{\sigma}(t))\psi^i(Y)-\varphi_i(Y)\psi^i(\dot{\sigma}(t))$ for any $Y\in\chi(E)$. But $\psi^(\dot{\sigma}(t))=0$, since along $\sigma$, $y^i=\dfrac{dx^i}{dt}$ and $\varphi_i(\dot{\sigma}(t))=0$ by virtue of the Euler-Lagrange eqs. Thus, $\dot{\sigma}(t)\cdot\theta=0$. \hfill$\blacksquare$

This theorem opens up a way in which contact geometry can come into the theory of Lagrangian systems [2]. We do not follow this way. Our geometrization is centered on a metrical structure derived from a regular time-dependent Lagrangian.

Finally, if we compare (3.4) with (2.7), we see that if $G^k$ is homogeneous of degree two with respect to $y$, a fact which is equivalent to $N^k_i y^i=2G^k$, it follows that the extremals of ${\cal L}$ coincide with the paths of the canonical nonlinear connection $N_L$ and with the integral curves of the semispray associated with $N_L$ as well.

\section{A metrical almost contact structure on $E$}
Let $R^n=(M,L(t,x,y))$ be a rheonomic Lagrange space. The canonical nonlinear connection produces a decomposition of the tangent bundle $TE$ as a direct sum $TE=HE\oplus VE$. Let $(\de_0,\de_i,\dot{\pp}_i)$ be the local frame adapted to this decomposition and $(dt,dx^i,\de y^i)$ its dual. Let us consider a linear map $F:T_uE\to T_uE$ given by $$F(\de_0)=0,\ \ F(\de_i)=-\pp_i,\ \ F(\dot{\pp}_i)=\de_i.\leqno(4.1)$$ Then $u\to F_u$, $u\in E$, defines an $(1,1)-$tensor field on $E$. It is obvious that rank$F=2n$ and an easy calculation gives $F^3+F=0$. Thus $F$ defines an $f(3,1)-$structure on $E$ [11].

{\bf Theorem 4.1.} {\it Let $RL^n=(M,L(t,x,y))$ be a rheonomic Lagrange space. Then the manifold $E=\R\times TM$ carries an almost contact structure $(F,\de_0,dt)$.}

\ms

\noindent{\it Proof.} We have $dt(\de_0)=1$ and equation (4.1) gives $F^2(\de_i)=-\de_i$, $F^2(\dot{\pp}_i)=-\dot{\pp}_i$. Thus it
follows that $F^2=-I+\de_0\times dt$.

The torsion tensor field [14] of the almost contact structure $(F,\de_0,dt)$ reduces to the Nijenhuis tensor $N_F$ of $F$. Thus, the almost contact structure $(F,\de_0,dt)$ is normal if and only if $N_F=0$.

Evaluating $N_F$ in the frame $(\de_0,\de_i,dt)$ one obtains the following theorem.

{\bf Theorem 4.2.} {\it The almost contact structure $(F,\de_0,dt)$ is normal if and only if

$(1)$ The canonical nonlinear connection $N_L=(N^k_0,N^k_i)$ is without curvature; and

$(2)$ $\dot{\pp}_i N^k_0=0$.}

As it is easy to check, the functions $(g_{ij})$ given by (3.1) are the components of a $d-$tensor of type $(0,2)$ on $E$. This will be called the metrical or fundamental tensor field of $RL^n$. Using it we can define the following $(0,2)-$tensor field on $E$: $$G=dt\otimes dt+g_{ij}dx^i\otimes dx^j+g_{ij}\de y^i\otimes\de y^j.\leqno(4.2)$$ We notice that if $(g_{ij})$ is positive definite, then $G$ is a Riemann metric on $E$. Otherwise, it is said it defines a metrical structure on $E$.

{\it Remark 4.1.} The horizontal and vertical distributions are orthogonal each to the other with respect to $G$.

A direct calculation in the frame $(\de_0,\de_i,\dot{\pp}_i)$ gives the following result.

{\bf Theorem 4.3.} {\it The metrical structure $G$ satisfies the following equations $$\begin{array}{l}G(FX,FY)=G(X,Y)-dt(X)dt(Y),\\ \\ dt(X)=G(\de_0,X),\ \ X,Y\in\chi(E).\end{array}\leqno(4.3)$$}

In other words, the previous theorem says that $(F,\de_0,dt,G)$ is a metrical almost contact structure on $E$. Recall that this metrical almost contact structure is completely determined by $L$. The particular form of $L$ could provide examples of structures which cover the classifications quoted in [11].

\section{$N-$linear connections on $E$}
Now we shall consider a class of linear connections on $E$ which are compatible with a nonlinear connection $N$, in particular with $N_L$, as well as with the almost contact structure associated with it. These will be called $N-$linear connections, recalling their compatibility with $N$.

The decomposition $T_uE=H_uE\oplus V_uE$ produced by a nonlinear connection $N$ induces a decomposition $$X=X^H+X^V,\ \ X\in\chi(E),\leqno(5.1)$$ where $X^H(X^V)$, is a vector field on $E$ taking its values in horizontal (vertical) distribution.

The decomposition (5.1) induces a decomposition of any tensor field on $E$ in horizontal and vertical parts. We denote also by $h$ and $v$ the horizontal and vertical projectors defined by (5.1), and then $P=h-v$ is an almost product structure on $E$.

{\bf Definition 5.1.} A linear connection $D:\chi(E)\times\chi(E)\to\chi(E)$, $(X,Y)\to D_X Y$ is said to be an $N-$linear connection if $$(a)\ \ D_XP=0,\ \ (b)\ \ D_XF=0,\ \ (c)\ \ D_X \de_0=0,$$ hold for any $X\in\chi(E)$.

Condition $(a)$ is equivalent to the fact that $D_X$ preserves by parallelism the horizontal and vertical distributions, i.e., $(D_X Y^H)^V=0$ and $(D_X Y^V)^H=0$, or $D_XY=(D_X Y^H)^H+(D_X Y^V)^V$. Now, if one sets $$D^H_X Y=D_{X^H}Y,\ \ D_{X}^H f=X^Hf,\ \ f\in{\cal F}(E)\leqno(5.2)$$ and extends $D^H_X$ to any $d-$tensor field on $E$ by the usual method, one obtains an operator called the $h-$covariant derivation in the algebra of $d-$tensor fields on $E$. Similarly, one may construct an operator for the $v-$covariant derivation, setting $$D^V_X Y=D_{X^V}Y,\ \ D_{X^V}f=X^V f,\ \ f\in{\cal F}(E).\leqno(5.3)$$ Now we state the following local characterization of an $N-$linear connection.

{\bf Theorem 5.1.} {\it To give an $N-$linear connection on $E$ is equivalent to give in every local chart on $E$, a set of functions $D\Gamma=(L^i_{j0},L^i_{jk},C^i_{jk})$ which satisfy on overlaps,} $$\begin{array}{c}\wt{L}^i_{j0}\dfrac{\pp\wt{x}^j}{\pp x^h}=\dfrac{\pp\wt{x}^i}{\pp x^k}L^k_{h0},\\ \\ \wt{L}^i_{jk}\dfrac{\pp\wt{x}^j}{\pp x^r}\dfrac{\pp\wt{x}^k}{\pp x^s}=\dfrac{\pp\wt{x}^i}{\pp x^h}L^h_{rs}-\dfrac{\pp^2\wt{x}^i}{\pp x^r\pp x^s},\\ \\ C^i_{jk}\dfrac{\pp\wt{x}^j}{\pp x^r}\dfrac{\pp\wt{x}^k}{\pp x^s}=\dfrac{\pp\wt{x}^i}{\pp x^h}C^h_{rs}.\end{array}\leqno(5.4)$$

\ms

\noindent{\it Proof.} If we express $D_X Y$ in a local chart, it comes out that it is well-defined by $$\begin{array}{ll}D_{\de_0}\de_j=L^0_{j0}\de_0+L^i_{j0}\de_i, & D_{\de_0}\dot{\pp}_j=M^i_{j0}\dot{\pp}_i,\\ \\
D_{\de_k}\de_j=L^0_{jk}\de_0+L^i_{jk}\de_i, & D_{\de_k}\dot{\pp}_j=M^i_{jk}\dot{\pp}_i,\\ \\
D_{\dot{\pp}_k}\de_j=Q^0_{jk}\de_0+Q^i_{jk}\de_i, & D_{\dot{\pp}_k}\dot{\pp}_j=C^i_{jk}\dot{\pp}_i,\end{array}$$ where (a) and (c)  from Definition 5.1 were taken into consideration. Taking into consideration (b) from the same definition, these equations reduce to
$$\begin{array}{lll}D_{\de_0}\de_j=L^i_{j0}\de_i, & D_{\de_k}\de_j=L^i_{jk}\de_i, & D_{\dot{\pp}_k}\de_j=C^i_{jk}\de_i,\\ \\ D_{\de_0}\dot{\pp}_j=L^i_{j0}\dot{\pp}_i, & D_{\de_k}\dot{\pp}_j=L^i_{jk}\de_i, & D_{\dot{\pp}_k}\dot{\pp}_j=C^i_{jk}\dot{\pp}_i,\end{array}\leqno(5.5)$$ and thus a set of functions $D\Gamma=(L^i_{j0},L^i_{jk},C^i_{jk})$ appears. If a transformation of coordinates on $E$ is performed, it turns out that these functions satisfy (5.4).

Conversely, given a set of functions $D\Gamma$, which on overlaps satisfy (5.4), by using (5.5), a well-defined linear connection on $E$ is obtained, and by a direct calculation one proves it satisfies (a), (b), and (c) from Definition 5.1, that is, it is an $N-$linear connection.\hfill$\blacksquare$

We notice that (5.4) shows that $(L^i_{j0})$ are the components of a $d-$tensor field of type $(1,1)$, $(C^i_{jk})$ are the components of a $d-$tensor field of type $(1,2)$ and $(L^i_{jk})$ define a Douglas connection.

Let $$T=t^{0...i...}_{...j...}\de_0\otimes...\times\de_i\otimes...\times \de y^j\otimes...$$ be a $d-$tensor field on $E$. By (5.5), one obtains the $h-$ and $w-$covariant derivative of $T$ as follows: $$\begin{array}{l}D_{\de_0}T=t^{0...i...}_{...j...|0}\de_0\otimes...\otimes\de_i\times...\otimes\de y^j\otimes...,\\ \\
D^H_X T=X^k t^{0...i...}_{...j...|k}\de_0\otimes...\times\de_i\otimes...\times\de y^j\otimes...,\\ \\ D^V_X T=\dot{X}^k t^{0...i...}_{...j...|k}\de_0\otimes...\times\de_i\otimes...\times\de y^j\otimes...,\end{array}$$ when $X=X^k\de_k+\dot{X}^k\dot{\pp}_k$, where we have put $$\begin{array}{l}t^{0...i...}_{...j...|0}=\de_0 t^{0...i...}_{...j...}+L^i_{k0}t^{0...k...}_{...j...}-L^k_{j0}t^{0...i...}_{...k...},\\ \\ t^{0...i...}_{...j...|h}=\de_h t^{0...i...}_{...j...}+L^i_{kh}t^{0...k...}_{...j...}-L^k_{jh}t^{0...i...}_{...k...},\\ \\ t^{0...i...}_{...j...|h}=\dot{\pp}_h t^{0...i...}_{...j...}+L^i_{kh}t^{0...k...}_{...j...}-C^k_{jh}t^{0...i...}_{...k...},\end{array}\leqno(5.6)$$ The torsion of an $N-$linear connection decomposes because of (5.1) into five $d-$tensor fields (the sixth identically vanishes) whose local components, in the adapted frame, are the following: $$\begin{array}{lll}T^i_{jk}=L^i_{jk}-L^i_{kj}, & R^i_{\a\beta}=\de_\beta N^i_\a-\de_\a N^i_\a, & C^i_{jk}, \\ \\ P^i_{0j}=\dot{\pp}_j N^i_0 - L^i_{j0}, & P^i_{kj}=\dot{\pp}_j N^i_k-L^i_{jk}, & S^i_{jk}=C^i_{jk}- C^i_{kj}.\end{array}\leqno(5.6')$$ All these functions will be called torsions of $D\Gamma$.

Analogously, one may prove that the curvature of an $N-$linear connection is locally determined by the following functions called curvatures of $D\Gamma$: $$\begin{array}{l}R^{\ i\ }_{j\ \a\beta}=\de_\beta L^i_{\ j\a}-\de_\a L^i_{\ j\beta}+L^i_{h\beta}L^h_{j\a}-L^h_{j\beta}L^i_{h\a}+C^i_{\ jh}R^h_{\a\beta},\\ \\ P^{\ i\ }_{j\ 0k}=\dot{\pp}_k L^i_{\ j0}-C^i_{\ jk|0}+C^i_{\ jh}P^h_{\ 0k},\\ \\ P^{\ i\ }_{j\ hk}=\dot{\pp}L^i_{\ jh}-C^i_{\ jh|k}+C^i_{\ js}P^s_{\ hk}, \\ \\ S^{\ i\ }_{j\ hk}=\dot{\pp}_k C^i_{\ jh}-\dot{\pp}_h C^i_{\ jk}+C^s_{\ jh}C^i_{\ sk}-C^s_{\ jk}C^i_{\ sh}.\end{array}\leqno(5.7)$$ We say $T^i_{jk}$ is the $h(hh)-$torsion of $D\Gamma$ and $S^i_{\ jk}$ is $v(vv)-$torsion of $D\Gamma$. The torsions and curvatures of $D\Gamma$ satisfy a number of a Bianchi identities. We do not write them here.

\section{Metrical $N-$linear connections}
Let $RL^n=(M,L)$ be a rheonomic Lagrange space and let us consider $N-$linear connections on $E$ which are compatible with the metrical structure $G$ defined by $L$.

{\bf Definition 6.1.} An $N-$linear connection $D$ on $E$ is said to be metrical if $D_X g=0$, for any $X\in\chi(E)$.

A direct calculation in local coordinates leads to the following result.

{\bf Theorem 6.1.} {\it An $N-$linear connection $D=(L^i_{j0},L^i_{jk}, C^i_{jk})$ is metrical if and only if $$g_{ij|0}=0,\ \ g_{ij|k}=0,\ \ g_{_{ij}|_k}=0.\leqno(6.1)$$}

As to the existence of metrical $N-$linear connections, we have the following theorem.

{\bf Theorem 6.2.} {\it Let $RL^n=(M,L(t,x,y))$ be a rheonomic Lagrange space. If $T^i_{\ jk}$ and $S^i_{\ jk}$ are two arbitrary skew-symmetric $d-$tensor fields on $E={\R}\times TM$, then there exists a set of metrical $N-$linear connections on $E$, such that each of them has $T^i_{\ jk}$ and $S^i_{\ jk}$ as $h(hh)-$ and $v(vv)-$torsions; respectively. The local coefficients of a connection from this set are given as follows: $$\begin{array}{l}L^k_{\ i0}=\dfrac12 g^{kh}\de_0 g_{ih}+O^{jk}_{ih}X^h_{j0},\\ \\ L^k_{\ ij}=\dfrac12 g^{kh}(\de_i g_{hj}+\de_j g_{ih}-\de_h g_{ij}+g_{is}T^s_{\ jk}+g_{js}T^s_{\ ih}+g_{hs}T^s_{\ ij}),\\ \\
C^k_{\ ij}=\dfrac12 g^{kh}(\dot{\pp}_i g_{hj}+\dot{\pp}_j g_{ih}-\dot{\pp}_h g_{ij}+g_{is}S^s_{\ jk}+g_{js}S^s_{\ ih}+g_{hs}S^s_{\ ij}),\end{array}\leqno(6.2)$$ where $X^h_{\ j0}$ is an arbitrary $d-$tensor field on $E$, and $O^{jk}_{ih}$ denotes the Obata operator $$O^{jk}_{ih}=\dfrac12(\de^j_i\de^k_{h}-g_{ih}g^{jk}).\leqno(6.3)$$}

\ms

\noindent{\it Proof.} The condition $g_{ij|k}=0$ is equivalent to $\de_k g_{ij}=L^h_{\ ik}g_{hj}+L^h_{\ jk}g_{ih}$. Permuting $(k,i,j)$ to $(i,j,k)$ and $(j,k,i)$ in this equality, adding two and subtracting one from the equalities thus obtained, and denoting $L^i_{\ jk}-L^i_{\ kj}=T^i_{\ jk}$, one obtains $L^k_{\ ij}$ in the form (6.2). One may proceed analogously in order to obtain $C^k_{\ ij}$ as in (6.2), using $g_{ij|k}=0$ and denoting $C^i_{\ jk}-C^i_{\ kj}=S^i_{\ jk}$.

Next, it is easy to check that $L^k_{\ i0}=\dfrac12 g^{kh}\de_0 g_{ih}$ are solutions of the equations $g_{ij|0}=\de_0g_{ij}-L^k_{\ i0}g_{kj}-L^k_{\ j0}g_{ik}=0$, in the unknowns $L^k_{\ i0}$.

Now, if $L^k_{\ i0}$ are any solutions of these equations, then $B^k_{\ i0}=L^k_{\ i0}-\dfrac12 g^{kh}\de_0 g_{ih}$ satisfy the equations $g_{ki}B^k_{\ j0}+g_{jk}B^k_{\ i0}=0$. The general solutions of these equations are $B^k_{\ j0}=O^{jk}_{ih}X^h_{j0}$, where $X^h_{j0}$ is an arbitrary $d-$tensor field. Thus $L^k_{\ i0}$ has the form given in (6.2).\hfill$\blacksquare$

Taking $X^i_{\ j0}=0$ in (6.2) one obtains the following corollary.

{\bf Corollary 6.1.} {\it Let $RL^n=(M,L(t,x,y))$ be a rheonomic Lagrange space and $T^i_{\ jk}$, $S^i_{\ jk}$ be two arbitrary skew-symmetric $d-$tensor fields on $E$. Then, there exists an unique metrical $N-$linear connection $D\Gamma=\(\dfrac12 g^{kh}\de_0 g_{ih},L^k_{\ ij},C^k_{\ ij}\)$, whose $h(hh)-$torsion is $T^i_{\ jk}$ and $v(vv)-$torsion is $S^i_{\ jk}$. The coefficients $L^i_{\ jk}$ and $C^i_{\ jk}$ are given by (6.2).}

\ms

\noindent{\it Proof.} The uniqueness of $L^i_{\ jk}$ and $C^i_{\ jk}$ from (6.2) follows by contradiction.\hfill$\blacksquare$

In particular, taking $T^i_{\ jk}=S^i_{\ jk}=0$, one obtains the next corollary.

{\bf Corollary 6.2.} {\it Let $RL^n= (M,L(t,x,y))$ be a rheonomic Lagrange space. Then, there exists a set of metrical $N-$linear connections, such that each of them has the vanishing $h(hh)-$ and $v(vv)-$torsion. The local coefficients of any connection of this set are given by $$\begin{array}{l}L^k_{\ i0}=\dfrac12 g^{kh}\de_0 g_{hi}+O^{jk}_{ih}X^h_{\ j0},\\ \\
L^k_{\ ij}=\dfrac12 g^{kh}(\de_i g_{hj}+\de_j g_{hi}-\de_h g_{ij}),\\ \\
C^k_{\ ij}=\dfrac12 g^{kh}(\dot{\pp}_i g_{hj}+\dot{\pp}_j g_{hi}-\dot{\pp}_h g_{ij}),\end{array}\leqno(6.4)$$ where $X^h_{\ j0}$ is an arbitrary $d-$tensor field on $E$.}

{\bf Definition 6.2.} The metrical $N-$linear connection whose local coefficients are given by (6.4), with $X^h_{\ j0}=0$ will be called the canonical metrical $N-$linear connection on $E$. It will be denoted by $\overset{c}{D}\Gamma$.

The $N-$linear connection $\overset{c}{D}\Gamma$ is completely determined by the time-de\-pen\-dent Lagrangian $L$. Thus $\overset{c}{D}\Gamma$ is similar to the connection $C\Gamma$ in Lagrange spaces [3].

The $h-$ and $v-$covariant derivatives of $C=y^i\dot{\pp}_i$, with respect to $\overset{c}{D}\Gamma$ lead us to introduce the following deflection tensors for $D$: $$D^i_{\ o}=y^i_{|0},\ \ D^i_{\ k}=y^i_{|k},\ \ d^i_{\ k}=y^i_{|_k}.\leqno(6.5)$$ Setting $D_{k0}=g_{ki}y^i_{|0}$, $D_{kj}=g_{ki}y^i_{|j}$, $d_{kj}=g_{ki}d^i_{\ j}$, and keeping in mind that $\overset{c}{D}\Gamma$ is metrical, one gets $$\begin{array}{l}D_{i0|k}-D_{ik|0}=R_{ji0k}y^j-d_{ih}R^h_{\ 0k},\\ \\
D_{hi|k}-D_{hk|i}=R_{jhik}y^j-d_{hs}R^s_{\ ik},\\ \\ D_{h0}|_k-d_{hk|0}=P_{jh0k}-d_{kj}P^j_{\ 0k},\\ \\ D_{hi}|_k-d_{hk|i}=P_{jhik}y^j-D_{hj}C^j_{\ ik}-d_{hj}P^j_{\ ik},\\ \\ d_{ik}|_h-d_{ih}|_k=S_{jikh}y^j.\end{array}\leqno(6.6)$$
We may also introduce the $h-$ and $v-$electromagnetic tensor fields, respectively, $$F_{ij}=\dfrac12(D_{ij}-D_{ji}),\ \ f_{ij}=\dfrac12(d_{ij}-d_{ji}).\leqno(6.7)$$ As it is easy to check, $f_{ij}=0$. As for $F_{ij}$, a direct calculation gives the following result.

{\bf Theorem 6.3.} {\it The tensor field $F_{ij}$, given by $(6.7)$ satisfies the following Maxwell equations $$\begin{array}{l}F_{ij|k}+F_{jk|i}+F_{ki|j}=-\dd\sum_{(i,j,k)}R^h_{\ jk}C_{ish}y^s,\\ \\ F_{ij}|_k+F_{jk}|_i+F_{ki}|_j=0,\end{array}\leqno(6.8)$$}

\section{Some time-dependent lagrangians}
Let $\wt{TM}$ be the manifold of nonvanishing vectors on $M$ and let $F:\R\times TM\to\R$ be a smooth function on $\R\times\wt{TM}$ and only continuous at the points $(t,x,0)$. Assume $F$ is positive on $\R\times\wt{TM}$ and homogeneous of degree one with respect to $y$.

A quadratic form is defined by $$h_{ij}(t,x,y)=\dfrac12\dot{\pp}_i\dot{\pp}_j F^2.\leqno(7.1)$$ If this is positive definite, $h_{ij}$ will be called a rheonomic Finsler metric on $M$ and the pair $RF^n=(M,F)$ will be called a rheonomic Finsler space. If we set $L=F^2$, it turns out that $(M,L)$ is a rheonomic Lagrange space. Thus, we may study the geometry of $RF^n$ regarding it as a rheonomic Lagrange space whose Lagrangian $L$ is positive, differentiable only on ${\R}\times\wt{TM}$ and homogeneous of degree two with respect to $y$.

Thus, the canonical nonlinear connection for $(M,L)$ will be called the Cartan nonlinear connection of $RF^n$, and the canonical metrical $N-$linear connection of $(M,L)$ will be called the Cartan metrical connection of $RF^n$, in such a way that the terminology corresponds to that from Finsler geometry [12]. By the Euler theorem on homogeneous functions, one finds $$L=F^2=h_{ij}(t,x,y)y^iy^j.\leqno(7.2)$$ Introducing the Cartan tensor fields of $RF^n$, $$C_{ijk}=\dfrac12\dot{\pp}_i h_{jk},\ \   C_{ij0}=\dfrac12\dfrac{\pp g_{ij}}{\pp t},\ \ C_{0ijk}=\pp_0 C_{ijk},\leqno(7.3)$$ where $\pp_0$ do stands for $\dfrac{\pp}{\pp t}$ the same theorem leads to the next proposition.

{\bf Proposition 7.1.} {\it The following identities hold: $$\begin{array}{l}y^i C_{ijk}=y^i C_{jik}=y^i C_{jki}=0,\\ \\ y^i C_{0ijk}=y^i C_{0jik}=y^i C_{0jki}=0.\end{array}\leqno(7.4)$$}

Introducing the usual Christoffel symbols, $$\gamma_{\ jk}^i=\dfrac12 h^{ir}(\pp_j h_{rk}+\pp_k h_{jr}-\pp_r h_{jk}),$$ we may state the following result.

{\bf Theorem 7.1.} {\it The local coefficients of the Cartan nonlinear connection are as follows $$N^i_0=\dfrac12 h^{ik}\pp_0 \dot{\pp}_k F^2,\ \ N^i_{\ k}=\dot{\pp}_k G^i,\ \textrm{where }\leqno(7.5)$$ $$G^i=\dfrac12\gamma^i_{\ jk}y^j y^k.\leqno(7.6)$$}

{\bf Theorem 7.2.} {\it The Cartan metrical connection $\overset{c}{F}T=(\overset{c}{F}^i_{\ j0},\overset{c}{F}^i_{\ jk},\overset{c}{C}^i_{\ jk})$ is as follows: $$\begin{array}{l}\overset{c}{F}^i_{\ j0}=\dfrac12 h^{ik}\de_0 h_{kj},\\ \\ \overset{c}F^i_{\ jk}=\dfrac12 h^{is}(\de_j h_{sk}+\de_k h_{js}-\de_s h_{jk}),\\ \\ \overset{c}{C}^i_{\ jk}=h^{is}C_{sjk},\end{array}\leqno(7.7)$$ where $\de_j$ is constructed by the help of $(7.5)$.}

The proofs are achieved by direct calculation. Also, by a direct calculation, one gets $$\begin{array}{l}N^i_0 y^i=\pp_0 F^2, \ \ \textrm{with\ }y_i=h_{is}y^s,\\ \\ y^i_{\ |k}=0,\ \ F^2_{\ |k}=0,\ \ y^i|_k=\de^i_{\ k},\ \ F^2|_k=2h_{ik}y^i.\end{array}\leqno(7.8)$$ By (7.8), the deflection tensors of an $RF^n$ space are $D_{ij}=0$, $d_{ij}=h_{ij}$, and thus the $h-$ and $v-$electromagnetic tensors identically vanishes. Thus, no rheonomic Finsler space supports a theory of electromagnetism.

It is clear that $h_{ij}$ from (7.1) is $0-$homogeneous with respect to $y$. This fact suggests that we consider rheonomic Lagrange spaces whose metrical tensor fields are $0-$homogeneous with respect to $y$, that is, the functions $g_{ij}$ given by (3.1) are $0-$homogeneous with respect to $y$. As to the general form of the Lagrangians of these spaces, we have the following theorem.

{\bf Theorem 7.3.} {\it If the metrical tensor field $(g_{ij})$ of a space $RL^n=(M,L)$ is $0-$homogeneous with respect to $y$, then $L$ has the general form $$L(t,x,y)=g_{ij}(t,x,y)y^i y^j+A_i(t,x)y^i+U(t,x),\leqno(7.9)$$ where $A_i$ is a covector field and $U$ is a real function on $\R\times M$.}

\ms

\noindent{\it Proof.} If we put $\overset{\circ}{L}(t,x,y)=g_{ij}(t,x,y)y^iy^j$, by the homogeneity of $g_{ij}$, it follows $\dot{\pp}_i\dot{\pp}_j(L— \overset{\circ}{L})=0$, which implies (7.9).\hfill$\blacksquare$

The following time-dependent Lagrangian, a particular form of (7.9), $$L(t,x,y)=a_{ij}(t,x)y^i y^j+A_i(t,x)y^i+U(t,x),\leqno(7.10)$$ where $(a_{ij})$ is a time-dependent Riemann metric, was used in treating some problems of dynamics [15-17].

Let's apply our previous theory to a rheonomic Lagrange space with the time-dependent Lagrangian (7.10). First, we note that its fundamental metric tensor field is just $(a_{ij}(t,x))$. The canonical nonlinear connection is given by $$\begin{array}{l}N^k_0(t,x,y)=a^{kh}(\pp_0 a_{hi})y^i+\pp_0 a(t,x),\\ \\ N^k_h(t,x,y)=a^k_{\ hi}(t,x)y^i-a^{kj}A_{hj},\end{array}\leqno(7.11)$$ where $$A_{jh}=\dfrac12 \(\dfrac{\pp A_j}{\pp x^h}-\dfrac{\pp A_h}{\pp x^j}\),\leqno(7.12)$$ and $a^k_{\ hi}$ denotes the Christoffel symbols constructed with $(a_{ij}(t,x))$.

The canonical metrical $N-$linear connection is $\overset{c}{D}\Gamma=(0,a^i_{\ jk}(t,x),0)$. The covariant deflection tensor fields are as follows: $D_{ij}=A_{ij}$, $d_{ij}=a_{ij}$. It comes out that $F_{ij}=A_{ij}$ and the term of $h-$electromagnetic tensor for $F_{ij}$ is supported by the form of $A_{ij}$. The Maxwell equations reduce to the classical ones.

Finally, we mention the possibility of ignoring that the metric tensor $g_{ij}(t,x,y)$ of a rheonomic Lagrange space is provided by a regular time-dependent Lagrangian, and study the geometry of the pair $(M,g_{ij}(t,x,y))$. Many results discussed in the above may be extended to this more general setting (cf. [5, Ch. XIII]).


\newpage

\runningauthor={M. ANASTASIEI}
\runningtitle={CERTAIN GENERALIZATIONS OF FINSLER METRICS}
\noindent
\baselineskip 8pt
\noindent{\footnotesize{Contemp. Math., 196, 161-169, Amer. Math. Soc.,\hfill\break Providence, RI, 1996}}
\vskip 2cm
\baselineskip 11.5pt plus .15pt
\centerline{\bf\Large CERTAIN GENERALIZATIONS}
\vskip .2cm
\centerline{\bf\Large OF FINSLER METRICS}
\vskip .5cm
\centerline{\bf by Mihai ANASTASIEI}
\vskip .5cm

\setcounter{section}{0}
\section{Introduction}
Scrutinizing the main body of results from Finsler geometry it
is observed that many of them depend on the Finsler metric
only and not on the fundamental Finsler function. Moreover, there are
many such results in which only some basic properties of
the Finsler metric are involved.

These facts led R. Miron ([9] [11]) to propose the study of
ge\-ne\-ra\-li\-zed Lagrange metrics, $GL$--metrics for
brevity, whose definition is tailored after the basic
properties of Finsler metrics.

The geometry of these metrics proved to be useful in the
Theory of Ge\-ne\-ral Relativity, Gauge Theory, and
Ecology (cf. [11] and references therein).

Certain problems from Mechanics and Theoretical Physics
require one, even at the Finslerian level, to study the geometry
of a $GL$--metric which, furthermore, depends on a
special variable analogous to  the physical time.
We  contributed to this study in [4].

The main objective of this paper is to review from our own
viewpoint the generalizations of Finsler metrics mentioned
above. We take this opportunity to cast a new light on
some well--known results and to add several new ones.
Some new examples are provided, too.

\section{Some properties of the Finsler metric}

Let $M$ be a real, smooth i.e. $C^\infty$, finite dimensional
manifold and $\tau:TM\to M$ its tangent bundle. Set
$\buildrel{\circ}\over{T}M=TM\setminus\{0_x\in T_xM$, $x\in
M\}$. Let $(U,(x^i))$ be a local chart on $M$. The indices
$i,j,k,...$ will run from $1$ to $n=\dim M$ and the
Einstein convention on summation will be implied. Associate
to $v\in \tau^{-1}(U)$ the coordinates $(x^i(\tau(u))$ and
$(y^i)$ provided by $v_{\tau(v)}=y^i\partial_i,$
$\partial_i:=\displaystyle\frac{\partial}{\partial x^i}$ and $TM$ becomes a smooth
orientable manifold. A change of coordinates
$(x^i,y^i)\to(x^{i'},y^{i'})$ on $TM$ is as follows:
$$x^{i'}=x^{i'}(x^1,...,x^n),\ y^{i'}=(\partial_jx^{i'})y^j,\
{\rm rank}\,(\partial_jx^{i'})=n.\leqno(2.1)$$

Let $F^n=(M,F)$ be a Finsler space and $\gamma_{ij}(x,y)$ the
local components of its Finsler metric. We list the following  known
properties of the Finsler metric, some of which are stated just
as the definition.\smallskip

P$_{\rm1}$. A change of coordinates (2.1) implies
$$\gamma_{ij}(x,y)=(\partial_ix^{i'})/(\partial_jx^{j'})
\gamma_{i'j'}(x',y').\leqno(2.2)$$

Thus $(\gamma_{ij}(x,y))$ are the components of a special,
{\it distinguished} tensor field on
$\buildrel{\circ}\over{T}M$ in the sense that their
transformation law (2.2) is similar with that of the
components of a tensor field on $M$. Throughout Finsler
geometry and  its generalizations one meets such
geometrical objects i.e. defined on
$\buildrel{\circ}\over{T}M$ or $TM$ but transforming under
(2.1), as being on $M$. We called them $d$--geometrical
objects, [11].\smallskip

P$_{\rm2}$. $\gamma_{ij}(x,y)=\gamma_{ji}(x,y)$\ \ \
(symmetry). \smallskip

P$_{\rm3}$. $\det(\gamma_{ij}(x,y))\ne0$\ \ \
(non--degeneracy). \smallskip

This property is usually postulated in a stronger form: the
quadratic form $\gamma_{ij}\xi^i\xi^j$, $(\xi^i)\in\R^n$
is positive definite.\smallskip

P$_{\rm4}$.
$\gamma_{ij}(x,y)=\displaystyle\frac{1}{2}\buildrel{\circ}\over{\partial}_i
\buildrel{\circ}\over{\partial}_jF^2,\
\buildrel{\circ}\over{\partial}_i:=\displaystyle\frac{\partial}{\partial y^i}.$\smallskip

P$_{\rm5}$. $\gamma_{ij}(x,y)$ are $p$
(positively)--homogeneous functions of zero degree  with respect to $(y^i)$.
Recall that $F$ is $p$--homogeneous of degree $1$ in $y^i$.\smallskip

P$_{\rm6}$. The function
$\buildrel{\circ}\over{C}_{ijk}=\displaystyle\frac{1}{2}
\buildrel{\circ}\over{\partial}_k\gamma_{ij}$ are the components
of a totally symmetric $d$--tensor field on
$\buildrel{\circ}\over{T}M$. Moreover,
$y^k\buildrel{\circ}\over{C}_{ijk}=0.$\smallskip

P$_{\rm7}$. Let $\buildrel{\circ}\over{\gamma}\!\!{}^i_{jk}(x,y)=
\displaystyle\frac{1}{2}\gamma^{ih}(\partial_j\gamma_{hk}+
\partial_k\gamma_{jh}-\partial_h\gamma_{jk})$
be the Christofell symbols derived from $(\gamma_{jk})$.
Then $\buildrel{\circ}\over{G}\!{}^i=\displaystyle\frac{1}{2}
\buildrel{\circ}\over{\gamma}\!\!{}^i_{jk} y^jy^k$
are the components of the (geodesic) spray
$\buildrel{\circ}\over{S}=y^i\partial_i+\buildrel{\circ}\over{G}\!{}^i
\buildrel{\circ}\over{\partial}_i$ on
$\buildrel{\circ}\over{T}M$ and
$\buildrel{\circ}\over{N}\!{}^i_j=
\buildrel{\circ}\over{\partial}_j\buildrel{\circ}\over{G}\!{}^i$
has the following  law of transformation under (2.1):
$$\buildrel{\circ}\over{N}\!{}^{i'}_{j'}(\partial_ix^{j'})=
(\partial_jx^{i'})\buildrel{\circ}\over{N}\!{}^j_i-
(\partial_i\partial_jx^{i'})y^j,\leqno(2.3)$$
that is, these functions are the coefficients of the nonlinear
Cartan connection.

Set
$\buildrel{\circ}\over{\delta}_i:=\partial_i\,-\buildrel{\circ}\over{N}\!{}^k_i
\buildrel{\circ}\over{\partial}_k$ and it results that
$\delta_i=(\partial_ix^{i'})\delta_{i'}$. For $v\in
\buildrel{\circ}\over{T}M,$ the linear space $H_v$ spanned
by $(\delta_i)_v$ is supplementary to
the vertical space $V_v={\rm Ker}(D\tau)_v$ spanned by
$(\buildrel{\circ}\over{\partial}_i)_v$, that is,
$$T_v\buildrel{\circ}\over{T}M=H_v\oplus V_v\ \ \ \
\mbox{(direct sum)}.\leqno(2.4)$$

P$_{\rm8}$. The function
$(\buildrel{\circ}\over{L}\!{}^i_{jk},\buildrel{\circ}\over{C}\!{}^i_{jk})$
given by
$$\begin{array}{lcl}
\buildrel{\circ}\over{L}\!{}^i_{jk}&=&
\frac{1}{2}\gamma^{ih}(\buildrel{\circ}\over{\delta}_j\gamma_{hk}+
\buildrel{\circ}\over{\delta}_k\gamma_{jh}-
\buildrel{\circ}\over{\delta}_h\gamma_{jk}),\\
\buildrel{\circ}\over{C}\!{}^i_{jk}&=&
\frac{1}{2}\gamma^{ih}(\buildrel{\circ}\over{\partial}_j\gamma_{hk}+
\buildrel{\circ}\over{\partial}_k\gamma_{jh}-
\buildrel{\circ}\over{\partial}_h\gamma_{jk})=
\gamma^{ih}\buildrel{\circ}\over{C}_{hjk},\end{array}\leqno(2.5)$$
are the local coefficients of the Cartan connection. This
connection is $h$--metrical
$(\gamma_{ij{\scriptstyle{\buildrel{\circ}\over{|}}}k}=0),$
$v$--metrical
$(\gamma^i_{jh}{\scriptstyle{\buildrel{\circ}\over{\big|}}}{}_{{}_k}=0)$,
$h$--symmetric
$(\buildrel{\circ}\over{L}\!{}^i_{jk}=\buildrel{\circ}\over{L}\!{}^i_{kj}),$
$v$--symmetric
$(\buildrel{\circ}\over{C}\!{}^i_{jk}=\buildrel{\circ}\over{C}\!{}^i_{kj})$
and is free of deflection
$(\buildrel{\circ}\over{D}\!{}^i_j=
y^i_{{\scriptstyle{\buildrel{\circ}\over{|}}}j}=0).$
In other words, it satisfies the well--known Matsumoto's
axioms. The list could be continued but these properties
are essential for developing  Finsler geometry.

\section{A generalization of the Finsler metrics: $GL$--metrics}

A collection of functions  $(g_{ij}(x,y))$ locally defined on
$TM$ and satisfying P$_{\rm1}$--P$_{\rm3}$ is called a
generalized Lagrange metric, shortly a $GL$--metric. As P$_{\rm7}$
cannot be recovered from P$_{\rm1}$--P$_{\rm3}$ only, we
introduce the assumption
\begin{itemize}
\item[(H$_{\rm1}$)] There exists a non--linear connection
on $TM$ i.e. a set of coefficients $(N^i_j(x,y))$ verifying
(2.3).\end{itemize}
This is always true if $M$ is paracompact. Notice that we
shall indicate the general case by deleting the superscript
``$\circ$'' from the entities previously introduced. Thus we
may consider $(\partial_i)$ and the decomposition (2.4) holds
for $v\in TM$. The functions  provided by (2.5) define a
connection with the first four properties of the Cartan
connection. Its deflection generally does not vanish. From
now on the torsions, the curvatures, the $h$-- and
$v$--paths and so on  can be introduced and studied.

The postulate (H$_{\rm1}$) is not very strong as the
hypothesis of paracompactness of $M$ is generally accepted.
But an arbitrary non--linear connection i.e. without any
relationship to $(g_{ij}(x,y))$ is far from useful.

Fortunately, in the most important examples there exists a
non--linear connection determined by or strongly related
to the given $GL$--metric.
\bigskip

\noindent {\bf Example.} For any positive functions  $a$ and $b$ on
$\buildrel{\circ}\over{T}M$ we set
$$g_{ij}(x,y)=a(x,y)\gamma_{ij}(x,y)+b(x,y)y_iy_j,\ y_i=\gamma_{ik}y^k.
\leqno(3.1)$$
This is a $GL$--metric. Indeed, it is easy to check that
$$g^{jk}=\frac{1}{a}\left(\gamma^{jk}-\frac{b}{a+bF^2}\
y^jy^k\right),\leqno(3.2)$$ verifies
$g_{ij}g^{jk}=\delta^k_i$.

In order to study it we have on hand the non--linear
connection $(\buildrel{\circ}\over{N}\!{}^i_j(x,y))$. We stress
that for various functions  $a$ and $b$ the $GL$--metric (3.1) supplies
all the $GL$-metrics treated in [11].

A $GL$-metric is said to be an $L$--metric if there exists a smooth function
$L:TM\to R$ such that
$$g_{ij}(x,y)=\frac{1}{2}
\buildrel{\circ}\over{\partial}_i
\buildrel{\circ}\over{\partial}_jL(x,y).\leqno(3.3)$$
Such a func\-tion, called a regular Lagrangian, exists if and only if
$(C_{ijk})$ is a totally symmetric $d$--tensor field. If
$L$ exists, it is not unique since $\tilde
L(x,y)=L(x,y)+\varphi_i(x)y^i+c$ is a new solution of (3.3).
Choosing such an $L$, the pair $(M,L)$ is called a Lagrange
space. In particular, $(M,F^2)$ is a Lagrange space.

For $L$--metrics, a canonical non--linear connection is
derived from the Euler-Lagrange equations provided by the
variational problem $\delta\!\displaystyle\int^{t_1}_{t_0}\!\!Ldt\!=\!0,$
by first considering $G^i=\displaystyle\frac{1}{4}g^{ik}[(
\buildrel{\circ}\over{\partial}_k\partial_jL)y^j-\partial_kL]$ (the
components of the ca\-no\-ni\-cal semi--spray) and then taking
$N^i_j=\buildrel{\circ}\over{\partial}_jG^i$.

If one requires that a $L$--metric be
$(m-2)-p-$homogeneous,
 then $L$ is uniquely
determined and is $m-(p)-$homogeneous. For such
$L$--metrics the functions $G^i$ and the connection
$(L^i_{jk},C^i_{jk})$ is deflection free. Thus these
$L$--metrics are closely related to Finsler metrics, [5], [6].

Coming--back to the Example, we notice that $(g_{ij}(x,y))$
is an $L$--metric if and only if
$$[(\buildrel{\circ}\over{\partial}_ka)\gamma_{ij}-
(\buildrel{\circ}\over{\partial}_ia)\gamma_{kj}]+
[(\buildrel{\circ}\over{\partial}_kb)y_i-
(\buildrel{\circ}\over{\partial}_ia)y_k]y_j+
b[y_i\gamma_{kj}-y_k\gamma_{ij}]=0.
\leqno(3.4)$$
Contracting this by $(\gamma^{ij})$ one gets
$$(n-1)(\buildrel{\circ}\over{\partial}_ka)-b(n-1)y_k+
(\buildrel{\circ}\over{\partial}_kb)F^2-
(\buildrel{\circ}\over{\partial}_ib)y^iy_k=0.
\leqno(3.5)$$\smallskip

\noindent {\bf Remark 3.1.} Even for simple functions $a$ and $b$, the
$GL$-metric  (3.1) does not reduce to an $L$--metric. For instance,
if $a$ and $b$ are positive constants, (3.5)
simplifies to $b(n-1)y_k=0,$ which does not hold for $n>1.$
So (3.4) fails.\bigskip

\noindent {\bf Remark 3.2.} Let $a=\alpha(F^2)$ and
$b=\beta(F^2)$ with $\alpha,\beta:{I\R}^*_+\to {\R}^*_+$. Then (3.5)
implies $\beta=2\alpha'$ and the condition $a+bF^2>0$
becomes $\alpha+2\alpha't>0,$ $t\to\alpha(t)$, $t\in{\R}^*_+$.
Set $\alpha=\varphi'$ with $\varphi:{\R}\to {\R}^*_+$, $\varphi'>0,$
$\varphi'(t)+2\varphi''t>0.$ One obtains the $\varphi$--Lagrange
metrics studied in [6].

\section{Almost Hermitian Model of a $GL$-metric}
Let $M$ be endowed with a $GL$-metric  $(g_{ij}(x,y))$ and a
non--linear connection $(N^i_j(x,y)).$ The decomposition
(2.4) implies a decomposition $X=hX+vY$ for every  vector field
(v.f.) $X$ on $TM$. Denote by $P$ the almost product structure
provided by the horizontal and vertical distributions
according to: $P(hX)=hX$, $P(vX)=-vX$. Consider also the
almost complex structure $F$ defined as follows:
$F(hX)=-vX,$ $F(vX)=hX.$ Next, setting
$G=g_{ij}(x,y)dx^i\otimes dx^j+g_{ij}\delta y^i\otimes\delta
y^j$, $\delta y^i=dy^i+N^i_kdx^k,$ one gets a metrical structure on
$TM$ which is a Riemannian structure if $(g_{ij})$ is positive
definite. It is easily seen that $(F,G)$ is an almost
Hermitian structure.

This simple construction has much more implications in the
geometry of $(g_{ij}(x,y))$ then it seems at a first insight.
Indeed, the coefficients $(L^i_{jk},C^i_{jk})$ supply a
linear connection $D$ on $TM$ given in the adapted basis
$(\delta_i,\buildrel{\circ}\over{\partial}_i)$ as follows:
$$\begin{array}{lclclcl}
D_{\delta_k}\delta_j&=&L^i_{jk}\delta_i,&\ \ \ \ &
D_{\buildrel{\circ}\over{\partial}_k}\delta_j&=&C^i_{jk}\delta_i,\\
D_{\delta_k}\buildrel{\circ}\over{\partial}_j&=
&L^i_{jk}\buildrel{\circ}\over{\partial}_i,&\ \ \ \ &
D_{\buildrel{\circ}\over{\partial}_k}\buildrel{\circ}\over{\partial}_j
&=&C^i_{jk}\buildrel{\circ}\over{\partial}_i.\end{array}\leqno(4.1)$$
This connection preserves the both distributions $(DP=0)$,
is almost complex $(DF=0)$ and is metrical $(DG=0)$.
Moreover, if its torsion $T$ is decomposed into vertical
and horizontal components, then $hT(h\cdot,h\cdot)=0$ and
$vT(v\cdot,v\cdot)=0$. Conversely, every linear connection
$D$ on $TM$, with the above properties, has
$(L^i_{jk},C^i_{jk})$ from (2.5) as local coefficients in
the adapted basis
$(\delta_i,\buildrel{\circ}\over{\partial}_j)$. Thus the
study of the connection (2.5) is equivalent to the study of
such a linear connection $D$ on $TM$ endowed with $(F,G)$.
This is a reason to call $(F,G)$ the almost Hermitian model
of $(g_{ij}(x,y))$. A first important application of this
model is due to R. Miron. He considered the Einstein
equations for $G$ and projecting them on horizontal and
vertical distributions, he arrived at a correct form of the
Einstein equations for $(g_{ij}(x,y))$, [9]. See also
[2],[3],[11],[12].

A simpler usage of the almost Hermitian model is as follows.
Looking for a meaning of the divergence of a $d$--vector
field $(X^i(x,y))$ we observe that it defines an horizontal
vector field $hX=X^i(x,y)\delta_i$ as well as a vertical vector field
$vX=X^i(x,y)\buildrel{\circ}\over{\partial}_i$. As $(TM,G)$ is an
orientable Riemannian manifold (the positiveness of
$(g_{ij})$ is implied), we define an
$h$--divergence $({\rm div}_hX)$ and a $v$--divergence
$({\rm div}_vX)$ according to ${\cal L}_{*X}dv=({\rm div}_*X)dv,$
$*=h,v$, where ${\cal L}$ means the Lie derivative and $dv$ is
the volume element associated to $G$.

Since $D$ has torsion, it comes out that the usual formula
for the divergence of a vector field $Z$ on $TM$ is ${\rm div}\,
Z=Trace(Y\to D_YZ+T(Z,Y))$. In the adapted basis one finds
${\rm div}_hX=X^i_{|i}-X^kP_k$, $P_k=P^i_{ki}$, $P^i_{jk}=
\buildrel{\circ}\over{\partial}_kN^i_j-L^i_{kj}$,
${\rm div}_vX=X^i_{ \big| i}+X^kC_k,$ $C_k=C^i_{ki}$. For the
$L$--metrics described in Remark 3.2, in particular, for Finsler
metrics we get ${\rm div}_h\, S=0$, a generalization of a Liouville
theorem from the Riemannian geometry. For any function $f$ on
$TM$ we have a $h$--gradient ${\rm
grad}_hf=(g^{ik}\delta_kf)\delta_i$ and a $v$--gradient:
${\rm grad}_vf=(g^{ik}\buildrel{\circ}\over{\partial}_kf)
\buildrel{\circ}\over{\partial}_i$.
Accordingly, we may define the $h$--Laplacean
$\Delta_hf={\rm div}_h({\rm grad}_hf)$ and the $v$--Laplacean
$\Delta_vf={\rm div}_v({\rm grad}_vf)$.

The function $\varepsilon=g_{ij}(x,y)y^iy^j$ is called the absolute
energy of the $GL$-metric  $(g_{ij}(x,y))$. For $L$--metrics
discussed in the Remark 3.2, the absolute energy is
$h$--harmonic i.e. $\Delta_h\varepsilon=0$.

Let $\C=y^i\buildrel{\circ}\over{\partial}_i$ be the
Liouville vector field on $TM$. The postulate (H$_{\rm1}$) is
clearly {\it implied} by the fol\-lowing one.
\begin{itemize}
\item[(H$_{\rm2}$)] There exists a linear connection
$\nabla$ in the vertical bundle which is re\-gu\-lar, that is,
the space $H_v=\{X_v\mid\nabla_{X_v}\C=0\}$ is
supplementary to $V_v$, $v\in TM$.\end{itemize}

Let $h_v$ be the inverse of the isomorphism $D_v\tau:H_v\to
T_{\tau(v)}M$ and $\delta_i=h_v(\partial_i)$. Then
$(D_v\tau)(\delta_i-\partial_i)=0.$ Thus
$\delta_i=\partial_i-N^k_i\buildrel{\circ}\over{\partial}_k$, where
the sign ``$-$'' is taken for the sake of convenience. The
functions  $(N^k_i)$ define a nonlinear connection.

Let the linear connection $\nabla$ be given as fol\-lows:
$\nabla_{\partial_k}\buildrel{\circ}\over{\partial}_j=
\Gamma^i_{jk}\buildrel{\circ}\over{\partial}_i$,
$D_{\buildrel{\circ}\over{\partial}_k}
\buildrel{\circ}\over{\partial}_j=B^i_{jk}\buildrel{\circ}\over{\partial}_i$.
The condition
$\nabla_{\delta_k}(y^j\buildrel{\circ}\over{\partial}_j)=0
\,\Longleftrightarrow\,(\delta^i_h+y^jA^i_{jh}) N^h_k=y^j\Gamma^i_{jk}$
shows that the regularity of $\nabla$ is equivalent to the
regularity of the matrix $(\delta^i_h+y^jA^i_{jh})$. The
triad
$(N^i_j,L^i_{jk}=\Gamma^i_{jk}-N^h_kB^i_{jh},B^i_{jk})$ is
an usual Finsler connection.

The postulate (H$_{\rm2}$) is involved in what we called
the vector bundle model of Finsler geometry (cf. [1]). This
model was recently used by D.Bao, S.S.Chern [7] and Z.Shen [13] for
solving some global problems in Finsler geometry. A variant
of it, usefull for Physics, was developed by J.G. Vargas
and D.Torr [14]. An essentially different model, called by us
the principal bundle model (cf. [1]), is due to M. Matsumoto [8].

\section{Finsler geometry of a vector bundle}
It is to be observed that the geometry of a $GL$--metric
essentially depends on a non--linear connection on $TM$.
The extension of this notion of connection to a vector
bundle $\pi:E\to M$ is quite natural. It is nothing but a
supplementary distribution to the vertical distribution
$u\to({\rm Ker}\,D\pi)_u$, $u\in E$. Notice that the
horizontal distribution is non--holonomic, so a study of
this pair of distributions is of interest.

If $E$ is endowed with a metrical structure $G$, we may take as
non--linear connection the orthogonal distribution to the
vertical distribution. Then $G$ takes the form
$$G=g_{ij}(x,y)dx^i\otimes dx^j+g_{ab}(x,y)\delta
y^a\otimes\delta y^b,\ \delta y^a=dy^a+N^a_kdx^k,\leqno(5.1)$$
where $(x^i,y^a),\ a,b,c,...=1,...,\ m=$ fibre dimension,
are the local coordinates on $E$ and $(N^a_i(x,y))$ are the
local coefficients of the non--linear connection defined by
$G$.

A change of coordinates $(x^i,y^a)\to(x^{i'},y^{a'})$ on
$E$ has the form
$$\begin{array}{ll}
x^{i'}=x^{i'}(x^1,...,x^n), & {\rm rank}(\partial_jx^{x'})=n,\\
y^{a'}=M^{a'}_b(x)y^b, & {\rm
rank}(M^{a'}_b)=m.\end{array}\leqno(5.2)$$
The coefficients $(N^a_i)$ of a non--linear connection have
the fol\-lowing  transformation law under (5.2):
$$N^{a'}_{i'}(\partial_ix^{i'})=M^{a'}_a(x)N^a_i-
(\partial_iM^{a'}_a(x))y^a.\leqno(5.3)$$
A geometrical study of the pair $(E,G)$ using the
above ingredients was performed by R. Miron [10]. Some
applications of his theory we have pointed out in [2],[3]
(see also [11]).

\section{Rheonomic $GL$--metrics}
Let us consider the functions  $(g_{ij}(t,x,y))$ with the
properties of a $GL$--metric. Assume $t$ remains unchanged
under (2.1), that is, $t$ is viewed as absolute time. We
call such a collection of functions a {\it rheonomic
$GL$--metric}, shortly a $RGL$--metric. It is clear that
this kind of $GL$--metric is living on $\R\times TM$, a
manifold which could be thought of as fibered in three
different ways projecting it on $\R$, $TM$ or ${\R}\times
TM$. Each of these fibrations has a certain value for
geometrizing problems from Mechanics or Calculus of
Variations. As more appropriate seems to be the fibration
$$\pi:{\R}\times TM\to{\R}\times M,\ \pi(t,v)=(t,\tau(v)),\
v\in TM.$$
Set $E=\R\times TM.$ The manifold $E$ is coordinizated by
$(t,x^i,y^i)$ and the $\pi$ takes the form
$(t,x^i,y^i)\to(t,x^i).$ It is convenient to put $x^0=t$
and to use the Greek indices $\alpha,\beta,\gamma,...$
ranging over $0,1,2,...,n.$ A non--linear connection can be
given by $(n+1)$ local vector fields, say $\delta_\alpha.$
Choosing $\delta_\alpha$ such that  they are projected to $\partial_\alpha$,
one gets
$$\delta_\alpha=\partial_\alpha-N^i_\alpha(t,x,y)
\buildrel{\circ}\over{\partial}_i.\leqno(6.1)$$

The invariance of the horizontal subspaces requires the
condition\break
$\delta_\alpha=(\partial_\alpha x^{\alpha'})\delta_{\alpha'}$, when
a change of coordinates on $\R\times TM$ is performed.
This implies the following law of transformation for
$(N^i_\alpha)$:
$$N^{i'}_{\alpha'}(\partial_\beta
x^{\alpha'})=(\partial_kx^{i'})N^k_\beta-
(\partial_\beta\partial_kx^{i'})y^k.\leqno(6.2)$$
If one rewrites (6.1) in the form
$$\delta_0=\partial_t-N_0^i(t,x,y)\buildrel{\circ}\over{\partial}_i,\
\delta_i=\partial_i-N_i^k(t,x,y)\buildrel{\circ}\over{\partial}_k,$$
it comes out from (6.2) that $(N^i_0(t,x,y))$ change like
the components of a $d$--vector field and $(N^i_j(t,x,y))$
change like the coefficients of a non--linear connection on
$TM$. Thus, we may identify a non--linear connection on $E$ with
the pair $(N^i_0,N^i_j)$.

Let $(\delta_0,\delta_i,\buildrel{\circ}\over{\partial}_i)$ be
the basis adapted to the decomposition $T_uE=H_uE\oplus
V_uE$ and $(dt,dx^i,\delta y^i)$ its dual. Denote by $P$
the almost product structure on $E$ associated as in \S4 to the
decomposition $T_uE=N_uE\oplus V_uE$ and define a tensor
field $\Phi$ of type $(1,2)$ on $E$ as fol\-lows:
$$\Phi(\delta_0)=0,\
\Phi(\delta_i)=-\buildrel{\circ}\over{\partial}_i,\
\Phi(\buildrel{\circ}\over{\partial}_i)=\delta_i.\leqno(6.3)$$
It easily comes out that $(\Phi,\delta_0,\delta t)$
is an almost contact structure on $E$. Using $(g_{ij}(t,x,y))$
the fol\-lowing metrical structure on $G$ is obtained
$$G=dt\otimes dt+g_{ij}dx^i\otimes dx^j+g_{ij}\delta
y^i\otimes\delta y^j.\leqno(6.4)$$
It is easy to see that $(\Phi,\delta_0,dt,G)$ is a metrical
almost contact struc\-tu\-re\ on $E$. This will be called the almost
contact model for $g_{ij}(t,x,y)$. As in the almost
Hermitian model it is quite natural to look for a linear
connection $D$ on $E$ with the properties:
$$DP=0,\ D\Phi=0,\ DG=0,\ D\delta_0=0,\
hT(h\cdot,h\cdot)=0,\ vT(v\cdot,v\cdot)=0.\leqno(6.5)$$
In the frame
$(\delta_0,\delta_i,\buildrel{\circ}\over{\partial}_i)$ this
connection has the coefficients
$(L^i_{j0},L^i_{jk},C^i_{jk})$, where the latter two are
similar with those from (2.5) while the first has the form
$$L^i_{j0}=\frac{1}{2}\,g^{ih}\delta_0g_{hj}+\frac{1}{2}
(\delta_j^k\delta^i_h-g_{jh}g^{ki})X^h_{k0},\leqno(6.6)$$
with $X^h_{k0}$ an arbitrary $d$--tensor field, cf. [4].

This set of connections allows us to develop the geometry
of the $RGL$--metric $g_{ij}(t,x,y)$.

As for $GL$--metrics, in the most important examples, the
non--linear connection is completely determined by
$(g_{ij})$. A $RGL$--metric will be called a rheonomic
$L$--metric if there exists  a smooth function $L:{\R}\times TM\to{\R}$
such that
$$g_{ij}(t,x,y)=\frac{1}{2}\,\buildrel{\circ}\over{\partial}_i
\buildrel{\circ}\over{\partial}_jL.\leqno(6.7)$$
If $L$ exists, it is not unique. Taking one $L$ as solution
of (6.6), the pair $(M,L)$ is called a rheonomic Lagrange
space. In particular, we arrive at the notion of theonomic
Finsler space as a pair $(M,F)$ with $F:R\times TM\to R$ a
positive function, smooth on $R\times\buildrel{\circ}\over{T}M$,
$p-$homogeneous of degree $1$ with respect to $(y^i)$ such that the functions
$g_{ij}(t,x,y)\!=\!\displaystyle\frac{1}{2}\!
\buildrel{\circ}\over{\partial}_{i}
\buildrel{\circ}\over{\partial}_j\!F^2$ satisfy
$\det(g_{ij}(t,x,y))\!\ne\!0.$ For a theory of
rheonomic Finsler and Lagrange spaces we refer to [4].

\section*{References}
\begin{description}
\item 1.
Anastasiei, M.,
{\it Models of Finsler and Lagrange geometry.}
Proc. of IV$^{\rm th}$ Nat. Sem. on Finsler and Lagrange
Spaces,
Bra\c sov, Romania, 1986, pp. 43--56.
\item 2.
Anastasiei, M.,
{\it Vector bundles. Einstein equations.}
Analele \c st. ale Univ. "Al.I.Cuza" Ia\c si,
T. XXXII, s.I-a, Mat., 1986, p. 17--24.
\item 3.
Anastasiei, M.,
{\it Conservation laws in the $\{V,H\}$--bundle model of
Relativity.}
Tensor N.S. 46 (1987), 323--328.
\item 4.
Anastasiei, M.,
{\it The geometry of Time--Dependent Lagrangians.}
Mathl. Comput. Modelling,
vol. 20, no. 415, 67--81, 1994.
\item 5.
Anastasiei, M., Antonelli P.L.,
{\it The Differential Geometry of Lagrangians which Generate
Sprays,} Lagrange and Finsler Geometry, Kluwer Academic Publishers, 1995.
\item 6.
Antonelli, P.-L., Hrimiuc, D.,
{\it A new class of spray--generating Lagrangians,} Lagrange and
Finsler Geometry, Kluwer Academic Pu\-bli\-shers, 1995.
\item 7.
Bao, D., Chern, S.S.,
{\it On a notable connection in Finsler geometry.}
Houston Journal of Mathematics, vol. 19 (1993), 135--180.
\item 8.
Matsumoto, M.,
{\it Foundations of Finsler Geometry and Special Finsler Spaces.}
Kaiseisha Press, Japan, 1986.
\item 9.
Miron, R.,
{\it A Lagrangian Theory of Relativity.}
Analele \c st. ale Univ.\break "Al.I.Cuza" Ia\c si,
T. XXXII, s.I-a, Mat., 1986, p. 37--62.
\item 10.
Miron, R.,
{\it Techniques of Finsler geometry in the theory of vector
bundles.} Acta Sci. Math. 49 (1985), 119--129.
\item 11.
Miron, R., Anastasiei, M.,
{\it The Geometry of Lagrange spaces: Theory and Applications.}
Kluwer Academic Publishers,
FTPH vol. 59, 1994, 285 p.
\item 12.
Miron, R., Ro\c sca, R., Anastasiei, M., Buchner, K.,
{\it New Aspects of Lagrangian Relativity.}
Foundations of Physics Letters, vol. 5 (1992), 141--171.
\item 13.
Shen, Z.,
{\it On a connection in Finsler Geometry.}
Houston Journal of Mathematics, vol. 20 (1994), 591--602.
\item 14.
Vargas, J., Torr, D.,
{\it Finslerian struc\-tu\-res: The Cartan--Clifton method of the moving
frame.} J. Math. Phys. 34 (10), 1993, 4898--4913.
\end{description}

\noindent{\footnotesize{Faculty of Mathematics\\
``Al.I.Cuza'' University of Ia\c si\\
R--6600 Ia\c si, Romania\\
E--mail: anastas@uaic.ro}}

\newpage





\def\G{\Gamma}

\runningauthor={M. ANASTASIEI}
\runningtitle={A HISTORICAL REMARK ON THE CONNECTIONS OF CHERN AND RUND}
\noindent
\baselineskip 8pt
\noindent{\footnotesize{Contemporary Mathematics\hfill\break
Volume 196,171-176, Amer. Math. Soc., Providence, RI, 1996}}
\vskip 2cm
\baselineskip 11.5pt plus .15pt
\centerline{\bf\Large A HISTORICAL REMARK ON}
\vskip .2cm
\centerline{\bf\Large THE CONNECTIONS OF CHERN AND RUND}
\vskip .5cm
\centerline{\bf \footnotesize{BY}}
\vskip .5cm
\centerline{\bf \footnotesize{M. ANASTASIEI}}
\vskip 1cm

\setcounter{section}{0}
\section{Introduction}
Let $M$ be a real, $n-$dimensional smooth (i.e. $C^\9$) manifold and $\tau:TM\to M$ its tangent bundle. In a local chart $(U,x^i)$ on $M$, a tangent vector $v\in T_pM$, $p\in M$, has the form $v=y^i\left.\dfrac{\pp}{\pp x^i}\right|_p$ and it is usual to take $(\tau^{-1}U,x^i\equiv x^i\circ\tau,y^i)$ as local coordinates on $TM$. Throughout the paper the indices run from 1 to $n$ and the Einstein convention on summation is implied.

A local change of coordinates $x^i\to\wt{x}^i$ on $M$ induces in turn a change of coordinates $(x^i,y^i)\to(\wt{x}^i,\wt{y}^i)$ on $TM$:
$$\begin{array}{l}\wt{x}^i=\wt{x}^i(x^1,...,x^n),\ {\rm{rank}}\(\dfrac{\pp\wt{x}^i}{\pp x^k}\)=n,\\ \\ \wt{y}^i=\dfrac{\pp\wt{x}^i}{\pp x^k}(x)y^k.\end{array}\leqno(1.1)$$ Set $(x,y):=(x^i,y^i)$ and $\overset{\circ}{TM}=TM\setminus\{(x,0)\}$.

A fundamental Finsler function is a function $F:TM\to\R$, $(x,y)\to F(x,y)$, with the properties

\noindent $(1.2)$ $F(x,y)\geq 0$ with equality if and only if $y=0$,

\noindent $(1.3)$ F is smooth on $\overset{\circ}{TM}$ and only continuous on $TM\setminus\overset{\circ}{TM}$,

\noindent $(1.4)$ $F(x,\lbb y)=\lbb F(x,y)$, $\lbb>0$,

\noindent $(1.5)$ $g_{ij}(x,y)\xi^i\xi^j\geq 0$ with equality if, and only if, $(\xi^i)=0$, where $$g_{ij}(x,y)=\dfrac12\dfrac{\pp^2 F^2}{\pp y^i\pp y^j}.\leqno(1.6)$$

The pair $F^n=(M,F)$ is called a Finsler space. Its geometry is called Finsler geometry.

The geometrical objects from Finsler geometry are in fact living on the sphere bundle $SM\to M$, $SM=TM/\sim$, $(x,y)\sim(x,\wt{y})$ if, and only if, there exists an $a>0$ such that $\wt{y}=ay$. However, for convenience we shall work with the slit tangent bundle $\overset{\circ}{TM}$ instead of $SM$.

The equivalence problem in Finsler geometry is to decide whether two fundamental Finsler functions $F$ and $\wt{F}$ will transform into each other under a diffeomorphism $(x,y)\to(\wt{x},\wt{y})$. In order to solve this problem using E. Cartan's equivalence method, S.S. Chern has introduced in 1948, \cite{1}, a remarkable connection in Finsler geometry by means of some connection $1-$forms. That connection remained outside of the mainstream of the development of Finsler geometry in the next decades. It was only briefly treated in the monograph by H. Rund, \cite{6}, and not at all in those of M. Matsumoto \cite{4}, R. Miron and M. Anastasiei \cite{5}. Chern came back to his connection in 1992, \cite{2}. Then, in a large joint paper with Bao, \cite{3}, its extraordinary usefulness in treating global problems in Finsler geometry was shown.

This fact appeared quite strange to us since along years of study of Finsler geometry and its generalizations we thought of and experienced a mechanism of producing Finsler connections. Thus we decided to see what is the place of Chern's connection among all Finsler connections.

Let $\tau^*TM\to\overset{\circ}{TM}$ be the pull-back of the tangent bundle by $\tau$. An interpretation of Chern's connection as a linear connection in this pull-back bundle has been sketched in \cite{BCh}. We have, however, chosen to relate it to the Cartan nonlinear connection associated to $F$. This allows us to view Chern's connection as a Finsler connection, \cite{M}, or in the terminology from \cite{MA} as a normal $d-$connection. Quite surprisingly we arrived at the Rund connection as defined in \cite{R}, \cite{M}, \cite{MA}. Thus in the famous diagram involving the four remarkable Finsler connection: Berwald, Rund, Cartan, Hashiguchi, \cite{M}, p. 120, Rund's name has to be replaced by Chern's who discovered the connection in question almost ten years earlier. In fact, Rund had a little bad luck with this connection (cf. Remark 18.1 in \cite{M}). These facts do not diminish at all the contribution of Rund and any history of Finsler geometry has to put his name on an outstanding place.

The structure of the paper is as follows. In \S2, we recall Chern's connection $1-$forms. Then in \S3, viewing Chern's connection as a Finsler connection we show that it coincides with the Rund connection.

\section{The Chern connection $1-$form}
We follow \cite{BCh} for recalling the definition and some properties of Chern's
connection $1-$forms.

Set $\pp_i:=\dfrac{\pp}{\pp x^i}$, $\overset{\circ}{\pp}_i:=\dfrac{\pp}{\pp y^i}$.

The homogeneity stipulation (1.4) implies
$$y^i\overset{\circ}{\pp}_i F=F,\leqno(2.1)$$ $$y^i\overset{\circ}{\pp}_i\overset{\circ}{\pp}_iF=0,\leqno(2.2)$$ $$y^i g_{ij}=\dfrac12\overset{\circ}{\pp}_j F^2,\leqno(2.3)$$ $$F^2=g_{ij}y^i y^j,\leqno(2.4)$$ $$y^iC_{ijk}=0,\ C_{ijk}=\dfrac12\overset{\circ}{\pp}_k\overset{\circ}{\pp}_i\overset{\circ}{\pp}_j F^2.\leqno(2.5)$$

By (1.1) and (1.5) it follows that, as a well-defined $(0,2)-$tensor field on $\overset{\circ}{TM}$, $$g=g_{ij}dx^i\otimes dx^j\leqno(2.6)$$ is symmetric and positive definite.

The sections of $\tau^*TM\to\overset{\circ}{TM}$ will be called $\tau-$vector fields or vector fields along $\tau$. Let $\chi(\tau)$ be the set of all $\tau-$vector fields. The fibre of $\tau^*TM\to\overset{\circ}{TM}$ at $u\in\overset{\circ}{TM}$ is $T_{\tau(u)}M$. It has a basis $\left(\dfrac{\pp}{\pp x^i}\right)_{\tau(u)}$ and an inner product given by (2.6). A $\tau-$vector field $\overline{X}\in\chi(\tau)$ is locally given as $\overline{X}=X^i(x,y)\(\dfrac{\pp}{\pp x^i}\)$, the components $(X^i(x,y))$ being smooth functions and transforming under (1.1) as follows $$\wt{X}^i(\wt{x},\wt{y})=\dfrac{\pp\wt{x}^i}{\pp x^k}=\dfrac{\pp\wt{x}^i}{\pp x^k}(x)X^k(x,y).\leqno(2.7)$$ This suggests that we take into consideration $\mathbb{T}=y^i\(\dfrac{\pp}{\pp x^i}\)_{\tau(u)}$ as a remarkable element of $\chi(\tau)$.

By (2.4) and (2.7) one gets $$g(\mathbb{T},\mathbb{T})=F^2\leqno(2.8)$$ i.e. the length of $\mathbb{T}$ is just $F$.

Let $\{e_i\}$ be a local orthonormal (with respect to $g$) frame field for the vector bundle $\tau^*TM\to\overset{\circ}{TM}$ such that $e_n=\dfrac{y^i}{F}\dfrac{\pp}{\pp x^i}$ and $\{w^i\}$ its dual $\oo-$frame. One finds that $\oo^n=(\overset{\circ}{\pp}_iF)dx^i$. Let us set $$\oo^i=v^i_k dx^k,\ e_j=u_j^h\pp_h\leqno(2.9)$$ $$dx^i=u_k^i \oo^k,\ \pp_j=v^h_je_h.\leqno(2.10)$$
These show that $\oo^i$ and $e_j$ can be regarded as $1-$forms and vector fields on $\overset{\circ}{TM}$, respectively.

According to \cite{BCh}, \S2, there exists a set of $1-$forms $\oo^i_j$ on $\overset{\circ}{TM}$ such that $$d\oo^i=\oo^j\wedge\oo^i_j,\leqno(2.11)$$ $$\oo_{ik}+\oo_{ki}:=\oo^j_i\de_{jk}+\oo^j_k\de_{ji}=-H_{ikj}\oo^j_n,\leqno(2.12)$$ $$H_{abc}v^a_iv^b_jv^c_k=F\overset{\circ}{\pp}_k g_{ij}.\leqno(2.13)$$

The $1-$forms $\oo^i_j$ define Chern's connection. We do not write down the fairly complicated expression of these $1-$forms in which the partial derivatives of $F$ are involved. We notice only the following combinations of these partial derivatives which will be used later. $$G_i:=\dfrac14(y^k\overset{\circ}{\pp}_i\pp_k F^2-\pp_iF^2),\leqno(2.14)$$ $$G^i=g^{ik}G_k,\leqno(2.15)$$ $$G^i_j:=\overset{\circ}{\pp}_jG^i.\leqno(2.16)$$ In the structure equation (2.11), $d$ means the exterior differentiation on $TM$.

Let $\Gamma^j_i$ be the representation of $\oo^j_i$ in the natural frame. One defines a covariant differentiation $\nabla$ by $$\nabla\pp_k=\Gamma^i_k\otimes\pp_i,\leqno(2.17)$$ and one proves that $$\Gamma^i_k=\Gamma^i_{kh}dx^h,\ \Gamma^i_{kh}=\Gamma^i_{hk},\leqno(2.18)$$ and with $\Gamma^i_{kh}=g^{ij}\Gamma_{jkh}$ one finds $$\Gamma_{jkh}=\dfrac12(\pp_h g_{jk}-\pp_j g_{kh}+\pp_k g_{hj})+\dfrac12(M_{jkh}-M_{khj}+M_{hjk}),\leqno(2.19)$$ where $$M_{jkh}=-G^t_h\overset{\circ}{\pp}_t g_{jk}.\leqno(2.20)$$

\section{Chern's connection as Finsler connection}

Recall that according to \cite{4}, Ch. I, III, a Finsler connection (a normal $d-$connection in \cite{5}, Ch. VII) is a triad $(N^i_j(x,y),F^i_{jk}(x,y),C^i_{jk}(x,y))$ where $(N^i_j(x,y))$ are the local coefficients of a nonlinear connection, $F^i_{jk}(x,y)$ behave like the coefficients of a linear connection and $C^i_{jk}(x,y)$ are the components of a tensor field. Such a connection is called $h-$metrical if $$g_{ij|k}:=\de_k g_{ij}-F^h_{ik}g_{hj}-F^h_{jk}g_{ih}=0,\leqno(3.1)$$ where $\de_k=\pp_k-N^i_k\overset{\circ}{\pp}$, and $v-$metrical if $$g_{ij}|_k:=\overset{\circ}{\pp}_k g_{ij}-C^h_{ik}g_{hj}-C^h_{jk}g_{ih}=0.\leqno(3.2)$$

Now we shall regard Chern's connection as a Finsler connection showing that it is a $h-$metrical one. First, we re-express $\Gamma^i_{kj}$ as follows. Considering $\de_i=\pp_i-G^j_i\overset{\circ}{\pp}_j$ we observe that $\de_i g_{kh}=\pp_i g_{kh}+M_{khi}$. Inserting this in (2.19) we find $$\Gamma_{jkh}=\dfrac12(\de_k g_{jh}+\de_h g_{jk}-\de_j g_{kh}).\leqno(3.3)$$
Now we must check that $(\Gamma^i_{kh})$ behave like $(F^i_{kh})$ under (1.1). But we can avoid this complicated calculation as we shall see below.

For the covariant differentiation of $g$ with respect to Chern's connection we have $(\nabla g)(\pp_i,\pp_j)=dg_{ij}-\Gamma_i^k g_{kj}-\Gamma^k_j g_{ik}=dg_{ij}-(\Gamma^k_{ih}g_{kj}+\Gamma^k_{jh}g_{ik})dx^h$ and using (3.1) we find $(\nabla g)(\pp_i,\pp_j)=(\overset{\circ}{\pp}_h g_{ij})(dy^h+G^h_s dx^s)$.

We put $\de y^h=dy^h+G^h_k dx^k$ and it is easily checked that $\de y^h(\de_k)=0$.

Going back to the above formulae we conclude that Chern's connection satisfies $$dg_{ij}=\Gamma^k_i g_{kj}+\Gamma^k_j g_{ik}+2C_{kij}\de y^k.\leqno(3.4)$$ We note also that from $(\nabla g)(\pp_i,\pp_j)=(\overset{\circ}{\pp}_h g_{ij})\de y^h$ it results that $\nabla$ is metrical only for those tangent vectors $v$ which verify $\de y^h(v)=0$. Recall that for $C_{ijk}=0$, $F^n$ reduces to a Riemannian space.

This fact motivates us to introduce the following:

{\bf Definition 3.1.} A tangent vector $X_n\in T_u\overset{\circ}{TM}$ is said to be horizontal if $\de y^h(X_u)=0$.

Thus $\nabla$ is metrical along horizontal vectors, in particular along the $\de_i$'s and on the subspace spanned by them, called the horizontal subspace of $T_u\overset{\circ}{TM}$.

The significance of (3.3) is underlined by

{\bf Proposition 3.1.} {\it There exists a unique set of $1-$forms $\{\Gamma^i_j\}$ on $\overset{\circ}{TM}$ satisfying

$(a)$ $d(dx^i)=dx^j\wedge\Gamma^i_j$,

$(b)$ $dg_{ij}=\Gamma^k_i g_{kj}+\Gamma^k_j g_{ik}+2C_{hij}\de y^h$.}

\noindent{\it Proof.} The existence was proved in the above. Let $\wt{\Gamma}^i_j=\wt{\Gamma}^i_{jk} dx^k+\widehat{\Gamma}^i_{jk}dy^k$ be $1-$forms satisfying (a) and (b). From $dx^j\wedge\wt{\Gamma}^i_j=0$ it follows that $\wt{\Gamma}^i_{jk}=\wt{\Gamma}^i_{kj}$ and $\widehat{\Gamma}^i_{jk}=0$. Subtracting member by member the equations (b) for $\Gamma$'s and $\widehat{\Gamma}$'s one obtains $(\wt{\Gamma}^s_{ik}-\Gamma^s_{ik})g_{sj}+(\wt{\Gamma}^s_{jk}-\Gamma^s_{jk})g_{sj}=0$. Permuting cyclicly the indices $i,j,k$ one gets two new equations which added and subtracting from the result the previous one gives $(\wt{\Gamma}^s_{ij}-\Gamma^s_{ij})g_{sk}=0$. Hence $\wt{\Gamma}^s_{ij}=\Gamma^s_{ij}$, q.e.d.

{\it Remark 3.1.} As $(\pp_i,\overset{\circ}{\pp}_i)$ is the natural frame on $\overset{\circ}{TM}$, (2.17) is equivalent to $$\begin{array}{l}\nabla_{\pp_j}\pp_i=\Gamma^k_{ij}\pp_k,\\ \\ \nabla_{\overset{\circ}{\pp}_j}\pp_i=0.\end{array}\leqno(3.5)$$

Calculating (3.4) for $(\overset{\circ}{\pp}_h)$ one finds $$g_{ij}|_h:=(\nabla_{\overset{\circ}{\pp}_h}g)(\pp_i,\pp_j)=2C_{ijh}.\leqno(3.6)$$

In Finsler geometry there exist four remarkable Finsler connections which have in common the Cartan nonlinear connection of coefficients $(\overset{c}{N}^i_j)$. Among them we have the Rund connection which has the form $(\overset{c}{N}^i_j(x,y)$, $F^i_{jk}(x,y),0)$ with the coefficients $F^i_{jk}(x,y)$ given by $$F^i_{jk}=\dfrac12 g^{ih}(\overset{c}{\de}_j g_{hk}+\overset{c}{\de}_h g_{jk}),\leqno(3.7)$$ where $\overset{c}{\de}_j=\pp_j-\overset{c}{N}^k_j(x,y)\overset{\circ}{\pp}_k$.

This connection is $h-$metrical but it is not $v-$metrical since by (3.2) we have $$g_{ij}|_k=2C_{ijk}\neq 0,\leqno(3.6')$$ except when $F^n$ is a Riemannian space.

Looking at Chern's connection we see that the $\Gamma$'s from (3.3) coincide with the $F$'s from (3.7) if the $(G^i_j)$ given by (2.16) are just the $(\overset{c}{N}^i_j(x,y))$ of Cartan.

This indeed holds as we now prove.

Let $\gamma^i_{jk}(x,y)$ be the ``Christoffel symbols'' $$\gamma^i_{jk}=\dfrac12 g^{ih}(\pp_j g_{hk}+\pp_k g_{jk}-\pp_h g_{jk}).$$

Then the coefficients of the Cartan nonlinear connection are $$\overset{c}{N}^i_j=\overset{\circ}{\pp}_j(\Gamma^i_{kh}y^k y^h).\leqno(3.8)$$ By (2.16) it is sufficient to check that $2G^i=\gamma^i_{kh}y^ky^h$. Equivalently, $$4G_i=(\pp_j g_{ik}+\pp_k g_{ji}-\pp g_{jk})y^jy^k.\leqno(3.9)$$ By (2.14), $4G_i=y^k\overset{\circ}{\pp}_i\pp_k(F^2)-\pp_iF^2$. Using (2.3) and (2.4), the righthand side of (3.9) becomes $$2\pp_j(g({ik}y^k)y^j-\pp_i(g_{jk}y^j y^k)=\pp_j(\overset{\circ}{\pp}_iF^2)y^j-F^2.$$ Hence (3.9) holds.

The equalities $G^i_j=\overset{c}{N}^i_j$, $\Gamma^i_{jk}=F^i_{jk}$, and (3.6) in conjunction with (3.6)', show that we may think of Chern's connection as the Finsler connection $(G^i_j,\Gamma^i_{jk},0)$ and furthermore it coincides with the Rund connection.

Since this Finsler connection was first introduced by Chern, it is quite natural that it bear his name. However, Chern has rather graciously suggested that it be called the Chern-Rund connection.


\newpage

\def\rr{\mathbb{R}}

\runningauthor={M. ANASTASIEI}
\runningtitle={FINSLER CONNECTIONS IN GENERALIZED LAGRANGE SPACES}
\noindent
\baselineskip 8pt
\noindent{\footnotesize{Balkan J. Geom. Appl.\hfill\break 1 (1996), no. 1, 1-9}}
\vskip 2cm
\baselineskip 11.5pt plus .15pt
\centerline{\bf\Large FINSLER CONNECTIONS}
\vskip .2cm
\centerline{\bf\Large IN GENERALIZED LAGRANGE SPACES}
\vskip .5cm
\centerline{\bf by Mihai ANASTASIEI}
\vskip 1cm

\begin{abstract}
The Chern--Rund connection  from Finsler geometry is settled in the generalized Lagrange
spaces. For the geometry of these spaces, we refer to [5].

{\bf AMS Subject Classification:} 53C60.

{\bf Key words and phrases:} Finsler connections, generalized
Lagrange spaces, Chern--Rund connection.
\end{abstract}

\setcounter{section}{0}
\section*{Introduction}
In a recent paper, \cite{1}, we showed that in a Finsler
space the connection introduced by S.S. Chern in 1948 is the
same with the connection proposed by H. Rund ten years later
and bearing his name. Accordingly, we proposed the name of
Rund be replaced with that of Chern, but several geometers
including S.S. Chern himself, suggested to call it from now
on a Chern--Rund connection.

As S.S. Chern and D. Bao showed in \cite{2}, the
Chern--Rund connection is very convenient in treating of many
global problems in Finsler geometry. This fact determined us to
come back to the subject.

The efforts made in defining a covariant derivative and
accordingly, a parallel displacement in Finsler space led
to a concept generically called a Finsler connection. Among the
Finsler connections there exist four, which are remarkable by
their properties named the Cartan, Berwald, Chern--Rund and
Hashiguchi connections, respectively. These are usually put
together in a nice commutative diagram (cf. \cite[Ch. III]{3}).

The most utilized is the Cartan connection because it is fully
metrical i.e. $h-$ and $v-$metrical, in spite of the fact
it has torsion.

But there are some problems involving the Berwald connection
which is by no means metrical or the Hashiguchi connection which is only $v-$metrical.

The Chern--Rund connection being $h-$metrical and free of
torsion is the nearest to the Levi--Civita connection a fact
which explains its adequacy for global problems in Finsler
geometry.

The Finsler connections are also suitable for the geometries more
general than the Finslerian one as the Lagrange geometry or generalized Lagrange geometry. Our purpose is to review Finsler
connections and to settle the Chern--Rund connection in this more general framework.

First, we give in \S1 a definition of Finsler connection by
local components and introduce its compatibility with a
generalized Lagrange metric. Then, in \S2, a Finsler
connection is defined as a pair $(N,\nabla)$, where $N$ is a
nonlinear connection on $TM$ and $\nabla$ is a linear
connection in the pull--back bundle $c^{-1}TM\longrightarrow TM$ with
$\tau:TM\longrightarrow TM$, the tangent bundle over a manifold $M$. These
definitions are equivalent. The four remarkable connections
mentioned above are characterized. A special attention is
paid to the possibility of determining $N$ from $\nabla$.\bigskip

\noindent{\bf Acknowledgement.} We are indebted to Prof. Dr. Radu
Miron who suggested us several improvements of the first
version of this paper.

\section{Finsler connections. A definition by local components}
\setcounter{equation}{0}

Let $M$ be a smooth i.e. $C^\9$ manifold of finite dimension $n$
and $\tau:TM\to M$ its tangent bundle. A local chart
$(U,(x^i))$ on $M$ induces a local chart $(\tau^{-1}(U),(x^i,y^i))$ on $TM$, where $x^i\equiv
x^i\circ\tau$ and $(y^i)$ are provided by $u=y^i\left.\dfrac{\pp}{\pp x^i}\right|_p,\ p=\tau(u)$.

A change of coordinates $(x^i,y^i)\longrightarrow(\tilde x^i,\tilde y^i)$ on $TM$ has the form
\begin{equation}
\begin{array}{l}
\tilde x^i=\tilde x^i(x^1,...,x^n),\ {\rm
rank}\,\left(\dd\frac{\pp\tilde x^i}{\pp x^j}\right)=n\vspace{1mm}\\
\tilde y^i=\dd\frac{\pp\tilde x^i}{\pp x^j}(x)y^j.
\end{array}
\end{equation}
The indices $i,j,k,...,$ will run from $1$ to $n$ and
Einstein's convention on summation is implied.

Let $L:TM\longrightarrow R$ be a scalar function on $TM$. Then $\wt L(\tilde
x(x),\tilde y(y))=L(x,y)$, from which, taking partial
derivatives and using (1.1), one gets
\begin{equation}
\dd\frac{\pp L}{\pp y^i}=\dd\frac{\pp\tilde x^k}{\pp x^i}\
\dd\frac{\pp L}{\pp \tilde y^k},
\end{equation}
\begin{equation}
\dd\frac{\pp L}{\pp x^i}=\dd\frac{\pp\tilde x^k}{\pp x^i}\
\dd\frac{\pp L}{\pp \tilde x^k}+\frac{\pp^2\tilde x^k}{\pp x^j\pp
x^i}y^j \frac{\pp L}{\pp \tilde y^k}\cdot
\end{equation}
According to (1.2), the set of functions $\left(\dd\frac{\pp
L}{\pp y^i}(x,y)\right)$ may be regarded as the components
of a covector field on $TM$. From (1.2), it follows that
$\left(\dd\frac{\pp^2L}{\pp y^i\pp y^j}(x,y)\right)$ may be
also viewed as the components of a (symmetric) tensor field
on $TM$. Thus on $TM$ there exist geometrical objects whose law of
transformation under (1.1) is the same as of the corresponding objects on $M$. These were called $d-$objects ($d$ is from
distinguished) in \cite{5}, Finsler objects in \cite{3} and
sometimes $M-$objects.

The geometry of $d-$objects is essentially involved in the
study of those metrical structures which are more general than
Riemannian structures i.e. Finsler structures, Lagrange structures, generalized
Lagrange structures (see \cite{5}).

Coming back to (1.3), we see that the behavior of the operators $\dd\frac{\pp}{\pp x^i}$ is drastically different from
that of $\dd\frac{\pp}{\pp y^i}\cdot$ Let us introduce a
correction of $\dd\frac{\pp}{\pp x^i}:=\pp_i$,
\begin{equation}
\delta_iL=\pp_iL+N^k_i(x,y)\dot\pp_k,\ \dot\pp_{\dot
k}:=\frac{\pp}{\pp y^{\dot k}},
\end{equation}
such that, with respect to (1.1):
\begin{equation}
\delta_iL=\frac{\pp\tilde x^k}{\pp x^i}\,\tilde\delta_kL,
\end{equation}
i.e. $(\delta_iL)$ to appear as the components of a covector field on $TM$. Then the functions
$(N^k_i(x,y))$ have to satisfy
\begin{equation}
\dfrac{\pp\tilde x^j}{\pp x^i}\,\tilde N^h_j=N^j_i\,\dfrac{\pp\tilde
x^h}{\pp x^j}+\dfrac{\pp^2\tilde x^h}{\pp x^i\pp
x^j}\,y^j.
\end{equation}
Note that $(N^j_i(x,y))$ are not the components of a
$(1,1)$--tensor field on $TM$ but the difference of two
sets of this type is so.

As it is well--known, when $M$ is paracompact, there exists on $M$
a linear connection, say of local coefficients $(\Gamma^i_{jk}(x))$. Then
$N^i_k(x,y)=\Gamma^i_{jk}(x)y^j$ verify (1.6). This example
assures also the existence of a nonlinear connection within a
generally accepted hypothesis on $M$.

The local vector fields $(\delta_i)$, $i=1,2,...,n$, given by
(1.4) are linearly independent and in a point $u\in TM$
they span an $n-$dimensional subspace $H_uTM$ of $T_uTM$.

Let $\tau_{*,u}$ be the tangent mapping (the Jacobian) of
$\tau$. Then $V_uTM=\ker\tau_{*,u}$ is called the vertical
subspace of $T_uTM$. A vertical vector is of the form
$X^k(x,y)\dot\pp_k$ such that under (1.1) one has
\begin{equation}
\wt X^k=\dfrac{\pp x^k}{\pp x^i}\,X^i.
\end{equation}
We immediately have
\begin{equation}
T_uTM=V_uTM\oplus H_uTM.
\end{equation}
Furthermore, $\tau_{*,u}$ restricted to $H_uTM$ gives an
isomorphism of it with $T_{\tau(u)}M$ such that
$\tau_{*,u}(\delta_i)=\pp_i\big|_{\tau(u)}$.

Conversely, if a supplement of $V_uTM$ in $T_uTM$ is
specified by a basis $(\delta_i)$, $i=1,2,...,n$, which is
carried by $\tau_*$ to $(\pp_i)$, then letting $\delta_i=\pp_i-N^k_i\dot\pp_i$, the condition
$\delta_i=\dd\frac{\pp\tilde x^k}{\pp x^i}\wt\delta_k$ implies (1.6)
for $(N^k_i)$. One says that $(N^k_i(x,y))$ are the coefficients
of a nonlinear connection.

A reason for this term is that when $(N^k_i)$ are linear
with respect to $(y)$ i.e. $N^k_i(x,y)=G^k_{ji}(x)y^j$,
then $(G^k_{ji})$ are the coefficients of a linear connection on $M$.

Summarizing the foregoing discussion we may formulate the
following two equivalent definitions for a nonlinear connections.

\ms

{\bf Definition 1.1.} A {\it nonlinear connection} is a set of functions
$(N^i_j(x,y))$ defined on each domain of local chart on
$TM$ such that an overlaps, (1.6) holds good.

\ms

{\bf Definition 1.2.} A {\it nonlinear connection} is a smooth distribution $u\longrightarrow H_uTM$ supplementary to the
vertical distribution $u\longrightarrow V_uTM$ i.e. (1.7) holds good for every $u\in TM$.

\ms

Let $(v^i(x,y))$ be the components of a $d-$vector field.
Then $\left(\dd\frac{\pp v^i}{\pp y^j}(x,y)\right)$ are the
components of a $d-$tensor field of type $(1,1)$. In other
words the partial derivatives with respect to $(y^i)$ are
covariant. However, in some circumstances, these have to be
replaced by
\begin{equation}
v^i{}_{|j}=\dfrac{\pp v^i}{\pp y^j}+C^i_{kj}(x,y)v^k,
\end{equation}
where $(C^i_{kj}(x,y))$ are the components of a $d-$tensor field. One
of them is as follows.

First, we introduce

\ms

{\bf Definition 1.3.} A $d-$tensor field of type $(0,2)$ of
components $(g_{ij}(x,y))$ which is

a) symmetric, i.e. $g_{ij}=g_{ji}$,

b) nondegenerate i.e. $\det(g_{ij}(x,y))\ne0$ and

c) the quadratic form $g_{ij}(x,y)\xi^ixi^j$ $(\xi\in\rr^n)$

\noindent has constant signature is called a {\it generalized Lagrange metric} ($GL-$metric
for brevity).

Extending (1.9), the covariant derivative of $(g_{ij})$ is given by
\begin{equation}
g_{ij\big|k}=\pp_jv^i-C^h_{ik}g_{hj}-C^h_{jk}g_{ih}.
\end{equation}
One says that the $GL-$metric $(g_{ij}(x,y))$ is $v-$covariant constant if $g_{ijg|k}=0$.

For the general $v^i_{|j}$, the condition $g_{ij|k}=0$ can be fulfilled with
\begin{equation}
\stackrel{c}{C}\!{}^h_{ij}=\dfrac{1}{2}g^{hk}(\dot\pp_i g_{kj}+\dot\pp_j g_{ik}-\dot\pp_k g_{ij}).
\end{equation}
The partial derivatives with respect to $(x^i)$ are far to be covariant derivatives. A correction of them could be
$\pp_jv^i+H^i_{kj}(x,y)v^k$, but $(H^i_{kj}(x,y))$ have a complicated law of transformation $A$ better one is
\begin{equation}
v^i_{|j}=\delta_jv^i+F^i_{kj}(x,y)v^k,
\end{equation}
since then $(F^i_{kj}(x,y))$ changes under (1.1) as the local coefficients of a linear connection on $M$. These derivatives
can be extended to any $d-$tensor field. For instance, the $v-$covariant derivative of $(g_{ij}(x,y))$ is given by
(1.10) and its $h-$covariant derivative is
\begin{equation}
g_{ij|k}=\delta_kg_{ij}-F^h_{ik}g_{hj}-F^h_{jk}g_{ih}.
\end{equation}

The $GL-$metric $(g_{ij}(x,y))$ is said to be $h-$covariant constant if $g_{ij|h}=0$. It is easy to check that the equation
$g_{ij|h}=0$ is satisfied with
\begin{equation}
\stackrel{c}{F}\!{}^k_{ij}=\dfrac12 g^{kh}(\delta_ig_{hj}+\delta_jg_{ih}-\delta_hg_{ij}).
\end{equation}

The foregoing discussions suggest

\ms

{\bf Definition 1.4} A Finsler connection is a triad $F\Gamma=(N^i_j(x,y)$, $F^i_{jk}(x,y)$, $C^i_{jk}(x,y))$, where
$N^i_j(x,y)$ are the coefficients of a nonlinear connection, $F^i_{jk}(x,y)$ are like the coefficients of a linear connection on $M$
and $C^i_{jk}(x,y)$ are the components of a $d-$tensor field.

\ms

We have also got a first example of Finsler connection $C\Gamma=(N^i_j(x,y)$, $\stackrel{c}{F}\!{}^i_{jk}(x,y)$,
$\stackrel{c}{C}\!{}^i_{jk}(x,y))$.

\ms

{\bf Definition 1.5} Let $F\Gamma$ be a Finsler connection and $(g_{ij}(x,y))$ a $GL-$metric. $F\Gamma$ is said to be
$h-$metrical if $g_{ij|h}=0$, $v-$metrical if $g_{ij|h}=0$ and metrical if the both equations hold.

\ms

In the above we have proved

\ms

{\bf Proposition 1.1} {\it The Finsler connection $C\Gamma$ is metrical.}

\ms

The following $d-$tensor fields are called the torsions of $F\Gamma$:
\begin{equation}
\begin{array}{ll}
T^i_{jk}=F^i_{jk}-F^i_{kj},& R^i_{jk}=\delta_k N^i_j-\delta_j N^i_k,C^i_{jk},\vspace{1mm}\\
P^i_{jk}=\dot\pp_kN^i_j-F^i_{kj},& S^i_{jk}=C^i_{jk}-C^i_{kj}.
\end{array}
\end{equation}

\noindent{\it Remark.} $R^i_{jk}$ is the integrability
tensor of the horizontal distribution. It measures also the
curvature of the nonlinear connection $N$.

The $d-$tensor fields
\begin{equation}
D^i_j=F^i_{kj} y^k-N^i_j,\ d^i_j=\delta^i_j+C^i_{kj} y^k,
\end{equation}
where $(\delta^i_j)$ is Kronecker' symbol, are called
$h-$deflection and $v-$deflection of $F\Gamma$, respectively.

From (1.6) we infer that $G^i_{jk}=\dot\pp_j N^i_k$ transform
under (1.1) as $F^i_{jk}$. Thus $B\Gamma=(N^i_j,G^i_{jk},0)$ is a
Finsler connection. It will be called the Berwald connection. This
connection is no $v-$metrical nor $h-$metrical and is free of
torsions if and only if $N$ is integrable $(R^i_{jk}=0)$ and symmetric
$(\dot\pp_jN^i_k=\dot\pp_kN^i_j)$.

The connection  $C\Gamma$ will be called the Cartan connection. It is
$h-$metrical, $h-$symmetric $(\stackrel{c}{F}\!{}^i_{jk}(x,y)=\stackrel{c}{F}^i_{kj}(x,y))$,
$v-$metrical and $v-$symmetric. The Finsler connection $H\Gamma=(N^i_j,G^i_{jk}(x,y)$, $\stackrel{c}{C}^i_{kj}(x,y))$
will be called the Hashiguchi connection. This is $v-$metrical, no $h-$metrical and has torsion. The Finsler connection
$CR\Gamma=(N^i_j,\stackrel{c}{F}^i_{jk}(x,y),0)$ will be called the Chern--Rund connection. This is $h-$metrical but not
$v-$metrical.

Summarizing, for a fixed nonlinear connection $N$ and a
$GL-$metric $(g_{ij}(x,$ $y))$ we have four typical Finsler
connections: $B\Gamma$, $C\Gamma$, $H\Gamma$ and $CR\Gamma$.

Let us replace $TM$ by $T_0M=TM\setminus0$.

A $GL-$metric $(g_{ij}(x,y))$ on $T_0M$ reduces to a Finsler
metric if there exists a fundamental Finsler function $F:T_0M\longrightarrow R_+$
such that $g_{ij}(x,y)=\dfrac12\dot\pp_i\dot\pp_jF^2(x,y)$. Taking as
$N$ the Cartan nonlinear connection of coefficients $\stackrel{c}{N}^i_j=\dfrac12\dot\pp_j\g^i_{oo}$,
$\g^i_{oo}=\g^i_{jk} y^jy^k,$ $\g^i_{jk}=\dfrac12g^{ih}(\pp_jg_{hk}+\pp_kg_{jh}-\pp_hg_{jk})$,
the afore mentioned Finsler connections reduce to the four remarkable connections in Finsler geometry (\cite[Ch. III]{3}).

The form of $D^i_j$ in (1.6) shows that one may associate
to any $F\Gamma$ a new Finsler connection $(F^i_{kj} y^k-D^i_j$, $F^i_{kj}$, $C^i_{kj})$ whose $h-$deflection is just
$D^i_j$, when this is prescribed. In particular, for $D^i_j=0$ a Finsler connection  without $h-$deflection is
obtained. In Finsler geometry $B\Gamma$, $C\Gamma$, $H\Gamma$ and $CR\Gamma$ are $h-$deflection free. So we have an explanation why the nonlinear connection was noted quite late in Finsler geometry.

\section{Another definition of Finsler connections}
\setcounter{equation}{0}

Let be $\tau^{-1}TM=\{(u,v)\in TM\times TM,\tau(u)=\tau(v)\}$
fibered over $TM$ by $\pi(u,v)=u$. The local fiber in
$(u,v)$ is $T_{\tau(u)}M$. A section in $(\tau^{-1}TM,\pi,TM)$
is locally of the form $\bar X=\bar X^i(x,y)\bar\pp_i$ with
$(\bar\pp_i)$ the natural basis in $T_{\tau(u)}M$. It
follows that under (1.1) we have
\begin{equation}
\wt{\bar X}{}^i=\dfrac{\pp\tilde x^i}{\pp x^k}\,\bar X^k.
\end{equation}
$\bar X$ will be called a $\tau-$vector field on $TM$. It
can be identified with the $d-$vector field $(\bar
X^i(x,y))$. More general, the tensorial algebra of the
pull--back bundle $\tau^{-1}TM$ can be thought of as algebra
of $d-$tensor fields on $TM$. There exists a remarkable $\tau-$vector
field ${\mathbb{C}}:u\longrightarrow(u,u)$, which locally is
$y^i\bar\pp_i$ and so it can be identified to the Liouville
vector field $\mathbb{C}=y^i\dot\pp_i$.

\ms

{\bf Theorem 2.1.} {\it There exists a one--to--one correspondence between the set of Finsler connections $F\Gamma$ and
the set of pairs $(N,\nabla)$ with $N$ a nonlinear connection on $TM$ and $\nabla$ a linear connection in the pull--back
bundle $\tau^{-1}TM$.}

\ms

\noindent{\it Proof.} If $F\Gamma$ is specified by $(N^i_j,F^i_{jk},C^i_{jk})$, we take
$N=(N^i_j)$ and define $\nabla$ by
\begin{equation}
\nabla_{\delta_k}\bar\pp_i=F^i_{jk}\bar\pp_i,\
\nabla_{\dot\pp_k}\bar\pp_j=C^i_{jk}\bar\pp_i.
\end{equation}
In the natural basis $\nabla$ takes the form
\begin{equation}
\nabla_{\pp_k}\bar\pp_j=\Gamma^i_{jk}\bar\pp_i,\ \
\nabla_{\dot\pp_k}\bar\pp_i=C^i_{jk}\bar\pp_i.
\end{equation}
\begin{equation}
\Gamma^i_{jk}=F^i_{jk}+N^h_kC^i_{jh}.
\end{equation}

Conversely, given $N=(N^i_j)$ and $\nabla$ specified by
(2.3) it results that $(N^i_j$, $F^i_{jk}$, $C^i_{jk})$ with $F^i_{jk}$
given by (2.4) is a Finsler connection. \hfill$\Box$

\ms

A $GL$--metric $(g_{ij}(x,y))$ defines a metrical structure $g$ in
the bundle $\tau^{-1}TM$:
\begin{equation}
g=gij(x,y)dx^i\otimes dx^j.
\end{equation}

Conversely, any metrical structure in the bundle $\tau^{-1}TM$
defines by (2.5) a $GL-$metric.

One easily checks

\ms

{\bf Theorem 2.2} {\it In the correspondence $F\Gamma\longleftrightarrow(N,\nabla)$ we have
\begin{itemize}
\item[{\rm a)}] $F\Gamma$ is h--metrical if and only if $\nabla_{hX}g=0,$
\item[{\rm b)}] $F\Gamma$ is v--metrical if and only if $\nabla_{vX}g=0,$
\item[{\rm c)}] $F\Gamma$ is metrical if and only if $\nabla_{X}g=0,$\\
for every $X\in{\cal X}(TM)$.
\end{itemize}}

\ms

Let $\rho:TTM\longrightarrow\tau^{-1}TM$ be the morphism of vector
bundles given by $\rho(X_u)=(u,\tau_{*,u}X_u)$, $X_u=T_uTM,$ $u\in TM$. It follows that $\ker\rho_u=V_uTM$
i.e. $\rho(\dot\pp_i)=0$ and $\rho(\delta_i)=\bar\pp_i$. Alternatively, we may define a morphism
$\sigma:TTM\longrightarrow\tau^{-1}TM$ on basis by $\sigma(\delta_i)=0$,
$\sigma(\dot\pp_i)=\bar\pp_i$. We say that
\begin{equation}
\begin{array}{l}
\mathbb{T}_\rho(X,Y)=\nabla_X\rho(Y)-\nabla_Y\rho(X)-\rho[X,Y],\vspace{1mm}\\
\mathbb{T}_\sigma(X,Y)=\nabla_X\sigma(Y)-\nabla_Y\sigma(X)-\sigma[X,Y],\
X,Y\in{\cal X}(TM),
\end{array}
\end{equation}
are torsions of $\nabla$.

The following characterizations of the Finsler connections
$B\Gamma$, $H\Gamma$, $CR\Gamma$ and $C\Gamma$ follow.

\ms

{\bf Theorem 2.3.} {\it In the correspondence $F\Gamma\longleftrightarrow(N,\nabla)$ we have
\begin{itemize}
\item[{\rm a)}] $B\Gamma\longleftrightarrow(N,\nabla)$ with
$\mathbb{T}_\sigma(hX,vY)=0,\ \mathbb{T}_\rho(hX,vY)=0$;
\item[{\rm b)}] $H\Gamma\longleftrightarrow(N,\nabla)$ with
$\mathbb{T}_\sigma(hX,vY)=0,\ \mathbb{T}_\sigma(vX,vY)=0,\ \nabla_{vX}g=0$;
\item[{\rm c)}] $CR\Gamma\longleftrightarrow(N,\nabla)$ with
$\mathbb{T}_\rho(hX,vY)=0,\ \mathbb{T}_\rho(hX,hY)=0,\ \nabla_{hX}g=0$;
\item[{\rm d)}] $C\Gamma\longleftrightarrow(N,\nabla)$ with
$\mathbb{T}_\rho(hX,vY)=0,\ \mathbb{T}_\sigma(vX,vY)=0,\
\nabla_{X}g=0$.
\end{itemize}}

\ms

\noindent{\it Proof.} The local expressions of $\mathbb{T}_\rho$ and $\mathbb{T}_\sigma$ in
conjunction with Theorem 2.2 give the desired
results.\bigskip

Now the following question appears. Which conditions should satisfy
$\nabla$ in order to determine $N$ such that the pair $(N,\nabla)$ to correspond to a Finsler connection. An answer is as
follows.

\ms

{\bf Definition 2.1.} A linear connection $\nabla$ in the
pull--back bundle $\tau^{-1}TM$ is said to be regular if the
subspace $\{X_u\mid\nabla_{X_u}{\mathbb{C}}=0$, $X\in{\cal X}(TM)\}$
of $T_uTM$ is supplementary to $V_uTM$ for every $u\in TM$.

\ms

By the definition, every regular connection $\nabla$ induces a
nonlinear connection $N$ on $TM$. The pair $(N,\nabla)$, as we have
seen before, corresponds to a Finsler connection $F\Gamma$. This
$F\Gamma$ has to be of a particular form. Indeed, one has

\ms

{\bf Theorem 2.4.} {\it There exists a bijection between
the set of regular connections in $\tau^{-1}TM$ and the set of
Finsler connections $F\Gamma=(N^i_j$, $F^i_{jk}$, $C^i_{jk})$ satisfying
$D^i_j=0$ and $\det(d^i_h)\ne0$.}

\ms

\noindent{\it Proof.} Let $\nabla$ be specified by (2.3). Using $N=(N^i_j)$
provided by the regularity of $\nabla$, we define $F^i_{jk}$ as
in (2.4). Then $0=\nabla_{\delta_k}{\mathbb{C}}=(y^jF^i_{jk}-N^i_k)\bar\pp_k$
implies $D^i_k=0$. Contracting (2.4) by $y^j$ we get
$N^h_k(d^i_h)=y^j\Gamma^i_{jk}$ and as $(N^h_k)$ is specified this
equation has to have an unique solution. Hence with necessity
$\det(d^i_h)\ne0$.

Conversely, let $(N,\nabla)$ be in correspondence with $F\Gamma$.
The condition $D^i_j=0$ assures that the subspace $\{X_u|\nabla_{X_u}{\mathbb{C}}=0$, $X\in{\cal X}(TM)$, $u\in
TM\}$ is contained in the horizontal subspace $H_uTM$ of
$N$. The condition $\det(d^i_k)\ne0$ implies that this subspace
is supplementary to $V_uTM$. Thus $\nabla$ is regular and the
nonlinear connection derived from it coincides with $N$.\hfill$\Box$

Let us assume that $(g_{ij})$ reduces to a Finsler metric on
$T_0M$. Then $C\Gamma$ is characterized by the following Matsumoto's
axioms:
$$T^i_{jk}=0,\ g_{ij|k}=0,\ S^i_{jk}=0,\ g_{ij\big|k}=0,\
D^i_j=0.\leqno(*)$$
It results $d^i_j=\delta^i_j$.

Combining these with Theorems 2.4 and 2.3, one
obtains

\ms

{\bf Theorem 2.5.} {\it Let $F^n=(M,F)$ be a Finsler space.
There exists a unique regular connection $\nabla$ in $\pi^{-1}T_0M$
satisfying the conditions:
$$\mathbb{T}_\rho(hX,hY)=0,\ \mathbb{T}_\sigma(vX,vY)=0,\ \nabla_Xg=0,\
X,Y\in{\cal X}(T_0M)$$ where $h$ and $v$ are projectors of $N$
induced by $\nabla$.}

\ms

We note that $\nabla$ is determined by $F$ only.

According to \cite{5} the Chern--Rund connection in a Finsler
space is cha\-rac\-te\-ri\-zed by the following axioms:
$$T^i_{jk}=0,\ g_{ij|k}=0,\ C^i_{jk}=0,\ D^i_j=0.$$
We have again $d^i_j=\delta^i_j$. By Theorems 2.3 and 2.4 we have

\ms

{\bf Theorem 2.6.} {\it Let $F^n=(M,F)$ be a Finsler space.
There exists a unique regular connection $\nabla$ in $\pi^{-1}T_0M$
satisfying the conditions:
$$\mathbb{T}_\rho(hX,hY)=0,\ \mathbb{T}_\rho(hX,vY)=0,\ \nabla_{hX}g=0,\
X,Y\in{\cal X}(T_0M)$$ where $h$ and $v$ are projectors of $N$
induced by $\nabla$.}

\ms

The systems of axioms for $H\Gamma$ and $B\Gamma$ discussed for
minimality in \cite{5} give similar results in view of
Theorems 2.3 and 2.4.

The Finsler connections may be viewed also as special liner
connections on $TM$ or in the Finsler bundle $\pi^{-1}LM$, where
$LM$ is the principal bundle of linear frames on $M$. We
refer to \cite{5} and \cite{3}, respectively.




\newpage

\runningauthor={M. ANASTASIEI}
\runningtitle={JACOBI FIELDS IN GENERALIZED LAGRANGE SPACES}
\noindent
\baselineskip 8pt
\noindent{\footnotesize{Rev. Roum. Math. Pures Appl.}}
\hfill\break
{\footnotesize{ 42 (1997), 9-10, 689-695}}
\hfill\break
{\footnotesize{ Collection of papers in honor of Academician Radu Miron on his 70th birthday}}
\vskip 2cm
\baselineskip 11.5pt plus .15pt
\centerline{\bf\Large JACOBI FIELDS IN}
\vskip .2cm
\centerline{\bf\Large GENERALIZED LAGRANGE SPACES}
\vskip .5cm
\centerline{\bf {\footnotesize{BY}}}
\vskip .5cm
\centerline{\bf {\footnotesize{M. ANASTASIEI and I. BUCATARU}}}
\vskip 1cm

\begin{abstract}
We consider the first variation of those curves on tangent manifold $TM$ which have property that are parallel with respect to the canonical metrical connection in a generalized Lagrange space. Accordingly we introduce and study the Jacobi fields on $TM$. Several particular cases are discussed.
\end{abstract}

\setcounter{section}{0}
\section{Introduction}

Among the notions introduced and studied by Prof.Radu Miron, very interesting and useful for applications is that of generalized Lagrange space, $GL$-space for brevity. This is a pair made up by a smooth manifold $M$ and a generalized Lagrange metric, shortly a $GL$-metric. Roughly speaking a $GL$-metric is a metrical structure in the vertical bundle over the manifold $TM$. Viewing in local coordinates one can see that its definition was tailored after the basic properties of a Finsler metric. Thus a $GL$-space appears as a very large generalization of a Finsler space. However, this notion preserves many properties of a Finsler space, the existence of a canonical metrical connection being an important one. The autoparallel curves of this connection are remarkable since in the Finslerian framework these are projecting on the geodesics of the Finsler metric. Calling then also geodesics, we consider their first variations, in Section 3, and accordingly we find a Jacobi equation whose solutions are called Jacobi fields. Some properties of these are found. Next, in Section 4, we consider horizontal and vertical Jacobi fields and we investigate some particular cases.
The Section 2 is devoted to some preliminaries and notations.

We express our hearty thanks to Prof.Radu Miron for his constant encouragements and valuable suggestions for our researches along many years.

\section{Generalized Lagrange spaces}

Let $M$ be a real, smooth i.e. $C^{\infty}$ manifold of finite dimension $n$
and $TM$ its tangent manifold projected to $M$ by the mapping $\tau$.
Set $\buildrel{\circ}\over{TM}=TM\setminus\{0_x\in T_xM,x\in M\}$.
Let $(U,(x^i))$ be a local chart on $M$. The indices $i,j,k,...$will run from 1
to $n$ and the Einstein convention on summation will be implied.
Associate to $v\in \tau^{-1}(U)$ the coordinates $x=(x^i(\tau(u)))$ and
$y=(y^i)$, provided by $V_{\tau(u)}=y^i\partial_i$, $\partial_i=
\dfrac{\partial}{\partial x^i}$, and $TM$ becomes a smooth orientable manifold.
A change of coordinates $(x^i,y^i) \mapsto (x^{i'},y^{i'})$ on $TM$ is as follows:
\begin{equation}
\label{(2.1)}
x^{i'}=x^{i'}(x^1,x^2,...,x^n),\ \ y^{i'}=(\partial_jx^{i'})y^j,\ \
\rm rank (\partial_jx^{i'})=n.
\end{equation}
\begin{definition}
{\rm A set of matrices $(g_{ij}(x,y))$ defined on $\tau^{-1}(U)$ for any open
set $U$ in a smooth atlas on $M$ is said to be a $GL$-metric if\\
1. $g_{ij}(x,y)=g_{ji}(x,y)$, \\
2. $g_{ij}(x,y)=(\partial_i x^{k'})(\partial_j x^{h'})g_{h'k'}(x'(x),y'(y))$
on $U\cap U'$, \\
3. det$(g_{ij}(x,y))\neq 0$, \\
4. The signature of the quadratic form $g_{ij}(x,y)\xi^i\xi^j$, $(\xi^i)\in
R^n$ is constant.}
\end{definition}
The simplest example of a $GL$-metric is a Riemannian one $\gamma_{ij}(x)$. This
is provided according to
\begin{equation}
\label{(2.2)}
\gamma_{ij}(x)=\dfrac{1}{2}\dfrac{\partial^2 L}{\partial y^i\partial y^j},\
\rm with \ L : {\buildrel{\circ}\over{TM}} \rightarrow R \ \rm given \ by
\end{equation}
\begin{equation}
\label{(2.3)}
L(x,y)=\sqrt{\gamma_{ij}(x)y^iy^j}.
\end{equation}
A little more general $GL$-metric is a Finslerian one which is provided by
(2.2) with a function $L=F: \stackrel{\circ}{TM} \rightarrow R_+$
which is positively homogeneous of degree 1 with respect to $y$ i.e.
\begin{equation}
\label{(2.4)}
F(x,\lambda y)=\lambda F(x,y),\ \ \lambda >0.
\end{equation}
A Lagrange metric is a $GL$-metric provided by (2.2) with any smooth
function $L: TM \rightarrow R$.

A large class of $GL$-metrics which are not reducible to the previous ones was
considered in \cite{[Ma2]}:
\begin{equation}
\label{(2.5)}
g_{ij}(x,y)=a(x,y)\gamma_{ij}(x,y)+b(x,y)y_iy_j,
\end{equation}
where $\gamma_{ij}(x,y)$ is a Finsler metric, $a$ and $b$ are smooth functions
on $\stackrel{\circ}{TM}$ such that $a(x,y)>0$, $b(x,y)\geq 0$ and
$y_i=\gamma_{ij}y^j$. Particular forms of these $GL$-metrics were studied in Chapters
11 and 12 of the monograph \cite{[RM]}.

\section{Jacobi fields}
\setcounter{equation}{0}

Let us consider, together with a $GL$-metric $(g_{ij}(x,y))$, a nonlinear
connection $(N^i_j(x,y))$. We have the decomposition
\begin{equation}
\label{(3.1)}
X=hX+vX\ \ \rm for\ \rm every\ \ X \in \chi(TM).
\end{equation}
Denote by $P$ the almost product structure provided by the horizontal and vertical distributions according to
\begin{equation}\label{(3.2)}
P(hX)=hX,\ \ P(vX)=-vX.
\end{equation}
Consider also the almost complex structure $F$ defined as follows:
\begin{equation}\label{(3.3)}
F(hx)=-vX,\ \ F(vX)=hX.
\end{equation}
Next, setting $G=g_{ij}(x,y)dx^i \times dx^j +g_{ij}(x,y)\delta y^i \otimes \delta y^j$,
$\delta y^i=dy^i+ N^i_k(x,y)dx^k$, one gets a metrical structure on $TM$ which is Riemannian if $(g_{ij})$ is positive definite.

\begin{theorem} [\cite{[RM]}]
There exists a unique linear connection $D$ on $TM$ with the properties:
$DP=O, \ \ DF=0,\ \ DG=0$ and $hT(h\cdot,h\cdot)=0,\ \ vT(v\cdot,v\cdot)=0,$
where $T$ denotes its torsion. In the basis $(\delta_i,
\stackrel{.}{\partial_i}),\ \
\delta_i=\partial_i-N^k_i\stackrel{.}{\partial_k},\ \ $ this connection is as
follows:
\begin{equation}\label{(3.3)}
\begin{array}{l}
D_{\delta_k}\delta_j=L^i_{jk}\delta_i, \ \ D_{\stackrel{.}\partial_k}\delta_j=C^i_{jk}\delta_i, \\ \\
D_{\delta_k}\stackrel{.}\partial_j=L^i_{jk}\stackrel{.}\partial_i,\ \ D_{\stackrel{.}\partial_k}\stackrel{.}\partial_j=C^i_{jk}\stackrel{.}\partial_i
\end{array}
\end{equation}
where
\begin{equation}\label{(3.4)}
\begin{array}{l}
L^i_{jk}(x,y)=\dfrac{1}{2}g^{ih}(\delta_jg_{hk}+\delta_kg_{jh}-\delta_hg_{jk}), \\ \\
C^i_{jk}(x,y)=\dfrac{1}{2}g^{ih}(\stackrel{.}\partial_jg_{hk}+\stackrel{.}\partial_kg_{jh}-\stackrel{.}\partial_hg_{jk}).
\end{array}
\end{equation}
\end{theorem}
We note that $D$ has torsion since the other three components of $T$ do not vanish. We have
\begin{equation}\label{(3.5)}
\begin{array}{l}
vT(\delta_k,\delta_j)=R^i_{jk}\delta_i,\ \ R^i_{jk}=\dfrac{\delta N^i_j}{\delta x^k}-\dfrac{\delta N^i_k}{\delta x^j} \\ \\
hT(\partial_k,\delta_j)=C^i_{jk}\delta_i,\\ \\
vT(\partial_k,\delta_j)=P^i_{jk}\partial_i,\ \ P^i_{jk}=\dfrac{\partial N^i_j}{\partial y^k}-L^i_{kj}.
\end{array}
\end{equation}

Thus connection is different from the Levi-Civita connection of $G$ since its torsion does not vanish. We call geodesics the autoparallel curves on $TM$ with respect to $D$. Consider a geodesic $c: [0,1]\rightarrow TM$ such that $c([0,1])\subset
\tau^{-1}(U)$, where $(U,x^i)$ is a local chart on $M$.
Thus the equation of $c$ is
\begin{equation}\label{(3.6)}
\left\{\begin{array}{l}
x^i=x^i(t)\\
y^i=y^i(t),\ \ t \in [0,1]
\end{array}\right.
\end{equation}
The tangent vector field is $\dot c(t)=\dfrac{dx^i}{dt}\dfrac{\partial}
{\partial x^i}+\dfrac{dy^i}{dt}\dfrac{\partial}{\partial y^i}.$
It can be written in the form
\begin{equation}\label{(3.6)}
\begin{array}{l}
\dot c(t)=\dfrac{dx^i}{dt}\delta_i+\(\dfrac{dy^i}{dt}+N^i_k(x(t),y(t))
\dfrac{dx^k}{dt}\)\dot{\partial_i}.
\end{array}
\end{equation}
It results that $\dot c(t)$ is a horizontal vector field if and only if
\begin{equation}\label{(3.7)}
\begin{array}{l}
\dfrac{\delta y^i}{dt}=\dfrac{dy^i}{dt}+N^i_k(x(t),y(t))\dfrac{dx^k}{dt}=0.
\end{array}
\end{equation}
When this condition holds at we say that $c$ is a horizontal geodesic.
By (3.8), $\dot c(t)$ is a vertical vector field if and only if $x^i=x^i_0$ (constant) i.e. the curve $c$ is in the tangent space $T_{p_0}M,\ p_0=(x^i_0).$
In this case we say that $c$ is a vertical geodesic.
The condition $D_{\dot c(t)}\dot c(t)=0$ takes locally the form:
\begin{equation}\label{(3.8)}
\left\{\begin{array}{l}
\dfrac{d^2x^k}{dt^2}+L^k_{ij}\dfrac{dx^i}{dt}\dfrac{dx^j}{dt}+C^k_{ij}
\dfrac{dx^i}{dt}\dfrac{\delta y^j}{dt}=0 \\ \\
\dfrac{\delta ^2x^k}{dt^2}+L^k_{ij}\dfrac{dx^i}{dt}\dfrac{\delta y^j}{dt}+C^k_{ij}
\dfrac{\delta y^i}{dt}\dfrac{\delta y^j}{dt}=0,
\end{array}\right.
\end{equation}
where we have put
\begin{equation}\label{(3.9)}
\begin{array}{l}
\dfrac{\delta^2y^k}{dt^2}=\dfrac{d^2 y^k}{dt^2}+N^k_h\dfrac{d^2x^h}{dt^2}+
\dfrac{dN^k_h}{dt}\dfrac{dx^h}{dt}=\dfrac{d}{dt}\(\dfrac{\delta y^k}{dt}\).
\end{array}
\end{equation}

{\bf Remark.}
{\rm The form of equations (3.10) is preserved by the affine transformation
$t \mapsto c_1t+c_0,\ \ c_0,c_1 \in R$ of parameter, only.}

{\bf Remark.}
{\rm If $c$ is a horizontal geodesic, (3.10) reduces to }
$$
\dfrac{d^2x^k}{dt^2}+L^k_{ij}\frac{dx^i}{dt}\frac{dx^j}{dt}=0\leqno(3.10')
$$
{\rm while if $c$ is a vertical geodesic it becomes}
$$
\dfrac{d^2y^k}{dt^2}+C^k_{ij}\frac{dy^i}{dt}\frac{dy^j}{dt}=0\leqno(3.10'')
$$

\begin{definition}
{\rm Let $c: I \rightarrow TM,\ \ I=[0,1]$ be a geodesic on $TM$. A first order
variation of it is a smooth mapping $\alpha:(-\varepsilon,\varepsilon)\times
I \rightarrow TM$ such that $\alpha(0,t)=c(t)$ and $\alpha_s(t)=\alpha(s,t)$ is
a geodesic for every $s\in(-\varepsilon,\varepsilon),\ \varepsilon \in R$
and $|\varepsilon|$ small.}
\end{definition}

Let $\(\dfrac{\partial}{\partial s}, \dfrac{\partial}{\partial t}\)$ be the natural basis
of the tangent space to $A=(-\varepsilon, \varepsilon)\times I$ in the point
$(s,t)$. We set
$$
\alpha_{\ast,(s,t)}\left.\(\dfrac{\partial}{\partial t}\)\right|_{s=0}=\tau(t),\ \
\alpha_{\ast,(s,t)}\left.\(\dfrac{\partial}{\partial s}\)\right|_{s=0}=V(t).
$$

The vector field $t\mapsto \tau(t)$ is in fact $\dot{c}(t)$, the tangent vector
field to the curve $c$ and the vector field $t \mapsto V(t)$ will be called
the variation vector field induced by $\alpha$.

As $\alpha_{\ast,(s,t)}\left[\dfrac{\partial}{\partial t},\dfrac{\partial}{\partial s}\right]
=[\tau(t),V(t)]$ we infer $[\tau,V]=0$. Thus $T(\tau,V)=D_{\tau}V-D_V\tau$ and
$R(\tau, V)\tau =D_{\tau}D_V\tau-D_VD_{\tau}\tau-D_{[\tau,V]}\tau=
D_{\tau}D_V\tau$ since $D_{\tau}\tau=0$ ($c$ is a geodesic).
Furthermore, $R(\tau, V)\tau =D_{\tau}(D_{\tau}V-T(\tau,V)=D^2_{\tau}V-
D_{\tau}T(\tau,V).$ Thus $V$ satisfies the following equation
$$
D^2_{\tau}V+R(V,\tau)\tau-D_{\tau}T(\tau,V)=0
$$
\begin{definition}
{\rm It is called Jacobi field along of a geodesic $c$ any vector field $X$
which is solution of the following Jacobi equation:}
\begin{equation}\label{(3.11)}
D^2_{\dot{c}(t)}X+R(X,\dot{c}(t))\dot{c}(t)-D_{\dot{c}(t)}T(\dot{c}(t),X)=0
\end{equation}
\end{definition}

As in the Riemannian framework one proves:
\begin{proposition}
$1)$ The solution $X$ of the Jacobi equation is uniquely determined by the
initial conditions $X(t_0)=X_0$ and $(D_{\dot{c}(t)}X)(t_0)=V_0,\ t_0\in I$.\\
$2)$ The set of Jacobi fields is a linear space of dimension $4n$. \\
$3)$ The vector fields $\tau : t \mapsto \tau(t)$ and $\widetilde{\tau} : t
\mapsto t\dot{c}(t)$ are Jacobi vector fields along the geodesic $c$. \\
$4)$ Any Jacobi vector field $X$ along $c$ is of the form $X=a\tau+
b\widetilde{\tau}+Y$, with $a$ and $b$ constants and $Y$ is a Jacobi vector
field which is ortogonal to $\tau$ with respect to $G$.
\end{proposition}

\section{Some particular cases}
\setcounter{equation}{0}

The following particular cases have to be considered: \\
a) $c$ is a horizontal geodesic and $X$ is horizontal. \\
b) $c$ is a vertical geodesic and $X$ is vertical.

In the case a) we have $T(X,\tau)=T(hX,h\tau)=hT(hX,h\tau)+vT(hX,h\tau)=
-v[X,\tau]=0$ since $[X,\tau]=0$.

Thus (3.12) reduces to
\begin{equation}\label{(4.1)}
D^2_{\tau}X+R(X,\tau)\tau=0.
\end{equation}

We notice that for a Finsler metric of a Finsler space $F^n$, $D$ is exactly
the Cartan connection of $F^n$. In (4.1), $R$ is the $(hh)h$- curvature
of $D$ which, coincides to the Chern-Rund connection (see \cite{[Ma2]}).
Thus for a Finsler metric the equation (4.1) is nothing but the equation
(4.13) in \cite{[BaCh]}. Taking $c$ as the lift $\(x^i,\dfrac{dx^i}{dt}\)$ of a curve
$x=x^i(t)$ on $M$ one may try a study similar to that from \cite{[BaCh]}
with some cautions regarding the parameter of the curve $x^i=x^i(t)$.

In case b) we have again $T(X,\tau)=0$ and (3.12) reduces to
\begin{equation}\label{(4.2)}
D^2_{\tau}X+R(X,\tau)\tau=0.
\end{equation}
In this case $R$ is the $(vv)v$-curvature of $D$, usually denoted by $S$.
Now the curve $c$ is entirely in $T_{p_0}M$, $p_0=(x^i_0)$. The space $T_{p_0}M$
has a pseudo-Riemannian structure given by $g_{ij}(x_0,y)$ whose curvature is $S$.
This pseudo-Riemannian structure is not flat except if the $GL$-metric $g_{ij}(x,y)$
is a Riemannian one. The equation (4.2) is exactly the Jacoby equation for
$(T_{p_0}M, g_{ij}(x_0,y))$ and when $g_{ij}(x_0,y)$ is positive defined one may
apply the theory from the Riemannian case. The geodesics in $(T_{p_0}M, g_{ij}(x_0,y))$ are sometimes called $v-$paths.

Let us consider the $GL$-metric
\begin{equation}
g_{ij}(x,y)=\gamma_{ij}(x)+b(x,y)y_iy_j,
\end{equation}
where $\gamma_{ij}(x)$is a Riemannian metric and $b : \stackrel{\circ}{TM}
\rightarrow R$ is a smooth function such that $b(x,y)>0$. Together with this $GL$-metric
we may consider the nonlinear connection $N^i_j(x,y)=\gamma^i_{kj}(x)y^k$,
where $\gamma^i_{kj}(x)$ are the Christoffel symbols derived from $(\gamma_{ij}(x))$.
It is not difficult to see that, with this choice the projection
$\tau (TM, G) \rightarrow (M, \gamma)$ is a Riemannian submersion. Thus the
general theory of submersion may be used in order to investigate (4.1).
It follows that if a curve is horizontal in a point it is horizontal at any points
and any horizontal curves is projected by $\tau$ on a geodesic of $(M,\gamma)$.
Furthermore, the Jacobi fields on $TM$ which are solutions of (4.1) are projected  by
$\tau_{\ast}$ on Jacobi fields on $(M, \gamma)$.

\medskip

\noindent{\footnotesize{\it Received April 10, 1997}}

\newpage

\def\ijk{^i_{jk}}
\def\bld#1#2{\buildrel{#1}\over{#2}}
\def\derp#1#2{\displaystyle\frac{\partial#1}{\partial#2}}
\def\derpp#1{\displaystyle\frac{\partial}{\partial#1}}

\runningauthor={M. ANASTASIEI}
\runningtitle={THE BEIL METRICS ASSOCIATED TO A FINSLER SPACE}
\noindent
\baselineskip 8pt
\noindent{\footnotesize{Balkan J. Geom. Appl.}}
\hfill\break
{\footnotesize{3, no. 2, 1--16, 1998}}
\vskip 2cm
\baselineskip 11.5pt plus .15pt
\centerline{\bf\Large THE BEIL METRICS}
\vskip .2cm
\centerline{\bf\Large ASSOCIATED TO A FINSLER SPACE}
\vskip .5cm
\centerline{\bf {\footnotesize{BY}}}
\vskip .5cm
\centerline{\bf {\footnotesize{M. ANASTASIEI and H. SHIMADA}}}
\vskip 1cm

\setcounter{section}{0}
\section{Introduction}
\setcounter{equation}{0}
Let $F^n=(M,F)$ be a Finsler space with $M$ a
smooth i.e. $C^\infty$ manifold and $F:TM\to R$, $(x,y)\to F(x,y)$.
Assume that $F^n$ is endowed with a Finsler $1$--form
$\beta_i(x,y)$ and set $\beta=\beta_i(x,y) y^i.$ Here $i,j,k,...$
will run from 1 to $n={\rm dim}\:M$ and the Einstein convention on
summation is implied. Then $^*F=L(F,\beta)$ in some conditions on
$L$ is so that $^*F^n=(M,^*F)$ is a new Finsler space. It is said
that $^*F^n$ is obtained from $F^n$ by a $\beta$--change
\cite{7},\cite{10}.

Typical for $^*F^n$ are the Randers and Kropina spaces which are
obtained from a Riemannian space by particular $\beta$--changes.

Let $g_{ij}(x,y)$ be the Finsler metric tensor of $F^n$. If one wishes
the construction of a new Finsler metric $^*g_{ij}$ which depends on
$g_{ij}(x,y)$, then because of the linear structure of the set of Finsler
tensor fields of a given type, the most general choice is
$$^*g_{ij}(x,y)=\rho(x,y)g_{ij}(x,y)+\sigma(x,y) B_{ij}(x,y),\leqno(1.1)$$
for $\rho$ and $\sigma$ two Finsler scalars and $B_{ij}(x,y)$ a
symmetric Finsler tensor field of type $(0,2)$. We may say that
$^*g_{ij}$ is obtained from $g_{ij}$ by a $B$--change.

It is clear that $^*g_{ij}$ is no longer a Finsler metric except if
some strong conditions on $\rho,\sigma$ and $B_{ij}$ are imposed. Metrics
similar to (1.1) appear in \cite{2} and \cite{5} from physical
considerations. See also \cite{11}.

In order to relax such conditions we do not ask $^*g_{ij}$ be a Finsler
metric but a generalized Lagrange metric in Miron' sense, shortly a
$GL$--metric. For the theory of the $GL$--metrics we refer to
\cite{9}, ch.X.

As such $(^*g_{ij})$ has to satisfy
\begin{itemize}
\item[a)] $\det(^*g_{ij})\ne0$ and
\item[b)] The quadratic from $^*g_{ij}(x,y)\xi^i\xi^j,\ (\xi^i)\in
\R^n,$ to be of constant signature.\end{itemize}

Even this minimal requirements are not easy to be fullfiled except
for some particular $\sigma,\rho$ and $B_{ij}.$

By our best knowledge the following two particular forms of the
$GL$--metric (1.1) were studied
$$^*{g}_{ij}(x,y)=e^{2\alpha(x,y)}g_{ij}(x,y).\leqno(1.2)$$
This class of $GL$--metrics contains the Miron--Tavakol metrics used by
them in General Relativity and the Antonelli metrics which were
introduced by P.L. Antonelli for some studies in Biology and
Ecology. For details see \cite{9}, ch.XI, and reference therein.
$$^*g_{ij}(x,y)=g_{ij}(x,y)+\sigma(x,y) y_iy_j,\ y_i=g_{ij}(x,y) y^j.\leqno(1.3)$$
Particular forms of the $GL$--metric (1.3) were used by R. Miron in
Relativistic Geometrical Optics. See also \cite{9}, ch.XII.

Some particular forms of the $GL$--metric
$$^*g_{ij}(x,y)=g_{ij}(x,y)+\sigma(x,y) B_i(x,y) B_j(x,y),\leqno(1.4)$$
with $B_i(x,y)=g_{ij}(x,y) B^j(x,y)$ for $B^j(x,y)$ a given Finsler vector
field were introduced by R.G. Beil in order to develop his interesting unified field theory (\cite{4}). These were
called Beil metrics. As such we refer to $^*g_{ij}$ in (1.4) as to the Beil metric, too. The following comment of R.G. Beil is illuminating on (1.4). ``Since in my unified theory the quantity $k$ which correspond to your $\sigma$
is related to the gravitational constant, this means that a possible physical
interpretation of your theory with a $y-$dependent $\sigma$ is that gravitation
itself is velocity dependent. This possibility is mentioned, for example, in
Section 40.8 of the famous book {\it Gravitation} by Misner, Thorne and Wheeler''.
See \cite{13}.

The particular form of (1.4) obtained for $\sigma=1$ and $B_{i}=\dfrac{\partial f}{\partial x^i}$,
$f:M\rightarrow\R$ was considered by C. Udri\c ste in \cite{14}. He proved that if $f$ is proper i.e.
$f^{-1} (K)$ is a compact set whenever K is compact, then the Finsler manifold
$(M,^*g_{ij}(x,y)) $ is complete. A Riemannian version of (1.1), that is,
was used by T. Aubin in order to prove that any compact Riemannian manifold
of dimension greater then 2 admits a metric whose scalar curvature is a negative constant.
See \cite{3} and for other connected results.

The geometry of the $GL$--metrics (1.4) was not investigated in a systematic
way. It is our purpose to fill this gap. After some preliminaries
in Section 2, we show in Section 3 that $(^*g_{ij})$ from (1.4) is a
$GL$--metric and we point out cases when it reduces to a Lagrange or to a
Finsler metric. In Section 4 we discuss possibilities for
introducing metrical connections for the $GL$--space $(M,^*g_{ij})$.
In Section 5 we digress on parallel and resp. concurrent Finsler vector
fields showing that the usual definitions for these notions are
also justified from the viewpoint of the almost Hermitian model of
a $GL$--space. For such a model see \cite{9}, ch. X. Section 6 is
devoted to the analysis of the $GL$--metric (1.4) for $B^i$ a
concurent Finsler vector field. For $\sigma$ a constant we
rediscover a modification of a Finsler function studied by M. Matsumoto
and K. Eguchi in \cite{8}. The case when $\sigma$ is a solution of
the so--called Tavakol--Van der Berg equation is investigated, too. In
Section 7 we treat a Beil metric associated to a Finsler space with
$(\alpha,\beta)$--metric. It is a future task to find properties of
the $GL$--metric (1.4) when $F^n$ is a particular Finsler space or its
dimension is low (2 or 3).

\section{Preliminaries}
\setcounter{equation}{0}

Let $M$ be a smooth i.e. $C^\infty$ manifold, paracompact and of dimension $n$, $TM$ its tangent manifold and
$\tau:TM\to M$ its tangent bundle. If $x=(x^i)$, $i,j,k,...=1,...,n$ are local coordinates on $M$, then the induced
coordinates on $TM$ will be $(x,y)=(x^i:x^i\circ\tau,y^i)$ with $(y^i)$ provided by $u_x=y^i\left.\derpp{x^i}\right|_x,$ $u\in T_xM$, $x\in M$. The change of coordinates $(x,y)\to(\tilde x,\tilde y)$ on $TM$
are as follows.
$$\begin{array}{l}
\tilde x^i=\tilde x^i(x^1,...,x^n),\ {\rm rank}\:\left(\derp{\tilde
x^i}{x^k}\right)=n\vspace*{1,4mm}\\ \tilde y^i=\derp{\tilde x^i}{x^k}(x)
y^k.\end{array}\leqno(2.1)$$

The geometrical objects on $TM$ whose local components change by (2.1)
as on $M$ i.e. ignoring their dependence on $y$, will be called
Finsler objects as in \cite{7} or $d$--objects as in \cite{9}.

We set $\partial_i:=\derpp{x^i},$ $\dot\partial_i:=\derpp{y^i}$ and notice
that the vertical subspace of $T_uTM$ i.e. $V_uTM={\rm Ker}\:(D\tau)_u,$
$u\in TM$, where $D\tau$ means the differential of $\tau$, is
spanned by $(\dot\partial_i).$ The $d$--objects can be expressed using
$(\dot\partial_i).$

A function $F:TM\to\R$ which is positive, smooth on $TM\setminus0$ and
only continuous in the rest, positively homogeneous of degree 1 with
respect to $y$ i.e. $F(x,\lambda y)=\lambda F(x,y),$ $\lambda>0$ and with the
quadratic form $g_{ij}(x,y)\xi^i\xi^j,$ $(\xi^i)\in\R^n$
nondegenerate and of constant signature, where
$$g_{ij}(x,y)=\displaystyle\frac12\dot\partial_i\dot\partial_jF^2,\leqno(2.2)$$
is called a fundamental Finsler function. The pair $F^n=(M,F)$ is called a Finsler
space.

The function $g_{ij}(x,y)$ are the components of a Finsler tensor field
called the Finsler metric of $F^n$.

A supplement $H_uTM$ of $V_uTM$ i.e. the decomposition in a direct
sum $T_uTM=H_uTM\oplus V_uTM$ holds, will be called the horizontal space
and the distribution $u\to H_uTM$ will be called a horizontal
distribution. A basis of it of the form
$\delta_i=\partial_i-N_i^k(x,y)\dot\partial_k,$ provides the functions $(N^k_i(x,y))$
called the local coefficients. These functions have a special rule of change
by (2.1) and in turn they completely determine the horizontal
distribution called also a nonlinear connection. Then
$(\delta_i,\dot\partial_i)$ is a basis adapted to the previous
decomposition of $T_uTM$. The Finsler objects may be also expressed
by using $(\delta_i).$ We notice that $(\delta_i)$ are Finsler vector
fields. For more details we refer to
\cite{7},\cite{9}.

\section{The Beil metric}
\setcounter{equation}{0}

Let $F^n=(M,F)$ be a Finsler space and
$g_{ij}(x,y)$ its Finsler metric. Assume that $F^n$ is endowed with a
Finsler vector field $B=B^i(x,y)\dot\partial_i$ and let $B_i(x,y) dx^i$
the Finsler 1-form with $B_i=g_{ik}B^k.$ The lowering and rising of
indices will be done with $(g_{ij})$ and $(g^{jk})$, where
$g^{jk}g_{ki}=\delta^j_i,$ respectively. Let $\sigma:TM\to\R,$
$(x,y)\to\sigma(x,y)$ a Finsler scalar. We set
$$^*g_{ij}(x,y)=g_{ij}(x,y)+\sigma(x,y) B_i(x,y) B_j(x,y).\leqno(3.1)$$

The functions $(^*g_{ij})$ from (3.1) define for $\sigma>0$ a positive
definite $GL$--metric called {\it the Beil metric.}

It is clear that  $(^*g_{ij})$ are the components of a symmetric $d$--tensor
field. We look for the inverse of the matrix
$(^*g_{ij})$ in the form $^*g^{jk}={}^*g^{jk}-{}^*\sigma B^jB^k$ with
$^*\sigma$ to be determined. From $^*g_{ij}{}^*g^{jk}=\delta^k_i$ it
follows that $^*\sigma=\displaystyle\frac\sigma{1+\sigma B^2},$ with
$B^2=B_iB^i=g_{ij} B^iB^j$ (the length of $B$ with respect to
$g_{ij}$). Thus we have
$$^*g^{jk}=g^{jk}-\displaystyle\frac\sigma{1+\sigma B^2}B^jB^k.\leqno(3.2)$$

Consequently, we have $\det(g_{ij})\ne0.$

The quadratic from
$\Phi(\xi)={}^*g_{ij} \xi^i\xi^j=g_{ij}\xi^i\xi^j+\sigma(B_k\xi^k)^2$ is
clear positive definite in our hypothesis.\hfill{\bf q.e.d.}\medskip

We notice that (3.2) holds in the weaker condition
$\sigma\ne-\displaystyle\frac1{B^2}$ and if $g_{ij}\xi^i\xi^j$ is only of
constant signature, the signature of $\Phi(\xi)$ will be constant
for some $\sigma$ and $(B^k)$ at least locally. \bigskip

{\it Remark 3.1.} The $GL$--metric (3.1) appears in papers by R.G. Beil
(\cite{4}) for $F^n$ a pseudo--Riemannian space or a
Minkowski space. It was called Beil's metric.\bigskip

We notice that for $B^i=y^i$ in (3.1) one obtains a general version
of the Synge metric which was used by R. Miron for a geometrical theory
of Relativistic Optics (cf. \cite{9}, ch.XI).

In the following we shall assume $B^i\ne y^i$ and use the ideas and
techniques from \cite{9}, ch.XI.

One says that $^*g_{ij}$ is reducible to a Lagrange metric, shortly
an $L$--metric if there exists a Lagrangian $L:TM\to\R$ such that
$^*g_{ij}=\displaystyle\frac12\dot\partial_i\dot\partial_jL.$ A necessary and sufficient condition
for $^*g_{ij}$ be reducible to an $L$--metric is the symmetry in all
indices of the Cartan tensor field
$^*C_{ijk}=\displaystyle\frac12\dot\partial_k{}^*g_{ij}$ i.e.
$$\dot\partial_k{}^*g_{ij}=\dot\partial_i{}^*g_{kj}.\leqno(3.3)$$
Using (3.1) this condition becomes
$$\begin{array}{c}
\dot\sigma_kB_iB_j-\dot\sigma_iB_kB_j+\sigma
(\dot\partial_kB_i\cdot B_j-\dot\partial_iB_k\cdot B_j)+\vspace*{1,4mm}\\
+\sigma(B_i\cdot \dot\partial_kB_j-B_k\cdot\dot\partial_iB_j)=0,\
\dot\sigma_k:=\dot\partial_k\sigma.\end{array}\leqno(3.4)$$
Multiplying it by $B^j$ one gets
$$\begin{array}{c}
B^2(\dot\sigma_kB_i-\dot\sigma_iB_k)+\sigma
B^2(\dot\partial_kB_i-\dot\partial_iB_k)+\vspace*{1,4mm}\\
+\sigma(B_i\cdot\dot\partial_kB_j\cdot B^j-
B_k\dot\partial_iB_j\cdot B^j)=0.\end{array}\leqno(3.5)$$
If (3.4) is an identity, then (3.5) should be an identity for any
$\sigma$ and $B_i.$ But for $B_i=B_i(x)$ and $\sigma=F^2$, (3.5)
reduces to $y_kB_i-y_iB_k=0$ which is not an identity for any $B_i.$
Thus in general $^*g_{ij}(x,y)$ is not reducible to an $L$--metric.

We have a case when $^*g_{ij}(x,y)$ is an $L$--metric as follows.\bigskip

{\bf Proposition 3.1.} {\it Assume $B_i=B_i(x)$. If $\sigma(x,y)=f(B_i(x)
y^i)$ for a smooth function $f:\R\to\R,$ then $^*g_{ij}$ is an
$L$--metric.}\bigskip

Indeed, it is easy to check that in these hypothesis (3.4)
identically holds. Notice that we do not know which is $L$ such
that $^*g_{ij}=\displaystyle\frac12\dot\partial_i\dot\partial_jL.$

It is said that $^*g_{ij}(x,y)$ is weakly regular if its absolute
energy $${\cal E}(x,y):={}^*g_{ij}(x,y)
y^iy^j=F^2(x,y)+\sigma(x,y)(B_iy^i)^2\leqno(3.6)$$ is a regular Lagrangian
i.e. the matrix with the entries
$$a_{kh}(x,y)=\displaystyle\frac12\dot\partial_k\dot\partial_h{\cal E},\leqno(3.7)$$
is of ${\rm rank}\: n.$

A direct calculation yields
$$a_{kh}=g_{kh}+\displaystyle\frac12\dot\sigma_{kh}\beta^2+
\beta(\dot\sigma_k\dot\beta_h+\dot\sigma_h\dot\beta_k)+
\sigma\dot\beta_k\dot\beta_h+\sigma\beta\dot\beta_{kh},\leqno(3.8)$$
$$\beta:=B_i(x,y) y^i,\dot\beta_k:=\dot\partial_k\beta,
\dot\beta_{kh}:=\dot\partial_k\dot\partial_h\beta,
\dot\sigma_{kh}:=\dot\partial_k\dot\partial_h\sigma,\
\dot\sigma_k:=\dot\partial_k\sigma\leqno(3.8)'$$
It is hopeless to decide if $a_{kh}$ is invertible or not. However
we have some interesting particular cases.\bigskip

{\bf Proposition 3.2 }\ \vspace*{-1,5mm}
\begin{itemize}
\item[a)] {\it If $B$ is orthogonal to the Liouville vector field
$\C=y^i\dot\partial_i,$ then $^*g_{ij}$ is weakly regular and
$a_{kh}(x,y)=g_{kh}(x,y).$}
\item[b)] {\it If $B_i=B_i(x)$ and $\sigma(x,y)=f(\beta)$ for some
smooth function $f:R\to R$, then $^*g_{ij}$ is weakly regular if and only if
$1+\varphi(\beta)B^2\ne0,$ where $2\varphi(\beta)=\beta^2f''+4\beta
f'+2f,$ $f'=\displaystyle\frac{df}{d\beta},$ $f''=\displaystyle\frac{d^2f}{d\beta^2}$ and we
have}\end{itemize}
$$a_{kh}(x,y)=g_{kh}(x,y)+\varphi(x,y) B_k(x) B_h(x).\leqno(3.9)$$

\medskip{\bf Proof.} a) The condition $B$ orthogonal to  $\C$ is equivalent to $\beta=0.$
Thus ${\cal E}(x,y)=F^2(x,y)$ and so $a_{kh}=g_{kh}.$

b) By a direct calculation one finds (3.9). Hence $(a_{kh})$ has
the same form as $^*g_{kh}$ with $\sigma$ replaced by $\varphi.$ The
conclusion follows.

We keep the hypothesis $B_i=B_i(x)$ and $\sigma=f(\beta)$,
$\beta\ne0.$ From (3.9) we see that we have again $a_{kh}=g_{kh}$
when $\varphi=0.$ The differential equation $\beta^2f''+4\beta f'+2f=0$
takes the form $(\beta^2f'+2\beta f)'=0$ and so its general solution
is $f(\beta)=\displaystyle\frac a\beta+\displaystyle\frac b{\beta^2},$ $a,b\in\R$. The
metric $^*g_{ij}$ becomes
$$^*g_{ij}=g_{ij}+\left(\displaystyle\frac a{B_i(x) y^i}+
\displaystyle\frac b{(B_s(x)y^s)^2}\right)B_i(x) B_j(x).\leqno(3.10)$$
Notice that although $^*g_{ij}$ is an $L$--metric, we do not yet
know the Lagrangian $L$.

The absolute energy of $^*g_{ij}$ is now ${\cal E}=F^2+a(F_i(x) y^i)+b$
and the Lagrange space $L^n=(M,{\cal E})$ is called an almost
Finslerian--Lagrange space (see Section 6, ch.IX of \cite{9}).

We may put (3.9) into the form
$$a_{kh}(x,y)={}^*g_{kh}+
\left(\displaystyle\frac12\beta^2f''+2\beta f'\right)B_kB_h.\leqno(3.9)'$$
Thus we see that $a_{kh}={}^*g_{kh}$ if and only if $f$ is a solution of the
differential equation
$$\displaystyle\frac12 f''\beta^2+2\beta f'=0 \mbox{ i.e. }f(\beta)=c-\displaystyle\frac
d{\beta^3},\ c,d\in\R.$$
We know that $^*g_{kh}$ is an $L$--metric (in previous
hypothesis). The condition $a_{kh}={}^*g_{kh}$ gives $L$ in the form
$L(x,y)={\cal E}(x,y)+A_i(x) y^i+\psi(x)$, where $A_i$ is a covector and
$\psi$ a scalar. Inserting here ${\cal E}$ we get
$$L(x,y)=F^2(x,y)+c(B_i(x,y)y^i)^2-\displaystyle\frac d{B_i(x) y^i}+A_i(x) y^i+\psi(x),\
c,d\in\R. \leqno(3.10)'$$
Therefore we found a case when $^*g_{ij}$ is an $L$--metric with $L$
of explicit form $(3.10)'$.\bigskip

{\bf Remark 3.2} In the hypothesis of a) in Proposition 3.2, $^*g_{ij}$ is
not necessarily an $L$--metric. If $\sigma(x,y)$ and $B_i(x,y)$ are
positively homogeneous of degree $0$, then $^*g_{ij}(x,y)$ is so and
$(M,^*g_{ij})$ is a generalized Finsler space in Izumi' sense (see
\cite{6}).\bigskip

{\bf Remark 3.3.} The condition $B$ orthogonal to $\C$ is equivalent
with the condition $B$ is tangent to the indicatrix bundle $I(M)\subset
TM.$\bigskip

{\bf Caution.} The conditions $\beta=0$ and $B_i=B_i(x)$ are
incompatible since they lead to  $B=0.$\bigskip

{\bf Remark 3.4.} If in (3.10) we take $d=0,$ $A_i=0,$ $\psi=0,$
$c>0$, then $^*F^2:=L(x,y)$ is positively homogeneous of degree $2$
and so $^*F^n=(M,^*F)$ becomes a Finsler space. Notice that $^*F$ is
getting from $F$ by a $\beta$--change and in this case $^*g_{ij}$
reduces to a Finsler metric.\bigskip

{\bf Remark 3.5.} An interesting Beil metric can be associated to
a Finsler space $F^n$ with an $(\alpha,\beta)$--metric. Here
$\alpha^2=a_{ij}(x) y^iy^j$ and $\beta=b_i(x) y^i,$ where $a_{ij}$ is a
Riemannian metric an $b_i$ a covector field on $M$. One may consider
$$^*g_{ij}(x,y)=a_{ij}(x)+\sigma(x,y) b_i(x) b_j(x),\leqno(3.12)$$
where $\sigma$ is a Finsler scalar such that $1+\sigma b^2\ne0$ for
$b^2=a^{ij}b_ib_j.$ This $GL$--metric is not reducible to an $L$--metric
or a Finsler metric. The previous discussion applies, too.

\section{Metrical connections for $GL\!=\!(M,^*g_{ij}(x,y))$}
\setcounter{equation}{0}

In Finsler geometry as well as in their generalizations, the nonlinear connections
play an important role. For instance these connections allow us to
work with $d$-- or Finsler objects and so to keep and check easily
the geometrical meaning of calculation in local coordinates.

A nonlinear connection always exists if $M$ is paracompact. But the
nonlinear connections derived from or associated in a way to a $GL$--metric are much
more useful. There are no possibilities to find nonlinear connections for any
$GL$--metric. But there are some classes of $GL$--metrics for which such
possibilities exist. One is that of weakly regular $GL$--metrics and as it
is well known there exist nonlinear connections canonically derived from a Lagrangian, a
Finslerian or a Riemannian metric. See \cite{9} for details.

We recall here the Cartan nonlinear connection for $F^n$. Set
$$\gamma^i_{jk}(x,y)=\displaystyle\frac12 g^{ih}(\partial_jg_{hk}+
\partial_kg_{hj}-\partial_hg_{jk}),\
\gamma^i_{00}:=\gamma^i_{jk} y^jy^k.\leqno(4.1)$$
Then $\bld{\circ}N{}\!^i_j=\displaystyle\frac12\dot\partial_j\gamma^i_{00}$ are the local
coefficients of the Cartan nonlinear connection.

For any Finsler connection $F\Gamma(N)$ we denote by ${\scriptstyle|} k$ and
${\big|}_k$ its $h$-- and $v$--covariant derivatives. Then $F\Gamma(N)$
is called $h$--metrical if $g_{ij{\scriptstyle|} k}=0$ and $v$--metrical if
$g_{ij}{\big|}{}_k=0.$

We consider
$$\begin{array}{l}
F^i_{jk}=\displaystyle\frac12 g^{ih}(\delta_j g_{hk}+\delta_kg_{jh}-\delta_hg_{jk}),\vspace*{1,4mm}\\
C^i_{jk}=\displaystyle\frac12 g^{ih}(\dot\partial_j
g_{hk}+\dot\partial_kg_{jh}-\dot\partial_hg_{jk}),\end{array}\leqno(4.2)$$
where $\delta_j=\partial_j-\bld{\circ}N{}\!^k_j\dot\partial_k.$
For $F^n$ we have four remarkable Finsler connections based on
$(\bld\circ N{}\!^i_j)$.

We mention here only the Cartan connection
$C\Gamma(\bld\circ N)=(\bld\circ N{}\!^i_j,F^i_{jk},C^i_{jk})$. This is
$v$-- and $h$--metrical and two torsions of it vanishes.

Let us come back to the $GL$--metric (3.1). We cannot derive a nonlinear connection
from it. But since it is constructed with $g_{ij}(x,y)$, we may take
into consideration the Cartan nonlinear connection $(\bld\circ N{}\!^i_j)$ and
then all possible nonlinear connections have the form $N^i_j=\bld\circ
N{}\!^i_j-A^i_j$ with $A^i_j(x,y)$ an arbitrary Finsler tensor field
of type $(1,1).$

Now we replace in the right side of (4.2) the metric $g_{ij}$ by
$^*g_{ij}$ and the operator $\delta_j$ by $^s\delta_j=\partial_j-\bld\circ
N{}\!^k_j\dot\partial_k+A^k_j\dot\partial_k$ and denote the results in the
left side by $^sF^i_{jk}$ and $^sC^i_{jk}$, respectively. Thus we get a
Finsler connection $^sC\Gamma(N)=(N^i_j,^sF^i_{jk},^sC^i_{jk})$ which we call
standard metrical connection of $GL.$

This connection is metrical i.e. $^*g_{ij\bld s{\scriptstyle|} k}=0,$
$^*g_{ij}\bld s{\big|}{}\!_k=0$ and its $h(hh)$-- torsion and
$v(vv)$--torsion vanish. It is clear that it depends on $A^i_j$ but
if $A^i_j$ is given apriori it is the unique Finsler connection with the
above properties. For $A^i_j=0$ we set $^*F:={}^sF$ and $^*C:={}^sC$.
Thus we have
$$\begin{array}{l}
^sF^i_{jk}=^*F^i_{jk}+\displaystyle\frac12{}^*g^{ih}(A^s_j\dot\partial_s{}^*g_{hk}+
A^s_k\dot\partial_s{}^*g_{hj}-A^s_h\dot\partial_s{}^*g_{jk})\vspace*{1,4mm}\\
^sC^i_{jk}={}^*C^i_{jk}.\end{array}\leqno(4.3)$$
The first equation in (4.3) takes also the form
$$\begin{array}{l}
^sF_{jik} = ^*F_{jik}+ ^*C_{kis}A^s_j+ ^*C_{jis}A^s_k- ^*g^{ih}A^l_h {}^*C_{jkl}.
\end{array}$$
{\bf Remark 4.1.} If $(^*g_{ij})$ reduces to an $L$--metric or to a
Finsler metric, (4.3) becomes
$$\begin{array}{l}
^sF{}\!^i_{jk}={}^*F^i_{jk}+C^i_{ks}A^s_j\vspace*{1,4mm}\\
^sC{}\!^i_{jk}={}^*C^i_{jk}.\end{array}\leqno(4.3)'$$
We notice the following possible choices of $A^i_j:$
$\lambda(x,y)\delta^i_j,$ $y^iy_j,$ $B^iy_j$, $y^iB_j,$ $B^iB_j$.

By (3.1) we find
$$\begin{array}{ll}
^*F^i_{jk}
&=B^i_sF^s_{jk}+\displaystyle\frac\sigma2\bld*g{}\!^{ih}
[\delta_j(B_hB_k)+\delta_k(B_hB_j)-\delta_h(B_jB_k)]+\vspace*{1,4mm}\\
&+\displaystyle\frac12{}^*g^{ih}(\sigma_jB_hB_k+
\sigma_kB_hB_j-\sigma_hB_jB_k),\vspace*{1,4mm}\\
^*C^i_{jk}
&=B^i_sC^s_{jk}+\displaystyle\frac\sigma2\bld*g{}\!^{ih}
[\dot\partial_j(B_hB_k)+\dot\partial_j(B_hB_j)-
\dot\partial_j(B_jB_k)]+\vspace*{1,4mm}\\
&+\displaystyle\frac12{}^*g^{ih}(\dot\sigma_jB_hB_k+
\dot\sigma_kB_hB_j-\dot\sigma_hB_jB_k),\
\ \mbox{with}\end{array}\leqno(4.4)$$
$$B^i_s=\dot\partial^i_s-{}^*\sigma B^iB_s,\ \sigma_k:=\delta_k\sigma,\
\dot\sigma_k:=\dot\partial_k\sigma,\ ^*\sigma=\sigma/(1+\sigma
B^2).\leqno(4.4)'$$ Now, $^sF{}\!^i_{jk}$ and $^sC{}\!^i_{jk}$
are easily deduced from (4.3).\bigskip

{\bf Remark 4.2.} The matrix $B^i_s$ is invertible. Its inverse is
$(B^{-1})^s_k=\delta^s_k+\sigma B^sB_k.$ As such from (4.4) we can find
$F$ and $C$ as depending on $^*F$ and $^*C$.\bigskip

In order to evaluate the torsions and curvatures of $^*C\Gamma(^cN)$ it
is more convenient to put (4.4) into the form
$$\begin{array}{l}
^*F^i_{jk}=F^i_{jk}+\Lambda^i_{jk},\vspace*{1,4mm}\\
^*C^i_{jk}=C^i_{jk}+\bld\circ\Lambda{}\!^i_{jk},\ \mbox{for}\end{array}\leqno(4.5)$$
$$\begin{array}{ll}
\Lambda^i_{jk}
&=\displaystyle\frac12{}^*g^{ih}[\delta_k(\sigma B_jB_h)+\delta_j(\sigma B_hB_k)-\delta_h(\sigma
B_jB_k)]+\vspace*{1,4mm}\\
&-{}^*\sigma B^iB^hF_{jhk}\vspace*{1,4mm}\\
\bld\circ\Lambda{}\!^i_{jk}
&=\displaystyle\frac12{}^*g^{ih}[\dot\partial_k(\sigma B_jB_h)+\dot\partial_j(\sigma B_hB_k)-
\dot \partial_h(\sigma B_jB_k)]+\vspace*{1,4mm}\\
&-{}^*\sigma B^iB^hC_{ijk}.\end{array}\leqno(4.5)'$$
The torsions of $^*C\Gamma(^cN)$ are as follows.
$$\begin{array}{l}
^*T^i_{jk}=0,\ ^*R^i_{jk}=R^i_{jk},\ ^*S^i_{jk}=0\vspace*{1,4mm}\\
^*P^i_{jk}=P^i_{jk}-\Lambda^i_{kj}\mbox{ and }^*C^i_{jk}\mbox{ from
(4.5)}.\end{array}\leqno(4.6)$$ As for the curvatures we have\medskip

\noindent$\!\!\!\begin{array}{lcl}
(4.7)&\ \
&^*S_j{}^i{}_{kh}=S_j{}^i{}_{kh}+
\bld\circ\Lambda{}\!_j{}^i{}_{kh}+(C^s_{jk}
\bld\circ\Lambda{}\!^i_{sh}+\bld
o\Lambda{}\!^s_{jk}C^i_{sh}-(k/h))\vspace*{1,4mm}\\
(4.7)'&&\bld\circ\Lambda{}\!_j{}^i{}_{kh}=
\dot\partial_h\bld\circ\Lambda{}\!^i_{jk}+\bld\circ\Lambda{}\!^s_{jk}\bld
o\Lambda{}\!^i_{sh}-(k/h),\end{array}$\medskip

\noindent where $-(k/h)$ means the substraction of the preceeding terms
with $k$ replaced by $h$.
$$\begin{array}{cll}
^*P{}\!_j{}^i{}_{kh}=P_j{}^i{}_{kh}+\dot\partial_h\Lambda^i_{jk}-
\bld\circ\Lambda{}\!^i_{jh{\scriptstyle|} k}-C^i_{jh{\scriptstyle\|}k}-
\bld\circ\Lambda{}\!^i_{jh{\scriptstyle\|}k}+\vspace*{1,4mm}\\
+\dot\partial_kC^i_{jh}+\dot\partial_k\bld
o\Lambda^i_{jh}-C^i_{js}\Lambda^s_{hk}+\bld
o\Lambda{}\!^i_{js}P^s_{hk}- \bld
o\Lambda{}\!^i_{js}\Lambda^s_{kh},\end{array}\leqno(4.8)$$
where ${ {\scriptstyle\|}k}$ denotes a covariant derivative constructed with
$\Lambda^i_{jk}$.
$$^*R_j{}^i{}_{kh}=R_j{}^i{}_{kh}+\Lambda_j{}^i{}_{kh}+
(F^s_{jk}\Lambda^i_{sh}+\Lambda^s_{jk}F^i_{sh}-(k/h))+
\bld\circ\Lambda{}\!^s_{js}R^s_{kh},\leqno(4.9)$$
where
$$\Lambda_j{}^i{}_{kh}=\delta_h\Lambda^i_{jk}+
\Lambda^s_{jk}\Lambda^i_{sh}-(k/h).\leqno(4.9)'$$

\section{Parallel and concurrent Finsler vector fields}
\setcounter{equation}{0}

Let $B^i(x,y)$ be a Finsler vector field and $F\Gamma(N)$ be a Finsler connection. Then it
is said that $(B^i)$ is parallel if
$$B^i_{|k}=0,\ B^i{\big|}_k=0\leqno(5.1)$$
and $(B^i)$ is concurrent if
$$B^i_{|k}=-\delta^i_k,\ B^i{\big|}_k=0.\leqno(5.2)$$
It is our purpose to confirm the correctness of these definitions
from the viewpoint of the almost K\"ahlerian model of a Finsler space (see
\cite{9}, ch.VII for details on this model). A different
confirmation of these definitions is given in \cite{8} using the
principal Finsler bundle model due to M. Matsumoto. The giving of
$N$ is equivalent to the decomposition
$$T_uTM=H_uTM\oplus V_uTM,\ u\in TM\mbox{ (Whitney'
sum)}.\leqno(5.3)$$ Accordingly we have two projectors $h$ and $v$ and
an almost product structure $P$ such that if we put $X=hX+vX$ for a
vector field $X$ on $TM$, then
$$P(hX)=hX,\ P(vX)=-vX.\leqno(5.5)$$
The set of Finsler connections is in a one--to--one correspondence
with the set of linear connections on $TM$ which verify
$$D_XP=0,\ D_XF=0\mbox{ for any vector field $X$ on $TM$}.\leqno(5.6)$$
By the very definition, a vector field $B$ on $TM$ is parallel
with respect to $D$ if
$$D_XB=0,\leqno(5.7)$$
and is concurrent if
$$D_XB=-X,\mbox{ for any vector field $X$ on $TM$}.\leqno(5.8)$$
Let $(\delta_i,\dot\partial_i)$ be the usual adapted basis for the
decomposition (5.3). The above mentioned one--to--one
correspondence is established by
$$\begin{array}{ll}
D_{\delta_k}\delta_j=L^i_{jk}\delta_i,
&D_{\dot\partial_k}\delta_j=V^i_{jk}\delta_i,\vspace*{1,4mm}\\
D_{\delta_k}\dot\partial_j=L^i_{jk}\dot\partial_i,
&D_{\dot\partial_k}\dot\partial_j=V^i_{jk}\dot\partial_i,\end{array}\leqno(5.9)$$
for $D\leftrightarrow F\Gamma(N)=(N^i_j,L^i_{jk},V^i_{jk}).$

It is obvious that (5.7) is equivalent to
$$D_{\delta_k}B=0,\ D_{\dot\partial_k}B=0,\leqno(5.7)'$$
and (5.8) is equivalent to
$$D_{\delta_k}B=-\delta_k,\ D_{\dot\partial_k}B=-\dot\partial_k.\leqno(5.8)'$$
Let now be $B=B^i(x,y)\delta_i+\hat B^i(x,y)\dot\partial_i.$ Then $(5.7)'$ is
equivalent by virtue of (5.9) with
$$B^i_{|k}=0,\ B^i{\big|}_k=0,\
\hat B^i_{|k}=0,\ \hat B^i{\big|}_k=0.\leqno(5.7)''$$
One may associate to $B^i(x,y)$ at least the following three vector
fields on $TM:$ $B^i\delta_i,$ $B^i\dot\partial_i,$
$B^i\delta_i+B^i\dot\partial_i$ and it is obvious by $(5.7)''$ that
$B^i(x,y)$ is parallel in the sense of (5.1) if and only if at least one from
these vector fields on $TM$ is parallel with respect to $D$. Thus
(5.1) is in agreement with the usual definition of parallelism.

Let us make a similar analysis for concurrent Finsler vector fields.
By (5.8), $B$ is concurrent on $TM$ if and only if
$$B^i_{|k}=-\delta^i_k,\ B^i_{\big|_k}=0,\ \hat B^i_{|k}=0,\
\tilde B^i{\big|}_k=-\delta^i_k.\leqno(5.10)$$
Now we assume that $D$ or $F\Gamma(N)$ is of Cartan type, that is,

$$y^i_{|k}=0,\ y^i{\big|}_k=\delta^i_k.\leqno(5.11)$$
The tensors $y^i_{|k}$ and $y^i{\big|}_k$ are called $h$--deflection
and $v$--deflection tensors, respectively. The equations (5.11) hold
for all four remarkable connections in Finsler spaces.

If moreover we assume that $\hat B^i$ is positively homogeneous of
degree 1, a natural assumption in Finslerian setting, writing $\hat
B^i{\big|}_k=-\delta^i_k$ in the form $\dot\partial_k\hat B^i+V^i_{jk}\hat
B^j=-\delta^i_k$ and contracting it by $y^k$ it results using (5.11)
that $y^k\dot\partial_k\hat B^i=-y^i.$ Thus by the Euler theorem, $\hat
B^i=-y^i$ and then $\hat B^i_{|k}=0$ reduces to $y^i_{|k}=0$ i.e.
the first equation in (5.11). Concluding, if we associate to the Finsler
vector field $B^i(x,y)$ the vector field
$B=B^i(x,y)\delta_i-y^i\dot\partial_i$ on $TM$, we find that $(B^i(x,y))$ is
concurrent in the sense of (5.2) if and only if $B$ is concurrent by the new
definition of concurrence on any manifold. In other words, the condition
(5.2) is in agreement with the notion of concurrence for vector
fields.

\section{The metric $^*g_{ij}$ with $B^i(x,y)$ a concurrent Finsler vector
field}
\setcounter{equation}{0}

In this section we are dealing with the $GL$--metric $^*g_{ij}$
given by (3.1) for $B^i(x,y)$ a concurrent Finsler vector field with
respect to the Cartan connection $C\Gamma$ of $F^n$ i.e.
$$B^i_{|j}=-\delta^i_j,\ B^i{\big|}_j=0.\leqno(6.1)$$
First we notice some results on concurrent Finsler vector fields due
to M. Matsumoto and K. Eguchi \cite{8}.

If $B^i(x,y)$ is concurrent we have with respect to $C\Gamma:$\bigskip
\begin{itemize}
\item[(6.2)] $B_{i|j}=-g_{ij},\ B_i{\big|}{}\!_j=0,$
\item[(6.3)] $B^hR_{hijk}=0,\ B^hP_{hijk}+C_{ijk}=0,\
B^hS_{hijk}=0,$
\item[(6.4)] $B^iC_{ijk}=C^s_{jk}B_s=0,$
\item[(6.5)] $B^i=B^i(x)$ and $B_i=B_i(x)$ i.e. $B^i$ and $B_i$ are
functions on position only,
\item[(6.6)] $\partial_iB_j=\partial_jB_i=F^s_{ij}B_s-g_{ij},\
\partial_kB^i=-\delta^i_k-F^i_{sk}B^k.$\end{itemize}
In these circumstancies a direct calculation yields
$$\begin{array}{l}
\Lambda\ijk=\displaystyle\frac{^*\sigma}{2\sigma}B^i(\sigma_kB_j+
\sigma_jB_k+\sigma(B^s\sigma_s)B_jB_k-2\sigma
g_{jk})-\frac12\sigma^iB_jB_k\\
\bld\circ\Lambda{}\!\ijk=\displaystyle\frac{^*\sigma}{2\sigma}B^i
(\dot\sigma_kB_j+\dot\sigma_jB_k+
\sigma(B^s\dot\sigma_s)B_jB_k-\frac12\dot\sigma^iB_jB_k,\mbox{ where}
\end{array}\leqno(6.7)$$
$$\sigma_k:=\delta_k\sigma,\ \dot\sigma_k:=\dot\partial_k\sigma,\
\sigma^i=g^{i_k}\sigma_k,\
\dot\sigma^i=g^{i_k}\dot\sigma_k.\leqno(6.7)'$$
Looking at (6.7) we see that the simplest case is given by
$$\sigma_k=0,\ \dot\sigma_k=0.\leqno(6.8)$$
From (6.8) it results that $\sigma$ is a constant $c$. And
$^*F^2:={}^*g_{ij} y^iy^j$ takes the form
$$^*F^2=F^2+c\beta^2,\ \beta=B_i(x) y^i.\leqno(6.9)$$
Thus, for $c>0,$ $^*F$ is a new Finsler function which is obtained from $F$
by a particular $\beta$--change.

The case $c=1$ is studied in \cite{8}.

Further on we have
$$^*F^i_{jk}=F^i_{jk}-{}^*\sigma B^ig_{jk},\
^*C^i_{jk}=C^i_{jk}.\leqno(6.10)$$

{\bf Remark 6.1.} The Cartan nonlinear connection of $^*F^n=(M,^*F)$ is given by
$N^i_j=\bld cN{}\!^i_j-\bld*\sigma B^iy_j$ i.e. the difference
tensor is $A^i_j=\bld*\sigma B^iy_j.$ Inserting it in $(4.3)'$ we
find $^sF^i_{jk}={}^*F^i_{jk}$. Therefore, in the geometry of $^*F^n$ we may
equally use $\bld c{N}\!^i_j$ or $N^i_j.$

By (6.10) we immediately get
$$^*S_{ijkh}=S_{ijkh}.\leqno(6.11)$$
Again by (6.10) but after a long calculation one finds
$$^*R_{ijkh}=R_{ijkh}+{}^*\sigma(g_{ik}g_{jh}-g_{ih}g_{jk}).\leqno(6.12)$$
This suggests us to take into consideration the case when $F^n$ is
$h$--isotropic i.e. there exists a constant $K$ such that
$R_{ijkh}=K(g_{ik}g_{jh}-g_{ih}g_{jk}).$ A contraction of this
last equation by $B^i$ gives for $K\ne0,$ $B_kg_{jh}-B_hg_{jk}=0$ in
virtue of (6.3). A new contraction by $B^k$ yields
$B^2g_{jh}=B_jB_h$ which contradicts the condition ${\rm rank}\:(g_{ij})=n>1.$
Thus we have\bigskip

{\bf Theorem 6.1.} {\it If $F^n$ is $h$--isotropic, then it does
not admit any con\-current Finsler vector field.}\bigskip

The proof of the following two theorems are the same as for $c=1$ (see
Theorems 14 and 15 in \cite{8}).\bigskip

{\bf Theorem 6.2.} {\it If $F^n$ admits a concurrent Finsler vector
field, then there is no a Finsler vector field which to be
concurrent with respect to $^*F$ given by $(6.9).$}

{\bf Theorem 6.3.} {\it If $F^n$ admits a concurrent Finsler vector
field and is $R3$--like, then $^*F^n=(M,^*F)$ with $^*F$ from
$(6.5)$ is also $R3$--like.}\bigskip

Now we consider a more complicated case
$$\sigma_k=0,\ \dot\sigma_k\ne0.\leqno(6.13)$$

{\bf Remark 6.2.} The equation $\sigma_k:=\derp{\sigma}{x^h}-\bld
cN{}\!^s_k\derp\sigma{y^s}=0$ is known as Tavakol--Van der Berg
equation. A solution of it is for instance $\sigma=aF^2$ for $a\in\R.$
For more details see \cite{12}.

Now (6.10) is replaced by
$$\begin{array}{l}
^*F\ijk=F\ijk-{}^*\sigma B^ig_{jk}\\
^*C\ijk=C\ijk+\displaystyle\frac{^*\sigma}{2\sigma}
B^i(\dot\sigma_kB_j+\dot\sigma_jB_k+\sigma
(B^s\dot\sigma_s)B_jB_k)-\frac12\dot\sigma^iB_jB_k.
\end{array}\leqno(6.14)$$
The Remark 6.1 is still valid for this case. Precisely, if we ask
for the vanishing of the $h$--deflection of $^*F\Gamma(\bld cN)$,
then $^*N^i_j=\bld cN{}\!^i_j-\bld*\sigma B^iy_j$ and so
$^sF\Gamma(\bld*N)$ coincides with $^*F\Gamma(\bld cN).$

Now we notice\medskip

\noindent$\!\!\!\begin{array}{lcl}
(6.15)&\ \ \
&^*C_j=C_j+\displaystyle\frac{^*\sigma B^2}{2\sigma}\dot\sigma_j,\
C_j:=C^i_{ji},\vspace*{1,4mm}\\
(6.16)&&^*C_{jik}=C_{jik}+\displaystyle\frac12(\dot\sigma_k
B_iB_j+\dot\sigma_jB_iB_k-\dot\sigma_iB_jB_k).\end{array}$\medskip

\noindent A long calculation yields
$$\begin{array}{c}
^*R_{jskh}=R_{jskh}+{}^*\sigma(g_{jk}g_{sh}-g_{jh}g_{sk})+\\
+\displaystyle\frac{^*\sigma}\sigma B_s(\partial_k\sigma\cdot
g_{jh}-\partial_h\sigma\cdot g_{jk})+\frac12
B_jB_sR^q_{kh}\dot\sigma_q.\end{array}\leqno(6.17)$$
Let us assume that $F^n$ is a locally Minkowski space. Then
$R_j{}^i{}_{kh}=0$ and $C^i_{jk|h}=0.$ In a local chart in which
$g_{ij}$ do not depend on $x$ we have $\bld cN{}\!^i_j=0$ and so
$\partial_k\sigma=\bld cN{}\!^p_j\dot\sigma_p=0$ i.e. $\sigma$ does not
depend on $x$.

The equation (6.17) reduces to
$$^*R_{jskh}={}^*\sigma(g_{jk}g_{sh}-g_{jh}g_{sk}).\leqno(6.18)$$
It takes also the form
$$\begin{array}{lcl}
^*R_{jskh}={}^*\sigma({}^*g_{jk}{}^*g_{sh}-
{}^*g_{jh}{}^*g_{sk})+
\sigma{}^*\sigma(B_jB_{hsk}+B_sB_{kjh})\mbox{ for}\\
B_{hsk}:=B_hg_{sk}-B_kg_{sh}.\end{array}\leqno(6.18)'$$
We notice that $B_{hsk}$ is never vanishing since otherwise a
contraction by $B^h$ gives a contradiction with
${\rm rank}\:(g_{ij})=n>1.$

\section{A Beil metric for a Finsler space with $(\alpha,\beta)$--metric}
\setcounter{equation}{0}

Here we consider again the Beil metric described in Remark 3.5.
Let $F^n$ be a Finsler space with an $(\alpha,\beta)$--metric. A
natural Beil metric is then
$$^*g_{ij}(x,y)=a_{ij}(x)+\sigma(x,y) b_i(x) b_j(x).\leqno(7.1)$$
Let $\gamma^i_{jk}$ be the Christoffel symbols for $a_{ij}(x).$ Then $\bld
cN{}\!^i_j=\gamma^i_{jk} y^k=:\gamma^i_{j0}$ and the triple
$\Gamma=(\gamma^i_{j0},\gamma^i_{jk},0)$ may be thought of as a
Finsler connection.

We have\bigskip

{\bf Theorem 7.1.} {\it If $b_i(x)$ is parallel and $\sigma$ is
covariant constant with respect to $\Gamma$, then $\Gamma$ is like
Chern--Rund connection for $(^*g_{ij}).$}

\medskip{\bf Proof.} Let $_{;k}$ denote the $h$--covariant
derivative with respect to $\Gamma.$ Notice that $v$--covariant
derivative is just the derivative with respect to $y$. Our
hypothesis read
$$b_{i;k}=0,\ \delta_k\sigma=0,\
\delta_k=\partial_k-\gamma^s_{k0}\dot\partial_s.\leqno(7.2)$$
Then we easily get
$$\begin{array}{l}
^*g_{i;jk}=(\delta_k\sigma)b_ib_j=0\vspace*{1,4mm}\\
^*g_{ij,k}=(\dot\partial_k\sigma)b_ib_j=2^*C_{ikj}.\end{array}\leqno(7.3)$$
Thus $\Gamma$ is $h$--metrical and no metrical for $^*g_{ij}.$ Hence
it is similar to the Chern--Rund connection from Finsler
geometry.\hfill{\bf q.e.d.}\medskip

The Chern--Rund connection is a remarkable one in Finsler geometry
(\cite{1}). Notice that its $h$--deflection vanishes.

From now on we assume $b_{i;k}=0$ and $\delta_k\sigma=0.$

A direct calculation yields
$$\begin{array}{l}
^*F\ijk=\gamma\ijk,\\
^*C\ijk=\displaystyle\frac{^*\sigma}{2\sigma}
b^i(\dot\sigma_kb_j+\dot\sigma_jb_k+\sigma
(b^h\dot\sigma_h)b_jb_k)-\frac12\dot\sigma^ib_jb_k.
\end{array}\leqno(7.4)$$
The first equation in (7.4) is important in many respects. For
instance using it we find the $h$--curvature of $^*F\Gamma(\bld cN)$ in
the form
$$^*R_h{}^i{}_{jh}=\gamma_h{}^i{}_{jh}+
\bld\circ\Lambda_h{}^i{}_sR^s_{jk},\leqno(7.5)$$
where $\gamma_h{}^i{}_{jh}$ is the curvature tensor of $a_{ij}(x)$ and
$R^i_{jk}=\gamma_0{}^i{}_{jk}$.  Here, as before, the index $0$
indicates the contraction by $y$. Consequently, (7.5) takes the
form
$$^*R_h{}^i{}_{jk}=(\delta^i_s\delta^r_h+
\bld\circ\Lambda{}\!^i{}_{hs}y^r)\gamma_r{}^s{}_{jk}.\leqno(7.6)$$
From Ricci identities we find $\gamma_i{}^s{}_{jk}b_s=0$ and from (7.5) we
deduce
$$^*R_{hijk}=\gamma_{hijk}+\displaystyle\frac12
b_hb_i\gamma_0{}^s{}_{jk}\dot\sigma_s.\leqno(7.7)$$
As for Ricci curvatures one finds
$$^*R_{ij}=r_{ij},\leqno(7.8)$$
where $r_{ij}$ is the Ricci curvature for $(a_{ij}(x)).$ From here it results
$$^*R=r,\leqno(7.9)$$
where $^*R$ and $r$ are the scalar curvatures for $(^*g_{ij})$ and
$(a_{ij}(x))$, respectively.

So, the $h$--Einstein tensor field of $^*g_{ij}$ i.e.
$^*E_{ij}=^*R_{ij}-\displaystyle\frac12{}^*R^*g_{ij}$
is related to the Einstein
tensor $E_{ij}$ of $a_{ij}(x)$ by
$$^*E_{ij}=E_{ij}+\displaystyle\frac{\sigma r}2b_ib_j.\leqno(7.10)$$
Consequently, the $h$--Einstein equation for $GL$ i.e.
$^*E_{ij}=\kappa^*\tau_{ij}$ with $\kappa\in\R$ reduces to
$$r_{ij}-\displaystyle\frac r2a_{ij}=\kappa\tau_{ij},\leqno(7.11)$$
where
$$\tau_{ij}={}^*\tau_{ij}-
\displaystyle\frac{\sigma r}{2\kappa}b_ib_j.\leqno(7.12)$$
The equation (7.11) is the
Einstein equation for $(M,a_{ij}(x))$ but with the
energy--momentum tensor influenced
by a field described by $b_i$. In the the unified theory of R.G. Beil the
term $b_ib_j$ in (7.12) is a "matter term" which could be the energy density of the self-field
of a charged object.

\noindent{\footnotesize{University "Al.I.Cuza" Ia\c si\\ Faculty of mathematics\\ 6600, Ia\c si, Romania}}\medskip

\noindent{\footnotesize{Hokkaido Tokai University\\ Minami-ku, Minami-sawa\\ 5-1, Sapporo 005, Japan}}

\newpage


\runningauthor={M. ANASTASIEI}
\runningtitle={LOCALLY CONFORMAL K\"{A}HLER STRUCTURES ON TANGENT MANIFOLD OF A SPACE FORM}
\noindent
\baselineskip 8pt
\noindent{\footnotesize{Libertas Math.}}
\hfill\break
{\footnotesize{19, 71--76, 1999}}
\vskip 2cm
\baselineskip 11.5pt plus .15pt
\centerline{\bf\Large LOCALLY CONFORMAL}
\vskip .2cm
\centerline{\bf\Large K\"{A}HLER STRUCTURES ON TANGENT}
\vskip .2cm
\centerline{\bf\Large MANIFOLD OF A SPACE FORM}
\vskip .5cm
\centerline{\bf by Mihai ANASTASIEI}
\vskip .5cm

\begin{abstract}
A set of locally conformal K\"{a}hler structures on tangent manifold $TM$ of a
space form $M$ is pointed out. This is found in a study of a type of Sasaki
metric whose second term is a special deformation of the first one.
Introducing an adequate almost complex structure we find at first a large
class of locally conformal almost K\"{a}hler structures on $TM$ for $M$ a
(pseudo)- Riemannian manifold. When $M$ is a space form, a subset of it is
made of locally conformal K\"{a}hler structures. One of them was found by R.
Miron in [3].
\end{abstract}

\setcounter{section}{0}
\section{Introduction}
\setcounter{equation}{0}
\setcounter{theorem}{0}
\setcounter{corollary}{0}

Let $(M,g)$ be a (pseudo)-Riemannian manifold and $\triangledown$ its
Levi-Civita connection. In a local chart $(U,(x^{i}))$ we set $g_{ij}=g(\partial_{i},\partial_{j})$, where $\partial_{i}:\dfrac{\partial}{\partial x^{i}}$ and we denote by $\gamma_{jk}^{i}(x)$ the Christoffel
symbols giving $\nabla$. Let $(x^{i},y^{i})\equiv(x,y)$ be the local
coordinates on the manifold $TM$ projected on $M$ by $\tau$. The indices $i,j,k...$ will run from $1$ to $n={\rm{dim}} M$.

The functions $N_{j}^{i}(x,y):=\gamma_{jk}^{i}(x)y^{k}$ are the local
coefficients of a nonlinear connection, that is the local vector fields $\delta_{i}=\partial_{i}-N_{i}^{k}(x,y)\stackrel{.}{\partial_{k}}$, where $\stackrel{.}{\partial_{k}}:\dfrac{\partial}{\partial y^{k}}$ span a distribution on $TM$ called horizontal which is supplementary to the
vertical distribution $u\rightarrow V_{u}TM=\ker\tau_{*,u},u\in TM$. Let
us denote by $u\rightarrow H_{u}TM$ the horizontal distribution and let $(\delta_{i},\stackrel{.}{\partial_{i}})$ be the basis adapted to the decomposition $T_{u}TM=H_{u}TM\oplus V_{u}TM,u\in TM$. The basis dual of it is $(dx^{i},\delta y^{i})$ with $\delta y^{i}=dy^{i}+N_{k}^{i}(x,y)dx^{k}$.

The Sasaki metric on $TM$ is as follows
\begin{equation}
G_{S}=g_{ij}(x)dx^{i}\otimes dx^{j}+g_{ij}(x)\delta y^{i}\otimes\delta
y^{j}.
\end{equation}

If in the second term of $G_{S}$ one replaces $g_{ij}(x)$ with the
components $h_{ij}(x,y)$ of a generalized Lagrange metric (see Ch. X in [4])
one gets a type of Sasaki metric
\begin{equation}
G(x,y)=g_{ij}(x)dx^{i}\otimes dx^{j}+h_{ij}(x,y)\delta y^{i}\otimes\delta
y^{j}.
\end{equation}

In particular, $h_{ij}(x,y)$ could be a deformation of $g_{ij}(x)$, a case
studied by the present author and H. Shimada in [1].

In this paper we are concerning with the metrical structure (1.2) in the
case when $h_{ij}(x,y)$ is the following special deformation of $g_{ij}(x)$
\begin{equation}
h_{ij}(x,y)=a(L^{2})g_{ij}(x)+b(L^{2})y_{i}y_{j},
\end{equation}
where $L^{2}=g_{ij}(x)y^{i}y^{j},y_{i}=g_{ij}(x)y^{j}$ and $a,b:{\rm{Im}}(L^{2})\subseteq \Bbb{R}_{+}\longrightarrow \Bbb{R}_{+}$ with $a>0,b\geq 0$.

For $b=0$ and $a=\dfrac{c^{2}}{L^{2}}$ for any constant $c$, the metrical
structure (1.2), (1.3) was studied by R. Miron in [3] as an homogeneous lift
of $g_{ij}(x)$ to $TM$.

In the following Section we introduce an almost complex structure which
paired with $G$ given by (1.2), (1.3) provides a large set of almost
Hermitian structures on $TM$. Then, in Section 3 we show that all these
structures are locally conformal almost K\"{a}hler structures. Finally, we find
in Section 4 that, when $(M,g)$ is of constant curvature, a part of them are
locally conformal K\"{a}hler structures.

\section{Some almost Hermitian structures on $TM$}
\setcounter{equation}{0}
\setcounter{theorem}{0}
\setcounter{corollary}{0}

Let $F_{S}$ be the almost complex structure on $TM$ given in the adapted
basis ($\delta_{i},\stackrel{.}{\partial_{i}})$ by
\begin{equation}
F_{S}(\delta_{i})=-\stackrel{.}{\partial_{i}},F_{S}(\partial_{i})=\delta_{i}.
\end{equation}

It is well known that the pair $(G_{S}$, $F_{S})$ is an almost K\"{a}hler
structure on $TM$, that is $G_{S}(F_{S}X,F_{S}Y)=G_{S}(X,Y)$ and the 2-form
$$\Omega (X,Y)=G_{S}(F_{S}(X),Y)\textrm{ is closed, }X,Y\in \chi (M).$$

The pair $(G,F_{S})$ with $G$ given by (1.2), (1.3) is no longer an
almost Hermitian structure. We look for a new almost complex structure which
paired with $G$ to provide an almost Hermitian structure. We modify $F_{S}$
to a linear map $F$ given in the basis ($\delta_{i},\stackrel{.}{\partial_{i}})$ as follows
\begin{equation}
F(\delta_{i})=(\alpha\delta_{i}^{k}+\beta y_{i}y^{k})\stackrel{.}{\partial_{k}},F(\stackrel{.}{\partial_{j}})=(\gamma \delta_{j}^{h}+\delta y_{j}y^{h})\delta_{h},
\end{equation}
where $\alpha,\beta,\gamma,\delta$ are functions on $TM$ to be
determined. The condition $F^{2}=-I$ (identity) leads to
\begin{equation}
\alpha\gamma=-1,\alpha\delta+\beta\gamma +\beta\delta L^{2}=0.
\end{equation}

Then the condition $G(F(X),F(Y))=G(X,Y)$ gives
\begin{equation}
a\alpha^{2}=1,\gamma^{2}=a,2\gamma\delta +\delta^{2}L^{2}=b,(2\alpha
\beta+\beta ^{2}L^{2})(a+bL^{2})+b\alpha^{2}=0
\end{equation}

The solution of the system of equations (2.3), (2.4) is
\begin{equation}
\alpha=-\dfrac{1}{\sqrt{a}},\beta=\dfrac{\sqrt{a}+\sqrt{a+bL^{2}}}{L^{2}\sqrt{a(a+bL^{2})}},\gamma =\sqrt{a},\delta =-\dfrac{\sqrt{a}+\sqrt{a+bL^{2}}}{L^{2}}.
\end{equation}

We notice that for $b=0$, besides the solution provided by (2.5), that is
\begin{equation}
\alpha=-\dfrac{1}{\sqrt{a}},\gamma =\sqrt{a},\beta =\dfrac{2}{L^{2}\sqrt{a}},\delta=-\dfrac{2\sqrt{a}}{L^{2}},
\end{equation}
there exists also the solution
\begin{equation}
\alpha =-\dfrac{1}{\sqrt{a}},\gamma=\sqrt{a},\beta=0,\delta=0.
\end{equation}

Let us make the substitution $a\longrightarrow\dfrac{a^{2}}{L^{2}}$, $b\longrightarrow \dfrac{b^{2}-a^{2}}{L^{4}}$.
Then (2.5) and (2.6) are unified to
\begin{equation}
\alpha =-\dfrac{L}{a},\beta =\dfrac{a+b}{abL},\gamma =\dfrac{a}{L},\delta=-%
\dfrac{a+b}{L^{3}},b\geq a>0
\end{equation}
and (2.7) modifies to
\begin{equation}
\alpha =-\dfrac{L}{a},\gamma =\dfrac{a}{L},\beta =\delta=0.
\end{equation}

The metric $G$ takes the form
\begin{equation}
\begin{array}{l}
G_{a,b}(x,y)=g_{ij}(x)dx^{i}\otimes dx^{j}+\left(\dfrac{a^{2}}{L^{2}}g_{ij}(x)+\dfrac{b^{2}-a^{2}}{L^{4}}y_{i}y_{j}\right)\delta y^{i}\otimes\delta y^{j},\\ \hfill b\geq a>0.
\end{array}
\end{equation}

Let $F_{a,b}$ be the almost complex structures given by (2.2), (2.8) and $F_{a}$ those given by (2.2), (2.9). Then the pairs $(G_{a,b},F_{a,b})$ and $(G_{a,a},F_{a})$ are almost Hermitian structures on $TM$.

For $a^{2}=\dfrac{L^{2}}{1+L^{2}},b=L^{2}$, the metric $G_{a,b}(x,y)$ is the
Cheeger-Gromoll metric, [5],[6]
\begin{equation}
G_{CG}(x,y)=g_{ij}(x)dx^{i}\otimes dx^{j}+\dfrac{1}{1+L^{2}}
(g_{ij}(x)+y_{i}y_{j})\delta y^{i}\otimes \delta y^{j}.
\end{equation}

If $a^{2}=\varphi^{\prime}L^{2},b^{2}=L^{2}(\varphi^{\prime}+2\varphi^{\prime\prime}L^{2})$ for $\varphi:\Bbb{R}_{+}\longrightarrow \Bbb{R}_{+}$ with $\varphi^{\prime}(t)\neq 0,t\in{\rm{Im}}(L^{2})$, one
obtains the Antonelli - Hrimiuc metrical structure, [2]
\begin{equation}
G_{AH}(x,y)=g_{ij}(x)dx^{i}\otimes dx^{j}+(\varphi^{\prime}g_{ij}(x)+2\varphi^{\prime\prime}y_{i}y_{j})\delta y^{i}\otimes\delta
y^{j}.
\end{equation}

\section{Locally conformal almost K\"{a}hler structures on TM}
\setcounter{equation}{0}
\setcounter{theorem}{0}
\setcounter{corollary}{0}

Let $\Omega (X,Y)=G_{a,b}(F_{a,b}X,Y),X,Y\in \chi(TM)$ be the 2-form
associated to the almost Hermitian structure $(G_{a,b},F_{a,b})$.
\begin{theorem}
The almost Hermitian structures $(G_{a,b},F_{a,b})$ are
locally conformal almost K\"{a}hlerian structures, that is
\begin{equation}
d\Omega=\Omega\wedge\theta,\theta=\dfrac{2a^{\prime}L+b}{aL}dL.
\end{equation}
\end{theorem}
\begin{proof}
We shall check (3.1) on the basis $(\delta_{i},\stackrel{.}{\partial_{i}})$. If we rewrite (2.2) in the form
\begin{equation}
F(\delta_{i})=A_{i}^{k}\stackrel{.}{\partial _{k}},F(\stackrel{.}{\partial_{i}})=B_{j}^{h}\delta_{h},
\end{equation}
we easily get
\begin{equation}
\Omega(\delta_{i},\delta_{j})=0,\Omega(\delta_{i},\stackrel{.}{\partial_{j}})=A_{i}^{k}h_{kj},\Omega(\stackrel{.}{\partial _{j},}\delta_{i})=B_{j}^{k}g_{ki},\Omega\stackrel{.}{(\partial_{i},}\stackrel{.}{\partial_{j},})=0,
\end{equation}
with $A_{i}^{k}h_{kj}+B_{j}^{k}g_{ki}=0$.

Thus $\Omega $ is completely determined by
\begin{equation}
\Omega_{ij}:=B_{j}^{k}g_{ki}=\gamma g_{ij}+\delta y_{i}y_{j};\Omega
=\Omega_{ij}\delta y^{i}\wedge dx^{j}.
\end{equation}

Next we have the following essential components of $d\Omega$:

$d\Omega(\delta_{i},\delta_{j},\stackrel{.}{\partial_{k}})=\delta_{j}\Omega_{ik}-\gamma_{ki}^{s}\Omega_{sj}-\delta_{i}\Omega_{jk}
-\gamma_{kj}^{s}\Omega_{si}$,

$d\Omega(\delta_{i},\stackrel{.}{\partial_{j}},\stackrel{.}{\partial_{k}})=\stackrel{.}{\partial_{j}}\Omega _{ik}-\stackrel{.}{\partial_{k}}\Omega_{ij}$.

Now we consider the Berwald connection $(N_{j}^{i}=\gamma_{kj}^{i}(x)y^{k}$, $\gamma_{kj}^{i}(x),0)$ on $TM$ (see Ch.8 in [4]) and
denote by $|k$ its $h-$covariant derivative. Thus because of $\Omega_{jk|i}=\delta_{i}\Omega_{jk}-\gamma_{ji}^{s}\Omega_{sk}-\gamma
_{ki}^{s}\Omega_{js}$, we have $d\Omega(\delta_{i},\delta_{j}\stackrel{.},{\partial_{k}})=\Omega_{ki|j}-\Omega_{kj|i}$.

The following formulae are verified by a direct calculation.
\begin{equation}
g_{ij|k}=0,y_{|k}^{j}=0,y_{i|k}=0,\delta_{k}L^{2}=0,\delta_{k}\psi
(L^{2})=0,
\end{equation}
$\stackrel{.}{\partial_{k}}y_{i}=g_{ik},\stackrel{.}{\partial_{k}} L^{2}=2y_{k},\stackrel{.}{\partial_{k}}\psi(L^{2})=2y_{k}\psi^{\prime}(L^{2})$,\\
for any $\psi:{\rm{Im}}(L^{2})\subseteq R_{+}\longrightarrow R_{+}$.

Using (3.5) it immediately results $\Omega_{kj|i}=0$ and so $d\Omega(\delta_{i},\delta_{j},\stackrel{.}{\partial_{k}})=0$. Consequently, $d\Omega$ is completely determined by $d\Omega(\delta_{i},\stackrel{.}{\partial_{j},}\stackrel{.}{\partial _{k}})=\stackrel{.}{(\partial_{j}}\gamma)g_{ik}-(\stackrel{.}{\partial _{k}}\gamma)g_{jk}+(\stackrel{.}{\partial
_{j}\delta)y_{k}y_{i}-(\stackrel{.}{\partial_{k}}}\delta)y_{j}y_{i}+\delta(g_{ij}y_{k}-g_{ik}y_{j})$.

Inserting here $\stackrel{.}{\partial_{j}}\gamma,stackrel{.}{\partial_{j}}\delta$ with $\gamma,\delta$ from (2.8) one arrives to
\begin{equation}
d\Omega(\delta_{i},\stackrel{.}{\partial_{j},}\stackrel{.}{\partial_{k}})=(2\gamma^{\prime}-\delta)(g_{ik}y_{j}-g_{ij}y_{k}) =\dfrac{2a^{\prime }L^{2}+b}{L^{3}}(g_{ik}y_{j}-g_{ij}y_{k}).
\end{equation}

Let be $\theta_{0}=dL^{2}=2y_{i}\delta y^{i}$. Thus $\theta_{0}(\delta_{i})=0$ and $\theta_{0}\stackrel{.}{(\partial_{j}})=2y_{j}$. Evaluating $\Omega\wedge\theta_{0}$ on the basis $(\delta_{i},\stackrel{.}{\partial_{i}})$ one finds the essential component
\begin{equation}
\Omega\wedge\theta_{0}(\delta_{i},\stackrel{.}{\partial_{j}},\stackrel{.}{\partial_{k}})=2(\Omega_{ik}y_{j}-\Omega _{ij}y_{k})=\dfrac{2a }{L}(g_{ik}y_{j}-g_{ij}y_{k}).
\end{equation}
Comparing (3.6) with (3.7) one obtains $d\Omega=\dfrac{2a^{\prime}L^{2}+b}{2aL^{2}}\Omega\wedge\theta_{0}$ which is just (3.1).
\end{proof}

Obviously $\theta$ is globally defined. Moreover, $\theta$ is closed. This
fact can be directly verified using (3.5) or by differentiating (3.1).

Looking at (3.6) we notice that contracting $g_{ik}y_{j}-g_{ij}y_{k}=0$ with
$g^{ik}$ one gets $(n-1)y_{j}=0$ which is a contradiction. Thus we have

\begin{theorem}
The almost Hermitian structures $(G_{a,b},F_{a,b})$ are almost K\"{a}hler structures if and only if
\begin{equation}
2a^{\prime }L^{2}+b=0,
\end{equation}
holds good.
\end{theorem}

We put $t=L^{2}$ and think (3.8) as a first order differential equation: $2ta^{\prime }(t)+b(t)=0.$ Its general solutions is $a(t)=c-\dfrac{1}{2}\dd\int\dfrac{b(t)}{t}dt$ for a constant $c$. Thus for various functions $b$ we find
a set of pairs $(a,b)$ for which (3.8) holds. Choosing among these pairs
those which verify $b\geq a>0$ we find a set of almost K\"{a}hler structures on
$TM$. For instance, if we take $b(t)=2t$ it results $a(t)=c-t$ and $b\geq
a>0$ holds if $\dfrac{c}{3}\leq L^{2}(x,y)<c$, for $c>0$. When $a=b$, the
equation (3.8) has the general solution $a(t)=\dfrac{c}{\sqrt{t}}$. It follows

\begin{corollary}
The almost Hermitian structures $(G_{a,a},F_{a,a})$ are almost K\"{a}hler structures if and only if $a(L^{2})=\dfrac{c}{\sqrt{L^{2}}},c>0$.
\end{corollary}

The almost Hermitian structures $(G_{a,a},F_{a})$ have to be separately
considered. Repeating for them the proof of Theorem 3.1 one obtains

\begin{theorem}
The almost Hermitian structures $(G_{a,a},F_{a})$ are locally conformal almost K\"{a}hler structures, that is
\begin{equation}
d\Omega=\Omega\wedge\theta,\theta=\dfrac{2a^{\prime}L-a}{aL}dL.
\end{equation}
\end{theorem}

The following corresponds to Theorem 3.2

\begin{theorem}
The almost Hermitian structures $(G_{a,a},F_{a})$ are almost K\"{a}hler structures if and only if $a=c\sqrt{L^{2}},c>0$.
\end{theorem}

\begin{proof}
The almost K\"{a}hler condition is now $2a^{\prime }L^{2}-a=0.$
Integrating the equation $2a^{\prime }t-a=0$ one gets $a=c\sqrt{t}$.
\end{proof}

\noindent{\it Remark.}
For $a=c\sqrt{L^{2}},c>0,G_{a,a}$ is very close to $G_{S}$ which is obtained for $c=1$.

\section{Locally conformal K\"{a}hler structures on TM}
\setcounter{equation}{0}
\setcounter{theorem}{0}
\setcounter{corollary}{0}

In order to find conditions that $(G_{a,b},F_{a,b})$ be a locally conformal
K\"{a}hler structure we have to put zero for the Nijenhuis tensor field of $F:=F_{a,b}$,
\begin{equation}
N_{F}=[FX,F]-F[FX,Y]-F[X,FY]-[X,Y],X,Y\in \chi (TM).
\end{equation}

As the evaluation of $N_{F}$ on the basis $(\delta_{i},\stackrel{.}{\partial_{i}})$ is in general very complicated we confine ourselves to the structures $(G_{a,a},F_{a})$. In this case, the conditions
\begin{equation}
N_{F}(\delta_{i},\delta_{j})=0,N_{F}(\delta_{i},\stackrel{.}{\partial_{j}})=0,N_{F}\stackrel{.}{(\partial_{j},}\stackrel{.}{\partial _{k}})=0,
\end{equation}
are equivalent with six equations. Three of them are identities because of $\delta_{i}\alpha=\delta_{i}\gamma=0$ and the other three are each one equivalent with
\begin{equation}
R_{ij}^{k}=\dfrac{2a^{\prime}L^{2}-a}{a^{3}}(y_{j}\delta_{i}^{k}-y_{i}\delta_{j}^{k}),
\end{equation}
where $R_{ij}^{k}=R_{sij}^{k}(x)y^{s}$ and $R_{sij}^{k}$ is the curvature
tensor of $\nabla$.

By a contraction with $g_{rk}$ the Eq. (4.3) reduces to
\begin{equation}
R_{srij}(x)y^{s}=\dfrac{2a^{\prime}L^{2}-a}{a^{3}}(g_{js}g_{ri}-g_{is}g_{rj})y^{s}.
\end{equation}

The Eq. (4.4) remember us the condition that $(M,g)$ is of constant curvature (space form). It suggests us to look for functions $a$ such that $\dfrac{2a^{\prime }L^{2}-a}{a^{3}}=k$, where $k$ is a constant. For $t=L^{2}$, solving the Bernoulli equation $a^{\prime }=\dfrac{1}{2t}a+\dfrac{k}{2t}a^{3}$ one gets $a(L^{2})=\sqrt{\dfrac{L^{2}}{c-kL^{2}}}$ for $c-kL^{2}>0$, where $c$ is a constant of integration. For these functions $a$, the Eq. (4.4) becomes
\begin{equation}
R_{srij}(x)y^{s}=-k(g_{js}g_{ri}-g_{is}g_{rj})y^{s},
\end{equation}
which says that $(M,g)$ is of constant curvature $-k$. Thus we have proved

\begin{theorem}
If the (pseudo)-Riemannian manifold $(M,g)$ is of constant curvature $k\in R$, for $a(L^{2})=\sqrt{\dfrac{L^{2}}{c+kL^{2}}}$ with $c$ a constant such that $c+kL^{2}>0$, the structures $(G_{a,a},F_{a})$ are locally conformal K\"{a}hler structures on $TM$.
\end{theorem}

The explicit form of these structures is as follows:
\begin{equation}
G_{a,a}(x,y)=g_{ij}(x)dx^{i}\otimes dx^{j}+\dfrac{1}{c+kL^{2}}
(g_{ij}(x))\delta y^{i}\otimes \delta y^{j}.
\end{equation}
\begin{equation}
F_{a}(\delta _{i})=-\sqrt{c+kL^{2}}\stackrel{.}{\partial _{i}},F_{a}(%
\stackrel{.}{\partial _{i}})=\dfrac{1}{\sqrt{c+kL^{2}}}\delta _{i},
\end{equation}

The 1-form $\theta$ is
\begin{equation}
\theta=\dfrac{kL}{c+kL^{2}}dL.
\end{equation}

\begin{corollary}
For $a(L^{2})=c_{0}\sqrt{L^{2}}$, with $c_{0}$ a strict positive constant,
the pairs $(G_{a,a},F_{a})$ are K\"{a}hler structures on $TM$ if and only if $(M,g)$ is flat.
\end{corollary}

\begin{proof}
If $(M,g)$ is flat, by the Theorem 4.1 for $a(L^{2})=c_{0}\sqrt{L^{2}},c_{0}=\dfrac{1}{\sqrt{c}}$, the pair $(G_{a,a},F_{a})$ is a locally conformal K\"{a}hler structure and by the Theorem 3.4 this is almost K\"{a}hler. Thus $(G_{a,a},F_{a})$ is a K\"{a}hler structure on $TM$. Conversely, if the pair $(G_{a,a},F_{a})$ with $a(L^{2})=c_{0}\sqrt{L^{2}}$ is a K\"{a}hler structure,
the Eq. (4.3) gives $R_{ij}^{k}=0$, equivalently $R_{srij}(x)=0$, that is $(M,g)$ is flat.
\end{proof}

Looking at (4.6) and (4.7) we see that the structures $(G_{a,a},F_{a})$ from
Corollary 4.1 are very close to $(G_{S},F_{S})$ which is obtained for $c=1$.
Thus the Corollary 4.1 covers a well-known result: $(G_{S},F_{S})$ is a
K\"{a}hlerian structure if and only if $(M,g)$ is flat.

Finally, we notice that for $c=0$ and $k\longrightarrow \dfrac{1}{k^{2}}$ in
(4.6) and (4.7) one obtains the locally conformal K\"{a}hler structure found by
R. Miron in [3].

\newpage

\runningauthor={M. ANASTASIEI}
\runningtitle={DEFORMATIONS OF FINSLER METRICS}
\noindent
\baselineskip 8pt
\noindent{\footnotesize{Finslerian geometries (Edmonton,AB, 1998),53-65}}
\hfill\break
{\footnotesize{Fund.Theories Phys., 109, Kluwer Acad. Publ. Dordrecht, 2000.}}
\vskip 2cm
\baselineskip 11.5pt plus .15pt
\centerline{\bf\Large DEFORMATIONS OF FINSLER METRICS}
\vskip .5cm
\centerline{\bf by Mihai ANASTASIEI and Hideo SHIMADA}
\vskip .5cm

\def\vsp{\vspace*{1,5mm}}
\def\trio#1#2#3#4#5{{#1}_{#2}{}^{#3}{}_{{#4}{#5}}}
\def\scri{\scriptstyle{|}}
\def\T{\mathbb{T}}
\def\bigi{\big|}
\def\G{{\Gamma}}
\def\g{{\gamma}}
\def\de{{\delta}}
\def\bld#1#2{{\buildrel{#1}\over{#2}}}
\def\bldo#1{\buildrel{\circ}\over{#1}{}\!}


\hyphenation{Ber-wald va-nish va-ni-shes pa-ral-lel res-pect}

\def\ijk{^i_{jk}}
\def\bld#1#2{\buildrel{#1}\over{#2}}
\def\derp#1#2{\displaystyle\frac{\partial#1}{\partial#2}}
\def\derpp#1{\displaystyle\frac{\partial}{\partial#1}}

\begin{abstract}
Let $F^n=(M,F(x,y))$ be a Finsler space and $g_{ij}(x,y)$ its Finsler metric. We consider a deformation of $g_{ij}(x,y)$ of
the form $$^*g_{ij}(x,y)=a(x,y)g_{ij}(x,y)+b(x,y) B_{i}(x,y)B_{j}(x,y),\leqno(1.1)$$ with two Finsler scalars $a>0$, $b\geq 0$ and $B_{i}(x,y)$ a Finsler co-vector. It follows that $^*g_{ij}$ is a generalized Lagrange metric in Miron'sense, briefly a GL--metric, see the monograph by R. Miron and M. Anastasiei \cite{8}. The metric $^*g_{ij}$ unifies the Antonelli metrics, the Miron--Tavakol
metrics, the Synge metrics (all treated in \cite{8}) as well as the Antonelli--Hrimiuc $\phi$-Lagrange metrics, \cite{2}, the Beil
metrics, \cite{4}, and the vertical part of the Cheeger--Gromoll metric, \cite{10}. We prove some general results on the geometry of
the GL- space $(M,^*g_{ij}(x,y)$. Next, the Levi-Civita connection and the curvature of a Riemannian metric on the
tangent manifold $TM$, induced by $g_{ij}$ and $^*g_{ij}$ are determined. These are used for the study of a Riemannian submersion involving the Cheeger--Gromoll metric.
\end{abstract}

\setcounter{section}{0}
\section{Deformations of Finsler metrics}

Let $F^n=(M,F)$ be a Finsler space with a smooth i.e. $C^\infty$ manifold $M$ and $F:TM\to R$, $(x,y)\mapsto F(x,y)$.
Here $x=(x^i)$ are coordinates on $M$ and $(x,y)=(x^i,y^i)$ are coordinates on the tangent manifold $TM$ projected on $M$ by
$\tau$. The indices $i,j,k,...$ will run from 1 to $n={\rm dim}\:M$ and the Einstein convention on summation is implied.
The geometrical objects on $TM$ whose local components change like on $M$ i.e. ignoring their dependence on $y$, will be called
Finsler objects as in \cite{7} or $d$--objects as in \cite{8}.

We set $\partial_i:=\derpp{x^i},$ $\dot\partial_i:=\derpp{y^i}$ and notice that the vertical subspace of $T_uTM$ i.e. $V_uTM={\rm Ker}\:(D\tau)_u$, $u\in TM$, where $D\tau$ means the differential of $\tau$, is spanned by $(\dot\partial_i)$. The $d$--objects can be expressed using $(\dot\partial_i)$.

The Finsler metric $g_{ij}(x,y)=\displaystyle\frac12\dot\partial_i\dot\partial_jF^2$ will be assumed positive definite. We have $F^2(x,y)=g_{ij}(x,y)y^iy^j$ and $F^2$ will be called the absolute energy of $F^n$. Assume that $F^n$ is endowed with a $d-$vector field or a Finsler vector field $B=B^i(x,y)\dot\partial_i$ and let $B_i(x,y) dx^i$ the Finsler 1-form with $B_i=g_{ik}B^k$. Set $B^2=B_iB^i$ and consider the following deformation of $g_{ij}(x,y)$: $$^*g_{ij}(x,y)=a(x,y)g_{ij}(x,y)+b(x,y) B_{i}(x,y)B_{j}(x,y),\leqno(1.1)$$ with two Finsler scalars $a>0$, $b\geq 0$. The metric $^*g_{ij}$ is no longer a Finsler metric but it is a positive definite generalized Lagrange metric in Miron'sense, briefly a GL--metric, see Ch. X in \cite{8}. It is easy to check that $^*g^{jk}=\displaystyle\frac1{a}g^{jk}- cB^jB^k$ is the inverse of $^*g_{ij}$ for $c=\displaystyle\frac{b}{a(a+bB^2)}$.

Various particular forms of $^*g_{ij}(x,y)$ were previously considered by some authors. The conformal case i.e. $b=0$, $a=exp(2\sigma(x,y))$ was studied and applied by R. Miron and R.K. Tavakol in {\it General Relativity}. The case $a=1$ and $B_i=y_i$
provides, for a convenient form of $b(x,y)$, a metric which generalizes the Synge metric from Relativistic Optics. This case was studied by R. Miron and T. Kawaguchi. For $b=0,$ $a=exp(2\sigma(x))$ and $g_{ij}(x,y)=g_{ij}(y)$ one gets the Antonelli metric which was used in Ecology. For the results on all these metrics we refer to the chapters XI and XII in \cite{8} and the references therein. The case $a=b=1$ and $B_i(x,y)= B_i(x)=\dfrac{\partial f}{\partial x^i},f: M \rightarrow\R$ was considered by C. Udri\c ste in \cite{11} for studying the completeness of a Finsler manifold. The Riemannian version of this case i.e. $g_{ij}(x,y)=g_{ij}(x)$
was intensively used by Th. Aubin in \cite{3}. The case $a=1$ and $g_{ij}(x,y)=g_{ij}(x)$ with various choices of $b$ and $B_i$
was introduced and studied by R. G. Beil for constructing a new unified field theory, \cite{5}.

One says that $^*g_{ij}$ is reducible to a Lagrange metric, shortly an $L$--metric if there exists a Lagrangian $L:TM\to\R$ such that $^*g_{ij}=\displaystyle\frac12\dot\partial_i\dot\partial_jL$. A necessary and sufficient condition for $^*g_{ij}$ be reducible to an $L$--metric is the symmetry in all indices of the Cartan tensor field $^*C_{ijk}=\displaystyle\frac12\dot\partial_k{}^*g_{ij}$ i.e.
$$\dot\partial_k{}^*g_{ij}=\dot\partial_i{}^*g_{kj}.\leqno(1.2)$$ Using (1.1) this condition becomes
$$\begin{array}{c}\dot a_kg_{ij}-\dot a_jg_{ik} + \dot b_kB_iB_j - \dot b_jB_iB_k +
b(\dot\partial_kB_i\cdot B_j-\dot\partial_iB_k\cdot B_j+\vspace*{1,4mm}\\
+B_i\cdot \dot\partial_kB_j-B_k\cdot\dot\partial_iB_j)=0,\
\dot a_k:=\dot\partial_k a,\dot b_k:=\dot\partial_k b.\end{array}\leqno(1.3)$$

Now we suppose that $a(x,y)=a(F^2)$ and $b(x,y)=b(F^2)$ assuming that the ranges of the real functions $a$ and $b$ from the right
hand are included in $Im(F^2)$. It results $\dot\partial_i a=2a'(F^2)y_i$ because of $\dot\partial_i F^2=2y_i$. Similarly,
$\dot\partial_i b=2b'(F^2)y_i$. We take $B_i=y_i$. For the GL-metric (1.1) subjected to the above conditions, (1.3) reduces to
$$(2a-b')(g_{ij}y_k-g_{ik}y_j)=0.\leqno(1.4)$$ Now if the equation $g_{ij}y_k-g_{ik}y_j=0$ is multiplied by $g^{ij}$ one gets $(n-1)y_k= 0$ which is a contradiction for $n\geq 1$. Thus we have

\medskip

{\bf Theorem 1.1.} {\it The GL-metric (1.1) with $B_i=y_i$, $a(x,y)=a(F^2),\\ b(x,y)=b(F^2)$ is an L-metric if and only if
$2a= b'$.}

\medskip

As always we may take $a= \phi'$, it comes out that the metric from Theorem 1.1 is essentially the $\phi$-Lagrange
metric of Antonelli-- Hrimiuc, \cite{2}, i.e. $$^*g_{ij}(x.y)= ag_{ij}(x,y)+ 2a'y_iy_j\leqno(1.5)$$
The Cheeger-Gromoll metric is a Riemannian metric on $TM$ of the form $$G_{CG}=g_{ij}dx^i\otimes dx^j+
\dd\frac{1}{1+F^2}(g_{ij}(x)+y_iy_j)\delta y^i\otimes \delta y^j,\leqno(1.6)$$ for $\delta y^i=dy^i+\gamma^i_{jk}y^jdx^k$,
where $\gamma^i_{jk}$ are the Christoffel symbols of $g_{ij}(x)$. This suggests considering the following GL-metric of type (1.1)
which generalizes the ``vertical part'' in (1.6): $$^*g_{ij}=\frac{1}{1+F^2}(g_{ij}(x)+y_iy_j),\leqno(1.7)$$
which we call a CGL--metric.

\medskip

{\bf Corolary 1.1.} {\it The $CGL$-metric $(1.7)$ is never reducible to a L-metric nor to a  Finsler metric.}

\medskip

\section{Metrical connection of the GL--space \hfill\break
$(M,{}^*g_{ij}(x,y))$}

The geometry of $^*g_{ij}(x,y)$ is naturally connected with the geometry of $F^n$. It is our purpose to express the geometrical
objects associated to $^*g_{ij}(x,y)$ using similar ones for $F^n$. If $\g^i_{jk}(x,y)$ are the generalized Christoffel symbols
for $g_{ij}(x,y)$ and we put $\g^i_{00}:=\g^i_{jk}y^jy^k$, then $\bldo N^i_j=\displaystyle\frac12\dot\partial_j\g^i_{00}$ are the local coefficients of the Cartan nonlinear connection. The Cartan connection for $F^n$ is $C\Gamma=(\bldo N^i_j,F^i_{jk},C^i_{jk}),$
where $$\begin{array}{l}F^i_{jk}=\displaystyle\frac12 g^{ih}(\de_jg_{hk}+\de_kg_{jh}-\de_hg_{jk}),\vsp\\
C^i_{jk}=\displaystyle\frac12 g^{ih}(\dot\partial_jg_{hk}+\dot\partial_kg_{jh}-\dot\partial_hg_{jk}),\end{array}\leqno(2.1)$$
for $\de_j=\partial_j-\bldo N^k_j\dot\partial_k$.

This connection is $h$--metrical, i.e. $g_{ij\bldo\scri k}=0$ and $v$--metrical, i.e. $g_{ij\bldo\bigi k}=0$. Here $\bldo\scri\scriptstyle k$ and $\bldo\bigi\scriptstyle k$ denote the $h$-- and $v$--covariant derivatives with respect to $C\G$. Moreover, two torsions of it vanish. We may consider a similar connection for $^*g_{ij}(x,y)$. Indeed, let
$^*C\G=(\bldo N^i_j,{}{}^*F^i_{jk},{}{}^*C^i_{jk})$ be the $d$--connection given by
$$\begin{array}{l}^*F^i_{jk}=\displaystyle\frac12{}^* g^{ih}(\de_j{}^*g_{hk}+\de_k{}^*g_{jh}-\de_h{}^*g_{jk}),\vsp\\
^*C^i_{jk}=\displaystyle\frac12{}^* g^{ih}(\dot\partial_j{}^*g_{hk}+\dot\partial_k{}^*g_{jh}-\dot\partial_h{}^*g_{jk}).\end{array}\leqno(2.2)$$
This $d$--connection is $h$--metrical i.e. $^*g_{ij\bld*\scri k}=0$ and $v$--metrical i.e. $^*g_{ij\bld*\bigi k}=0$ and the
torsions $^*T^i_{jk}:={}^*F^i_{jk}-{}^*F^i_{kj}=0$, $^*S^i_{jk}:={}^*C^i_{jk}-{}^*C^i_{kj}=0$. Moreover, when $\bldo N^i_j(x,y)$ is fixed, $^*C\G$ is the unique $d$--connection with these properties. It will be called the canonical metrical connection
of $^*g_{ij}(x,y)$. Using (1.1) in (2.2), after some calculation one gets

\medskip

{\bf Proposition 2.1.} {\it The metrical connection $^*C\G$ is given by $$^*F^i_{jk}=F^i_{jk}+\Phi^i_{jk},\ \ ^*C^i_{jk}=C^i_{jk}+\Lambda^i_{jk},\leqno(2.3)$$
$$\begin{array}{r}
\Phi^i_{jk}=\displaystyle\frac12{}^*g^{ih}[a_jg_{hk}+a_kg_{jk}-a_hg_{jk}
+\de_j(bB_kB_h)+\vsp\\
+\de_k(bB_jB_h)-\de_k(bB_jB_k)]-ac
B^iB^hF_{jhk}\end{array}\leqno(2.4)$$
$$\begin{array}{r}
\Lambda^i_{jk}=\displaystyle\frac12{}^*g^{ih}[\dot a_jg_{hk}+\dot ag_{jh}
-\dot a_hg_{jk}+\dot\partial_j(bB_kB_h)+\vsp\\
+\dot\partial_k(bB_jB_h)-\dot\partial_h(bB_jB_k)]-
acB^iB^hC_{ihk}\end{array}\leqno(2.5)$$
with the notations
$$\begin{array}{ll}
a_k=\de_ka,\ \dot a_k=\dot\partial_ka,
&F_{jhk}=\displaystyle\frac12(\de_jg_{hk}+\de_kg_{jh}-\de_hg_{jk}),\vsp\\
&C_{jhk}=\displaystyle\frac12(\dot\partial_jg_{kh}+\dot\partial_kg_{jh}-\dot\partial_hg_{jk}).
\end{array}\leqno(2.6)$$}

\medskip

{\bf Proposition 2.2.} {\it The torsions of $^*C\G$ are as follows:
$$\begin{array}{l}
^*T^i_{jk}=0,\ ^*R^i_{jk}=R^i_{jk}:=\de_k\bldo N^i_j-\de_j\bldo N^i_k,\
^*S^i_{jk}=0\vsp\\
^*P^i_{jk}=P^i_{jk}-\Phi^i_{jk}\mbox{ where }
P^i_{jk}=\dot\partial_k N^i_j-F^i_{jk}\mbox{ and }^*C^i_{jk}
\mbox{ from }(2.3).\end{array}\leqno(2.7)$$}

\medskip

{\bf Proposition 2.3.} {\it The curvatures of $^*C\G$ are as follows: $$^*\trio Sjikh=\trio Sjikh+\trio\Lambda jikh+
(C^s_{jk}\Lambda^i_{sh}+\Lambda^s_{jc}C^i_{sh}-k/h),\leqno(2.8)$$ $$\trio\Lambda jikh=\dot\partial_h\Lambda^i_{jk}+
\Lambda^s_{jk}\Lambda^i_{sh}-k/h,\leqno(2.8)'$$ where $-k/h$ means the subtraction of the preceding terms with
$k$ replaced by $h$.}

$$^*\trio Rjikh=\trio Rjikh+\trio\Phi jikh+
(F^s_{jk}\Phi^i_{sh}+\Phi^s_{jk}F^i_{sh}-k/h)+
\Lambda^i_{js}R^s_{kh},\leqno(2.9)$$
$$\trio\Phi jikh=\de_h\Phi^i_{jk}+\Phi^s_{jk}\Phi^i_{sh}-k/h,\leqno(2.9)'$$
$$^*\trio Pjikh=\trio Pjikh+\Phi^i_{jk\bldo\bigi
h}-\Lambda^i_{jh\bldo\scri
k}+\Lambda^i_{js}P^s_{kh}+C^s_{kh}\Phi^i_{sj}+
\Phi^s_{jk}\Lambda^i_{sh}-\Phi^i_{sk}\Lambda^s_{jh}.\leqno(2.10)$$

\section{On a Riemannian metric on $TM$}

Let $TM$ be the tangent manifold to $M$ endowed with the fundamental Finsler function $F$ and the Finsler metric $g_{ij}(x,y)$. Consider the Cartan nonlinear connection $(\bldo N^a_j(x,y))$ and then $(\de_i=\partial_i-\bldo N^a_i\dot\partial_a,\dot\partial_a)$ is a local frame on $TM$ adapted to the decomposition of $T_uTM$ into a direct sum of vertical and horizontal subspaces. From now on we shall use two types of indices: $a,b,c,...$ will indicate  vertical components and $i,j,k,...$ will
indicate horizontal ones. All have the same range $\{1,2,...,n\}$. Let be $h_{ab}(x,y)=\de^i_a\de^j_b{}^*g_{ij}(x,y)$, where $\de^i_a$ is the Kronecker symbol, and $$G(x,y)=g_{ij}(x,y) dx^i\otimes dx^j+h_{ab}(x,y)\de y^a\otimes\de
y^b,\leqno(3.1)$$ where $\de y^a=dy^a+N^a_k(x,y) dx^k$.

Then $(TM,G(x,y))$ is an oriented Riemannian manifold. The horizontal and vertical distributions are mutually orthogonal with
respect to $G$. It is our purpose to study the Riemannian metric $G$. First, we compute the coefficients of the Levi--Civita
connection $D$ of $G$ in the frame $(\de_i,\dot\partial_a)$. We set
$$\begin{array}{ll}
D_{\de_k}\de_j=F^i_{jk}\de_i+A^a_{jk}\dot\partial_a,
&D_{\dot\partial_b}\de_j=\widetilde
C^i_{jb}\de_i+E^a_{jb}\dot\partial_a, \vsp\\
D_{\de_k}\dot\partial_b=L^a_{k}\dot\partial_a+D^i_{bk}\de_i,
&D_{\dot\partial_b}\dot\partial_c=C^a_{cb}\dot\partial_a+B^i_{cb}\de_i\end{array}\leqno(3.2)$$
Let $\T$ be the torsion of $D$ i.e. $\T(X,Y)=D_XY-D_YX-[X,Y]$ for
$X,Y$ vector fields on $TM$. The condition $D$ is torsion--free is
equivalent to
$$\T(\de_i,\de_j)=\T(\de_i,\dot\partial_a)=\T(\dot\partial_a,\dot\partial_b)=0.\leqno(3.3)$$
Using the following equations
$$[\de_i,\de_j]=R^a_{ij}\dot\partial_a,\
[\de_j,\dot\partial_b]=(\dot\partial_bN^a_j)\dot\partial_a,\
[\dot\partial_a,\dot\partial_b]=0\leqno(3.4)$$ where
$R^a_{ij}=\de_jN^a_i-\de_iN^a_j,$ one finds that (3.3) is
equivalent to
$$\begin{array}{ll}
F^k_{ij}=F^k_{ji},&A^a_{ij}-A^a_{ji}=-R^a_{ij}\vsp\\
D^k_{ai}=\widetilde C^k_{ia},&L^b_{ai}=\dot\partial_aN^b_i+E^b_{ia}\vsp\\
C^a_{bc}=C^a_{cb},&B^i_{bc}=B^i_{cb}.\end{array}\leqno(3.5)$$
The condition that $D$ is metrical, that is, $XG(X,Y)=G(D_XY,Z)+G(Y,D_XZ)$, written in the frame
$(\de_i,\dot\partial_a)$ gives
$$\begin{array}{ll}
F^h_{ji}g_{hk}+F^h_{ki}g_{hj}=\de_ig_{jk},
&\widetilde C^i_{ja}g_{ik}+\widetilde C^i_{ka}g_{ij}=
\dot\partial_ag_{jk},\vsp\\
A^c_{ji}h_{ca}+D^k_{ai}g_{kj}=0,
&E^c_{ja}h_{cb}+B^k_{ba}g_{kj}=0,\vsp\\
L^c_{ai}h_{cb}+L^c_{bi}h_{ca}=\de_ih_{ab},
&C^e_{ba}h_{ec}+C^e_{ca}h_{eb}=\dot\partial_ah_{bc}.\end{array}\leqno(3.6)$$
The systems (3.5) and (3.6) have the unique solution
$$\begin{array}{l}
F^k_{ij}=\displaystyle\frac12 g^{kh}(\de_ig_{hj}+\de_jg_{hi}-\de_hg_{ij}),\
A^a_{jk}=\displaystyle\frac12(-R^a_{jk}-h^{ab}\dot\partial_bg_{jk}),\vsp\\
\widetilde C^i_{jb}=\displaystyle\frac12g^{ih}(\dot\partial_bg_{jh}+
h_{bc}R^c_{hj)}=D^i_{bj},\vsp\\
E^a_{ib}=\displaystyle\frac12h^{ac}h_{bc\scriptstyle{\|}i},\
L^a_{bi}=\dot\partial_bN^a_i+\displaystyle\frac12
h^{ac}h_{bc\scriptstyle{\|}i},\vsp\\
B^k_{ab}=-\displaystyle\frac12g^{kj}h_{ab\scriptstyle{\|}j},\
C^a_{bc}=\displaystyle\frac12h^{ad}(\dot\partial_bh_{dc}+
\dot\partial_ch_{bd}-\dot\partial_dh_{bc}).
\end{array}\leqno(3.7)$$
Here $h_{bc\scriptstyle{\|}i}$ denotes the $h$--covariant derivative of $h_{bc}$ with respect to the Berwald connection
$B\G=(\bldo N^a_i,\dot\partial_b N^a_i,0)$. Now we shall compute the components of the curvature of $D$ in the same frame. To this aim we shall consider an intermediate linear connection $\nabla$ on
$TM$:
$$\begin{array}{ll}
\nabla_{\de_j}\de_k=F^i_{jk}\de_i,
&\nabla_{\dot\partial_b}\de_j=D^i_{bj}\de_i\vsp\\
\nabla_{\de_k}\dot\partial_b=L^a_{bk}\dot\partial_a,
&\nabla_{\dot\partial_b}\dot\partial_c=C^a_{cb}\dot\partial_a.\end{array}\leqno(3.8)$$
This connection is metrical with respect to $G$ i.e. $\nabla_XG=0$, it preserves the horizontal and vertical
distributions and it has three non-vanishing torsions: $R^a_{jk},D^i_{bj},$ $P^a_{jb}=\displaystyle\frac12
h^{ac}h_{bc\scriptstyle{\|}j}$.

The curvature of $\nabla$ has six components in the form (see p. 48 of \cite{8}):
$$\begin{array}{lcl}\widehat{\trio Rhijk}=\de_kF^i_{hj}+F^m_{hj}F^i_{mk}-j/k+
D^i_{ah}R^a_{jk},\vsp\\
\widetilde{\trio
Rbajk}=\de_kL^a_{bj}+L^c_{bj}L^a_{ck}-j/k+C^a_{bc}R^c_{jk},\vsp\\
\widetilde{\trio Pjika}=\dot\partial_a F^i_{jk}-D^i_{aj\scri
k}+D^i_{bj}P^b_{ka},\vsp\\
\trio Pbakc=\dot\partial_cL^a_{bk}-C^a_{bc\scri
k}+C^a_{bd}P^d_{kc},\vsp\\
\widetilde{\trio Sjibc}=\dot\partial_c D^i_{bj}+D^h_{bj}D^i_{ch}-b/c,\vsp\\
\trio Sbacd=\dot\partial_dC^a_{bc}+C^e_{bc}C^a_{ed}-c/d.
\end{array}\leqno(3.9)$$
Here and in the following $\scriptstyle{|k}$ and $\bigi\scriptstyle{a}$ will denote $h$-- and $v$--covariant
derivatives with respect to $\nabla$.

{\it Remark 3.1.} $\trio Sbacd$ is nothing but $^*\trio Sjikh$. And the other tensors in (3.9) can be expressed with $\trio Rjikh,\trio Pjikh,\trio Sjikh$ or with their $*$--counterparts. For instance, $\widehat R_h{}^i{}_{jk}=\trio Rhijk+\displaystyle\frac12 g^{is}h_{ac}R^c_{sh}R^a_{jk}$.

Let $K$ be the curvature tensor field of the Levi--Civita connection $D$. We shall denote its components by the same letter
$K$ indexed with two types of indices with the understanding that different indices means different components. There will be twelve
components of $K$. After calculation one finds
$$\begin{array}{l}
K(\dot\partial_b,\dot\partial_c)\dot\partial_d:=\trio Kdacb\dot\partial_a+\trio
Kdicb\de_i,\vsp\\
\trio Kdacb=\trio Sdacb+B^i_{cd}E^a_{ib}-B^i_{db}E^a_{ic},\
\trio Kdicb=B^i_{cd\bigi b}-B^i_{bd\bigi c},\vsp\\
K_{abdc}=S_{abdc}+\displaystyle\frac12(B^i_{ad}h_{bc\scriptstyle{\|}i}-
B^i_{ac}h_{bd\scriptstyle{\|}i},
\end{array}\leqno(3.10)$$
$$\begin{array}{l}
K(\dot\partial_b,\dot\partial_c)\de_j=\trio Kjacb\dot\partial_a+\trio
Kjicb\de_i,\vsp\\
\trio Kjicb=\widetilde{\trio Sjicb}+E^d_{jc}B^i_{db}-E^d_{jb}B^i_{dc},\
\trio Kjacb=E^a_{jc\bigi b}-E^a_{jb\bigi c},
\end{array}\leqno(3.11)$$
$$\begin{array}{l}K(\dot\partial_b,\de_j)\dot\partial_c:=\trio
Kcajb\dot\partial_a+\trio Kcijb\de_i,\vsp\\
\trio Kcajb=\trio Pcajb-B^k_{cb}A^a_{kj}+D^k_{cj}E^a_{kb},\vsp\\
\trio Kcijb=D^i_{cj\bigi b}-B^i_{bc\scri
j}-P^d_{jb}B^i_{dc}+D^k_{bj}D^i_{ck},
\end{array}\leqno(3.12)$$
$$\begin{array}{l}
K(\dot\partial_b,\de_k)\de_j:=\trio Kjakb\dot\partial_a+\trio
Kjikb\de_i,\vsp\\
\trio Kjakb=A^a_{jk\bigi b}-E^a_{jb\scri
k}+D^h_{bk}A^a_{jh}+P^c_{kb}E^a_{jc},\vsp\\
\trio Kjikb=\widetilde{\trio Pjikb}+A^c_{jk}B^i_{cb}-E^c_{jb}D^i_{ck},\vsp\\
K_{jakb}=A_{ajk\bigi b}-E_{ajb\scri k}+A_{ajh}D^h_{bk}+E_{ajc}P^c_{kb},
\end{array}\leqno(3.13)$$
$$\begin{array}{l}
K(\de_j,\de_k)\dot\partial_b:=\trio Kbakj\dot\partial_a+\trio
Kbikj\de_i,\vsp\\
\trio Kbakj=\widetilde{\trio Rbakj}+D^h_{bk}A^a_{hj}-D^h_{bj}A^a_{hk},\vsp\\
\trio Kbikj=D^i_{bk\scri j}-D^i_{bj\scri k}-R^c_{jk}B^i_{bc},
\end{array}\leqno(3.14)$$
$$\begin{array}{l}
K(\de_j,\de_k)\de_h:=\trio Khakj\dot\partial_a+\trio Khikj\de_i,\vsp\\
\trio Khikj=\widehat{\trio
Rhikj}+A^b_{hk}D^i_{bj}-A^b_{hj}D^i_{bk},\vsp\\ \trio
Khakj=A^a_{hk\scri j}-A^a_{hj\scri k}+R^c_{kj}E^a_{hc},\vsp\\
K_{hikj}=R_{hikj}+D_{ibj}A^b_{hk}-D_{ibk}A^b_{hj}.
\end{array}\leqno(3.15)$$
Now easily follows

{\bf Proposition 3.1.} {\it The sectional curvatures of $D$ are as follows:
$$\begin{array}{l}
K_{ab}=[S_{abab}+\displaystyle\frac12(B^i_{aa}h_{bb\scriptstyle{\|i}}-B^i_{ab}
h_{ab\scriptstyle{\|i}})]/(h_{aa}h_{bb}-h^2_{ab}),\vsp\\
K_{ja}=(A_{ajj\bigi a}-E_{aja\scri
k}+A_{ajh}D^h_{aj}+E_{ajc}P^c_{ja})/g_{jj} g_{aa}\vsp\\
K_{ji}=(R_{jiji}+D_{ibi}A^b_{jj}-
D_{ibj}A^b_{ji})/(g_{ii}g_{jj}-g^2_{ij}).
\end{array}\leqno(3.16)$$}

In the following we assume that $F^n$ reduces to a Riemannian space i.e. $g_{ij}(x,y)=g_{ij}(x).$ The Cartan nonlinear connection
reduces to $\bldo N^i_j(x,y)=\g^i_{jk}(x)y^k,$ where $(\g^i_{jk}(x))$ are the Christoffel symbols of the metric $g=(g_{ij}(x))$. We consider the corresponding Riemannian metric $G$ given by (3.1) and we have

{\bf Proposition 3.2.} {\it The mapping $\tau:(TM,G)\to(M,g)$ is a Riemannian submersion.}

Indeed, $\tau$ is of maximal rank $n$ and its differential $D\tau$ preserves the lengths of horizontal vectors as it follows from $G(\de_i,\de_j)=g_{ij}(x)$.

Let $h$ and $v$ denote the projections of $T_uTM$ onto the subspaces of horizontal and vertical vectors, respectively.
Following B. O'Neil, \cite{9}, the fundamental tensor fields of the Riemannian submersion $\tau$ are as follows: \bigskip

\noindent $\begin{array}{lcl}
(3.17)&\hspace*{15mm}&S(X,Y)=hD_{vX}Y+vD_{vX}hY,\vsp\\
(3.18)&&N(X,Y)=vD_{hX}hY+hD_{hX}vY,\ \ X,Y\in{\cal X}(TM).\end{array}$\bigskip

In the frame $(\de_i,\dot\partial_a)$ we have
$$S(\de_i,\de_j)=0,S(\de_i,\dot\partial_a)=0,S(\dot\partial_a,\de_i)=E^j_{ia}\de_j,
S(\dot\partial_a,\dot\partial_b)=B^i_{ab}\de_i.\leqno(3.19), $$
$$N(\de_i,\de_j)=\displaystyle\frac12 R^a_{ij}\dot\partial_a,
N(\de_i,\dot\partial_a)=D^i_{ai}\de_j,
N(\dot\partial_a,\de_i)=0,N(\dot\partial_a,\dot\partial_b)=0\leqno(3.20).$$
By (3.19) and (3.7) it follows

{\bf Proposition 3.3.} {\it The Riemannian submersion $\tau:(TM,G)\to(M,g)$ is
totally geodesics, i.e. $S=0$ if and only if
$$^*g_{ij\scriptstyle{\|}k}=0,\leqno(3.21)$$
where $\scriptstyle{\|k}$ denotes the $h$--covariant derivative
with respect to the Berwald connection
$(\g^i_{jk}(x)y^k,\g^i_{jk}(x),0)$.}

{\bf Proposition 3.4.} {\it The tensor field $N$ vanishes if and only if the
Riemannian metric $g$ is flat.}

\section{Deformations of Riemannian metrics}

The geometrical objects associated to $^*g_{ij}(x,y)$ are generally complicated. Some simplifications appear for particular choices
of $a$, $b$ and $B_i$. We studied in a previous paper, \cite{1}, the case $a=1$ and a concurrent d-vector field $B^i(x,y)$ while M. Kitayama studied the case $a=1$ and a parallel d-vector field $B^i(x,y)$, \cite{6}.
Here we selected for a detailed analysis the following deformation of a Riemannian metric $g=(g_{ij}(x))$:
$$^*g_{ij}(x,y)=a(F^2)g_{ij}(x)+b(F^2)y_iy_j,\leqno(4.1)$$ where $F^2(x,y)=g_{ij}(x)y^iy^j,y_i=g_{ij}(x)y^j$.

Accordingly, we consider the Riemannian submersion $\tau:(TM,G)\to(M,g)$, where
$$G(x,y)=g_{ij}(x,y) dx^i\otimes dx^j+(a(F^2)g_{ij}(x)+b(F^2)y_iy_j)\de y^a\otimes\de
y^b,\leqno(4.2)$$

The $GL$--metric (4.1) contains as a particular case the $\phi$--Lagrange metric associated to a
Riemannian space while $G$ generalizes the Cheeger--Gromoll metric studied by Sekizawa \cite{10}.
The Cartan connection for $(M,g_{ij}(x))$ reduces to $C\G=(\g^i_{jk}(x)y^j,\g^i_{jk},0)$. The $v$--covariant derivative
$\bldo\bigi\scriptstyle{k}$ coincides with the partial derivative with respect to $(y^k)$. The $h$--covariant derivative
$\bldo\scri\scriptstyle{k}$ reduces to the usual covariant derivative for the objects which do not depend on $(y^i)$ and coincides with $\scriptstyle{\|k}$ for the others.

We notice for the later use the following formulae
$$\de_kF^2=0,\ y^i_{\bldo\scri k}=0,\ y_{i\bldo\scri k}=0,\
y^i_{\bldo\bigi k}=\de^i_k,\
y_{i\bldo\bigi k}=g_{ik}(x)\leqno(4.3)$$
$$\begin{array}{ll}
\de_ka=0,&\de_kb=0\vsp\\
\dot\partial_ka=2a'y_k,&\dot\partial_kb=2b'y_k.\end{array}\leqno(4.4)$$
By a direct calculation one proves

{\bf Proposition 4.1.} {\it The $d$--connection $C\G$ of the $GL$--metric $(4.1)$
is given by
$$\begin{array}{ll}
^*F^i_{jk}&=\gamma^i_{jk}(x)  i.e.\ \ \ \ \ \Phi^i_{jk}=0\vsp\\
^*C^i_{jk}(x,y)&=\Lambda^i_{jk}(x,y)=
\displaystyle\frac{a'}a(\de^i_ky_j+\de^i_jy_k)+\vsp\\
&+\displaystyle\frac{b-a'}{a+bF^2}y^ig_{jk}+
\frac{ab'-2a'b}{a(a+bF^2)}y^iy_jy_k.\end{array}\leqno(4.5)$$}

From (4.3) and (4.4) it results
$$^*g_{ij\bldo\scri k}=0,\ ^*g_{ij\bldo\bigi k}=2a'g_{ij} y_k+
b(g_{ik}y_j+g_{jk}y_i)+2b'y_iy_jy_k.\leqno(4.6)$$
Thus $^*g_{ij}$ is $h$--metrical and not $v$--metrical with respect
to $C\G$. The torsions of $^*C\G$ of the $GL$--metric $(4.1)$
are vanishing excepting $^*R^i_{jk}=\g^i_{hjk}(x)y^h$ and $^*C^i_{jk}$
from $(4.5)$. As for its curvatures we find
$$^*\trio Rjikh=\trio rjikh(x)+\Lambda^i_{js}  R^s_{kh},\leqno(4.7)$$
$$^*\trio Pjikh=0\mbox{ because of $\Lambda^i_{jh\bldo\scri
k}=0$},\leqno(4.8)$$ $$^*\trio Sjikh=\Lambda^i_{jkh}\mbox{ from
}(2.8)',\leqno(4.9)$$ where $r^ i_{j kh}$ is the curvature tensor of
$(g_{ij}(x))$.

Using $y_sR^s_{ kh}=y_s\trio rspkhy^p=r_{pikh}y^py^i=0,$ one gets
$$^*\trio Rjikh=\trio rjikh(x)+\displaystyle\frac{a'}ay_jR^i_{ kh}+
\frac{b-a'}{a+bF^2}y^iR_{jkh},\leqno(4.7)'$$
$$^*\trio R0ikh=\left(1+\displaystyle\frac{a'F^2}a\right)R^i_{ kh},
\leqno(4.7)''$$ where ``$0$'' denotes the contraction by $(y^j)$.

Now we consider the Riemannian metric $G$ given by (4.2). The
Levi--Civita connection of it has the local coefficients
$$\begin{array}{l}
F^k_{ij}=\g^k_{ij}(x),\
A^a_{jk}=-\displaystyle\frac12\trio r0ajk,\vsp\\
D^i_{bj}=\displaystyle\frac a2\trio rjib0=\widetilde C^i_{jb},\vsp\\
E^a_{ib}=0=B^k_{ab},\
L^a_{bi}=\g^a_{bi}(x),\
C^a_{bc}=\Lambda^a_{bc}.\end{array}\leqno(4.10)$$
The curvature of $\nabla$ from (3.9) reduces to
$$\begin{array}{lcl}
\trio{\widehat R}hijk=\trio rhijk(x)+\displaystyle\frac a2\trio
rhia0\cdot\trio r0ajk,\vsp\\
\widetilde{\trio Rbajk}={}^*\trio Rbajk,\vsp\\
\widetilde{\trio Pjika}=-\displaystyle\frac a2\trio rjia{0;k},\vsp\\
\trio Pbakc=0\mbox{ because of }\Lambda^a_{bc\scri k}=0,\vsp\\
\widetilde{\trio Sjibc}=a\trio rjibc+\left(a'y_c\trio
rjib0+\displaystyle\frac{a^2}4\trio
rjhb0r^i_{hc0}-b/c\right),\vsp\\
\trio Sbacd=\Lambda^ a_{b cd}.\end{array}\leqno(4.11)$$
The curvature of the Levi--Civita connection $D$ are given by
$$\begin{array}{l}
\trio Kdabc=\trio\Lambda dabc,\ \trio Kdicb=0,\vsp\\
\trio Kjibc=\trio{\widetilde S}jibc,\ \trio Kjabc=0,\vsp\\
\trio Kcajb=0,\ \trio Kcijb=\displaystyle\frac a2\trio rjicb-
\displaystyle\frac{a'}2y_b\trio rjic0-\displaystyle\frac{a'}2y_c\trio
rjib0+ \displaystyle\frac{a^2}4\trio rsib0\trio
rjsc0,\vsp\\
\trio Kjakb=-\displaystyle\frac12\trio rbajk+\displaystyle\frac a4\trio r0aik-
\displaystyle\frac{a'}{2a}y_b\trio r0ajk-
\frac{b-a'}{2(a+bF^2)}y^ar_{0bjk},\vsp\\
\trio Kjikb=-\displaystyle\frac a2\trio rjib{0;k},\vsp\\

\trio Kbakj=\trio rbakj+\displaystyle\frac{a'}ay_b\trio r0akj+
\displaystyle\frac{b-a'}{a+bF^2}y^ar_{0bkj}-\displaystyle\frac a4\trio
rkhb0\trio r0ahk+\vsp\\
\hfill+\displaystyle\frac a4\trio rjhb0\trio
r0ajk,\\
\trio Kbikj=\displaystyle\frac a2(\trio rkib{0;j}-\trio
rjib{0;k}),\vsp\\
\trio Khikj=\trio rhikj+\displaystyle\frac a2\trio rhia0\trio
r0akj- \displaystyle\frac a4\trio r0ahj\trio rkia0,\vsp\\
\trio Khakj=\displaystyle\frac12(\trio r0ah{j;k}-\trio r0rah{k;j}).
\end{array}\leqno(4.12)$$
An inspection of (4.11) and (4.12) gives

{\bf Theorem 4.1.} {\it If $(M,g)$ is flat, then $(TM,G)$ is flat if and only
if $\trio\Lambda jikh=0$.}

This theorem shows that $G$ is less ``rigid'' than the Sasaki metric
of $(g_{ij}(x))$ which is locally flat if and only if $(g_{ij}(x))$ is
locally flat.

Now if we fix $x=x_0,$ then $^*g_{ij}(x_0,y)$ is a Riemannian metric
in the fibre $T_{x_0}M$ and $\trio\Lambda jikh$ is just its curvature tensor field. Thus we may reformulate Theorem
\ref{t4.1} in the form \bigskip

{\bf Theorem 4.1$'$.} {\it If $(M,g)$ is flat, then
$(TM,G)$ is flat if and only if $(T_{x_0}(M)$, $^*g_{ij}(x_0,y))$ is a flat
Riemannian manifold for every $x_0\in M.$}\bigskip

For the conformal case i.e. $b=0$ one finds
$$\Lambda^i_{jk}=\frac{a'}a(\de^i_ky_j-\de^i_jy_k-y^ig_{jk})\leqno(4.13)$$
$$\begin{array}{ll}
\trio\Lambda jikh&=
\left[2\left(\displaystyle\frac{a'}a\right)'-
\displaystyle\frac{a'}a^2\right]
\left(\de^i_ky_jy_h+y^iy_kg_{jh}-h/k\right)+\vsp\\
&+\displaystyle\frac{a'}a^2F^2(\de^i_kg_{jh}-\de^i_hg_{jk}).
\end{array}\leqno(4.14)$$
It follows

{\bf Proposition 4.2.} {\it $\trio\Lambda jikh=0\iff a={\rm constant.}$}

From Theorem \ref{p3.3} and (4.6) one deduces

{\bf Proposition 4.3.} {\it The Riemannian submersion $\tau:(TM,G)\to (M,g)$ with
$G$ given by $(4.2)$ is totally geodesics.}

The other consequences of the previous formulae will be presented
elsewhere.

\noindent {\footnotesize{University ``Al. I. Cuza'' Ia\c si\\
Faculty of Mathematics\\ 6600, Ia\c si, Romania}}

\medskip

\noindent {\footnotesize{Hokkaido Tokai University\\
Minami-ku, Minami-sawa\\ 5-1, Sapporo 005, Japan}}


\newpage

\def\rr{\mathbb{R}}
\def\ov#1{\overline{#1}}

\runningauthor={M. ANASTASIEI}
\runningtitle={A FRAMED $f-$STRUCTURE}
\noindent
\baselineskip 8pt
\noindent{\footnotesize{An.Univ.Bucure\c sti Mat. Inform.\hfill\break 49 (2000), no. 2, 3-9}}
\vskip 2cm
\baselineskip 11.5pt plus .15pt
\centerline{\bf\Large A FRAMED $f-$STRUCTURE}
\vskip .2cm
\centerline{\bf\Large ON TANGENT MANIFOLD}
\vskip .2cm
\centerline{\bf\Large OF A FINSLER SPACE}
\vskip .5cm
\centerline{\bf by Mihai ANASTASIEI}
\vskip 1cm

\begin{abstract}
It is shown that the slit tangent manifold $\buildrel{{\circ}}\over{TM}$ of a Finsler
space $F^n=(M,L)$ carries a natural framed $f-$structure of
corank 2. When this is restricted to the indicatrix bundle $TM$
defined by $L=1$, one gets a {\it contact Riemannian structure on
$TM$} which is Sasakian iff $F^n$ is of constant curvature 1.

{\bf Mathematics Subject Classifications 2000:} 53C60.

{\bf Key words and phrases:} tangent manifold, framed
$f-$structure, indicatrix bundle.
\end{abstract}

\setcounter{section}{0}
\section{Introduction}
Let $M$ be a smooth i.e. $C^\infty $ manifold of dimension $n$ and
$\tau:TM\to M$ its tangent bundle. If $(x^i)$, $i,j,k...=1,...,n$
are local coordinates on $M$, the induced local coordinates on
$TM$ will be denoted by $(x,y)\equiv(x^i=x^i\circ\tau,y^i)$, where
$(y^i)$ are the components of a vector from $T_pM$, $p(x^i)$, in
the natural basis $\left(\partial_i:=\displaystyle\frac\partial{\partial x^i}\right)$.

Let $F^n=(M,L)$ be a Finsler space. The function
$L:\buildrel{{\circ}}\over{TM}:=TM\setminus\{(x,0)\}\to{\mathbb{R}}_+$ is smooth, positively
homogeneous of degree 1 with respect to $(y^i)$ and the matrix
with the entries $g_{ij}(x,y)=\dfrac12\dfrac{\pp^2L^2}{\pp y^i\partial y^j}$ is of
rank $n$. From the homogeneity of $L$ it follows $L^2(x,y)=g_{ij}(x,y)
y^iy^j=y^iy_i$ for $y_i=g_{ij} y^j$.

If $\gamma^i_{jk}(x,y)$ are the ``generalized'' Christoffel symbols constructed using $g_{ij}(x,y)$,
and $\gamma^i_{00}(x,y)=\g^i_{jk}(x,y) y^iy^j$, then the functions
$N^i_j(x,y)=\displaystyle\frac12\dot\partial_j(\gamma^i_{00})$, $\dot\partial_j:=\displaystyle\frac{\partial}{\pp y^j}$
are the local coefficients of the nonlinear Cartan connection of $F^n$.
For details see Ch. VIII in [7].

Using them, a new local basis $(\delta_i,\dot\partial_i)$, where
$\delta_k=\partial_i-N_i^k\dot\partial_k$, on $\buildrel{{\circ}}\over{TM}$ is introduced.
The dual of this basis is $(dx^i,\delta y^i=dy^i+N^i_kdx^k)$. If the
quadratic form of matrix $(g_{ij}(x,y))$ is positive defined, then
$G_S=g_{ij}(x,y) dx^i\otimes dx^j+g_{ij}(x,y)\delta y^i\otimes\delta y^j$ is a
Riemannian metric on the tangent manifold $\buildrel{{\circ}}\over{TM}$, called the
Sasaki--Matsumoto lift of $(g_{ij})$ to $\buildrel{{\circ}}\over{TM}$. The linear operator
$F$ given in the local basis $(\delta_i,\dot\partial_i)$ as follows:
$F(\delta_i)=-\dot\partial_i$, $F(\dot\partial_i)=\delta_i$, defines an almost
complex structure on $\buildrel{{\circ}}\over{TM}$ and the pair $(F,G_S)$ is an almost
K\"ahlerian structure on $\buildrel{{\circ}}\over{TM}$.

On $\buildrel{{\circ}}\over{TM}$ there exist two remarkable vector fields: $C=y^i\dot\partial_i$, called
the Liouville vector field and $S=y^i\dot\partial_i$, which is the
geodesic spray of $F^n$.

A framed $f-$structure is a natural generalization of an almost
contact structure. It was introduced by S.I. Goldberg and K. Yano
[3]. We recall its definition following [5, p.47].

Let $N$ be a $(2n+s)-$dimensional manifold endowed with an
endomorphism $f$ of rank $2n$, of the tangent bundle, satisfying
$f^3+f=0$. If there exist on $N$ the vector fields $(\xi_\alpha)$ and the
1-forms $y^\alpha$ $(\alpha=1,2,...,s)$ such that
$\eta^\alpha(\xi_\beta)=\delta^\alpha_\beta,$ $f(\xi_\alpha)=0$,
$\eta^\alpha\circ f=0$, $f^2=-I+\displaystyle\sum_\alpha\eta^\alpha\otimes\xi_\alpha$, where $I$ is
the identity automorphism of the tangent bundle, then $N$ is said
to be a framed $f-$manifold.

\section{A framed $f-$structure on $\buildrel{{\circ}}\over{TM}$}
\setcounter{equation}{0}

Denote $\xi_1=y^i\delta_i=S$ and $\xi_2=y^i\dot\partial_j=C$. From the
definition of $F$ it follows

\ms

{\bf Lemma 2.1.} $F(\xi_1)=-\xi_2,\ F(\xi_2)=\xi_1$.

\ms

We introduce the 1-forms $\eta^1=\displaystyle\frac{y_i}{L^2}dx^i$ and
$\eta^2=\displaystyle\frac{y_i}{L^2}=\delta y^i$. These are globally defined
on $\buildrel{{\circ}}\over{TM}$. By a direct calculation one gets

\ms

{\bf Lemma 2.2.} $\eta^1\circ F=\eta^2,\ \eta^2\circ F=-\eta^1$.

\ms

Let be $G=\displaystyle\frac1{L^2}G_S$. One easily verifies

\ms

{\bf Lemma 2.3.} {\it $\eta^1(X)=G(X,\xi_1),\
\eta^2(X)=G(X,\xi_2)$, for every $X\in{\cal X}(\buildrel{{\circ}}\over{TM})$, the
module of vector fields on $\buildrel{{\circ}}\over{TM}$.}

\ms

We define a tensor field of type $(1,1)$ on $\buildrel{{\circ}}\over{TM}$ by
\begin{equation}
f(X)=F(X)+\eta^1(X)\xi_2-\eta^2(X)\xi_1,\ X\in{\cal X}(\buildrel{{\circ}}\over{TM}).
\end{equation}

\ms

{\bf Theorem 2.1.} {\it The ensemble $(f,(\xi_a),(\eta^b))$
$a,b,...=1,2$ provides a framed $f-$structure on $\buildrel{{\circ}}\over{TM}$,
that is the followings hold:
\begin{itemize}
\item[{\rm(i)}] $f$ is of ${\rm rank}\,2n-2$ and $f^3+f=0$,
\item[{\rm(ii)}] $\eta^a(\xi_b)=\delta^a_b,\ f(\xi_a)=0,\ \eta^a\circ
f=0$,
\item[{\rm(iii)}] $f^2(X)=-X+\eta^1(X)\xi_1+\eta^2(X)\xi_2,\
X\in{\cal X}(\buildrel{{\circ}}\over{TM})$.
\end{itemize}}

\ms

\noindent{\it Proof.} Using (2.1) and the Lemmas 2.1 and 2.2 one easily
checks (ii) and (iii). Applying $f$ on the equality (ii) one
obtains the second part in (i). From the second equations in (ii),
we see that ${\rm span}\{\xi_1,\xi_2\}\subseteq{\rm Ker}\,f$. If
$X=X^k\delta_k+\dot X^k\dot\partial_k$ belongs to ${\rm Ker}\,f$ and it is
not in ${\rm span}\{\xi_1,\xi_2\}$, we hve $X^iy_i=0$ and $\dot
X^iy_i=0$ on $\buildrel{{\circ}}\over{TM}$, hence $X=0.$ Therefore, ${\rm
Ker}\,f={\rm span}\{\xi_1,\xi_2\}$ and ${\rm rank}\,f=2n-2$,
q.e.d.

\ms

{\bf Theorem 2.2.} {\it The Riemannian metric $G=\displaystyle\frac1{L^2}G_S$
verifies
\begin{equation}
G(fX,fY)=G(X,Y)-\eta^1(X)\eta^1(Y)-\eta^2(X)\eta^2(Y),\
X,Y\in{\cal X}(\buildrel{{\circ}}\over{TM}).
\end{equation}}

\ms

\noindent{\it Proof.} From (2.1) the following local expression of $f$ is
obtained
\begin{equation}
\begin{array}{l}
f(\delta_i)=\left(-\delta^j_i+\displaystyle\frac1{L^2}y_iy^j\right)\dot\partial_j,\vspace*{1,2mm}\\
f(\dot\partial_i)=\left(\delta_i^j-\displaystyle\frac1{L^2}y_iy^j\right)\delta_j.\end{array}
\end{equation}
Using (2.3) one finds
$$\begin{array}{l}
G(f(\delta_i),f(\delta_j))=\displaystyle\frac1{L^2}\left(g_{ij}-\displaystyle\frac1{L^2}y_iy_j\right),\vspace*{1,2mm}\\
G(f,\delta_i),f(\dot\partial_j)=0,\vspace*{1,2mm}\\
G(f(\dot\partial_i),f(\dot\partial_j))=
\displaystyle\frac1{L^2}\left(g_{ij}-\displaystyle\frac1{L^2}y_iy_j\right).\end{array}$$
From here easily follows (2.2). As (2.3) shows the operator $f$
appears as a deformation of $F$ similar with that studied in [1].

\ms

\noindent{\it Remark 2.1.} The metric $G$ is homogeneous on the fibres
of $\buildrel{{\circ}}\over{TM}$ while $G_S$ is not. See [6].

Let us set $\phi(X,Y)=G(fX,Y)$ for $X,Y\in{\cal X}(\buildrel{{\circ}}\over{TM})$. Using
Theorems 2.1, 2.2 one verifies
\begin{equation}
\phi(Y,X)=-\phi(X,Y).
\end{equation}
Thus $\phi$ is a 2-form on $\buildrel{{\circ}}\over{TM}$.

Theorem 2.1 shows that the annihilator of $\phi$ is ${\rm
span}\{\xi_1,\xi_2\}$. A direct calculation gives
$[\xi_2,\xi_1]=\xi_1.$ Hence the distribution ${\rm
span}\{\xi_1,\xi_2\}$ is integrable even if $\phi$ is not closed.
(The annihilator of a closed 2-form is always integrable.) A
calculation in local coordinates leads to
\begin{equation}
\phi=d\eta^1+\varphi,\mbox{ where }\varphi=\displaystyle\frac1{L^4}y_iy_jdx^i\wedge\delta y^j.
\end{equation}
Thus $\phi$ is closed if and only if $\varphi$ is closed and this
happens under strong restrictions on the curvatures of the Cartan
connection. Concluding, $\phi$ is in general an almost
presymplectic structure on $\buildrel{{\circ}}\over{TM}$. Notice that $d\eta^1$ is a
symplectic structure on $\buildrel{{\circ}}\over{TM}$. It appears as a deformation of the
symplectic structure $\phi_S(X,Y)=G_S(FX,Y)$, $X,Y\in{\cal X}(\buildrel{{\circ}}\over{TM})$ since we have
$d\eta^1={\displaystyle\frac{1}{L^2}}\phi_S+2\varphi$.

\section{An almost contact structure on the indicatrix bundle of $F^n$}
\setcounter{equation}{0}

The set $IM=\{(x,y)\in\buildrel{{\circ}}\over{TM}\mid L(x,y)=1\}$ is called the
indicatrix bundle of $F^n$. This set is a submanifold of dimension
$2n-1$ of $\buildrel{{\circ}}\over{TM}$. We show that the framed $f-$structure on
$\buildrel{{\circ}}\over{TM}$ given by Theorem 2.2 induces an almost contact
structure on $\buildrel{{\circ}}\over{TM}$. (This has to be compared with that from
[4].)

It is well--known that $\xi_2=y^i\dot\partial_i$ is normal to $IM$. We
notice that it has the length 1 with respect to $G$. Thus the
vector fields tangent to $IM$ verify $G(X,\xi_2)=0$. Let us
restrict to $IM$ the notions introduced above. Denote the
restrictions putting a bar over that symbol. For $X,Y,...$ vector
fields which are tangent to $IM$ we have:
\begin{itemize}
\item[--] $\overline{\xi}_1=\xi_1$ since $\xi_1$ is tangent to $IM$,
\item[--] $\eta^2\equiv0$ on $IM$ since $\eta^2(X)=G(X,\xi_2)$,
\item[--] $\overline{G}=G_S$ because $L^2=1$ on $IM$,
\item[--] $\overline{f}(X)=F(X)+\overline{\eta}^2(X)\xi_2$ is an endomorphism of
the tangent bundle of $IM$ since $G(\overline{f}(X),\xi_2)=0$.
\end{itemize}
We put $\overline{\xi}=\overline{\xi}_1,\ \overline{\eta}=\overline{\eta}^1$.

\ms

{\bf Theorem 3.1.} {\it The ensemble $(\overline{f},\overline{\xi},\overline{\eta})$
provides an almost contact structure on $IM$, that is the followings
hold:
\begin{itemize}
\item[{\rm(i)}] $\overline{f}^3+\overline{f}=0,\ {\rm rank}\,\overline{f}=2n-2=(2n-1)-1$
\item[{\rm(ii)}] $\overline{\eta}(\overline{\xi})=1,\ \overline{f}(\overline{\xi})=0,\
\overline{\eta}\circ\overline{f}=0$
\item[{\rm(iii)}] $\ov f^2(X)=-X+\ov\eta(X)\ov\xi,$ for $X$ a
vector tangent to $IM$.
\end{itemize}}

\noindent{\it Proof.} All questions follows from those proved in Theorem
2.2 by virtue of the above considerations on the restrictions
to $IM$ of the ensemble $(f,(\xi_a),(\eta^a)),$ $a=1,2.$

From Theorem 2.2 it follows

\ms

{\bf Theorem 3.2.} {\it The Riemannian metric $G_S$ verifies
\begin{equation}
G_S(\ov fX,\ov fY)=G_S(X,Y)-\ov\eta(X)\ov\eta(Y),
\end{equation}
for $X,Y$ vectors tangent to $IM$.}

\ms

One checks that $(\delta_i,\ov f\delta_j)$, $j=1,...,n-1$, is a local
frame on a neighborhood with $y^n\ne0$ on $IM.$ As
the points $(x,0)$ are outside of $IM$ one always may consider
such a local frame.

Let $\Omega(X,Y)=G_S(\ov fX,Y)$ be the 2-form usually associated
to an almost contact structure.

By a direct calculation one gets
\begin{equation}
\begin{array}{l}
d\ov\eta(\delta_i,\delta_j)=0=\Omega(\delta_i,\delta_j)\vspace*{1mm}\\
d\ov\eta(\delta_i,\ov f\delta_j)=g_{ij}-y_iy_j=\Omega(\delta_i,f\delta_j)\vspace*{1mm}\\
d\ov\eta(\ov f\delta_i,\ov f\delta_j)=0=\Omega(\ov f\delta_i,\ov f\delta_j),\end{array}
\end{equation}
in other words $\Omega=d\ov\eta$.

In all these calculation we have used the Cartan connection of the
Finsler space $F^n$.

Thus we have

\ms

{\bf Theorem 3.3.} {\it Let $F^n$ be endowed with the Cartan
connection. Then the structure $(\ov f,\ov\xi,\ov\eta,G_S)$ is a
contact Riemannian structure on $IM$.}

\ms

The structure $(\ov f,\ov\xi,\ov\eta)$ is called {\it normal} if
$N=N_{\ov f}+d\ov\eta\otimes\ov\xi=0,$ where $N_{\ov f}$ is the
Nijenhuis tensor field of $\ov f$. And it is said to be Sasakian
if it is normal and $\Omega=d\ov\eta.$ Again by a direct
calculation one find
$N(\delta_i,\delta_j)=(y_i\delta^h_j-y_j\delta^h_i-R^h{}_{ij})\dot\partial_h$, and
the vanishing of this term implies the vanishing of $N(\ov
f\delta_i,\delta_j)$ and $N(\ov f\delta_i,\ov f\delta_j)$. But
$N(\delta_i,\delta_j)=0$ is equivalent with
\begin{equation}
R^h{}_{ij}=y_i\delta^h_j-y_j\delta^h_i,
\end{equation}
where $R^h{}_{ij}=R_k{}^h{}_{ij}y^k$ and $R_k{}^h{}_{ij}$ is the
$hh-$curvature of the Cartan connection.

The equality (3.3) takes also the form
\begin{equation}
R_{ihk}=g_{ik}y_h-g_{ih}y_k,
\end{equation}
which says that $F^n$ is of constant curvature 1.

Thus we have

\ms

{\bf Theorem 3.4.} {\it Let $F^n$ endowed with Cartan
connection. Then the structure $(\ov f,\ov\xi,\ov\eta, G_S)$ on
$IM$ is a Sasakian structure if and only if $F^n$ is of constant
curvature $1$.}

\ms

For Finsler spaces of constant curvature we refer to [2]. It seems
that the almost contact structure given in Theorem 3.1 is very
close with that obtained in [4]. They have the same properties
(Theorems 3.3 and 3.4).

\bigskip

\hspace*{8.1cm}{\footnotesize{{\it Faculty of Mathematics}}}\\ \hspace*{8cm}{\footnotesize{\it ``Al.~I.~Cuza'' University of Iasi}}


\newpage

\def\scrii{\scriptstyle{\|}}
\def\vsp{\vspace*{1,5mm}\\ }
\def\vspp{\vspace*{3mm}\\ }
\def\vi{\raise 2pt\hbox{,}}
\def\st#1#2{{\mathrel{\mathop{#2}\limits_{#1}}{}\!}}

\runningauthor={M. ANASTASIEI}
\runningtitle={SOME RIEMANNIAN ALMOST PRODUCT STRUCTURES ON TANGENT MANIFOLD}
\noindent
\baselineskip 8pt
\noindent{\footnotesize{Algebras Groups Geom.}}
\hfill\break
{\footnotesize{17, no. 3, 253--262, 2000}}
\vskip 2cm
\baselineskip 11.5pt plus .15pt
\centerline{\bf\Large SOME RIEMANNIAN ALMOST PRODUCT}
\vskip .2cm
\centerline{\bf\Large STRUCTURES ON TANGENT MANIFOLD}
\vskip .5cm
\centerline{\bf {\footnotesize{BY}}}
\vskip .5cm
\centerline{\bf {\footnotesize{M. ANASTASIEI\footnote{Partially supported by CNCSIS Bucure\c sti, Romania}}}}
\vskip 1cm

\begin{abstract}
The tangent manifold $TM$ of a smooth i.e. $C^\infty,$ paracompact
manifold $M$, fibered over $M$ by the natural projection $\tau$,
carries an integrable distribution ${\rm Ker}\,\tau_*$, called {\it
vertical distribution}. If one takes a supplementary distribution
of it, called {\it horizontal}, an almost\\ product structure $P$ on
$TM$ appears. One endows the vertical distribution with a Riemannian
metric $g$. Then $g$ can be prolonged to a Riemannian metric $G$ on
$TM$ in such a way that the pair $(P,G)$ becomes a Riemannian almost
product structure. In this paper we propose a deformation of $P$ suggested the
almost complex case, \cite{1}.
This produces six new Riemannian almost product structures. Some
properties of these structures are pointed out. The particular
case when $g$ is the vertical lift of a Riemannian metric on $M$ is
considered.

{\bf MSC2000 : 53 C  15}
\end{abstract}

\setcounter{section}{0}
\section{A standard Riemannian almost product structure on $TM$}
Let $M$ be a smooth i.e. $C^\infty$ paracompact manifold of dimension
$n$ with local coordinates $(x^i)$, $i,j,k...=1,...,n$. Denote by
$TM$ its tangent manifold with local coordinates $(x^i,y^i)$ and
projection $\tau:TM\to M.$ It is known that $TM$ is also
paracompact. Let $V_uTM={\rm Ker}\,\tau_{*,u}$ for $u\in TM$. Then
$u\to V_uTM$ is an integrable distribution on $TM$, called {\it
vertical distribution} and $VTM=\st u\cup V_uTM$ is the vertical
bundle over $TM.$

Let $HTM$ be a vector bundle over $TM$ which is supplementary to
$VTM$. Such a vector bundle, called {\it horizontal}, always
exists since $TM$ is paracompact. It is said also that it defines
a nonlinear connection on $TM$. Thus we have the decomposition
$$T_uTM=H_uTM\oplus V_uTM.\leqno(1.1)$$
The projectors $h$ and $v$ produced by the direct sum (1.1)
provide an almost product structure $P$ ($AP$--structure for
brevity), given by
$P=h-v.$ Thus we notice
$$P^2=I,\ h=\frac12(I+P),\ v=\frac12(I-P).\leqno(1.2)$$
The horizontal and vertical subspaces in $T_uTM$ are eigenspaces
of $P$ corresponding to the eigenvalues $+1$ and $-1$,
respectively.

The vertical distribution is locally spanned by
$\dot\partial_i:=\displaystyle\frac\partial{\partial y^i}\cdot$ Looking for a basis $(\delta_i)$
in $H_uTM$ in such a way that $\tau_*(\delta_i)=\partial_i:=\displaystyle\frac\partial{\partial
x^i}\vi$ one finds that
$$\delta_i=\partial_i-N_i{}^j(x,y)\dot\partial_j$$ (the sign "$-$" is for
convenience), where the functions $(N_i{}^j(x,y))$ transform by a
change of coordinates $(x^i,y^i)\to(\widetilde x^i,\widetilde y^i)$ as follows:
$$\widetilde N_i{}^k\partial_j\widetilde x^i=\partial_h\widetilde x^k\cdot N_j{}^h-\partial_j\partial_h\widetilde
x^k.y^h.\leqno(1.3)$$ In terminology from \cite{5}, $(N_i{}^j)$ define a
nonlinear connection.

The basis $(\delta_i,\dot\partial_i)$ is adapted to the decomposition
(1.1). Its dual is $(dx^i,\delta y^i)$ for $\delta y^i=
dy^i+N_k{}^i(x,y)dx^k.$ In the adapted basis $P$ takes the form
$$P(\delta_i)=\delta_i,\ P(\dot\partial_i)=-\dot\partial_i,\leqno(1.4)$$
i.e. it has the matrix
$\left(\begin{array}{cc}
I_n&0\\ 0&-I_n\end{array}\right).$

In general, the horizontal distribution is not integrable. The
nonintegrability is measured by the functions $R^k{}_{ij}(x,y)$ from
$$[\delta_i,\delta_j]=R^k{}_{ij}(x,y)\dot\partial_k.\leqno(1.5)$$
The functions $(R^k{}_{ij})$ behave like the components of a tensor of
$M$ i.e. they define a $d$--tensor field. These functions are regarded
as the curvature of the nonlinear connection $(N_j{}^i)$. We
notice for the later use
$$[\delta_i,\dot\partial_j]=\dot\partial_j(N_i{}^k)\dot\partial_k.\leqno(1.6)$$
One says that $P$ is integrable if the horizontal and vertical
distributions are integrable. Thus we have

\bigskip{\bf Theorem 1.1.} {\it The $AP$-structure $P$ is integrable
if and only if $R^k{}_{ij}(x,y)\!=\!0$,\break equivalently the nonlinear connection
$(N_j{}^i)$ is without curvature.}\bigskip

Now let us endow the vertical bundle over $TM$ with a Riemannian
me\-tric $g$. We may do this since $M$ is paracompact. The local
components of $g$, given by $g_{ij}(x,y)=g(u)(\dot\partial_i,\dot\partial_j)$,
$u\in TM$ define a $d$--tensor field of type $(0,2)$, symmetric
and positive defined. In fact $g$ is nothing but a generalized
Lagrange metric introduced by R. Miron and studied by him and his co\-workers, see \cite[5, Ch. X--XII]{}.
The Riemannian metric $g$ may be extended to a Riemannian
metric $G$ on $TM$ given in the form
$$G(u)=g_{ij}(x,y)dx^i\otimes dx^j+g_{ij}(x,y)\delta y^i\otimes \delta y^i,\
u=(x,y)\in TM.\leqno(1.7)$$
It is clear that $H_uTM$ and $V_uTM$ are orthogonal with respect
to $G$. One easily see that $G(PX,PY)=G(X,Y)$ for any vector
fields $X,Y$ on $TM$. Thus we have

\bigskip{\bf Theorem 1.2.} {\it The pair $(P,G)$ is a Riemannian
$AP$--structure on $TM$.}\bigskip

Let $D$ be a linear connection on $TM$ with the torsion $T.$

If $P$ is parallel with respect to $D$ i.e. $D_XP=0,$ then the
Nijenhuis tensor field associated to $P$ takes the form
$$N_P(X,Y)=T(X,Y)+PT(X,PY(-PT(PX,Y)-T(PX,PY)$$ for $X,Y$ vector
fields on $TM$. This form proves

\bigskip{\bf Theorem 1.3.} {\it If the Levi--Civita connection of $G$
makes $P$ parallel, then the $AP$--structure $P$ is integrable.}

\section{Deformations of the Riemannian $AP$--structure $P$}

We set $y_i=g_{ij}(x,y) y^j$ and consider the following deformations of
$P$ $$\begin{array}{l}
P_d(\delta_i)=(\alpha\delta_i^k+\beta y_iy^k)\delta_k,\vsp
P_d(\dot\partial_i)=(\gamma\delta_i^k+\delta y_iy^k)\dot\partial_k,\end{array}\leqno(2.1)$$
for $\alpha,\beta,\gamma,\delta$ functions on $TM$, to be determined in such a way
that $P^2_d\!=\!I$ and $G(P_d\cdot,P_d\cdot)=G(\cdot,\cdot)$. This deformation is suggested
by the almost complex case, see \cite{1}.
The condition $P^2_d=I$ shows that $\alpha,\beta,\gamma,\delta$
have to be so\-lu\-tions of the
following system of equations
$$\begin{array}{ll}
\alpha^2=1,&\beta(2\alpha+\beta F^2)=0,\\
\gamma^2=1,&\delta(2\gamma+\delta F^2)=0,\end{array}\leqno(2.2)$$
for $F^2=g_{ij}(x,y) y^iy^j.$\medskip

This system of equations has sixteen solutions. Inserting them in (2.1)
one finds, leaving aside the trivial $AP$--structures $\pm I$,
fourteen $AP$--structures from which seven are essential, the
other seven differing by a sign from the previous ones. If we put
$A=(A_i^j)=\left(\delta_i^j-\displaystyle\frac2{F^2}y_iy^j\right),$ these
$AP$--structures are given in matrix form as follows\newpage
$$\begin{array}{c}
P_0\equiv P=\left(\begin{array}{ccc}I_n&0\cr 0&-I_n\end{array}\right),\
P_1=\left(\begin{array}{ccc}I_n&0\cr 0&A\end{array}\right),\
P_2=\left(\begin{array}{ccc}I_n&0\cr 0&-A\end{array}\right),\vsp
P_3=\left(\begin{array}{ccc}A&0\cr 0&I_n\end{array}\right),\
P_4=\left(\begin{array}{ccc}A&0\cr 0&-I_n\end{array}\right),\vsp
P_6=\left(\begin{array}{ccc}A&0\cr 0&A\end{array}\right),\
P_6=\left(\begin{array}{ccc}A&0\cr 0&-A\end{array}\right).\end{array}\leqno(2.3)$$

The condition $G(P_d\cdot,P_d\cdot)=G(\cdot,\cdot)$ does not impose any restriction on
$\alpha,\beta,\gamma,\delta$ previously determined. Thus we have

\bigskip{\bf Theorem 2.1.} {\it The pairs $(P_\alpha,G)$, $\alpha=0,...,6$
are seven Riemannian $AP$--structures on $TM$.}

\bigskip{\it Remark 2.1.} The set $\{I,P_0,P_1,...,P_6\}$ has a group
structure given by the table

$$\begin{array}{l|llllllll}
&I&P_0&P_1&P_2&P_3&P_4&P_5&P_6\\\hline
I&I&P_0&P_1&P_2&P_3&P_4&P_5&P_6\\
P_0&P_0&P_1&P_2&P_1&P_4&P_3&P_6&P_5\\

P_1&P_1&P_2&I&P_0&P_5&P_6&P_3&P_4\\
P_2&P_2&P_1&P_0&I&P_6&P_5&P_4&P_3\\
P_3&P_3&P_4&P_5&P_6&I&P_0&P_1&P_2\\
P_4&P_4&P_3&P_6&P_5&P_0&I&P_2&P_1\\
P_5&P_5&P_6&P_3&P_4&P_1&P_2&I&P_0\\
P_6&P_6&P_5&P_4&P_3&P_2&P_1&P_0&I\end{array}$$

This group is commutative. Its proper subgroups are
$\{I,P_{\alpha}\}$ for $\alpha = 0,1,...,6$,
and $\{I,P_0,P_1,P_2\},$ $ \{I,P_0,P_3,P_4\},
\{I,P_1,P_4,P_6\},$ $ \{I,P_0,P_5,P_6\}, $
$\{I,$ $P_1,$ $P_3,$ $P_5\},$
$\{I,P_2,P_3,P_6\}, $ $\{I,P_2,P_4,P_5\}.$
The last seven are isomorphic with the Klein group.
The group can be also seen as a Burnside group $B(2,3)$ generated
by $\{P_0,P_3,P_5\}$

Let $h_\alpha=\displaystyle\frac12(I+P_\alpha)$ and
$V_\alpha=\displaystyle\frac12(I-P_\alpha)$ be the projectors
defined by $P_\alpha$ and let us set $H_\alpha={\rm Ker}\, v_\alpha,$
$v_\alpha={\rm Ker}\, h_\alpha$, for $\alpha=0,1,...,6$ with $h_0=h,$
$v_0=v,$ $H_0=H,$ $V_0=V.$

For identifying the distributions $H_\alpha$ and $V_\alpha$,
$\alpha=1,...,6$, we consider the vector fields
$C=y^i\dot\partial_i$ and $S=y^i\delta_i$ and denote by the same
letters the $1$-dimensional distributions defined by them.

Furthermore, we denote by $C^\bot$ the orthocomplement of $C$ in $V$,
that is, $C^\bot=\{A^i\dot\partial_i\mid g_{ij}y^iA^j=A^jy_j=0\}$
and by $S^\bot$ the orthocomplement of $S$ in $H$, that is,
$S^\bot=\{X^i\delta_i\mid g_{ij}y^iX^j=X^jy_i=0\}$. With this
notations the following result holds.

\bigskip{\bf Theorem 2.2.} {\it The distributions
defining $P_\alpha$ are as follows:}
$$\left\{\begin{array}{l}
H_0=H\\
V_0=V\end{array}\right.\
\left\{\begin{array}{l}
H_1=H\oplus C^\bot\\
V_1=C\end{array}\right.\
\left\{\begin{array}{l}
H_2=H\oplus C\\
V_2=C^\bot\end{array}\right.\
\left\{\begin{array}{l}
H_3=V\oplus S^\bot\\
V_3=S\end{array}\right.$$
$$\left\{\begin{array}{l}
H_4=S^\bot\\
V_4=V\oplus S\end{array}\right.\
\left\{\begin{array}{l}
H_5=[S]^\bot\oplus[C]^\bot\\
V_5=[S]\oplus[C]\end{array}\right.\
\left\{\begin{array}{l}
H_6=S^\bot\oplus C^\bot\\
V_6=S\oplus C\end{array}\right.$$

\noindent{\it Proof.} For $\alpha=1,$ we have
$h_1(\delta_i)=\delta_i$,
$h_1(\dot\partial_i)\!=\!\left(\delta_i^k-
\displaystyle\frac1{F^2}y_iy^k\right)\dot\partial_k,\,
v_i(\delta_i)\!=\!0,$ $v(\dot\partial)=
\displaystyle\frac1{F^2}y_iy^k\dot\partial_k$. From these
equations one gets
$$\begin{array}{ll}
V_1&={\rm Ker}\, h_1=\left\{X^i\delta_i+A^i\dot\partial_i\mid
X^i\delta_i+ \left(\delta_i^k-
\displaystyle\frac1{F^2}y_iy^k\right)A^i\dot\partial_k=0\right\}=\vsp
&=\left\{X^i\delta_i+A^i\dot\partial_i\mid X^i =0,\
A^k=\displaystyle\frac1{F^2}(A^iy_i)y^k\right\}=C\mbox{ and}\vsp
H_1&={\rm Ker}\, v_1=\left\{X^i\delta_i+A^i\dot\partial_i\mid
A^iy_i=0\right\}=H\oplus C^\bot.\end{array}$$
Similarly, one finds the other distributions.

A study of the integrability of the distributions $V_\alpha$,
$\alpha=0,1,...,6$ gives

\bigskip{\bf Theorem 2.3.}
\begin{itemize}
\item[1)] {\it The distributions $V_0=V,V_1,V_3$ are always
integrable.}
\item[2)] {\it The distribution $V_2$ is integrable if there
exists a real function $L$ on $TM$ such that
$g_{ij}(x,y)=\displaystyle\frac12\dot\partial_i\dot\partial_jL.$}
\item[3)] {\it The distribution $V_5$ is integrable if the
nonlinear connection $(N^i_j)$ is positively homogeneous of degree
$1$.}
\item[4)] {\it The distributions $V_4$ and $V_6$ are never integrable.}
\end{itemize}

\noindent{\it Proof.} 1) We have noticed that $V$ is integrable.
The distributions $V_1$ and $V_3$ are
1-dimensional, hence integrable.\medskip

2) Let $A=A^i\dot\partial_i$ and $B=B^j\dot\partial_j$ in $C^\bot$
i.e. $A^iy_i=0,$ $B^iy_j=0$. Then $[A,B]\in C^\bot$ if and only if
$A^i\dot\partial_i(B^j)y_j-B^i\dot\partial_i(A^j)y_j=0,$
equivalently
$A^jB^i\dot\partial_i(y_j)-A^iB^j\dot\partial_i(y_j)=0.$ This
condition identically holds if
$\dot\partial_ig_{jk}=\dot\partial_jg_{ik}.$\medskip

3) It is clear that the distribution $S\oplus C$ is integrable if
$[C,S]$ belongs to it. We hve
$[C,S]=S+y^j(N^k_j-y^i\dot\partial_i(N^k_j))\dot\partial_k=S$, if
the functions $(N^k_j(x,y))$ are positively homogeneous of degree
1 with respect to $(y^i)$.\medskip

4) A direct calculation.

\bigskip{\it Remark 2.2.} If $g_{ij}(x,y)$ is the metric
tensor of a Finsler space and $(N^i_j(x,y))$ is the Cartan
nonlinear connection, the hypothesis in 2) and 3) of Theorem 2.3
are satisfied and so in this case the distributions
$V_0,V_1,V_3,V_5$ are integrable.\bigskip

The distributions $H_\alpha,$ $\alpha=0,1,...,6$ are not
integrable or they are so in very strong conditions. We renounce
to write down such conditions. A Riemannian $AP$--structure is
integrable if the both distributions defining it are integrable.
From the above it follows the Riemannian $AP$--structures $P_4$
and $P_6$ are never integrable. The others are integrable only
under some strong conditions on $(g_{ij})$ and $(N^i_j)$.

The Riemannian $AP$--structures were classified by A.M. Naveira
\cite{6}. Modulo a duality there exists thirty--six different
classes described by conditions on $\nabla h_\alpha$, where
$\nabla$ denotes the Levi--Civita connection of $G$. See also
\cite{4}.

From \cite{2} it follows that the Levi--Civita connection $\nabla$ can be
taken in the form
$$\begin{array}{ll}
\nabla_{\delta_k}\delta_j=F^i_{jk}\delta_i+A^a_{jk}\dot\partial_a,
&\nabla_{\dot\partial_b}\delta_j=\widetilde C_j{}^i{}_b\delta_i+E_j{}^a{}_b\dot\partial_a,\vsp
\nabla_{\delta_k}\dot\partial_b=L^a_{bk}\dot\partial_a+D^i_{bk}\delta_i,
&\nabla_{\dot\partial_b}\dot\partial_c=C^a_{cb}\dot\partial_a+B^i_{cb}\delta_i,
\end{array}\leqno(2.1)$$
with the connection coefficients given by
$$\begin{array}{l}
F^i_{jk}=\displaystyle\frac12 g^{ih}(\delta_jg_{hk}-\delta_kg_{jh}-\delta_hg_{jk}),\vsp
A^a_{jk}=\displaystyle\frac12(-R^a_{jk}-g^{ab}\dot\partial_bg_{jk}),\vsp
\widetilde C_j{}^i{}_b=D^i_{bj}=\displaystyle\frac12g^{ih}(\dot\partial_bg_{jh}+g_{bc}R^c{}_{hj}),\vsp
E^a_{ib}=\displaystyle\frac12g^{ac}g_{bc\scrii i},\ L^a_{bi}=\dot\partial_bN^a_i+
\displaystyle\frac12g^{ac}g_{bc\scrii i},\vsp
B^k_{ab}=-\displaystyle\frac12g^{kj}g_{ab\scrii j},\ C^a_{bc}=
\displaystyle\frac12g^{ad}(\dot\partial_bg_{dc}+\dot\partial_cg_{bd}-\dot\partial_dg_{bc}).
\end{array}\leqno(2.2)$$
Here $g_{bc\scrii i}$ denotes the $h$--covariant derivative
with respect to the Berwald connection i.e.
$$g_{bc\scrii i}=\delta_ig_{bc}-(\dot\partial_bN_i^d)g_{dc}
-(\dot\partial_cN_i^d)g_{bd}.$$
We do not classify $P_\alpha$ here but we remark that $P_4$ and
$P_6$ cannot be in Naveira's class $\cal P.$ Indeed, the
conditions $\nabla P_\alpha=0$ characterizing $\cal P$ implies in
virtue of Theorem 1.3 that $P_\alpha$ should be integrable. But
$P_4$ and $P_6$ are never integrable.

A distribution ${\cal D}$ is geodesically
invariant if all geodesics with initial vector in ${\cal D}$ remain
tangent to $D$ for all time. As it was proved
in \cite{4}, a distribution ${\cal D}$ is geodesically invariant if and only if for
any sections $X,Y$ of ${\cal D}$, the symmetric product
$X:Y=\nabla_XY+\nabla_YX$ is again a section of $D$. See also \cite{1}.
Using (2.1) and (2.2) one gets

\bigskip{\bf Theorem 2.4.}
\begin{itemize}
\item[1)] {\it The distribution $H$ is geodesically invariant if
and only if $g_{ij\scrii k}=0.$}
\item[2)] {\it The distribution $V$ is geodesically invariant if
and only if $\dot\partial_kg_{ij}=0.$}
\end{itemize}

The first condition in Theorem 2.4 tells us that we have to assume no
dependence on $y=(y^i)$ in $(g_{ij})$. In this case $g(g_{ij})$
reduces to a Riemannian metric on $M$.  If we continue to work with an arbitrary nonlinear
connection $(N^i_j)$, the second condition in Theorem 2.4 is not verified. But
if we take $N^i_j(x,y)=\gamma^i_{jk}(x)y^k$, where $(\gamma^i_{jk})$ are
Christoffel symbols derived from $(g_{ij}(x))$, then the condition
$g_{jk\scrii h}=0$ reduces to
$\partial_hg_{jk}-\gamma^i_{hj}g_{ik}-\gamma^i_{hk}g_{ji}=0$ which obviously
holds. We consider this case in what
follows.
Thus we may state

\bigskip{\bf Corollary 2.1.} {\it Let $(M,g)$ be a Riemannian manifold.
One considers $TM$ endowed with the nonlinear connection
$N^i_j(x,y)=\gamma^i_{jk}(x)y^k$ and one defines the Riemannian metric $G$
with this nonlinear connection $(G$ becomes the Sasaki metric of
$g)$. Then the distributions $V$ and $H$ are geodesically in\-variant.}\bigskip

We associate to every $P_\alpha,$ $\alpha=0,1,...,6$, the symmetric tensor
field $\phi_\alpha$ defined by
$$\phi_\alpha(X,Y)=G(P_\alpha X,Y),\ X,Y\in{\cal X}(TM).\leqno(2.6)$$
Using the matrix form of $P_\alpha$ one obtains
$$\begin{array}{l}
\phi_1(u)=g_{ij} dx^i\otimes dx^j-g_{ij}\delta y^i\otimes \delta y^j,\\
\phi_1(u)=g_{ij} dx^i\otimes dx^j+\left(g_{ij}-\displaystyle\frac2{F^2}y_iy_j\right)\delta
y^i\otimes \delta y^j,\vsp
\phi_2(u)=g_{ij} dx^i\otimes dx^j-\left(g_{ij}-\displaystyle\frac2{F^2}y_iy_j\right)\delta
y^i\otimes \delta y^j,\end{array}$$
and so on.

As $\left(g_{ij}-\displaystyle\frac2{F^2}y_iy_j\right)$ is an invertible matrix,
with the inverse\break $\left(g^{jk}-\displaystyle\frac2{F^2}y^iy^j\right)$, it follows
that $\phi_\alpha$, $\alpha=0,1,...,6$ are pseudo--Riemannian metrics on $TM$.

If we look at the first terms in $\phi_0,\phi_1,\phi_2,$ it
appears as obvious

\bigskip{\bf Theorem 2.5.} {\it The maps $\tau:(TM,\phi_\alpha)\to(M,g)$,
$\alpha=0,1,2$ are Riemannian submersions.}\bigskip

A fundamental tensor field of the submersions $\tau_\alpha$ is
$$S_\alpha(X,Y)=h_\alpha\nabla_{v_\alpha X}v_\alpha
Y+v_\alpha\nabla_{v_\alpha X}h_\alpha Y,\ X,Y\in{\cal X}(TM)$$ and the submersion $\tau_\alpha$ is called {\it
totally geodesic} if $S_\alpha$ identically vanishes. Here $\alpha=0,1,2.$

A tedious calculation in which the identity $y_iy^k\dot\partial
\left(\displaystyle\frac1{F^2}y_iy^h\right)=0$ is used proves

\bigskip{\bf Theorem 2.6.} {\it The Riemannian submersions
$\tau_\alpha:(TM,\phi_\alpha)\to(M,g)$ are totally geodesic,
$\alpha=0,1,2$.}

\newpage

\def\llr{\Longleftrightarrow }
\def\bigi{\big|}
\def\scri{\scriptstyle{|}}
\def\bld#1#2{{\buildrel{#1}\over{#2}}}
\def\({\left(}
\def\){\right)}
\def\io{if and only if }
\def\tmd{$TM$}
\def\ij{_{ij}}
\def\ji{_{ji}}
\def\ijk{^i_{jk}}
\def\oo{{\omega}}
\def\ooo{{\Omega}}
\def\pp{{\partial}}
\def\de{{\delta}}
\def\G{{\Gamma}}
\def\g{{\gamma}}
\def\9{{\infty}}
\def\fcs{functions}
\def\coefs{coefficients}
\def\fwg{following}
\def\0{\leqno}
\def\lagn{Lagrangian}
\def\frax#1#2{\displaystyle\frac{#1}{#2}}
\def\1{^{-1}}
\def\w{\wedge }
\def\bk{\bigskip }
\def\n{\noindent }
\def\cds{conditions}
\def\barr{\begin{array}}
\def\earr{\end{array}}

\runningauthor={M. ANASTASIEI}
\runningtitle={SYMPLECTIC CONNECTIONS IN LAGRANGE GEOMETRY}
\noindent
\baselineskip 8pt
\noindent{\footnotesize{Finsler and Lagrange geometries(Ia\c si,2001)}}
\hfill\break
{\footnotesize{9–16, Kluwer Acad. Publ., Dordrecht, 2003}}
\vskip 2cm
\baselineskip 11.5pt plus .15pt
\centerline{\bf\Large SYMPLECTIC CONNECTIONS}
\vskip .2cm
\centerline{\bf\Large IN LAGRANGE GEOMETRY}
\vskip .5cm
\centerline{\bf {\footnotesize{BY}}}
\vskip .5cm
\centerline{\bf {\footnotesize{M. ANASTASIEI}}}
\vskip 1cm

\begin{abstract}
A Lagrangian structure and, in particular, a Finslerian or a Riemannian structure on a
manifold $M$ induces a symplectic structure on $TM$. We investigate the linear
connections on $TM$ depending on the Lagrangian structure only, which are compatible
with this symplectic structure and have no torsion.\\ {\bf MSC 2000:} 53C60,53D99
\end{abstract}

\setcounter{section}{0}
\section*{Introduction}
On $TM$ we have the vertical bundle as the kernel of the differential of the
projection $TM\to M.$ We take a supplement of it, that is, a horizontal bundle or a
nonlinear connection and we consider the natural almost complex product $F$ on $TM$
associated to these bundles. We show in the first section of this paper that the
symplectic structures on $TM$ having the horizontal and vertical bundles as Lagrangian
subbundles and being compatible with $F$ are essentially induced by a Lagrangian
structure on $M$. In the second section we state precisely the symplectic structure
$\ooo_L$ induced by a Lagrangian structure $L$. In the third section we are interested
in linear connections on \tmd\ which are compatible with $\ooo_L$ and have no torsion,
called symplectic connections. First, we notice that the liner Cartan connection of
$L$ is compatible with $\ooo_L$ but has torsion. By a deformation of it we find a set
of symplectic connections from which set one depending on $L$ only is single out. In
the case that the horizontal distribution is integrable, a set of symplectic
connections preserving by parallelism the vertical and horizontal bundles is found. I.
Vaisman discovered \cite{3} three classes of symplectic connections: flat, Ricci flat
and with reducible curvature. Among the symplectic connections that we have found, two
flat symplectic connections and a Ricci flat one are pointed out in section 4. We
refer to the monograph \cite{2} for notations and terminology.

\section{Symplectic structures on \tmd}
We shall work in the category of real, smooth, i.e. $C^\9$ and
finite dimensional manifolds. Let $M$ be a manifold of dimension
$n$ and $\tau:TM\to M$ its tangent bundle. Let $(U,(x^i))$,
$i=1,2,...,n$ be a coordinate chart on $M$. Then
$(\tau\1(U),(x^i\circ\tau,y^i)),$ where $(y^i)$ are the components
of a tangent vector $v_x,$ $x\in U,$ in the natural basis
$\pp_i:=\frax\pp{\pp x^i}$ of $T_xM$, is a coordinate chart on
\tmd. The indices $i,j,k...$ will range from 1 to $n$ and the
Einstein convention on summation will be used.

Let $\tau_*:TTM\to TM$ be the differential of $\tau.$ The union of
$V_uTM:=\ker\tau_{*,n}$ for $u\in TM$ defines the vertical bundle
over \tmd. We may thought it as a distribution on \tmd\ called
vertical distribution. This is locally spanned by
$\dot\pp_i:=\frax\pp{\pp y^i}$, hence it is integrable. Thus we
may speak about vertical foliation whose leaves are $T_xM$, $x\in
M.$ A {\it non--linear connection} $N$ is a subbundle $HTM$ of
$TTM$, called horizontal, that is supplementary to the vertical
bundle, i.e. the \fwg\ decomposition holds
$$TTM=VTN\oplus HTM\mbox{ (Whitney's sum)}\0(1.1)$$
We also view the horizontal subbundle as a distribution $u\to
H_uTM$ called the {\it horizontal distribution} on \tmd.

Locally, we shall use the {\it adapted bases} $(\de_i,\dot\pp_i)$,
where
$$\de_i=\pp_i-N^k_i(x,y)\dot\pp_km\0(1.2)$$
span the horizontal distribution, and their dual cobases
$(dx^i,\de y^i)$, where $$\de y^i=dy^i+N^i_k(x,y)dx^k.\0(1.3)$$
The functions $(N^i_k)$ are called the local \coefs\ of the
non--linear connection $N$. If these \fcs\ are linear with respect
to $(y^i)$, that is, $N^i_k(x,y)=\Gamma^i_{kj}(x)y^j$, it comes
out that $(\Gamma^i_{kj}(x))$ are the local \coefs\ of a linear
connection on $M$.

The tensor fields on \tmd\ get a natural multiple grading induced by (1.1). When this
is made explicit by the use of the adapted bases and their dual cobases, the \coefs\
of the components are \fcs\ depending on $(x,y)$ but transform under a change of
coordinates on \tmd\ as tensors on $M$. it is said in \cite{2} that these components
or their \coefs are $d$-tensor fields on \tmd. here $d$ is for ``distinguished''. In
particular, for the spaces of differential forms we have
$$\wedge^k(TM)=\oplus_{p+q=k}\wedge^{pq}(TM),\0(1.4)$$ where $p$ is the $V$-degree and
$q$ is the $H$-degree. Thus any 2-form $\ooo$ on \tmd\ can be written as
$$\ooo=\frac12b_{ij}(x,y)dx^i\w dx^j+a\ij(x,y)dx^i\w\de y^i+ \frac12 c\ij(x,y)\de
y^i\w \de y^j,\0(1.5)$$ with $b\ij=-b\ji,\ c\ij=-c\ji.$ Each term in (1.5) is a
distinguished 2-form on \tmd. The \coefs\ $a\ij(x,y)$, $b\ij(x,y)$, $c\ij(x,y)$
transform under a change of coordinates on \tmd\ as the components of tensors on $M$,
the last two being skew symmetric.

Let us suppose that $\ooo$ given by (1.5) defines a symplectic structure o \tmd. From
$$\barr{ll} \ooo(\de_i,\de_j)=b\ij,&\ooo(\de_i,\dot\pp_j)=a\ij,\\
\ooo(\dot\pp_i,\de_j)=-a\ji,&\ooo(\dot\pp_i,\dot\pp_j)=c\ij,\earr\0(1.6)$$ it comes
out that the vertical (horizontal) bundle is a Lagrangian subbundle with respect to
$\ooo$ \io $c\ij=0$ $(b\ij=0)$. In the sequel we shall be interested only in
symplectic structures on \tmd\ that make the vertical and horizontal bundles the
Lagrangian subbundles of $TTM$. Thus we consider only the symplectic structures on
$TM$ given by the 2-forms $$\ooo=a\ij(x,y)dx^i\wedge\de y^i,\0(1.7)$$ satisfying the
\cds $$\det(a\ij(x,y))\ne0\llr\ooo\mbox{ is nondegenerate},\0(1.8)$$ $$\barr{l}
\sum_{(ijk)}a_{ih}R^h{}_{jk}=0,\
\de_ia_{jk}+a_{ih}\dot\pp_kN^h_j=\de_ja_{ik}+a_{jk}\dot\pp_kN^h_i,\\
\dot\pp_ka\ij=\dot\pp_ja_{ik},\earr\0(1.9)$$ where
$$[\de_j,\de_k]=R^h{}_{jk}\dot\pp_h,\ R^h{}_{jk}=\de_kN^h_j-\de_jN^h_k.\0(1.10)$$ The
eqs. (1.9) are equivalent with $d\ooo=0$. The \fcs\ $(R^h{}_{jk}(x,y))$ define a
$d$-tensor of type $(1,2)$. It vanishes \io the horizontal distribution is integrable.

Now we consider the almost complex structure $F$ on \tmd\ defined
by
$$F(\de_i)=-\dot\pp_i,\ F(\dot\pp_i)=\de_i.\0(1.11)$$
Let $\chi(TM)$ the set of vector fields on \tmd. It is easy to
check

\bk\n{\bf Proposition 1.1.} {\it For $X,Y\in\chi(TM)$ we have
$$\ooo(FX,FY)=\ooo(X,Y),\0(1.12)$$ \io $a\ij=a\ji.$}\bk

We confine ourselves to the case when $\ooo$ from $(1.7)$
satisfies $(1.12)$. We put $a\ij=-g\ij$ with $g\ij=g\ji$ and we
write $\ooo$ in the form
$$\ooo=g\ij(x,y)\de y^i\wedge dx^j.\0(1.13)$$

The $d$-tensor field $g=g\ij(x,y)\de y^i\otimes\de y^j$ with $\det(g\ij)\ne0$ and such
that the quadratic form $g\ij\xi^i\xi^j$, $(\xi^i)\in\rr^n,$ has constant signature,
is called a generalized Lagrange metric, shortly a $GL$-metric, \cite{2}. One may
consider also the $d$-tensor field $g\ij9x,y)dx^i\otimes dx^j$ which summed with $g$
gives a metrical structure on \tmd: $$G=g\ij(x,y)dx^i\otimes dx^j+g\ij(x,y)\de
y^i\otimes \de y^j.\0(1.14)$$ One easily verifies

\bk\n{\bf Proposition 1.2} {\it For every $X,Y\in\chi(TM)$ one has
$$G(X,Y)=\ooo(X,FY),\0(1.15)$$ $$G(FX,FY)=G(X,Y).\0(1.16)$$}

Thus the pair $(F,G)$ is an almost Hermitian structure o \tmd\ and
$\ooo$ appears as its fundamental 2-form. As $d\ooo=0,$ we have
that $(F,G)$ is an almost K\"ahler structure. It reduces to a
K\"ahler structure \io $R^h{}_{jk}=0$ and
$\dot\pp_kn^i_h=\dot\pp_hN^i_k$, cf. \cite{2}, Ch.7.

The \fcs\ $(g\ij(x,y))$ have to satisfy the \cds
$$\barr{l}
\sum_{(ijk)}R_{ijk}=0,\
\de_ig_{jk}+g_{ih}\dot\pp_kN^h_j=\de_jg_{ik}+g_{jh}\dot\pp_kN^h_i,\\
dot\pp_kg\ij=\dot\pp_jg_{ik},\mbox{ for
}R_{ijk}:=g_{ih}R^h{}_{jk}, \earr\0(1.9)'$$
in order that $d\ooo=0$ for $\ooo$ given by (1.13).

The third equality in $(1.9)'$ holds \io
$$g\ij(x,y)=\frax12\dot\pp_i\dot\pp_jL(x,y),\mbox{ for some
function }L\mbox{ on }TM.\0(1.17)$$
We shall take the assumption (1.17) for the rest of this paper.

\section{Lagrangian symplectic structures on \tmd}

We call a Lagrangian structure on $M$ a regular Lagrangian on \tmd, that is a function
$L:TM\to R$ such that the matrix $(g\ij(x,y))$ given by (1.17) has $\det(g\ij)\ne0$
and the quadratic form $g\ij(x,y)\xi^i\xi^j$, $\xi\in\rr^n$, is of constant signature
on \tmd. The pair $(M,L)$ is called a Lagrange manifold. We send to the monograph
\cite{2} for the geometry of these manifolds. It is known that a Lagrangian structure
determines a non--linear connection. This can be constructed as follows,
\cite[Ch.IX]{2}. The Euler--Lagrange equations for $L$ take the form
$$\frax{d^2x^i}{dt^2}+2G^i\(x,\frax{dx}{dt}\)=0,\0(2.1)$$ The \fcs\ $(-2G^i9x,y))$ are
the \coefs\ of a semispray (second order differential equation) on $M$ and one proves
that $$N^i_j(x,y)=\dot\pp_jG^i(x,y),\0(2.3)$$ are the local \coefs\ of a non--linear
connection $N_L\circ$ \tmd.

Now we may consider the adapted bases and their dual cobases with
respect to $N_L.$ We keep the notations from the first section but
we refer now to $N_L$ only.

Thus for the Lagrange manifold $(M,L)$ we have $g\ij(x,y)$ given by (1.17) and
$(N^i_j(x,y))$ given by (2.3). The symplectic structure $$\ooo_L=g\ij(x,y)\de y^i\w
dx^j,\0(2.4)$$ will be called a Lagrangian symplectic structure.

That $\ooo_L$ is indeed a symplectic structure and not only an
almost symplectic one it follows from

\bk\n{\bf Proposition 2.1.} {\it $\ooo_L=d\oo_L$ for
$\oo_L=\frax12(\dot\pp_jL)dx^j.$}\bk

\ms\n{\it Proof.} We have $d\oo_L=\frax12(\pp_i\dot\pp_jL)dx^i\w
dx^j+g\ij dy^i\w dx^j.$ Inserting here $dy^i=\de y^i-N^i_kdx^k,$
one gets
$$d\oo_L=\(\frax12\pp_i\dot\pp_jL-g_{kj}N^k_i\)dx^i\w dx^j+
g\ij\de y^i\w dx^j=\ooo_L$$
because of the symmetry in the indices $i,j$ of
$A\ij=\frax12\pp_i\dot\pp_jL-g_{kj}N^k_i.$ Indeed, a direct
calculation gives
$$4A\ij=(\pp_i\dot\pp_jL+\dot\pp_i\pp_jL)-2y^s\pp_sg\ij+
4G^k\dot\pp_kg\ij,\ \mbox{q.e.d.}$$

On the other hand the condition $d\ooo_L=0$ is equivalent with
$(1.9)'$ written for $(g\ij)$ given by (1.17) and $(N^j_i)$ given
by (2.3). By Proposition 2.1 the \cds\ $(1.9)'$ become identities.

\section{Symplectic connections for $(TM,\ooo_L)$}

A linear connection $\nabla$ on \tmd\ endowed with $\ooo_L$ will be called almost
symplectic if $\nabla\ooo_L=0$. The term of symplectic connection is reserved for
those almost symplectic connections which have no torsion. It is known that any
symplectic manifold admits infinitely many almost symplectic connections and
infinitely many symplectic connections. See \cite{3} for a clear review of this
matter.

For a Lagrange manifold $(M,L)$, besides the non--linear connection $N_L$ we have a
canonical linear connection determined by $L$ only, called the linear Cartan
connection. We recall it following  \cite{2} and \cite{1}.

Locally, this has the form
 $$\barr{ll} \bld
cD_{\de_k}\de_j=F^i_{jk}\de_i, &\bld cD_{\dot\pp_k}\de_j=C^i_{jk}\de_i,\\ \bld
cD_{\de_k}\dot\pp_j=F^i_{jk}\dot\pp_i,& \bld
cD_{\dot\pp_k}\dot\pp_j=C^i_{jk}\dot\pp_i. \earr\0(3.1)$$
 where
$$F\ijk=F^i_{kj},\ C\ijk=C^i_{kj},\0(3.2)$$ and the condition that $\bld cD$ is
metrical with respect to $G$  is equivalent to $$\barr{l} g_{ij\scri
k}:=\de_kg\ij-F^h_{ik}g_{jh}-F^h_{jk}g_{ih}=0,\\ g_{ij}|_k
=\dot\pp_kg\ij-C^h_{ik}g_{jh}-C^h_{jk}g_{ih}=0.\earr\0(3.3)$$ Then, from (3.2) and
(3.3), $F\ijk$ and $C\ijk$ are uniquely determined in the form $$\barr{l}
F\ijk=\frax12g^{ih}(\de_jg_{hk}+\de_kg_{hj}-\de_hg_{jk}),\\
C\ijk=\frax12g^{ih}(\dot\pp_jg_{hk}+\dot\pp_kg_{hj}-\dot\pp_hg_{jk})
=\frax12g^{ih}\dot\pp_hg_{jk}\earr\0(3.4)$$ Notice that although $\bld cD$ is a
metrical connection, it does not coincide with the Levi--Civita connection of $G$
since it has torsion. Indeed, we have $$T(\de_k,\de_j)=R\ijk\dot\pp_i,\
T(\dot\pp_k,\de_j)=C\ijk\de_i+P\ijk\dot\pp_i,\0(3.5)$$ where
$P\ijk=\dot\pp_k(N^i_j)-F^i_{kj}.$ The $d$-tensor fields $R\ijk,$ $C\ijk$, $P\ijk$
vanish only for very particular Lagrangians. For instance, if we consider a Riemannian
metric $(\g\ij(x))$ on$M$ and we put $$L(x,y)=\g\ij(x)y^iy^j,\0(3.6)$$ we obtain a
Lagrangian for which $P\ijk=C\ijk=0$ but $R\ijk\ne0$ unless if the Riemannian metric
$(\g\ij(x))$ is flat. Now we can prove a simple but important result.

\bk\n{\bf Theorem 3.1.} {\it The linear Cartan connection of the
symplectic manifold $(TM,\ooo_L)$ is an almost symplectic
connection.}\bk

\ms\n{\it Proof.} Using $\ooo_L(X,FY)=G(X,Y)$ and $\bld cDF=0$ in the form $\bld
cD_XFY=F\bld cD_XY$, one easily obtains $$(D_X\ooo_L)(Y,Z)=(D_X G)(FY,Z)=0,\mbox{ for
every } X,Y,Z\in\chi(TM),\mbox{ q.e.d.}$$

We notice that $\bld cD$ is completely determined by $L$. For the Lagrangian (3.6), it
reduces to the Levi--Civita connection of the Riemannian metric $(\g\ij(x)).$

In the proof of Theorem 3.2 we have used the both \cds\ $\bld cDF=0$ and $\bld cDG=0.$
We may ask if there exists an almost symplectic connection on \tmd preserving the
horizontal and vertical distributions i.e. a distinguished linear connection
satisfying only one or none from these two \cds.

Let $D$ be a distinguished linear connection, shortly a
$d$-connection,on \tmd. We follow the theory from \cite[Ch.3]{2}
where such connections are considered on the total space of a
vector bundle.

Using the projectors $h$ and $v$, we have the following decomposition of $D$:
$$D_XY=hD_{hX}hY+vD_{vX}vY+hD_{vX}hY+vD_{hX}vY,\ X,Y\in\chi(TM).\0(3.7)$$ When we take
$$hD_{vX}hY=h[vX,hY],\ vD_{hX}vY=v[hX,vY],\0(3.8)$$ the definition of a connection is
respected. The first two terms in the right hand of (3.7) will be determined from the
condition $D\ooo_L=0,$ which gives $$\barr{l} \ooo_L(hD_{hX}hY,hZ)=0\\
\ooo_L(hD_{hX}hY,vZ)=(hX)\ooo_L(hY,vZ)-\ooo_L(hY,v[hX,vZ]), \earr\0(3.9)$$ $$\barr{l}
\ooo_L(vD_{vX}vY,vZ)=0\\ \ooo_L(vD_{vX}vY,hZ)=(vX)\ooo_L(hY,vZ)-\ooo_L(h[hX,vZ],vY).
\earr\0(3.10)$$ With $hD_{hX}hY$ and $vD_{vX}vY$ uniquely determined from (3.9) and
(3.10), respectively, and (3.8), the condition $D\ooo_L=0$ holds.

The $d$-connection $D$ is not an almost complex one i.e. $DF\ne0$ nor a metrical one
i.e. $DG\ne0.$

From (3.7)--(3.10) and $d\ooo_L=0$ it follows that the torsion of$D$ vanishes \io
$v[hX,hY]=0,$ that is, the horizontal distribution is integrable. Thus $D$ becomes a
symplectic connection only in a very restrictive condition. Now we wish to avoid it
and in order to do so we have to renounce  to the condition that $D$ is a
$d$-connection. We shall determine a set of symplectic connections for $(TM,\ooo_L)$
as a subset of all linear connections on \tmd\ and we single out one which is
completely determined by $L$.

\bk\n{\bf Theorem 3.2.} {\it There exists a linear connection $\nabla$ on \tmd\ which
is almost symplectic with respect to $\ooo_L$, is without torsion and depends on $L$
only.}\bk

\ms\n{\it Proof.} Having the linear Cartan connection $\bld cD$, we set
$\nabla_XY=\bld cD_XY+A(X,Y)$ for some tensor field $A$ of type $(1,2)$ and $X,Y\in
\chi(TM).$ The condition that the torsion of $\nabla$ vanishes reads
$$T(X,Y)+A(X,Y)-A(Y,X)=0,\0(3.11)$$ and since $\bld cD$ is an almost symplectic
connection, $\nabla$ is almost symplectic connection \io
$$\ooo_L(A(X,Y),Z)+\ooo_L(Y,A(X,Y))=0,\ X,Y,Z\in\chi(TM).\0(3.12)$$ Locally, we put
$$\barr{ll} A(\de_k,\de_j)=A\ijk\dot\pp_i,& A(\dot\pp_k,\de_j)=E\ijk\dot\pp_i,\\
A(\de_k\dot\pp_j)=D\ijk\de_i,& A(\dot\pp_k,\dot\pp_j)=B\ijk\de_i. \earr\0(3.13)$$ Thus
we already took a particular form of $A$. Then (3.11) is equivalent with
$$R\ijk=A\ijk-A^i_{kj},\ E\ijk=-P^i_{kj},\ D\ijk=C\ijk,\ B\ijk=B^i_{kj},\0(3.14)$$ and
(3.12) is equivalent to $$\barr{ll} A^h_{ik}g_{hj}-A^h_{jk}g_{hi}=0,&
D^h_{ik}g_{hj}-D^h_{jk}g_{hi}=0,\\ E^h_{ik}g_{hj}-E^h_{jk}g_{hi}=0,&
B^h_{ik}g_{hj}-B^h_{jk}g_{hi}=0. \earr\0(3.15)$$ The tensorial eqs. (3.15) could be
solved using the Obata operators associated to $(g\ij(x,y)).$ For brevity, we shall
not introduce them here. Instead, we check that the system of eqs. (3.14) and (3.15)
has the solution $$\barr{l} A\ijk=\frax13(R\ijk+g^{ih}R_{jhk}),\ E\ijk=-P^i_{kj},\\
D\ijk=C\ijk,\ B\ijk=\frax12(X\ijk+g^{ih}g_{\ell j}X^\ell_{hk}), \earr\0(3.16)$$ for
some $X$ which satisfies $$X\ijk=X^i_{kj},\ g_{sj}X^s_{hk}=g_{sk}X^s_{hj},\mbox{
otherwise arbitrary}.\0(3.17)$$ Indeed, the first eq. in (3.14) is verified by virtue
of (2.5). The others are clearly verified. In (3.15), the first, the second and the
fourth equations become identities by virtue of (3.16). The third is equivalent to
$P^h_{ik}g_{hj}=P^h_{jk}g_{hi}.$ Inserting $P^h_{ik},$ after some calculation we find
that this equation is equivalent to $g_{ij{\scriptstyle||k}}=g_{ij{\scriptstyle||k}}$,
which by (2.5) is an identity. Notice that the condition (2.5), that is $d\ooo_L=0$ is
essential in the solving of (3.14) and of (3.15) as well.

The connection $\nabla$ is determined by $\bld cD$, the torsion $R\ijk$ and $P\ijk$ as
well as by the unknown $d$-tensor field $X\ijk$ satisfying (3.17). We may take
$X\ijk=C\ijk$ since $C_{ijk}=g_{is}C^s_{jk}$ is a completely symmetric $d$-tensor.
This choice singles out a symplectic connection $\bld s\nabla$ that depends on $L$
only. The theorem is proved.

We remark that the choice which we have made is not unique. Thus$\bld s\nabla$ is not
canonical in any way. However we shall treat only it in the following. And we denote
it simply by $\nabla.$

\section{Symplectic curvature tensor field of the symplectic
connection $\nabla$}

In\cite{3}, I. Vaisman established the decomposition of the space of tensors which
have the symmetries of the curvature of a symplectic connection into ${\rm
Sp}(n)$--irreducible components. Accordingly, he discovered three classes of
symplectic connections: flat, Ricci flat and with reducible curvature. A natural
question is to which class our symplectic connection $\nabla$ belongs. For obtaining
an answer we have to compute the symplectic curvature tensor of $\nabla$. We shall do
this in the adapted bases $(\de_i,\dot\pp_i).$

Let $^\nabla R$ and $^DR$ be the curvature tensor of type $(1,3)$ of $\nabla$ and $D$,
respectively. Using $\nabla=\bld cD+A$ we get $$\barr{l} ^\nabla R(X,Y)Z={}^DR(X,Y)Z+
(D_XA)(Y,Z)-(D_YA)(X,Z)+\\ +A(T(X,Y),Z)+A(X,A(Y,Z))-A(Y,A(X,Z)),\ \hfill X,Y,Z\in \chi
(TM). \earr\0(4.1)$$ where $T$ is the torsion of $D$ locally given by (3.5). With the
notations (3.9), the local components of $(D_XA)(Y,Z)$ are given by $$\barr{ll}
(D_{\de_k}A)(\de_j,\de_h)=A^i_{hj\scri k}\de_i,
&(D_{\de_k}A)(\de_j,\dot\pp_h)=D^i_{hj\scri k}\de_i,\\
(D_{\de_k}A)(\dot\pp_j,\de_h)=E^i_{hj\scri k}\dot\pp_i,
&(D_{\de_k}A)(\dot\pp_j,\dot\pp_h)=B^i_{hj\scri k}\de_i,\\
(D_{\dot\pp_k}A)(\de_j,\dot\pp_h)=A^i_{hj}|_k\dot\pp_i,
&(D_{\dot\pp_k}A)(\de_j,\dot\pp_h)=D^i_{hj\scri k}\de_i,\\
(D_{\dot\pp_k}A)(\dot\pp_j,\de_h)=E^i_{hj}|_k\dot\pp_i,
&(D_{\dot\pp_k}A)(\dot\pp_j,\dot\pp_h)=B^i_{hj}|_k\de_i, \earr\0(4.2)$$ where $(\scri
k)$ and $(|k)$ denote the $h$- and $v$-covariant derivatives with respect to $D$. The
curvature operator $^DR(X,Y)$ carries horizontal vector fields to horizontals and the
vertical vector fields to verticals. Its action on horizontals is as follows
$$\barr{l} ^DR(\de_k,\de_j)\de_h=R_h{}^i{}_{jk}\de_i,\\
^DR(\dot\pp_k,\de_j)\de_h=P_h{}^i{}_{jk}\de_i,\\
^DR(\dot\pp_k,\dot\pp_j)\de_h=S_h{}^i{}_{jk}\de_i, \earr\0(4.3)$$ and its action on
verticals is similarly determined by the same $d$-tensors
$R_h{}^i{}_{jk}$, $P_h{}^i{}_{jk}$, $S_h{}^i{}_{jk}$, given by $$\barr{l}
R_h{}^i{}_{jk}=\de_kF^i_{hj}+F^s_{hj}F^i_{sk}- (j/k)+C^{}_{hs}R^s_{jk},\\
P_h{}^i{}_{jk}=\dot\pp_kF^i_{hj}- C^i_{hk\scri j}+C^i_{hs}P^s_{jk},\\
S_h{}^i{}_{jk}=\dot\pp_kC^i_{hj}+C^s_{hj}C^i_{sk}-(j/k),\earr\0(4.4)$$ where $(j/k)$
means the preceding terms with $k$ changed to $j$ and $j$ changed to $k$.

The curvature operator $^\nabla R(X,Y)$ does not preserve
the horizontal and vertical distributions. As such $^\nabla
R$ has twelve components. We put
$$\barr{l}
^\nabla R(\de_k,\de_j)\de_h=
{}^\nabla R_h{}^i{}_{jk}\de_i+K_h{}^i{}_{jk}\dot\pp_i,\\
^\nabla R(\dot\pp_k,\de_j)\de_h=
{}^\nabla P_h{}^i{}_{jk}\de_i+K_h{}^i{}_{jk}\dot\pp_i,\\
^\nabla R(\dot\pp_k,\dot\pp_j)\de_h=
{}^\nabla S_h{}^i{}_{jk}\de_i+M_h{}^i{}_{jk}\dot\pp_i,\\
^\nabla R(\de_k,\de_j)\dot\pp_h=
\widetilde K_h{}^i{}_{jk}\de_i+\widetilde R_h{}^i{}_{jk}\dot\pp_i,\\
^\nabla R(\dot\pp_k,\de_j)\dot\pp_h=
\widetilde L_h{}^i{}_{jk}\de_i+\widetilde P_h{}^i{}_{jk}\dot\pp_i,\\
^\nabla R(\dot\pp_k,\dot\pp_j)\dot\pp_h=
\widetilde M_h{}^i{}_{jk}\de_i+\widetilde S_h{}^i{}_{jk}\dot\pp_i.
\earr\0(4.5)$$
Using (3.5) and (4.2) an explicit form of these components
is obtained as follows
$$\barr{l}
^\nabla R_h{}^i{}_{jk}=R_h{}^i{}_{jk}+
(A^s_{hj}D^i_{sk}-(k/j)),\\
K_h{}^i{}_{jk}=A^i_{hj\scri k}-
(k/j)+R^s_{jk}E^i_{hs}.
\earr\0(4.6)$$
$$\barr{l}
^\nabla P_h{}^i{}_{jk}=P_h{}^i{}_{jk}-
E^s_{hk}D^i_{sj}+A^s_{hj}B^i_{sk},\\
L_h{}^i{}_{jk}=A^i_{hj\scri k}-E^i_{hk\scri j}+
C^s_{jk}A^i_{hs}+P^s_{jk}E^i_{hs}.
\earr\0(4.7)$$
$$\barr{l}
^\nabla S_h{}^i{}_{jk}=S_h{}^i{}_{jk}+
(E^s_{hj}B^i_{sk}-(k/j)),\\
M_h{}^i{}_{jk}=E^i_{hj}\bigi k-(k/j).
\earr\0(4.8)$$
$$\barr{l}
\widetilde K_h{}^i{}_{jk}=D^i{hj\scri k}-(j/k)+
R^s_{jk}B^i_{hs},\\
\widetilde R_h{}^i{}_{jk}=R_h{}^i{}_{jk}+(D^s_{hj}A^i_{sk}-
(k/j)).
\earr\0(4.9)$$
$$\barr{l}
\widetilde L_h{}^i{}_{jk}=D^i{hj}\bigi_k+
C^s_{jk}D^i_{hs}+P^s_{jk}B^i_{hs},\\
\widetilde P_h{}^i{}_{jk}=P_h{}^i{}_{jk}+D^s_{hj}E^i_{sk}-
B^s_{hk}A^i_{sj}.
\earr\0(4.10)$$
$$\barr{l}
\widetilde M_h{}^i{}_{jk}=B^i{hj}\bigi k-(j/k)\\
\widetilde S_h{}^i{}_{jk}=S_h{}^i{}_{jk}+(B^s_{hj}E^i_{sk}-
(k/j)).
\earr\0(4.11)$$

Then, if we take in (4.6)--(4.11),
$A^i_{jk}=\frax13(R\ijk+g^{ih}R_{jhk}),$
$E\ijk=-P\ijk$ and $D\ijk=B\ijk=C\ijk,$ we obtain the
twelve components of the curvature tensor of type $(1,3)$
of $\nabla.$

The symplectic curvature tensor is defined by $$S(X_2,X_2,X_3,X_4)=\ooo_L ({}^\nabla
R(X_3,X_4)X_2,X_1),\ X_1,...X_4\in\chi(TM).\0(4.12)$$ This is skew symmetric in the
last two arguments and symmetric with respect to the first two arguments. Moreover,
the cyclic sum over the last three arguments vanishes, \cite{3}. It is locally
determined by the twelve tensors $^\nabla R_{hijk},...,\widetilde S_{hijk}$ from
(4.6)--(4.11) with the upper index brought down with $(g^{hs})$ on the second place.

The Ricci curvature tensor of $\nabla$ is defined by the
usual formula
$$\sigma(X,Y)={\rm Tr}(V\to{}^\nabla R(V,X)Y).\0(4.13)$$
A direct calculation gives
$$\barr{ll}
\sigma(\de_j,\de_h)={}^\nabla
R_h{}^i{}_{jki}+L_h{}^i{}_{ji},\\
\sigma(\de_j,\dot\pp_h)=\widetilde K_h{}^i{}_{ji}+
\widetilde P_h{}^i{}_{ji},\\
\sigma(\dot\pp_j,\dot\pp_h)=-\widetilde L_k{}^i{}_{ji}+
\widetilde S_h{}^i{}_{ji}.
\earr\0(4.14)$$

The components of curvature tensors just found are quite complicated. Thus it is quite
sure that for a general Lagrangian $L$, the symplectic connection $\nabla$ is a
general one, too. The above formula simplify for the Lagrangian defined by a
Riemannian metric $\g\ij(x)$ as in (3.6).

Let $\g\ijk(x)$ be the Christoffel symbols of $\g\ij(x)$
and $r_h{}^i{}_{jk}(x)$ its curvature tensor. Then
$\ooo_L=\g\ij(x)\de y^i\w dx^j$ and $\nabla$ takes the form
$$\barr{l}
\nabla_{\de_k}\de_j=\g\ijk\de_i+A\ijk\dot\pp_i,\
\nabla_{\dot\pp_k}\de_j=0,\\
\nabla_{\de_k}\dot\pp_j=\g\ijk\dot\pp_i,\
\nabla_{\dot\pp_k}\dot\pp_j=0,\earr\0(4.15)$$
where $A\ijk$ is given by (3.12) with
$R\ijk=r_h{}^i{}_{jk}(x)y^h.$ An inspection of (4.4),
(4.6)--(4.11) shows that the nonzero components of $^\nabla
R$ are the following ones
$$^\nabla R_h{}^i{}_{jk}=r_j{}^i{}_{hk}(x),\0(4.16)$$
$$K_h{}^i{}_{jk}=A^i_{hj\scri k}-(k/j)=
\frax13(r_q{}^i{}_{hj;k}+r_j{}^i{}_{hq;k})y^q-(k/j),\0(4.17)$$
$$L_h{}^i{}_{jk}=A^i_{hj}|_k=
\frax13(r_k{}^i{}_{hj}+r_j{}^i{}_{hk}),\0(4.18)$$
where $(;k)$ means the covariant derivative with respect to
the Levi--Civita connection of $(\g\ij).$
From (4.16)--(4.18) it follows

\bk\n{\bf Theorem 4.1.} {\it Let be the symplectic manifold $(TM,\ooo_L)$ for
$L(x,y)=\g\ij(x)y^iy^j$ and $(\g\ij(x))$ a Riemannian metric. The symplectic
connection $\nabla$ is flat \io $(\g\ij)$ is flat.}\bk

From (4.14) it results that the only non--zero component of
the Ricci curvature tensor of $\nabla$ is
$$\sigma(\de_j,\de_h)=\frax23 r_{jh}(x),\0(4.19)$$
where $r_{jh}(x)=r_j{}^i{}_{hi}(x)$ is the Ricci tensor of
$\g\ij(x).$ Thus we have

\bk\n{\bf Theorem 4.2.} {\it The same hypothesis as in Theorem {\rm4.1}. The
symplectic connection $\nabla$ is Ricci flat \io $\g\ij(x)$ is Ricci flat.}\bk

\begin{flushright}
{\footnotesize{June, 2001}}
\end{flushright}



\newpage

\def\llr{\Longleftrightarrow }
\def\frax#1#2{\displaystyle\frac{#1}{#2}}
\def\derp#1#2{\displaystyle\frac{\partial#1}{\partial#2}}
\def\derpp#1{\displaystyle\frac{\partial}{\partial#1}}
\def\bld#1#2{{\buildrel{#1}\over{#2}}}
\def\scri{\scriptstyle{|}}
\def\a{{\alpha}}
\def\b{{\beta}}
\def\de{{\delta}}
\def\G{{\Gamma}}
\def\g{{\gamma}}
\def\pp{{\partial}}
\def\wt{\widetilde}
\def\dd{\displaystyle}
\def\vsp{\vspace*{1,5mm}\\ }
\def\n{\noindent}

\runningauthor={M. ANASTASIEI}
\runningtitle={METRIZABLE LINEAR CONNECTIONS IN VECTOR BUNDLES}
\noindent
\baselineskip 8pt
\noindent{\footnotesize{Publ. Math. Debrecen}}
\hfill\break
{\footnotesize{62, no. 3-4, 277--287, 2003}}
\vskip 2cm
\baselineskip 11.5pt plus .15pt
\centerline{\bf\Large METRIZABLE LINEAR CONNECTIONS}
\vskip .2cm
\centerline{\bf\Large IN VECTOR BUNDLES}
\vskip .5cm
\centerline{\bf {\footnotesize{BY}}}
\vskip .5cm
\centerline{\bf {\footnotesize{M. ANASTASIEI}}}
\vskip 1cm

\begin{center}
\footnotesize{Dedicated to Professor Dr. Lajos Tam\'{a}ssy at his 80th anniversary}
\end{center}

\begin{abstract} A linear connection $\nabla$ in a vector bundle is said to be
metrizable if the vector bundle admits a Riemannian metric $h$ with the
property $\nabla h=0.$ Sufficient conditions for the linear
connection $\nabla$ to be metrizable are provided.

Mathematics Subject Classification: Primary 53C60; Secondary 53C05.\\
Key words and phrases: vector bundles, linear connections, metrizability
\end{abstract}

\setcounter{section}{0}
\section*{Introduction}
The problem of the metrizability of a linear connection was
treated by many authors in various contexts (see the paper \cite{7}
by L. Tamassy and the references therein). When a linear
connection $\nabla$ in a vector bundle $ \xi = (E,p,M) $ is metrizable, its parallel translations are
isometries with respect to any Riemannian metric $h$ in $\xi $ with $\nabla
h=0$. Using a
local chart around a point $x$ in $M$, the holonomy group  $\phi(x)$ may be identifed with a subgroup
of $GL(m,\R)$, where $m$ is the dimension of fibre. With this
identification, a necessary condition for $\nabla$ be metrizable
is that the holonomy group to be contained in the orthogonal group $O(m).$ We prove two
versions of the converse of this fact (Theorems 3.1 and 3.2).
Then, we are dealing with the same problem when the vector bundle $\xi $ is
endowed with a Finsler function. The linear connection $\nabla$
induces a nonlinear connection on $E$ and a linear connection $D$
in the vertical vector bundle over $E$. The Finsler function $F$
defines a Riemannian metric $g$ in the vertical vector bundle over $E$. We show
that if $g$ is covariant constant on horizontal directions, then
$\nabla$ is metrizable (Theorem 4.2). When the tangent bundle of a
manifold $M$ is endowed with a Finsler function $F$ one says that
$(M,F)$ is a Finsler manifold. In this case our result has to be compared with the one due to Z. Szab\'o, ( \cite{6}) regarding the
metrizability of the Berwald connection.

If the cotangent bundle of a manifold $M$ is endowed with a
Finsler function $K$, then the pair $(M,K)$ is called a {\it
Cartan space}. This notion was introduced and studied by R. Miron
in \cite{3}. In this case Theorem 4.1 has to be compared with our previous results
on the metrizability of Berwald--Cartan connection \cite{1}.

The first two sections of the paper are devoted to some
preliminaries from the theory of vector bundles and linear
connections in vector bundles.

\section{Vector bundles}

Let $\xi=(E,p,M)$ be a vector bundle of rank $m$. Here $E$ and $M$
are smooth i.e. $C^{\infty}$ manifolds with ${\rm dim}\,M=n$, ${\rm
dim}\,E=n+m$ and $p:E\to M$ is a smooth submersion. The fibres
$E_x=p^{-1}(x),$ $x\in M$ are linear spaces of dimension $m$ which
are isomorphic with the type fibre $\R^m$.

Let $\{(U_\a,\psi_\a)\}_{\a\in A}$ be an atlas on $M$. A vector
bundle atlas is $\{(U_\a,\varphi_\a$, $\R^m)\}_{\a\in A}$ with the
bijections $\varphi_\a:p^{-1}(U_\a)\to U_\a\times\R^m$ in the form
$\varphi_{\alpha}=(p(u), \varphi_{\alpha , p(u)} (u))$, where $\varphi_{\alpha , p(u)} : E_p(u) \to \R^m $ is a bijection.
The given atlas on $M$ and a vector bundle atlas provide an atlas $(p^{-1}(U_\a),\phi_\a)_{\a\in
A}$ on $E$. Here $\phi_\a:p^{-1}(U_\a)\to\varphi_\a(U_\a)\times\R^m$ is
the bijection given by $\phi_\a(u)=(\psi_\a(p(u)),\varphi_{\a,p(u)}(u)).$ For $x\in M,$ we
put $\psi_\a(x)=(x^i)\in\R^m$ and we take $(x^i,y^a)$ as local
coordinates on $E$. If $(U_\b,\psi_\b)$ is such that $x\in
U_\a\cap U_\b\ne\emptyset$ and $\psi_\b(x)=(\wt x^i),$ then
$\psi_\b\circ\psi^{-1}_\a$ has the form
$$\wt x^i=\wt
x^i(x^1,...,x^n),\ {\rm rank}\left(\frax{\pp\wt x^i}{\pp
x^j}\right)=n.\leqno(1.1)$$

Let $(e_a)$ be the canonical basis of $\R^m$.
Then $\varphi^{-1}_{\a,x}(e_a)=\varepsilon_a(x)$ is a basis of $E_x$ and $u\in
E_x$ takes the form $u=y^a\varepsilon_a(x).$ We put $\wt y^a=M^a_b(x)y^b$
with ${\rm rank}(M^a_b(x))=m.$ Then $\phi_\b\circ\phi^{-1}_\a$ has
the form
$$\begin{array}{ll}
\wt x^i=\wt x^i(x^1,...,x^n),
&{\rm rank}\left(\frax{\pp\wt x^i}{\pp x^j}\right)=n\vsp
\wt y^a=M^a_b(x)y^b,
&{\rm rank}(M^a_b(x))=m.\end{array}\leqno(1.2)$$
The indices $i,j,k,...,a,b,c,..$ will take the values $1,2,...,n$
and $1,2,...,m,$ respectively. The Einstein convention on
summation will be used.

We denote by ${\cal F}(M),{\cal F}(E)$ the ring of real functions on $M$
and $E$, respectively and by ${\cal X}(M),$ resp. $\Gamma(E)$,
${\cal X}(E)$ the module of sections of the tangent bundle of $M$,
resp. of the bundle $\xi$ and of the tangent bundle of $E$. On
$U_\a$, the vector fields $\left(\pp_k:=\dfrac{\pp}{\pp x^k}\right)$ provide a
local basis for ${\cal X}(U_\a).$ The sections $\varepsilon_a:U_\a\to
p^{-1}(U_\a)$ given by $\varepsilon_a(x)=\varphi^{-1}_{\a,x}(e_a)$ will be taken as
canonical basis for $\Gamma(p^{-1}(U_\a))$ and a section $A:U_\a\to
p^{-1}(U_\a)$ will take the form $A(x)=A^a(x)\varepsilon_a(x).$

Let $\xi^*=(E^*,p^*,M)$ be the dual of the vector bundle $\xi.$ We
take as local basis of $\Gamma(E^*)$ on $U_\a$, the sections
$\theta^a:U_\a\to p^{*-1}(U_\a),$ $x\to\theta^a(x)\in E^*_x$ such
that $\theta^a(\varepsilon_b(x))=\delta^a_b$.

Next, we may consider the tensor bundle of type
$(r,s){\cal T}^r_s(E):=E\underbrace{\otimes\cdots\otimes}_r E$ $\otimes
E^*\underbrace{\otimes\cdots\otimes}_sE^*$ over $M$ and its
sections. For $g\in\Gamma(E^*\otimes E^*)$ we have the local
representation $g=g_{ab}(x)\theta^a\otimes\theta^b$.
As $E^*\otimes E^*\cong L_2(E,\R),$ we may regard $g$ as a smooth
mapping $x\to g(x):E_x\times E_x\to\R$ with $g(x)$ a bilinear
mapping given by $g(x)(s_a,s_b)=g_{ab}(x).$

If the mapping $g(x)$ is symmetric i.e. $g_{ab}=g_{ba}$ and
positive--definite i.e. $g_{ab}(x)\zeta^a\zeta^b>0$ for every
$0\ne(\zeta^a)\in\R^m,$ one says that $g$ defines a Riemannian metric
in the vector bundle $\xi.$

The sets of sections $\G(T^r_s(E))$ are ${\cal F}(M)$-modules for
every natural numbers $r,s$. On the sum
$\dd\bigoplus_{r,s}\G(T^r_s(E))$ a tensor product can be defined
and one gets a tensorial algebra ${\cal T}(E)$. For the vector bundle
$(TM,\tau,M)$ this reduces to tensorial algebra of the manifold
$M$.

\section{Linear connections in a vector bundle}

\bigskip\n{\bf Definition 2.1.} A linear connection in the vector
bundle $\xi=(E,p,M)$ is a mapping $\nabla:{\cal X}(M)\times\G(E)\to\G(E)$,
$(X,A)\to\nabla_XA$ which is ${\cal F}(M)$-linear in the first
argument, additive in the second and
$$\nabla_X(fA)=X(f)A+f\nabla_XA,\ f\in{\cal F}(M).\leqno(2.1)$$
For $X=X^k(x)\pp_k$ and $A=A^a(x)\varepsilon_a(x),$ we get
$$\nabla_XA=X^k(\pp_kA^a+\G^a_{bk}(x)A^b)\varepsilon_a(x),\leqno(2.2)$$
where the local coefficients $\G^a_{bk}(x)$ are defined by
$$\nabla_{\pp_k}\varepsilon_b=\G^a_{bk}\varepsilon_a.\leqno(2.3)$$
If $\wt\G^c_{dj}$ are the local coefficients of $\nabla$ on $U_\b$
such that $U_\a\cap U_\b\ne\emptyset$, then we have
$$\wt\G^c_{dj}(\wt x(x))=M^c_a(x)(M^{-1})^b_d\derp{x^k}{\wt
x^j}\G^a_{bk}(x)- \derp{M^c_b}{x^k}\ \derp{x^k}{\wt x^j}(M^{-1})^b_d.\leqno(2.4)$$

A section $A$ of $\xi$ is called {\it parallel} if $\nabla_XA=0$
for every $X\in{\cal X}(M).$

The linear connection $\nabla$ induces operators of covariant
derivative $\nabla_k$ in the tensorial algebra ${\cal T}(E)$ taking
$\nabla_kf=\pp_kf$, $\nabla_k\b_a=\pp_k\b_a-\G^c_{ak}\b_c$ and
requiring that $\nabla_k$ to satisfy the Newton--Leibniz rule
with respect to the tensor product and to commute with the all
contractions.

Let $c:[0,1]\to M$ be a curve on $M$ and $A:t\to A(t):=A(c(t))$ a
section of $\xi$ along the curve $c$.  Then $\nabla_{\dot
c(t)}A=:\frax{\nabla A}{dt}$ is called the covariant derivative of
$A$ along $c$.

On $U_\a\cap c[0,1]$ if one puts $c(t)=(x^i(t)),$ we get
$$\frax{\nabla A}{dt}=\left(\frax{dA^a}{dt}+\G^a_{bk}(x(t))A^b\frax{dx^k}{dt}\right)\varepsilon_a. \leqno(2.5)$$
The section $t\to A(t)$ is said to be {\it parallel} on $c$ if
$\frax{\nabla A}{dt}=0.$ This means that the functions $(A^a(t))$
have to be solutions of the following system of ordinary linear
differential equations
$$\frax{dA^a}{dt}+\G^a_{bk}(x)A^b\frax{dx^k}{dt}=0.\leqno(2.6)$$
For given initial conditions $A^a(0)=(u^a)\in E_{c(0)}$ the system
(2.6) admits a unique solution that can be prolonged beyond $U_\a$
providing a parallel section $A$ along $c$. If we associates to
$(u^a)=A^a(0)$ the element $(v^a)=A^a(1)\in E_{c(1)}$ one gets a
linear isomorphism $P_c:E_{c(0)}\to E_{c(1)},$ called the {\it
parallel translation} of $E_{c(0)}$ to $E_{c(1)}$ along $c$. The
parallel translations can be defined along any curve or segment of
curve providing linear isomorphisms between fibres in various
point of curves on $M$. In particular, if one considers the loops
with the origin in $x\in M$, the corresponding parallel
translations as linear isomorphisms $E_x\to E_x$ can be composed
and a group $\phi(x)$ called the holonomy group in $x\in M$ is
obtained.

When $M$ is connected, the holonomy groups $\phi(x),$ $x\in M$,
are isomorphic and one speaks about the holonomy group $\phi$
associated to or defined by $\nabla$.

The covariant derivative along $c$ can be recovered from parallel
translations according to the following known

\bigskip\n{\bf Lemma 2.1.} {\it Let $A$ be a section of $\xi$ along a
curve on $M$, $c$: $t\to c(t)$, $t\in\R$, starting from $x=c(0)$. Then
$$(\nabla_{\dot c(0)}A)(x)=\lim_{t\to0}\frax1t(P_c(A(t))-A(0)),\leqno(2.7)$$
where $P_c:E_{c(t)}\to E_x$ is the parallel translation along $c$.}\medskip

\section{A sufficient condition for $\nabla$ be metrizable}

Let $\nabla$ be a linear connection in the vector bundle $\xi=(E,p,M)$. Assume that the manifold $M$ is connected. One says that $\nabla$ is {\it metrizable} if there exists a Riemannian metric $g$ in $\xi$ such that $\nabla g=0.$ When $\nabla$ is metrizable all parallel translations $P_c:(E_x,g_x)\to(E_y,g_y)$ for any curve $c$ and any points $x,y$ joining them in $M$ are isometries. In particular, the holonomy group $\phi(x)$ is a subgroup of the orthogonal group of $(E_x,g_x).$ These facts follow from

\bigskip\n{\bf Lemma 3.1.} {\it Let $g$ be any Riemannian metric in the vector bundle $\xi$ and $c:t\to c(t),$ $t\in\R$, a curve in $M$ with $c(0)=x.$ Then
$$\left(\nabla_{\dot c(0)}g\right)(A,B)=\dd\lim_{t\to0}\frax1t(g_{c(t)}(P_cA,P_cB)-g_x(A,B)),\leqno(3.1)$$
where $A,B\in E_x$ and $P_c:E_x\to E_{c(t)}$ is the parallel translation along $c$.}\medskip

\n{\it Proof.} Let $\wt A,\wt B$ be sections of $\xi$ which are parallel on $c$ such that $\wt A(0)=A,$ $\wt B(0)=B$. Then $P_cA=\wt A(t)$ and $P_c(B)=\wt B(t).$ By the Taylor theorem and using the condition that $\wt A$ and $\wt B$ are parallel sections on $c$, in the natural basis $(\varepsilon_a)$ we get $(P_cA)^a=\wt A^a(t)=A^a+\frax{d\wt A}{dt}(\tau)t=A^a-\G^a_{ck}(x(\tau))\wt A^c(\tau)\frax{dx^k}{dt}t$ and a similar formula for $(P_cB)^b,$ $a,b=1,2,...,m.$ Then, using again the Taylor theorem, omitting the terms which contain $t^2$, we may write:
$$\begin{array}{c}
g_{ab}(t)(P_cA)^a(P_cB)^b-g_{ab}(x)A^aB^b=\\ \\
\left(g_{ab}(x)+\frax{dg_{ab}}{dt}(\theta)t\right)(P_cA)^a(P_cB)^b-\vsp
-g_{ab}(x)A^aB^b=
\left(\frax{dg_{ab}}{dt}-g_{ac}\G^c_{bk}\frax{dx^k}{dt}-g_{cb}\G^c_{ak}\frax{dx^k}{dt}\right)t,\end{array}\leqno(3.2)$$
where the terms in the last paranthesis are computed for $\tau,\tau',\theta\in(0,t).$

Dividing in (3.2) by $t$ and taking $t\to0$, one obtains (3.1).

By Lemma 3.1 we have also that if all parallel translations of $\nabla$ are isometries with respect to $g$, then $\nabla g=0$. Thus, in order to prove that $\nabla$ is metrizable we need to find a Riemannian metric $g$ such that all parallel translations of $\nabla$ to be isometries with respect to $g$. Taking an arbitrary bundle chart $(U_\a,\varphi_\a,\R^m)$, using the linear isomorphism $\varphi_{\a,x}:E_x\to\R^m,$ we may identify $\phi (x)$, $x\in U_\a,$ with a subgroup of $GL(\R^m)$. When $\nabla$ is metrizable, by Lemma 3.1 it follows that this subgroup is contained in the orthogonal group $O(m).$ Therefore, a necessary condition for $\nabla$ be metrizable is that its holonomy group to be contained in $O(m).$ We show two versions of the converse.

\bigskip\n{\bf Theorem 3.1.} {\it Let $\nabla$ be a linear connection in the vector bundle $\xi=(E,p,M)$ with $M$ connected. Assume that there exists a point $x_0\in M$ such that the holonomy group $\phi(x_0)$ is contained in the orthogonal group of $E_{x_0}$ when $E_{x_0}$ is regarded as being isomorphic with the Euclidean space $(R^m,<,>)$ via a fixed bundle chart. Then $\nabla$ is metrizable.}\medskip

\n{\it Proof.} Let $h_0$ be the inner product on $E_{x_0}$ induced by $<,>$ via the bundle chart $(U_\a,\varphi_\a,\R^m)$, $x_0\in U_\a,$ that is,
$$h_0(u,v)=\,<\varphi_{\a,x_0}u,\varphi_{\a,x_0}v>.\leqno(*)$$
By hypothesis this inner product is invariant under the group $\phi(x_0).$ Let be any $x\in M.$ We join $x$ with $x_0$ using a curve $c:[0,1]\to M,$ $c(0)=x,$ $c(1)=x_0,$ consider the parallel translation $P_c:E_x\to E_{x_0}$ and define an inner product $h_x$ in $E_x$ by
$$h_x(A,B)=h_0(P_cA,P_cB),\ A,B\in E_x.\leqno(3.3)$$

\bigskip\n{\bf Lemma 3.2.} {\it The inner product $h_x$ does not depend
on the curve $c$.}

Indeed, if $\wt c$ is another curve joining $x$ with $x_0$, we
consider the reverse $c_-$ of $c$ and the loop $\wt c\circ c_-$ in
$x_0$. It follows that $h_0\left(P_{\wt c\circ c_-}u,P_{\wt c\circ
c_-}v\right)=h_0(u,v)$, $u,v\in E_{x_0}$. Inserting here $u=P_cA$ and
$v=P_cB$ and taking into account (3.3), the Lemma follows.

The mapping $x\to h_x$ is smooth since $P_c$ smoothly depends  on $x$ according to the general theory of differential equations. Thus we obtain a Riemannian metric $h$ in $\xi$. The parallel translations of $\nabla$ are isometrics with respect to $h$. Indeed, for $y$ a point of $M$ different from $x$, any parallel translation from $E_x$ to $E_y$ has the form $P_{\sigma_-\circ c}=P_{\sigma_-}\circ P_c$ for $\sigma_-$ the reverse of a curve $\sigma$ joining $y$ with $x_0.$ This is an isometry as a product of isometries. Therefore, we may conclude using Lemma 3.1, that $\nabla h=0.$ q.e.d.

The following version of the Theorem 3.1 extends to the vector bundle setting a result of B.G. Schmidt \cite{5}.

\bigskip\n{\bf Theorem 3.2.} {\it Let $\nabla$ be a linear connection
in the vector bundle $\xi=(E,p,M)$ with $M$ connected. Assume that
for a fixed $x_0\in M$, the holonomy group $\phi(x_0)$ leaves
invariant a given positive--definite quadratic form $h_0$ on
$E_{x_0}$. Then there exists a Riemannian metric $h$ in $\xi$ such
that $\nabla h=0.$}\medskip

\n{\it Proof.} Let denote by the same letter $h_0$ the inner product in
$E_{x_0}$ defined by the quadratic form $h_0$. This inner product
could be obtained by transferring one from $\R^m$ using a bundle
chart. By hypothesis the inner product $h_0$ is invariant under
$\phi(x_0)$. From now on the reasoning proving Theorem 3.1 can be
entirely repeated in order to find $h$ such that $\nabla h=0.$

\bigskip\n{\bf Remark 3.1.}  The Riemannian metric $h$ found in Theorem 3.1
is not unique and is not canonical in any way. The same applies
for $h$ found in Theorem 3.2.

\section{Another condition for $\nabla$ be metrizable}

We are dealing with the problem of the metrizability of a linear
connection $\nabla$ in a vector bundle endowed with a Finsler
function.

\bigskip\n{\bf Definition 4.1.} Let $\xi=(E,p,M)$ be a vector bundle of rank $m$. A {\it Finsler function} on $E$ is a nonnegative real function $F$ on $E$ with the properties
\begin{itemize}
\item[1)] $F$ is smooth on $E\setminus\{(x,0),x\in M\},$
\item[2)] $F(x,\lambda y)=\lambda F(x,y)$ for all $\lambda>0,$
\item[3)] The matrix with the entries $g_{ab}(x,y)=\frax12\ \derp{^2F^2}{y^a\pp y^b}$ is positive definite.
\end{itemize}
On the manifold $E$ we have the vertical distribution $u\to V_uE=\ker p_{*,u}$ where $p_*$ denotes the differential of $p$. This is spanned by $\dot\pp_a:=\derpp{y^a}$. A distribution $u\to H_uE$ which is supplementary to the vertical distribution is called a {\it horizontal} distribution or a {\it nonlinear connection} on $E$. This is usually taken as spanned by $\delta_i=\pp_i-N^a_i(x,y)\dot\pp_a,$ where the functions $(N^a_i(x,y))$ are called the {\it coefficients} of the given nonlinear connection. Under a change of coordinates they behave as follows:
$$\wt N{}^a_j\ \derp{\wt x^j}{x^k}=M^a_b(x)N^b_k(x,y)-\derp{M^a_b}{x^k}y^b,\leqno(4.1)$$
a fact which is equivalent with
$$\delta_i=\derp{\wt x^k}{x^i}\wt\delta_k.\leqno(4.1)'$$
Introducing the horizontal distribution we have
$$T_uE=H_uE\oplus V_uE,\ u\in E.\leqno(4.2)$$
It is convenient to decompose the geometrical objects on $E$ according to (4.2) using the adapted basis $(\delta_i,\dot\pp_a)$ and its dual $(dx^i,\delta y^a=dy^a+N^a_i(x,y)dx^i).$

The linear connection $\nabla$ in $\xi$ defines a nonlinear connection on $E$ taking $N^a_i(x,y)=\G^a_{bi}(x)y^b.$ Indeed, using (2.4) it is easy to check that these functions satisfy (4.1). From now on we shall use only the decomposition (4.2) provided by these functions.

Furthermore, the linear connection $\nabla$ induces a linear connection $D$ in the vertical bundle over $E$ as follows: $D:{\cal X}(E)\times\G(VE)\to \G(VE),$ $(X,Z)\to D_XZ$ is given for $Z=Z^a\dot\pp_a$ by
$$D_{\delta_k}\dot\pp_a=\G^a_{bk}(x)\dot\pp_a,\ D_{\dot\pp_b}\dot\pp_a=0.\leqno(4.3)$$
We call $D$ the vertical lift of $\nabla$ and we use $D_{\delta_k}$ for defining a {\it horizontal} covariant derivative operator in the tensor algebra of the vertical bundle, denoted by $\scri k$, setting
$$\begin{array}{l}
f_{\scri k}=\delta_kf\mbox{ for any function on  }E\vsp
X^a_{\scri k}=\delta _kX^a+\G^a_{bk}(x)X^b.\end{array}\leqno(4.4)$$

For a fixed $x\in E$, the pair $(E_x,F_x)$ is a Minkowski space. Here $F_x$ denotes the restriction of $F$ to $E_x$ and it is obvious that this is a Minkowski norm on  $E_x.$

Now we show that under certain conditions the parallel translations of $\nabla$ are isometries of Minkowski spaces.

\bigskip\n{\bf Theorem 4.1.} {\it Let $\xi=(E,p,M)$ be a vector bundle
of rank $m$ with $M$ connected, endowed with a Finsler function
$F$ and with a linear connection $\nabla$ as well. Let $\scri k$
be the horizontal covariant derivative operator defined by the
vertical lift $D$ of $\nabla$. If $F_{\scri k}=0,$ then the
parallel translation defined by $\nabla$,
$P_c:(E_x,F_x)\to(E_y,F_y)$ is an isometry of Minkowski spaces for
any points $x,y\in M$ and any curve $c:[0,1]\to M$ joining them.}\medskip

\n{\it Proof.} Let be $u\in E_x$ and $t\to A(t)$, $t\in[0,1]$ a section of
$\xi$ which is parallel along $c$ and $A(0)=u.$ Its local
components $A^a$ are solutions of the system of differential
equations (2.6). And $P_c(u)=A(1):=v.$

We know already that $P_c$ is a linear isomorphism. Let us write the condition $F_{\scri k}=0$ for the points $(x(t),A(t))$ of $E$ where $t\to x(t)$ is the local representation of the curve $c$. We obtain:
$$0=\left(\derp F{x^k}-A^b\G^a_{bk}\derp F{y^a}\right)\frax{dx^k}{dt}\bld{(2.6)}{=\!=\!=}
\derp F{x^k}\ \frax{dx^k}{dt}+\derp F{y^a}\ \frax{dA^a}{dt}=\frax{dF(x(t),A(t))}{dt}\cdot$$
Thus the function $F(x(t),A(t))$ is constant. It follows $F(x,u)=F(y,P_cu)$, that is, $F_x(u)=F_y(P_cu).$ In other words, $P_c$ is an isometry of Minkowski spaces $(E_x,F_x)$ and $(E_y,F_y)$. q.e.d.

\bigskip\n{\bf Corollary 4.1.} {\it In the hypothesis of Theorem
{\rm4.1}, the holonomy group $\phi(x)$ consists of isometries of
the Minkowski space $(E_x,F_x)$.}\bigskip

The functions $g_{ab}(x,y)$ define a Riemannian metric in the vertical bundle over $E$ by $g=g_{ab}(x,y)\delta y^a\otimes\delta y^b$. We call $(g_{ab}(x,y))$ the Finsler metric associated with $F$.

The condition $F_{\scri k}=0$ from the hypothesis of Theorem 4.1 can be replaced with
$g_{ab\scri k}=0,$ because of

\bigskip\n{\bf Lemma 4.1.} {\it $F_{\scri k}=0$ is equivalent with
$g_{ab\scri k}=0.$}\medskip

\n{\it Proof.} The homogeneity of $F$ implies $F^2(x,y)=g_{ab}(x,y)y^ay^b.$ Then $F^2_{\scri k}=2FF_{\scri k}=g_{ab\scri k}y^ay^b+2g_{ab}y^a_{\scri k}y^b=g_{ab\scri k}y^ay^b$ since $y^a_{\scri k}=0.$ Thus if $g_{ab\scri k}=0,$ then $F_{\scri k}=0.$ In order to prove the converse, we notice that $\dot\pp_a(H_{\scri k})=(\dot\pp_aH)_{\scri k}$ for any function $H$ on $E$. This follows by a direct calculation taking care that $\dot\pp_aH$ is a vertical $1$-form. Using this ``commutation'' formula we get $g_{ab\scri k}=\frax12\dot\pp_a\dot\pp_b(F^2_{\scri k})=\dot\pp_a\dot\pp_b(FF_{\scri k})=0.$ q.e.d.

Now we are ready to prove the main result of this section.

\bigskip\n{\bf Theorem 4.2.} {\it Let $\nabla$ be a linear connection
in the vector bundle $\xi=(E,p,M)$ with $M$ connected. Suppose
that $E$ is endowed with a Finsler function $F$ with the
associated Finsler metric $g_{ab}(x,y).$ Let $\scri k$ be the
$h$-covariant derivative operator induced by $\nabla$. If
$g_{ab\scri k}=0$, then $\nabla$ is metrizable.}\medskip

\n{\it Proof.} For a fixed $x_0\in M$ we have the Minkowski space
$(E_{x_0},F_{x_0})$. Let $G$ be the group of all linear
isomorphisms of $E_{x_0}$ which preserve the set $S_{x_0}=\{u\in
E_{x_0},F_{x_0}(u)=1\}$. This $G$ is a compact Lie group since $S_{x_0}$
is compact. In our hypothesis, according to Lemma 4.1 and
Corollary 4.1, the holonomy group $\phi(x_0)$ is a Lie subgroup of
$G$. Let $<,>$ be any inner product on $E_{x_0}$. Define a new
inner product on $E_{x_0}$ by
$$h_{x_0}(u,v)=\frax1{{\rm vol}(G)}\int_G<gu,gv>\mu_G,\leqno(4.5)$$
for $u,v\in E_{x_0},\ g\in G$ and $\mu_G$ the bi-invariant Haar
measure on $G$.

It follows that for every $a\in G$ we have
$$h_{x_0}(au,av)=h_{x_0}(u,v),\ u,v\in E_{x_0}.\leqno(4.6)$$
In particular, (4.6) holds for any element of $\phi(x_0)\subset
G$. Thus $\phi(x_0)$ leaves invariant the inner product $h_{x_0}$
in $E_{x_0}.$ The inner product $h_{x_0}$ is extended by parallel
translations to a Riemannian metric $h$ in $\xi$. Furthermore, this
metric verifies $\nabla h=0$ since all parallel translations of
$\nabla$ become isometries with respect to $h$.  Thus $\nabla$ is
metrizable. q.e.d.

\bigskip\n{\bf Remark 4.1.} The Riemannian metric $h$ is not unique and it
is not canonical in any way.


\newpage

\def\bld#1#2{{\buildrel{#1}\over{#2}}}
\def\scri{\scriptstyle{|}}
\def\bigi{\big|}
\def\de{{\delta}}
\def\wh#1{\widehat{#1}}
\def\wt#1{\widetilde{#1}}
\def\derp#1#2{\displaystyle\frac{\partial#1}{\partial#2}}
\def\derpp#1{\displaystyle\frac{\partial}{\partial#1}}
\def\N{\mathbb{N}}

\runningauthor={M. ANASTASIEI}
\runningtitle={GEOMETRY OF BERWALD VECTOR BUNDLES}
\noindent
\baselineskip 8pt
\noindent{\footnotesize{Algebras Groups Geom.\hfill\break 21 (2004), no. 3, 251-262}}
\vskip 2cm
\baselineskip 11.5pt plus .15pt
\centerline{\bf\Large GEOMETRY OF BERWALD}
\vskip .2cm
\centerline{\bf\Large VECTOR BUNDLES}
\vskip .5cm
\centerline{\bf by Mihai ANASTASIEI}
\vskip .5cm
\centerline{\footnotesize{Dedicated to Prof. Dr. Radu MIRON on the occasion of his 75th birthday}}

\begin{abstract}
Let $\xi$ be a vector bundle endowed with a nonlinear connection
$N$. It is called a Berwald vector bundle if the local
coefficients of the Berwald linear connection defined by $N$ do
not depend on the variables $y$ in fibres of $\xi$. Thus they
define a linear connection $\nabla$ in $\xi$. One endows $\xi$
with a regular Lagrangian $L$. A compatibility condition
between $L$ and $N$ is introduced and consequences of it on the
holonomy group of $\nabla$ are derived. Assuming that $L$ is
homogeneous of degree two in $y$, one proves that $\nabla$ is
metrizable. Some particular cases and examples are discussed.\\
{\bf MCS 2000:} 53C07, 53C60.\\
{\bf Kewwords and phrases:} vector bundles, nonlinear connections,
Berwald connections.\end{abstract}

\section*{Introduction}
In the very influential paper \cite{8}, R. Miron develops the
geometry of the total space of a vector bundle using ideas and
techniques from Finsler geometry. He considers on the total space
$E$ of the vector bundle $\xi=(E,p,M)$ a distribution that is
supplementary to the vertical distribution i.e. a nonlinear
connection and decomposes all geometrical objects on $E$ with respect to these
distributions. On this way he
proposes an elegant treatment of the linear connections and of
metrical structures on $E$. From a nonlinear connection $N$ a
linear connection in the vertical bundle over $E$ is easily
derived. This is called the Berwald connection associated to $N$.
When it happens that the local coefficients of this connection do
not depend on the variables $y$ in fibres, they define a linear
connection $\nabla$ in the vector bundle $\xi$ and the pair
$(\xi,N)$ will be called a Berwald bundle. Some properties of the
pairs $(\xi,N)$ are given in Section 1. Then, in Section 2, we
endow $E$ with a regular Lagrangian $L$ and introduce a natural
condition of compatibility between $N$ and $L$. Some direct
consequences of this compatibility are given in Proposition 2.1.
Then we consider the parallel translations defined by $\nabla$ and
we show in Theorem 2.1 that these are compatible with the
structures induced by the Lagrangian on the fibres of $\xi$. In
particular, the holonomy group $\phi(x)$, $x\in M$, of $\nabla$
preserves the indicatrix defined by $L$. The differentials of the
elements of the holonomy group $\phi(x)$, provide a group of
linear isomorphism of the vertical subspace $V_uE$, $p(u)=x$. We
show in Theorem 2.2 that the elements of this group are also
isometries with respect to the pseudo--Riemannian metric induced
by $L$ in the vertical bundle over $E$. In Section 3 we treat the
case $L=F^2$, where $F$ is a Finsler function. In this case we
prove that $\nabla$ is metrizable, that is there exists a
Riemannian metric $h$ in $\xi$ such that $\nabla h=0$. Some
particular cases and examples are discussed in Section 4. The
notations and terminology are those from \cite{9} and \cite{5}.

\setcounter{section}{0}
\section{Berwald vector bundles}

Let $\xi=(E,p,M),\ p:E\to M$, be a vector bundle of rank $m$. Here
$M$ is a smooth i.e. $C^\infty$ manifold of dimension $n$. The type
fibre is $\R^m$ and $E$ is a smooth manifold of dimension $n+m$.
The projection $p$ is a smooth submersion. Let $(U,(x^i))$ be a
local chart on $M$ and let $\varepsilon_a(x)$, $x\in U$, be a field of
local sections of $\xi$ over $U$. Then every section $A$ of $\xi$
over $U$ takes the form $A=A^a(x)\varepsilon_a(x)$, $x\in U$, and an
element $u\in p^{-1}(x):=E_x$ can be written as $u=y^a\varepsilon_a(x)$,
$(y^a)\in\R^m$. The indices $i,j,k,...$ will range over
$\{1,2,...,n\}$ and the indices $a,b,c,...$ will take their values
in $\{1,2,...,m\}$. The convention on summation over repeated
indices of the same kind will be used.

The local coordinates on $p^{-1}(U)$ will be $(x^i,y^a)$ and a change
of coordinates $(x^i,y^a)\to(\wt x^i,\wt y^a)$ on $U\cap \wt U\ne\emptyset$ has the form
\begin{equation}
\begin{array}{l}
\wt x^i=\wt x^i(x^1,...,x^n),\ {\rm rank}\left(\derp{\wt x^i}{x^j}\right)=n,\vspace*{1,5mm}\\
\wt y^a=M^a_b(x)y^b,\ {\rm rank}(M^a_b(x))=m,\ \ \forall  x\in U\cap\wt U.
\end{array}
\end{equation}
On $E$ we have the vertical distribution $u\to V_uE={\rm
Ker}\,p_{x,u},$ where $p_{*}$ denotes the differential of $p$.
This consists of vectors which are tangent to fibres and it is
locally spanned by $\left(\dot\partial_a:=\derpp{y^a}\right)$. We shall regard
also the vertical distribution as a vector subbundle $VE:=\displaystyle\bigcup_{u\in E}V_uE\to E$ of $TE\to E$. Its sections will
be called vertical vector fields of $E$. The tensorial algebra
${\cal T}(VE)=\oplus{\cal T}^p_q(VE),\ p,q\in\N$ of this subbundle will
be used. Its elements will be indicated by the word ``vertical''.

\begin{definition} A nonlinear connection $N$ on $E$ is a
distribution $N:u\to N_uE$, $u\in E$, on $E$, which is supplementary
to the vertical distribution on $E$.
\end{definition}

We take the distribution $N$ as being locally spanned by
$\de_k=\partial_k-N^a_k(x,y)\dot\partial_a,$ for $\partial_k:=\derpp{x^k}$. By a
change of coordinates (1.1), the condition $\de_k=\derp{\wt x^i}{x^k}\ \wt\de_i$ is equivalent with
\begin{equation}
\wt N^a_j\partial_k\wt x^j=M^a_b(x)N^b_k(x,y)-\partial_k(M^a_b(x))y^b.
\end{equation}
It is important to notice that from (1.2) it follows that the set of functions
$F^a_{bk}(x,y)=\dot\partial_bN^a_k(x,y)$ behaves under a change of
coordinates (1.1) as the local coefficients of a linear connection
in the vertical bundle over $\xi$, that is
\begin{equation}
\wt F{}^a_{bk}(\wt x(x),\wt y(x,y))=M^a_c(x)\wt M^d_b(\wt x(x))
\derp{x^i}{\wt x^k}\ F^c_{di}(x,y)-\partial_i(M^a_c(x))\derp{x^i}{\wt x^k}\,y^c,
\end{equation}
where $\left(\derp{x^i}{\wt x^k}\right)$ is the inverse matrix of $\left(\derp{\wt
x^k}{x^j}\right)$ and $(\wt M^d_b)$ denotes the inverse matrix of $(M^b_c)$.

We should like to construct a linear connection $D$ in the
vertical bundle $VE\to E$. In order to do this it suffices to
provide $D_{\de_k}\dot\partial_a$ and $D_{\dot\partial_a}\dot\partial_b$. Using
(1.3) we have the possibility
\begin{equation}
D_{\de_k}\dot\partial_a=F^b_{ak}(x,y)\dot\partial_b,\
D_{\dot\partial_b}\dot\partial_c=V^a_{bc}(x,y)\dot\partial_a,\tag{1.3$^\circ$}
\end{equation}
where necessarily $(V^a_{bc}(x,y))$ behave like the components of
a vertical tensor field of type $(1,2)$.

In particular, we may take $V^a_{bc}=0$ and introduce

\begin{definition} The linear connection $D$ in the
vertical bundle $VE\to E$ given by
\begin{equation}
D_{\de_k}\dot\partial_a=F^b_{ak}(x,y)\dot\partial_b,\ \
D_{\dot\partial_a}\dot\partial_b=0,
\end{equation}
is called the {\it Berwald connection} associated to $N$.
\end{definition}

\begin{definition} We call the pair $(\xi,N)$ a Berwald
bundle if the functions $F^a_{bk}(x,$ $y)=\dot\partial_bN^a_b(x,y)$
depends on $x$ only.
\end{definition}

When $(\xi,N)$ is a Berwald bundle, the functions
$F^a_{bk}(x,y)=F^a_{bk}(x)$ define a linear connection $\nabla$ in
$\xi$ by
\begin{equation}
\nabla_{\partial_k}\varepsilon_b=F^a_{bk}(x)\varepsilon_a,
\end{equation}
for $(\varepsilon_a)$ a basis of local sections in $\xi$.

Conversely, if $\xi$ is endowed with a linear connection
of local coefficients $F^a_{bk}(x)$, then the functions
\begin{equation}
N^a_k(x,y)=F^a_{bk}(x,y)y^b,
\end{equation}
define by setting $\de_{\dot k}=\partial_{\dot k}-N^a_k(x,y)\dot\partial_a$
a nonlinear connection on $E$ such that $(\xi,N)$ becomes a
Berwald bundle. In other words, any vector bundle endowed with a
linear connection is a Berwald bundle.

We notice that the nonlinear connection (1.6) is positively
homogeneous of degree 1 in $y=(y^a)$. This suggests us to confine
ourselves to the pairs $(\xi,N)$ with the functions $(N^a_k(x,y))$
positively homogeneous of degree 1 in $y$. The examples to be
given later will fall in this category. This assumption requires
to eliminate from $E$ the image of the null section as we shall
do in the following.

It is well known that, see \cite{8}, \cite{9}, the Berwald connection
induces a covariant derivative in the tensorial algebra of the
vertical bundle. This splits in two operators of covariant
derivative. The first one is called $h$--covariant derivative and
is defined on functions and vertical vector fields as follows:
\begin{equation}
f_{\scri k}=\de_kf,\ X^a_{\scri k}=\de_kX^a+F^a_{bk}(x,y)X^b.
\end{equation}
It is extended by usual rules to any vertical tensor field. The second, called the
$v$-covariant derivative, is simply the partial derivative with
respect to $y$
\begin{equation}
f\bigi_a=\dot\partial_a f,\ X^a\bigi_b=\dot\partial_bX^a,
\end{equation}
since we have chosen $V^a_{bc}=0$.

We use the notation $\scri k$ and $\bigi_a$ for denoting the $h$-
and $v$-covariant derivatives of any vertical tensor field.

\begin{lemma} Let $\xi$ be endowed with a positively
$1$-homogeneous nonlinear connection $N$ and $\scri k$ the
$h$-covariant derivative defined by the Berwald connection
associated to it. Then
\begin{equation}
y^a_{\scri k}=0,
\end{equation}
holds.
\end{lemma}

\begin{proof} $y^a_{\scri k}=\de_ky^a+F^a_{bk}(x,y)y^b=F^a_{bk}(x,y)y^b-N^a_k(x,y)=0$
because of Euler theorem on homogeneous functions.
\end{proof}

\begin{lemma} Let $(\xi,N)$ be a Berwald bundle. Then
for any vertical tensor field $T$ of local coefficients
$T^{a_1...a_r}_{b_1...b_s}(x,y)$ we have
\begin{equation}
T^{a_1...a_r}_{b_1...b_s}{}_{\scri k}\bigi_a=T^{a_1...a_r}_{b_1...b_s}\bigi_{a\scri k}.
\end{equation}
\end{lemma}

\begin{proof} One verifies (1.10) by a direct calculation keeping in mind
that $F^a_{bk}=\dot\partial_aN^a_k$ do not depend on $y$.
\end{proof}

\section{Berwald bundles endowed with regular Lagrangians}
\setcounter{equation}{0}

We recall that in $\xi=(E,p,M)$, $E$ means in fact $E\setminus\{(x,0),\ x\in M\}$.

\begin{definition} A smooth function $L:E\to\R$ is called a regular Lagrangian on $E$ if
\begin{itemize}
\item[{\rm(i)}] the matrix with the entries
$g_{ab}(x,y)=\displaystyle\frac12\dot\partial_a\dot\partial_bL$ is nondegenerate,
\item[{\rm(ii)}] the quadratic form $g_{ab}(x,y)\zeta^a\zeta^b,\
(\zeta^a)\in\R^m$, is of rank constant.
\end{itemize}
\end{definition}
A regular Lagrangian $L$ induces a pseudo--Riemannian metric $g$
in the vertical bundle over $E$, given locally by
\begin{equation}
g(\dot\partial_a,\dot\partial_b)=g_{ab}(x,y).
\end{equation}
It provides also a set of vertical tensor fields by successively
deriving it with respect to $(y^a)$
\begin{equation}
C_{abc}(x,y)=\frac14\dot\partial_a\dot\partial_b\dot\partial_cL,\
D_{abcd}(x,y)=\frac18\dot\partial_a\dot\partial_b\dot\partial_c\dot\partial_dL,\
\mbox{etc.}
\end{equation}

\begin{definition} Let $\xi$ be endowed with a positively
$1$-homogeneous nonlinear connection $N$ and with a regular
Lagrangian $L$. We say that $N$ is compatible with $L$ if
\begin{equation}
L_{\scri k}:=\de_kL=0.
\end{equation}
\end{definition}

\begin{definition} If $(\xi,N)$ is a Berwald bundle with a
regular Lagrangian $L$ such that $(2.3)$ holds, the pair $(N,L)$
will be called a Berwald Lagrange structure, shortly a $BL$
structure for $\xi$.
\end{definition}

\begin{proposition} If $\xi$ has a $BL$ structure, then
\begin{itemize}
\item[{\rm(i)}] $g_{ab\scri k}=0,\ C_{abc\scri k}=0,\ D_{abcd\scri
k}=0$ etc.
\item[{\rm(ii)}] $g^{ab}{}_{\scri k}=0,\ y_{a\scri k}=0\
(y_a=g_{ab}y^b),\ C^a_{bc\scri k}=0\ (C^a_{ab}=g^{ae}C_{ebc})$.
\end{itemize}
\end{proposition}

\begin{proof} Easy consequences of (2.3) and of the commutation formulae
(1.10).
Assume that $\xi$ has a $BL$ structure.

Let be $c:[0,1]\to M$, $t\to c(t),\ t\in[0,1]$ a smooth curve on
$E$. A section $A$ of $\xi$ along $c$ given as $A(t)=A^a(t)\varepsilon_a$
is said to be {\it parallel} with respect to linear
connection $\nabla$ given by $(F^a_{bk}(x))$ if in a local chart on
$M$,
\begin{equation}
\frac{dA^a}{dt}+F^a_{bk}(c(t))A^b(t)\ \frac{dc^k}{dt}=0,
\end{equation}
holds.
\end{proof}

For the initial conditions $c(0)=x$ and $A^a(0)=A^a_0$, the system
of differential equations (2.4) admits a unique solution
$A^a(x(t))$ and if one assigns to $(A^a_0)\in E_x$ the element
$A^a(x(1))\in E_{c(1)=z}$ one obtains an application $P_c:E_x\to
E_z$ called {\it parallel translation} along $c$.

Moreover, from the linearity of the system (2.4) it follows that
$P_c$ is a linear isomorphism. Now if one considers all loops on
$M$ in $x\in M$, the corresponding parallel translations as linear
isomorphisms $E_x\to E_x$ provide a group with respect to their
composition, called the holonomy group $\phi(x)$ of $\nabla$ in
$x\in M$. When $M$ is connected, all these groups are isomorphic and
one speaks about the holonomy group $\phi$ of $\nabla$.

Let $L_x$ be the restriction of $L$ to the fibre $E_x$. We call
$L$-map a linear isomorphism $f:(E_x,L_x)\to(E_z,L_z)$ with the
property $L_x(u)=L_z(f(u))$ for every $u\in E_x$.

\begin{theorem} If $\xi$ admits a $BL$ structure,
then all parallel translations of $\nabla$ are $L$-maps. In
particular, the holonomy groups $\phi(x),$ $x\in M$, consists of
$L$-maps.
\end{theorem}

\begin{proof} Let $c:[0,1]\to M$ be a curve joining the points $x=c(0)$ and
$z=c(1)$ of $M$. Consider a parallel section $A(t):=A(c(t))$,
$t\in[0,1]$, of $\xi$ along $c$. We show that the function $f:t\to
L(x(t),A(t))$, $t\in[0,1],$ is constant. Indeed,
$$\frac{dL((x,y),A(t))}{dt}=(\partial_k)\ \frac{dx^k}{dt}+(\dot\partial_aL)\
\frac{dA^a}{dt}\bld{(2.4)}{=\!\!=\!\!=}L_{\scri k}\ \frac{dx^k}{dt}=0.$$
Consider $A_0\in E_x$ and $A(t)$ the unique solution of (2.4) with
the initial condition $A_0$. Then$P_c(A_0)=A_1$, where $A_1=A(1)$
and since $f$ is constant, we get
$L_x(A_0)=L_z(A_1)=L_z(P_c(A_0))$, q.e.d.
\end{proof}

The subset $I_x=\{A\in E_x\mid L_x(A)=1\}$ of $E_x$ is called the
indicatrix of $L.$ Let $G(I_x)$ be the group of all linear
isomorphisms of $E_x$ which leave invariant the indicatrix $I_x$.
From Theorem 2.1, it easily follows

\begin{corollary} The holonomy group $\phi(x)$ is a
subgroup of $G(I_x)$.
\end{corollary}

Let us continue to consider a parallel translation $P_c:E_x\to
E_z$. Its differential $(P_c)_{*,u}$, $u\in E$ is a linear
isomorphism $V_uE\to V_vE$ for $v=P_c(u)$ since $P_c$ itself is a
linear isomorphism and $T_u(E_x)$ is nothing but $V_uE$. We denote
it by $P_c^v$.

In particular, the differentials of the elements of $\phi(x)$ are
linear isomorphisms of $V_uE$ with $p(u)=x$ and these provide a
group $\phi^v(u)$ that is a subgroup of $GL(V_uE)$.

We call $\phi^v(u)$ the vertical lift of $\phi(x)$. For every
$u\in E$, $(V_uE,g_u)$ is a pseudo--Euclidean space.

\begin{theorem} The mapping $P^v_c:V_uE\to V_vE$,
$v\in P_c(u)$, are linear isometries of pseudo--Euclidean spaces.
In particular, the group $\phi^v(u)$ is a subgroups of the
isometries of $(V_uE,g_u)$.
\end{theorem}

\begin{proof} We fix the curve $c$ joining $x,z$ in $M$ and denote by
$(P^a_b)$ the matrix of $P_c:E_x\to E_z$ in the basis $(\varepsilon_a(x))$
and $\varepsilon_a(z)$. Here we tacitly assumed that $c$ is in a domain
$U$ of a local chart on $M$. If it is not so we divide $c$ in
segments. The matrix of $P^v_c$ is the same $P^a_b$ in the basis
$\dot\partial_a\bigi_u$ and $\dot\partial_a\bigi_v$. As $P_c$ is an $L$-map,
we have $L(x,u^a)=L(y,P^a_bu^b)$. We derive this equality two
times with respect to $(u^a)$ and we obtain
$$\derp{^2L}{u^a\partial u^b}=\derp{^2L}{y^c\partial y^d}\ P^c_aP^c_b,$$
that is $g_{ab}(u)=g_{cd}(v)P^c_aP^c_b$. This exactly means that
$P^v_c$ is an isometry of the pseudo--Euclidean spaces
$(V_uE,g_u)$ and $(V_vE,g_v)$. q.e.d.
\end{proof}

\section{Berwald bundles endowed with Finsler functions}
\setcounter{equation}{0}
\setcounter{definition}{0}
\setcounter{theorem}{0}
\setcounter{proposition}{0}
\setcounter{corollary}{0}

Let $\xi$ be a vector bundle.

\begin{definition} A smooth function
$F:E:=E\setminus0\to\R,$ $(x,y)\to F(x,y)$ is called a Finsler
function if
\begin{itemize}
\item[{\rm(i)}] $F(x,y)\ge0$,
\item[{\rm(ii)}] $F(x,\lambda y)=\lambda F(x,y),\ \forall  \lambda>0$,
\item[{\rm(iii)}] the matrix with the entries $g_{ab}(x,y)=\frac12\
\dot\partial_a\dot\partial_bF^2$ is positive definite
$(g_{ab}(x,$ $y)\zeta^a\zeta^a>0$ for $(\zeta^a)\ne0)$.
\end{itemize}
When $\xi$ is endowed with a Finsler function $F$ we call it a vector
Finsler bundle. If $(\xi,N)$ is a Berwald bundle, the pair $(N,F)$
will be called a Berwald Finsler structure, shortly a $BF$
structure for $\xi$ if $F_{\scri k}:=\de_kF=0$.
\end{definition}

If we put $L=F^2$, we obtain a regular Lagrangian. Thus any $BF$
structure is a $BL$ structure. As such, the properties of $BL$
structures proved in the previous section are valid for $BF$
structures. We show new properties for $BF$ structures.

\begin{proposition} If $\xi$ has a $BF$ structure
then $F_{\scri k}=0$ if and only if $g_{ab\scri k}=0$.
\end{proposition}

\begin{proof} $F_{\scri k}=0$ implies $L_{\scri k}=F^2{}_{\scri k}=0$ and by
Proposition 2.1, one gets $g_{ab\scri k}=0$. Conversely, applying
the Euler theorem to $F^2$ one obtains
$F^2(x,y)=g_{ab}(x,y)y^ay^b$. And the $h$-covariant derivation
yields $F^2{}_{\scri k}=2FF_{\scri k}=g_{ab\scri
k}y^ay^b+2g_{ab}y^ay^b{}_{\scri k}=0$ since $y^b_{\scri k}=0$.
Hence $F_{\scri k}=0.$ q.e.d.
\end{proof}

The pairs $(E_x,F_x)$ are called Minkowski spaces and $F_x$ is
called a Min\-kowski norm on $E_x$. The reason is that $F_x$,
besides the conditions (i)--(iii) from Definition 3.1 satisfies
also (see \cite{5} p.6; (iv) $F_x(y)>0$ whenever $y\ne0$; (v)
$F_x(y_1+y_2)\le F_x(y-1)+F_x(y_2)$.

The linear isomorphisms of $E_x$ keeping $F_x$ will be called
isometries.

We already know by Theorem 2.1 that if $\xi$ has a $BF$ structure,
all parallel translations defined by $\nabla$ are isometries.

In particular, the elements of $\phi(x)$ are isometries of the
Minkowski space $(E_x,F_x)$. And $\phi(x)$ is a subgroup of the
$G(I_x)$, the group of all linear isomorphism which leave
invariant the indicatrix $I_x$.

These facts are basic in the proof of the main result of this
section.

\begin{theorem} If $\xi$ has a $BF$ structure, the
linear connection $\nabla$ is metrizable, that is, there exists a
Riemannian metric $h$ in $\xi$ such that $\nabla h=0$.
\end{theorem}

\begin{proof} Let be $x_0\in M$ and the Minkowski space $(E_{x_0},F_{x_0})$.
The indicatrix $I_x$ is compact. It follows that the group
$G:=G(I_x)$ is a compact Lie group. We know that $G$ contains
$\phi(x)$ as a Lie subgroup but in general $\phi(x)$ is not
compact. Let $\langle\cdot\rangle$ be an arbitrary inner product in
$E_{x_0}$. Define a new inner product on $E_{x_0}$ by
$$h_{x_0}(u,v)=\frac1{{\rm Vol}(G)}\int_G<gu,gv>\,\mu_G,\ \mbox{ for }g\in G,\ u,v\in E_{x_0},$$
where $\mu_G$ denotes the bi--invariant Haar measure on $G$. It
follows that $h_{x_0}$ is $G$--invariant and, in particular, it is
$\phi(x_0)$--invariant, i.e., $h_{x_0}(Pu,Pv)=h_{x_0}(u,v)$ for
any $P\in\phi(x_0)$. Now we transfer $h_{x_0}$ to all the points
of $M$. For any point $x\in M$, we consider a curve $c$ joining
$x$ with $x_0$ $(c(0)=x$, $c(1)=x_0)$.
\end{proof}

Define $h_x(A,B)=h_{x_0}(P_cA,P_cB)$, $A,B\in E_x$. The property
that $h_{x_0}$ is $\phi(x_0)$--invariant assures that $h_x$ does
not depend on the curve $c$.

The mapping $h:x\longrightarrow h_x:E_x\times E_x\to R$ is smooth since $P_c$
smoothly depends on $x$ by a general result about dependence of
solutions of an ordinary differential equation on initial data.
Thus a Riemannian metric $h$ in $\xi$ is obtained. The proof is
ended with the help of

\begin{lemma} Let $h$ be a Riemannian metric in $\xi$
and $t\to c(t)$, $t\in\R$, a curve with $c(0)=x\in M$. Then
\begin{equation}
\lim_{t\to0}\ \frac1t\left(h_{c(t)}(P_cA,P_cB)-h_x(A,B)\right)=
\left(\nabla_{\dot c(0)}h\right)(A,B)(x),
\end{equation}
where $A,B\in E_x$ and $P_c:E_x\to
E_{c(t)}$ is the parallel translation along $c$.
\end{lemma}

Indeed, by the definition of $h$, the term in the left side of
$(3.1)$ vanishes.

For the proof of Lemma 3.1 we refer to \cite{2}.

\section{Particular cases}
\setcounter{equation}{0}
\setcounter{definition}{0}
\setcounter{theorem}{0}
\setcounter{proposition}{0}
\setcounter{corollary}{0}

\noindent{\bf4.1.} Let $\xi=\tau_M=(TM,\tau,M)$ be the tangent bundle of
$M$. If $\tau_M$ is endowed with a Finsler function $F$, the pair
$(M,F)$ is called a Finsler manifold. For the geometry of these
manifolds we refer to \cite{7}, \cite{5}.

The Finsler function $F$ induces the Cartan nonlinear connection
$\bld\circ N{}^i_j(x,y)=\gamma^i_{j0}-C^i_{jk}\gamma^j_{00}$,
where
$2\gamma^i_{jk}=g^{ih}(\partial_jg_{kh}+\partial_kg_{jh}-\partial_hg_{jk})$,
$2C^i_{jk}=g^{ih}\dot\partial_hg_{jk}$, $\gamma^i_{j0}=\gamma^i_{jk}y^k$
and $\gamma^i_{00}=\gamma^i_{jk}(x,y)y^jy^k$. Of course,
$g_{jk}=\displaystyle\frac12\ \dot\partial_j\dot\partial_kF^2$ denotes the Finsler
metric. This nonlinear connection is $p$--homogeneous of degree 1
in $y$ and is compatible with $F$, that is, $F_{\scri k}=0$. If
the local coefficients $\bld\circ
G{}^i_{jk}(x,y)=\dot\partial_j\bld\circ N{}^i_k(x,y)$ of the Berwald
connection associated to $(\bld\circ N{}^i_j)$ depend on $x$ only,
the Finsler manifold $(M,F)$ is called a Berwald space. In
\cite[p.263--64]{5} there are given five properties
characterizing the Berwald spaces. Among them we notice the
condition $C_{ijk\scri h}=0$. Thus,
if $\tau_M$  is endowed with a Finsler function
$F$ for wich $C_{ijk\scri h}=0$, the pair $(\bld\circ N,F)$ is a Berwald Finsler structure. Particularizing our
results from Section 3 the results previously proved by Y. Ichijyo
\cite{6} and Z. Szab\'o \cite{10} are obtained.\bigskip

\noindent{\bf4.2.} Let $\xi=\tau_M$ be endowed with a regular Lagrangian
$L$. Then the pair $(M,L)$ is called a Lagrange manifold. For the
geometry of Lagrange manifolds we refer to \cite{9}. The
Lagrangian defines a nonlinear connection
$N^i_j(x,t)=\dot\partial_jG^i$, where
$4G^i=g^{ik}(y^h\dot\partial_k\partial_hL-\partial_kL)$ but, in general, this is
not $p$-homogeneous nor compatible with $L$. We notice that $N$ is
provided by the semi--spray $(G^i(x,y))$ that in turn is derived
from $L$. A question is whether there exist Lagrangians which to
generate sprays, that is the functions $(G^i(x,y))$ to be
$p$--homogeneous of degree 2 in $y$.

A first example was given and studied in \cite{3}. A larger class
of such Lagrangians called $\varphi$--Lagrangians is proposed and
studied in \cite{4}.

Let $\tau_M$ be endowed with a Finsler function $F$. We eliminate
the image of null section $\{0_x,\ x\in M\}$ from $TM$. Let
$\varphi:\R_+\to\R$ be a smooth function. Then $L=\varphi(F^2)$ is a
Lagrangian and one proves (\cite{4}) that if $\varphi'(t)\ne0$ and
$\varphi'(t)+t\varphi''(t)\ne0$ for any $t\in{\rm Im}(F^2)$, then $L$ is a
regular Lagrangian, called a $\varphi$-Lagrangian. For a
$\varphi$-Lagrangian, the functions $G^i(x,y)$ are $p$-homogeneous of
degree 2 in $y$. Moreover, $G^i(x,y)=\gamma^i_{00}$ and the
nonlinear connection $N$ provided by a $\varphi$-Lagrangian coincides
with the Cartan nonlinear connection $\bld\circ N$ of $(M,F)$, cf.
\cite{4}. It follows easily that $N$ is compatible with $L$.
Furthermore, the pair $(\tau_M,N)$ is a Berwald bundle if and only
if $(\tau_M,\bld\circ N)$ is a Berwald bundle. It follows that
$(N,L)$ is a Berwald Lagrange structure for $\tau_M$ if and only
if $(\bld\circ N,F)$ is a Berwald Finsler structure for $\tau_M$.
The connection $\nabla$ is the same for these structures and by
Theorem 3.1 it is metrizable.

\newpage

\def\rr{\mathbb{R}}
\def\to{\rightarrow}
\def\9{\infty}
\def\a{{\alpha}}
\def\b{{\beta}}
\def\calx{{\cal X}}
\def\calf{{\cal F}}
\def\calh{{\cal H}}
\def\p{\partial}

\runningauthor={M. ANASTASIEI}
\runningtitle={MINKOWSKIAN $G-$STRUCTURES IN VECTOR BUNDLES}
\noindent
\baselineskip 8pt
\noindent{\footnotesize{Lagrange and Hamilton Geometries and Applications\hfill\break  Ed Fair Partners, Bucharest, 2004, 1-10}}
\vskip 2cm
\baselineskip 11.5pt plus .15pt
\centerline{\bf\Large MINKOWSKIAN $G-$STRUCTURES}
\vskip .2cm
\centerline{\bf\Large IN VECTOR BUNDLES}
\vskip .5cm
\centerline{\bf by Mihai ANASTASIEI}
\vskip .5cm

\begin{center}
    {\footnotesize{Dedicated to the Memory of Grigorios TSAGAS
(1935-2003),\par President of Balkan Society of Geometers (1997-2003)}}
\end{center}

\begin{abstract}
A natural generalization of the usual inner product on $\rr^m$ is
the so-called Minkowski norm. For a smooth vector bundle $\xi$
with the type fibre $\rr^m$ endowed with a Minkowski norm, a
$G-$structure for $\xi$, generalizing the $O(m)-$structures, is
defined and called a Minkowskian $G-$structure. Several properties
of these structures are pointed out in Theorem A and Theorem B.
Some of them extend to the vector bundles the results given by Y.
Ichijio ([2], [3]) for tangent bundle. If applied to the cotangent
bundle our results enrich the geometry of Cartan spaces presented
in the monograph [5] by R. Miron et al. on the line of our paper
[1].

{\bf MSC2000: 53C60, 53C10.}

{\bf Keywords and phrases:} vector bundles, Minkowskian norms,
$G$-struc\-tu\-res.
\end{abstract}

\setcounter{section}{0}
\section*{Introduction}

Let $\xi=(E,\pi,M)$ be a smooth i.e. $C^\infty$vector bundle of
rank $m$. Assume that $M$ is connected. The type fibre of $\xi$ is
$\rr^m$ and its structural group is $GL(m,\rr)$. The linear space
$\rr^m$ has a natural inner product $<$ , $>$ and it is well known
that this can be transferred with the help of the bundle charts to
a Riemannian structure $g$ in $\xi$ if and only if $\xi$ admits an
$O(m)-$structure. Moreover, it is also known that if $\xi$ admits
an $O(m)-$structure, then there exists a linear connection
$\nabla$ in $\xi$ that is metrical with respect to $g$ $(\nabla
g=0)$. Then $\nabla$ can be also seen as a principal connection in
the principal bundle of the frames of $\xi$ having the property
that its connection 1-form takes the values in the Lie algebra of
$O(m)$.

The condition $\nabla g=0$ is equivalent with the fact that all
parallel translations defined by $\nabla$ are isometries.

Consequently, the fibres $(E_x,g_x),x\in M$, of $\xi$ are all
congruent i.e. linearly isometrically isomorphic.



We prove the following theorems.

\medskip
{\bf Theorem A.} {\it If $\xi=(E,\pi,M)$ admits a Minkowskian
$G_f-$structure, then
\begin{itemize}
    \item[{\rm{i)}}] Each fibre $E_x$, $x\in M$ becomes a Minkowski
    space,
    \item[{\rm{ii)}}] A Finsler function $F(x,y)=f(\mu^a_b (x)y^b)$, $\mu^a_b (x)\in
    GL(m,\rr)$ is defined on $E$.

    Denote by $(g_{ab}(x,y))$ the Finsler metric associated to
    $F$.
    \item[{\rm{iii)}}] Let $\nabla$ be a linear $G_f-$connection and $|k$
    be the horizontal covariant derivative defined by its vertical
    lift to $E$. Then $F_{|k}=0$ and $g_{ab|k}=0$.
    \item[{\rm{iv)}}] the fibres $E_x,x\in M$ are all congruent each
    others as Minkowski spaces.
\end{itemize}}

\medskip

A pair $(F,\nabla)$ with $F$ a Finsler function on $E$ and
$\nabla$ a linear connection in $\xi$ such that $g_{ab|k}=0$ will
be called a $(F,\nabla)-${\bf structure} for $\xi$.

Theorem A says that if $\xi $ admits a $G_f-$ structure, then $\xi
$ admits a $(F,\nabla)$- structure. The converse holds, too.

\medskip

{\bf Theorem B.} {\it Let $\xi=(E,\pi,M)$ be a vector bundle of
rank $m$. Assume that it admits $(F,\nabla)-$structure. Then $F$
induces a Minkowski norm $f$ on $\rr^m$ and $\xi$ admits a
$G_f-$structure such that the $(F,\nabla)-$structure induced by it
is just that initially  given.}

\medskip

The notions entering in the contents of these theorems will be
explained below in the appropriate places.

Our results extend to any vector bundle some of the results due to
Y. Ichijio for tangent bundle, [2], [3].

\section{Vector bundles. Minkowskian $G-$structures}
Let $\xi=(E,\pi,M)$ be a smooth vector bundle of rank $m$.

Assume that $M$ is connected and its dimension is $n$. Then $E$ is
a smooth manifold of dimension $n+m$.

Let $\{(U_\a,\psi_\a)\}_{\a\in A}$ be a smooth atlas on $M$. A
vector bundle atlas is then $\{(U_\a,\varphi_\a,\rr^m)\}_{\a\in
A}$, where $\varphi_\a:\pi^{-1}(U_\a)\to U_\a\times\rr^m$ are
diffeomorphisms of the form
$\varphi_\a(u)=(\pi(u),\varphi_{\a,\pi(u)}(u))$,
$u\in\pi^{-1}(U_\a)$ such that for every $x\in U_\a\cap U_\b\neq
\phi$, $\varphi_{\b,x}\circ \varphi_{\a,x}^{-1}$ belongs to
$GL(m,\rr)$.

The manifold structure of $E$ is defined by the atlas
$\{(\pi^{-1}(U_\a),\phi_\a)\}$ with $\phi_\a:\pi^{-1}(U_\a)\to
\psi_\a(U_\a)\times\rr^m$ given by
$\phi_\a(u)=(\psi_\a(\pi(u)),\varphi_{\a,\pi(u)}(u))$. The
mappings
$(\phi_\b\circ\phi_\a^{-1})(u)=((\psi_\b\circ\psi^{-1}_\a)(\pi(u))$,
$(\varphi_{\b,\pi(u)}\circ\varphi^{-1}_{\a,\pi(u)})(u))$
are smooth.

Let $(e_a),a=1,2,...,m$ be the canonical basis of $\rr^m$. The
mappings $\varepsilon_{\a,a}:U_\a\to\pi^{-1}(U_\a)$,
$\varepsilon_{\a,a}(x)=\varphi^{-1}_\a(x,e_a)$ are $m$ linearly
independent local sections of $\xi$, that is
$(\varepsilon_{\a,a}(x))$ is a basis of the fibre $E_x$, $x\in M$.

If we put $\psi_\a(x)=(x^i)$, $\psi_\b(x)=(\widetilde{x}^j)$,
$i,j,k...=1,2,...,n$, then $\psi_\b\circ\psi^{-1}_\a$ has the form
\begin{itemize}
    \item[(1.1)] $\qquad \widetilde{x}^j=\widetilde{x}^j(x^1,...,x^n)$,
    rank$\left(\dfrac{\partial\widetilde{x}^j}{\partial
    x^i}\right)=n$.
\end{itemize}

For $u\in E$ with $\pi(u)=x$, we can write
$u=y^a\varepsilon_{\a,a}(x)=\widetilde{y}^b\varepsilon_{\b,b}(x)$.
If we put $(\varphi_{\b,x}\circ\varphi^{-1}_{\a,x})(e_a)=M^b_a
(x)e_b$, then $s_{\a,a}(x)=M^b_a (x)s_{\a,b}$ and
$\widetilde{y}^b=M^b_a (x)y^a$.

The local coordinates on $E$ will be $(x^i,y^a)$ and a change of
coordinates $(x^i,y^a)\to(\widetilde{x}^i,\widetilde{y}^b)$ has
the form
\begin{itemize}
    \item[(1.2)] $\qquad \widetilde{x}^i=\widetilde{x}^i(x^1,...,x^n)$,
    rank$\left(\dfrac{\partial\widetilde{x}^i}{\partial
    x^j}\right)=n$,

    $\qquad \widetilde{y}^a=M^a_b(x)y^b$, rank$(M^a_b(x))=m$.
\end{itemize}

The Einstein convention on summation will be applied for the
indices $i,j,k,...=1,...,n$ as well as for the indices
$a,b,c...=1,...,m$.

Let $G$ be a Lie subgroup of $GL(m,\rr)$. One says that $\xi$
admits a $G$-{\bf structure} if there exists a vector bundle atlas
$\{(U_\a,\varphi_\a,\rr^m)\}_{\a\in A}$ such that the mapping
$U_\a\cap U_\b\to GL(m,\rr)$, $x\to\varphi_{\b,x}\circ
\varphi^{-1}_{\a,x}$ take their values in $G$.

The existence of a $G-$structure in $\xi$ is equivalent with the
existence of an open covering $(U_\a)_{\a\in A}$ of $M$ and of $m$
sections $s_{\a,a}:U_\a\to\pi^{-1}(U_\a)$ for every $\a\in A$ such
that
\begin{itemize}
    \item[i)] $(s_{\a,a}(x))$ is a basis in $E_x$, $x\in U_\a$,
    \item[ii)] $s_{\a,a}(x)=M^b_a (x)s_{\b,b}(x)$ with $M^b_a (x)\in
    G$, for every $\a,\b\in A$ with $U_\a\cap U_\b\neq\phi$.
\end{itemize}

This means that a $G-$structure in $\xi$ is a {\bf reduction to
$G$} of the principal bundle of frames of $\xi$ (see [4]).

The basis $(s_{\a,a}(x))$ are called frames {\bf adapted} to the
given $G-$structure.

A {\bf Minkowski norm} $f$ on $\rr^m$ is a non-negative real
function on $\rr^m$ with the properties:
\begin{enumerate}
    \item $f$ is smooth on $\rr^m\setminus 0$,
    \item $f(\lambda y)=\lambda f(y)$ for all $\lambda>0$
    \item The matrix with the entries $g_{ij}(y)=\dfrac{1}{2}\dfrac{\partial^2 f^2}{\partial y^i \partial
    y^j}$ is positive definite.
\end{enumerate}

The pair $(\rr^m,f)$ is called a {\bf Minkowski space}.

If it happens that $f(-y)=f(y)$, then $f(\lambda y)=|\lambda|f(y)$
and one says that $f$ is an absolutely homogeneous Minkowski norm.
One proves [BCS, p.6] that any absolutely homogeneous Minkowski
norm is a norm on $\rr^m$.

Let $\rr^m$ be endowed with a Minkowski norm $f$ and let be
$G_f=\{T\in GL(m,\rr)| f(Ty)=f(y),\forall y\in \rr^m\}$.

Then $G_f$ is a closed subgroup of $GL(m,\rr)$. Indeed, if
$(T_n)_{n\geq 0}$ is a sequence in $G_f$ that converges to $T_0$,
making $n\to \9$ in the equality $f(T_n y)=f(y)$ $\forall y$ we
get $f(T_0 y)=f(y)$ $\forall y$, that is $T_0\in G_f$. It follows
that $G_f$ is a Lie subgroup of $GL(m,\rr)$. A $G_f-$structure in
$\xi$ will be called a {\bf Minkowskian structure in $\xi$}.

\section{Finsler vector bundle. Connections}

Let $\xi=(E,\pi,M)$ be a smooth vector bundle of rank $m$.

A {\bf Finsler function} on $E$ is a non-negative real function
$F$ on $E$ with the properties:
\begin{enumerate}
    \item $F$ is smooth on $E\setminus\{(x,0),x\in M\}$ and only
    continuous on the set $\{(x,0),x\in M\}$,
    \item $F(x,\lambda y)=\lambda F(x,y)$ for all $\lambda>0$,
    \item The matrix with the entries $g_{ab}(x,y)=\dfrac{1}{2}\dfrac{\partial^2 F^2}{\partial y^a\partial
    y^b}$is positive definite.
\end{enumerate}

The pair $(\xi,F)$ is called a {\bf Finsler vector bundle}.

For every $x\in M$, the function $F_x:E_x\to R$ given by
$F_x(u)=F(x,u)\forall u \in E_x$ is a Minkowski norm on $E_x$.
Thus the fibres of a Finsler vector bundle are all Minkowski
spaces.

On $E$ we have the vertical distribution $u\to V_u E=\ker
\pi_{*,u}$ made by the vectors which are tangent to fibres.

A distribution $u\to H_u E$ which is supplementary to it is called
a horizontal distribution or a nonlinear connection for $\xi$. The
vertical distribution is spanned by
$\left(\dfrac{\partial}{\partial y^a}\right)$. As a local basis
for the horizontal distribution it is usually taken
$\delta_i=\dfrac{\partial}{\partial
x^i}-N^a_i(x,y)\dfrac{\partial}{\partial y^a}$, where the
functions $(N^a_i(x,y))$ are called the local coefficients of a
given nonlinear connection. By a change of coordinates on $E$
these functions behave in such a way that the transformation law
\begin{itemize}
    \item[(2.1)] $\qquad \delta_i=\dfrac{\partial \widetilde{x}^j}{\partial
    x^i}\widetilde{\delta}_j$,
\end{itemize}
is assured.

When a nonlinear connection is considered we have the
decomposition
\begin{itemize}
    \item[(2.2)] $\qquad T_u E=H_u E\oplus V_u E$, $u\in E$.
\end{itemize}

Then all the geometric objects on $E$ can be decomposed
accordingly.

If the functions $(N^a_i(x,y))$ are linear in $(y^a)$, that is
$N^a_i(x,y)=\Gamma^a_{bi}(x)y^b$, the nonlinear connection becomes
a linear one. In this case we may define an operator of covariant
derivative $\nabla$: $\calx(M)\times\Gamma(E)\to\Gamma(E)$,
$(X,\sigma)\to \nabla_X\sigma=X^j(x)\left(\dfrac{\partial
\sigma^a}{\partial
x^j}+\Gamma^a_{bj}(x)\sigma^b\right)\varepsilon_a$, where
$X=X^j\dfrac{\partial}{\partial x^j}$,
$\sigma=\sigma^a(x)\varepsilon_a$ and
$\nabla_{\dfrac{\partial}{\partial
x^k}}\varepsilon_a=\Gamma^b_{ak}(x)\varepsilon_b$.

We denoted by $\Gamma(E)$ the $\calf(M)-$module of sections in
$\xi$ and $\calx(M):=\Gamma(TM)$.

A linear connection $\nabla$ in $\xi$ induces a linear connection
$D$ in the vertical bundle over $E$ as follows:

The operator $D$: $\calx(E)\times \Gamma(VE)\to\Gamma(VE)$ is
defined by the equations
\begin{itemize}
    \item[(2.3)] $\qquad D_{\delta_k}\dot\partial_a=\Gamma^b_{ak}(x)\dot\p_b$,
    $D_{\dot\p_b}\dot\p_a=V^c_{ab}(x,y)\dot\p_c$,
\end{itemize}
where $\dot\p_a:=\dfrac{\p}{\p y^a}$ and $C^c_{ab}(x,y)$ are the
components of a vertical tensor field
$V^c_{ab}\dot\p_c\otimes\delta y^a\otimes\delta y^b$, $\delta
y^a=dy^a+\Gamma^a_{b_i}(x)y^b dx^i$.

The functions $C^c_{ab}$ can be taken zero. In such a case, $D$
will be called the vertical lift of $\nabla$.

We shall use $D_{\delta_k}$ for defining a {\bf horizontal
covariant derivative} in the tensor algebra of the vertical bundle
over $E$. It will be denoted by $|k$ and it is obtained by the
usual extension procedure starting with
\begin{itemize}
    \item[(2.4)] $\qquad f_{|k}=\delta_k f$ for every real function $f$ on $E$,

    $\qquad A^a_{|k}=\delta_k A^a+\Gamma^a_{bk}(x)A^b$, for
    $A=A^a\dot\p_a$ a section in the vertical bundle.
\end{itemize}

The vector field $C=y^a\dot\p_a$ is called the Liouville vector
field on $E$ and $D_{\delta_k}C=y^a_{|k}\dot\p_a$ is called the
deflexion tensor field of $D$. We have

\medskip

{\bf Lemma 2.1} {\it $D_{\delta_k}C=0$.}

\medskip

Indeed, $y^a_{|k}=\delta_k
y^a+\Gamma^a_{b_k}(x)y^b=-\Gamma^a_{b_k}(x)y^b+\Gamma^a_{b_k}(x)y^b=0$.

\medskip

{\bf Lemma 2.2} {\it For every real function $H$ on $E$ we have
$$\dot\p_a(H_{|k})=(\dot\p_a H)_{|k}.$$}


\noindent{\it Proof.} A direct calculation keeping in mind that $\dot\p_a
H$ are he coefficients of a vertical 2-form and so $$(\dot\p_a
H)_{|k}=\delta_k(\dot\p_a H)-\Gamma^b_{ak}(x)\dot\p_b H.$$

\section{Proof of Theorem A} Let be $\xi$ with the type fibre the
Minkowski space $(\rr^m,f)$. Assume that $\xi$ admits a
$G_f-$structure. Let $(s_{\a,a}(x))$ be a frame in $E_x$ adapted
to this $G_f-$structure. For $u\in E_x$ we have $u=y^a
\varepsilon_{\a,a}(x)=z^a s_{\a,a}(x)$. We define
$F_\a:E_x\to[0,\9)$ by $F_\a (u)=f(z^a)$. For $x\in U_\a\cap U_\b$
we have also $F_\b(u)=f(\widetilde{z}^b)$, where
$(\widetilde{z}^b)$ are given by $u=\widetilde{z}^b s_{\b,b}(x)$.
It follows that $\widetilde{z}^b=M^b_a(x)z^a$ with $(M^b_a(x))\in
G_f$. Consequently, $f(\widetilde{z}^b)=f(z^a)$ and
$F_\a(u)=F_\b(u)$. In the other words, the function $F$ defined by
$F(u)=f(z^a)$ does not depend on the chosen local chart. It is
clear that for every $x\in M$, this $F$ is a Minkowski norm on
$E_x$. Thus i) of Theorem A is proved.

If we put
$s_{\a,a}(x)=\lambda^b_a(x)\varepsilon_{\a,b}(x),(\lambda^b_a(x))\in
GL(m,\rr)$, it results $z^a=\mu^a_b(x)y^b$ with
$(\mu^a_b)=(\lambda^a_b)^{-1}$.

The function $F:(x,y)\to F(x,y)=f(\mu^a_b(x)y^b)$ is a Finsler
function on $E$. Indeed, this is smooth on $E\setminus\{(x,0),x\in
M\}$ since $f$ is smooth on $\rr^m\setminus 0$ and the functions
$(\mu^a_b)$ as the entries of the inverse matrix of a matrix whose
entries are smooth are also smooth.

Note that $y^a=0$ if and only if $z^a=0$. The function $F$ is
positively homogeneous of degree 1 in $(y^a)$ because $f$ is
homogeneous of degree 1. The matrix $(g_{ab}(x,y))$ has in this
case the entries $g_{ab}(x,y)=\dfrac{1}{2}\dfrac{\p^2 f^2}{\p y^c
\p y^d}\mu^c_a(x)\mu^d_b(x)$ and for any $(\varsigma^a)\in\rr^m$
we get
$$g_{ab}(x,y)\varsigma^a\varsigma^b=\dfrac{1}{2}\dfrac{\p^2 f^2}{\p y^c \p
y^d}\sigma^c\varsigma^d\textrm{ for
}\varsigma^c=\mu^c_a(x)\xi^a.$$ It follows that the matrix
$(g_{ab})$ is positive definite since $f$ is a Minkowski norm.
Thus ii) of Theorem A is proved. Note that ii) implies i) because
the Finsler function $F$ will provide by restriction a Minkowski
norm in each fibre of $\xi$.

We fix an open set $U_\a$ in $M$ and take a field of frames
$(s_a)$ of local sections adapted to the Minkowskian structure
$G_f$. A $G_f-$connection is a connection in the principal bundle
of frames of $\xi$ whose connection 1-form $\theta$ has values in
the Lie algebra $g_f$ of $G_f$. The 1-form $\theta$ is completely
determined by a matrix $(\theta^a_b(x))$ of 1-forms on $M$ using a
fixed basis of $g_f$. The operator of covariant derivative is
$\nabla^*_X s_a=\theta^b_a(X)s_b$. If one sets
$\theta^b_a=\Gamma^{*b}_{ak}(x)dx^k$, then
$\theta^b_a(X)=X^k\Gamma^{*b}_{ak}(x)$ for $X=X^k\p_k$ and so
$\nabla^*_X s_a=\Gamma^{*b}_{ak}(x)X^k s_b$ with the matrix
$(\Gamma^{*b}_{ak}(x)X^k)$ in $g_f$. The following Lemma gives a
characterization of the elements of $g_f$.

\medskip

{\bf Lemma 3.1} {\it A matrix $A=(A^a_b)\in g_f$ if and only if
\begin{itemize}
    \item[{\rm(3.1)}] $\qquad \dfrac{\p f}{\p z^a}A^a_b z^b=0$ for
    every $(z^a)\in\rr^m$.
\end{itemize}}

\medskip

\noindent{\it Proof.} If $A\in g_f$ then $\exp tA\in G_f$, hence $f((\exp
tA)z)=f(z)$ $\forall z\in\rr^m$. This means that
$\dfrac{d}{dt}f((\exp tA)z)|_{t=0}=0$, a equation that is
equivalent with (3.1).

Note that for $f=\sqrt{\langle z,z\rangle}$, the equation (3.1)
reduces to the skew symmetry of $A$.

Let $(\varepsilon_a)$ be the natural frame (corresponding to the
canonical basis $(e_a)$ in $\rr^m$) on $U_\a$ so that
$u=y^a\varepsilon_a=z^b s_b$. It result $z^b=\mu_c^b y^c$ for
$\mu^b_c=(\lambda^a_b)^{-1}$, where as before
$s_a(x)=\lambda^b_a(x)\varepsilon_b(x)$. We put $\nabla^*
\varepsilon_a=X^k\Gamma^b_{ak}(x)\varepsilon_b$.

Then $\nabla^*_X s_a=\nabla^*_X(\lambda^b_a
\varepsilon_b)=X^k(\p_k\lambda^c_a+\Gamma^c_{bk}(x)\lambda^b_a)\varepsilon_c$.

If we think $(\lambda^b_a)$ as a set of $m$ vector fields, the
last parenthesis is $\nabla_k \lambda^c_a$, hence $\nabla^k_X
\varepsilon_a=X^k(\nabla_k \lambda^c_a)\varepsilon_c$.

On the other hand, $\nabla_X^* s_a=X^k\Gamma^{*b}_{ak}\lambda^c_b
\varepsilon_c$ and by a comparison we get $\nabla_k
\lambda^c_a=\Gamma^{*b}_{ak}\lambda^c_b$ or
\begin{itemize}
    \item[(3.2)] $\qquad\Gamma^{*b}_{ak}=(\nabla_k \lambda^c_a)\mu^b_c$.
\end{itemize}

Lemma 3.1 applied for $(X^k\Gamma^{*b}_{ak})$ says that
\begin{itemize}
    \item[(3.3)] $\qquad\dfrac{\p f}{\p z^b}(X^k
    \Gamma^{*b}_{ak})z^a=0$ $\forall (z^a)\in\rr^m$.
\end{itemize}

Inserting here $\Gamma^{*b}_{ak}$ given by (3.2) one gets
\begin{itemize}
    \item[(3.4)] $\qquad\dfrac{\p f}{\p z^b}(\nabla_k \lambda^c_a)\mu^b_c z^a=0$.
\end{itemize}

Now we consider the nonlinear connection
$N^a_k(x,y)=\Gamma^a_{bk}(x)y^b$ and the vertical lift of linear
connection $(\Gamma^a_{bk}(x))$ denoting by $|i$ the corresponding
horizontal covariant derivative.

Recall that $F(x,y)=f(\mu^a_b (x)y^b)$ and compute $F_{|k}$.

We have $$F_{|k}=\dfrac{\p F}{\p x^k}-N^a_k\dfrac{\p F}{\p
y^a}=\dfrac{\p f}{\p z^c}\dfrac{\p \mu^c_b}{\p x^k}
y^b-\Gamma^a_{bk}y^b\dfrac{\p f}{\p z^c}\mu^c_a=\dfrac{\p f}{\p
z^c}y^b\left(\dfrac{\p \mu^c_b}{\p
x^k}-\Gamma^a_{bk}\mu^c_a\right).$$

From $\mu^c_b \lambda^a_c=\delta^a_b$ it follows $$\dfrac{\p
\mu^c_b}{\p x^k}\lambda^a_c=-\mu^c_b\dfrac{\p \lambda^a_c}{\p
x^k}\textrm{ or }\dfrac{\p \mu^c_b}{\p x^k}=-\mu^d_b\dfrac{\p
\lambda^a_d}{\p x^k}\mu^c_a.$$

Inserting this in the last form of $F_{|k}$ we get
$$F_{|k}=-\dfrac{\p f}{\p z^c}y^b\mu^d_b\left(\dfrac{\p
\lambda^a_d}{\p
x^k}+\Gamma^a_{ek}\lambda^e_d\right)\mu^c_a=-\dfrac{\p f}{\p
z^c}z^d(\nabla_k \lambda^a_d)\mu^c_a=0$$ by (3.4).

Using Lemma 2.2, we compute
$$g_{ab|k}=\dfrac{1}{2}(\dot\p_a\dot\p_b
F^2)_{|k}=\dfrac{1}{2}\dot\p_a\dot\p_b(F^2_{|k})=0$$ since
$F_{|k}=0\Leftrightarrow F^2_{|k}=0$.

Thus the point iii) of Theorem A is proved.

Let us consider a smooth curve $c:[0,1]\to M; t\to c(t)$, joining
the points $x=c(0)$ and $y=c(1)$ and let us denote by $P_c:E_x\to
E_y$ the parallel translation along $c$ defined by a linear
$G_f-$connection $\nabla$ in $\xi$.

It associates to an element $u=A(0)\in E_x$ the unique element
$A(1)$ from $E_y$, where $t\to A(t)$ is a section in $\xi$ along
$c$ which is parallel along $c$, that is its components $(A^a(t))$
are solutions of the system of differential equations
\begin{itemize}
    \item[(3.5)]
    $\qquad\dfrac{dA^a}{dt}+\Gamma^a_{bk}(x(t))A^b(t)\dfrac{dx^k}{dt}=0$.
\end{itemize}

Consider $F$ in the points $(x(t),A(t))$ and compute
$$\dfrac{dF(x(t),A(t))}{dt}=\dfrac{\p F}{\p
x^k}\dfrac{dx^k}{dt}+\dfrac{\p F}{\p
y^a}\dfrac{dA^a}{dt}\buildrel{(3.5)}\over{=}\dfrac{dx^k}{dt}\left(\dfrac{\p
F}{\p x^k}-\Gamma^a_{bk}A^b\dfrac{\p F}{\p y^a}\right)=0$$ because
of $F_{|k}=0$.

Thus the function $t\to F(x(t),A(t))$ is constant. Hence
$F_x(u)=F_y(P_c u)$. In the other words, the linear isomorphism
$P_c$ preserves the Minkowskian norms. This proves the point (iv)
and thus Theorem A is completely proved.

\section{Proof of Theorem B}

Let be a pair $(F,\nabla)$ with $F$ a Finsler function on $E$ and
$\nabla$ a linear connection in $\xi$ such that $g_{ab|k}=0$.

\medskip

{\bf Lemma 4.1} {\it $g_{ab|k}=0$ implies $F_{|k}=0$.}

\medskip

\noindent{\it Proof.} The homogeneity of $F$ implies by a repeated use of
the Euler theorem that $F^2(x,y)=g_{ab}(x,y)y^a y^b$.

Then $F^2_{|k}=g_{ab|k}y^a y^b+2g_{ab}y^a_{|k}y^b=0$, by
hypothesis and Lemma 2.1. Hence $F_{|k}=0$, q.e.d.

We have proved in the end of Section 3 that if $F_{|k}=0$, then
all parallel translations of $\nabla$ are isometries of Minkowski
spaces.

In particular, the holonomy group, let say $H$, of $\nabla$ is
made of isometries of Minkowski spaces.

Let $\calh$ be the Lie algebra of $H$ and an element
$A=(A^a_b)\in\calh$. Then $\exp tA\in H$ and we have
\begin{itemize}
    \item[(4.1)] $\qquad F(x, (\exp tA)y)=F(x,y),\forall x\in M,\forall y\in
    E_x$.
\end{itemize}

This is equivalent with
\begin{itemize}
    \item[(4.1')] $\qquad \dfrac{d}{dt}F(x,(\exp tA)y)|_{t=0}=0$.
\end{itemize}

The linear connection $\nabla$ in $\xi$ corresponds to an
infinitesimal connection in the principal bundle of linear frames
of $\xi$.

By the Holonomy Theorem ([4]) this principal bundle admits an $H-$
structure (a reduction to the Lie subgroup $H$) such that $\Gamma$
becomes an $H-$ connection.

Correspondingly, $\xi$ admits a reduction to $H$ such that
$\nabla$ is an $H-$ connection.

Let be $s_a=\lambda^b_a(x)\varepsilon_a$ a field of frames on the
open set $U_\a$ containing $x$. We think (4.1) and (4,1') in this
frame taking $y=y^a\varepsilon_a=\varepsilon^a s_a$. Thus
$\dfrac{d}{dt}((\exp tA)y)|_{t=0}=A_y=(A^a_b\xi^b)s_a=A^a_b\mu^b_c
y^c\lambda^d_a\varepsilon_d$.

Expanding (4.1') we find
\begin{itemize}
    \item[(4.2)] $\qquad (\dot\p_d F)\lambda^d_a A^a_b\mu^b_c
    y^c=0$,
\end{itemize}
where $\dot\p_d F$, $d=1,...,m$ mean the partial derivatives with
the second set of $m$ variables of $F(\cdot,\cdot)$.

When we put $\nabla_X s_a=X^k \Gamma^{*b}_{ak}s_b$, we necessarily
have $(X^k\Gamma^{*b}_{ak})\in\calh$.

In the natural frame we set
$\nabla_X\varepsilon_a=X^k\Gamma^b_{ak}\varepsilon_b(x)$ and as
before we get
\begin{itemize}
    \item[(4.3)] $\qquad \nabla_X s_a=X^k(\nabla_k \lambda^c_a)\mu^b_c
    s_b(x)$.
\end{itemize}

Thus by comparison it follows (3.2).

Now we write (4.2) for the matrix $(X^k\Gamma^{*a}_{bk})$. We get
\begin{itemize}
    \item[(4.4)] $\qquad (\dot\p_d F)(\nabla_k \lambda^d_a)\mu^a_c
    y^c=0$.
\end{itemize}

Let be $F(x^i,y^a)=F(x^i,\lambda^a_b(x)\xi^b):=f(x,\xi)$.

We show that $f$ does not depend on $x$.

We compute
$$\p_k F=\p_k F+(\dot\p_a F)\p_k (\lambda^a_b (x))\xi^b \buildrel{F_{|k}=0}\over{=}y^b\Gamma^a_{bk}\dot\p_a F+(\dot\p_a
F)\p_k(\lambda^a_b(x))\xi^b=$$ $$=(\dot\p_a
F)(\p_k\lambda^a_c+\Gamma^a_{bk}\lambda^b_c)\xi^c=(\dot\p_a
F)(\nabla_k\lambda^a_c)\mu^c_e y^e=0,$$ by (4.4).

Thus $F(x^i,y^i)=F(x,\lambda^a_b(x)\xi^b)=f(\xi^a)$.

We regard $f$ as a function on $\rr^m$ and it obvious that $f$ is
a Minkowski norm.

Now we show that the holonomy group $H\subset G_f$. Let $T\in H$
with $(T^{*a}_{\ \ b})$ its matrix in the frame $(s_a)$ and
$(T^a_b)$ its matrix in the frame $(\varepsilon_a)$. Then
$T^{*a}_{\ \ b}=\mu^a_c T^c_d\lambda^d_b$. We have
$f(\xi^a)=f(\mu^a_b y^b)=F(x^i,y^a)=F(x^i,T^a_b y^b)=f(\mu^c_a
T^a_b y^b)=f(\mu^c_a\lambda^a_e T^{*e}_d\mu^d_b
y^b)=f(T^{*c}_d\xi^d)$. Thus $T\in G_f$.

As $\xi$ admits an $H-$structure, we may say that it admits also a
$G_f-$ structure. If one reviews the proof of Theorem A it comes
out that the $(\nabla,F)-$ structure induced by this
$G_f-$structure is just that initially given.

\vspace{5mm}

{\bf References}

\begin{enumerate}
    \item Anastasiei, M., {\it The geometry of Berwald - Cartan
    spaces.} To appear. Communicated at the Workshop: Differential
    Geometry and its Applications, June 24-30, Thessaloniki,
    Greece.
    \item Ichijy\~{o}, Y., {\it Finsler manifolds modelled on a Minkowski
    space}. J. Math. Kyoto Univ. (JMKYAZ) 16-3 (1976), 639-652.
    \item Ichijy\~{o}, Y., {\it Finsler Manifolds with a Linear
    Connection}. J. Math. Tokushima Univ., vol. 10 (1976), 1-11.
    \item Kobayashi, S., Nomizu, K., {\it Foundations of Differential
    Geometry}, vol. I, Intersciences Publishers, 1963.
    \item Miron, R., Hrimiuc, D., Shimada, H., Sab\u{a}u, V.S., {\it The Geometry of Hamilton and Lagrange
    Spaces}. Kluwer Academic Publishers, FTPH 118, 2001.
\end{enumerate}

\vspace{5mm}





\newpage

\def\bigi{\big|}
\def\scri{\scriptstyle{|}}

\runningauthor={M. ANASTASIEI}
\runningtitle={SOME NEW PROPERTIES OF BERWALD - CARTAN SPACES}
\noindent
\baselineskip 8pt
\noindent{\footnotesize{Libertas Mathematica\hfill\break 24, 3-10 (2004)}}
\vskip 2cm
\baselineskip 11.5pt plus .15pt
\centerline{\bf\Large SOME NEW PROPERTIES OF}
\vskip .2cm
\centerline{\bf\Large BERWALD - CARTAN SPACES}
\vskip .5cm
\centerline{\bf by Mihai ANASTASIEI{\footnote[1]{This paper was partially
supported by CNCSIS-Romania}}}
\vskip .5cm

\begin{abstract}
A manifold endowed with a regular Hamiltonian which is
2--homogeneous in momenta was called a {\it Cartan space}.The
geometry of regular Hamiltonians as smooth functions on the
cotangent bundle is mainly due to R. Miron and it is now
systematically described in the monograph \cite{MHSS}. An
interesting particular class of Cartan spaces is given by the
so--called Berwald--Cartan spaces. In this paper some new
properties of the Berwald--Cartan
spaces are proved.\\
{\bf 2000 Mathematics Subject Classification : 53C60}\\
{\bf Phrases and key words:} cotangent bundle, homegeneous Hamiltonians.
\end{abstract}

\setcounter{section}{0}
\section*{Introduction}

Analytical Mechanics and some theories in Physics brought into
discussion regular Lagrangians and their geometry, \cite{MA}. A
regular Lagrangian which is 2-homogeneous in velocities is nothing
but the square of a fundamental Finsler function and its geometry
is Finsler geometry. This geometry was developed since 1918 by P.
Finsler, E. Cartan, L. Berwald and many others, see \cite{M} and
the most recent graduate text \cite{BCS}. But in Mechanics and
Physics there exists also regular Hamiltonians whose geometry is
also useful. This geometry is mainly due to R. Miron, \cite{Mi}, and it is now
systematically presented in the monograph \cite{MHSS}. A manifold
endowed with a regular Hamiltonian which is 2-homogeneous in
momenta was called a Cartan space. The notion of Cartan space was introduced
by R. Miron in \cite{Mi}. A particular and interesting
class of Cartan spaces is given by the so--called Berwald--Cartan
spaces, shortly $BC$-spaces. The geometry of the $BC$-spaces can
be found in \cite{MHSS}, Chs. 6-7. Our purpose is to prove some
new properties of these spaces. A Cartan space is a pair $(M,K)$
for $M$ a smooth manifold and $K$ a regular Hamiltonian which is
$2$-homogeneous in momenta. A $BC$ space is defined as a Cartan
space whose Chern--Rund connection coefficients of the canonical
metrical connection do not depend on momenta, that is,
$H^i_{jk}(x,p)=H^i_{jk}(x).$ For a Cartan space the pair
$(T^*_xM,K(x,p))$ for any fixed $x\in M$ is a Minkowski space. We
prove (Theorem 3.2) that for $BC$ spaces the Minkowski spaces
$(T^*_xM,K(x,p))$ are all linearly isometric to each other.
Noticing that the functions $H^i_{jk}(x)$ defines a symmetric
linear connection $\nabla$ on $M$ we prove (Theorem 3.3) that
$\nabla$ is metrizable, that is, there exists a Riemannian metric
on $M$ whose Levi--Civita connection is $\nabla$. These proofs are
presented in Section 3. Some preliminaries from the geometry of
cotangent bundle are given in Section 1, and Section 2 contains
necessary facts from the geometry of Cartan spaces.

\section{Preliminaries}

Let $M$ be an $n$-dimensional $C^\infty$ manifold and $\tau^*:T^*M\to
M$ its cotangent bundle. If $(x^i)$ are local coordinates on $M$,
then $(x^i,p_i)$ will be taken as local coordinates on $T^*M$ with
the momenta $(p_i)$ provided by $p=p_idx^i$ where $p\in T^*_xM$,
$x=(x^i)$ and $(dx^i)$ is the natural basis of $T^*_xM$. The
indices $i,j,k...$ will run from 1 to $n$ and the Einstein
convention on summation will be used. A change of coordinates
$(x^i,p_i)\to(\widetilde x^i,\widetilde p_i)$ on $T^*M$ has the form
$$\begin{array}{ll}
\widetilde x^i=\widetilde x^i(x^1,...,x^n),&{\rm rank}\left(
\displaystyle\frac{{\partial}\widetilde x^i}{{\partial} x^j}\right)=n\vspace*{1,5mm}\\
\widetilde p_i=\displaystyle\frac{{\partial} x^j}{{\partial} \widetilde x^i}(\widetilde x)p_j,\end{array}\leqno(1.1)$$
where $\left(\displaystyle\frac{{\partial} x^j}{{\partial} \widetilde x^i}\right)$ is the inverse of the
Jacobian matrix $\left(\displaystyle\frac{{\partial}\widetilde x^j}{{\partial} x^k}\right).$

Let $\left({\partial}_i:=\displaystyle\frac{\partial}{{\partial} x^i},\ {\partial}^i:=\displaystyle\frac{\partial}{{\partial} p_i}\right)$ be
the natural basis in $T_{(x,p)}T^*M$. The change of coordinates
(1.1) produces
$$\begin{array}{l}
{\partial}_i=({\partial}_i\widetilde x^j)\widetilde{\partial}_j+({\partial}_i\widetilde p_j)\widetilde{\partial}^j,\vspace*{1,5mm}\\
\widetilde{\partial}^i=({\partial}_j\widetilde x^i){\partial}^j.\end{array}\leqno(1.2)$$
The natural cobasis $(dx^i,dp_i)$ from $T^*_{(x,p)}T^*M$
transforms as follows.
$$d\widetilde x^i=({\partial}_j\widetilde x^i)dx^j,\ d\widetilde p_i=\displaystyle\frac{{\partial} x^j}{{\partial} \widetilde
x^i}\ dp_j+\displaystyle\frac{{\partial}^2x^j}{{\partial}\widetilde x^i{\partial}\widetilde x^k}\
p_j\,dx^k.\leqno(1.3)$$

The kernel $V_{(x,p)}$ of the differential
$d\tau^*:T_{(x,p)}T^*M\to T_xM$ is called the {\it vertical}
subspace of $T_{(x,p)}T^*M$ and the mapping $(x,p)\to V_{(x,p)}$
is a regular distribution on $T^*M$ called the {\it vertical
distribution}. This is integrable with the leaves $T^*_xM,$ $x\in
M$ and is locally spanned by $({\partial}^i).$ The vector field
$C^*=p_i{\partial}^i$ is called the Liouville vector field and
${\omega}=p_idx^i$ is called the Liouville $1$-form on $T^*M$. Then
$d{\omega}$ is the canonical symplectic structure on $T^*M$. For an
easier handling of the geometrical objects on $T^*M$ it is usual
to consider a supplementary distribution to the vertical
distribution, $(x,p)\to N_{(x,p)}$, called the {\it horizontal
distribution} and to report all geometrical objects on $T^*M$ to
the decomposition
$$T_{(x,p)}T^*M=N_{(x,p)}\oplus V_{(x,p)}.\leqno(1.4)$$
The pieces produced by the decomposition (1.4) are called
$d$--geometrical objects $(d$ is for distinguished) since their
local components behave like geometrical objects on $M$, although
they depend on $x=(x^i)$ and momenta  $p=(p_i).$

The horizontal distribution is taken as being locally spanned by
the local vector fields
$$\delta_i:={\partial}_i+N_{ij}(x,p){\partial}^j,\leqno(1.5)$$
and for a change of coordinates (1.1), the condition
$$\delta_i=({\partial}_i\widetilde x^j)\widetilde\delta_j\mbox{ for
}\widetilde\delta_j:=\widetilde{\partial}_j+\widetilde N_{jk}(\widetilde x,\widetilde p)\widetilde{\partial}^k,\leqno(1.6)$$
is equivalent with
$$\widetilde N_{ij}(\widetilde x,\widetilde p)=\displaystyle\frac{{\partial} x^s}{{\partial}\widetilde x^i}\
\displaystyle\frac{{\partial} x^r}{{\partial}\widetilde x^j}\ N_{sr}(x,p)+
\displaystyle\frac{{\partial}^2x^r}{{\partial}\widetilde x^i{\partial}\widetilde x^r}\ p_r.\leqno(1.7)$$
The horizontal distribution is called also a {\it nonlinear
connection} on $T^*M$ and the functions $(N_{ij})$ are called the
local coefficients of this nonlinear connection. It is important
to note that any regular hamiltonian on $T^*M$ determines a
nonlinear connection whose local coefficients verify
$N_{ij}=N_{ji}.$ The basis $(\delta_i,{\partial}^i)$ is adapted to the
decomposition (1.4). The dual of it is $(dx^i,\delta p_i)$, for
$\delta p_i=dp_i-N_{ji}dx^j$ and then $\delta\widetilde p_i=\displaystyle\frac{{\partial}
x^j}{{\partial}\widetilde x^i}\ \delta p_j.$

\section{Cartan spaces}
A {\it Cartan structure} on $M$ is a function $K:T^*M\to[0,\infty)$
with the following properties:
\begin{enumerate}
\item $K$ is $C^\infty$ on $T^*M\setminus0$ for $0=\{(x,0),\ x\in
M\},$
\item $K(x,\lambda p)=\lambda K(x,p)$ for all $\lambda>0,$
\item The $n\times n$ matrix $(g^{ij})$, where
$g^{ij}(x,p)=\displaystyle\frac12\ {\partial}^i{\partial}^jK^2(x,p),$ is positive--definite
at all points of $T^*M\setminus0.$
\end{enumerate}
We notice that in fact $K(x,p)>0$, whenever $p\ne0.$

\medskip\noindent{\bf Definition 2.1.} The pair $(M,K)$ is called a
{\it Cartan space}.

\medskip\noindent{\it Example.} Let $(\gamma_{ij}(x))$ be the matrix of the
local coefficients of a Riemannian metric on $M$ and
$(\gamma^{ij}(x))$ its inverse. Then
$K(x,p)=\sqrt{\gamma^{ij}(x)p_ip_j}$ gives a Cartan structure.
Thus any Riemannian manifold can be regarded as a Cartan space.
More examples can be found in Ch. 6 of \cite{MHSS}.

We put $p^i=\displaystyle\frac12\ {\partial}^iK^2$ and $C^{ijk}=-\displaystyle\frac14\
{\partial}^i{\partial}^j{\partial}^kK^2.$ The properties of $K$ imply
$$\begin{array}{l}
p^i=g^{ij}p_j,\ p_i=g_{ij}p^j,\ K^2=g^{ij}p_ip_j=p_ip^j,\vspace*{1,5mm}\\
C^{ijk}p_k=C^{ikj}p_k=C^{kij}p_k=0.\end{array}\leqno(2.1)$$
One considers the {\it formal Christoffel symbols}
$$\gamma^i_{jk}(x,p):=\displaystyle\frac12\
g^{is}({\partial}_kg_{js}+{\partial}_jg_{sk}-{\partial}_sg_{jk})\leqno(2.2)$$ and the
contractions $\gamma^\circ_{jk}(x,p):=\gamma^i_{jk}(x,p)p_i,$
$\gamma^\circ_{j\circ}:=\gamma^i_{jk}p_ip^k.$ Then the functions
$$N_{ij}(x,p)=\gamma^\circ_{ij}(x,p)-\displaystyle\frac12\
\gamma^\circ_{h\circ}(x,p){\partial}^hg_{ij}(x,p),\leqno(2.3)$$ verify (1.7).
In other words, these functions define a nonlinear connection on
$T^*M$. This nonlinear connection was discovered by R. Miron, \cite{Mi}. Thus a decomposition (1.4) holds. From now on we shall use
only the nonlinear connection given by (2.3).

A linear connection $D$ on $T^*M$ is said to be an
$N$--{\it linear connection} if
\begin{itemize}
\item[$1^\circ$] $D$ preserves by parallelism the distributions
$N$ and $V$,
\item[$2^\circ$] $D\theta=0,$ for $\theta=\delta p_i\wedge dx^i.$
\end{itemize}
One proves that an $N$-linear connection can be represented in the
adapted basis $(\delta_i,{\partial}^i)$ in the form
$$\begin{array}{ll}
D_{\delta_j}\delta_i=H^k_{ij}\delta_j,&
D_{\delta_j}{\partial}^i=-H^i_{kj}{\partial}^k,\vspace*{1,5mm}\\
D_{\partial^j}\delta_i=V_i^{kj}\delta_k,&
D_{\partial^j}{\partial}^i=-V_k^{ij}\delta^k,\end{array}\leqno(2.4)$$
where $V_i^{kj}$ is a $d$--tensor field and $H^k_{ij}(x,p)$ behave
like the coefficients of a linear connection on $M$. The functions
$H^k_{ij}$ and $V_i^{kj}$ define operators of $h$--covariant and
$v$-covariant derivatives in the algebra of $d$-tensor fields,
denoted by $\scri k$ and $\bigi^k$, respectively.
For $g^{ij}$ these are given by
$$\begin{array}{l}
g^{ij}{}_{\scri k}=\delta_kg^{ij}+g^{sj}H^i_{sk}+g^{is}H^j_{sk},\vspace*{1,5mm}\\
g^{ij}\bigi^k={\partial}^kg^{ij}+g^{sj}V_s^{ik}+g^{is}V_s^{jk}.\end{array}\leqno(2.5)$$
An $N$-linear connection given in the adapted basis $(\delta_i,{\partial}^j)$
as $D\Gamma(N)=(H^i_{jk},V_j^{ik})$ is called {\it metrical} if
$$g^{ij}{}_{\scri k}=0,\ \ g^{ij}\bigi^k=0.\leqno(2.6)$$
One verifies that the $N$-linear connection
$C\Gamma(N)=(H^i_{jk},C^{jk}_i)$ with
$$\begin{array}{l}
H^i_{jk}=\displaystyle\frac12\ g^{is}(\delta_jg_{sk}+\delta_kg_{js}-\delta_sg_{jk}),\vspace*{1,5mm}\\
C^{jk}_i=-\displaystyle\frac12\
g_{is}({\partial}^jg^{sk}+{\partial}^kg^{sj}-{\partial}^sg^{jk})=g_{is}C^{sjk},\end{array}\leqno(2.7)$$
is metrical and its $h$-torsion $T^i_{jk}:=H^i_{jk}-H^i_{kj}=0,$
$v$-torsion $S^{jk}_i:=C^{jk}_i-C^{kj}_i=0$ and the deflection
tensor $\Delta_{ij}=N_{ij}-p_kH^k_{ij}=0.$ Moreover, it is unique
with these properties. This is called the canonical metrical
connection of the Cartan space $(M,K).$ It has also the following
properties:
$$\begin{array}{llll}
K_{\scri j}=0,&K\bigi^j=\displaystyle\frac{p^j}K,&
K^2{}_{\scri j}=0,&K^2\bigi^j=2p^j,\vspace*{1,5mm}\\
p_{i\scri j}=0,&p_i\bigi^j=\delta^j_i,
&p^i{}_{\scri i}=0,&p^i\bigi^j=g^{ij}.\end{array}\leqno(2.8)$$
Besides $C\Gamma(N)$ one may consider on $T^*M$ three other
important $N$-linear connection which are partially or not at all
metrical: Chern--Rund connection $CR\Gamma(N)=(H^i_{jk},0)$, the Hashiguchi connection
$H\Gamma(N)=({\partial}^iN_{jk},C^{kj}_i)$ and
 the Berwald connection $B\Gamma(N)=({\partial}^iN_{jk},0)$.

\section{Berwald--Cartan spaces}

Let $C^n=(M,K)$ be a Cartan space with the canonical metrical
connection $C\Gamma(N)=(H^i_{jk},C^{jk}_i)$ given by (2.7).

\medskip\noindent{\bf Definition 3.1.} The Cartan space $C^n$ is called a {\it
Berwald--Cartan space}, shortly a $BC$ space, if the connection
coefficients $H^i_{jk}$ do not depend on momenta, that is,
$H^i_{jk}(x,p)=H^i_{jk}(x).$\bigskip

In \cite{MHSS}, by direct methods or using the duality between
Finsler and Cartan spaces given by the Legendre map, one proves

\medskip\noindent{\bf Theorem 3.1.} {\it The following assertions are
equivalent:
\begin{itemize}
\item[$1^\circ$] The Cartan space $C^n$ is a $BC$ space,
\item[$2^\circ$] The coefficients $B^i_{jk}={\partial}^iN_{jk}$ of the
Berwald connection are functions of position only, that is $B^i_{jk}(x,p)=B^i_{jk}(x),$
\item[$3^\circ$] The curvature $P_j{}^i{}_k{}^h:=\dot{\partial}^hB^i_{jk}$
of the Berwald connection vanishes.
\item[$4^\circ$] $C^{ijk}{}_{\scri h}=0.$
\end{itemize}}

For the Cartan space $C^n=(M,K)$, the function
$K_x:=K(x,\cdot):T^*_xM\to\rr$ is a {\it Minkowski norm} for every
$x\in M$. Thus we have the Minkowski spaces $(T^*_xM,K_x)$, $x\in
M.$ For $BC$ spaces, the following theorem holds.

\medskip\noindent{\bf Theorem 3.2.} {\it Let $(M,K)$ be a $BC$ space. Whenever
$M$ is connected the Minkowski spaces $(T^*_xM,K_x)$ are all
linearly isometric to each other.}

\medskip\noindent{\it Proof.} Let ${\omega}={\omega}_idx^i$ an $1$-form and $v=v^j{\partial}_j$
a vector field on $M$. Using the connection coefficients
$H^i_{jk}(x)$ we may define a covariant derivative of ${\omega}$ in the
direction of $v$ as follows:
$\nabla_v{\omega}=v^k({\partial}_k{\omega}_i-H^j_{ik}{\omega}_j)dx^i.$

We restrict ${\omega}$ to a curve $c:t\to x(t),$ $t\in\R,$ on $M$,
define the covariant derivative of ${\omega}$ along $c$ by
$\displaystyle\frac{\nabla{\omega}}{dt}=\left[\displaystyle\frac{d{\omega}_i}{dt}-H^j_{ik}{\omega}_j\
\displaystyle\frac{dx^k}{dt}\right]dx^i$ and we say that ${\omega}$ is parallel
along $c$ if $\displaystyle\frac{\nabla{\omega}}{dt}=0$. Let us estimate
$\displaystyle\frac{dK^2(x(t),{\omega}(t))}{dt}.$ We write the equality
$K^2(x,p)=g^{ij}(x,p)p_jp_j$ for $(x(t),{\omega}(t))$ and we obtain
that along the curve $c$:
$\displaystyle\frac{dK^2}{dt}=\displaystyle\frac{dg^{ij}}{dt}{\omega}_i{\omega}_j+2g^{ij}{\omega}_i\
\displaystyle\frac{d{\omega}_j}{dt}.$ But $\displaystyle\frac
d{dt}(g^{ij})=(\delta_kg^{ij})\displaystyle\frac{dx^k}{dt}+({\partial}^kg^{ij})\displaystyle\frac{\delta p_k}{dt}$
and using $g^{ij}{}_{\scri k}=0$ as well as the last equation
(2.1) one gets:
$$\displaystyle\frac{dK^2}{dt}=2g^{ij}{\omega}_i\left(\displaystyle\frac{d{\omega}_j}{dt}-H^s_{jk}{\omega}_s\
\displaystyle\frac{dx^k}{dt}\right).$$
From here we read

\medskip\noindent{\bf Lemma 3.1.} {\it If the $1$-form ${\omega}$ is parallel along
the curve $c:t\to x(t)$, then the function $K(t):=K(x(t),{\omega}(t))$ is constant
along the curve $c$.}\bigskip

Let $x,y$ be points of $M$ joined by a curve $c:[0,1]\to M$ such that
$c(0)=x,$ $c(1)=y.$ Let be $\alpha\in T^*_xM$. We consider the
unique solution ${\omega}=({\omega}_i)$ of the system of linear ordinary
differential equations $\displaystyle\frac{d{\omega}_i}{dt}-H^j_{ik}{\omega}_j\
\displaystyle\frac{dx^k}{dt}=0$ with the initial condition ${\omega}(0)=\alpha$ and
we associate to $\alpha$ the element $\alpha'={\omega}(1)$ of
$T^*_y M.$ The mapping $T^*_xM\to T^*_yM$ given by
$\alpha\to\alpha'$ is a linear isomorphism. By Lemma 3.1,
$K(x(t),{\omega}(t))$ has the same values at $t=0$. Hence
$K_x(\alpha)=K_y(\alpha').$ This means that the Minkowski spaces
$(T^*_xM,K_x)$ and $(T^*_yM,K_y)$ are linearly isometric for every
$x,y\in M,$ q.e.d.

Another interesting property of $BC$ spaces is as follows.

The connection coefficients $H^i_{jk}(x,p)=H^i_{jk}(x)$ define a
symmetric linear connection $\nabla$ on $M$ and it happens that
this is {\it metrizable}, that is, there exists on $M$ a
Riemannian metric $h$ such that $\nabla$ is the Levi--Civita
connection associated to it. This $h$ is not unique.

We prove this fact by adapting an idea of Z.I. Szab\'o \cite{S}.
The duality with Finsler spaces is not used.

\bigskip\noindent{\bf Theorem 3.3.} {\it Let $C^n=(M,K)$ be a $BC$ space with
$M$ connected and $\nabla$ the symmetric linear connection on $M$
of local coefficients $H^i_{jk}(x,p)=H^i_{jk}(x).$ Then there exists a
Riemannian metric $h$ on $M$ such that $\nabla$ is the
Levi--Civita connection of it.}

\medskip\noindent{\it Proof.} Let be the Minkowski space
$(T^*_{x_0}M,K_{x_0})$ for a fixed $x_0\in M.$ Then
$S_{x_0}=\{{\omega}\mid K_{x_0}({\omega})=1\}$ is a compact subset of
$T^*_{x_0}M.$ Let $G$ be the group of all linear isomorphisms of
$T^*_{x_0}M$ that preserve $S_{x_0}.$ This $G$ is a compact Lie
group. It contains as a subgroup the holonomy group $H_{x_0}$
defined by $(H^i_{jk}(x))$ according to Lemma 3.1. In general,
$H_{x_0}$ is not compact.

Let $<,>$ be {\it any} inner product in $T^*_{x_0}M.$ Define a new
inner product on $T^*_{x_0}M$ by
$$h_{x_0}({\varphi},{\omega})=\displaystyle\frac1{{\rm vol}(G)}\int_G<a{\varphi},a{\omega}>\mu_G,\
{\varphi},{\omega}\in T^*_{x_0}M,\leqno(3.1)$$
for $a\in G$, where $\mu_G$ denotes the bi--invariant Haar measure
on $G$. It results $h_{x_0}(b{\varphi},b{\omega})=h_{x_0}({\varphi},{\omega})$ for every
$b\in G$ (from the properties of $\mu_G)$, that is $h_{x_0}$ is
$G$-invariant. In particular, $h_{x_0}$ is $H_{x_0}$-invariant.

Let now any $x\in M$ and a curve $c:t\to c(t)$ joining $x$ with
$x_0$, $c(0)=x,$ $c(1)=x_0.$ Denote by $P_c:T^*_xM\to T^*_{x_0}M$
the parallel transport of covectors defined by $H^i_{jk}(x).$ For
every ${\varphi}\in T^*_xM$, $P_c({\varphi})={\omega}(1)\in T^*_{x_0}M,$ where
${\omega}=({\omega}_i)$ is the unique solution of the system of linear
differential equations
$$\displaystyle\frac{d{\omega}_i}{dt}-H^i_{jk}{\omega}_j\ \displaystyle\frac{dx^k}{dt}=0,\mbox{
with }{\omega}(0)={\varphi}.\leqno(3.2)$$
In the proof of Theorem 3.2 we have seen that $P_c$ is a linear
isometry of Minkowski spaces. We define an inner product on
$T^*_xM$ by
$$h_x({\varphi},\psi)=h_{x_0}(P_c{\varphi},P_c\psi),\ {\varphi},\psi\in
T^*_{x_0}M.\leqno(3.3)$$

\medskip\noindent{\bf Lemma 3.2.} {\it $h_x$ does not depend on the curve
$c$.}\bigskip

Indeed, if $\widetilde c$ is another curve joining $x$ and $x_0$, denote
by $c_-$ the reverse of $c$ and consider the loop $\widetilde c\circ
c_-$. Then $P_{\widetilde c\circ c_-}\in H_{x_0}$ and from the
$H_{x_0}$-invariance of $h_{x_0}$, that is,
$h_{x_0}(P_{\widetilde c\circ c_-}{\varphi}$, $P_{\widetilde c\circ
c_-}\psi)=h_{x_0}({\varphi},\psi)$ we get
$h_{x_0}(P_{\widetilde c}{\varphi},P_{\widetilde c}\psi)=h_{x_0}(P_c{\varphi},P_c\psi)$ as
we claimed.

The mapping $x\to h_x:T^*_xM\times T^*_xM\to R$ is smooth since
$P_c$ smoothly depends on $x$, according to a general result
regarding the dependence of solution of system of differential
equations by initial data. Thus we have constructed a Riemannian
metric $h$ in the cotangent bundle of $M$.

The connection coefficients $(H^i_{jk}(x))$ define a linear
connection $\nabla$ in the cotangent bundle as follows:
$$\nabla:{\cal X}(M)\times\Gamma(T^*M)\to\Gamma(T^*M),\
(X,{\omega})\to\nabla_X{\omega}=X^k\left(
\displaystyle\frac{{\partial}{\omega}_i}{{\partial} x^k}-H^j_{ik}{\omega}_j\right)dx^i$$
and the operator $\nabla_X,$ $X\in{\cal X}(M)$, extends to the
tensorial algebra of the cotangent bundle. For instance, if we
regard $h$ as a section in the vector bundle $L^s_2(T^*M,\R)$,
then we have
$$(\nabla_Xh)({\varphi},\psi)=X(h({\varphi},\psi))-h(\nabla_X{\varphi},\psi)-h({\varphi},\nabla_X\psi).\leqno(3.4)$$

\medskip\noindent{\bf Lemma 3.3.} {\it $\nabla_Xh=0,\ X\in{\cal X}(M).$}

\medskip\noindent{\it Proof.} We choose a basis $({\varphi}_i(x))$ in $T^*_xM.$ It
suffices to show that\break $(\nabla_Xh)({\varphi}_i(x),{\varphi}_j(x))=0.$ Let be the vector
$X=\left.\displaystyle\frac{dc}{dt}\right|_\circ$ tangent to a curve $c$ starting from
$x\in M$ at $t=0.$ We parallel translate ${\varphi}_i(x)$ along $c$ and
we obtain a field of basis ${\varphi}_i(t)$ along $c$. The general
formula
$$\displaystyle\frac{\nabla
h}{dt}({\varphi},w)=\displaystyle\frac{dh({\varphi},\psi)}{dt}-h\left(
\displaystyle\frac{\nabla{\varphi}}{dt},\psi\right)-h\left({\varphi},\displaystyle\frac{\nabla\psi}{dt}\right),$$
gives
$$\displaystyle\frac{\nabla h}{dt}({\varphi}_i(x),{\varphi}_j(x))=
\left.\displaystyle\frac{dh({\varphi}_i,{\varphi}_j)}{dt}\right|_{t=0}$$
because of $\displaystyle\frac{\nabla{\varphi}_i}{dt}=0.$

Now we show that $h({\varphi}_i(t),{\varphi}_j(t))$ does not depend on $t$.

Indeed,
$h_{c(t)}({\varphi}_i(t),{\varphi}_j(t))=h_{x_0}(P_{{\varphi}_i},P_{{\varphi}_j}),$ where
$P$ is the parallel translation from $T^*_{c(t)}M$ to $T_{x_0}M.$
This $P$ may be thought as the composition of a parallel
translation $P_2$ from $T^*_{c(t)}M$ to $T^*_xM$ and of a parallel
translation $P_1$ from $T^*_xM$ to $T^*_{x_0}M.$
We have
$h_{c(t)}({\varphi}_i(t),{\varphi}_j(t))=
h_{x_0}((P_2\circ P_1){\varphi}_i,(P_2\circ P_1){\varphi}_j)=
h_{x_0}(P_1{\varphi}_i,P_2{\varphi}_j)=h_x({\varphi}_i(x),{\varphi}_j(x)).$ Hence
$h_{c(t)}({\varphi}_i(t),{\varphi}_j(t))$ does not depend on $t$, as we
claimed.

This fact ends the proof of Lemma 3.3.

To end the proof of Theorem, we take the covariant part of $h$ as
a section in the vector bundle $L^s_2(TM,\R)$ and so we get a
Riemannian metric on $M$, denoted with the same letter $h$. The
operator $\nabla_X$ acts also on vector fields on $M$ by the rule
$\nabla_XY=X^k\left(\displaystyle\frac{{\partial} Y^i}{{\partial} x^k}+H^i_{jk} Y^j\right)$ for $Y=Y^i\
\displaystyle\frac{\partial}{{\partial} x^i}$ and $(X,Y)\to\nabla_XY$ gives a linear
connection on $M$ such that $\nabla_Xh=0.$ As $\nabla$ has no
torsion, it coincides with the Levi--Civita connection of $h$, q.e.d.

\medskip\noindent{\it Remark.}  An alternative way to prove Lemma 3.3 is
to prove first that
$\displaystyle\frac{\nabla h}{dt}({\varphi},\psi)=\displaystyle\lim_{t\to0}\
\displaystyle\frac{h(P_c\varphi,P_c\psi)-h({\varphi},\psi)}t,$
where $P_c$ is the parallel translation from $c(0)$ to $c(t)$.

\bigskip\noindent{\bf Acknowledgements.} The author is grateful to Professor Radu Miron for stimulating discussions during the preparation of this work.

\noindent{\footnotesize{Faculty of Mathematics,\\
University ``Al.I. Cuza'' Ia\c si,\\
6600, Ia\c si, Romania}}

\newpage

\def\llr{\Longleftrightarrow }
\def\frax#1#2{\displaystyle\frac{#1}{#2}}
\def\derp#1#2{\displaystyle\frac{\partial#1}{\partial#2}}
\def\derpp#1{\displaystyle\frac{\partial}{\partial#1}}
\def\bld#1#2{{\buildrel{#1}\over{#2}}}
\def\scri{\scriptstyle{|}}
\def\a{{\alpha}}
\def\b{{\beta}}
\def\de{{\delta}}
\def\G{{\Gamma}}
\def\g{{\gamma}}
\def\pp{{\partial}}
\def\wt{\widetilde}
\def\dd{\displaystyle}
\def\n{\noindent }
\def\bigi{\big|}

\def\wh#1{\widehat{#1}}
\def\wt#1{\widetilde{#1}}

\runningauthor={M. ANASTASIEI}
\runningtitle={FINSLER VECTOR BUNDLES. METRIZABLE CONNECTIONS}
\noindent
\baselineskip 8pt
\noindent{\footnotesize{Period. Math. Hungar.}}
\hfill\break
{\footnotesize{48, no. 1-2, 83–-91, 2004}}
\vskip 2cm
\baselineskip 11.5pt plus .15pt
\centerline{\bf\Large FINSLER VECTOR BUNDLES.}
\vskip .2cm
\centerline{\bf\Large METRIZABLE CONNECTIONS}
\vskip .5cm
\centerline{\bf {\footnotesize{BY}}}
\vskip .5cm
\centerline{\bf {\footnotesize{M. ANASTASIEI}\footnote{Lecture given at the Workshop on Finsler Geometry and its Applications, August 11-15,2003, Debrecen, Hungary}}}
\vskip 1cm

\centerline{\footnotesize{Dedicated to Prof. Dr. Lajos Tam\'{a}ssy at his 80th anniversary}}

\begin {abstract}
A vector bundle $\xi = (E,\pi,M)$ of rank $m$ is called a
Finsler vector bundle if $E$ is endowed with a continuous,
positive function $F$ which is smooth on $ E\backslash 0$,
positively homogeneous of degree $1$ in fibre variables and whose
Hessian is positive definite. Then the fibres $E_x,x\in M,$ of
$\xi$ are Minkowski spaces with the Minkowski norm $F(x, )$.

A nonlinear connection $N$ in $\xi$ induces a linear connection in the vertical bundle over $E$ (Berwald connection) and an operator $|_k$ of $h-$covariant derivative. We say that $N$ is compatible with $F$ if $F_{|k}=0$ and in this case we show that the parallel translations of $N$ preserve the norms $F(x, )$. Next we consider the case when the coefficients of the Berwald connection do not depend of the fibre variables and we prove that the linear connection in $\xi $ defined by these coefficients is metrizable. As a corollary a metrizability condition for any linear connection in the Finsler vector bundle $\xi $ is provided.  \\

\noindent{\bf Mathematics Subject Classification} : Primary 53C60; Secondary 53C05.\\
{\bf Key words and phrases} : Finsler vector bundles, linear connections, metrizability
\end{abstract}

\setcounter{section}{0}
\section*{Introduction}
The notion of Finsler function can be considered not only for tangent bundles but also for any vector bundle and the notion of
Finsler vector bundle is obtained. A vector bundle $\xi=(E,\pi,M)$ of rank $m$ is called a {\bf Finsler vector bundle} if $E$ is
endowed with a continuous, positive function $F$ which is smooth on $ E\backslash 0$, positively homogeneous of degree $1$ in fibre
variables and whose Hessian is positive definte. Any Riemannian metric in $\xi $ defines a Finsler function and Finsler functions
of Randers type can be considered. When $M$ is a paracompact manifold, the vector bundle $\xi$ can be endowed with a nonlinear
connection $N$.This defines a linear connection in the vertical bundle over $E$ called the Berwald connection associated to $N$.
We use it in Section 2 in order to define two kinds of compatibility between $F$ and $N$ that coincide when the Berwald
connection does not depend on variables from fibres. In this case the Berwald connection may be thought as a linear connection
$\nabla$ in $\xi$ and in Section 3 we show that $\nabla$ is a metrizable connection, that is there exists a Riemannian metric
$h$ in $\xi$ such that $\nabla h = 0$. As a corollary we point out a metrizability condition for any linear connection in the
Finsler vector bundle $\xi$. For the problem of metrizability of linear connections we refer to the paper \cite{5}, \cite{6} by L.
Tamassy as well as to our papers \cite{1} and \cite{2}.

\section{Finsler vector bundles}

Let $\xi=(E,p,M),\ p:E\to M$, be a vector bundle of rank $m$. Here $M$ is a smooth i.e. $C^\infty$ manifold of dimension $n$. The type fibre is $\R^m$ and $E$ is a smooth manifold of dimension $n+m$. The projection $p$ is a smooth submersion. Let $(U,(x^i))$ be a local chart on $M$ and let $\varepsilon_a(x),$ $x\in U$, be a field of local sections of $\xi$ over $U$. Then every section $A$ of $\xi$ over $U$ takes the form $A=A^a(x)\varepsilon_a(x),$ $x\in U$, and an element $u\in p^{-1}(x):=E_x$ can be written as $u=y^a\varepsilon_a(x)$, $(y^a)\in\R^m.$ The indices $i,j,k,...$ will range over $\{1,2,...,n\}$ and the indices $a,b,c,...$ will take their values in $\{1,2,...,m\}$. The convention on summation over repeated indices of the same kind will be used.

The local coordinates on $p^{-1}(U)$ will be $(x^i,y^a)$ and a change of coordinates $(x^i,y^a)\to(\wt x^i,\wt y^a)$ on $U\cap
\wt U\ne\emptyset$ has the form
$$\begin{array}{l}
\wt x^i=\wt x^i(x^1,...,x^n),\ {\rm rank}\left(\derp{\wt x^i}{x^j}\right)=n,\vspace*{1,5mm}\\
\wt y^a=M^a_b(x)y^b,\ {\rm rank}(M^a_b(x))=m,\ \ \forall  x\in U\cap\wt U.\end{array}\leqno(1.1)$$

We denote by ${\cal F}(M),{\cal F}(E)$ the ring of real functions on $M$ and $E$, respectively and by ${\cal X}(M),$ resp. $\Gamma (E)$, ${\cal X}(E)$ the module of sections of the tangent bundle of $M$, resp. of the bundle $\xi$ and of the tangent bundle of $E$. On $U$, the vector fields $(\partial_k:=\frax{\partial}{\partial x^k})$ provide a local basis for ${\cal X}(U).$

Let $\xi^*=(E^*,p^*,M)$ be the dual of the vector bundle $\xi$. We take as local basis of $\Gamma(E^*)$ on $U_\a$, the sections
$\theta^a:U\to p^{*-1}(U),$ $x\to\theta^a(x)\in E^*_x$ such that $\theta^a(\varepsilon_b(x))=\delta^a_b$. A section $\beta$ of $\xi^*=(E^*,p^*,M)$ will take the form $\beta(x)=\beta_{a}\theta^a$.

Next, we may consider the tensor bundle of type $(r,s)$, denoted as ${\cal T}^r_s(E):=E\underbrace{\otimes\cdots\otimes}_r E\otimes
E^*\underbrace{\otimes\cdots\otimes}_sE^*$ over $M$ and its sections. For $g\in\Gamma(E^*\otimes E^*)$ we have the local
representation $g=g_{ab}(x)\theta^a\otimes\theta^b$. As $E^*\otimes E^*\cong L_2(E,\R),$ we may regard $g$ as a smooth
mapping $x\to g(x):E_x\times E_x\to\R$ with $g(x)$ a bilinear mapping given by $g(x)(s_a,s_b)=g_{ab}(x)$.

If the mapping $g(x)$ is symmetric i.e. $g_{ab}=g_{ba}$ and positive-definite i.e. $g_{ab}(x)\zeta^a\zeta^b>0$ for every
$0\ne(\zeta^a)\in\R^m,$ one says that $g$ defines a Riemannian metric in the vector bundle $\xi$.

The sets of sections $\G(T^r_s(E))$ are ${\cal F}(M)$-modules for every natural numbers $r,s$. On the sum
$\dd\bigoplus_{r,s}\G(T^r_s(E))$ a tensor product can be defined and one gets a tensorial algebra ${\cal T}(E)$. For the vector bundle $(TM,\tau,M)$ this reduces to the tensorial algebra of the manifold $M$.

A vector bundle $\xi=(E,p,M)$ is called {\bf a Finsler vector bundle} if it is endowed with a Finsler function defined as follows.

\bigskip\n{\bf Definition 1.1.} Let $\xi=(E,p,M)$ be a vector bundle of rank $m$. A {\it Finsler function} on $E$ is a nonnegative real function $F$ on $E$ with the properties
\begin{itemize}
\item[1)] $F$ is smooth on $E\setminus\{(x,0),x\in M\},$
\item[2)] $F(x,\lambda y)=\lambda F(x,y)$ for all $\lambda>0,$
\item[3)] The matrix with the entries $g_{ab}(x,y)=\frax12\ \derp{^2F^2}{y^a\pp y^b}$ is positive definite.
\end{itemize}
On $E$ we have the vertical distribution $u\to V_uE={\rm
Ker}\,p_{x,u},$ where $p_{*}$ denotes the differential of $p$.
This consists of vectors which are tangent to fibres and it is
locally spanned by $\left(\dot\partial_a:=\derpp{y^a}\right)$. We shall regard
also the vertical distribution as a vector subbundle
$VE:=\displaystyle\bigcup_{u\in E}V_uE\to E$ of $TE\to E$. Its sections will
be called vertical vector fields of $E$. The tensorial algebra
${\cal T}(VE)=\oplus{\cal T}^p_q(VE),\ p,q\in{\mathbb{N}}$ of this subbundle will
be used. Its elements will be indicated by the word ``vertical''.

A Finsler function $F$ on $E$ induces a Riemannian metric $g$
in the vertical bundle over $E$, given locally by
$$g(\dot\partial_a,\dot\partial_b)=g_{ab}(x,y).\leqno(1.1)$$
It provides also a set of vertical tensor fields by successively
deriving it with respect to $(y^a)$
$$C_{abc}(x,y)=\frac14\dot\partial_a\dot\partial_b\dot\partial_cL,\
D_{abcd}(x,y)=\frac18\dot\partial_a\dot\partial_b\dot\partial_c\dot\partial_dL,\
\mbox{etc.}\leqno(1.2)$$

The homogeneity of $F$ implies that the functions $g_{ab}(x)$ are positively homogenous of degree $0$ in $y^a$ and the components of vertical tensor fields from (1.2) are positively homogeneous in $y^a$ of degree $-1,-2,...$ etc. When the Euler theorem on homogeneous functions is applied to $F$ one gets $$F^2(x,y)=g_{ab}(x,y)y^ay^b.\leqno(1.3)$$
If the functions $g_{ab} $ do not depend on $y$ we obtain the simplest example of Finsler function on $E$. We may put this differently. Let $h_{ab}(x)$ be a Riemannian metric in the vector bundle $\xi$. Then $F$ given by $F^2(x,y)=h_{ab}(x)y^ay^b$ is a Finsler function on $E$. Thus any Riemannian vector bundle is a particular Finsler vector bundle.
On using the Riemannian metric $h_{ab}(x)$ as well as the components $\beta_a(x) $ of a section $\beta$ in $\xi^*$ and assuming that $h^{ab}\beta_a\beta_b<1$ one may construct a  Finsler function of Randers type on $E$ as follows
$$F(x,y)=\sqrt{h_{ab}(x)y^ay^b}+\beta_a(x)y^a.\leqno(1.4)$$
If we set $\alpha=\sqrt{h_{ab}(x)y^ay^b}$ and $\beta=\beta_a(x)y^a$ a  Finsler function on $E$ can be given as
$$F(x,y)=L(\alpha,\beta).\leqno(1.5)$$ for $L$ a homogeneous of degree one function in the both variables.

\section{Finsler vector bundles with nonlinear connections}

Let $\xi=(E,\pi,M)$ be a Finsler vector bundle of rank $m$ endowed  with the Finsler function $F$.

\medskip\noindent{\bf Definition 2.1} A nonlinear connection $N$ on $E$ is a
distribution $N:u\to N_uE$, $u\in E$, on $E$, which is supplementary to the vertical distribution $u \longrightarrow V_u$
on $E$.\medskip

We take the distribution $N$ as being locally spanned by $\de_k=\partial_k-N^a_k(x,y)\dot\partial_a$. By a
change of coordinates (1.1), the condition $\de_k=\derp{\wt x^i}{x^k}\ \wt\de_i$ is equivalent with
$$\wt N^a_j\partial_k\wt x^j=M^a_b(x)N^b_k(x,y)-\partial_k(M^a_b(x))y^b\leqno(2.1)$$
It is important to notice that from (2.1) it follows that the set of functions
$F^a_{bk}(x,y)=\dot\partial_bN^a_k(x,y)$ behaves under a change of coordinates (1.1) as the local coefficients of a linear connection in the vertical bundle over $\xi$, that is $$\wt F{}^a_{bk}(\wt x(x),\wt y(x,y))=M^a_c(x)\wt M^d_b(\wt x(x))
\derp{x^i}{\wt x^k}\ F^c_{di}(x,y)-\partial_i(M^a_c(x))\derp{x^i}{\wt x^k}\,y^c,\leqno(2.2)$$
where $\left(\derp{x^i}{\wt x^k}\right)$ is the inverse matrix of $\left(\derp{\wt x^k}{x^j}\right)$ and $(\wt M^d_b)$ denotes the inverse matrix of $(M^b_c).$

We should like to construct a linear connection $D$ in the vertical bundle $VE\to E$. In order to do this it suffices to
provide $D_{\de_k}\dot\partial_a$ and $D_{\dot\partial_a}\dot\partial_b.$ Using (2.2) we have the possibility
$$D_{\de_k}\dot\partial_a=F^b_{ak}(x,y)\dot\partial_b,\ D_{\dot\partial_b}\dot\partial_c=V^a_{bc}(x,y)\dot\partial_a,\leqno(2.3)^\circ$$ where necessarily $(V^a_{bc}(x,y))$ behave like the components of a vertical tensor field of type $(1,2)$.

In particular, we may take $V^a_{bc}=0$ and introduce

\medskip\noindent{\bf Definition 2.2.} The linear connection $D$ in the vertical bundle $VE\to E$ given by
$$D_{\de_k}\dot\partial_a=F^b_{ak}(x,y)\dot\partial_b,\ \ D_{\dot\partial_a}\dot\partial_b=0,\leqno(2.3)$$
is called the {\it Berwald connection} associated to $N$.\medskip

\medskip\noindent{\bf Definition 2.3.} We call the pair $(\xi,N)$ a Berwald
bundle if the functions $F^a_{bk}(x,y)=\dot\partial_bN^a_b(x,y)$ depend on $x$ only.\medskip

When $(\xi,N)$ is a Berwald bundle, the functions $F^a_{bk}(x,y)=F^a_{bk}(x)$ define a linear connection $\nabla$ in
$\xi$ by $$\nabla_{\partial_k}\varepsilon_b=F^a_{bk}(x)\varepsilon_a,\leqno(2.4)$$
for $(\varepsilon_a)$ a basis of local sections in $\xi$.

Conversely, if $\xi$ is endowed with a linear connection of local coefficients $\G^a_{bk}(x)$, then the functions
$$N^a_k(x,y)=\G^a_{bk}(x)y^b,\leqno(2.5)$$ define by setting $\de_{\dot k}=\partial_{\dot k}-N^a_k(x,y)\dot\partial_a$
a nonlinear connection on $E$ such that $(\xi,N)$ becomes a Berwald bundle. In other words, any vector bundle endowed with a
linear connection is a Berwald bundle.

We notice that the nonlinear connection (2.5) is positively homogeneous of degree 1 in $y=(y^a)$. This suggests us to confine
ourselves to the pairs $(\xi,N)$ with the functions $(N^a_k(x,y))$ positively homogeneous of degree 1 in $y$. The examples to be
given later will fall in this category. This assumption requires to eliminate from $E$ the image of the null section as we shall
do in the following.

It is well known that, see R.Miron \cite{4}, R. Miron and M. Apostasies \cite{5}, the Berwald connection induces a covariant
derivative in the tensorial algebra of the vertical bundle. This splits in two operators of covariant derivative. The first one is
called $h$--covariant derivative and is defined on functions and vertical vector fields as follows:
$$f_{\scri k}=\de_kf,\ X^a_{\scri k}=\de_kX^a+F^a_{bk}(x,y)X^b.\leqno(2.6)$$ It is extended by usual
rules to any vertical tensor field. The second, called the $v$-covariant derivative, is simply the partial derivative with
respect to $y$ $$f\bigi_a=\dot\partial_a f,\ X^a\bigi_b=\dot\partial_bX^a,\leqno(2.7)$$
since we have chosen $V^a_{bc}=0$.

We use the notation $\scri k$ and $\bigi_a$ for denoting the $h$- and $v$-covariant derivatives of any vertical tensor field.

\medskip\noindent{\bf Lemma 2.1.} {\it Let $\xi$ be endowed with a positively $1$-homogeneous nonlinear connection $N$ and $\scri k$ the $h$-covariant derivative defined by the Berwald connection associated to it. Then $$y^a_{\scri k}=0,\leqno(2.8)$$ holds.}

\medskip\noindent{\it Proof.} $y^a_{\scri k}=\de_ky^a+F^a_{bk}(x,y)y^b=F^a_{bk}(x,y)y^b-N^a_k(x,y)=0$
because of Euler theorem on homogeneous functions.

\medskip\noindent{\bf Lemma 2.2.} {\it Let $(\xi,N)$ be a Berwald bundle. Then for any vertical tensor field $T$ of local coefficients $T^{a_1...a_r}_{b_1...b_s}(x,y)$ we have $$T^{a_1...a_r}_{b_1...b_s}{}_{\scri k}\bigi
_a=T^{a_1...a_r}_{b_1...b_s}\bigi_{a\scri k}.\leqno(2.9)$$}

\medskip\noindent{\bf Proof.} One verifies (2.9) by a direct calculation keeping in mind
that the functions $F^a_{bk}=\dot\partial_aN^a_k$ do not depend on $y$.

We recall that in $\xi=(E,p,M),$ $E$ means in fact $E\setminus\{(x,0),\ x\in M\}.$

\medskip \noindent {\bf Definition 2.2.} Let $(\xi,F)$ be a Finsler vector bundle endowed with a positively
$1$-homogeneous nonlinear connection $N$. We say that $N$ is weakly compatible with $F$ if
$$F_{\scri k}:=\de_kF=0.\leqno(2.10)$$

In the following $N(N^a_i)$ will denote a positively $1$-homogeneous nonlinear connection. Given $N$ we may consider the Berwald connection $(\dot{\partial}_bN^a_i,0)$ and we may speak about $g_{ab\scri k}$.

\bigskip\noindent{\bf Definition 2.3.} Let $(\xi,F)$ be a Finsler vector bundle endowed with a positively
$1$-homogeneous nonlinear connection $N$. We say that $N$ is strongly compatible with $F$ if
$$g_{ab\scri k}=0.\leqno(2.11)$$

The terminology just introduced is explained by

\medskip \noindent {\bf Lemma 2.3.} {\it Let $(\xi,F)$ be a Finsler vector bundle endowed with a positively
$1$-homogeneous nonlinear connection $N$.Then $g_{ab\scri k}=0$ implies $F_{\scri k}= 0$.
The converse holds if the functions $( \dot{\partial}_bN^a_i)$ depends on $x$ only.}

\medskip \noindent {\it Proof.} We covariantly derive in the equality (1.3) and we get $F^2_{\scri k} = g_{ab\scri k}y^ay^b +2g_{ab}(x,y)y^ay^b_{\scri k}= 0 $ by (2.11) and the Lemma 2.1. For the converse, we covariantly derive in the equality defining $g_{ab}$. If the functions $(\dot{\partial}_bN^a_i)$ do not depend on $y$, the Lemma 2.2 applies in order to get $$g_{ab\scri k} = \frax{1}{2}\frax{\partial ^2(F^2_{\scri k})}{\partial y^a \partial y^b} =0, $$ q.e.d.

Let be $c:[0,1]\to M$, $t\to c(t),\ t\in[0,1]$ a smooth curve on
$E$. A section $A$ of $\xi$ along $c$ given as $A(t)=A^a(t)\varepsilon_a$
is said to be {\it parallel} with respect to the nonlinear connection $N$ if $ A_*(\dot c )$ are horizontal vectors. Here $A_*$ means the differential of the section  $A:M \to E$. A direct calculation shows that the section $A$ is parallel along the curve $c$ if and only if in any local chart on $M$, we have $$\frac{dA^a}{dt}+N^a_{k}(c(t),A(t)\
\frac{dx^k}{dt}=0,\leqno(2.12)$$ where $t \to x^k(t))$ are the local equations of the curve $c$.

For the initial conditions $c(0)=x$ and $A^a(0)=A^a_0$, the system of differential equations (2.12) admits an unique solution
$A^a(x(t))$ and if one assigns to $(A^a_0)\in E_x$ the element $A^a(x(1))\in E_{c(1)=z}$ one obtains an application $P_c:E_x\to
E_z$ called {\it parallel translation} along $c$, defined by $N$. We notice that because of the homogeneity of the functions $N^a_i$ the solutions of (2.12) are defined on $[0,1]$. The application $P_c:E_x\to E_z$ is a bijection and in general is not a linear map since the system (2.12) is not a linear one.

Now if one considers all loops on $M$ in $x\in M,$ the corresponding parallel translations as bijections from  $E_x\to
E_x$ provide a group with respect to their composition, called the helotomy group $\phi(x)$ of $N$ in $x\in M$. This is not a linear
group.

Let $F_x$ be the restriction of $F$ to the fibre $E_x.$ We call $F$-map a bijection $f:(E_x,F_x)\to(E_z,F_z)$ with the
property $F_x(u)=F_z(f(u))$ for every $u\in E_x.$

\medskip\noindent{\bf Theorem 2.1.} {\it Let the Finsler vector bundle $(\xi ,F)$ be endowed with a nonlinear connection $N$ which is weakly compatible with $F$. Then all parallel translations of $\nabla$ are $F$-maps. In particular, the holonomy groups $\phi(x),$ $x\in M$, consists of $F$-maps.}\medskip

\medskip\noindent{\it Proof.} Let $c:[0,1]\to M$ be a curve joining the points $x=c(0)$ and $z=c(1)$ of $M$. Consider a parallel section $A(t):=A(c(t))$, $t\in[0,1]$, of $\xi$ along $c$. We show that the function $f:t\to
F(x(t),A(t))$, $t\in[0,1]$, is constant. Indeed, $$\frac{dF(x(t),A(t))}{dt}=(\partial_k)\ \frac{dx^k}{dt}+(\dot\partial_aF)\
\frac{dA^a}{dt}\bld{(2.4)}{=}F_{\scri k}\ \frac{dx^k}{dt}=0.$$ Consider $A_0\in E_x$ and $A(t)$ the unique solution of (2.4) with
the initial condition $A_0$. Then $P_c(A_0)=A_1$, where $A_1=A(1)$ and since $f$ is constant, we get
$F_x(A_0)=F_z(A_1)=L_z(F_c(A_0)),$ q.e.d.

\section{Metrizability of Berwald connection}

Let the Finsler vector bundle $(\xi,F)$ be endowed with a nonlinear connection $N$ which is weakly compatible with $F$ and such that $(\xi,N)$ is a Berwald bundle. Then by Theorem 2.1, all parallel translations defined by $\nabla$ are isometries, that is, linear $F$- maps.

In particular, the elements of $\phi(x)$ are isometries of the Minkowski space $(E_x,F_x)$. And $\phi(x)$ is a subgroup of the
$G(I_x)$, the group of all linear isomorphisms which leave invariant the indicatrix $I_x$.

These facts are basic in the proof of the main result of this section.

\bigskip\noindent{\bf Theorem 3.1.} {\it If $\xi,F$ is endowed with a nonlinear connection $N$ which is weakly compatible with $F$ and $\xi,N$ is a Berwald bundle, then the linear connection $\nabla$ is metrizable, that is, there exists a
Riemannian metric $h$ in $\xi$ such that $\nabla h=0.$}\bigskip

\bigskip\noindent{\it Proof.} Let be $x_0\in M$  and the Minkowski space $(E_{x_0},F_{x_0}).$ The indicatrix $I_x$ is compact. It follows that the group $G:=G(I_x)$ is a compact Lie group. We know that $G$ contains $\phi(x)$ as a Lie subgroup but in general $\phi(x)$ is not compact. Let $<\,\cdot\,>$ be an arbitrary inner product in $E_{x_0}$. Define a new inner product on $E_{x_0}$ by
$$h_{x_0}(u,v)=\frac1{{\rm Vol}(G)}\int_G<gu,gv>\,\mu_G,\ \mbox{ for }g\in G,\ u,v\in E_{x_0},$$ where $\mu_G$ denotes the bi--invariant Haar measure on $G$. It follows that $h_{x_0}$ is $G$--invariant and, in particular, it is $\phi(x_0)$--invariant, i.e., $h_{x_0}(Pu,Pv)=h_{x_0}(u,v)$ for any $P\in\phi(x_0).$ Now we transfer $h_{x_0}$ to all the points of $M$. For any point $x\in M$, we consider a curve $c$ joining $x$ with $x_0$ $(c(0)=x,\ c(1)=x_0).$

Define $h_x(A,B)=h_{x_0}(P_cA,P_cB),$ $A,B\in E_x.$ The property that $h_{x_0}$ is $\phi(x_0)$--invariant assures that $h_x$ does
not depend on the curve $c$.

The mapping $h:x\longrightarrow  h_x:E_x\times E_x\to R$ is smooth since $P_c$ smoothly depends on $x$ by a general result about dependence of solutions of an ordinary differential equation on initial data. Thus a Riemannian metric $h$ in $\xi$ is obtained. The proof is ended with the help of

\bigskip\noindent{\bf Lemma 3.1.} {\it Let $h$ be a Riemannian metric in $\xi$ and $t\to c(t)$, $t\in\R$, a curve with $c(0)=x\in M.$ Then $$\lim_{t\to0}\ \frac1t\left(h_{c(t)}(P_cA,P_cB)-h_x(A,B)\right)= \left(\nabla_{\dot c(0)}h\right)(A,B)(x),\leqno(3.1)$$ where $A,B\in E_x$ and $P_c:E_x\to E_{c(t)}$ is the parallel translation along $c$.}\bigskip

Indeed, by the definition of $h$, the term in the left side of $(3.1)$ vanishes.

For the proof of Lemma 3.1 we refer to \cite{1}.

\bigskip\noindent{\bf Corollary  3.1.}{\it Let $\Gamma$ be a linear connection in the vector bundle $\xi=(E,p,M)$ with $M$ connected. Suppose that $E$ is endowed with a Finsler function $F$ with the associated Finsler metric $g_{ab}(x,y).$ Let $\scri k$ be the
$h$-covariant derivative operator induced by $\Gamma.$ If $g_{ab\scri k}=0$, then $\Gamma$ is metrizable.}\medskip

\bigskip\noindent{\it Proof.} The linear connection $\Gamma$ induces an $h$- covariant derivative operator and if $g_{ab\scri k}=0$ the Theorem 3.1 applies to get that $\Gamma$ is metrizable.


\newpage

\runningauthor={M. ANASTASIEI}
\runningtitle={GEOMETRY OF LAGRANGIANS AND SEMISPRAYS ON LIE ALGEBROIDS}
\noindent
\baselineskip 8pt
\noindent{\footnotesize{Geom. Balkan Press, Bucharest}}
\hfill\break
{\footnotesize{BSG Proc., 13,10-17, 2006}}
\vskip 2cm
\baselineskip 11.5pt plus .15pt
\centerline{\bf\Large GEOMETRY OF LAGRANGIANS AND}
\vskip .2cm
\centerline{\bf\Large SEMISPRAYS ON LIE ALGEBROIDS}
\vskip .5cm
\centerline{\bf {\footnotesize{BY}}}
\vskip .5cm
\centerline{\bf {\footnotesize{M. ANASTASIEI}}}
\vskip 1cm

\begin{abstract}
One considers a regular Lagrangian $L$ on the total space of a Lie
algebroid and one associates to it a semispray suggested by the
form of the Euler-Lagrange equations established by A. Weinstein,
\cite {5}. Some properties of this semispray are pointed out.

2000 Mathematics Subject Classifications : 53C60, 53C07

Key words: regular Lagrangian, Euler- Lagrange equations,
semisprays, Lie algebroids
\end{abstract}

\setcounter{section}{0}
\section{Introduction} In a paper appeared in 1996, \cite {5}, Alan Weinstein proposed a Lagrangian formalism for Lie algebroids.
This is general enough to include several Lagrangian formalisms as those on tangent bundles, on tangent subbundles and on Lie
algebras. He obtains the Euler - Lagrange equations using the Poisson structure on the dual of the given Lie algebroid and the
Legendre transformation defined by a regular Lagrangian on it. He also defines a notion of semispray. Later on, E. Martinez, \cite
{3}, develops a Lagrangian formalism for Lie algebroids that is similar to Klein's formalism, \cite{2}. He mainly uses a vector
bundle which replaces the double tangent bundle from the usual case. A notion of semispray appears in this setting, too.

In this paper we are mainly dealing with the notion of semipray in A. Weinstein' sense. In Section 2 we recall necessary facts from
the theory of vector bundles and establish the notations following the monograph \cite {4}.

Section 3 is devoted to semisprays on Lie algebroids. We give a definition that is a direct generalization of the one used in
tangent bundle case and we prove that this is equivalent with the definition given by A. Weinstein, \cite {5}. A local characterization is also provided. Three invariants are associated to any semispray.

In Section 4 we show that any regular Lagrangian on a Lie algebroid induces a semispray. This is done on a direct way: the
Euler - Lagrange equations obtained by A. Weinstein suggest the form of the local coefficients of a semispray and by a direct
calculation we checked that those coefficients are the appropriate ones. Some examples are pointed out.

\section{Vector bundles}

Let $\xi=(E,\pi,M)$ be a vector bundle of rank $m$. Here $E$ and $M$ are smooth i.e. $C^\infty$ manifolds with dim$M=n$,
dim$E=n+m$, and $\pi:E\to M$ is a smooth submersion. The fibres $E_x=\pi^{-1}(x)$, $x\in M$ are linear spaces of dimension $m$
which are isomorphic with the type fibre $\mathbb{R}^m$.

Let $\{(U_\alpha,\psi_\alpha)\}_{\alpha\in A}$ be an atlas on $M$. A vector bundle atlas is
$\{(U_\alpha$, $\varphi_\alpha$, $\mathbb{R}^m)\}_{\alpha\in A}$ with the bijections $\varphi_\alpha: \pi^{-1}(U_\alpha)\to
U_\alpha\times\mathbb{R}^m$ in the form $\varphi _\alpha(u)=(\pi(u),\varphi_{\alpha,\pi (u)})$, where
$\varphi_{\alpha,\pi(u)}:E_{\pi(u)}\to\mathbb{R}^m$ is a bijection. The given atlas on $M$ and a vector bundle atlas
provide an atlas $\{(\pi^{-1}(U_\alpha)$, $\Phi_\alpha)\}_{\alpha\in A}$ on $E$.

Here $\Phi_\alpha:\pi^{-1}(U_\alpha)\to\phi_\alpha(U_\alpha)\times\mathbb{R}^m$ is the bijection given by
$\phi_\alpha(u)=$ $(\psi_\alpha(\pi(u))$, $\varphi_{\alpha,\pi (u)}(u))$. For $x\in M$, we put
$\psi_\alpha(x)=(x^i)\in\mathbb{R}^n$ and if $(U_\beta,\psi_\beta)$ is another local chart such that $x\in
U_\alpha\cap U_\beta\neq\phi$, we set $\psi_\beta(x)=\widetilde{x}^i$ and then
$\psi_\beta\circ\psi_\alpha^{-1}$ has the form
\begin{equation}
\widetilde{x}^i=\widetilde{x}^i(x^1,\cdots,x^n),\
{\rm{rank}}\left(\dfrac{\partial \widetilde{x}^i}{\partial x^j
}\right)=n.\tag{1.1}
\end{equation}

Let $(e_a)$ be the canonical basis of $\mathbb{R}^m$. Then $\varphi^{-1}_{\alpha,x}(e_a):=\varepsilon_a(x)$ is a basis of
$E_x$ and $u\in E_x$ has the form $u=y^a\varepsilon_a (x)$.

We take $(x^i,y^a)$ as coordinates on $E$. For the bundle chart $(U_\beta,\varphi_\beta,\mathbb{R}^m)$ we put
$\varphi^{-1}_{\beta,x}(e_a)=\widetilde{\varepsilon}_a(x)$ and then $u=\widetilde{y}^a\widetilde{\varepsilon}_a(x)$. If we set
$\varepsilon_a(x)=M^b_a(x)\widetilde{\varepsilon}_b$ with rank$(M_a^b(x))=m$ it follows that $\widetilde{y}^a=M^a_b(x)y^b$.
Thus the mapping $\phi_\beta\circ\phi_\alpha^{-1}$ has the form
\begin{equation}
\begin{array}{l}\widetilde{x}^i=\widetilde{x}^i(x^1,\cdots,x^n), {\rm{rank}}
\left(\dfrac{\partial\widetilde{x}^i}{\partial x^j}\right)=n\vspace{1.5mm}\\ \widetilde{y}^a=M^a_b(x)y^b,
{\rm{rank}}(M^a_b(x))=m.
\end{array}
\tag{1.2}
\end{equation}

The indices $i,j,k,...$ and $a,b,c ...$ will take the values $1,2,... n$ and $1,2,... m$, respectively. The Einstein convention
on summation will be used.

We denote by ${\cal F}(M)$, ${\cal F}(E)$ the ring of real functions on $M$ and $E$ respectively, and by $\chi(M)$, respectively
$\Gamma(E)$, $\chi(E)$ the module of sections of the tangent bundle of $M$, respectively of the bundle $\xi$ and of the tangent bundle of $E$. On $U_\alpha$, the vector fields $\left(\partial_k:=\dfrac{\partial}{\partial x^k}\right)$ provide
a local basis for $\chi(U_\alpha)$. The sections $\varepsilon_a:U_a\to \pi^{-1}(U_\alpha)$,
$\varepsilon_a(x)=\varphi^{-1}_{\alpha,x}(e_a)$ provide a basis for $\Gamma(\pi^{-1}(U_\alpha))$ and a section $A:U_\alpha\to
\pi^{-1}(U_\alpha)$ will take the form $A(x)=A^a(x)\varepsilon_a(x)$, $x\in U_\alpha$.

Let $\xi^*=(E^*,p^*,M)$ be the dual of the vector bundle $\xi$. We may also consider the tensor bundle $T^r_s(E)$ over $E$. The set
of sections $\Gamma(T^r_s(E))$ are ${\cal F}(M)-$modules for any natural numbers $r,s$. On the sum $\oplus_{r,s}\Gamma(T^r_s(E))$ a
tensor product can be defined and one gets a tensor algebra $T(E)$. For the tangent bundle $(TM,\tau,M)$ this reduces to the
tensor algebra of the manifold $M$. The tensor algebra of the manifold $E$ could be also involved. Its elements are sections in
${\cal T}^r_s(TE)$. The tensorial algebra of $E$ contains the subset of $d-$tensor fields on $E$. For a general definition of
these tensor fields we refer to \cite {4}, Ch. III. Shortly, these tensor fields are defined by components depending on $(x^i,y^a)$
and transforming by a change of coordinates as tensors but with the matrices $\left(\dfrac{\partial\widetilde{x}^i}{\partial
x^j}\right)$ and $(M^a_b(x))$ and their inverses, only. Notice that in the law of transformation of a tensor field on $E$ could
appear also the matrix $\left(\dfrac{\partial M^a_b(x)}{\partial x^i}y^b\right)$.

A large class of examples is provided by the sections in the vertical bundle over $E$. We recall that the vertical bundle
$VE\to E$ is the union of the fibres $V_uE=\ker\pi_{*,u}$ over $u\in E$, where $\pi_{*,u}$ is the differential of $\pi$. A basis
of local section of $VE\to E$ is given by $\left(\left.\dfrac{\partial}{\partial y^a}\right|_u\right)$ and
its dual is $dy^a|_u$. The local components of any element in $\Gamma(T^r_s(VE))$, transform under a change of coordinates on
$E$ with the matrix $(M^a_b(x))$ and its inverse $(W^a_b)$. We call such an element a vertical tensor field.

Now if $L:E\to M$ is a smooth function on $E$ (called usually a Lagrangian) then it is easy to check that functions
$\dfrac{\partial L}{\partial y^a}$, $g_{ab}=\dfrac{1}{2}$ $\dfrac{\partial^2 L}{\partial y^a\partial y^b}$,
$C_{abc}=\dfrac{1}{2}\dfrac{\partial g_{ab}}{\partial y^c}$ define vertical tensor fields of covariance indicated by the position and number of indices.

\section{Semisprays for Lie algebroids}

A vector bundle $\xi=(E,\pi,M)$ is called a Lie algebroid if it has the following properties:

\begin{enumerate}
\item The space of sections $\Gamma(\xi)$ is endowed with a Lie algebra structure $[,]$;
\item There exists a bundle map $\rho:E\to TM$ (called the {\it anchor map}) which induces a Lie algebra homomorphism (also denoted by $\rho$) from $\Gamma(\xi)$ to $\chi(M)$.
\item For any smooth functions $f$ on $M$ and any sections $s_1,s_2\in\Gamma(\xi)$ the following identity is satisfied
$$[s_1,fs_2]=f[s_1,s_2]+(\rho(s_1)f)s_2.$$
\end{enumerate}

Locally, we set
\begin{equation}
\rho(s_a)=\rho^i_a\dfrac{\partial}{\partial x^i},\
[\varepsilon_a,\varepsilon_b]=L^c_{ab}s_c,\tag{3.1}
\end{equation}

A change of local charts implies
\begin{equation}
\widetilde{\rho}^i_a=W^b_a\rho^j_b\dfrac{\partial\widetilde{x}^i}{\partial
x^j},\tag{3.2}
\end{equation}
where $W^b_a$ is the inverse of the matrix $(M^a_b)$.

Examples of Lie algebroids: the tangent bundle $\tau:TM\to M$ with
$\rho=$identity, any integrable subbundle of $TM$ with the
inclusion as anchor map, $TP/G$ for $P(M,G)$ a $G-$principal
bundle, see \cite {5}.

For a function $f$ on $M$ one defines its vertical lift $f^v$ on $E$ by $f^v(u)=f(\pi(u))$ and its complete lift $f^c$ on $E$
by $f^c(u)= \rho_a^iy^a\dfrac{\partial f}{\partial x^i}$ for $u=(x,y)$ in $E$. If $A=A^a(x)\varepsilon_a$ is a section in
$\xi$, the vertical lift $A^v$ is a vector field on $E$ defined by $A^v(x,y)= A^a(x)\dfrac{\partial }{\partial y^a}$ and the complete
lift $A^c$ is a vector field on $E$ defined by $$A^c(x,y)=A^a\rho^i_a\dfrac{\partial}{\partial x^i}+\(\rho^i_b\dfrac{\partial A^a}{\partial x^i}-A^dL^a_{db}\)y^b\dfrac{\partial}{\partial y^a}.$$ In particular, $\varepsilon^v_a=\dfrac{\partial }{\partial y^a},\varepsilon_a^c=\rho^i_a\dfrac{\partial}{\partial x^i}-L^d_{ab}y^b\dfrac{\partial}{\partial y^d}$.

A semispray $S$ for the tangent bundle $\tau:TM\to M$ is a vector field on $TM$ which at the same time is a section in the vector
bundle $\tau_*:TTM\to TM$, that is we have $\tau_{TM}(S(u))=u$ and $\tau_{*,u}(S(u))=u$, $\forall u\in TM$, where $\tau_{TM}$ is the
vector bundle projection $TTM\to TM$. It follows that $\tau_{*,u}(S(u))=\tau_{TM}(S(u))$, $\forall u\in TM$.

This equation suggests the following

{\bf Definition 3.1.} {\it Let $\xi=(E,\rho,M)$ be a Lie algebroid with the anchor $\rho$. A vector field $S$ on $E$ will be called a semispray if
\begin{equation}
\pi_{*,u}(S(u))=(\rho\circ\tau_E)(S(u)),\ \forall u\in E\tag{3.3}
\end{equation} where $\tau_E:TE\to E$ is the natural projection.}

Let $c:I\to M$, $I\subseteq\mathbb{R}$ be a curve on $M$ and let $\widetilde{c}:I\to E$ be any curve on $E$ such that
$\pi\circ\widetilde{c}=c$. Denote by $\dot{\widetilde{c}}$ the vector field that is tangent to $\widetilde{c}$.

{\bf Definition 3.2.} {\it We say that $\widetilde{c}$ is {\bf admissible} if $$\pi_*(\dot{\widetilde{c}})=\rho(\widetilde{c}).$$}

In local charts on $M$ and $E$, we have $c(t)=(x^i(t))$, $\widetilde{c}(t)=(x^i(t),y^a(t))$ and
$\dot{\widetilde{c}}(t)=\dfrac{dx^i}{dt}\dfrac{\partial}{\partial x^i}+\dfrac{dy^a}{dt}\dfrac{\partial}{\partial y^a}$, $t\in I$.

It results

{\bf Lemma 3.1.} {\it The curve $\widetilde{c}$ is admissible if and only if
\begin{equation}
\dfrac{dx^i}{dt}(t)=\rho^i_a(x(t))y^a(t),\ \forall t\in
I.\tag{3.4}
\end{equation}}

Again in local charts, let be $S=X^i\dfrac{\partial}{\partial x^i}+Y^a\dfrac{\partial}{\partial y^a}$ a vector field on $E$.

This is a semispray if and only if
\begin{equation}
X^i(x,y)=\rho^i_a(x)y^a.\tag{3.5}
\end{equation}

Thus the coordinates $(Y^a(x,y))$ are not determined. We set for convenience $Y^a=-2G^a$. Furthermore, under a change of
coordinates $(x^i,y^u)\to(\widetilde{x}^i,\widetilde{y}^a)$, the coordinates $(X^i),(G^a)$ have to change as follows:
\begin{equation}
\widetilde{X}^i=\dfrac{\partial\widetilde{x}^i}{\partial
x^j}(x)X^j,\tag{3.6}
\end{equation}
\begin{equation}
\widetilde{G}^a=M^a_b G^b-\dfrac{1}{2}\dfrac{\partial
M^a_b}{\partial x^i}y^b\rho^i_c y^c.\tag{3.7}
\end{equation}

Using (3.2) one easily sees that the coordinates $(X^i(x,y))$ given by (3.5) verify (3.6).

Concluding, we have

{\bf Theorem 3.1.} {\it A vector field $S=(\rho^i_a y^a)\dfrac{\partial}{\partial x^i}-2G^a\dfrac{\partial}{\partial
y^a}$ on $E$ is a semispray if and only if the coordinates $(G^a)$ transform by $(3.7)$.}

The integral curves of $S$ are given by the system of differential equations
\begin{equation}
\dfrac{dx^i}{dt}=\rho^i_a(x)y^a,\
\dfrac{dy^a}{dt}+2G^a(x,y)=0.\tag{3.8}
\end{equation}

It comes out these curves are all admissible. The converse is also true, that is we have

{\bf Theorem 3.2.} {\it A vector field on $E$ is a semispray if and only if all its integral curves are admissible.}

{\it Remark 3.1.} The characterization of a semispray provided by the Theorem 3.2 was taken by A. Weinstein, \cite{5}, as definition
for a semispray on $E$.

{\it Remark 3.2.} (i) Let us assume that $\rho=0$. Then the admissible curves are all curves from the fibre $E_{x_0}$,
$x_0(x^i_0)\in M$. The integral curves of a semispray $S$ are given by the equations $\dfrac{dy^a}{dt}+2G^a(x_0,y)=0$.

(ii) The system of equations (3.8) is no longer equivalent with a second order differential equations as it happens for $TM$. Thus
the term of ``second order differential equations'' used sometimes for a semispray is no longer appropriate.

(iii) Let $D$ a distribution on $M$. We regard it as a subbundle of $TM$ and so we may view it as a Lie algebroid with the natural
inclusion as anchor map. Using a local basis on $D$ one can see that the admissible curves are those that are tangent to the
distribution $D$. For details we refer to \cite {1}.

Let $\widehat{S}$ be another semispray on $E$. Then $\widehat{S}=(\rho^i_a y^a)\dfrac{\partial}{\partial
x^i}-2\widehat{G}^a\dfrac{\partial}{\partial y^a}$, where the functions $(\widehat{G}^a(x,y))$ have to satisfy (3.7) under a
change of coordinates on $E$. It follows that $\widehat{S}-S=2(G^a-\widehat{G}^a)\dfrac{\partial}{\partial y^a}$
and the functions $D^a=G^a-\widehat{G}^a$ transform by the rule
\begin{equation}
\widehat{D}^a=M^a_b D^b.\tag{3.9}
\end{equation}

So we have proved

{\bf Theorem 3.3.} {\it Any two semisprays on $E$ differ by a vertical vector field on $E$.}

A different point of view on semisprays for algebroids was proposed by E.Martinez,\cite {3}. It can be shortly described as
follows.

Let ${\cal L}^\pi E $ be the subset of $E\times TE$ defined by ${\cal L}^\pi E=\{(u,z)| \rho(u)=\pi_*(z)\}$ and denote by
$\pi_L :{\cal L}^\pi E \longrightarrow E$ the mapping given by $\pi_L(u,z)=\tau_E(z)$. Then $({\cal L}^\pi E,\pi_L,E)$ is a
vector bundle over $E$ of rank $2m$. One proves that this vector bundle is also a Lie algebroid.

One associates to a section $A$ of $\xi$ the vertical lift $A^V$ and the complete lift $A^C$ as sections of $\pi_L:{\cal L}^\pi E
\longrightarrow E$ given by $$A^V(u)=(0,A^v(u)),\ A^C(u)=(A(\pi(u)),\ A^c(u)),\ \ u\in E.$$

If $\{s_a\}$ is a local basis of $\Gamma (E))$, then $\{s_a^V,s_s^C\}$ is a local basis for $\Gamma({\cal L}^\pi E)$.

The vector bundle $({\cal L}^\pi E,\pi_L, E)$ admits a canonical section $C$ called the {\it Liouville or Euler section} defined
by $C(u)=\(o,y^a\dfrac{\partial}{\partial y^a}\)$ for $u=y^a\varepsilon_a\in E$. A section $J$ of the vector bundle ${\cal L}^\pi E \bigoplus ({\cal L}^\pi E)^*\longrightarrow E$ characterized by the conditions $J(A^V)=0$, $J(A^C)=A^V$, $A\in\Gamma E$ is called the {\it vertical endomorphism.} We have that $J^2=0$. A section $S$ of the vector bundle $({\cal L}^\pi E,\pi_L,E)$ is said to be a semispray if it satisfies the condition $JS=C$. This definition is equivalent with the preceding one.
Indeed, in local coordinates if we set $S=A^a\varepsilon_a^C+S^a\varepsilon_a^V$, the condition $JS=C$ gives $A^a=y^a$ and
so $S=y^a\(\rho^i_a\dfrac{\partial}{\partial x^i}-L_{ab}^c y^b\dfrac{\partial}{\partial y^c}\)+S^a\dfrac{\partial}{\partial
y^a}=y^a\rho^i_a\dfrac{\partial}{\partial x^i}+S^a\dfrac{\partial}{\partial y^a}$ since $L^c_{ab}y^ay^b=0$.

For a semispray on $TM$, a case when this is equivalent with a system of second order differential equations (SODE), there exists
a way to find geometric invariants that to determine, up to a change of coordinates, the solutions of the system.

This way led to a KCC-theory named so as after Kosambi, Cartan and Chern.

The KCC-theory apparently does not work for semisprays on Lie algebroids. However, at least formally we can associate to a
semispray $S=(\rho^i_a y^a)\dfrac{\partial}{\partial x^i}-2G^a(x,y)\dfrac{\partial}{\partial y^a}$, the following invariants:
\begin{equation}
\zeta^a=2G^a-\dfrac{\partial G^a}{\partial y^b}y^b,\tag{3.10}
\end{equation}
\begin{equation}
\Xi^a=\dfrac{\partial G^a}{\partial y^b}-\dfrac{\partial
G^a}{\partial y^b\partial y^c}y^c,\tag{3.11}
\end{equation}
\begin{equation}
\Gamma^a=2G^a-2\dfrac{\partial G^a}{\partial
y^b}y^b+\dfrac{\partial G^a}{\partial y^b\partial y^c}y^b y^c.
\tag{3.12}
\end{equation}

Indeed, it is not difficult to check that all these sets of functions define vertical vector fields on $E$.

To find a complete list of such invariants could be a future task.

\section{A semispray derived from a regular Lagrangian}

Let $L:E\to R$ be a regular Lagrangian on the Lie algebroid $(E,[,],\rho)$, that is $L$ is a smooth functions such that the
matrix with the entries
\begin{equation}
g_{ab}(x,y)=\dfrac{1}{2}\dfrac{\partial^2 L}{\partial y^a\partial
y^b},\tag{4.1}
\end{equation}
is of rank $m$.

In \cite{5}, one associates to $L$ the Euler - Lagrange equations
\begin{equation}
\dfrac{d}{dt}\left(\dfrac{\partial L}{\partial y^a
}\right)=\rho^i_a\dfrac{\partial L}{\partial
x^i}+L^c_{ba}y^b\dfrac{\partial L}{\partial y^c},\tag{4.2}
\end{equation}
for $c(t)=(x^i(t),y^a(t))$ an admissible curve.

Expanding the derivative in (4.2), using (4.1) and (3.4), we may
put (4.2) in the form
\begin{equation}
\dfrac{dy^a}{dt}+2G^a_L(x,y)=0,\tag{4.3}
\end{equation}
with the notation
\begin{equation}
G^a_L=\dfrac{1}{4}g^{ab}\left(\dfrac{\partial^2 L}{\partial
y^b\partial x^i}\rho^i_c y^c-\rho^i_b\dfrac{\partial L}{\partial
x^j}-L^c_{bd}y^d\dfrac{\partial L}{\partial y^c}\right).\tag{4.4}
\end{equation}

We show that the function $(G^a_L)$ verifies (3.7) under a change
of coordinates on $E$.

We set
\begin{equation}
E_a=4g_{ab}G^b,\tag{4.5}
\end{equation}
where
\begin{equation}
E_a=\dfrac{\partial^2 L}{\partial y^a\partial x^i}\rho^i_b
y^b-\rho^i_a\dfrac{\partial L}{\partial
x^i}-L^c_{ba}y^b\dfrac{\partial L}{\partial y^c}.\tag{4.6}
\end{equation}

Then we use (3.2) as well as the following equations:
$$\dfrac{\partial L}{\partial x^i}=\dfrac{\partial L}{\partial \widetilde{x}^j}\dfrac{\partial
\widetilde{x}^j}{\partial x^i}+\dfrac{\partial L}{\partial
\widetilde{y}^a}\dfrac{\partial M^a_c}{\partial x^i}y^c$$
$$\dfrac{\partial^2 L}{\partial y^a\partial x^i}=M^b_a\left(\dfrac{\partial^2 L}{\partial
y^b\partial\widetilde{x}^j}\dfrac{\partial\widetilde{x}^j}{\partial
x^i}+2\widetilde{g}_{db}\dfrac{\partial M^d_c}{\partial
x^i}y^c\right)+\dfrac{\partial
L}{\partial\widetilde{y}^d}\dfrac{\partial M^d_a}{\partial x^i}$$
$$L^c_{ab}M^e_c=M^c_a M^d_b\widetilde{L}^e_{cd}+\rho^k_a\dfrac{\partial M^e_b}{\partial
x^k}-\rho^k_b\dfrac{\partial M^e_a}{\partial x^k}$$ in order to
derive
\begin{equation}
E_a=M^b_a\widetilde{E}_b+2M^b_a\widetilde{g}_{bd}\dfrac{\partial
M^d_c}{\partial x^i}y^c\rho^i_d y^d.\tag{4.7}
\end{equation}

Using this in (4.5) one shows that $\widetilde{G}^a_L$ is related
to $G^a_L$ as in (3.7).

Thus we have proved

{\bf Theorem 4.1.} {\it Let $L$ be a regular Lagrangian on the Lie
algebroid $(E,[,],\rho)$. Then $L$ defines a semispray
$S_L=(\rho^i_a y^a)\dfrac{\partial}{\partial
x^i}-2G^a_L(x,y)\dfrac{\partial}{\partial y^a}$, where the
function $G^a_L$ are given by (4.4).}

{\bf Example 4.1.} Let $g_{ab}(x)$ be the coefficients of a
Riemannian metric in the Lie algebroid $(E,[,],\rho)$. Then
\begin{equation}
L(x,y)=g_{ab}(x)y^a y^b\tag{4.8}
\end{equation}
is a regular Lagrangian on $E$. The semispray associated to it is
determined by the functions
\begin{equation}
G^a=\dfrac{1}{2}g^{ab}\left(\dfrac{\partial g_{bc}}{\partial
x^i}\rho^i_d-\dfrac{1}{2}\dfrac{\partial g_{cd}}{\partial
x^i}\rho^i_b-L^e_{db}g_{ec}\right)y^c y^d.\tag{4.9}
\end{equation}

{\bf Example 4.2.} A more general example is provided by the
regular Lagrangians which are homogeneous of degree 2 in $(y^a)$.
By the Euler theorem one obtains
\begin{equation}
L(x,y)=g_{ab}(x,y)y^a y^b,\tag{4.10}
\end{equation}
where $(g_{ab}(x,y))$ are homogeneous functions of degree 0.

As $\dfrac{\partial}{\partial y^a}$ are homogeneous functions of
degree 1 and the derivative with respect to $(x^j)$ does not
affect the degree of homogeneity, it results that the coefficients
$(G^a)$ from (4.4) are homogeneous of degree 2 in $(y^a)$. This
fact is equivalent with $\zeta^a=0$ and so we have a meaning of
the invariant $\zeta^a$. The corresponding semispray is called a
spray.


\newpage

\runningauthor={M. ANASTASIEI}
\runningtitle={MECHANICAL SYSTEMS ON LIE ALGEBROIDS}
\noindent
\baselineskip 8pt
\noindent{\footnotesize{Algebras Groups Geom.}}
\hfill\break
{\footnotesize{23, no. 3, 235-245, 2006}}
\vskip 2cm
\baselineskip 11.5pt plus .15pt
\centerline{\bf\Large MECHANICAL SYSTEMS}
\vskip .2cm
\centerline{\bf\Large ON LIE ALGEBROIDS}
\vskip .5cm
\centerline{\bf {\footnotesize{BY}}}
\vskip .5cm
\centerline{\bf {\footnotesize{M. ANASTASIEI}}}
\vskip 1cm

\centerline{\footnotesize{\it Dedicated to the 70th birthday of Professor
Ruggero Maria Santilli}}

\setcounter{section}{0}
\section{Introduction}
The simplest mathematical model for a mechanical system is made of
a Riemannian manifold $(M,g)$ with $M$ a smooth manifold of states
$x=(x^i)$ and $g=(g_{ij}(x))$ a Riemannian metric provided by the
kinetic energy $\dfrac{1}{2}g_{ij}(x)\dot x^i \dot x^j $ of the
system. The difference between kinetic energy and the potential
energy defines the Lagrangian $L$ of the system and the solution
curves of the Euler-Lagrange equations written for $L$ are the
evolution curves of the system. The regular Lagrangian $L$ is
living on the tangent manifold $TM$ and thus a new space (phase
space) is coming into play.

In many cases the mechanical systems involve external forces that
are not of gradient type. These forces are modelled by a covector
field $F=(F_i(x))$ or equivalently by a vector field of components
$(g^{ij}F_j)$ on the manifold $M$ and then the second Newton's law
of dynamics takes the form $ \dfrac{d}{dt}\(\dfrac{\partial
L}{\partial \dot x^i}\)-\dfrac{\partial L}{\partial x^i} = F_i(x)$
and it gives also the evolution curves of the system. Thus a
mechanical system with external forces is defined as a triple
$(M,L,F)$ for $L$ a regular Lagrangian and $F$ a covector(vector)
field. The theory of these systems was extensively and clearly
presented in an excellent book by R.M. Santilli, [5].

But there exist cases when the external forces depend also on
velocity, that is $F$ is living on $TM$. The corresponding theory
was developed by R. Miron and C. Frigioiu, [2] and Munoz-Lecanda
M.C. et al., [4].

On the other hand, A. Weinstein constructed in [6] a Lagrangian
formalism on a Lie algebroid. A Lie algebroid is a vector bundle
$(E,\pi,M)$ that is endowed with a Lie bracket $[,]$ for its
sections and is anchored to the tangent bundle with a bundle
morphism $\rho:E\longrightarrow TM$ that induces on sections a
Lie algebra homomorphism denoted also by $\rho$ such that for any
two sections $A,B$ and any function $f$ on $M$ we have $[A,fB]=
f[A,B] +\rho (A)f.B$. The formalism of A. Weinstein contains the
Euler - Lagrange equations for a Lagrangian on $E$. Thus it is
open a way for approaching the theory of mechanical systems with
external forces (not of gradient type) on Lie algebroids. This is
the aim of this paper. Moreover, for enlarging the applicability
of our theory we assume that the external forces depend on the
fibre variables. These variables may have various meaning (velocities in tangent bundle case).

Our main result says that if the system is dissipative then its
energy is decreasing on the evolution curves. The energy of the
system is also used for constructing a Lyapunov function for an
equilibrium point of the system. These are presented in Section 4.
The preceding sections are devoted to necessary facts from the
theory of vector bundles and of Lie algebroids.

\section{Vector bundles}

Let $\xi=(E,\pi,M)$ be a vector bundle of rank $m$. Here $E$ and
$M$ are smooth i.e. $C^\infty$ manifolds with dim$M=n$,
dim$E=n+m$, and $\pi:E\to M$ is a smooth submersion. The fibres
$E_x=\pi^{-1}(x)$, $x\in M$ are linear spaces of dimension $m$
which are isomorphic with the type fibre $\mathbb{R}^m$.

Let $\{(U_\alpha$, $\psi_\alpha)\}_{\alpha\in A}$ be an atlas on $M$.
A vector bundle atlas is
$\{(U_\alpha$, $\varphi_\alpha$, $\mathbb{R}^m)\}_{\alpha\in A}$ with
the bijections $\varphi_\alpha:\pi^{-1}(U_{\alpha})\to
U_\alpha\times\mathbb{R}^m$ in the form $\varphi_\alpha (u)=(\pi
(u),\varphi_{\alpha,\pi (u)}(u))$, where
$\varphi_{\alpha,\pi(u)}:E_{\pi(u)}\to\mathbb{R}^m$ is a
bijection. The given atlas on $M$ and a vector bundle atlas
provide an atlas $\{(\pi^{-1}(U_\alpha)$, $\Phi_\alpha)\}_{\alpha\in A}$ on $E$.
Here $\Phi_\alpha:\pi^{-1}(U_\alpha)\to\phi_\alpha(U_\alpha)\times\mathbb{R}^m$
is the bijection given by $\phi_\alpha(u)=$ $(\psi_\alpha(\pi(u))$, $\varphi_{\alpha,\pi
(u)}(u))$. For $x\in M$, we put
$\psi_\alpha(x)=(x^i)\in\mathbb{R}^n$ and if
$(U_\beta,\psi_\beta)$ is another local chart such that $x\in
U_\alpha\cap U_\beta\neq\phi$, we set
$\psi_\beta(x)=(\widetilde{x}^i)$ and then
$\psi_\beta\circ\psi_\alpha^{-1}$ has the form
\begin{equation}
\widetilde{x}^i=\widetilde{x}^i(x^1,\cdots,x^n),\
{\rm{rank}}\left(\dfrac{\partial \widetilde{x}^i}{\partial x^j
}\right)=n.\tag{2.1}
\end{equation}

Let $(e_a)$ be the canonical basis of $\mathbb{R}^m$. Then
$\varphi^{-1}_{\alpha,x}(e_a):=\varepsilon_a(x)$ is a basis of
$E_x$ and $u\in E_x$ has the form $u=y^a\varepsilon_a (x)$.

We take $(x^i,y^a)$ as coordinates on $E$. For the bundle chart
$(U_\beta,\Psi_\beta,\mathbb{R}^m)$ we put
$\varphi^{-1}_{\beta,x}(e_a)=\widetilde{\varepsilon}_a(x)$ and
then $u=\widetilde{y}^a\widetilde{\varepsilon}_a(x)$. If we set
$\varepsilon_a(x)=M^b_a(x)\widetilde{\varepsilon}_b$ with
rank$(M_a^b(x))=m$ it follows that $\widetilde{y}^a=M^a_b(x)y^b$.
Thus $\Phi_\beta\circ\Phi_\alpha^{-1}$ has the form
\begin{equation}
\begin{array}{l}\widetilde{x}^i=\widetilde{x}^i(x^1,\cdots,x^n), {\rm{rank}}
\left(\dfrac{\partial\widetilde{x}^i}{\partial
x^j}\right)=n\vspace{1.5mm}\\ \widetilde{y}^a=M^a_b(x)y^b,
{\rm{rank}}(M^a_b(x))=m.
\end{array}
\tag{2.2}
\end{equation}

The indices $i,j,k,... a,b,c ...$ will take the values $1,2,... n$
and $1,2,... m$, respectively. The Einstein convention on
summation will be used.

We denote by ${\cal F}(M), {\cal F}(E)$ the ring of real functions
on $M$ and $E$ respectively, and by $\chi(M)$, respectively
$\Gamma(E),\chi(E)$ the module of sections of the tangent bundle
of $M$, respectively of the bundle $\xi$ and of the tangent bundle
of $E$. On $U_\alpha$, the vector fields
$\left(\partial_k:=\dfrac{\partial}{\partial x^k}\right)$ provide
a local basis for $\chi(U_\alpha)$. The sections
$\varepsilon_a:U_a\to p^{-1}(U_\alpha)$,
$\varepsilon_a(x)=\varphi^{-1}_{\alpha,x}(e_a)$ provide a basis
for $\Gamma(p^{-1}(U_\alpha))$ and a section $A:U_\alpha\to
p^{-1}(U_\alpha)$ will take the form
$A(x)=A^a(x)\varepsilon_a(x)$, $x\in U_\alpha$.

Let $\xi^*=(E^*,p^*,M)$ be the dual of the vector bundle $\xi$. We
may also consider the tensor bundle $T^r_s(E)$ over $E$. The set
of sections $\Gamma(T^r_s(E))$ are ${\cal F}(M)-$modules for any
natural numbers $r,s$. On the sum $\oplus_{r,s}\Gamma(T^r_s(E))$ a
tensor product can be defined and one gets a tensor algebra
$T(E)$. For the tangent bundle $(TM,\tau,M)$ this reduces to the
tensor algebra of the manifold $M$. The tensor algebra of the
manifold $E$ could be also involved. Its elements are sections in
${\cal T}^r_s(TE)$. The tensorial algebra of $E$ contains the
subset of $d-$tensor fields on $E$. For a general definition of
these tensor fields we refer to [3], Ch. III. Shortly, these
tensor fields are defined by components depending on $(x^i,y^a)$
and transforming tensorially by a change of coordinates but with
the matrices $\left(\dfrac{\partial\widetilde{x}^i}{\partial
x^j}\right)$ and $(M^a_b(x))$ and their inverses, only. Notice
that in the law of transformation of a tensor field on $E$ could
appear also the matrix $\left(\dfrac{\partial M^a_b(x)}{\partial
x^i}y^b\right)$.

A large class of examples is provided by the sections in the
vertical bundle over $E$. We recall that the vertical bundle
$VE\to E$ is the union of the fibres $V_uE=\ker\pi_{*,u}$ over
$u\in E$, where $\pi_{*,u}$ is the differential of $\pi$. A basis
of local section of $VE\to E$ is given by
$\left(\left.\dfrac{\partial}{\partial y^a}\right|_u\right)$ and
its dual is $dy^a|_u$. The local components of any element in
$\Gamma(T^r_s(VE))$, transform under a change of coordinates on
$E$ with the matrix $(M^a_b(x))$ and its inverse $(W^a_b)$. We
call such an element a vertical tensor field.

Now if $L:E\to M$ is a smooth function on $E$ (called usually a
Lagrangian) then it is easy to check that functions
$\dfrac{\partial L}{\partial y^a}$, $g_{ab}=\dfrac{1}{2}$
$\dfrac{\partial^2 L}{\partial y^a\partial y^b}$,
$C_{abc}=\dfrac{1}{2}\dfrac{\partial g_{ab}}{\partial y^c}$ define
vertical tensor fields of covariance indicated by the position and
number of indices.

\section{Lagrangians on a Lie algebroid. Associated semispray}

A vector bundle $\xi=(E,\pi,M)$ is called a Lie algebroid if it
has the following properties:

\begin{enumerate}
\item The space of sections $\Gamma(\xi)$ is endowed with a Lie
algebra structure $[,]$;
\item There exists a bundle map $\rho:E\to TM$ (called the {\it anchor
map}) which induces a Lie algebra homomorphism (also denoted by
$\rho$) from $\Gamma(\xi)$ to $\chi(M)$.
\item For any smooth functions $f$ on $M$ and any sections
$s_1,s_2\in\Gamma(\xi)$ the following identity is satisfied
$$[s_1,fs_2]=f[s_1,s_2]+(\rho(s_1)f)s_2.$$
\end{enumerate}

Locally, we set
\begin{equation}
\rho(\varepsilon_a)=\rho^i_a\dfrac{\partial}{\partial x^i},\
[\varepsilon_a,\varepsilon_b]=L^c_{ab}\varepsilon_c,\tag{3.1}
\end{equation}

A change of local charts implies
\begin{equation}
\widetilde{\rho}^i_a=W^b_a\rho^j_b\dfrac{\partial\widetilde{x}^i}{\partial
x^j}.\tag{3.2}
\end{equation}

Examples of Lie algebroids: the tangent bundle $\tau:TM\to M$ with
$\rho=$identity, any integrable subbundle of $TM$ with the
inclusion as anchor map, $TP/G$ for $P(M,G)$ a $G-$principal
bundle, see [6].

Let $L:E\to R$ be a regular Lagrangian on the Lie algebroid
$(E,[,],\rho)$, that is $L$ is a smooth functions such that the
matrix with the entries
\begin{equation}
g_{ab}(x,y)=\dfrac{1}{2}\dfrac{\partial^2 L}{\partial y^a\partial
y^b},\tag{3.3}
\end{equation}
is of rank $m$.
Let $c:I\to M$, $I\subseteq\mathbb{R}$ be a curve
on $M$ and let $\widetilde{c}:I\to E$ be any curve on $E$ such
that $\pi\circ\widetilde{c}=c$. Denote by $\dot{\widetilde{c}}$
the vector field that is tangent to $\widetilde{c}$.

{\bf Definition 3.1.} {\it We say that $\widetilde{c}$ is {\bf
admissible} if
$$\pi_*(\dot{\widetilde{c}})=\rho(\widetilde{c}).$$}

In local charts on $M$ and $E$, we have $c(t)=(x^i(t))$,
$\widetilde{c}(t)=(x^i(t),y^a(t))$ and
$\dot{\widetilde{c}}(t)=\dfrac{dx^i}{dt}\dfrac{\partial}{\partial
x^i}+\dfrac{dy^a}{dt}\dfrac{\partial}{\partial y^a}$, $t\in I$.

It results

{\bf Lemma 3.1.} {\it The curve $\widetilde{c}$ is admissible if
and only if \begin{equation}
\dfrac{dx^i}{dt}(t)=\rho^i_a(x(t))y^a(t),\ \forall t\in
I.\tag{3.4} \end{equation}}

In [6], one associates to $L$ the Euler - Lagrange equations
\begin{equation}
\dfrac{d}{dt}\left(\dfrac{\partial L}{\partial y^a
}\right)=\rho^i_a\dfrac{\partial L}{\partial
x^i}+L^c_{ba}y^b\dfrac{\partial L}{\partial y^c},\tag{3.5}
\end{equation}
for $c(t)=(x^i(t),y^a(t))$ an admissible curve.

Expanding the derivative, using (3.3) and (3.4), we may put (3.5)
in the form
\begin{equation}
\dfrac{dy^a}{dt}+2G^a_L(x,y)=0,\tag{3.6}
\end{equation}
with the notation
\begin{equation}
G^a_L=\dfrac{1}{4}g^{ab}\left(\dfrac{\partial^2 L}{\partial
y^b\partial x^i}\rho^i_c y^c-\rho^i_b\dfrac{\partial L}{\partial
x^j}-L^c_{bd}y^d\dfrac{\partial L}{\partial y^c}\right).\tag{3.7}
\end{equation}

Let be $S=\rho^i_a(x)y^a\dfrac{\partial}{\partial
x^i}+Y^a\dfrac{\partial}{\partial y^a}$ a vector field on $E$,
where the coordinates $(Y^a(x,y))$ are not determined. We set for
convenience $Y^a=-2G^a$. Furthermore, under a change of
coordinates $(x^i,y^u)\to(\widetilde{x}^i,\widetilde{y}^a)$, the
coordinates $(X^i(x,y)=\rho^i_a(x)y^a),(G^a)$ have to change as
follows
\begin{equation}
\widetilde{X}^i=\dfrac{\partial\widetilde{x}^i}{\partial
x^j}(x)X^j,\tag{3.8}
\end{equation}
\begin{equation}
\widetilde{G}^a=M^a_b G^b-\dfrac{1}{2}\dfrac{\partial
M^a_b}{\partial x^i}y^b\rho^i_c y^c.\tag{3.9}
\end{equation}

Using (3.2) one easily sees that the coordinates $(X^i(x,y))$
 verify (3.8).

We say that the vector field $S$ as above is a semispray on $E$.
For more details on semisprays on $E$ we refer to [1].

Now we show that the function $(G^a_L)$ verifies (3.9) under a
change of coordinates on $E$.

We set
\begin{equation}
E_a=4g_{ab}G^b,\tag{3.10}
\end{equation}

 where
\begin{equation}
E_a=\dfrac{\partial^2 L}{\partial y^a\partial x^i}\rho^i_b
y^b-\rho^i_a\dfrac{\partial L}{\partial
x^i}-L^c_{ba}y^b\dfrac{\partial L}{\partial y^c}.\tag{3.11}
\end{equation}

Then we use (3.2) as well as the following equations:
$$\dfrac{\partial L}{\partial x^i}=\dfrac{\partial L}{\partial \widetilde{x}^j}\dfrac{\partial
\widetilde{x}^j}{\partial x^i}+\dfrac{\partial L}{\partial
\widetilde{y}^a}\dfrac{\partial M^a_c}{\partial x^i}y^c$$
$$\dfrac{\partial^2 L}{\partial y^a\partial x^i}=M^b_a\left(\dfrac{\partial^2 L}{\partial
y^b\partial\widetilde{x}^j}\dfrac{\partial\widetilde{x}^j}{\partial
x^i}+2\widetilde{g}_{db}\dfrac{\partial M^d_c}{\partial
x^i}y^c\right)+\dfrac{\partial
L}{\partial\widetilde{y}^d}\dfrac{\partial M^d_a}{\partial x^i}$$
$$L^c_{ab}M^e_c=M^c_a M^d_b\widetilde{L}^e_{cd}+\rho^k_a\dfrac{\partial M^e_b}{\partial
x^k}-\rho^k_b\dfrac{\partial M^e_a}{\partial x^k}$$ in order to
derive
\begin{equation}
E_a=M^b_a\widetilde{E}_b+2M^b_a\widetilde{g}_{bd}\dfrac{\partial
M^d_c}{\partial x^i}y^c\rho^i_d y^d.\tag{3.12}
\end{equation}

Using this in (3.10) one shows that $\widetilde{G}^a_L$ is related
to $G^a_L$ as in (3.9).

Thus we have proved

{\bf Theorem 3.1.} {\it Let $L$ be a regular Lagrangian on the Lie
algebroid $(E,[,],\rho)$. Then $L$ defines a semispray
$S_L=(\rho^i_a y^a)\dfrac{\partial}{\partial
x^i}-2G^a_L(x,y)\dfrac{\partial}{\partial y^a}$, where the
function $G^a_L$ are given by (3.7).}

{\bf Example 3.1.} Let $g_{ab}(x)$ be the coefficients of a
Riemannian metric in the Lie algebroid $(E,[,],\rho)$. Then
\begin{equation}
L(x,y)=g_{ab}(x)y^a y^b\tag{3.13}
\end{equation}
is a regular Lagrangian on $E$. The semispray associated to it is
determined by the functions
\begin{equation}
G^a=\dfrac{1}{2}g^{ab}\left(\dfrac{\partial g_{bc}}{\partial
x^i}\rho^i_d-\dfrac{1}{2}\dfrac{\partial g_{cd}}{\partial
x^i}\rho^i_b-L^e_{db}g_{ec}\right)y^c y^d.\tag{3.14}
\end{equation}

{\bf Example 3.2.} A more general example is provided by the
regular Lagrangians which are homogeneous of degree 2 in $(y^a)$.
By the Euler theorem one obtains
\begin{equation}
L(x,y)=g_{ab}(x,y)y^a y^b,\tag{3.15}
\end{equation}
where $(g_{ab}(x,y))$ are homogeneous functions of degree 0.

As $\dfrac{\partial}{\partial y^a}$ are homogeneous functions of
degree 1 and the derivative with respect to $(x^j)$ does not
affect the degree of homogeneity, it results that the coefficients
$(G^a)$ from (3.4) are homogeneous of degree 2 in $(y^a)$.  The
corresponding semispray is called a spray.

\section{Mechanical Lagrangian systems on algebroids}

Let $(E, [,],\rho)$ be a Lie algebroid.

\ms

{\bf Definition 4.1.} {\it A mechanical Lagrangian system with
external forces on the Lie algebroid $(E,[,],\rho)$ is
$\sum=(E,L,F)$ with $L$ a regular Lagrangian on $E$ and
$F=(F_a(x,y))$ a vertical covector field.}

Let be the functions
\begin{equation}
{\cal L}_a:=\dfrac{d}{dt}\left(\dfrac{\partial L}{\partial
y^a}\right)-\rho^i_a\dfrac{\partial L}{\partial
x^j}-L^c_{ba}y^b\dfrac{\partial L}{\partial y^c}\tag{4.1}
\end{equation}
defined on admissible curves on $E$.

Then the equalities ${\cal L}_a=0$ represent the Euler - Lagrange
equations associated to $L$.

We assume that the evolution equations of the system $\sum$ are as
follows:
\begin{equation}
{\cal L}_a(x(t),y(t))=F_a(x(t),y(t)),\tag{4.2}
\end{equation}
for $\widetilde{c}(t)=(x(t),y(t))$ an admissible curve on $E$.

The equations (4.2) after some arrangements take the form
\begin{equation}
\dfrac{dy^a}{dt}+2G^a(x,y)=\dfrac{1}{2}F^a(x,y),\tag{4.3}
\end{equation}
where the functions $(G^a)$ are given by (3.7), $F^a=g^{ab}F_b$,
and the equations $\dfrac{dx^i}{dt}=\rho^i_a(x)y^a$ hold.

Thus the evolution equations of the system $\sum$ become
\begin{equation}
\begin{array}{l}\dfrac{dx^i}{dt}=\rho^i_a(x)y^a,\vspace{1.5mm}\\
\dfrac{dy^a}{dt}=-2\left(G^a-\dfrac{1}{4}F^a\right).\end{array}\tag{4.4}
\end{equation}

The solutions of this system may be regarded as the integral
curves of a semispray
\begin{equation}
S^*=\rho^i_a(x)y^a\dfrac{\partial}{\partial
x^i}-2G^*(x,y)\dfrac{\partial}{\partial y^a},\
G^{*a}=G^a-\dfrac{1}{4}F^a.\tag{4.5}
\end{equation}

Indeed, $S^*$ is a semispray because it differs by the semispray
$S$ derived from $L$ by a vertical vector field.

\ms

{\bf Definition 4.2.} {\it We say that the mechanical Lagrangian
system $\sum$ is dissipative if $F_a(x,y)y^a\leq 0$ and that it is
strictly dissipative if $F_a(x,y)y^a\leq -\alpha y_a y^a$ with
$\alpha>0$ a constant and $y_a=g_{ab}y^b$.}

\ms

{\bf Theorem 4.1.} {\it Let be the mechanical Lagrangian system
$\sum$ with the evolution equations $(4.4)$. If it is dissipative
then its energy $E=y^a\dfrac{\partial L}{\partial y^a}-L$
decreases on the curves that are solutions of $(4.4)$. If
furthermore it is strictly dissipative its energy is strictly
decreasing on the curves solutions of $(4.4)$, assuming that these
have no singularities.}

\ms

\noindent{\it Proof.} Let be $\widetilde{c}(t)=(x^i(t),y^a(t))$ a curve
that is a solution of (4.4). Along this curve we have
$$\dfrac{dE}{dt}=\dfrac{dy^a}{dt}\dfrac{\partial L}{\partial
y^a}+y^a\dfrac{d}{dt}\left(\dfrac{\partial L}{\partial
y^a}\right)-\dfrac{\partial L}{\partial
x^i}\dfrac{dx^i}{dt}-\dfrac{\partial L}{\partial
y^a}\dfrac{dy^a}{dt}=$$ $$=y^a{\cal L}_a(x,y)=y^a F_a(x,y).$$

The last equality is based on (4.2) and to obtain the previous one
the equations
\begin{equation}
L^c_{ab}y^a y^b=0,\tag{4.6}
\end{equation}
have been used.

If the system $\sum$ is dissipative we have $\dfrac{dE}{dt}\leq 0$
and if it is strictly dissipative we have
$\dfrac{dE}{dt}\leq-\alpha y_a y^a<0$, q.e.d.

Now, we show that if $\sum$ is dissipative we can associate to it
a Lyapunov function.

Let $(x^i_0,y^a_0)$ be an equilibrium point of $S^*$.

If $\rho$ is injective this has the form $(x^i_0,0)$ with
$G^{*a}(x^i,0)=0$, a condition that is verified if $S^*$ is a
spray.

Assume that $(x^i_0, y^a_0)$ is a minimum point for the energy $E$
and set $\widetilde{E}(x,y)=E(x,y)-E(x_0,y_0)$.

We have
\begin{equation}
\widetilde{E}(x_0,y_0)=0,\ \widetilde{E}(x,y)>0.\tag{4.7}
\end{equation}

Let us denote by ${\cal L}_{S^*}$ the Lie derivative with respect
to $S^*$.

We have: ${\cal L}_{S^*}(E)=\rho^i_a y^a\dfrac{\partial
E}{\partial x^i}-2G^a\dfrac{\partial E}{\partial
y^a}+\dfrac{1}{2}F^a\dfrac{\partial E}{\partial y^a}$.

But $\dfrac{\partial E}{\partial y^a}=2g_{ab}y^b:=2y_a$. Hence
${\cal L}_{S^*}(E)=y^a E_a4G^a y_a+y_a F^a$, where $E_a$ was
defined in (3.11). Again (4.6) was used.

Looking at the connection between $E_a$ and $G^a$ it comes out
that the first two terms in the expression of ${\cal L}_{S^*}(E)$
cancel and so we have
\begin{equation}
{\cal L}_{S^*}(E)=y_a F^a\leq 0,\tag{4.8}
\end{equation}
since $\sum$ is dissipative.

Thus the function $\widetilde{E}$ is a Lyapunov function for $S^*$
in the equilibrium point $(x^i_0,y^a_0)$ but we can not conclude
that this point is stable.

In order to do so we need to introduce a Riemannian metric on $E$
and to prove that $S^*$ is complete with respect to that metric.
For details see [4].

For $E=TM$ endowed with a regular Lagrangian a Sasaki type metric
can be considered but that construction does not work except if
the algebroid $(E,[,],\rho)$ is endowed with a nonlinear
connection.

\newpage

\def\v{\varphi}
\def\o{\omega}
\def\g{\gamma}
\def\a{\alpha}
\def\b{\beta}
\def\l{\lambda}
\def\d{\delta}
\def\p{\partial}
\def\t{\theta}
\def\vp{\varepsilon}
\def\s{\sigma}
\def\D{\Delta}
\def\k{\kappa}
\def\G{\Gamma}
\def\O{\Omega}
\def\L{\Lambda}
\def\cal#1{\mathcal{#1}}

\def\spa{\vspace*{-.3cm}\\ }

\def\dd{\displaystyle}
\def\bk{\bigskip}

\runningauthor={M. ANASTASIEI}
\runningtitle={A GENERALIZATION OF MYERS THEOREM}
\noindent
\baselineskip 8pt
\noindent{\footnotesize{Analele {\c S}tiin{\c t}ifice ale Universit{\u a}{\c t}ii ``Al.~I.~Cuza'' Ia{\c s}i}}
\hfill\break
{\footnotesize{Mat. (N.S.) 53 (2007), suppl. 1, 33--40}}
\vskip 2cm
\baselineskip 11.5pt plus .15pt
\centerline{\bf\Large A GENERALIZATION OF}
\vskip .2cm
\centerline{\bf\Large MYERS THEOREM}
\vskip .5cm
\centerline{\bf {\footnotesize{BY}}}
\vskip .5cm
\centerline{\bf {\footnotesize{M. ANASTASIEI}\footnote{This work was partially supported by grant CNCSIS 1158/2007,
Romania}}}
\vskip 1cm

{\footnotesize{\it Dedicated to Academician Radu Miron
at his 80th anniversary}}

\begin{abstract} The Myers theorem extracts some
topological properties of a Riemannian manifold $(M,g)$ from the
assumptions that its Ricci curvature is uniformly bounded below by
a positive constant. The theorem was extended to Finsler
manifolds. Proofs of it can be seen in [1], Ch. 7, [3] Ch.7. In
1979, {\sc Galloway} ([2]) obtains the same topological properties
of $(M,g)$ assuming a weaker boundedness hypothesis on the Ricci
curvature.

In this paper we show that the version of Myers theorem due to
Galloway holds also for Finsler manifolds. So, a positive
answer to a problem posed by B. Suceav\u{a} in a private
communication is provided.

We mention that B. Suceav\u{a} proved a Myers type theorem in the
spirit of [2] for almost Hermitian manifolds [4].

Our proof is obtained by modifying some points in the proof from
[1] and by checking that some facts proved in [2] for Riemannian
manifolds hold also for Finsler manifolds.

{\bf Mathematics Subject Classification 2000:} 53C60.

{\bf Key words:} Finsler manifolds, Ricci scalar, Myers
theorem.
\end{abstract}

\setcounter{section}{0}

\section{Preliminaries} We shall use the notations and
the terminology from [1] without comments.

Let $(M,F)$ be a Finsler manifold. The Finsler structure $F$ is a
function $F:TM\to [0,\infty)$, $(x,y)\to F(x,y)$ which is
$C^{\infty}$ on the slit tangent bundle $TM\backslash 0$,
positively homogeneous in $y$ and whose Hessian matrix $g_{ij}:=
\dfrac{1}{2}\dfrac{\p^{2}F^{2}}{\p y^{i}\p y^{i}}$ is
positive-definite at every point of $TM\backslash 0.$

The Chern connection is a linear connection in the pull-back
bundle $\pi^{*}TM$ over $TM\backslash 0$, where $\pi:TM\to M$ is
the natural projection. It is only $h$-metrical and it has two
curvatures $R_{j} \ ^{i}\ _{kh}$, $P_{j} \ ^{i}\  _{kh}$.

Let be $y$ a non zero element of $T_{x}M$. Then $g(x,y) =
g_{ij}(x,y)dx^{i}\otimes dx^{j}$ is an inner product which is used
to measure lengths and angles in $T_{x}M$.One calls $y$ a flagpole
of the flag (a plane in $T_{x}M$) spanned by $l=\dfrac{y}{F(x,y)}$,
and another unit vector $V$ which is orthogonal to the flagpole.

The flag curvature is then given as
$$
K(x,y,l\wedge V):= V^{i}(l^{j}R_{jikh}l^{h})V^{k}=:
V^{i}R_{ik}V^{k}. \leqno(1.1)
$$
The raising and lowering of indices is made by using $g^{ij}$ and
$g_{ij}$, respectively. Sometimes, the flag curvature is denoted
simply $K(l,V)$. If $V$ is not a unit vector, then we have
$g_{(x,y)}(V,V)K(l,V) = V^{i}R_{ik}V^{i}$. Let $\{l, e_{\a},
\a=1,\ldots,n-1\}$ be a $g$-orthonormal basis for the fiber of
$\pi^{*}TM$ over the point $(x,y)\in TM\backslash 0.$ With respect
to it one has $K(x,y,l\wedge e_{\a})=R_{\a\a}$. The Ricci scalar
denoted by Ric is
$$
Ric:= \sum_{\a=1}^{n-1}K(x,y,l\wedge e_{\a}) =
\sum_{\a=1}^{n-1}R_{\a\a}. \leqno(1.2)
$$
In any basis one gets
$$
Ric = g^{ik}R_{ik}. \leqno(1.3)
$$
The Ricci tensor is defined as follows
$$
Ric_{jk} = \frac{1}{2}\frac{\p^{2}(F^{2}Ric)}{\p y^{j}\p y^{k}}
\leqno(1.4)
$$
and one shows that
$$
Ric = l^{j}l^{k}Ric_{jk}. \leqno(1.5)
$$
Equivalently,
$$
Ric(x,y) = \frac{1}{F^{2}(x,y)}[y^{i}y^{j}Ric_{ij}]. \leqno(1.6)
$$

If $(M,F)$ has constant flag curvature $c$, then
$$
Ric = (n-1)c,\;\; Ric_{jk} = (n-1)cg_{jk}. \leqno(1.7)
$$

Let $\s (t), 0\leq t \leq L,$  be a unit geodesic with velocity
field $T$. One abbreviates $g_{(\s,T)}$ by $g_{T}$.

For a vector field $W(t):= W^{i}(t)\dfrac{\p}{\p x^{i}}$ along
$\s$, the expression,
$$
D_{T}W =
\left[\dfrac{dW^{i}}{dt}+W^{j}T^{k}(\Gamma_{jk}^{i}(G,T))\right]\dfrac{\p}{\p
x^{i}}\leqno(1.8)
$$
is called covariant derivative with reference vector $T$. The
formula 1.8 can be stated for any curve but for geodesics one has
$$
\frac{d}{dt}g_{T}(V,W) = g_{T}(D_{T}V, W) +g_{T}(V,D_{T}W)
\leqno(1.9)
$$
for any vector fields $V,W$ along $\s$.

Note that (1.9) holds for any curve if $V$ or $W$ is proportional
to $T$.

The constant speed geodesics are solutions of $D_{T}T=0,$ with
reference vector $T$.

One says that $W$ is parallel long $\s$ if $D_{T}W,$ with
reference vector $T$. Parallel transport (with reference vector
$T$) one defines on the standard way. By (1.9) the parallel
transport preserves $g_{T}$-lengths and angles.

For two continuous and piecewise $C^{\infty}$ vector fields $V$
and $W$ along $\s$ the index form is
$$
I(V,W) = \int_{0}^{L}[g_{T}(D_{T}V, D_{T} W) -
g_{T}(R(V,T)T,W)]dt. \leqno(1.10)
$$
Here all $D_{T}$ are calculated with reference vector $T$ and
$$R(V,T)T:=(T^{j}R_{jkh}^{i}T^{h})V^{k}\frac{\p}{\p x^{i}}$$ is
evaluated at the point $(\s,T)$.

The index form is bilinear and symmetric. We quote from [1] the
following facts

{\bf Proposition 1.1} [1, p. 174] {\it Let $\s(t) = \exp_{p}(tT),$
$0\leq t \leq r$ be a constant speed geodesic from $p = \s(0)$ to
$q = \s(r)$.}

{\it The following five statements are mutually equivalent:

$(a)$ The point $q$ is not conjugate to $p$ along $\s$.

$(b)$ Any Jacobi field that vanishes as both points $p$ and $q$
must be identically zero along $\s$.

$(c)$ Take the variation field of any variation of $\s$ by
geodesics. If it vanishes at $p$ and $q$, then it must be
identically zero along $\s$.

$(d)$ Given any $v \in T_{p}M$ and $w\in T_{q}M$, there exits a
unique Jacobi field $J$ along $\s$ that equals $v$ at $p$ and $w$
at $q$.

$(e)$ The derivative $\exp_{p*}$ of the exponent map $\exp_{p}$ is
nonsingular at the location $rT$ in $T_{p}M$.}

{\bf Proposition 1.2} [1, p. 182] {\it Let $\s(t), 0 \leq t\leq
r$ be a geodesic in a Finsler manifold $(M,F)$. Suppose no point
$\s(t), 0 < t\leq r$ is conjugate to $p:=\s(0)$. Let $W$ be any
piecewise $C^{\infty}$ vector field along $\s$ and let $J$ denote
the unique Jacobi field along $\s$ that has the same boundary
values as $W$. That is, $J(0) = W(0)$ and $J(r) = W(r)$. Then}
$$
I(W,W)>I(J,J). \leqno(1.11)
$$

Equality holds if and only if $W$ is actually a Jacobi field, in
which case the said $J$ coincides with $W$.

We close this Section by quoting, for the sake of comparison, the
Bonnet-Myers theorem from [1], p. 194:

\medskip
{\it Let $(M,F)$ be a forward geodesically complete connected
Finsler manifold of dimension $n$. Suppose its Ricci scalar has
the following uniform positive lower bound
$$
Ric\geq (n-1)\l>0.
$$
Equivalently, suppose its Ricci tensor satisfies
$y^{i}y^{i}Ric_{ij}(x,y)\geq (n-1)\l F^{2}(x,\\y)$ with $\l>0$.
Then:

\begin{itemize}
\item[$(1)$] Along every geodesic the distance between any two
successive conjugate points is at most $\dfrac{\pi}{\sqrt{\l}}$. In
other words, every geodesic with length $\dfrac{\pi}{\sqrt{\l}}$ or
longer must contain conjugate points.

\item[$(2)$] The diameter of $M$ is at most $\dfrac{\pi}{\sqrt{\l}}$.

\item[$(3)$] $M$ is in fact compact.

\item[$(4)$] The fundamental group $\pi(M,x)$ is
finite.\end{itemize}}

\section{A generalization of Bonnet - Myers theorem} Looking over
the proof of Bonnet-Myers theorem given in [1], p. 194-198 it comes
out that essential is a proof of its first statement.

Thus we give a more general form of this statement as follows:

{\bf Lemma 1.} {\it Let $\s(t)$, $0\leq t \leq L$ be a unit speed
geodesic with velocity field $T$. If
$$
Ric(T,T)\geq a+\dfrac{df}{dt},\; \mbox{for\; a\; constant}\; a>0
\leqno(2.1)
$$
and some function $f$ with $|f(t)|\leq C, C \geq 0,$ and
$$
L\geq \dfrac{\pi}{a}(c+\sqrt{c^{2} + a(n-1)}), \leqno(2.2)
$$
then $\s$ must contain conjugate points.}

{\it Remarks.}

(i) For $c=0$ and $a= (n-1)\l$, Lemma 2.1 reduces to the assertion
(1) of the Bonnet-Myers theorem.

(ii) The condition (2.1) on Ricci allows and negative values of
$Ric(T,T)$ along $\s$.

\bigskip
{\bf Proof.} Using the parallel transport with reference vector
$T$ one construct a moving frame $\{e_{i}(t)\}$ along $\s$ such
that

(i) Each $e_{i}$ is parallel along $\s$, that is $D_{T}e_{i} = 0,$

(ii) $\{e_{i}(t)\}$ is a $g_{T}$-ortonormal frame,

(iii) $e_{n}= T.$

Define $W_{\a}(t) = f_{\a}(t)e_{\a}(t)$ for some smooth functions
$f_{\a}$, $\a=1,2$, ..., $n-1$.

Fix a positive $r\leq L$ and consider the index from $I$ for
$\s(t), 0 \leq t \leq r$. By (1.10) we have
$$
I(W_{\a}, W_{\a}) = \int_{0}^{r}[\|D_{T}W_{\a}\|^{2} -
\|W_{\a}\|K(T,W_{\a})] dt,
$$
where the abbreviation $\|V\|:= g_{T}(V,V)$ was used and
$K(T,W_{\a})$ is the flag curvature evaluated at the point
$(\s(t),T)\in TM\backslash 0.$

As $D_{T}W_{\a}=\dfrac{df_{\a}}{dt}e_{\a}$, it results $\|D_{T}
W_{\a}|^{2} = |f_{\a}(t)|^{2}$. It is known that the flag
curvature does not depend on vectors spanning the flag. Thus we
have $K(T,W_{\a}) = K(T,e_{\a})$.

Using these facts, $I(W_{\a},W_{\a})$ takes the form
$$
I(W_{\a},W_{\a}) =
\int_{0}^{r}\left[\left(\frac{df_{\a}}{dt}\right)^{2} -
f_{\a}^{2}K(T,e_{\a})\right]dt.
$$
We take $f_{\a}(t) = \sin\dfrac{\pi t}{r}$ and we get
$$
I(W_{\a},W_{\a}) = \frac{\pi^{2}}{2r} -
\int_{0}^{L}\sin^{2}\frac{\pi t}{r}K(T,e_{\a})dt.
$$

Summing over $\a$ one obtains
$$
\sum_{\a}I(W_{\a},W_{\a}) = (n-1)\frac{\pi^{2}}{2r} -
\int_{0}^{r}Ric (T,T)\sin^{2}\frac{\pi t}{r}dt.
$$
By hypotheses, $-Ric (T,T)\leq -a-\dfrac{df}{dt}$. Hence
$$
\sum_{\a}I(W_{\a},W_{\a})\leq (n-1)\frac{\pi^{2}}{2r} -
\int_{0}^{r}\(a+\frac{df}{dt}\)\sin^{2}\frac{\pi t}{r}dt.
$$
An integration by parts gives first
$$
\sum_{\a}I(W_{\a}, W_{\a}) \leq (n-1)\frac{\pi^{2}}{2r} -
\frac{ar}{2}+\frac{\pi}{r}\int_{0}^{r}f(t)\sin \frac{2\pi t}{r}dt,
$$
and then using $|\sin u|\leq 1$ and $\|f(t)|\leq c,$ one finds
$$
\sum_{\a}I(W_{\a},W_{\a})\leq (n+1)\frac{\pi^{2}}{2r} -
\frac{ar}{2} +\pi c
$$
and we have $\sum_{\a}I(W_{\a},W_{\a})\leq o$ if $r\geq
\frac{\pi}{a}(c+\sqrt{c^{2}+a(n-1)})$ an inequality that holds for
$r=L$ by hypothesis. It follows that some $I(W_{\a},W_{\a})$ must
be nonpositive and let denote that $W_{\a}$ by $W$.

We proceed by contradiction. Suppose that $\s(t), 0 \leq t\leq
r=\frac{\pi}{a}(c+\sqrt{c^{2} + a(n-1)})$ contains no conjugate
points.

By Proposition 1.1, the vector field $W$, with $W(0) = W(r) = 0,$
can not be a Jacobi field since is nowhere zero on $(0,r)$. And by
the same Proposition 1.1 the unique Jacobi field which vanishes at
the endpoints of $\s(t)$, $0\leq t\leq r$ is identically zero
field. By Proposition 1.2 we have $0=I(J,J) <I(W,W)\leq 0$, which
is a contradiction and lemma is proved. In combination with
Theorem 7.5.1 from [1], Lemma 1 tell us that the said geodesic
$\s$ minimizes arc length among ``nearly'' piecewise $C^{\infty}$
curves from $p=\s(0)$ to $q = \s(r)$, $r=\dfrac{\pi}{a}(c+\sqrt{c^{2} +a(n-1)})$. The following two
consequences of this Lemma cover the content of the Bonnet-Myers
theorem.

{\bf Theorem 1.} {\it Let $(M,F)$ be a forward geodesically
complete connected Finsler manifold. Suppose there exists
constants $a>0$ and $c\geq 0$ such that for every pair of points
in $M$ and minimal geodesic $\s$ joining those points having unit
tangent vector $T$,  the Ricci curvature satisfies
$$
Ric(T,T) \geq a+\frac{df}{dt}\; \mbox{along}\;\s
$$
where $f$ is some function of arclength $t$ satisfying $|f(t)|\leq
c$ along $\s$. Then $M$ is compact and its $\mbox{diam}\, (M)\leq
\frac{\pi}{a}(c+\sqrt{c^{2}+a(n-1)})$}.

\ms

\noindent{\it Proof.} Since $M$ is forward geodesically complete, by the
Hopf-Rinow theorem any pair of points in $M$ can be joined by a
minimal geodesic. By Lemma 1, such a geodesic must have the length
less than or equal with $\dfrac{\pi}{a}(c+\sqrt{c^{2}+a(n-1)})$.
Thus $\mbox{diam}\, (M)\leq \dfrac{\pi}{a}(c+\sqrt{c^{2}+ a(n-1)})$
and so $M$ is bounded. Using again the Hopf-Rinow theorem one
deduces that $M$ is compact. $\hfill\square$

{\bf Theorem 2.} {\it Let $(M,F)$ be a forward geodesically
complete connected Finsler manifold. Suppose there exist constants
$a>0$ and $c\geq 0$ such that for every pair of points in $M$
{\rm(}not necessarily distinct{\rm)} and geodesic $\s$ with unit
tangent vector $T$ joining these points, the Ricci curvature
satisfies (2.1) where $f$ is some function of the arclength $t$
satisfying $|f(t)|\leq c$ along $\s.$ Then the universal covering
manifold of $M$ is compact, with diameter bounded by
$\dfrac{\pi}{a}(c+\sqrt{c^{2} + a(n-1)})$, and hence the
fundamental group of $M$ is finite.}

\ms
\noindent{\it Proof.} Let $\widetilde{M}$ be the universal covering
manifold of $M$ with the universal covering map
$p:\widetilde{M}\to M$. In [1] p. 197 one proves that $p$ endows
$\widetilde{M}$ with the same local geometry as $M$. Repeating
word by word the proof of Theorem 1.3 from [2] it comes out that
$\widetilde{M}$ satisfies the hypothesis of Theorem 2.1, hence it
is compact. It follows its closed subset $p^{-1}(x)$ is compact
and being discrete is finite. Since $\pi_{1}(M,x)$ is bijective
with $p^{-1}(x)$  it is itself finite.$\hfill\square$

\noindent{\footnotesize{\it Received: 15.X.2007}}

\newpage

\runningauthor={M. ANASTASIEI}
\runningtitle={SEMISPRAYS ON LIE ALGEBROIDS. APPLICATIONS}
\noindent
\baselineskip 8pt
\noindent{\footnotesize{Tensor (N.S.)}}
\hfill\break
{\footnotesize{69, 190--198, 2008}}
\vskip 2cm
\baselineskip 11.5pt plus .15pt
\centerline{\bf\Large SEMISPRAYS ON LIE ALGEBROIDS.}
\vskip .2cm
\centerline{\bf\Large APPLICATIONS$^*$}
\vskip .5cm
\centerline{\bf {\footnotesize{BY}}}
\vskip .5cm
\centerline{\bf {\footnotesize{M. ANASTASIEI}\footnote{Communicated at The 8th Conference of Tensor Society on Differential Geometry, Functional and Complex Analysis,
Informatics and their Applications, held at Varna, Bulgaria, August, 22-26, 2005}}}
\vskip 1cm

{\footnotesize{Dedicated to Prof. Dr. Tomoaki Kawaguchi at his 70th anniversary}}

\setcounter{section}{0}
\section*{Introduction}

In any Lagrangian formalism for Lie algebroids (see A. Weinstein,
[9]\footnote{Numbers in brackets refer to the references at the
end of the paper}, E. Martinez, [5]) , the notion of semispray on
a Lie algebroid has a central place. If one looks at various
definitions of a semispray on a Lie algebroid (see M. Anastasiei,
[1]) it comes out that in defining a semispray the anchor map only
is used. In the other words, as it will be shown in this paper
(Section 2) the notion of semispray can be considered also on the
anchored vector bundles. Moreover, we will show in Section 3 that
the set of the anchored vector bundle is the largest with this
property. Of course, this set includes the set of Lie algebroids
and on a Lie algebroid the assertion that any regular Lagrangian
on it induces a semispray holds as in the tangent bundle case. We
will prove it in Section 4 (see also M. Anastasiei, [ 1]). We close
the paper with an application of semisprays to the mechanical
systems on a Lie algebroid (see M. Anastasiei, [3]). The first
Section is devoted to some preliminaries on vector bundles.

\section{Preliminaries on vector bundles}

Let $\xi=(E,\pi,M)$ be a vector bundle of rank $m$. Here $E$ and
$M$ are smooth i.e. $C^\infty$ manifolds with dim$M=n$,
dim$E=n+m$, and $\pi:E\to M$ is a smooth submersion. The fibres
$E_x=\pi^{-1}(x)$, $x\in M$ are linear spaces of dimension $m$
which are isomorphic with the type fibre $\mathbb{R}^m$.

Let $\{(U_\alpha,\psi_\alpha)\}_{\alpha\in A}$ be an atlas on $M$.
A vector bundle atlas is
$\{(U_\alpha$, $\varphi_\alpha$, $\mathbb{R}^m)\}_{\alpha\in A}$ with
the bijections $\varphi_\alpha: \pi^{-1}(U_\alpha)\to
U_\alpha\times\mathbb{R}^m$ in the form $\varphi _\alpha(u)=(\pi
(u),\varphi_{\alpha,\pi (u)})$, where
$\varphi_{\alpha,\pi(u)}:E_{\pi(u)}\to\mathbb{R}^m$ is a
bijection. The given atlas on $M$ and a vector bundle atlas
provide an atlas
$\{(\pi^{-1}(U_\alpha)$, $\Phi_\alpha)\}_{\alpha\in A}$ on $E$.\\
Here $\Phi_\alpha:\pi^{-1}(U_\alpha)\to\phi_\alpha(U_\alpha)\times\mathbb{R}^m$
is the bijection given by
$\phi_\alpha(u)=$ $(\psi_\alpha(\pi(u))$, $\varphi_{\alpha,\pi
(u)}(u))$. For $x\in M$, we put
$\psi_\alpha(x)=(x^i)\in\mathbb{R}^n$ and if
$(U_\beta,\psi_\beta)$ is another local chart such that $x\in
U_\alpha\cap U_\beta\neq\phi$, we set
$\psi_\beta(x)=\widetilde{x}^i$ and then
$\psi_\beta\circ\psi_\alpha^{-1}$ has the form
\begin{equation}
\widetilde{x}^i=\widetilde{x}^i(x^1,\cdots,x^n),\
{\rm{rank}}\left(\dfrac{\partial \widetilde{x}^i}{\partial x^j
}\right)=n.\tag{1.1}
\end{equation}

Let $(e_a)$ be the canonical basis of $\mathbb{R}^m$. Then
$\varphi^{-1}_{\alpha,x}(e_a):=\varepsilon_a(x)$ is a basis of
$E_x$ and $u\in E_x$ has the form $u=y^a\varepsilon_a(x)$.

We take $(x^i,y^a)$ as coordinates on $E$. For the bundle chart
$(U_\beta,\varphi_\beta,\mathbb{R}^m)$ we put
$\varphi^{-1}_{\beta,x}(e_a)=\widetilde{\varepsilon}_a(x)$ and
then $u=\widetilde{y}^a\widetilde{\varepsilon}_a(x)$. If we set
$\varepsilon_a(x)=M^b_a(x)\widetilde{\varepsilon}_b$ with
rank$(M_a^b(x))=m$ it follows that $\widetilde{y}^a=M^a_b(x)y^b$.
Thus $\phi_\beta\circ\phi_\alpha^{-1}$ has the form
\begin{equation}
\begin{array}{l}\widetilde{x}^i=\widetilde{x}^i(x^1,\cdots,x^n), {\rm{rank}}
\left(\dfrac{\partial\widetilde{x}^i}{\partial
x^j}\right)=n\vspace{1.5mm}\\ \widetilde{y}^a=M^a_b(x)y^b,
{\rm{rank}}(M^a_b(x))=m.
\end{array}
\tag{1.2}
\end{equation}

The indices $i,j,k,...$ and $a,b,c...$ will take the values
$1,2,...n$ and $1,2,...m$, respectively. The Einstein convention
on summation will be used.

We denote by ${\cal F}(M), {\cal F}(E)$ the ring of real functions
on $M$ and $E$ respectively, and by $\chi(M)$, respectively
$\Gamma(E),\chi(E)$ the module of sections of the tangent bundle
of $M$, respectively of the bundle $\xi$ and of the tangent bundle
of $E$. On $U_\alpha$, the vector fields
$\left(\partial_k:=\dfrac{\partial}{\partial x^k}\right)$ provide
a local basis for $\chi(U_\alpha)$. The sections
$\varepsilon_a:U_a\to \pi^{-1}(U_\alpha)$,
$\varepsilon_a(x)=\varphi^{-1}_{\alpha,x}(e_a)$ provide a basis
for $\Gamma(\pi^{-1}(U_\alpha))$ and a section $A:U_\alpha\to
\pi^{-1}(U_\alpha)$ will take the form
$A(x)=A^a(x)\varepsilon_a(x)$, $x\in U_\alpha$.

Let $\xi^*=(E^*,p^*,M)$ be the dual of the vector bundle $\xi$. We
may also consider the tensor bundle $T^r_s(E)$ over $E$. The set
of sections $\Gamma(T^r_s(E))$ are ${\cal F}(M)-$modules for any
natural numbers $r,s$. On the sum $\oplus_{r,s}\Gamma(T^r_s(E))$ a
tensor product can be defined and one gets a tensor algebra
$T(E)$. For the tangent bundle $(TM,\tau,M)$ this reduces to the
tensor algebra of the manifold $M$. The tensor algebra of the
manifold $E$ could be also involved. Its elements are sections in
${\cal T}^r_s(TE)$. The tensorial algebra of $E$ contains the
subset of $d-$tensor fields on $E$. For a general definition of
these tensor fields we refer to [6], Ch. III. Shortly, these
tensor fields are defined by components depending on $(x^i,y^a)$
and transforming as tensors by a change of coordinates but with
the matrices $\left(\dfrac{\partial\widetilde{x}^i}{\partial
x^j}\right)$ and $(M^a_b(x))$ and their inverses, only. Notice
that in the law of transformation of a tensor field on $E$ could
appear also the matrix $\left(\dfrac{\partial M^a_b(x)}{\partial
x^i}y^b\right)$.

A large class of examples is provided by the sections in the
vertical bundle over $E$. We recall that the vertical bundle
$VE\to E$ is the union of the fibres $V_uE=\ker\pi_{*,u}$ over
$u\in E$, where $\pi_{*,u}$ is the differential of $\pi$. A basis
of local section of $VE\to E$ is given by
$\left(\left.\dfrac{\partial}{\partial y^a}\right|_u\right)$ and
its dual is $dy^a|_u$. The local components of any element in
$\Gamma(T^r_s(VE))$, transform under a change of coordinates on
$E$ with the matrix $(M^a_b(x))$ and its inverse $(W^a_b)$. We
call such an element a vertical tensor field.

\section{Semisprays for anchored vector bundles}

A vector bundle $\xi=(E,\pi,M)$ is called anchored (with the
tangent bundle $TM$) if there exists a v.b. morphism $\rho:E
\mapsto M$ called the anchor map.

The v.b. morphism $\rho$ induces a ${\cal F}(M)$ - module
homomorphism from $\Gamma(E)\mapsto\chi(M)$ denoted also by
$\rho$.

Locally, we set
\begin{equation}
\rho(\varepsilon_a)=\rho^i_a\dfrac{\partial}{\partial x^i}.
\tag{2.1}
\end{equation}

A change of local charts implies
\begin{equation}
\widetilde{\rho}^i_a=W^b_a\rho^j_b\dfrac{\partial\widetilde{x}^i}{\partial
x^j},\tag{2.2}
\end{equation}
where $W^b_a$ is the inverse of the matrix $(M^a_b)$.

\medskip

{\bf Examples.}
\begin{enumerate}
\item A trivial example of anchored v.b. is the tangent bundle
itself with the identity mapping as anchor.

\item A less trivial example is a provided by a subbundle of the
tangent bundle i.e. a distribution $D$ on $M$ with the inclusion
mapping as anchor. Let be $dim D = m <n $ and $(X_1,...,X_m)$ a
base of local sections of $D$. Then we may write $X_a =
X_a^i\dfrac{\partial}{\partial x^i}$ with $rank(X_a^i)=m$.
The anchor is given by
\begin{equation}
\rho (X_a)=X_a^i \dfrac{\partial }{\partial x^i},\tag{2.3}
\end{equation}

\item Let $P$ be a principal $G-$ bundle of projection $p$ over
$M$. Then $TP/G$ is a vector bundle over $M$ whose sections are
the $G-$ invariant vector fields on $P$. The derivative $p_* : TP
\mapsto TM$ passes to a mapping from $TP/G \mapsto TM$ which is
the anchor.
\end{enumerate}

We recall that a semispray $S$ for the tangent bundle $\tau:TM\to
M$ is a vector field on $TM$ which at the same time is a section
in the vector bundle $\tau_*:TTM\to TM$, that is we have
$\tau_{TM}(S(u))=u$ and $\tau_{*,u}(S(u))=u$, $\forall u\in TM$,
where $\tau_{TM}$ is the vector bundle projection $TTM\to TM$. It
follows that $\tau_{*,u}(S(u))=\tau_{TM}(S(u))$, $\forall u\in
TM$.

This equation suggests the following

{\bf Definition 2.1.} {\it Let $\xi=(E,\rho,M)$ be a an anchored
v.b. with the anchor $\rho$. A vector field $S$ on $E$ will be
called a semispray if
\begin{equation}
\pi_{*,u}(S(u))=(\rho\circ\tau_E)(S(u)),\ \forall u\in E\tag{2.4}
\end{equation}
where $\tau_E:TE\to E$ is the natural projection.}

Let $c:I\to M$, $I\subseteq\mathbb{R}$ be a curve on $M$ and let
$\widetilde{c}:I\to E$ be any curve on $E$ such that
$\pi\circ\widetilde{c}=c$. Denote by $\dot{\widetilde{c}}$ the
vector field that is tangent to $\widetilde{c}$.

{\bf Definition 2.2.} {\it We say that $\widetilde{c}$ is {\bf
admissible} if
$$\pi_*(\dot{\widetilde{c}})=\rho(\widetilde{c}).$$}

In local charts on $M$ and $E$, we have $c(t)=(x^i(t))$,
$\widetilde{c}(t)=(x^i(t),y^a(t))$ and
$\dot{\widetilde{c}}(t)=\dfrac{dx^i}{dt}\dfrac{\partial}{\partial
x^i}+\dfrac{dy^a}{dt}\dfrac{\partial}{\partial y^a}$, $t\in I$.

It results

{\bf Lemma 2.1.} {\it The curve $\widetilde{c}$ is admissible if
and only if
\begin{equation}
\dfrac{dx^i}{dt}(t)=\rho^i_a(x(t))y^a(t),\ \forall t\in
I.\tag{2.5}
\end{equation}}

Again in local charts, let be $S=X^i\dfrac{\partial}{\partial x^i}+Y^a\dfrac{\partial}{\partial y^a}$ a vector field on $E$.

This is a semispray if and only if
\begin{equation}
X^i(x,y)=\rho^i_a(x)y^a.\tag{2.6}
\end{equation}

Thus the coordinates $(Y^a(x,y))$ are not determined. We set for convenience $Y^a=-2G^a$. Furthermore, under a change of
coordinates $(x^i,y^u)\to(\widetilde{x}^i,\widetilde{y}^a)$, the coordinates $(X^i),(G^a)$ have to change as follows:
\begin{equation}
\widetilde{X}^i=\dfrac{\partial\widetilde{x}^i}{\partial
x^j}(x)X^j,\tag{2.7}
\end{equation}
\begin{equation}
\widetilde{G}^a=M^a_b G^b-\dfrac{1}{2}\dfrac{\partial
M^a_b}{\partial x^i}y^b\rho^i_c y^c.\tag{2.8}
\end{equation}

Using (2.2) one easily sees that the coordinates $(X^i(x,y))$ given by (2.6) verify (2.7).

Concluding, we have

{\bf Theorem 2.1.} {\it A vector field $S=(\rho^i_a y^a)\dfrac{\partial}{\partial x^i}-2G^a\dfrac{\partial}{\partial
y^a}$ on $E$ is a semispray if and only if the coordinates $(G^a)$ transform by (2.8).}

The integral curves of $S$ are given by the system of differential
equations
\begin{equation}
\dfrac{dx^i}{dt}=\rho^i_a(x)y^a,\
\dfrac{dy^a}{dt}+2G^a(x,y)=0.\tag{2.9}
\end{equation}

It comes out that these curves are all admissible. The converse is
also true, that is we have

{\bf Theorem 2.2.} {\it A vector field on $E$ is a semispray if
and only if all its integral curves are admissible.}

{\bf Remark 2.1.} The characterization of a semispray provided by
the Theorem 3.2 was taken by A. Weinstein, [9], as definition for
a semispray on Lie algebroids.

{\bf Remark 2.2.}

(i) Let us assume that $\rho=0$. Then the admissible curves are
all curves from the fibre $E_{x_0}$, $x_0(x^i_0)\in M$. The
integral curves of a semispray $S$ are given by the equations
$\dfrac{dy^a}{dt}+2G^a(x_0,y)=0$.

(ii) For a distribution $D$ on $M$ the condition (2.5) tell us
that the tangent vector field $\dfrac {dc}{dt}= y^a(t)X_{a}(c(t),$
that is $\dfrac {dc}{dt}$ is a section in the vector subbundle
$D$. In other words the admissible curves are in this case all the
curves that are tangent to the distribution $D$. See also M.
Anastasiei, [4].

Let $\widehat{S}$ be another semispray on $E$. Then
$\widehat{S}=(\rho^i_a y^a)\dfrac{\partial}{\partial
x^i}-2\widehat{G}^a\dfrac{\partial}{\partial y^a}$, where the
functions $(\widehat{G}^a(x,y))$ have to satisfy (2.8) under a
change of coordinates on $E$. It follows that
$\widehat{S}-S=2(G^a-\widehat{G}^a)\dfrac{\partial}{\partial y^a}$
and the functions $D^a=G^a-\widehat{G}^a$ transform by the rule
\begin{equation}
\widehat{D}^a=M^a_b D^b.\tag{2.10}
\end{equation}

By (2.10) we have that $D^a\dfrac{\partial}{\partial y^a}$ is a
vertical vector field.

So we have proved

{\bf Theorem 2.3.} {\it Any two semisprays on $E$ differ by a
vertical vector field on $E$.}

\section{Homogeneous semisprays(sprays)on anchored vector bundles}

For every real member $c > 0$ let $h_c$ denote the homothety $E\to
E,$ given by $u\to cu, u\in E$. A semispray $S$ on $E$ is called
a spray if

\begin{equation}
S(h_{c}(u)) = ch_{c,*} S(u). \tag{$H$}
\end{equation}

Locally,$h_c :(x^i,y^a) \mapsto (x^i,cy^a)$ and the condition
$(H)$ is equivalent with
\begin{equation}
G^a(x,cy) = c^2G^a(x,y).\tag{$H_0$}
\end{equation}

Let be $C=y^a\dfrac{\partial}{\partial y^a}$ the Liouville vector
field on $E$.

Using the Euler theorem on homogeneous functions one verifies that
$(H_0)$ is equivalent with
\begin{equation}
[C,S] = S. \tag{3.1}
\end{equation}

We notice that if we assume that $S$ is smooth on $E$ the condition $(H_0)$ reduces to the assertion that  $G^a$ are
homogeneous polynomials of degree 2 in $y^a$  because of

{\bf Lemma 3.1 ([8]).} {\it Let $V$ and $V'$ be linear spaces and $f:V\mapsto V'$ a mapping that is at least $r>0$ times
differentiable at $0\in V$ and positively homogeneous of degree $r$. Then $f$ is a homogeneous polynomial of degree $r$.}

When $S$ is smooth only on $E\setminus\{0\}$ the condition $(H_0)$ is in use.

As we have seen till now, given an anchored v.b. we may find in principle a semispray by pointing out a set of functions $(G^a)$
subject to (2.8). If someone tries to define a semispray on any vector bundle it is reasonable to try to define first a spray
since this has a simpler form. Thus he will start with a vector field  $S_0$ on $E$ that verifies the condition $(H_0)$.

If $S_0 = X^i(x,y)\dfrac{\partial}{\partial x^i}+Y^a(x,y)\dfrac{\partial}{\partial y^a}$, it will result that
$(X^i(x,y))$ are linear functions in $y^a$, that is $X^i=\rho_a^i(x)y^a$ and $(Y^a(x,y))$ are homogeneous polynomials of
degree 2 in $y^a$. The map $\pi_* \circ S_0$ carries a section $y^a\varepsilon _a$ to $\rho_a^i(x)y^a\dfrac{\partial }{\partial
x^i} $ i.e. it defines a morphism $E\mapsto TM$. As $\tau_E\circ S_0=id_E$ holds, the condition (2.1) is fulfilled, i.e. $S_0$
is a spray.

Concluding, if one wishes the extension of the notion of semispray to vector bundles, one has to assume that vector bundle is
anchored. In the other word, the class of anchored v.b. is the largest in which the notion of semispray can be considered. It
contains the class of Lie algebroids.

\section{A  semispray derived from a Lagrangian on a Lie algebroid}

A vector bundle $\xi=(E,\pi,M)$ is called a Lie algebroid if it has the following properties:

\begin{enumerate}
\item The space of sections $\Gamma(\xi)$ is endowed with a Lie
algebra structure $[,]$;
\item There exists a bundle map $\rho:E\to TM$ (called the {\it anchor
map}) which induces a Lie algebra homomorphism (also denoted by
$\rho$) from $\Gamma(\xi)$ to $\chi(M)$.
\item For any smooth functions $f$ on $M$ and any sections
$s_1,s_2\in\Gamma(\xi)$ the following identity is satisfied
$$[s_1,fs_2]=f[s_1,s_2]+(\rho(s_1)f)s_2.$$
\end{enumerate}

Locally, we set
$$\rho(\varepsilon_a)=\rho^i_a\dfrac{\partial}{\partial x^i},\
[\varepsilon_a,\varepsilon_b]=L^c_{ab}\varepsilon_c.$$

Let $L:E\to R$ be a regular Lagrangian on the Lie algebroid
$(E,[,],\rho)$, that is $L$ is a smooth functions such that the
matrix with the entries
\begin{equation}
g_{ab}(x,y)=\dfrac{1}{2}\dfrac{\partial^2 L}{\partial y^a\partial
y^b},\tag{4.1}
\end{equation}
is of rank $m$.

In [9], one associates to $L$ the Euler - Lagrange equations
\begin{equation}
\dfrac{d}{dt}\left(\dfrac{\partial L}{\partial y^a
}\right)=\rho^i_a\dfrac{\partial L}{\partial
x^i}+L^c_{ba}y^b\dfrac{\partial L}{\partial y^c},\tag{4.2}
\end{equation}
for $c(t)=(x^i(t),y^a(t))$ an admissible curve.

Expanding the derivative, using (4.1) and (3.4), we may put (4.2)
in the form
\begin{equation}
\dfrac{dy^a}{dt}+2G^a_L(x,y)=0,\tag{4.3}
\end{equation}
with the notation
\begin{equation}
G^a_L=\dfrac{1}{4}g^{ab}\left(\dfrac{\partial^2 L}{\partial
y^b\partial x^i}\rho^i_c y^c-\rho^i_b\dfrac{\partial L}{\partial
x^j}-L^c_{bd}y^d\dfrac{\partial L}{\partial y^c}\right).\tag{4.4}
\end{equation}

In [1]we have shown by a direct calculation that the function
$(G_L^a)$ verify (2.8) under a change of coordinates.

In the other words  we have proved

{\bf Theorem 4.1.} {\it Let $L$ be a regular Lagrangian on the Lie
algebroid $(E,[,],\rho)$. Then $L$ defines a semispray
$S_L=(\rho^i_a y^a)\dfrac{\partial}{\partial
x^i}-2G^a_L(x,y)\dfrac{\partial}{\partial y^a}$, where the
function $G^a_L$ are given by $(4.4)$.}

{\bf Example 4.1.} Let $g_{ab}(x)$ be the coefficients of a
Riemannian metric in the Lie algebroid $(E,[,],\rho)$. Then
\begin{equation}
L(x,y)=g_{ab}(x)y^a y^b\tag{4.5}
\end{equation}
is a regular Lagrangian on $E$. The semispray associated to it is
determined by the functions
\begin{equation}
G^a=\dfrac{1}{2}g^{ab}\left(\dfrac{\partial g_{bc}}{\partial
x^i}\rho^i_d-\dfrac{1}{2}\dfrac{\partial g_{cd}}{\partial
x^i}\rho^i_b-L^e_{db}g_{ec}\right)y^c y^d.\tag{4.6}
\end{equation}

{\bf Example 4.2.} A more general example is provided by the
regular Lagrangians which are homogeneous of degree 2 in $(y^a)$.
By the Euler theorem one obtains
\begin{equation}
L(x,y)=g_{ab}(x,y)y^a y^b,\tag{4.7}
\end{equation}
where $(g_{ab}(x,y))$ are homogeneous functions of degree 0.

As $\dfrac{\partial}{\partial y^a}$ are homogeneous functions of
degree 1 and the derivative with respect to $(x^j)$ does not
affect the degree of homogeneity, it results that the coefficients
$(G^a)$ from (4.4) are homogeneous of degree 2 in $(y^a)$. The
corresponding semispray is nothing but a spray.

\section{Mechanical Lagrangian systems on Lie algebroids}

Let $(E, [,],\rho)$ be a Lie algebroid.

{\bf Definition 5.1.} {\it A mechanical Lagrangian system with
external forces on the Lie algebroid $(E,[,],\rho)$ is
$\sum=(E,L,F)$ with $L$ a regular Lagrangian on $E$ and
$F=(F_a(x,y))$ a vertical covector field.}

Let be the functions
\begin{equation}
{\cal L}_a:=\dfrac{d}{dt}\left(\dfrac{\partial L}{\partial
y^a}\right)-\rho^i_a\dfrac{\partial L}{\partial
x^j}-L^c_{ba}y^b\dfrac{\partial L}{\partial y^c}\tag{5.1}
\end{equation}
defined on admissible curves on $E$.

Then the equalities ${\cal L}_a=0$ represent the Euler - Lagrange
equations associated to $L$.

We assume that the evolution equations of the system $\sum$ are as
follows:
\begin{equation}
{\cal L}_a(x(t),y(t))=F_a(x(t),y(t)),\tag{5.2}
\end{equation}
for $\widetilde{c}(t)=(x(t),y(t))$ an admissible curve on $E$.

The equations (5.2) after some arrangements take the form
\begin{equation}
\dfrac{dy^a}{dt}+2G^a(x,y)=\dfrac{1}{2}F^a(x,y),\tag{5.3}
\end{equation}
where the functions $(G^a)$ are given by (4.4), $F^a=g^{ab}F_b$,
and the equations $\dfrac{dx^i}{dt}=\rho^i_a(x)y^a$ hold.

Thus the evolution equations of the system $\sum$ become
\begin{equation}
\begin{array}{l}\dfrac{dx^i}{dt}=\rho^i_a(x)y^a,\vspace{1.5mm}\\
\dfrac{dy^a}{dt}=-2\left(G^a-\dfrac{1}{4}F^a\right).\end{array}\tag{5.4}
\end{equation}

The solutions of this system may be regarded as the integral
curves of a semispray
\begin{equation}
S^*=\rho^i_a(x)y^a\dfrac{\partial}{\partial
x^i}-2G^*(x,y)\dfrac{\partial}{\partial y^a},\
G^{*a}=G^a-\dfrac{1}{4}F^a.\tag{5.5}
\end{equation}

Indeed, $S^*$ is a semispray because it differs by the semispray
$S$ derived from $L$ by a vertical vector field.

{\bf Definition 5.2.} {\it We say that the mechanical Lagrangian
system $\sum$ is dissipative if $F_a(x,y)y^a\leq 0$ and that it is
strictly dissipative if $F_a(x,y)y^a\leq -\alpha y_a y^a$ with
$\alpha>0$ a constant and $y_a=g_{ab}y^b$.}

{\bf Theorem 5.1.} {\it Let be the mechanical Lagrangian system
$\sum$ with the evolution equations $(5.4)$. If it is dissipative
then its energy $E=y^a\dfrac{\partial L}{\partial y^a}-L$
decreases on the curves that are solutions of $(5.4)$. If
furthermore it is strictly dissipative its energy is strictly
decreasing on the curves solutions of $(5.4)$, assuming that these
have no singularities.}

\noindent{\it Proof.} Let be $\widetilde{c}(t)=(x^i(t),y^a(t))$ a curve
that is a solution of (5.4). Along this curve we have
$$\dfrac{dE}{dt}=\dfrac{dy^a}{dt}\dfrac{\partial L}{\partial
y^a}+y^a\dfrac{d}{dt}\left(\dfrac{\partial L}{\partial
y^a}\right)-\dfrac{\partial L}{\partial
x^i}\dfrac{dx^i}{dt}-\dfrac{\partial L}{\partial
y^a}\dfrac{dy^a}{dt}=$$ $$=y^a{\cal L}_a(x,y)=y^a F_a(x,y).$$

The last equality is based on (5.2) and to obtain the previous one
the equations
\begin{equation}
L^c_{ab}y^a y^b=0,\tag{5.6}
\end{equation}
have been used.

If the system $\sum$ is dissipative we have $\dfrac{dE}{dt}\leq 0$
and if it is strictly dissipative we have
$\dfrac{dE}{dt}\leq-\alpha y_a y^a<0$, q.e.d.

Now, we show that if $\sum$ is dissipative we can associate to it
a Lyapunov function.

Let $(x^i_0,y^a_0)$ be an equilibrium point of $S^*$.

If $\rho$ is injective this has the form $(x^i_0,0)$ with
$G^{*a}(x^i,0)=0$, a condition that is verified if $S^*$ is a
spray.

Assume that $(x^i_0, y^a_0)$ is a minimum point for the energy $E$
and set $\widetilde{E}(x,y)=E(x,y)-E(x_0,y_0)$.

We have
\begin{equation}
\widetilde{E}(x_0,y_0)=0,\ \widetilde{E}(x,y)>0.\tag{5.7}
\end{equation}

Let us denote by ${\cal L}_{S^*}$ the Lie derivative with respect
to $S^*$.

We have: ${\cal L}_{S^*}(E)=\rho^i_a y^a\dfrac{\partial
E}{\partial x^i}-2G^a\dfrac{\partial E}{\partial
y^a}+\dfrac{1}{2}F^a\dfrac{\partial E}{\partial y^a}$.

Expanding this and using again (5.6) we get
\begin{equation}
{\cal L}_{S^*}(E)=y_a F^a\leq 0,\tag{5.8}
\end{equation}
since $\sum$ is dissipative.

Thus the function $\widetilde{E}$ is a Lyapunov function for $S^*$
in the equilibrium point $(x^i_0,y^a_0)$ but we can not conclude
that this point is stable.

In order to do so we need to introduce a Riemannian metric on $E$
and to prove that $S^*$ is complete with respect to that metric.
For details see [7].

\newpage

\def\v{\varphi}
\def\g{\gamma}
\def\a{\alpha}
\def\b{\beta}
\def\l{\lambda}
\def\d{\delta}
\def\p{\partial}
\def\o{\omega}
\def\t{\theta}
\def\vp{\varepsilon}
\def\s{\sigma}
\def\D{\Delta}
\def\k{\kappa}
\def\G{\Gamma}
\def\O{\Omega}
\def\L{\Lambda}
\def\rr{\mathbb{R}}
\def\spa{\vspace*{-.3cm}\\ }
\def\cal#1{\mathcal{#1}}
\def\mb#1{\mathbf{#1}}
\def\bd#1{\boldsymbol{#1}}
\def\dd{\displaystyle}

\runningauthor={M. ANASTASIEI}
\runningtitle={BANACH LIE ALGEBROIDS}
\noindent
\baselineskip 8pt
\noindent{\footnotesize{Analele {\c S}tiin{\c t}ifice ale Universit{\u a}{\c t}ii ``Al.~I.~Cuza'' Ia{\c s}i}}
\hfill\break
{\footnotesize{Mat. (S.N.), Tom. LVII, 2011, f. 2}}
\vskip 2cm
\baselineskip 11.5pt plus .15pt
\centerline{\bf\Large BANACH LIE ALGEBROIDS}
\vskip .5cm
\centerline{\bf {\footnotesize{BY}}}
\vskip .5cm
\centerline{\bf {\footnotesize{M. ANASTASIEI}}}
\vskip 1cm

\begin{abstract} First, we extend the notion of second order differential equations
(SODE) on a smooth manifold to anchored Banach vector bundles.
Then we define the Banach Lie algebroids as Lie algebroids
structures modeled on anchored Banach vector bundles and prove
that they form a category.

{\bf Mathematics Subject Classification 2000:} 58B20, 58A99.

{\bf Key words:} Banach vector bundles, anchor, Second Order Differential Equations, Lie algebroids.
\end{abstract}

\setcounter{section}{0}
\section*{Introduction}

Lie algebroids are related to many areas of geometry ([2], [7]) and has recently become
 an object of extensive studies. See [6] for basic
definitions, ex\-am\-ples and  references. In 1996, Weinstein [8]
proposed some applications of the Lie algebroids in Analytical
Mechanics. New theoretical developments followed. See the survey
[5] by de Leon, Marrero and Martinez about Mechanics on
Lie algebroids.

In [1], we gave a construction of a semispray associated to a
regular Lagrangian on a Lie algebroid.

In this paper, we consider the notion of Lie algebroid in the
category of Banach vector bundles, that is vector bundles over
smooth Banach ma\-ni\-folds whose type fibres are Banach spaces. Such
a Banach vector bundle over base $M$ is called anchored if there
exists a morphism from it to the tangent bundle $TM$. First, we
extend the usual notion of second order differential equations
(SODE) to anchored Banach vector bundles and we show that if a
Banach vector bundle admits a homogeneous SODE it is necessarily
anchored. Then we define the Banach Lie algebroids as Lie
algebroid structures modeled on anchored Banach vector bundles.
In our setting only one from three equivalent definitions of a
morphism of Lie algebroids is working. Using it we show that the
Banach Lie algebroids form a category.

\section{Anchored Banach vector bundles}
\setcounter{equation}{0}

Let $M$ be a smooth i.e. $C^{\infty}$ Banach manifold modeled on a
Banach space $\mathbf{M}$ and let $\pi:E\to M$ be a Banach vector
bundle whose type fiber is a Banach space $\mathbf{E}$. We denote
by $\tau:TM\to M$ the tangent bundle of $M$.
\begin{definition}
We say that the vector bundle $\pi:E\to M$ is an anchored vector
bundle if there exists a vector bundle  morphism  $\rho:E\to TM$.
The morphism $\rho$ will be called the anchor map.
\end{definition}

Let $\cal{F}(M)$ be the ring of smooth real functions on $M$.
We denote  by $\G(E)$ the $\cal{F}(M)$-module of smooth sections
in the vector bundle $(E,\pi,M)$ and by $\cal{X}(M)$ the module of
smooth sections in the tangent bundle of $M$ (vector fields on
$M$).

The vector bundle morphism $\rho$ induces an $\cal{F}(M)$-module
morphism which will be denoted also by $\rho:\G(E)\to \cal{X}(M)$,
$\rho(s)(x) = \rho(s(x)),$ $x\in M, s\in\G(E)$.

Let $\{(U,\v), (V,\psi), \ldots\}$ be an atlas on $M$. Restricting
$U,V$ if necessary we may choose a vector bundle atlas
$\{(\pi^{-1}(U),\overline{\v}), (\pi^{-1}(V), \overline{\psi}),
\ldots \}$ with $\overline{\v}:\pi^{-1}(U)\to U\times \mathbb{E}$
given by $\overline{\v}(u) = (\pi(u), \overline{\v}_{\pi(u)}),$
where $\overline{\v}_{\pi(u)}:E_{\pi(u)}\to \mathbb{E}$ is a
toplinear isomorphism. Here $E_{\pi(u)}$ is the fiber of
$(E,\pi,M)$ in $u\in E.$ The given atlas on $M$ together with a
vector bundle atlas induce a smooth atlas $\{(\pi^{-1}(U),\phi),
(\pi^{-1}(U),\psi),\ldots \}$ on $E$ such that $E$ becomes a
Banach manifold modeled on the Banach space $\mathbb{M}\times
\mathbb{E}$. The map $\phi:\pi^{-1}(U)\to \v(U)\times \mathbb{E}$
is given by
$$
\phi(u) =
(\v(\pi(u)),\overline{\v}_{\pi(u)}(u)),\;\; u \in E.
$$
For a section $s:U \to \pi^{-1}(U)$, its local representation
$\phi\circ s\circ \v^{-1}:\v(U)\to \v(U)\times \mathbb{E}$ given
by $(\phi\circ s\circ \v^{-1})(\v(x)) = (\v \pi(s(x)),
\overline{\v}_{\pi(s(x))}(s(x))=(\v(x), \overline{\v}_{x}(s(x)))$
is completely determined by the map $s_{\v}:\v(U)\to \mathbb{E}$
given by $s_{\v}(\v(x)) = \overline{\v}_{x}(s(x))$ which will be
called the local representative (shortly l.r.) of $s$. On $U\cap
V$ we may speak also of the l.r. $s_{\psi}$ of a section $s:U\cup
V\to \pi^{-1}(U\cap V)$ given by $s_{\psi}(\psi(x)) =
\overline{\psi}_{x}(s(x))$. It is clear that we have
\begin{equation}
s_{\psi}(\psi(x)) = \overline{\psi}_{x}\circ
\overline{\phi}_{x}^{-1}(s_{\v}(\v(x))),\;\; x\in U\cup V.
\end{equation}
For a vector field $X:U \to \tau^{-1}(U)$ we have a l.r.
$X_{\v}:\v(U)\to \mathbb{M}$ and on $U\cap V$ we have also a l.r.
$X_{\psi}$ and one holds
\begin{equation}
X_{\psi}(\psi(x)) = d(\psi\circ
\v^{-1})(\v(x))(X_{\v}(\v(x))),\;\; x\in U\cap V,
\end{equation}
where $d$ means Frechet differentiation.

Locally, $\rho$ reduces to a morphism $U\times \mathbb{E}\to
U\times \mathbb{M}$, $(x,v)\to (x,\rho_{U}(x)v)$ with
$\rho_{U}(x)\in L(\mathbb{E},\mathbb{M})$, the space of continuous
linear maps from $\mathbb{E}$ to $\mathbb{M}$. We call
$\rho_{U}(x)$ the l.r. of $\rho$. On overlaps of local charts one
easily gets
\begin{equation}
\rho_{V}(x)\circ \overline{\psi}_{x}\circ \overline{\v}_{x}^{-1} =
d(\psi\circ \v^{-1})(\v(x))\circ \rho_{U}(x),\;\; x\in U\cap V
\end{equation}

\medskip

{\bf Example.}

1. The tangent bundle of $M$ is trivially anchored vector
bundle with $\rho = I$ (identity).

2. Let $A$ be a tensor field of type $(1,1)$ on $M$. It is
regarded as a section of the bundle of linear mappings
$L(TM,TM)\to M$ and also as a morphism $A :TM\to TM$. In  other
words, $A$ may be thought as an anchor map.

3. Any subbundle of $TM$ is an anchored vector bundle with
the anchor the inclusion map in $TM$.

4. Let $\pi:E \to M$ be only a submersion. The subspaces
$V_uE = \pi ^{-1}(x),\pi(u) =x$ of $TE$ over $E$ denoted by $VE$
form a subbundle  called the vertical subbundle. By Example 3)
this is an anchored Banach vector bundle.

The anchored vector bundles over the same base $M$ form a
category. The objects are the pairs $(E,\rho_{E})$ with $\rho_{E}$
the anchor of $E$ and a morphism $f: (E,\rho_{E})\to (F,
\rho_{F})$ is a vector bundle morphism $f:E\to F$ which verifies
the condition $\rho_{F}\circ f = \rho_{E}$.

\section{Semisprays in an anchored vector bundle}
\setcounter{equation}{0}

Let $(E,\pi, M)$ be an anchored vector bundle with the anchor map $\rho$ and let $\pi_{*} :TE\to TM$ be the differential (tangent map) of $\pi$.

We denote by $\tau_{E}:TE \to E$ the tangent bundle of $E$.

\begin{definition}
A section $S:E\to TE$ will be called a semispray if

\begin{itemize}
\item[$(i)$] $\tau_{E}\circ S =$ identity on $E$,
\item[$(ii)$] $\pi_{\ast}\circ S = \rho.$
\end{itemize}
\end{definition}

The condition (i) says that $S$ is a vector field on $E$. The
condition (ii) can be written also in the form
$$
\pi_{*,u}(S(u)) =\rho(u) = (\rho \circ \tau_{E})(S(u)),\;\; u \in
E.$$
When $E=TM$ and $\rho=$ identity on $TM$, $S$ is simultaneously
a vector field on $TM$ and a section in the vector bundle
$\pi_{*}: TTM\to TM$ i.e. it is a second-order vector field on $M$
in terminology from [3, p.96]. Such a vector field is frequently
called a second order differential equation (SODE) on $M$ or a
semispray.

As we will see below, in our context $S$ is no more related to a
second order differential equation on $M$ and so the corresponding
terminology is inadequate.

Let $c:J \to E$ for $\circ \in J \subset \mathbb{R}$ a curve on
$E$. The differential of $c$ is $c_{*}:J \times \mathbb{R}\to TE$
and using $\imath:J \to J \times \mathbb{R},$ $t\to (t,1),$ $t\in
J$ we set $c^{\prime}(t) =c_{*}\circ \imath.$
Then in general $\pi\circ c$ is a curve on $M$ and we have that
$(\pi\circ c)^{\prime}(t) = \pi_{*,c(t)}\circ c^{\prime}(t).$

\begin{definition}
A curve $c$ on $E$ will be called admissible if $(\pi\circ
c)^{\prime}(t) = \rho (c(t)),$ $\forall t\in J.$
\end{definition}

Locally, if $c:J \to \v(U)\times E,$ $t\to (x(t), w(t))$ then
$\pi\circ c:J \to \v(U)$ is $t\to x(t)$, $t\in J$ and it follows
that $c$ is an admissible curve if and only if
\begin{equation}
\frac{dx}{dt} = \rho_{U}(x(t))w(t),\;\; t\in J
\end{equation}

\begin{theorem}
A vector field $S$ on $E$ is a semispray if and only if all its
integral curves are admissible curves.
\end{theorem}

\noindent{\it Proof.}
Let $S$ be a semispray. A curve $c:J \to E$ is an integral curve
of $S$ if $c^{\prime}(t) = S(c(t))$. It follows $\pi_{*}\circ
c^{\prime}(t) = (\pi_{*}\circ S)(c(t))$ or $(\pi \circ
c)^{\prime}(t) = \rho(c(t))$, that is $c$ is an admissible curve.
Conversely, let $S$ be a vector field on $E$ whose integral curves
are admissible. For every $u\in E$ there exists an unique integral
curve $c:J \to E$ of $S$ such that $c(0) = u$ and $c^{\prime}(0) =
S(u)$. We have $\pi_{*}\circ c^{\prime}(0) = (\pi_{*}\circ S)(u)$,
$(\pi \circ c)^{\prime}(0) = (\pi_{*}\circ S)(u)$ and
$\pi_{*}\circ S = \rho(u)$ since $c$ is admissible.$\hfill\square$

We restrict to a local chart $(U,\v)$ on $M.$ Then $TU
\simeq \v(U)\times \mathbb{M}$, $E_{|U} \simeq \v(U)\times
\mathbb{E}$ and $TE_{|U} \simeq (\v (U)\times \mathbb{E})\times
\mathbb{M} \times \mathbb{E}$.

The l.r. of a vector field on $E$ is $S_{\v}:\v(U)\times
\mathbb{E} \to \v(U)\times \mathbb{E} \times \mathbb{M}\times
\mathbb{E}$, $S_{\v}(x,u) = (x,u,S_{\v}^{1}(x,u),
S_{\v}^{2}(x,u))$. As l.r. of $\pi_{*}$ is $\v(U)\times
\mathbb{E}\times \mathbb{M}\times \mathbb{E}\to \v(U)\times
\mathbb{M}$, $(x,u,y,v)\to (x,y)$ the condition $\pi_{*}\circ S =
\rho$ translates to $S_{\v}^{1}(x,u) = (x,\rho_{U}(x)u)$. We set
for convenience $S_{\v}^{2}(x,u) = -2 G_{\v}(x,u)$ and so the l.r.
of a semispray for the anchored vector bundle $(E,\pi, M)$ with
the anchor $\rho$ is given as follows:
\begin{equation}
S_{\v}(x,u) = (x,u,\rho_{U}(x)u, - 2G_{\v}(x,u)).
\end{equation}

Let $(V,\psi)$ be another local chart and let us set $h = \psi
\circ \v^{-1}: \v(U\cap V)\to \psi(U\cap V)$. Then $h_{*}:\v
(U\cup V)\times \mathbb{M}\to \psi(U\cap V)\times \mathbb{M}$ is
given by $(x,v)\to (x,dh(x)(v))$, $x\in \v(U\cup V), v\in
\mathbb{M}$.

Let us denote by $H:\v (U\cap V)\times \mathbb{E}\to \psi(U\cap
V)\times \mathbb{E}$ the map given by $H(x,u) = (h(x), M(x)u)$,
where $M(x) = \overline{\psi}_{x}\circ \overline{\v}_{x}^{-1}\in
L(\mathbb{E},\mathbb{E})$. Then $H_{*}$ is locally given as the
pair $(H,H^{\prime})$: $\v(U\cup V)\times \mathbb{E}\times
\mathbb{M}\times \mathbb{E}\to \psi(U\cup V)\times
\mathbb{E}\times \mathbb{M}\times \mathbb{E}$, where the
derivative $H^{\prime}(x,u)$ is given by the Jacobian matrix
operating on the column vector $^{t}(y,w)$ with $y\in \mathbb{M}$
and $w\in \mathbb{E}$. Thus $(H,H^{\prime})$ takes the form
$(x,u,y,v)\to (h(x), M(x)u, h^{\prime}(x)y, M^{\prime}(x)(y)(u) +
M(x)v)$ with prime being denoted the Frechet derivative.

If $S_{\psi}$ is l.r. of $S$ in the chart $(V,\psi)$, necessarily
we have $(H,H^{\prime})\circ S_{\v} = S_{\psi}$ with
$S_{\psi}(x,u) = (h(x), M(x)u, \rho_{U}(h(x))M(x)u, -
2G_{\psi}(h(x), M(x)u))$.

Computing $(H,H^{\prime})\circ S_{\v}$ and identifying with
$S_{\psi}$ one finds
\begin{eqnarray}
\rho_{V}(h(x))M(x)(u) &=& h^{\prime}(x)\rho_{U}(x)(u)\nonumber\\
G_{\psi}(h(x), M(x)u)&=& M(x)G_{\v}(x,u) -
\frac{1}{2}M^{\prime}(x)(\rho_{U}(x)u)u.
\end{eqnarray}
The first equation (2.3) is just (1.3) and the second provides the
connection between the l.r. $G_{\v}$ and $G_{\psi}$ on overlaps.
We have

\begin{theorem}
A vector field $S$ on $E$ is a semispray if and only if it has
l.r. $S_{\v}$ in the form $(2.2)$ and the functions involved in
$(2.2)$ satisfy $(2.3)$ on overlaps of local charts.
\end{theorem}

\noindent{\it Proof.}
The ``if'' part was proved in the above. The converse is
obvious.$\hfill\square$

We denote by $h_{\l}:E \to E$, $h_{\l}(u_{x}) =\l u_{x}$, $\l \in
\mathbb{R}$, $\l >0, x\in M,$ the homothety of factor $\l.$

\begin{definition}
We say that a semispray $S$ is a spray if the following equality
holds
\begin{equation}
S\circ h_{\l} = \l (h_{\l})_{*}\circ S.
\end{equation}
Locally, (2.4) is equivalent to
\begin{equation}
G_{\v}(x, \l v) = \l^{2}G_{\v}(x,v),\;\; (x,v)\in U \times
\mathbb{E}.
\end{equation}
\end{definition}
\noindent Indeed, $(S\circ h_{\l})(u) = S(\l u) = (x, \l v, \rho_{U}(\l v),
- 2G_{\v}(x,\l v)$ and $\l(h_{\l})_{*}S(u) =$ $(x, \l v, \l
\rho_{U}(v), - 2\l^{2}G_{\v}(x,\l v))$. Since $\rho_{U}$ is a
linear mapping, (2.4) implies (2.5) and conversely.
We look at (2.5). If we fix $x\in U$ and omit the index $\v$ we
get a mapping $G:\mathbb{E}\to \mathbb{E}$ that verifies $G(\l v)
= \l^{r}G(v)$ for all $\l>0$ and $r=2.$ We say that such a map is
positively homogeneous of degree $r$.

For such mapping the following Euler type theorem holds.

\begin{theorem}
Suppose that a mapping $G: \mathbb{E}\to \mathbb{E}$ is
differentiable away from the origin of $\mathbb{E}$. Then the
following two statements are equivalent:\vspace*{-.3cm}
\begin{itemize}
\item[$(i)$] $G$ is positively homogeneous of degree $r$,\vspace*{-.3cm}

\item[$(ii)$] $dG_{v} (v) = r G(v)$, for all $v\in
\mathbb{E}\backslash \{0\}$.
\end{itemize}
\end{theorem}

\noindent{\it Proof.}
Suppose (i) holds. Fix $v\in \mathbb{E}$ and differentiate the
equation $G(\l y) = \l^{r}G(v)$ with respect to the parameter
$\l.$ We get $dG_{\l v} (\l v) = r\l^{r-1}G(v)$ and for $\l=1,$
$dG_{v}(v) = rG(v),$ that is (ii) holds.

Conversely, suppose (ii), fix $v$ and consider the mapping $\l \to
G(\l v)$ with $\l>0.$ By the chain rule, we have $\dfrac{dG(\l v)}{d\l} = dG_{\l v}(v) = \dfrac{1}{\l}dG_{\l v}(\l v) =
\dfrac{r}{\l}G(\l v)$, that is the mapping $\l\to G(\l v)$ is a
solution of the differential equation
$
\frac{d}{d\l}G(\l v) - \frac{r}{\l}G(\l v) = 0.
$
The integrating factor $\frac{1}{\l^{r}}$ then gives $G(\l v) =
\l^{r}C,$ where the integrating constant $C$ is depending on our fixed $v$. Setting $\l=1,$
we get $C = G(v)$ and so $G(\l v) = \l^{r}G(v),$ that is (i)
holds, q.e.d.$\hfill\square$

The proof of Theorem 2.6 shows also that if $G:\mathbb{E}\to
\mathbb{E}$ is of class $C^{1}$ on $\mathbb{E}$ and positively
homogeneous of degree 1, then it is linear and $G(v) = dG_{v}(v)$.
Moreover, if $G$ is $C^{2}$ on $\mathbb{E}$ and is positively
homogeneous of degree 2, then it is quadratic, that is $2G(v) =
d_{v}^{2}G(v,v).$

Returning to the (2.5) we note that if $G_{\v}$ is of class
$C^{2}$ in the points $(x,0)$, then it is quadratic in $v$.
Thus $S$ satisfying (2.4) reduces to a quadratic spray. In order to avoid
this reduction we have to delete from $E$ the image of
the null section in the vector bundle $\pi: E \to M.$

Now, we show that if for a vector bundle $E\to M$ there exists a
vector field $S_{0}$ on $E$ that satisfies (2.4) then $\pi:E \to
M$ is an anchored vector bundle and $S_{0}$ is a spray.

Let be $S_{0}(x,v)=(x,v,S_{01}(x,v)$, $S_{0,2}(x,v))$ in a local chart on $E$. Then $S_{0}(h_{\l}u)=$ $S_{0}(x,\l v)=$ $(x,\l v,S_{01}(x,\l v)$, $S_{02}(x,\l v))$ and $(h_{\l})_{*}S_{0}(u)=$ $(x,\l v$, $S_{01}(x,v)$, $\l S_{02}(x,v))$. The condition (2.4) implies $S_{01}(x,\l v)=$ $\l S_{01}(x,v)$ and $S_{02}(x,\l v)=\l^{2}S_{02}(x,v)$. It follows that $S_{01}$ is a linear map with
respect to $v$. Hence we may put $S_{01}(x,v) = \rho_{U}(x)v,$
$\rho_{U}(x)\in L(\mathbb{E}, \mathbb{M})$. Using $\{\rho_{U}(x),
x\in M \}$ one defines a morphism $\rho:E \to TM.$ Thus $E\to M$
is an anchored vector bundle. As $(\pi_{*}\circ S_{0})(u) = (x,
S_{01}(x,v)) = (x, \rho_{U}(x)v)$ we have $\pi_{*}\circ S_{0}
=\rho$ and as $\tau_{E}\circ S_{0}$=indentity automatically holds
it follows that $S_{0}$ is a spray.

\section{Category of Banach Lie algebroids}
\setcounter{equation}{0}\setcounter{definition}{0}

Let $\pi:E \to M$ be an anchored Banach vector bundle with the
anchor $\rho_{E}:E \to TM$ and the induced morphism
$\rho_{E}:\G(E)\to \mathcal{X}(M)$.

Assume there exists defined a bracket $[,]_{E}$ on the space
$\G(E)$ that provides a structure of real Lie algebra on $\G(E)$.

\begin{definition}
The triplet $(E, \rho_{E}, [,]_{E})$ is called a Banach Lie
algebroid if
\begin{itemize}
\item[$(i)$] $\rho: (\G(E), [,]_{E})\to (\mathcal{X}(M), [,])$ is a
Lie algebra homomorphism and
\item[$(ii)$] $[s_{1}, fs_{2}]_{E} = f[s_{1}, s_{2}]_{E}
+\rho_{E}(s_{1})(f)s_{2}$, for every $f\in \mathcal{F}(M)$ and
$s_{1}, s_{2}\in \G(E)$.
\end{itemize}
\end{definition}

{\bf Example.}

1. The tangent bundle $\tau: TM\to M$ is a Banach Lie
algebroid with the anchor the identity map and the usual Lie
bracket of vector fields on $M$.

2. For any submersion $\pi:E \to M$, the vertical bundle
$VE$ over $E$ is an anchored Banach vector bundle. As the Lie
bracket of two vertical vector fields is again a vertical vector
field it follows that $(VE, i, [,]_{VE})$, where $i:VE \to TE$ is
the inclusion map, is a Banach Lie algebroid. This applies, in
particular, to any Banach vector bundle $\pi:E \to M.$


Let $\Omega^{q}(E):= \G(\L^{q}E^{*})$ be the $\mathcal{F}(M)-$
module of differential forms of degree $q$. In particular,
$\O^{q}(TM)$ will be denoted by $\O^{q}(M)$. The differential
operator $d_{E}:\O^{q}(E)\to \O^{q+1}(E)$ is given by the formula
\begin{eqnarray}
&&(d_{E} \o)(s_{0}, \ldots, s_{q})  =
\sum_{i=0,\ldots,q}(-1)^{i}\rho_{E}(s_{i})\o(s_{0}, \ldots,
\widehat{s}_{i}, \ldots, s_{q})\nonumber\\&&+ \sum_{0\leq i< j
\leq q} (-1)^{i+j}\o ([s_{i},s_{j}]_{E}, s_{0}, \ldots
\widehat{s}_{i}, \ldots , \widehat{s}_{j}, \ldots, s_{q})
\end{eqnarray}
for $s_{1}, \ldots, s_{q}\in \G(E)$, where hat over a symbol means
that symbol must be deleted.

For Lie algebroids constructed on vector bundles with finite
dimensional fibres there exist three different but equivalent
notions of morphisms.

For Banach Lie algebroids only one of them is working. We give it
here. For a detailed discussion on Lie algebroids morphisms see
[4]. Let $(E^{\prime}, \pi^{\prime}, M)$ be a Banach vector bundle
and $(E^{\prime}, \rho_{E^{\prime}}, [,]_{E^{\prime}})$ a Banach
Lie algebroids based on it.

\begin{definition}
A vector bundle morphism $f:E \to E^{\prime}$ over $f_{0}:M\to
M^{\prime}$ is a morphism of the Banach Lie algebroids $(E,
\rho_{E}, [,]_{E})$ and $(E^{\prime}, \rho_{E^{\prime}},
[,]_{E^{\prime}})$ if the map induced on forms
$f^{*}:\O^{q}(E^{\prime})\to \O^{q}(E)$ defined by
$(f^{*}\o^{\prime})_{x}(s_{1}, \ldots, s_{q}) =
\o^{\prime}_{f_{0}(x)}(fs_{1}, \ldots, fs_{q})$, $s_{1}, \ldots,
s_{2}\in \G(E)$ commutes with the differential i.e.
\begin{equation}
d_E \circ f^{*} = f^{*}\circ d_{E^{\prime}}.
\end{equation}
\end{definition}

Using this definition it is easy to prove

\begin{theorem}
The Banach Lie algebroids with the morphisms defined in the above,
form a category.
\end{theorem}

\vspace{.2cm} {\footnotesize{\it\noindent Received: 15.X.2009}}


\begin{thebibliography}{99}
\item P. Flaschel and W. Klingenberg (1972), {\it Riemannsche Hilbertmannigfaltigkeiten. Periodische Geod\"{a}tische}, Lecture Notes in Mathematics, 282.
\item S. Kobayashi (1965), {\it A Theorem on the Affine group of a Riemannian manifold}, Nagoya Math. Jour., 9, 30-41.
\item S. Lang (1962), {\it Introduction to differentiable manifolds}, Interscience NY.
\item L. Maxim (1972), {\it Connection compatibles with Fredholm structures on Banach manifolds}, An. St. Univ. Iasi, Mat., 18(2), 389-399.
\item J.P. Penot (1967), {\it De submersion en fibrations}, Seminaire de g\'{e}om\'{e}trie diff\'{e}rentielle de Melle P. Libermann, Paris.
\end{thebibliography}

\begin{thebibliography}{99}
  \bibitem{1} Bourbaki N., {\it Vari\'{e}t\'{e}s diff\'{e}rentialles et analytiques}, Fascicule de r\'{e}sultats, Hermann,
Paris, 1967.
  \bibitem{2} Elworthy K.D. and Tromba A.J., {\it Differential structures and Fredholm maps on Banach manifolds}, Proc. Symp. in Pure Math. t. XV., A.M.S., (1970),  p. 45-94.
  \bibitem{3} de la Harpe P., {\it The Clifford algebra and the spinor group of a Hilbert space}, Compositio Math., 25 (1972), p. 245—261.
  \bibitem{4} de la Harpe P., {\it Sur le groupe spinoriel associ\'{e} \`{a} un espace de Hilbert}, C.R. Acad. Sci. Paris, t. 274 (1972), p. 1484-1488.
  \bibitem{5} Hirzebruch F., {\it Topological methods in algebraic geometry},  Springer,  1966.
  \bibitem{6} Koschorke U., {\it Infinite dimensional $K-$theory}, Proc. Symp. in Pure Math., t. XV, A.M.S., (1970), p. 95-133.
  \bibitem{7} Lichnerowicz A., {\it Champs spinoriels et propagateurs en relativit\'{e} g\'{e}n\'{e}rale}, Bull. Soc. Math. France, 92 (1964), p. 11-100.
  \bibitem{8} Milnor J., {\it Remarks concerning spin manifolds}, Differential and combinatorial topology, In honor of M. Morse, Princeton Univ. Press (1965), p. 55-62.
  \bibitem{9} Plymen R.J. and Streater R. F., {\it A  model of the universal covering group of $SO(H)_2$}, Bull. London Math. Soc, 7 (1975), p. 283-289.
  \bibitem{10} Popovici I., {\it Consid\'{e}ration sur les structures spinorielles}, Rend Circ. Mat. Palermo, t. 23 (1974), p. 113-134.
  \bibitem{11} Serre J.-P., {\it Homologie singuli\`{e}re des espaces fibr\'{e}s}, Ann. of Math., 54 (1951), p. 456.
\end{thebibliography}

\begin{thebibliography}{99}

\bibitem{1} R. Hermann, {\it Differential Geometry and the Calculus of Variations}, Academic-Press, New York, (1968).

\bibitem{2} P. Libermann and Ch.-M. Marie, {\it Symplectic Geometry and Analytical Mechanics}, D. Reidel, Dordrecht, (1987).
\bibitem{3} R. Miron, {\it Lagrange geometry}, (this issue).
\bibitem{4} R. Miron and M. Anastasiei, {\it Vector Bundles, Lagrange Spaces, Applications to Relativity}, (in Romanian), Ed. Academiei Romane, Bucuresti, (1987).
\bibitem{5} R. Miron and M. Anastasiei, {\it Lagrange Spaces, Theory and Applications}, Kluwer Academic Publisher, Dordrecht, (1994).
\bibitem{6} R. Abraham and J. Marsden, {\it Foundations of Mechanics}, Benjamin Inc., New York, (1967).
\bibitem{7} P.L. Antonelli, Editor, {\it Mathematical Essays on Growth and the Emergence of Form}, The Univ. of Alberta Press, (1985).
\bibitem{8} P.L. Antonelli, {\it Finsler Volterra-Hamilton systems in ecology}, Tensor N.S. 50, 22-31 (1991).
\bibitem{9} R. Miron, {\it Techniques of Finsler geometry in the theory of vector bundles}, Acta Sci. Math. 49, 119-129 (1985).
\bibitem{10} M. Anastasiei and H. Kawaguchi, {\it A geometrical theory of time-dependent Lagrangians}, Tensor N.S. 48, 273-293 (1989); Tensor N.S. 49, 296-304 (1990).
\bibitem{11} M. de Leon and P.R. Rodrigues, {\it Methods of Differential Geometry in Analytical Mechanics}, North-Holland, Amsterdam, (1989).
\bibitem{12} M. Matsumoto, {\it Foundations of Finsler Geometry and Special Finsler Spaces}, Kaiseisha Press, Otsu, Japan, (1986).
\bibitem{13} S. Kobayashi and K. Nomizu, {\it Foundations of Differential Geometry}, Interscience, New York, (1963).
\bibitem{14} S. Sasaki, {\it Almost contact manifolds}, Tohoku Institute, 1965-1968.
\bibitem{15} N. Oshima, {\it Survey of rheonomic aspects in problems of dynamics}, RAAG Memoirs I (B-II), 62-69 (1955).
\bibitem{16} V.I. Arnold, {\it Mathematical Methods of Classical Mechanics}, Springer-Verlag, New York, (1978).
\bibitem{17} G.S. Asanov, {\it Finsler Geometry, Relativity and Gauge Theories}, D. Reidel, Dordrecht, (1985).

\end{thebibliography}

\begin{thebibliography}{99}
\bibitem{1}
Chern, S.S., {\it Local equivalence and Euclidean connection in Finsler spaces}, Sci. Rep. Nat. Tsing Hua Univ. Ser. A 5(1948), 95-121.
\bibitem{2}
Chern, S.S., {\it On Finsler geometry}, C.R. Acad. Sci. Paris 314(1)(1992), 757-761.
\bibitem{3}
Bao, D., Chern, S.S., {\it On a notable connection in Finsler geometry}, Houston J. Math. 19(1)(1993), 137-182.
\bibitem{4}
Matsumoto, M., {\it Foundations of Finsler geometry and Special Finsler Spaces}, Kaiseisha Press, Japan, 1986.
\bibitem{5}
Miron, R., Anastasiei, M., {\it The Geometry of Lagrange Spaces: Theory and Application}, Kluwer Academic Publishers, Dordrecht, Holland 1994.
\bibitem{6}
Rund, H., {\it The Differential Geometry of Finsler Spaces}, Springer, 1959.
\end{thebibliography}

\begin{thebibliography}{5}

\bibitem{1} Anastasiei M., {\it On Chern's connection in
Finsler geometry}, to appear.

\bibitem{2} Bao, D., Chern, S.S., {\it On a notable connection
in Finsler geometry}, Houston J. Math. 19(1), 1993, 137--182.

\bibitem{3} Matsumoto, M., {\it Foundations of Finsler geometry
and special Finsler spaces}, Kaiseisha Press, Japan, 1986.

\bibitem{5} Miron, R., Anastasiei, M., {\it The geometry of
Lagrange spaces: Theory and Applications.} Kluwer Academic
Publishers, FTPH 59, 1994.

\bibitem{4} Miron, R., Aikou, T., Hashiguchi, M., {\it On
minimality of axiomatic systems of remarkable Finsler
connections,} Rep. Fac. Sci. Kagoshima Univ. (Math. Phys.
\& Chem.) 26, 41--51, 1993.

\end{thebibliography}

\begin{thebibliography}{10}
\bibitem{[Ma1]} M. Anastasiei, {\it A Historical Remark on the Connection of Chern
and Rund.} Contemporary Mathematics, {\bf 196}(1996), 171-176.
\bibitem{[Ma2]} M. Anastasiei, {\it A class of generalized Lagrange spaces.}
An. \c St. Univ. ``Al.I.Cuza'' Ia\c si, Supliment, s I a, Matematic\u a,
{\bf 42}(1996), 259-264.
\bibitem{[BaCh]} D. Bao, S.S. Chern, {\it On a notable connection in Finsler
geometry.} Houston Journal of Mathematics, {\bf 19}, 1(1993), 135-180.
\bibitem{[RM]} R. Miron, M. Anastasiei {\it The Geometry of Lagrange spaces:
Theory and Applications.} Kluwer Academic Publishers, FTPH 59, 1994.
\bibitem{[ON]} B. O'Neil {\it The fundamental equations of a submersion.}
Michigan Math.J., {\bf 13}(1966), 459-469.
\bibitem{[ON]} B. O'Neil {\it Submersion and geodesic.}
Duke Math.J., {\bf 34}(1967), 363-374.
\end{thebibliography}

\begin{thebibliography}{n}
\bibitem{1} Anastasiei, M., {\it A historical remark on the connections
of Chern and Rund.} Contemporary Mathematics, vol.196, 1996,
171--176.

\bibitem{2} Aringazin, A.K., Asanov, G.S., {\it Problems of Finslerian
theory of gauge fields and gravitation.} Reports on Mathematical
Physics 25(1988), 183--241.

\bibitem{3} Aubin, T., {\it M$\acute{e}$triques riemanniennes et courbure.} J.Differential Geometry,
4(1970), 383--424.

\bibitem{4} Beil, R.G., {\it Comparison of unified field
theories.} Tensor N.S., 56(1995), 175--183.

\bibitem{5} Ikeda, S., {\it Advanced Studies in Applied Geometry.}
Seizansha, 1995, Japan.

\bibitem{6} Izumi, H., {\it On the geometry of generalized metric
spaces.} I. {\it Connections and identities}. Publ. Math. Debrecen
39/1-2(1991), 113--134.

\bibitem{7} Matsumoto, M., {\it Foundations of Finsler Geometry and
Special Finsler Spaces,} Kaiseisha Press, Saikawa, \= Otsu, Japan,
1986.

\bibitem{8} Matsumoto, M., Eguchi, K., {\it Finsler spaces admitting
a concurrent vector field}, Tensor, N.S. Vol. 28(1974), 239--249.

\bibitem{9} Miron, M., Anastasiei, M., {\it The Geometry of Lagrange
Spaces: Theory and Applications.} Kluwer, FTPH 59, 1994.

\bibitem{10} Shibata, C., {\it On invariant tensors of
$\beta$--changes of Finsler metrics,}
J. Math. Kyoto Univ. 24--1(1984), 163--188.

\bibitem{11} Shimada, H., {\it Cartan--like connections of special
generalized Finsler spaces. Differential Geometry and Its
Applications.} World Scientific, Singapore, 1990, 270--275.

\bibitem{12} Tavakol, R.K., Van der Bergh, N., {\it Viability
criteria for the theories of gravity and Finsler spaces.} Gen. Rel.
Grav. 18(1986), 849--859.

\bibitem{13} Misner,Ch.W.,Thorne, K.S., Wheeler, A.J., {\it Gravitation.} W.H. Freeman abd Company,
San Francisco, 1980

\bibitem{14} Udri\c ste C., Completeness of Finsler manifolds. Publ. Math. Debrecen, 42/ 1-2(1993), 45--50.
\end{thebibliography}

\begin{thebibliography}{9}
\bibitem{1}  Anastasiei M., Shimada H., {\it Deformations of Finsler metrics}, to
appear.

\bibitem{2}  Antonelli P.L., Hrimiuc D., {\it A new class of Spray-generating
Lagrangians}, in the volume {\bf Lagrange and Finsler geometry. Application to
Physics and Biology}, Eds. P.L. Antonelli and R. Miron, Kluwer Academic
Publishers, TFTPH 76, 1996, 81-92.

\bibitem{3}  Miron R., {\it The homogeneous lift of a Riemannian metric}, to
appear.

\bibitem{4}  Miron R., Anastasiei M., {\it The geometry of Lagrange spaces:
Theory and Applications}, Kluwer Academic Publishers. FTPH 59, 1994.

\bibitem{5}  Musso E., Tricerri F., {\it Riemannian metrics on tangent bundle},
Ann. Mat. Pura Appl.(4), 150, 1988, 1-20.
\end{thebibliography}

\begin{thebibliography}{n}
\bibitem{1} Anastasiei M., Shimada H., {\it The Beil metrics associated to a Finsler space}, to appear in Balkan
Journal of Geometry and Its Applications.

\bibitem{2} Antonelli P.L., Hrimiuc D., {\it A new class of spray generating Lagrangians}, in Kluwer Academic Publishers,
FTPH vol. 76, 1996, 81--92.

\bibitem{3} Aubin T., {\it M$\acute{e}$triques riemanniennes et courbure}, J. Differential Geometry,
4(1970), 383--424.

\bibitem{4} Beil R. G., {\it New Class of Finsler Metrics}, Int. Jour. Theor. Phys. vol.28, 1989, 659--667.

\bibitem{5} Beil R.G., {\it Finsler geometry and a unified field theory}, Contemporary Mathematics, vol. 196, 1996,
265-271.

\bibitem{6} Kitayama M., {\it Generalized Finsler spaces admitting a parallel Finsler vector field}, to appear in
Algebras, Groups and Geometries, 1989, Hadronic Press.

\bibitem{7} Matsumoto M., {\it Foundations of Finsler Geometry and Special Finsler Spaces}, Kaiseisha Press, Saikawa, \= Otsu, Japan, 1986.

\bibitem{8} Miron M., Anastasiei, M., {\it The Geometry of Lagrange Spaces: Theory and Applications}, Kluwer, FTPH 59, 1994.

\bibitem{9} O'Neil B., {\it The fundamental equations of a submersion}, Michigan Math. J. 134(1966), 459--469.

\bibitem{10} Sekizawa M., {\it Curvatures of Tangent Bundles with Cheeger-Gromoll Metric}, Tokyo J. Math. Vol. 14(2),
1991, 407--417.

\bibitem{11} Udri\c ste C., {\it Completeness of Finsler manifolds}, Publ. Math. Debrecen, 42/ 1-2(1993), 45--50.

\end{thebibliography}

\begin{thebibliography}{7}
\bibitem{1} Anastasiei, M., {\it Locally conformal Kaehler structures on tangent manifold of a space form}, Libertas Mathematica, vol.XIX, 71-76.
\bibitem{2} Bao, D., Chern, S.-S, Shen, Z., {\it An introduction to Riemann--Finsler Geometry}, Springer--Verlag, 2000.
\bibitem{3} Goldberg, S.I., Yano, K., {\it On normal globally framed $f$--manifolds}, T\^ohoku Math. J. 22 (1970), 362--370.
\bibitem{4} Hasegawa, I., Yamaguchi, K., Shimada, H., {\it Sasakian structures on Finsler manifolds} in P.L. Antonelli and R. Miron (eds.), Lagrange and Finsler Geometry, Kluwer Academic Publishers, 1996, p. 75--80.
\bibitem{5} Mihai, I., Ro\c sca, R., Verstraelen, L., {\it Some aspects of the differential geometry of vector fields}, {\sc padge}, katholieke Universiteit Leuven, vol.2, 1996.
\bibitem{6} Miron, R., {\it Homogeneous lift of a Finsler space}, to appear
\bibitem{7} Miron, R., Anastasiei, M., {\it The geometry of Lagrange spaces: theory and applications}, Kluwer Academic Publ., FTPH 59, 1994.

\end{thebibliography}

\bigskip


\begin{thebibliography}{n}
\bibitem{1} Anastasiei, M., {\it Locally conformal Kaehler structures on tangent manifold of a space form}, Libertas Matematica, Vol. XIX (1999), 71--77.
\bibitem{2} Anastasiei, M., Shimada, H., {\it Deformations
of Finsler Metrics}, in vol. {\it Finslerian Geometries}, P.L. Antonelli ed.,
Kluwer Acad. Publ. FTPH 109, Dordrecht, 2000, 53 --65.
\bibitem{3} Lewis, D. A.,{\it  A symmetric product for vector fields and its geometric meaning.} To appear in Mathematische Zeitschrift.
\bibitem{4} Mihai, I, Ro\c sca, R., Verstraelen, L. {\it Some aspects of the differential geometry of vector fields}. Centre for Pure and Applied Differential Geometry (PADGE), Katholieke Universiteit Leuven and Brussel, vol. 2, 1996.
\bibitem{5} Miron, R., Anastasiei, M., {\it The geometry of
Lagrange spaces: theory and applications}, Kluwer Acad. Publ. FTPH
59, Dordrecht, 1994.
\bibitem{6} Naveira, A. M., {\it A classification of Riemannian almost product structures}. Rend. Mat. Roma 3(1983), 577--592
\bibitem{7} Wood, C.M., {\it A class of harmonic almost--product
structures}, Journal of Geometry and Physics 14(1994), 25--42.
\end{thebibliography}

\begin{thebibliography}{n}
\bibitem{1}
Matsumoto, M., {\it Foundations of Finsler geometry and special
Finsler space}, Kaiseisha Press, Japan, 1986.
\bibitem{2}
Miron, R., Anastasiei, M., {\it The geometry of Lagrange spaces:
theory and applications}, Kluwer Academic Publishers, FTPH 59, 1994.
\bibitem{3}
Vaisman, Izu, {\it Symplectic curvature tensors}, Mh.Math. 100,
1985, 299--327.

\end{thebibliography}

\begin{thebibliography}{n}
\bibitem{1}
Anastasiei, M., {\it Geometry of Berwald-Cartan spaces}, to appear.
\bibitem{2}
Miron, R., {\it Techniques of Finsler geometry in the theory of vector bundles}, Acta Sci. Math. 49(1985), 119--129.
\bibitem{3}
Miron, R., {\it Hamilton geometry}, Univ. Timi\c soara(Romania), Sem. Mecanica 3(1987), 54 p.
\bibitem{4}
Miron, R., Anastasiei, {\it The Geometry of Lagrange Spaces: Theory and Applications}, Kluwer Academic Publishers, FTPH 59, 1994.
\bibitem{5}
Schmidt, B.G., {\it Conditions on a Connection to be a Metric Connection}, Commun. Math. Phys. 29(1973), 55-59.
\bibitem{6}
Szab\'o, Z.I., {\it Positive Definite Berwald Spaces}, Tensor N.S., Vol. 35(1981), 25--39.
\bibitem{7}
Tamassy, L., {\it Metrizability of Affine Connections}, Balkan J. of Geometry and Its Applications. Vol. 1(1996), no. 1, 83-90
\end{thebibliography}

\begin{thebibliography}{n}
\bibitem{1}
Anastasiei, M., {\it The geometry of Berwald Cartan spaces}, to
appear.
\bibitem{2}
Anastasiei, M., {\it Metrizable linear connections in vector
bundles}, to appear.
\bibitem{3}
Anastasiei, M., Antonelli, P.L., {\it Lagrangian which generate
sprays} in vol. {\it Lagrange and Finsler Geometry. Applications
to Physics and Biology}, edited by P.L. Antonelli and R. Miron.
Kluwer Academic Publishers, FTPH 76, 1996, p. 1--14.
\bibitem{4}
Antonelli, P.L., Hrimiuc, D., {\it A new class of
spray--generating Lagrangians}, idem, p. 81--92.
\bibitem{5}
Bao, D., Chern, S.-S., Shen, Z., {\it An Introduction to
Riemann--Finsler Geometry}, Springer--Verlag New York, Inc., 2000.
\bibitem{6}
Ichijy\_o, Y., {\it Finsler Manifolds with a Linear Connection},
J. Math. Tokushima Univ., vol. 10(1976), 1--11.
\bibitem{7}
Matsumoto, M., {\it Foundations of Finsler Geometry and Special
Finsler Spaces}, Kaiseisha Press, Otsushi, Japan.
\bibitem{8}
Miron, R., {\it Techniques of Finsler geometry in the theory of
vector bundles}, Acta Sci. Math., 49(1985), 119--129.
\bibitem{9}
Miron, R., Anastasiei, M., {\it The Geometry of Lagrange Spaces:
Theory and Applications}, Kluwer Academic Publishers. FTPH 59, 1994.
\bibitem{10}
Szab\'o, Z.I., {\it Positive Definite Berwald Spaces}, Tensor
N.S., vol. 35(1981), 25--39.
\end{thebibliography}

\begin{thebibliography}{n}
\bibitem{BCS}
Bao, D., Chern, S.-S., Shen, Z.,
An Introduction to Riemann--Finsler Geometry,
Graduate Texts in Mathematics 200, Springer--Verlag, 2000.
\bibitem{M}
Matsumoto, M.,
Foundations of Finsler Geometry and Special Finsler Spaces,
Kaiseisha Press, Otsu, 1986.
\bibitem{Mi}
Miron, R., Hamilton geometry. Univ. Timisoara(Romania), Sem. Mecanica 3(1987), 54 p.
\bibitem{MHSS}
Miron, R., Hrimiuc, D., Shimada, H., Sab\u au, V.-S.,
The geometry of Hamilton and Lagrange Spaces,
Kluwer Academic Publishers, FTPH 118, 2001.
\bibitem{MA}
Miron, R., Anastasiei, M.,
the Geometry of Lagrange Spaces: Theory and Applications,
Kluwer Academic Publishers, FTPH 59, 1994.
\bibitem{S}
Szab\'o, Z.I.,
Positive Definite Berwald Spaces,
tensor N.S., Vol. 35(1981), 25--39.
\end{thebibliography}

\begin{thebibliography}{n}
\bibitem{1}
Anastasiei, M., {\it Metrizable linear connections in vector bundles}, Publ. Math. Debrecen 62/3-4(2003),277-287.
\bibitem{2}
Anastasiei, M.,{\it Geometry of Berwald vector bundles}. To appear.
 \bibitem{3}
Miron, R., {\it Techniques of Finsler geometry in the theory of vector bundles}, Acta Sci. Math., 49(1985), 119--129.
\bibitem{4}
Miron, R., Anastasiei, {\it The Geometry of Lagrange Spaces: Theory and Applications}, Kluwer Academic Publishers, FTPH 59, 1994.
\bibitem{5}
Tamassy, L., {\it Metrizability of Affine Connections}, Balkan J. of Geometry and Its Applications. Vol. 1(1996), no. 1, 83-90
\bibitem{6}
Tamassy, L., {\it Geometry of the point Finsler spaces}, lecture presented on July 10, 2002 in Budapest at the J. Bolyai Conference to commemorate the 200th anniversary of his birth.

\end{thebibliography}

\begin{thebibliography}{6}
\bibitem{1} Anastasiei M., {\it Distributions on spray spaces}, BJGA,
6(2001) 1-6.
\bibitem{2} Klein J., {\it Espaces variationnels et m\'{e}canique}, Ann.
Inst. Fourier (Grenoble) 12(1962), 1-124.
\bibitem{3} Martinez E., {\it Lagrangian mechanics on Lie algebroids},
Acta Applicandae Mathematicae, 67(2001), 295-320.
\bibitem{4} Miron R., Anastasiei M., {\it Geometry of Lagrange spaces: theory and
applications}, FTPH 59, Kluwer Academic Publishers, 1994.
\bibitem{5} Weinstein A., {\it Lagrangian Mechanics and Grupoids},
Fields Institute Communications, vol. 7, 1996, p. 207-231.
\end{thebibliography}

\begin{thebibliography}{6}

\item Anastasiei M., {\it Semisprays on Lie algebroids}, to
appear.

\item Miron R. Frigioiu C., {\it Finslerian Mechanical systems}, Algebras, Groups and Geometries, Vol. 22(2), 2005,
151-168.

\item Miron R., Anastasiei M., {\it Geometry of Lagrange spaces: Theory and applications}, FTPH 59, Kluwer Academic Publishers, 1994.

\item Munoz - Lecanda M.C., Yaniz - Fernandes J.F., {\it Dissipative control of mechanical systems: a geometrical
approach}, SIAM J. Control Optim., 40(5), 1505-1516, 2002.

\item Santilli R.M., {\it Foundations of Theoretical Mechanics I}, Text and Monograph in Physics, Springer - Verlag, New-York, 1978

\item Weinstein A., {\it Lagrangian Mechanics and Grupoids}, Fields Institute Communications, vol. 7, 1996, p. 207-231.
\end{thebibliography}

\newpage

\begin{thebibliography}{11}

\bibitem{1} Bao, D.; Chen, S.S.; Shen, Z. -- {\it An Introduction to Riemannian-Finsler
Geometry}, Springer, Grodinate, ... Mathematics, nr. 200, 2000.

\bibitem{2} Galloway, J.G. -- {\it A generalization of Myers theorem and an application relativistic
cosmology,} J. Differential Geometry, 14(1979), 105-116.

\bibitem{3} Matsumoto, M. -- {\it Foundations of Finsler geometry and special Finsler
spaces}, Kaiseisha Preess, 1986.

\bibitem{4} Suceav\u{a}, B. -- {\it A Myers type theorem for almost Hermitian
mnaifolds}, Algebra, Geometry and their Applications, Seminar Proceedins, 2002, p. 1-4.

\end{thebibliography}

\begin{thebibliography}{12}

\bibitem{[1]} Anastasiei M., {\it Geometry of Lagrangians and
semisprays on Lie algebroids}. Proceedings of The 5-th Conference
of Balkan Society of Geometers, August 29 - September 2, 2005,
Mangalia, Romania. BSG Proceedings 13, Geometry Balkan Press,
2006, p.10-17
\bibitem{[2]} Anastasiei M., {\it Semisprays on anchored vector
bundles}. Proceedings of The 4th Annual Symposium on Mathematics
applied in Biology \& Biophysics. U.A.S.V.M. Ia\c si, Romania,
p.195-202
\bibitem{[3]} Anastasiei M., {\it Mechanical systems on Lie
algebroids}. To appear.
\bibitem{[4]} Anastasiei M., {\it Distributions on spray spaces},
BJGA, 6(2001) 1-6.
\bibitem{[5]} Martinez E., {\it Lagrangian mechanics on Lie
algebroids}, Acta Applicandae Mathematicae, 67(2001), 295-320.
\bibitem{[6]} Miron R.; Anastasiei, M., {\it The geometry of
Lagrange spaces: theory and applications}, Kluwer Academic
Publishers, FTPH 59 (1994).
\bibitem{[7]} Munoz - Lecanda M.C., Yaniz - Fernandes J.F., {\it
Dissipative control of mechanical systems: a geometrical
approach}, SIAM J. Control Optim., 40(5), 1505-1516, 2002.
\bibitem{[8]} Reckziegel H., {\it Generalized spray and the
theorem of Ambrose-Palais- Singer}, Dillen, Franki (ed.) et al.,
Geometry and topology of submanifolds, V. Proceedings of the
conferences on ``Differential geometry and vision'' and ``Theory of
submanifolds'' held in combination at Leuven and Brussels, Belgium,
July 1992. Singapore: World Scientific. 242-248 (1993).
\bibitem{[9]} Weinstein A., {\it Lagrangian Mechanics and
Grupoids}, Fields Institute Communications, vol. 7, 1996, p.
207-231.

\end{thebibliography}

\begin{thebibliography}{11}
\bibitem{1} Anastasiei, M., {\it Geometry of Lagrangians and semisprays on Lie algebroids}, Proceedings of the 5th Conference of Balkan Society of Geometers, 10--17, BSG Proc., 13, Geom. Balkan Press, Bucharest, 2006.
\bibitem{2} Dru\c{t}\u{a}, S.L., {\it The holomorphic sectional curvature of general natural K\"{a}hler structures on cotangent bundles}, An. \c{S}tiin\c{t}. Univ. ``Al. I. Cuza'' Ia\c{s}i. Mat. (N.S.), 56 (2010), 113--130.
\bibitem{3} Lang, S., {\it Fundamentals of Differential Geometry}, Graduate Texts in Mathematics, 191, Springer-Verlag, New York, 1999.
\bibitem{4} Higgins, P.J.; Mackenzie, K., {\it Algebraic constructions in the category of Lie algebroids}, J. Algebra, 129 (1990), 194--230.
\bibitem{5} de Le\'{o}n, M.; Marrero, J.C.; Mart\'{i}nez, E., {\it Lagrangian submanifolds and dynamics on Lie algebroids}, J. Phys. A, 38 (2005), R241--R308.
\bibitem{6} Mackenzie, Kirill C.H., {\it General Theory of Lie Groupoids and Lie Algebroids}, London Mathematical Society Lecture Note Series, 213, Cambridge University Press, Cambridge, 2005.
\bibitem{7} Neagu, M., {\it Jet geometrical objects depending on a relativistic time}, An. \c{S}tiin\c{t}. Univ. ``Al. I. Cuza'' Ia\c{s}i. Mat. (N.S.), 56 (2010), 407--428.
\bibitem{8} Weinstein, A., {\it Lagrangian mechanics and groupoids}, Mechanics day (Waterloo, ON, 1992), 207--231, Fields Inst. Commun., 7, Amer. Math. Soc., Providence, RI, 1996.
\end{thebibliography}
\end{document}